\documentclass{article}
\usepackage{latexsym}
\usepackage{graphics}
\usepackage{amsmath}
\usepackage{amsfonts}
\begin{document}
\bibliographystyle{plain}
\title{The decomposition of global conformal invariants VI: 
The proof of the proposition on local Riemannian invariants.}
\author{Spyros Alexakis\thanks{University of Toronto, alexakis@math.toronto.edu.
\newline
This work has absorbed
the best part of the author's energy over many years. 
This research was partially conducted during 
the period the author served as a Clay Research Fellow, 
an MSRI postdoctoral fellow,
a Clay Liftoff fellow and a Procter Fellow.  
\newline
The author is immensely indebted to Charles
Fefferman for devoting twelve long months to the meticulous
proof-reading of the present paper. He also wishes to express his
gratitude to the Mathematics Department of Princeton University
for its support during his work on this project.}} \date{}
\maketitle
\newtheorem{proposition}{Proposition}
\newtheorem{theorem}{Theorem}
\newcommand{\Sum}{\sum}
\newtheorem{lemma}{Lemma}
\newtheorem{observation}{Observation}
\newtheorem{formulation}{Formulation}
\newtheorem{definition}{Definition}
\newtheorem{conjecture}{Conjecture}
\newtheorem{corollary}{Corollary}
\numberwithin{equation}{section}
\numberwithin{lemma}{section}
\numberwithin{theorem}{section}
\numberwithin{definition}{section}
\numberwithin{proposition}{section}

\begin{abstract}
This is the last in a series of papers where we  prove a conjecture of Deser and Schwimmer
regarding the algebraic structure of ``global conformal
invariants''; these  are defined to
be conformally invariant integrals of geometric scalars.
 The conjecture asserts that the integrand of
any such  integral can be expressed as a linear
combination of a local conformal invariant, a divergence and of
the Chern-Gauss-Bonnet integrand.

 The present paper, jointly with \cite{alexakis4,alexakis5} gives 
a proof of an algebraic Proposition regarding local Riemannian invariants, 
which lies at the heart  of our resolution of the Deser-Schwimmer conjecture.
This algebraic Propositon may be of independent interest, applicable to related problems. 
\end{abstract}
\tableofcontents
\section{Introduction}

This paper is the sixth in the series of papers \cite{alexakis1}--\cite{alexakis6}, where we confirm 
a conjecture of Deser and Schwimmer on the algebraic structure of ``global
conformal invariants''. 
For the reader's convenience, we briefly review the Deser-Schwimmer 
conjecture:\footnote{We refer the reader to the introduction in \cite{alexakis1} 
for a detailed discussion of the notions of 
``Riemannian invariant'', ``local conformal invariant''.}

\begin{definition}
\label{globconfinv} Consider a Riemannian invariant $P(g)$ of weight $-n$
($n$ even). We will say that the integral $\int_{M^n}P(g)dV_g$ is 
a ``global conformal invariant'' if 
 the value of $\int_{M^n} P(g)dV_g$ remains invariant
under conformal re-scalings of the metric $g$.

In other words, $\int_{M^n}P(g)dV_g$ is 
a ``global conformal invariant'' if  for any
 $\phi\in C^\infty (M^n)$ we have 
 $\int_{M^n} P(e^{2\phi}g)dV_{e^{2\phi}g}=\int_{M^n}P(g)dV_g$.
\end{definition}

\par The Deser-Schwimmer conjecture \cite{ds:gccaad} asserts:

\begin{conjecture}
\label{desschwim} Let $P(g)$ be a Riemannian invariant of weight $-n$ 
such that the integral $\int_{M^n}P(g)dV_g$ 
is a global conformal invariant. Then
there exists a local conformal invariant $W(g)$, a Riemannian
vector field $T^i(g)$ and a constant $(Const)$ so that $P(g)$ can
be expressed in the form:
\begin{equation} \label{dseqn} P(g)=W(g)+div_i T^i(g)+(Const)\cdot
\operatorname{Pfaff}(R_{ijkl}).
\end{equation}
\end{conjecture}

We prove: 

\begin{theorem}
\label{thetheorem} Conjecture \ref{desschwim} is true.
\end{theorem}

For the reader's convenience, we recall in brief the main
results obtained in the earlier papers in this series; we then 
(again briefly) discuss the relationship of our work 
(and this paper in particular) with work 
related to Riemannian and conformal invariants, a 
subject largely inspired by Fefferman's program to
understand the singularities in the Bergman and 
Szeg\"o kernels of strictly pseudo-convex CR manifolds. 
Finally, we give a proper statement of the results we will be proving in the 
present paper and an outline of their proof. 
\newline

In \cite{alexakis1,alexakis2,alexakis3}
 we proved that the Deser-Schwimmer conjecture 
holds, {\it provided} one can show certain 
``Main algebraic propositions'', namely Proposition 5.2 in 
\cite{alexakis1} and Propositions 3.1, 3.2  in \cite{alexakis2}. 
The next three papers, \cite{alexakis4}--\cite{alexakis6} 
(including the present one) are devoted to proving these ``Main algebraic Propositons''. 

In \cite{alexakis4} we set up a multiple induction by which 
we will prove these Propositions. In particular, we presented 
the ``fundamemental Proposition'' 2.1 in \cite{alexakis4}
 (which we reproduce here, see Proposition \ref{giade} below) 
which is a generalization of 
the Main algebraic Propositions, and which depends on certain
 parameters (more on this below);  we explained that we would prove Proposition \ref{giade}
 by an induction on these parameters. 
Since Proposition 2.1 is rather complicated to even write out, 
we first recall Proposition 5.2 in \cite{alexakis1}:
\newline

{\it A simplified description of the main algebraic Proposition 5.2 in \cite{alexakis1}:}
Given a Riemannian metric $g$ over an $n$-dimensional 
 manifold $M^n$ and auxilliary $C^\infty$ scalar-valued functions 
$\Omega_1,\dots,\Omega_p$ defined over $M^n$, the 
objects of study are linear combinations of tensor fields 
$\sum_{l\in L} a_l C^{l,i_1\dots i_\alpha}_g$, where each 
$C^{l,i_1\dots i_\alpha}_g$ is a {\it partial contraction} 
with $\alpha$ free indices, in the form: 
\begin{equation}
\label{gen.form1}
 pcontr(\nabla^{(m)}R\otimes\dots\otimes\nabla^{(m_s)}R\otimes
\nabla^{(b_1)}\Omega_1\otimes\dots\otimes\nabla^{(b_m)}\Omega_p);
\end{equation}
here $\nabla^{(m)}R$ stands for the $m^{th}$ covariant derivative of 
the curvature tesnor $R$,\footnote{In other words it is an $(m+4)$-tensor; 
 if we write out its free 
indices it would be in the form 
$\nabla^{(m)}_{r_1\dots r_m}R_{ijkl}$.} and 
$\nabla^{(b)}\Omega_h$ stands for the $b^{th}$ covariant 
derivative of the function $\Omega_h$. A {\it partial contraction} 
means that we have list of pairs of indices $({}_a,{}_b),\dots, ({}_c,{}_d)$ 
in (\ref{gen.form1}), 
which are contracted against each other using the  metric $g^{ij}$. The remaining 
 indices (which are not contracted against another index in (\ref{gen.form1})) are the 
{\it free indices} ${}_{i_1},\dots, {}_{i_\alpha}$. 

The ``main algebraic Proposition'' of \cite{alexakis1} (roughly) asserts the following: Let 
$\sum_{l\in L_\mu} a_l C^{l,i_1\dots i_\mu}_g$ stand for a 
linear combination of partial contractions in the form (\ref{gen.form1}), 
where each $C^{l,i_1\dots i_\mu}_g$ has a given number $\sigma_1$ of factors and a given number 
$p$ of factor $\nabla^{(b)}\Omega_h$. Assume also that $\sigma_1+p\ge 3$, 
each $b_i\ge 2$,\footnote{This means that 
each function $\Omega_h$ is differentiated at least twice.} and that 
for each pair of contracting 
indices $({}_a,{}_b)$ in any given $C^{l,i_1\dots i_\mu}_g$, 
the indices ${}_a,{}_b$ do not belong to the 
same factor. Assume also the rank 
$\mu>0$ is fixed and each partial contraction $C^{l,i_1\dots i_\mu}_g, l\in L_\mu$ 
has a given {\it weight} $-n+\mu$.\footnote{See \cite{alexakis1} 
for a precise definition of weight.} 
Let also $\sum_{l\in L_{>\mu}} a_l C^{l,i_1\dots i_{y_l}}_g$ stand for 
a (formal) linear combination of partial contractions of 
weight $-n+y_l$, with all the properties of the 
terms indexed in $L_\mu$, {\it except} that now all the partial 
contractions have a different rank $y_l$, and each $y_l>\mu$. 

The assumption of the ``main algebraic Proposition'' is 
a local equation in the form:

\begin{equation}
\label{assumption.1} 
\sum_{l\in L_\mu} a_l Xdiv_{i_1}\dots Xdiv_{i_\mu}C^{l,i_1\dots i_\mu}_g+
\sum_{l\in L_{>\mu}} a_l Xdiv_{i_1}\dots Xdiv_{i_{y_l}}C^{l,i_1\dots i_{y_l}}_g=0,
\end{equation}
which
is assumed to hold {\it modulo} complete contractions with $\sigma+1$ factors.
Here given a partial contraction $C^{l,i_1\dots i_\alpha}_g$ 
in the form (\ref{gen.form1}) $Xdiv_{i_s}[C^{l,i_1\dots i_\alpha}_g]$
stands for sum of $\sigma-1$ terms in $div_{i_s}[C^{l,i_1\dots i_\alpha}_g]$
where the derivative $\nabla^{i_s}$ is {\it not} 
allowed to hit the factor to which the free 
index ${}_{i_s}$ belongs.\footnote{Recall that given a partial contraction
$C^{l,i_1\dots i_\alpha}_g$ in the form (\ref{gen.form1}) 
with $\sigma$ factors, $div_{i_s}C^{l,i_1\dots i_\alpha}_g$ 
is a sum of $\sigma$ partial contractions of rank $\alpha-1$. 
The first summand arises by adding a derivative $\nabla^{i_s}$ 
onto the first factor $T_1$ and then contracting the upper index ${}^{i_s}$
against the free index ${}_{i_s}$; the second summand 
arises by adding a derivative $\nabla^{i_s}$ 
onto the second factor $T_2$ and then contracting the upper index ${}^{i_s}$
against the free index ${}_{i_s}$ etc.}  

The main algebraic Proposition in \cite{alexakis1} 
then claims that there will exist a linear combination of partial 
contactions in the form (\ref{gen.form1}), 
$\sum_{h\in H} a_h C^{h,i_1\dots i_{\mu+1}}_g$ with all the properties 
of the terms indexed in $L_{>\mu}$, and all with rank $(\mu+1)$, so that:

\begin{equation}
\label{conclusion.1} 
\sum_{l\in L_\mu} a_l C^{l,(i_1\dots i_\mu)}_g+
\sum_{h\in H} a_h  Xdiv_{i_{\mu+1}}
C^{l,(i_1\dots i_\mu)i_{\mu+1}}_g=0;
\end{equation}
the above holds modulo terms of length $\sigma+1$. Also the symbol $(\dots)$ 
means that we are {\it symmetrizing} over the indices between parentheses. 
\newline

{\bf Local Invariants and Fefferman's program on the Bergman and Szeg\"o kernels.} 
 The theory of {\it local} invariants of Riemannian structures 
(and indeed, of more general geometries,
e.g.~conformal, projective, or CR)  has a long history. As stated above, the original foundations of this 
field were laid in the work of Hermann Weyl and \'Elie Cartan, see \cite{w:cg, cartan}. 
The task of writing out local invariants of a given geometry is intimately connected
with understanding polynomials in a space of tensors with  given symmetries, 
which remain invariant under the action of a Lie group. 
In particular, the problem of writing down all 
 local Riemannian invariants reduces to understanding 
the invariants of the orthogonal group. 

 In more recent times, a major program was laid out by C.~Fefferman in \cite{f:ma}
aimed at finding all scalar local invariants in CR geometry. This was motivated 
by the problem of understanding the  
local invariants which appear in the asymptotic expansions of the 
Bergman and Szeg\"o kernels of strictly pseudo-convex CR manifolds,
 in a similar way to which Riemannian invariants appear in the asymptotic expansion  
of the heat kernel; the study of the local invariants
in the singularities of these kernels led to important breakthroughs 
in \cite{beg:itccg} and more recently by Hirachi in \cite{hirachi1}.
It is worth noting that an analogous problem arises in
the context of understanding the asymptotic expansion of the
Szeg\"o kernel of strictly pseudo-convex domains in $\mathbb{C}^n$
(or alternatively of abstract CR-manifolds). In particular, the
leading term of the logarithmic singularity of the Szeg\"o kernel
exhibits a global invariance which is very similar to the one we
discuss here, see \cite{h:lsskgissd}.

\par This program was later extended  to conformal geometry in \cite{fg:ci}. 
Both these geometries belong to a 
broader class of structures, the
{\it parabolic geometries}; these admit a principal bundle with 
structure group a parabolic subgroup $P$ of a semi-simple 
Lie group $G$, and a Cartan connection on that principle bundle 
(see the introduction in \cite{cg1}). 
An important question in the study of these structures 
is the problem of constructing all their local invariants, which 
can be thought of as the {\it natural, intrinsic} scalars of these structures.

  In the context of conformal geometry, the first (modern) landmark 
in understanding {\it local conformal invariants} was the work of Fefferman 
and Graham in 1985 \cite{fg:ci},
where they introduced the {\it ambient metric}. This allows one to 
construct local conformal invariants of any order in odd 
dimensions, and up to order $\frac{n}{2}$ in even dimensions. 
The question is then whether {\it all} invariants arise via this construction. 

The subsequent work of Bailey-Eastwood-Graham \cite{beg:itccg} proved that 
 indeed in odd dimensions all conformal invariants 
arise via this construction; in even dimensions, 
they proved that the result holds  
when the weight (in absolute value) is bounded by the dimension. The ambient metric construction 
in even dimensions was recently extended by Graham-Hirachi, \cite{grhir}; this enables them to 
indentify in a satisfactory way {\it all} local conformal invariants, 
even when the weight (in absolute value) exceeds the dimension.  

 An alternative 
construction of local conformal invariants can be obtained via the {\it tractor calculus} 
introduced by Bailey-Eastwood-Gover in \cite{bego}. This construction bears a strong 
resemblance to the Cartan conformal connection, and to 
the work of T.Y.~Thomas in 1934, \cite{thomas}. The tractor 
calculus has proven to be very universal; 
tractor bundles have been constructed \cite{cg1} for an entire class of parabolic geometries. 
The relation betweeen the conformal tractor calculus and the Fefferman-Graham 
ambient metric  has been elucidated in \cite{cg2}.  
\newline

{\it Broad discussion on the main algebraic Proposition:} 
The present work, while pertaining to the questions above
(given that it ultimately deals with the algebraic form of local 
{\it Riemannian} and {\it conformal} invariants\footnote{Indeed, the prior work on 
 local {\it conformal} 
invariants played a central role in this endeavor, in \cite{alexakis2,alexakis3}.}),
 nonetheless addresses a different 
{\it type} of problem:  We here consider Riemannian invariants $P(g)$ for 
which the {\it integral} $\int_{M^n}P(g)dV_g$ remains invariant 
under conformal changes of the underlying metric; we then seek to understand 
the possible algebraic form of the {\it integrand} $P(g)$, 
ultimately proving that it can be de-composed 
in the way that Deser and Schwimmer asserted.

 Now, Proposition \ref{giade}\footnote{We recall that this is a generalization of
the ``main algebraic Propositions'' in \cite{alexakis1,alexakis2}, 
which we outlined above.} is purely a statement on local 
Riemannian invariants (in this case intrinsic scalar-valued functions
which depend both on the metric and on certain 
auxilliary functions $\Omega_1,\dots,\Omega_p$).
Thus, this is a proposition regarding the algebraic properties 
of the {\it classical} local Riemannian invariants.
While the author was led to led to the ``main algebraic Propositions'' 
out of the strategy that he felt was necessary to 
solve the Deser-Schwimmer conjecture, they can 
be thought of as results of independent interest. 
The {\it proof} of  Proposition \ref{giade}, presented
 in \cite{alexakis4,alexakis5,alexakis6} is in fact 
not particularily intuitive. It is proven by an induction
(this is perhaps natural); the proof of the inductive step 
relies on studying the conformal variation of 
the assumption of Proposition \ref{giade} below. 
This is perhaps unexpected: Proposition \ref{giade} 
deals {\it purely} with Riemannian invariants; accordingly (\ref{hypothese2}) 
holds for {\it all} Riemannian metrics. From that point of view, 
it is not obvious why restricting attention to the conformal variation of 
the equation (\ref{hypothese2}) should provide useful 
information on the underlying algebraic form of the terms in (\ref{hypothese2}). 
It is the author's 
sincere hope that deeper insight will be obtained 
in the future as to  {\it why} the algebraic 
 Propositions 5.2, 3.1,3.2 in \cite{alexakis1,alexakis2} hold. 
\newline

 Let us now recall Proposition 2.1 in \cite{alexakis4}: 

  This claim  
(reproduced as Proposition \ref{giade} below) 
 deals with tensor fields in the form:

\begin{equation}
\label{form2}
\begin{split}
&pcontr(\nabla^{(m_1)}R_{ijkl}\otimes\dots\otimes\nabla^{(m_{\sigma_1})}
R_{ijkl}\otimes \\&S_{*}\nabla^{(\nu_1)}R_{ijkl}\otimes\dots\otimes
S_{*}\nabla^{(\nu_t)} R_{ijkl}\otimes
\\& \nabla^{(b_1)}\Omega_1\otimes\dots\otimes
\nabla^{(b_p)}\Omega_p\otimes
\\& \nabla\phi_{z_1}\dots \otimes\nabla\phi_{z_w}\otimes\nabla
\phi'_{z_{w+1}}\otimes
\dots\otimes\nabla\phi'_{z_{w+d}}\otimes\dots \otimes
\nabla\tilde{\phi}_{z_{w+d+1}}\otimes\dots\otimes\nabla\tilde{\phi}_{z_{w+d+y}}).
\end{split}
\end{equation}
(See the introduction in \cite{alexakis4} for a detailed 
description of the above form). We recall that a 
(complete or partial) contraction in the above form is called ``acceptable''
$b_i\ge 2$ for every $1\le i\le p$. (In other words, 
we require that each of the functions $\Omega_i$ is differentiated at least twice).
\newline

 The claim of Proposition 2.1 in \cite{alexakis4} which we reproduce here  is a generalization of the 
``main algebraic Propositions'' in \cite{alexakis1,alexakis2}:

\begin{proposition}
\label{giade}
Consider two linear combinations of acceptable tensor fields in
the form (\ref{form2}):

$$\Sum_{l\in L_\mu} a_l
C^{l,i_1\dots i_{\mu}}_{g} (\Omega_1,\dots
,\Omega_p,\phi_1,\dots ,\phi_u),$$

$$\Sum_{l\in L_{>\mu}} a_l
C^{l,i_1\dots i_{\beta_l}}_{g} (\Omega_1,\dots
,\Omega_p,\phi_1,\dots ,\phi_u),$$
 where each tensor field above has real length $\sigma\ge 3$ and a given
simple character $\vec{\kappa}_{simp}$. We assume  that for each
$l\in L_{>\mu}$,  $\beta_l\ge \mu+1$. We also assume 
 that none of the tensor fields of maximal refined double
character in $L_\mu$ are ``forbidden'' (see Definition 2.12 in \cite{alexakis4}).

 We denote by
$$\Sum_{j\in J} a_j C^j_{g}(\Omega_1,\dots ,\Omega_p,
\phi_1,\dots ,\phi_u)$$ a generic linear combination of complete
contractions (not necessarily acceptable) in the form
(\ref{form1}) below, that are simply subsequent to
$\vec{\kappa}_{simp}$.\footnote{Of course if $Def(\vec{\kappa}_{simp})=\emptyset$ then  by
definition $\Sum_{j\in J} \dots=0$.} We assume that:

\begin{equation}
\label{hypothese2}
\begin{split}
&\Sum_{l\in L_\mu} a_l Xdiv_{i_1}\dots Xdiv_{i_{\mu}}
C^{l,i_1\dots i_{\mu}}_{g} (\Omega_1,\dots
,\Omega_p,\phi_1,\dots ,\phi_u)+
\\&\Sum_{l\in L_{>\mu}} a_l Xdiv_{i_1}\dots Xdiv_{i_{\beta_l}}
C^{l,i_1\dots i_{\beta_l}}_{g} (\Omega_1,\dots
,\Omega_p,\phi_1,\dots ,\phi_u)+
\\& \Sum_{j\in J} a_j
C^j_{g}(\Omega_1,\dots ,\Omega_p,\phi_1,\dots ,\phi_u)=0.
\end{split}
\end{equation}

\par We draw our conclusion with a little more notation: We break the index set
$L_\mu$ into subsets $L^z, z\in Z$, ($Z$ is finite)
with the rule that each $L^z$  indexes tensor fields
with the same refined double character, and conversely two tensor
fields with the same refined double character must be indexed in
the same $L^z$. For each index set $L^z$, we denote the
refined double character in question by $\vec{L}^z$. Consider
the subsets $L^z$ that index the tensor fields of {\it maximal} refined double
 character.\footnote{Note that in any set $S$ of $\mu$-refined double characters
 with the same simple character there is going to be a subset $S'$
 consisting of the maximal refined double characters.}
 We assume that the index set of those $z$'s
is $Z_{Max}\subset Z$.

We claim that for each $z\in Z_{Max}$ there is some linear
combination of acceptable $(\mu +1)$-tensor fields,

$$\Sum_{r\in R^z} a_r C^{r,i_1\dots i_{\alpha +1}}_{g}(\Omega_1,
\dots ,\Omega_p,\phi_1,\dots ,\phi_u),$$ where 
 each $C^{r,i_1\dots i_{\mu +1}}_{g}(\Omega_1, \dots
,\Omega_p,\phi_1,\dots ,\phi_u)$ has a $\mu$-double
character $\vec{L^z_1}$ and also the same set of factors 
$S_{*}\nabla^{(\nu)}R_{ijkl}$ as in $\vec{L}^z$ contain 
special free indices, so that:

\begin{equation}
\label{bengreen}
\begin{split}
& \Sum_{l\in L^z} a_l C^{l,i_1\dots i_\mu}_{g}
(\Omega_1,\dots ,\Omega_p,\phi_1,\dots
,\phi_u)\nabla_{i_1}\upsilon\dots\nabla_{i_\mu}\upsilon -
\\&\Sum_{r\in R^z} a_r X div_{i_{\mu +1}}
C^{r,i_1\dots i_{\mu +1}}_{g}(\Omega_1,\dots
,\Omega_p,\phi_1,\dots ,\phi_u)\nabla_{i_1}\upsilon\dots
\nabla_{i_\mu}\upsilon=
\\& \Sum_{t\in T_1} a_t
C^{t,i_1\dots i_\mu}_{g}(\Omega_1,\dots ,\Omega_p,,\phi_1,\dots
,\phi_u)\nabla_{i_1}\upsilon\dots \nabla_{i_\mu}\upsilon,
\end{split}
\end{equation}
modulo complete contractions of length $\ge\sigma +u+\mu +1$.
Here each
$$C^{t,i_1\dots i_\mu}_{g}(\Omega_1,\dots
,\Omega_p,\phi_1,\dots ,\phi_u)$$ is acceptable and is either
simply or doubly subsequent to $\vec{L}^z$.\footnote{Recall that
``simply subsequent'' means that the simple character of
$C^{t,i_1\dots i_\mu}_{g}$ is subsequent to $Simp(\vec{L}^z)$.}
\end{proposition}

 (See the first section in \cite{alexakis4} for a description 
 of the notions of {\it real length, acceptable tensor fields, 
simple character, 
refined double character, maximal refined double character, 
simply subsequent, strongly doubly subsequent}. 

 Proposition \ref{giade} is proven 
by an induction on four parameters, which we now recall:

{\it The induction:}  Denote
the left hand side of equation (\ref{hypothese2}) by
\\$L_{g}(\Omega_1,\dots ,\Omega_p,\phi_1,\dots ,\phi_u)$ or just
$L_{g}$ for short. We recall that for the complete contractions in $L_{g}$,
$\sigma_1$  stands for the number of factors
$\nabla^{(m)}R_{ijkl}$ and $\sigma_2$ stands for the number of
factors $S_{*}\nabla^{(\nu)} R_{ijkl}$. Also $\Phi$ stands for
the total number of factors
$\nabla\phi,\nabla\tilde{\phi},\nabla\phi'$ and $-n$ stands
for the weight of the complete contractions involved.
\newline

\begin{enumerate}

\item{We assume that Proposition \ref{giade} is true for all
linear combinations $L_{g^{n'}}$ with weight $-n'$,
$n'<n$, $n'$ even, that satisfy the
hypotheses of our Proposition.}

\item{We assume that Proposition \ref{giade} is true for all
linear combinations $L_{g}$ of weight $-n$ and real length
$\sigma'<\sigma$, that satisfy the
hypotheses of our Proposition.}

\item{We assume that Proposition \ref{giade} is true for all
linear combinations $L_{g}$ of weight $-n$ and real length
$\sigma$, with $\Phi'>\Phi$ factors $\nabla\phi,\nabla\tilde{\phi},\nabla\phi'$, 
that satisfy the hypotheses of our Proposition.}

\item{We assume that Proposition \ref{giade} is true for all
linear combinations $L_{g}$ of weight $-n$ and real length
$\sigma$, $\Phi$ factors $\nabla\phi,\nabla\tilde{\phi},\nabla\phi'$
 and with {\it fewer than $\sigma_1+\sigma_2$ curvature
factors} $\nabla^{(m)}R_{ijkl},S_{*}\nabla^{(\nu)}R_{ijkl}$,
provided $L_g$  satisfies the 
hypotheses of our Proposition.}
\end{enumerate}

\par We will then prove Proposition \ref{giade} for the linear
 combinations $L_{g}$ with weight $-n$, real length $\sigma$, $\Phi$  factors
 $\nabla\phi,\nabla\phi',\nabla\tilde{\phi}$ and with
$\sigma_1+\sigma_2$ curvature factors $\nabla^{(m)}R_{ijkl}$,
$S_{*}\nabla^{(\nu)}R_{ijkl}$. So we are proving our Proposition
by a multiple induction on the parameters $n,\sigma, \Phi,\sigma_1+\sigma_2$
of the linear combination $L_{g}$.
\newline

In \cite{alexakis4} we reduced the inductive step of
 Proposition \ref{giade} to three 
Lemmas 3.1, 3.2, 3.5;\footnote{Lemma 3.5 in \cite{alexakis4} depends on 
two {\it preparatory Lemmas}, 3.3, 3.4 in \cite{alexakis4}.} 
(in particular we distinguished cases I,II,III on 
Proposition \ref{giade} by examiniming 
the tensor fields appearing in (\ref{hypothese2}) and these three Lemmas corresponded to 
the three cases). In \cite{alexakis4} and \cite{alexakis5} we proved that 
 these three Lemmas  in \cite{alexakis4}, imply the inductive 
step of Proposition \ref{giade}. 

In the present paper we prove Lemmas 3.1, 3.2, 3,3, 3,4, 3.5. Lemmas 3.1, 3.2 
in \cite{alexakis4}. Lemmas 3.1, 3.2  will be derived in part A of 
the present paper, which consists of sections \ref{dasproof}, \ref{theanalysis}, \ref{easyproof}.
These two Lemmas are simpler to prove than Lemma 3.5; the analysis performed 
in part A will lay the groundwork for the proof of 
Lemma 3.5  in \cite{alexakis4} (and Lemmas 3.3, 3.4 in \cite{alexakis4}), 
in part B of the present paper, which consists of all the remaining sections.
\newline

For the reader's convenience, we will reproduce 
here the statements of Lemmas 3.1, 3.2, 3.5  
from \cite{alexakis4}, which will be 
proven in the present paper.\footnote{These Lemmas 
are reproduced as Lemmas \ref{zetajones}, \ref{pool2}, \ref{pskovb} 
in the present paper.} There will be separate discussions 
in the beginning of parts A and B outlining the ideas
and arguments that come into in the proofs of these Lemmas.  
 For the reader's convenience,
 however, we will first provide a very
 simple sketch of the claims of these two Lemmas. We do this in order 
to present the {\it gist} of their claims,  freed from the many notational 
conventions needed for the precise statement:
\newline

{\bf A simplified formulation of Lemmas 3.1, 3.2, 3.5 from \cite{alexakis4}:}
The assumption of our Lemma is the equation (\ref{hypothese2}). We
recall that all the tensor fields appearing in that equation 
have the same {\it $u$-simple character}, which
(in simple language) means that the factors $\nabla\phi_h$, 
$1\le h\le u$ contract against 
the different factors $\nabla^{(m)}R_{ijkl}, S_{*}\nabla^{(\nu)}R_{jkl}, \nabla^{(A)}\Omega_h$ according 
to the {\it same pattern}; for example, if the factor $\nabla\phi_1$ contracts against the index ${}_i$ of a
 factor $S_{*}\nabla^{(\nu)}_{r_1\dots r_\nu}R_{ir_{\nu+1}kl}$ and the factor $\nabla\phi_4$ 
contracts against one of the indices ${}_{r_1},\dots, {}_{r_\nu},{}_{r_{\nu+1}}$ 
for a tensor field $C^{l_1, i_1\dots i_a}_g$ in (\ref{hypothese2}), then 
the factors $\nabla\phi_1,\nabla\phi_4$ contract according to 
that rule in {\it all} the tensor fields in (\ref{hypothese2}). 

{\it The notion of $Xdiv$:} For a tensor field 
(i.e. a {\it partial contraction}) $C^{l,i_1\dots i_a}_g$ 
in the form (\ref{form2}), 
given a free index ${}_{i_s}$ which belongs to a factor $T$, the regular
 divergence $div_{i_s}C^{l,i_1\dots i_a}_g$  equals a sum of 
$\sigma+u$ $(a-1)$-tensor fields: The sum arises when we hit any of 
the $\sigma+u$ factors in $C^{l,i_1\dots i_a}_g$ 
by  a derivative $\nabla^{i_s}$,\footnote{(which contracts against the free index ${}_{i_s}$).} 
and then sum over all the resulting $(\mu-1)$-tensor fields. Now, $Xdiv_{i_s}C^{l,i_1\dots i_a}_g$ 
stands for the sum of $\sigma-1$ terms in the sum  $div_{i_s}C^{l,i_1\dots i_a}_g$ where 
we {\it only} consider the $\sigma-1$ terms where the derivative $\nabla^{i_s}$ has hit a factor in 
one of the forms $\nabla^{(m)}R_{ijkl},S_{*}\nabla^{(\nu)}R_{ijkl}$ or $\nabla^{(a)}\Omega_h$,\footnote{In 
other words $\nabla^{i_s}$ is not allowed to hit one of the 
$u$ factors $\nabla\phi_h$, $1\le h\le \sigma$.} but {\it not} the factor $T$
 to which the free index ${}_{i_s}$ belongs.
Thus, given any tensor field $C^{l,i_1\dots i_\alpha}_g(\Omega_1,\dots, \Omega_p,\phi_1,\dots,\phi_u)$, we can think of 
$Xdiv_{i_1}\dots Xdiv_{i_\alpha}C^{l,i_1\dots i_\alpha}_g(\Omega_1,\dots, \Omega_p,\phi_1,\dots,\phi_u)$ 
as a linear combination of complete contractions in the form: 
\begin{equation}
\label{form1}
\begin{split}
&pcontr(\nabla^{(m_1)}R_{ijkl}\otimes\dots\otimes\nabla^{(m_s)}R_{ijkl}
\otimes
\\& \nabla^{(b_1)}\Omega_1\otimes\dots\otimes \nabla^{(b_p)}\Omega_p
\otimes\nabla\phi_1\otimes\dots \otimes\nabla\phi_u);
\end{split}
\end{equation}

{\bf The (simplified)  statement of Lemma \ref{zetajones}:} This Lemma 
applies when there are tensor fields of rank $\mu$
in (\ref{hypothese2}) with special free indices in factors 
$S_{*}\nabla^{(\nu)}R_{ijkl}$.\footnote{Recall that a  free index 
in a factor $S_{*}\nabla^{(\nu)}R_{ijkl}$ is {\it special} when it is one of the indices ${}_k,{}_l$.} 
In that case, our Lemma picks out a particular subset of the tensor 
fields of rank $\mu$, all of which have a special free index 
in  a factor $S_{*}\nabla^{(\nu)}R_{ijkl}$ (denote the index 
set of these tensor fields by $L_\mu^*\subset L_\mu$); it also picks out
 one of those special free indices--say the index ${}_{i_1}$, which will occupy the position 
 ${}_k$ of the factor $S_{*}\nabla^{(\nu)}R_{ijkl}\nabla^i\tilde{\phi}_1$. 
 
 We then consider the $(\mu-1)$-tensor fields 
$C^{l,i_1\dots i_{\mu}}_g\nabla_{i_1}\phi_{u+1}$, $l\in L_\mu^*$.\footnote{These 
 $(\mu-1)$-tensor fields arise from $C^{l,i_1\dots i_\mu}_g$ by just 
 contracting the free index ${}_{i_1}$ against a new factor $\nabla\phi_{u+1}$.}
 The claim of Lemma \ref{zetajones} is (schematicaly) 
that there exists a linear combination of $(\mu+1)$-tensor 
 fields, $\sum_{h\in H} a_h C^{h,i_1\dots i_{\mu+1}}_g\nabla_{i_1}\phi_{u+1}$,  
where each tensor field $C^{h,i_1\dots i_{\mu+1}}_g$ is a partial contraction in the form (\ref{form2}), 
with the same $u$-simple character $\vec{\kappa}_{simp}$, and with the index ${}_{i_1}$
  occupying the position ${}_k$ in the factor $S_{*}\nabla^{(\nu)}R_{ijkl}$, such that:

\begin{equation}
\label{schemclaimzetajones} 
\begin{split}
&\Sum_{l\in L_\mu^*} a_l 
Xdiv_{i_2}\dots Xdiv_{i_\alpha}C^{l,i_1\dots i_{\alpha}}_{g} (\Omega_1,\dots
,\Omega_p,\phi_1,\dots ,\phi_u)\nabla_{i_1}\phi_{u+1}=
\\&\sum_{h\in H} a_h Xdiv_{i_2}\dots Xdiv_{i_{\alpha+1}}C^{h,i_1\dots i_{\alpha+1}}_g\nabla_{i_1}\phi_{u+1}
\\&+\sum_{j\in J} a_j C^{j,i_1}_g (\Omega_1,\dots
,\Omega_p,\phi_1,\dots ,\phi_u)\nabla_{i_1}\phi_{u+1};
\end{split}
\end{equation}
here the terms indexed in $J$ are ``junk terms'': they have length $\sigma+u$ (like the 
tensor felds indexed in $L_1$ and $H$)
and are in the general form (\ref{form1}). They are ``junk terms'' because 
they have one of the two following features: either the 
index ${}_{i_1}$ (which contracts against the factor $\nabla\phi_{u+1}$) 
belongs to some factor $S_{*}\nabla^{(\nu)}R_{ijkl}$ but is {\it not} a special index, {\it or}
one of the factors $\nabla\phi_h, 1\le h\le u$) which are supposed to contract against the
 index ${}_i$ in some factor $S_*\nabla^{(\nu)}R_{ijkl}$ for $\vec{\kappa}_{simp}$
  now contracts against a derivative index of some factor $\nabla^{(m)}R_{ijkl}$.\footnote{In 
  the formal language of Lemma \ref{zetajones}, introduced in \cite{alexakis4}, in this second scenario
  we would say that $C^{j,i_1}_g(\Omega_1,\dots
,\Omega_p,\phi_1,\dots ,\phi_u)$ is ``simply subequent' to
 the simple character $\vec{\kappa}_{simp}$. }
 \newline

{\bf The (simplified)  statement of Lemma \ref{pool2}:} 
This Lemma applies when no  tensor fields of rank $\mu$
in (\ref{hypothese2}) have special free indices in factors 
$S_{*}\nabla^{(\nu)}R_{ijkl}$,\footnote{Recall that a  free index 
in a factor $S_{*}\nabla^{(\nu)}R_{ijkl}$ is {\it special} when it is one of the indices ${}_k,{}_l$.} 
but there are tensor fields of rank $\alpha$ that have special free indices 
in $\nabla^{(m)}R_{ijkl}$.\footnote{Recall that a  free index 
in a factor $\nabla^{(m)}R_{ijkl}$ is {\it special} when it is one of the indices ${}_i,{}_j,{}_k,{}_l$.} 
In that case, our Lemma picks out a particular subset of the tensor 
fields of rank $\mu$, all of which have a special free index 
in  a factor $\nabla^{(m)}R_{ijkl}$ (denote the index set 
of these tensor fields by $L_\mu^*\subset L_\mu$); it also picks out
 one of those special free indices--say the index ${}_{i_1}$, which will occupy the position 
 ${}_i$ of the factor $\nabla^{(m)}R_{ijkl}$. 
 
 We then consider the $(\mu-1)$-tensor fields 
$C^{l,i_1\dots i_{\mu}}_g\nabla_{i_1}\phi_{u+1}$, $l\in L_\mu^*$.\footnote{These 
 $(\mu-1)$-tensor fields arise from $C^{l,i_1\dots i_\mu}_g$ by just 
 contracting the free index ${}_{i_1}$ against a new factor $\nabla\phi_{u+1}$.}
 The claim of Lemma \ref{zetajones} is (schematicaly) that there will exist a linear combination of $(\mu+1)$-tensor 
 fields, $\sum_{h\in H} a_h C^{h,i_1\dots i_{\mu+1}}_g\nabla_{i_1}\phi_{u+1}$,  
where each tensor field $C^{h,i_1\dots i_{\mu+1}}_g$ is a partial contraction in the form (\ref{form2}), 
with the same $u$-simple character $\vec{\kappa}_{simp}$, and with the index ${}_{i_1}$
  occupying the position ${}_i$ in the factor $\nabla^{(\nu)}R_{ijkl}$, such that:

\begin{equation}
\label{schemclaimzetajones} 
\begin{split}
&\Sum_{l\in L^*_\mu} a_l 
Xdiv_{i_2}\dots Xdiv_{i_\alpha}C^{l,i_1\dots i_{\alpha}}_{g} (\Omega_1,\dots
,\Omega_p,\phi_1,\dots ,\phi_u)\nabla_{i_1}\phi_{u+1}=
\\&\sum_{h\in H} a_h Xdiv_{i_2}\dots Xdiv_{i_{\alpha+1}}C^{h,i_1\dots i_{\alpha+1}}_g\nabla_{i_1}\phi_{u+1}
\\&+\sum_{j\in J} a_j C^{j,i_1}_g (\Omega_1,\dots
,\Omega_p,\phi_1,\dots ,\phi_u)\nabla_{i_1}\phi_{u+1};
\end{split}
\end{equation}
here the terms indexed in $J$ are ``junk terms'': they have length $\sigma+u$ (like the 
tensor felds indexed in $L_1$ and $H$)
and are in the general form (\ref{form1}). They are ``junk terms'' because 
they have one of the two following features: either the 
index ${}_{i_1}$ (which contracts against the factor $\nabla\phi_{u+1}$) 
is a derivative index in some factor $\nabla^{(m)}R_{ijkl}$, {\it or}
one of the factors $\nabla\phi_h, 1\le h\le u$ which are supposed to contract against the
 index ${}_i$ in some factor $S_*\nabla^{(\nu)}R_{ijkl}$ for $\vec{\kappa}_{simp}$ 
  now contracts against a derivative index of some factor $\nabla^{(m)}R_{ijkl}$.\footnote{In 
  the formal language of Lemma \ref{zetajones}, 
introduced in \cite{alexakis4}, in this second scenario
  we would say that $C^{j,i_1}_g(\Omega_1,\dots
,\Omega_p,\phi_1,\dots ,\phi_u)$ is ``simply subequent'' to
 the simple character $\vec{\kappa}_{simp}$. }
 \newline

{\bf The (simplified)  statement of Lemma \ref{pskovb}:}  Lemma 
\ref{pskovb} applies when all tensor fields of minimum rank $\mu$
in (\ref{hypothese2}) have no special free indices in factors
$S_{*}\nabla^{(\nu)}R_{ijkl}$ or $\nabla^{(m)}R_{ijkl}$.\footnote{Recall that a  free index 
in a factor $S_{*}\nabla^{(\nu)}R_{ijkl}$ is {\it special} when it is one of the indices ${}_k,{}_l$;
 a free index in a factor $\nabla^{(m)}R_{ijkl}$ is {\it special} 
when it is one of the indices ${}_i,{}_j,{}_k,{}_l$.} 
 In order to distinguish cases A and B of Lemma \ref{pskovb}, we must recall some 
 facts about the notion of {\it refined double character} 
 of $\mu$-tensor fields in the form (\ref{form2}) 
 and {\it the maximal refined double character} among the 
 $\mu$-tensor fields appearing  in (\ref{hypothese2}):\footnote{The reader
  is refered to \cite{alexakis4} for precise definitions of these notions.}
\newline

{\bf The notion of (refined) double character, and the comparison 
between different refined double characters:}
We recall that for a tensor field $C^{l,i_1\dots i_\mu}_g$
in the form (\ref{form2}) with no special free indices, its refined double character (which coincides with 
the {\it double character} in this case) encodes the pattern of
 distribution of the $\mu$ free indices among the different factors.  
 
\par Furthermore, in \cite{alexakis4} we introduced 
a weak ordering among refined double characters: Given two 
tensor fields $C^{l,i_1\dots i_\mu}_g, C^{r,i_1\dots i_\mu}_g$
 (with the same simple character, say $\vec{\kappa}_{simp}$), 
we have introduced a comparison between their (refined) 
double characaters: We formed a list of the numbers of free 
indices that belong to  the different factors, say $List_l=(s_1,\dots s_\sigma)$ 
and $List_r=(t_1,\dots t_\sigma)$, and considered the decreasing 
rearrangements of these lists, say $RList_l,RList_r$. We then decreed
 $C^{l,i_1\dots i_\mu}_g$ to be ``doubly subsequent''
to $C^{r,i_1\dots i_\mu}_g$ if $RList_r$ is lexicographically greater than $RList_l$. 
We also defined two {\it different} ``refined double characters'' for which 
neither one is doubly subsequent to the other to be {\it equipolent}.

Now, cases A and B for Lemma \ref{pskovb} are distinuished as follows: 
Let us consider the different $\mu$-tensor fields of maximal refined double character in 
(\ref{hypothese2}). Let us suppose that their
corresponding  lists of distributions of free indices (in decreasing 
rearrangement),\footnote{As defined in the previous paragraph.} 
is in the form  $(M,s_1,\dots ,s_{\sigma-1})$. 
Case A is when $s_1\ge 2$. Case B is when $s_1\le 1$. 
\newline

{\it A rough description of Lemma \ref{pskovb} in case A:} 
We canonicaly pick out a particular subset of the $\mu$-tensor 
fields of maximal refined double character in (\ref{hypothese2});
we denote the index set of these tensor fields by 
 $L^z,z\in Z'_{Max}$ ($\mu$-tensor fields with the same refined double character 
 are indexed in the same index set $L^z$ and vice versa).  
 
 For each  $C^{l,i_1\dots i_\mu}_g$, $l\in\bigcup_{z\in Z'_{Max}}L^z$, we denote by
$\dot{C}^{l,i_1\dots \hat{i}_{t\alpha +1}\dots i_\mu,i_{*}}_{g}$
the tensor field that arises from $C^{l,i_1\dots i_\mu}_{g}$
 by erasing a certain particular index ${}_{i_{t\alpha+1}}$ and adding a free derivative
 index ${}_{i_{*}}$ onto a particular other factor(s).\footnote{As noted 
in \cite{alexakis4}, this operation is well-defined.}
 
 The claim of Lemma \ref{pskovb} is (schematicaly) that there 
will exist a linear combination of $(\mu+1)$-tensor 
 fields, $\sum_{h\in H} a_h C^{h,i_1\dots i_{\mu+1}}_g\nabla_{i_1}\phi_{u+1}$,  
where each tensor field $C^{h,i_1\dots i_{\mu+1}}_g$ 
is a partial contraction in the form (\ref{form2}), 
with the same $u$-simple character $\vec{\kappa}_{simp}$, and with the index ${}_{i_1}$
  being a non-special index in the crucial factor, such that:

\begin{equation}
\label{schemclaimpskovb} 
\begin{split}
&\sum_{z\in Z'_{Max}} \sum_{l\in L^z} a_l 
Xdiv_{i_2}\dots Xdiv_{i_\mu}C^{l,i_1\dots i_{\mu}}_{g} (\Omega_1,\dots
,\Omega_p,\phi_1,\dots ,\phi_u)\nabla_{i_1}\phi_{u+1}=
\\&\sum_{l\in \tilde{L}} a_l 
Xdiv_{i_2}\dots Xdiv_{i_\mu}C^{l,i_1\dots i_{\mu}}_{g} (\Omega_1,\dots
,\Omega_p,\phi_1,\dots ,\phi_u)\nabla_{i_1}\phi_{u+1}+
\\&\sum_{h\in H} a_h Xdiv_{i_2}\dots 
Xdiv_{i_{\mu+1}}C^{h,i_1\dots i_{\mu+1}}_g\nabla_{i_1}\phi_{u+1}
\\&+\sum_{j\in J} a_j C^{j,i_1}_g (\Omega_1,\dots
,\Omega_p,\phi_1,\dots ,\phi_u)\nabla_{i_1}\phi_{u+1}.
\end{split}
\end{equation}
Here the $(\mu-1)$-tensor fields indexed in $\tilde{L}$ are acceptable in the 
form (\ref{form2}) and also have length $\sigma+u$
 (like the ones indexed in  each $L^z$), but they are doubly 
 subsequent to the $(\mu-1)$-tensor fields in the first line. 
The terms indexed in $J$ are ``junk terms''; they have length $\sigma+u$ (like the 
tensor felds indexed in $L^z$ and $H$)
and are in the general form (\ref{form1}). 
They are ``junk terms'' because 
one of the factors $\nabla\phi_h, 1\le h\le u$) which are supposed to contract against the
 index ${}_i$ in some factor $S_*\nabla^{(\nu)}R_{ijkl}$  for the $u$-simple character $\vec{\kappa}_{simp}$
  now contracts against a derivative index of some factor 
$\nabla^{(m)}R_{ijkl}$.\footnote{In the formal language of Lemma \ref{pskovb},
 introduced in \cite{alexakis4}, 
  we would say that $C^{j,i_1}_g(\Omega_1,\dots
,\Omega_p,\phi_1,\dots ,\phi_u)$ is {\it simply subsequent} to
 the simple character $\vec{\kappa}_{simp}$.}
\newline

We note that we proved in \cite{alexakis4} how Lemma \ref{pskovb} (in case A) implies the
 inductive step of Proposition \ref{giade}.  

{\it A rough description of Lemma \ref{pskovb} in case B:} In this case 
the claim of Lemma 3.5 in \cite{alexakis4}
coincides with that of Proposition \ref{giade}. 
\newline

{\bf The rigorous statement of Lemmas 3.1, 3.2 in \cite{alexakis4}:}

For both Lemmas 3.1, 3.2 in \cite{alexakis4}, we canonicaly pick out a particular subset of the 
$\mu$-tensor fields of maximal refined double character
 in (\ref{hypothese2});\footnote{We refer the reader to 
the discussion above Proposition 2.1 in 
\cite{alexakis4} for a rigorous definition of this notion.}
we denote the index set of these tensor fields by 
 $L^z,z\in Z'_{Max}$ ($\mu$-tensor fields with the same refined double character 
 are indexed in the same index set $L^z$ and vice versa).  

Then, Lemma 3.1 in \cite{alexakis4} asserts the following:

\begin{lemma}
\label{zetajones}
\par Assume (\ref{hypothese2}), with weight $-n$, real length $\sigma$,
 $u=\Phi$ and $\sigma_1+\sigma_2$ factors
$\nabla^{(m)}R_{ijkl},S_{*}\nabla^{(\nu)}R_{ijkl}$--assume also
that the tensor fields of maximal refined double character are not
 ``forbidden'' (see Definition 2.12 in \cite{alexakis4}).
 Suppose that
there are $\mu$-tensor fields in (\ref{hypothese2}) with at least one 
 special free index in  a factor $S_{*}\nabla^{(\nu)}R_{ijkl}$.
 We then claim that there is a linear combination of
acceptable tensor fields,
$$\Sum_{p\in P} a_p C^{p,i_1\dots i_{b}}_{g}
(\Omega_1,\dots ,\Omega_p,\phi_1,\dots ,\phi_u),$$ each with
$b\ge\mu+1$, with a simple character $\vec{\kappa}_{simp}$ and
where each \\$C^{p,i_1\dots i_b}_{g}(\Omega_1,\dots
,\Omega_p,\phi_1,\dots ,\phi_u)$ has the property that the free
index ${}_{i_1}$ is the index ${}_k$ in the critical factor
$S_{*}\nabla^{(\nu)}R_{ijkl}$ against which
$\nabla\tilde{\phi}_{Min}$ is contracting, so that modulo complete
contractions of length $\ge\sigma +u+2$:

\begin{equation}
\label{narod}
\begin{split}
&\Sum_{z\in Z'_{Max}}\Sum_{l\in L^z} a_l Xdiv_{i_2}\dots
Xdiv_{i_\mu}C^{l,i_1\dots i_\mu}_{g}(\Omega_1,\dots
,\Omega_p,\phi_1,\dots ,\phi_u)\nabla_{i_1}\phi_{u+1}+
\\&\Sum_{\nu\in N} a_\nu Xdiv_{i_2}\dots
Xdiv_{i_\mu}C^{\nu,i_1\dots i_\mu}_{g}(\Omega_1,\dots
,\Omega_p,\phi_1,\dots ,\phi_u)\nabla_{i_1}\phi_{u+1} -
\\&\Sum_{p\in P} a_p Xdiv_{i_2}\dots
Xdiv_{i_{b}}C^{p,i_1\dots i_{b}}_{g} (\Omega_1,\dots
,\Omega_p,\phi_1,\dots ,\phi_u)\nabla_{i_1} \phi_{u+1}=
\\& \Sum_{t\in T} a_t C^{t,i_{*}}_{g}
(\Omega_1,\dots , \Omega_p,\phi_1,\dots
,\phi_u)\nabla_{i_{*}}\phi_{u+1}.
\end{split}
\end{equation}
 Here each $C^{\nu,i_1\dots i_\mu}_{g}(\Omega_1,\dots
,\Omega_p,\phi_1,\dots ,\phi_u)\nabla_{i_1}\phi_{u+1}$ is
acceptable and has a simple character $\vec{\kappa}_{simp}$
 (and ${}_{i_1}$ is again the index ${}_k$ in the critical factor
  $S_{*}\nabla^{(\nu)}R_{ijkl}$), but
also has either strictly fewer than $M$ free indices in the critical factor
or is doubly subsequent to each $\vec{L}^z, z\in Z'_{Max}$.
 Each  $C^{t,i_{*}}_{g} (\Omega_1,\dots , \Omega_p,\phi_1,\dots
,\phi_u)\nabla_{i_{*}}\phi_{u+1}$ is either simply subsequent to
$\vec{\kappa}_{simp}$, {\it or}
 $C^{t,i_{*}}_{g} (\Omega_1,\dots ,
\Omega_p,\phi_1,\dots ,\phi_u)$ has a  $u$-simple character
$\vec{\kappa}_{simp}$ but the index ${}_{i_{*}}$ is not a special index.
 All complete contractions have the same weak $(u+1)$-simple character.
\end{lemma}

In order to state Lemma 3.2 in \cite{alexakis4}, 
we recall that it applies in the case where the tensor fields of
 maximal refined double character in (\ref{hypothese2}) have special free indices in some 
 factors $\nabla^{(m)}R_{ijkl}$, but no special free indices 
in any factor $S_{*}\nabla^{(\nu)}R_{ijkl}$. We recall also that for each
 $C^{l,i_1\dots i_\mu}_g,l\in L^z, z \in Z'_{Max}$ the set $I_{*,l}$ stands
 for the index set of special free indices.

We also recall that for 
 each $l\in L^z, z\in Z'_{Max}$ and each  ${}_{i_h}\in I_{*,l}$ (we may assume
 with no loss of generality that ${}_{i_h}$ is the index ${}_i$
in some factor $\nabla^{(m)}R_{ijkl}$), we denote by
$\tilde{C}^{l,i_1\dots i_\mu}_{g}(\Omega_1,\dots
,\Omega_p,\phi_1,\dots ,\phi_u)\nabla_{i_h}\phi_{u+1}$ the tensor
field that arises from $C^{l,i_1\dots i_\mu}_{g}(\Omega_1,\dots
,\Omega_p,\phi_1,\dots ,\phi_u)\nabla_{i_h}\phi_{u+1}$ by
replacing the expression $\nabla^{(m)}_{r_1\dots r_m}R_{i_hjkl}\nabla^{i_h}\phi_{u+1}$
by an expression $S_{*}\nabla^{(m)}_{r_1\dots r_m}R_{i_hjkl}\nabla^{i_h}\phi_{u+1}$. 
Then, Lemma 3.2 in \cite{alexakis4} asserts:

\begin{lemma}
\label{pool2} Assume (\ref{hypothese2}) with weight $-n$, real length $\sigma$,
 $u=\Phi$ and $\sigma_1+\sigma_2$ factors
$\nabla^{(m)}R_{ijkl},S_{*}\nabla^{(\nu)}R_{ijkl}$.
Suppose that no $\mu$-tensor fields
 have special free indices in factors $S_{*}\nabla^{(\nu)}R_{ijkl}$,
but some have special free indices in factors $\nabla^{(m)}R_{ijkl}$.
In the notation above we claim that
 there exists a linear combination $\Sum_{d\in D} a_d
C^{d,i_1,\dots ,i_b}_{g}(\Omega_1,\dots ,\Omega_p,\phi_1,\dots
,\phi_u,\phi_{u+1})$ of acceptable $b$-tensor fields (in the form (\ref{form2}) 
and $(u+1)$ factors $\nabla\phi$ and length $\sigma+u+1$) with a
$(u+1)$-simple character $\vec{\kappa}'_{simp}$ and $b\ge \mu$, so that:

\begin{equation}
\label{vecFskillb}
\begin{split}
&\Sum_{z\in Z'_{Max}}\Sum_{l\in L^z} a_l \Sum_{i_h\in
I_{*,l}}Xdiv_{i_1}\dots\hat{Xdiv}_{i_h}\dots
 Xdiv_{i_\mu}\tilde{C}^{l,i_1\dots
i_\mu}_{g}(\Omega_1,\dots ,\Omega_p,\phi_1,\dots
,\phi_u)\\&\nabla_{i_h}\phi_{u+1}
+\Sum_{\nu\in N} a_\nu Xdiv_{i_2}\dots Xdiv_{i_{\mu}} 
C^{\nu,i_1\dots i_{\mu}}_{g}(\Omega_1,\dots
,\Omega_p,\phi_1,\dots ,\phi_u)\nabla_{i_1}\phi_{u+1} -
\\& \Sum_{d\in D} a_d Xdiv_{i_1}\dots Xdiv_{i_b} C^{d,i_1,\dots
,i_b}_{g}(\Omega_1,\dots ,\Omega_p,\phi_1,\dots
,\phi_u,\phi_{u+1})=
\\& \Sum_{t\in T} a_t C^{t,i_{*}}_{g}
(\Omega_1,\dots ,\Omega_p,\phi_1,\dots
,\phi_u,\phi_u)\nabla_{i_{*}}\phi_{u+1},
\end{split}
\end{equation}
where the $(\mu-1)$-tensor fields $C^{\nu,i_1\dots
i_{\mu}}_{g}(\Omega_1,\dots ,\Omega_p,\phi_1,\dots
,\phi_u)\nabla_{i_1}\phi_{u+1}$ are acceptable, have
$(u+1)$-simple character $\vec{\kappa}'_{simp}$ but also either
 have fewer than $M$ free indices in the factor against which
  $\nabla_{i_h}\phi_{u+1}$ contracts,\footnote{``Fewer than $M$ free indices''
  where we also count the free index ${}_{i_h}$.} or are
doubly subsequent
 to all the refined double characters $\vec{\kappa}^z$,
$z\in Z'_{Max}$. Moreover we require that each $C^{\nu,i_1\dots
,i_\mu}_{g}$ has the property that at least one of the indices
${}_{r_1},\dots ,{}_{r_\nu},{}_j$ in the factor $S_{*}\nabla^{(\nu)}_{r_1\dots
r_\nu}R_{ijkl}$ is neither free nor contracting against  a factor
$\nabla\phi_h'$, $h\le u$. The complete contractions
$C^{t,i_{*}}_{g} (\Omega_1,\dots ,\Omega_p,\phi_1,\dots
,\phi_u,\phi_u)\nabla_{i_{*}} \phi_{u+1}$ are simply subsequent to
$\vec{\kappa}'_{simp}$.
\end{lemma}

{\bf The rigorous statement of Lemma 3.5 in \cite{alexakis4}:} 

\begin{lemma}
\label{pskovb}
\par Assume (\ref{hypothese2}) with weight $-n$, real length $\sigma$,
 $u=\Phi$ and $\sigma_1+\sigma_2$ factors
$\nabla^{(m)}R_{ijkl},S_{*}\nabla^{(\nu)}R_{ijkl}$, and additionally assume
 that no $\mu$-tensor field in (\ref{hypothese2}) has special free indices;
assume also that $L^{*}_\mu\bigcup L^{+}_\mu\bigcup 
L''_{+}=\emptyset$ (see the statement of Lemma 3.5 in 
\cite{alexakis4} and the 
discussion above it). 
 Recall the case A that we have distinguished above.

\par Consider case A: Let $k$ stand
 for the (universal) number of second critical factors
among the tensor fields indexed in $\bigcup_{z\in Z'_{Max}}L^z$.
Let also $\alpha$ be the number
of free indices in the (each) second critical factor in each $C^{l,i_1\dots i_\mu}_g$, for each 
$z\in Z'_{Max}$. We claim that:

\begin{equation}
\label{esmen}
\begin{split}
&{\alpha\choose{2}}\Sum_{z\in Z'_{Max}}\Sum_{l\in L^z} a_l
\Sum_{r=0}^{k-1}Xdiv_{i_2}\dots Xdiv_{i_{*}}\dot{C}^{l,i_1\dots
\hat{i}_{r\alpha+1}\dots i_\mu,i_{*}}_{g}
(\Omega_1,\dots,\Omega_p,\phi_1,\dots ,\phi_u)
\\&\nabla_{i_{r\alpha+2}}\phi_{u+1}+\Sum_{\nu\in N} a_\nu Xdiv_{i_2}\dots
Xdiv_{i_\mu}C^{\nu,i_1\dots
,i_\mu}_{g}(\Omega_1,\dots,\Omega_p,\phi_1,\dots
,\phi_u)\nabla_{i_1} \phi_{u+1}+
\\& \Sum_{t\in T_1} a_t Xdiv_{i_1}\dots Xdiv_{i_{z_t}}
C^{t, i_1\dots i_{z_t}}_{g}(\Omega_1,\dots ,\Omega_p,\phi_1,\dots
,\phi_{u+1})+
\\& \Sum_{t\in T_2}  a_t Xdiv_{i_2}\dots Xdiv_{i_{z_t}}
C^{t, i_1\dots i_{z_t}}_{g}(\Omega_1,\dots ,\Omega_p,\phi_1,\dots
,\phi_u)\nabla_{i_1}\phi_{u+1}+
\\& \Sum_{t\in T_3}  a_t Xdiv_{i_1}\dots Xdiv_{i_{z_t}}
C^{t, i_1\dots i_{z_t}}_{g}(\Omega_1,\dots ,\Omega_p,\phi_1,\dots
,\phi_{u+1})
\\& \big{(}+\Sum_{t\in T_4}  a_t Xdiv_{i_1}\dots Xdiv_{i_{z_t}}
C^{t, i_1\dots i_{z_t}}_{g}(\Omega_1,\dots ,\Omega_p,\phi_1,\dots
,\phi_{u+1})\big{)}=
\\& \Sum_{j\in J}  a_j C^j_{g}(\Omega_1,\dots
,\Omega_p,\phi_1,\dots ,\phi_{u+1})=0,
\end{split}
\end{equation}
modulo complete contractions of length $\ge\sigma +u+2$. Here each
$C^{\nu,i_1\dots i_\mu}_{g}$ is acceptable and has a simple
character $\vec{\kappa}^{+}_{simp}$ and a double character that is
 doubly subsequent to each $\vec{L}^{z,\sharp}, z\in Z'_{Max}$.\footnote{$\vec{L}^{z,\sharp}$ 
 is the refined $(u+1,\mu-1)$-double character of the 
 tensor fields $\dot{C}^{l,i_1\dots
\hat{i}_{r\alpha+1}\dots i_\mu,i_{*}}_{g}
(\Omega_1,\dots,\Omega_p,\phi_1,\dots ,\phi_u)$, $z\in Z'_{Max}$.}
$$\Sum_{t\in T_1} a_t
C^{t, i_1\dots i_{z_t}}_{g}(\Omega_1,\dots ,\Omega_p,\phi_1,\dots
,\phi_{u+1})$$
 is a generic linear combination of acceptable tensor fields
 with  a $(u+1)$-simple character $\vec{\kappa}^{+}_{simp})$,
 and with $z_t\ge \mu$.

$$\Sum_{t\in T_2} a_t
C^{t, i_1\dots i_{z_t}}_{g}(\Omega_1,\dots ,\Omega_p,\phi_1,\dots
,\phi_{u})$$
 ($z_t\ge \mu +1$) is a generic linear combination of acceptable
 tensor fields with a $u$-simple character $\vec{\kappa}_{simp}$, with the additional
restriction that the free index ${}_{i_1}$ that belongs to the
(a) crucial factor\footnote{I.e.~the second critical factor, in this
case.} is a special free index.\footnote{Recall that a special free
index is either an index ${}_k,{}_l$ in a factor
$S_{*}\nabla^{(\nu)}R_{ijkl}$ or an internal index in a factor
$\nabla^{(m)}R_{ijkl}$.}

\par Now, $t\in T_3$ means that there is one unacceptable factor
$\nabla\Omega_h$ (and it is not contracting against any factor
$\nabla\phi_t$) and moreover the tensor fields indexed in $T_3$
have $(u+1)$-simple character $\vec{\kappa}^{+}_{simp}$ and
$z_t\ge \mu$.\footnote{If $z_t=\mu$ then we additionally claim that
$\nabla\phi_{u+1}$ is contracting against a derivative index, and
if it is contracting against a factor $\nabla^{(B)}\Omega_h$ then
$B\ge 3$; moreover, in this case $C^{t,i_1\dots i_\mu}_g$ will
contain no special free indices.}

\par The sublinear combination $\Sum_{t\in T_4}\dots$ appears
only if the second critical factor is of the form
$\nabla^{(B)}\Omega_k$, for some $k$. In that case, $t\in T_4$
means that there is one unacceptable factor $\nabla\Omega_k$, and
it is contracting against a factor $\nabla\phi_r$:
$\nabla_i\Omega_k\nabla^i\phi_r$, and moreover if $z_t=\mu$ then
one of the free indices ${}_{i_1},\dots,{}_{i_{\mu}}$ is a derivative index, 
 and moreover if it belongs to a factor $\nabla^{(B)}\Omega_h$ then $B\ge 3$.

Finally,
$$\Sum_{j\in J}  a_j C^j_{g}(\Omega_1,\dots
,\Omega_p,\phi_1,\dots ,\phi_{u+1})$$ stands for a generic linear
combination of complete contractions that are $u$-simply subsequent to
$\vec{\kappa}_{simp}$.
\newline

In case B, we just claim that Proposition \ref{giade} is true. 
\end{lemma}

\subsection{Outline of Part A: The main strategy.} 
In part A of this paper we prove Lemmas \ref{zetajones} and \ref{pool2}, and 
set up the groundwork for the proof of Lemma \ref{pskovb} in part B. 

\par The starting point of this proof will be the analysis of one local equation: 
We denote the assumption of these Lemmas (the equation (\ref{hypothese2})) by
 $L_g(\Omega_1,\dots,\Omega_p,\phi_1,\dots,\phi_u)=0$, or just $L_g=0$ for short.
 The point of departure of our analysis will be to study the {\it first conformal variation}
  of this equation, $Image^1_{\phi_{u+1}}[L_g]=0$.\footnote{Since the equation $L_g=0$ is 
  assumed to hold for all Riemannian metrics $g$, we may consider its first variation
  under conformal deformations of $g$; i.e. $Image^1_{\phi_{u+1}}[L_g]=0$ is the new local equation 
  $\frac{d}{dt}|_{t=0}L_{e^{2t\phi_{u+1}}g}=0$.} 
  
Now, our first result here is to pick out a specific sublinear combination, 
$Image^{1,+}_{\phi_{u+1}}[L_g]$ in  $Image^1_{\phi_{u+1}}[L_g]$ 
and to prove that it must vanish separately, 
modulo junk terms that we do not care about; this is the content of Lemma \ref{lemtsabes} and is  
done in subsection \ref{otlinemi}. Roughly speaking, the sublinear combination $Image^{1,+}_{\phi_{u+1}}[L_g]$  consists
 of the terms in $Image^1_{\phi_{u+1}}[L_g]$  which have one of two properties: 
 
 \begin{enumerate}
\item  {\it Either} they have $\sigma+u+1$ factors 
($u$ of them in the form $\nabla\phi_1,\dots,\nabla\phi_u$ 
 and a new one $\nabla\phi_{u+1}$),\footnote{Recall that 
the terms in $L_g$ have $\sigma+u$ factors, $u$ of
  them in the form $\nabla\phi_1,\dots,\nabla\phi_u$).} 
 have the $u$-weak character 
$Weak(\vec{\kappa}_{simp})$,\footnote{This is the
same $u$-weak character as for all terms in (\ref{hypothese2}).}
In rough terms, this means that the factors $\nabla\phi_1,\dots,\nabla\phi_u$ 
contract against the different factors 
$\nabla^{(m)}R_{ijkl},\nabla^{(p)}\Omega_h$ according to
the same {\it pattern}, 
and also the new factor 
$\nabla\phi_{u+1}$ contracts against the ``correct'' factor.

\item {\it Or}, they have $\sigma+u$ factors, $u$ of them  
in the form $\nabla\phi_1,\dots,\nabla\phi_u$, and a new one 
in the form $\nabla^{(A+2)}\phi_{u+1}$; this new factor has replaced one of the factors $\nabla^{(A)}R_{ijkl}$ 
or $S_{*}\nabla^{(A)}R_{ijkl}$ in $L_g$, by virtue of the transformation law (\ref{curvtrans}).
 In this case, we additionnaly require that if we formally replace the
  factor $\nabla^{(A)}_{r_1\dots r_{A-2}r_{A-1}r_A}\phi_{u+1}$  by 
factor $\nabla^{(A-3)}_{r_1\dots r_{A-3}}R_{r_{A-2}r_{A-1}sr_A}\nabla^s\phi_{u+1}$
(by virtue of applying the curvature identity 
to the indices ${}_{r_{A-2}}, {}_{r_{A-1}}$), then the resulting term 
would satisfy the first property above. 
\end{enumerate}

We next naturally break up $Image^{1,+}_{\phi_{u+1}}[L_g]$ into three sublinear 
combinations, see (\ref{heidegger}) below.
In the rather technical subsection \ref{prelwork} we ``get rid'' of a specific sublinear combination in 
$Image^{1,+}_{\phi_{u+1}}[L_g]$  which would otherwise cause us trouble.

 In section \ref{theanalysis} we consider the terms in $Image^{1,+}_{\phi_{u+1}}[L_g]$ which have 
$\sigma+u$ factors. This sublinear combination is denoted by $CurvTrans[L_g]$. Our 
aim is to ``get rid'' of these terms (since the Lemmas we are proving
assert claims about linear combinations with $\sigma+u+1$ factors), 
by introducing correction terms with $\sigma+u+1$ factors in total, 
{\it which we can control}.\footnote{This phrase {\it mostly} means that the 
corrections terms will be generic linear combinations which are allowed in the RHSs of
the Lemmas we are proving.} In order to obtain correction terms which we can control, we 
 argue as follows: We establish that the sublinear combination 
 $CurvTrans[L_g]$ vanishes separately, modulo longer correction 
terms (which apriori we can {\it not control}). 
 Moreover (as we check after multiple calculations in
 section \ref{anal3sub}), $CurvTrans[L_g]$  {\it retains} a lot 
 of the algebraic structure of the terms in $L_g$; in particular, 
 it can be expressed a  linear combination of $Xdiv$'s of high 
 order, plus terms which are ``simply subsequent'', in an appropriate sense. 
Thus, we iteratively  apply {\it the inductive assumption} of Proposition \ref{giade}, to derive 
that we can write $CurvTrans[L_g]=(Correction.Terms)$, 
where the correction terms have length $\sigma+u+1$
and {\it also} retain the algebraic structure that we want. 
We note that the analysis of section \ref{anal3sub} applies to Lemmas 
\ref{zetajones}, \ref{pool2} {\it and} to Lemma \ref{pskovb}.

\par Finally, in section \ref{easyproof} we deal directly with the terms in $Image^{1,+}_{\phi_{u+1}}[L_g]$ 
which are ``born'' with $\sigma+u+1$ factors. 
Our analysis in part A applies only to Lemmas \ref{zetajones}, \ref{pool2}. 
After long calculations, we find that these terms have the algebraic features
we would like them to; so if we add that sublinear combination to the correction terms we obtained from 
$CurvTrans[L_g]$, we straightforwardly derive our Lemmas \ref{zetajones}, \ref{pool2}. 

A note is in order here: As explained in the simplified version of Lemmas \ref{zetajones}, \ref{pool2}, 
terms with $\sigma+u+1$ factors for which the factor $\nabla\phi_{u+1}$ 
contracts against a derivative index in
 a curvature factor $\nabla^{(m)}R_{ijkl}$ are considered ``junk terms'' in the
 statements of those Lemmas, and are thus allowed in the RHSs of our claims;
we do not have to worry about their algebraic form. In other words, we are allowed to ``throw away'' many terms 
 that apear in $Image^{1,+}_{\phi_{u+1}}[L_g]$, since tehy are allowed
in the RHSs of the Lemmas \ref{zetajones}, \ref{pool2}. This simplifies the 
 task of proving the Lemmas \ref{zetajones}, \ref{pool2}, and indeed we are
 able to derive them by this long analysis of the {\it single} equation 
 $Image^{1,+}_{\phi_{u+1}}[L_g]=0$. We will see in part B of this paper
 that the corresponding task for Lemma  3.5 in \cite{alexakis4} will be more arduous. 
\newline

{\it Technical remark:} In the setting of Lemma \ref{pool2} there  are certain special cases where 
the Lemma can not be derived from the analysis below; in those cases, Lemma \ref{pool2} will be
 derived directly, in a ``Mini-Appendix'' at 
the end of pat A.\footnote{In fact, in that case the claims of 
Lemma \ref{pool2} and Proposition \ref{giade} coincide.} These special cases are when the 
tensor fieds of maximal refined double character in (\ref{hypothese2}) 
satisfy:
\begin{enumerate}
\item Any factor $\nabla^{(m)}R_{ijkl}$ in one of the forms 
$\nabla^{(m)}_{r_1\dots r_m}R_{\sharp\sharp\sharp\sharp}$ or 
 $\nabla^{(m)}_{r_1\dots r_m}R_{(free)\sharp\sharp\sharp}$, where each 
of the indices ${}_{r_1},\dots ,{}_{r_m}$ contracts against a factor $\nabla\phi_h$;
  each index ${}_{\sharp}$ is contracting against 
an index in another factor in $C^{l,i_1\dots i_\mu}_g$. 

\item All the other factors in $C^{l,i_1\dots i_\mu}_g$ 
are simple factors in the form $S_{*}R_{ijkl}$, or  factor 
$\nabla^{(2)}\Omega_h$ which are either simple and contain 
at most one free index or contract against exactly one $\nabla\phi_h$ and 
contain no free index.\footnote{Recall that a factor $\nabla^{(A)}\Omega_h$ 
 is called ``simple'' if it is not contracting against any factor $\nabla\phi_h$. 
A factor $S_{*}\nabla^{(\nu)}R_{ijkl}$ is called ``simple'' 
if its indices ${}_{r_1},\dots,{}_{r_\nu},{}_j$ 
are not contracting against any factor $\nabla\phi'_h$.}
\end{enumerate}

 So when we deal with the 
 Lemma \ref{pool2} below, we will be assuming that  the tensor fields of maximal
  refined double character in (\ref{hypothese2}) are 
{\it not} in the forms described above.

\section{Proof of Lemmas \ref{zetajones} and 
\ref{pool2}: Notation and preliminary results.} \label{dasproof}

\subsection{Codification of the assumptions:}

\par Recall that the main assumption of Proposition \ref{giade} is
the equation (\ref{hypothese2}). This equation is also the main
assumption for each of the Lemmas \ref{zetajones}, \ref{pool2} and
Lemma \ref{pskovb}. Recall that we are seeking to prove Lemmas
\ref{zetajones}, \ref{pool2} and Lemma \ref{pskovb} for weight $-n$,
where all the {\it tensor fields} in (\ref{hypothese2}) have  a
given simple character $\vec{\kappa}_{simp}$ (with $u$ factors
$\nabla\phi$, and $p$ factors $\nabla^{(y)}\Omega_h$).

\par In the remainder of this paper, we will prove Lemmas \ref{zetajones} 
and \ref{pool2} and set the groundwork for the proof of Lemma \ref{pskovb} in part B. 
\newline

\subsection{Introduction: Some technical tools.}
\label{otlinemi}

\par The main tool in proving the Lemmas \ref{zetajones}, \ref{pool2} in the present paper,
will be a careful analysis of one equation. This analysis will also be one of the main ingredients 
of the proof of Lemma \ref{pskovb} in 
part B of this paper.
The starting point for this analysis will be the first 
conformal variation of the hypotheses of
 Lemmas \ref{zetajones} and \ref{pool2}. Recall that 
we denote the hypothesis of Lemmas \ref{zetajones} and \ref{pool2}
 by $L_g(\Omega_1,\dots,\Omega_p,\phi_1,\dots,\phi_u)=0$, for short. 
 Then, the first conformal variation is the equation 
 $$Image^1_{{\phi}_{u+1}}[L_g(\Omega_1,\dots,\Omega_p,\phi_1,\dots,\phi_u)]=0,$$
where we recall that  $Image^1_{\phi_{u+1}}[L_g]$ 
is defined via the formula $Image^1_{\phi_{u+1}}[L_g]=\frac{d}{dt}|_{t=0} L_{e^{2t\phi_{u+1}}}(\Omega_1,\dots,\Omega_p,\phi_1,\dots,\phi_u)$. 
Given that $L_g$ consists of complete contractions 
involving factors $\nabla^{(m)}R_{ijkl}$, 
$\nabla^{(p)}\Omega_h$, $\nabla\phi_h$, it will be useful to recall 
the transformation law (under the conformal change 
${\hat{g}}_{ij}(x)=e^{2\phi(x)} g_{ij}(x)$) of 
the tensor $R_{ijkl}$ and the Levi-Civita connection $\nabla$:

\begin{equation}
\label{curvtrans}
\begin{split}
&R_{ijkl}(e^{2\phi}g)=e^{2\phi}[R_{ijkl}(g)+ \nabla^{(2)}_{il}\phi
g_{jk}+\nabla^{(2)}_{jk}\phi g_{il}-\nabla^{(2)}_{ik}\phi g_{jl}-\nabla^{(2)}_{jl}\phi g_{ik}
\\&+\nabla_i\phi\nabla_k{\phi}g_{jl}+\nabla_j{\phi}\nabla_l{\phi}g_{ik}-\nabla_i{\phi}\nabla_l {\phi}
g_{jk} -\nabla_j{\phi}\nabla_k{\phi}g_{il}
\\&+|\nabla\phi|^2g_{il}g_{jk}- |\nabla\phi|^2g_{ik}g_{lj}],
\end{split}
\end{equation}

\begin{equation}
\label{levicivita} {\nabla}_k {\eta}_l(e^{2\phi}g)=
{\nabla}_k{\eta}_l(g) -\nabla_k{\phi} {\eta}_l -\nabla_l{\phi} {\eta}_k
+\nabla^s{\phi} {\eta}_s g_{kl}.
\end{equation}

\par Recall equation (\ref{hypothese2}), which we re-write:

\begin{equation}
\label{hypothesegen}
\begin{split}
&L_{g}(\Omega_1,\dots ,\Omega_p,\phi_1,\dots ,\phi_u) =\Sum_{l\in
L_\mu} a_l Xdiv_{i_1}\dots Xdiv_{i_\mu} C^{l,i_1\dots ,i_\mu}_{g}
(\Omega_1,\dots ,\Omega_p,\\&\phi_1,\dots ,\phi_u)
+\Sum_{l\in L'} a_l Xdiv_{i_1}\dots Xdiv_{i_a}
C^{l,i_1\dots ,i_a}_{g} (\Omega_1,\dots ,\Omega_p,\phi_1,\dots
,\phi_u)+
\\&\Sum_{j\in J\bigcup J^v} a_j C^j_{g}(\Omega_1,\dots ,\Omega_p,
\phi_1,\dots ,\phi_u)=0;
\end{split}
\end{equation}
the complete contractions indexed in $J^v$ have length $\ge\sigma
+u+1$ (the ones indexed in $J$ have length $\sigma+u$). All {\it
tensor fields} above are {\it acceptable},\footnote{Recall that this means that all 
functions $\Omega_h$, $1\le h\le p$ are differentiated at least twice.} 
and have a given $u$-simple
character $\vec{\kappa}_{simp}$.\footnote{See the definition 2.5 in  \cite{alexakis4} 
for the precise definition of this notion.} 
The tensor fields indexed in
$L_\mu$ have rank $\mu$, and the ones indexed in $L'$ have
rank strictly
greater than $\mu$.\footnote{The rank in question is not necessarily
the same for each $l\in L'$.} The complete contractions $C^j$ are
 simply subsequent to $\vec{\kappa}_{simp}$. The above equation holds
  perfectly--not modulo longer complete contractions.

\par Clearly, the above equation implies that:

\begin{equation}
\label{mniaaa} Image^1_{\phi_{u+1}}[L_g(\Omega_1,\dots
,\Omega_p,\phi_1,\dots ,\phi_u)]=0.
\end{equation}

In fact, we will not be using the equation (\ref{mniaaa}) itself to prove
Lemmas \ref{zetajones} and \ref{pool2}, but  a slight variant of
it: 

\par  We inquire which of the factors $\nabla^{(t_1)}\Omega_1,\dots
,\nabla^{(t_p)}\Omega_p$ in $\vec{\kappa}_{simp}$ are not
contracting against any factor $\nabla\phi_f$. With no loss of
generality, we assume they are the factors
$\nabla^{(t_1)}\Omega_1,\dots ,\nabla^{(t_Y)}\Omega_Y$.

 We will then be considering the equation:

\begin{equation}
\label{veryimportant}
\begin{split}
&S_{g}(\Omega_1,\dots ,\Omega_p,\phi_1,\dots ,\phi_{u+1})=
Image^1_{\phi_{u+1}}[L_{g}(\Omega_1,\dots ,\Omega_p,\phi_1,\dots ,
\phi_u)]
\\&+L_{g}(\Omega_1\cdot \phi_{u+1},\Omega_2,\dots ,\Omega_p,
\phi_1,\dots ,\phi_u)+\dots 
\\&+L_{g}(\Omega_1,\dots , \Omega_Y\cdot
\phi_{u+1},\dots ,\Omega_p,\phi_1,\dots ,\phi_u) =0,
\end{split}
\end{equation}
which holds perfectly, i.e. without correction terms.

\par Lemmas \ref{zetajones}, \ref{pool2} and \ref{pskovb} will be proven by
carefully analyzing this equation.  For now, we start by recalling
some facts regarding ``simple characters''.

\par We recall that for all the tensor fields and complete
contractions appearing in (\ref{hypothesegen}) and for each factor
$\nabla\phi_f$, $f, 1\le f\le u$, there is a unique factor
$\nabla^{(m)} R_{ijkl}$ or  $S_{*}\nabla^{(\nu)} R_{ijkl}$ or
$\nabla^{(B)}\Omega_h$ against which $\nabla\phi_f$ is
contracting. Therefore, for each $f,1\le f\le u$ we may
unambiguously speak of {\it the} factor against which
$\nabla\phi_f$ is contracting in each of the tensor fields and
contractions in (\ref{hypothese2}).

\par On the other hand, we may have factors $\nabla^{(m)} R_{ijkl}$ in
$\vec{\kappa}_{simp}$ that are not contracting against any factor
$\nabla\phi_h$. We recall that there is the same number of such
factors for all tensor fields with  the given simple character
$\vec{\kappa}_{simp}$. We will sometimes refer to such factors as
``generic factors of the form $\nabla^{(m)} R_{ijkl}$''.
\newline

{\bf Notational conventions:} Examine the conclusions of Lemmas
\ref{zetajones}, \ref{pool2}. Focus on the first
lines of those conclusions. For Lemma \ref{pool2}
 all the tensor fields in the first line have the
same $(u+1)$-simple character which we will denote by
$\vec{\kappa}^{+}_{simp}$. For Lemma \ref{zetajones} we observe that
all tensor fields in the first line of the conclusion have the
same $u$-simple character $\vec{\kappa}_{simp}$, and furthermore
the factor $\nabla\phi_{u+1}$ is contracting against the index
${}_k$ in a specified factor $T=S_{*}\nabla^{(\nu)}R_{ijkl}$ in
$\vec{\kappa}_{simp}$. While strictly speaking these tensor fields
are not in the form (\ref{form2}),\footnote{Because the factor $\nabla\phi_{u+1}$
is contracting against a special index in a factor $S_{*}R_{ijkl}$.} and hence we cannot speak of a
simple character, we will abuse language and define that in the
context of the proof of Lemma \ref{zetajones}, any tensor field
that is contracting against $(u+1)$ factors $\nabla\phi_h$ has a
$(u+1)$-simple character $\vec{\kappa}^{+}_{simp}$ if it satisfies
the properties explained in the previous sentence.

\par Now, recall that there is a well-defined notion of the ``crucial factor(s)'' in the
context of Lemmas \ref{zetajones}, 
\ref{pool2}, introduced in \cite{alexakis4}. In particular, for
each of the three Lemmas above, we may unambiguously speak of the
set of factors $\nabla\phi_h$ in $\vec{\kappa}_{simp}$ that are
contracting against the crucial factor. (Recall that this set may
also be empty--in that case we have a {\it set} of crucial
factors, which are all the factors $\nabla^{(m)}R_{ijkl}$ that are
not contracting against any $\nabla\phi_h$).

\par We denote by $(\vec{\kappa}_{simp})_1$ the set of numbers $h$ for which
$\nabla\phi_h$ is contracting against the crucial factor.

\par We will now introduce some further notation, for future
reference. We will formally construct a new $u$-simple character
{\it for contractions with $\sigma+u$ factors} (whereas
$\vec{\kappa}_{simp}^{+}$ corresponds to contractions with
$\sigma+u+1$ factors in total). (The definition that follows is highly un-intuitive,
 but the discussion below provides some intuition).

\begin{definition}
\label{laborat} Consider the simple character
$\vec{\kappa}_{simp}$. Pick a tensor field (in the form (\ref{form2})) with a $u$-simple
character $\vec{\kappa}_{simp}$. In the case where the crucial
factor(s) is (are) of the form
$S_{*}\nabla^{(\nu)}R_{ijkl}$ or $\nabla^{(m)}R_{ijkl}$, we will define a new $u$-simple
character $pre\vec{\kappa}_{simp}^{+}$ as follows:

\par Recall that $\vec{\kappa}_{simp}$ is a list of sets.
Consider the entry in $\vec{\kappa}_{simp}$ that corresponds to
the crucial factor. That entry will either be a set of numbers
$S_h$ or a set in the form $(\{\alpha\},S_h)$, where $\alpha$ is a
number and $S_h$ a set of numbers. Respectively, the entry will
either belong to the list $L_1$ or the list $L_2$. Then
$pre\vec{\kappa}_{simp}^{+}$ arises from the simple character
$\vec{\kappa}_{simp}$ by erasing this entry $(\{\alpha\},S_h)$ or $S_h$
 from $L_1$ or $L_2$ and adding an entry
$S_{p+1}=S_h\bigcup \{\alpha\}$ or $S_{p+1}=S_h$, respectively, in
$L_3$.
\end{definition}

\par  A more intuitive description of $pre\vec{\kappa}_{simp}^{+}$
is the following: Consider any complete contraction $C_g$ with a
$u$-simple character $\vec{\kappa}_{simp}$. Consider a crucial
factor $T$ in $C_g$,\footnote{For this definition the crucial factor must be a curvature factor.}
 along with its indices, $T_{r_1\dots
r_{m+4}}$. Assume that two of the indices in $T$ (say
${}_{r_{m+2}},{}_{r_{m+4}}$) are not contracting against factors
$\nabla\phi$ (this can always be done, i.e. we can always find a
contraction $C_g$ that satisfies this requirement-by adding
derivative indices onto $T$ if necessary).

Formally replace $T_{r_1\dots r_{m+4}}$ by a new factor
$\nabla^{(m+4)}_{r_1\dots r_{m+1}r_{m+3}}\Omega_{p+1}$. The indices that
contracted against ${}_{r_{m+2}},{}_{r_{m+4}}$ now become free.
$pre\vec{\kappa}_{simp}^{+}$ is then the $u$-simple character of
this new partial contraction.
\newline

\par Now, in the setting of Lemma \ref{pskovb} 
we will slightly generalize the above notions.
Although we have defined a notion of ``crucial factor'' (which is
either the
 ``critical'' or the ``second critical'' factor), we wish to allow ourselves
 some extra freedom, and thus we will be picking another factor (set of factors)
in $\vec{\kappa}_{simp}$ and we will call it (them) the {\it selected factor(s)}.

\par Our choice of selected factor(s) is entirely free: We can either pick any
well-defined factor in $\vec{\kappa}_{simp}$ and call it the
selected factor (recall a {\it well-defined} factor is either a
curvature factor that is contracting against some $\nabla\phi$ or
a factor $\nabla^{(p)}\Omega_h$), or we can pick the {\it set} of
factors $\nabla^{(m)}R_{ijkl}$ in $\vec{\kappa}_{simp}$ that are
not contracting against any factors $\nabla\phi$ and call all
those factors the {\it selected factors}.

\par Having chosen a (set of) selected factor(s) in $\vec{\kappa}_{simp}$, say $F_s$,
we then formally construct the $(u+1)$-simple character $\vec{\kappa}_{simp}^{+}$
by just adding a derivative $\nabla_a$ onto
the (one of the) factor(s) $F_s$ and just contracting
 ${}_a$ against a new factor $\nabla^a\phi_{u+1}$ (and $S_{*}$-symmetrizing if $F_s$
is of the form $S_{*}\nabla^{(\nu)}R_{ijkl}$). Then, if $F_s$ is a curvature factor,
 we define the $(u+1)$-simple character
 $pre\vec{\kappa}_{simp}^{+}$ as in
 Definition \ref{laborat}. (If $F_s$ is not a curvature term,
 then $pre\vec{\kappa}_{simp}^{+}$ is undefined).
\newline

\par A note: In the rest of this section, we will be referring
exclusively to the ``selected'' factor. In the setting of Lemmas
 \ref{zetajones}, \ref{pool2} and in 
 case A of Lemma \ref{pskovb}, we will take the selected factor(s) to
 be the crucial factor(s). Therefore, in most circumstances
the two notions coincide.
\newline

\par Now, we define sublinear combinations in
$Image^1_{\phi_{u+1}}[L_{g}]$,
 where $L_{g}$ is the left hand side of our hypothesis (\ref{hypothesegen}).
We note that these definitions still make sense for {\it any}
linear combination $L_{g}$ consisting of complete contractions
with a given weak character $\vec{\kappa}$ and a given
 (set of) selected factor(s) $T$ in $\vec{\kappa}$.

\begin{definition}
\label{plus}
 We denote by
$Image^{1,+}_{\phi_{u+1}}[L_{g}(\Omega_1,\dots ,
\Omega_p,\phi_1,\dots ,\phi_u)]$ the sublinear combination in
$Image^{1}_{\phi_{u+1}}[L_{g}(\Omega_1,\dots ,
\Omega_p,\phi_1,\dots ,\phi_u)]$ that consists of complete
contractions $C_{g}(\Omega_1,\dots , \Omega_p,\phi_1,\dots
,\phi_u)$ with the following properties:
\begin{enumerate}
\item{$C_{g}(\Omega_1,\dots , \Omega_p,\phi_1,\dots
,\phi_u,\phi_{u+1})$ must have no internal contractions.}

\item{If $C_{g}(\Omega_1,\dots , \dots ,\Omega_p,\phi_1,\dots
,\phi_u,\phi_{u+1})$  has length $\sigma +u+1$ then it must have a
factor $\nabla\phi_{u+1}$ (with only one derivative) and a weak
character Weak($\vec{\kappa}^{+}_{simp})$.}

\item{If $C_{g}(\Omega_1,\dots , \Omega_p,\phi_1,\dots
,\phi_u,\phi_{u+1})$  has length $\sigma +u$, then it has a factor
$\nabla^{(A)}\phi_{u+1}$, $A\ge 2$. Moreover, it must have a weak
$u$-character \\$Weak(pre\vec{\kappa}^{+}_{simp})$.}
\end{enumerate}
\end{definition}

\begin{definition}
\label{newdefinition} Define
$Image^{1,\alpha}_{\phi_{u+1}}[L_{g}(\Omega_1,\dots ,
\Omega_p,\phi_1,\dots ,\phi_u)]$ to be the sublinear combination
in $Image^1_{\phi_{u+1}}[L_{g}(\Omega_1,\dots ,
\Omega_p,\phi_1,\dots ,\phi_u)]$ that consists of complete
contractions $C_{g}(\Omega_1,\dots , \Omega_p,\phi_1,\dots
,\phi_u)$ with length $\sigma +u+1$ and a factor
$\nabla^{(A)}\phi_{u+1}$ with $A\ge 2$. We also denote a generic
linear combination of such complete contractions by $\Sum_{z\in Z}
a_z C^z_{g}(\Omega_1,\dots ,\Omega_p, \phi_1,\dots ,\phi_{u+1})$.

\par We define $Image^{1,\beta}_{\phi_{u+1}}[L_{g}(\Omega_1,
\dots ,\Omega_p,\phi_1,\dots ,\phi_u)]$ to stand for the
 sublinear combination of complete contractions
$C_{g}(\Omega_1,\dots ,\Omega_p,\phi_1, \dots ,\phi_{u+1})$ in
\\ $Image^1_{\phi_{u+1}}[L_{g}(\Omega_1,\dots ,
\Omega_p,\phi_1,\dots ,\phi_u)]$ with the
 following properties:

\begin{enumerate}
\item{$C_{g}(\Omega_1,\dots ,\Omega_p,\phi_1, \dots ,\phi_{u+1})$
must have precisely one internal contraction.}

\item{Either $C_{g}(\Omega_1,\dots ,\Omega_p,\phi_1, \dots
,\phi_{u+1})$ has length $\sigma +u+1$ and a factor
$\nabla\phi_{u+1}$.}

\item{Or $C_{g}(\Omega_1,\dots ,\Omega_p,\phi_1, \dots
,\phi_{u+1})$ has length $\sigma +u$ and a factor
$\nabla^{(A)}\phi_{u+1}$ ($A\ge 2$) (which does not contain the internal contraction).}

\end{enumerate}

\par Finally, we define
$Image^{1,\gamma}_{\phi_{u+1}}[L_{g}(\Omega_1, \dots
,\Omega_p,\phi_1,\dots ,\phi_u)]$ to stand for the
 sublinear combination of complete contractions
$C_{g}(\Omega_1,\dots ,\Omega_p,\phi_1, \dots ,\phi_{u+1})$ in
$Image^1_{\phi_{u+1}}[L_{g}(\Omega_1,\dots , \Omega_p,\phi_1,\dots
,\phi_u)]$ with one of the properties:

\begin{enumerate}
\item{Either $C_{g}(\Omega_1,\dots ,\Omega_p,\phi_1, \dots
,\phi_{u+1})$ has length $\sigma+u$, a factor $\nabla^{(A)}\phi_{u+1}$ 
 and a weak character which is not $Weak(pre\vec{\kappa}_{simp}^{+})$,}

\item{Or $C_{g}(\Omega_1,\dots ,\Omega_p,\phi_1, \dots
,\phi_{u+1})$ has length $\sigma+u+1$ and a factor
$\nabla\phi_{u+1}$ but its
 weak character is {\it not}
$Weak(\vec{\kappa}^{+}_{simp})$.}
\end{enumerate}
\end{definition}

\par We make note of the fact that we are {\it not}
imposing any restriction to the weak character of the complete
contractions that belong to \\$Image^{1,\beta}_{\phi_{u+1}}
[L_{g}(\Omega_1,\dots ,\Omega_p,\phi_1,\dots ,\phi_u)]$.

\par Now, assuming any equation of the form $L_{g}(\Omega_1, \dots
,\Omega_p,\phi_1,\dots ,\phi_u)=0$ (this equation is assumed to hold perfectly--here
$L_{g}$ consists of complete contractions with length $\ge\sigma+u$) we observe that:

\begin{equation}
\label{halevy}
\begin{split}
&Image^1_{\phi_{u+1}}[L_{g}(\Omega_1, \dots
,\Omega_p,\phi_1,\dots ,\phi_u)]=
Image^{1,+}_{\phi_{u+1}}[L_{g}(\Omega_1, \dots
,\Omega_p,\phi_1,\dots ,\phi_u)]+
\\& Image^{1,\alpha}_{\phi_{u+1}}[L_{g}(\Omega_1,
\dots ,\Omega_p,\phi_1,\dots ,\phi_u)]+
Image^{1,\beta}_{\phi_{u+1}}[L_{g}(\Omega_1, \dots
,\Omega_p,\phi_1,\dots ,\phi_u)]+
\\& Image^{1,\gamma}_{\phi_{u+1}}[L_{g}(\Omega_1,
\dots ,\Omega_p,\phi_1,\dots ,\phi_u)](=0),
\end{split}
\end{equation}
(modulo complete contractions of length $\ge\sigma+u+2$).

This follows by the definitions above and the transformation laws
 (\ref{curvtrans}) and (\ref{levicivita}).
\newline

\par Our next claim will be used frequently in the future, so we present
 it in somewhat  general notation:

\begin{lemma}
\label{lemtsabes}
\par Assume that $L_{g}(\Omega_1,\dots ,
\Omega_p,\phi_1,\dots ,\phi_u)$ is a  linear combination of complete
 contractions with no internal contractions and with a given weak character $Weak(\vec{\kappa})$
(where $\vec{\kappa}$ is any chosen simple character consisting of
$\sigma+u$ factors).
 Assume that

 \begin{equation}
 \label{gkiourntani}
 L_{g}(\Omega_1,\dots ,
\Omega_p,\phi_1,\dots ,\phi_u)=0,
\end{equation}
 modulo complete
contractions with length $\ge\sigma+u+1$.  Then, in the notation of
Definition \ref{plus}, we claim:

\begin{equation}
\label{tsabes}
 Image^{1,+}_{\phi_{u+1}}[L_{g}(\Omega_1,\dots ,
\Omega_p,\phi_1,\dots ,\phi_u)]+ \Sum_{z\in Z} a_z C^z_{g}
(\Omega_1,\dots ,\Omega_p,\phi_1,\dots ,\phi_{u+1})=0,
\end{equation}
modulo complete contractions of length $\ge \sigma+u+2$. Here
$\sum_{z\in Z}$ stands for a generic linear combination of
contractions with $\sigma+u+1$ factors, one of which is in the form
$\nabla^{(c)}\phi_{u+1}$, $c\ge 2$.
\end{lemma}

{\it Proof of Lemma \ref{lemtsabes}:} Our point of departure is equation (\ref{halevy}).
In view of that equation, we notice that it would suffice to show
 that:

\begin{equation}
\label{chrysophgh1}
Image^{1,\beta}_{\phi_{u+1}}[L_{g}(\Omega_1,\dots ,
\Omega_p,\phi_1,\dots ,\phi_u)]=\Sum_{z\in Z} a_z C^z_{g}
(\Omega_1,\dots ,\Omega_p,\phi_1,\dots ,\phi_{u+1}),
\end{equation}

\begin{equation}
\label{chryspohgh2}
Image^{1,\gamma}_{\phi_{u+1}}[L_{g}(\Omega_1,\dots ,
\Omega_p,\phi_1,\dots ,\phi_u)]=\Sum_{z\in Z} a_z C^z_{g}
(\Omega_1,\dots ,\Omega_p,\phi_1,\dots ,\phi_{u+1}),
\end{equation}
modulo complete contractions of length $\ge\sigma +u+2$. If we can
show the above equations, then by virtue of (\ref{halevy}) we will
immediately deduce our claim.
\newline

\par We will denote by $Image^{1,\beta,\sigma +u}_{\phi_{u+1}}
[L_{g}(\Omega_1,\dots , \Omega_p,\phi_1,\dots ,\phi_u)]$,
\\$Image^{1,\gamma,\sigma +u}_{\phi_{u+1}}[L_{g}(\Omega_1,\dots
,\Omega_p,\phi_1,\dots ,\phi_u)]$ the sublinear combinations of
complete
 contractions with length $\sigma +u$ in
\\ $Image^{1,\beta}_{\phi_{u+1}}
[L_{g}(\Omega_1,\dots ,\Omega_p,\phi_1,\dots ,\phi_u)]$,
 $Image^{1,\gamma}_{\phi_{u+1}}[L_{g}(\Omega_1,\dots ,\Omega_p,\phi_1,\dots ,\phi_u)]$, 
respectively. We will also denote by $\Sum_{w\in W^\beta} a_w
C^w_{g}(\Omega_1,\dots ,\Omega_p,\phi_1,\dots ,\phi_u)$,\\
$\Sum_{w\in W^\gamma} a_w C^w_{g}(\Omega_1,\dots
,\Omega_p,\phi_1,\dots ,\phi_u)$ generic linear combinations of
complete
 contractions as the ones that belong to the linear combinations
\\$Image^{1,\beta}_{\phi_{u+1}} [L_{g}(\Omega_1,\dots
,\Omega_p,\phi_1,\dots ,\phi_u)]$,
 $Image^{1,\gamma}_{\phi_{u+1}}[L_{g}(\Omega_1,\dots ,
 \Omega_p,\phi_1,\dots ,\phi_u)]$, respectively, which in
 addition have length $\sigma +u+1$.

\par We claim that:

\begin{equation}
\label{halevy3}
\begin{split}
&Image^{1,\beta,\sigma +u}_{\phi_{u+1}}[L_{g}(\Omega_1,\dots ,
\Omega_p,\phi_1,\dots ,\phi_u)]= \Sum_{w\in W^\beta} a_w
C^w_{g}(\Omega_1,\dots , \Omega_p,\phi_1,\dots ,\phi_u)
\\& +\Sum_{z\in Z} a_z C^z_{g}(\Omega_1,\dots ,
\Omega_p,\phi_1,\dots ,\phi_{u+1}),
\end{split}
\end{equation}

\begin{equation}
\label{halevy4}
\begin{split}
&Image^{1,\gamma,\sigma+ u}_{\phi_{u+1}}[L_{g}(\Omega_1,\dots ,
\Omega_p,\phi_1,\dots ,\phi_u)]= \Sum_{w\in W^\gamma} a_w
C^w_{g}(\Omega_1,\dots , \Omega_p,\phi_1,\dots ,\phi_u)
\\& +\Sum_{z\in Z} a_z C^z_{g}(\Omega_1,\dots ,
\Omega_p,\phi_1,\dots ,\phi_{u+1}).
\end{split}
\end{equation}

These equations are proven by the usual argument: Recall that
(\ref{halevy}) holds formally. Then, we observe that if we pick
out the sublinear combinations in those equations that consist of
complete contractions with $\sigma+u$ factors (denote those sublinear
combinations by $Z_g$), then $lin\{Z_g\}=0$ formally. (Recall that $lin\{Z_g\}$
stands for the {\it linearization} of the linear combination $Z_g$--see 
\cite{alexakis1}). Now, since both the weak
character and the number of internal
 contractions are {\it invariant} under the permutations that make $lin\{Z_g\}$ formally zero,
we derive that the linearizations of the left hand sides of
(\ref{halevy3}), (\ref{halevy4}) must vanish formally. Hence, by
repeating the permutations in the non-linear setting, we obtain
the right hand sides in (\ref{halevy3}), (\ref{halevy4}) as
 correction terms.

\par Thus, substituting the above two equations into (\ref{halevy}),
we obtain an equation in place of (\ref{halevy}), where the sums analogous
to $Image^{1,\beta}_{\phi_{u+1}}[L_{g}(\Omega_1,\dots ,
\Omega_p,\phi_1,\dots ,\phi_u)]$,
$Image^{1,\gamma}_{\phi_{u+1}}[L_{g}(\Omega_1,\dots ,
\Omega_p,\phi_1,\dots ,\phi_u)]$ contain no
 complete contractions of length $\sigma+u$.
 So if we denote by
$Image^{1,+,\sigma +u}_{\phi_{u+1}}[L_{g}(\Omega_1,\dots ,
\Omega_p,\phi_1,\dots ,\phi_u)]$ the sublinear combination in
$Image^{1,+}_{\phi_{u+1}}[L_{g}(\Omega_1,\dots ,
\Omega_p,\phi_1,\dots ,\phi_u)]$ that consists of complete
 contractions with $\sigma+u$ factors, we obtain that
\\$Image^{1,+,\sigma +u}_{\phi_{u+1}}[L_{g}(\Omega_1,\dots ,
\Omega_p,\phi_1,\dots ,\phi_u)]=0$, modulo longer complete
contractions.
 But then, since this equation must hold formally at the linearized level,
 by just repeating the permutations that make the
 linearized linear combination formally zero we derive, modulo
  complete contractions of length $\ge\sigma+u+2$:

\begin{equation}
\label{asnd}
\begin{split}
&Image^{1,+,\sigma +u}_{\phi_{u+1}}[L_{g}(\Omega_1,\dots ,
\Omega_p,\phi_1,\dots ,\phi_u)]= \Sum_{h\in H} a_h
C^h_{g}(\Omega_1,\dots , \Omega_p,\phi_1,\dots ,\phi_{u+1})
\\&+\Sum_{z\in Z} a_z C^z_{g}(\Omega_1,\dots ,
\Omega_p,\phi_1,\dots ,\phi_{u+1}),
\end{split}
\end{equation}
where each $C^h_{g}(\Omega_1, \dots ,\Omega_p,\phi_1,\dots
,\phi_{u+1})$ has length $\sigma +u+1$ and a factor
$\nabla\phi_{u+1}$ and no internal contractions and also has
 a weak character $Weak(\vec{\kappa}^{+})$.

\par By replacing this also into (\ref{halevy}), we derive:

\begin{equation}
\label{maryland}
\begin{split}
&\Sum_{h\in H} a_h C^h_{g}(\Omega_1,\dots , \Omega_p,\phi_1,\dots
,\phi_{u+1})+ \Sum_{w\in W^\beta} a_w C^w_{g}(\Omega_1,\dots ,
\Omega_p,\phi_1,\dots ,\phi_{u+1})+
\\&\Sum_{w\in W^\gamma} a_w C^w_{g}(\Omega_1,\dots ,
\Omega_p,\phi_1,\dots ,\phi_{u+1})+ \Sum_{z\in Z} a_z
C^z_{g}(\Omega_1,\dots , \Omega_p,\phi_1,\dots ,\phi_{u+1})=0.
\end{split}
\end{equation}

\par Now, since the above must hold formally (and the weak
 character as well as the number of internal contractions are
invariant under the permutations that make the linearizations of
the left hand side formally zero), we derive that:

$$\Sum_{w\in W^\beta} a_w C^w_{g}(\Omega_1,\dots ,
\Omega_p,\phi_1,\dots ,\phi_{u+1})=0,$$
$$\Sum_{w\in W^\gamma} a_w C^w_{g}(\Omega_1,\dots ,
\Omega_p,\phi_1,\dots ,\phi_{u+1})=0,$$
modulo complete contractions of length $\ge\sigma +u+2$.
\newline

Thus, (\ref{halevy3}), (\ref{halevy4}) combined with the above two equation  complete 
the proof of Lemma \ref{lemtsabes}. $\Box$
\newline

{\bf Further break-up of
$Image^{1,+}_{\phi_{u+1}}[L_{g}(\Omega_1,\dots
,\Omega_p,\phi_1,\dots ,\phi_u)]$:}
\newline

 We break
$$Image^{1,+}_{\phi_{u+1}}[L_{g}(\Omega_1,\dots ,\Omega_p,\phi_1,\dots ,\phi_u)]$$
into three sublinear combinations: We define
$CurvTrans[L_{g}(\Omega_1,\dots ,\Omega_p,\phi_1,\dots ,\phi_u)]$
to stand for the sublinear combination that arises by applying the
transformation law (\ref{curvtrans}) to a
 factor $\nabla^{(m)}R_{ijkl}$ or $S_{*}\nabla^{(\nu)}R_{ijkl}$, in some complete contraction in
$L_{g}$. We observe that the complete contractions in
$CurvTrans[L_{g}(\Omega_1,\dots ,\Omega_p,\phi_1,\dots ,\phi_u)]$
have length $\sigma +u$. They will each be in the form:

\begin{equation}
\label{plusou}
\begin{split}
&contr(\nabla^{(m_1)}R_{ijkl}\otimes\dots\otimes\nabla^{(m_s)}R_{ijkl}
\otimes S_{*}\nabla^{(\nu_1)}R_{ijkl}\otimes\dots\otimes
S_{*}\nabla^{(\nu_\tau)}R_{ijkl}
\\&\nabla^{(m)}\phi_{u+1}\otimes\nabla^{(f_1)}\Omega_1\otimes\dots\otimes\nabla^{(f_p)}
\Omega_p \otimes\nabla\phi_1\otimes\dots\otimes\nabla\phi_u).
\end{split}
\end{equation}

 We denote by $LC[L_{g}(\Omega_1,\dots
,\Omega_p,\phi_1,\dots ,\phi_u)]$ the sublinear combination in
$Image^{1,+}_{\phi_{u+1}}[L_{g}(\Omega_1,\dots ,\Omega_p,\phi_1,\dots
,\phi_u)]$ that arises by applying the transformation
 law (\ref{levicivita}). Finally, we denote by
$W[L_{g}(\Omega_1,\dots ,\Omega_p,\phi_1,\dots ,\phi_u)]$ the
sublinear
 combination in $Image^{1,+}_{\phi_{u+1}}[L_{g}(\Omega_1,\dots ,
\Omega_p,\phi_1,\dots ,\phi_u)]$ that arises by applying the
 transformation $R_{ijkl}\longrightarrow e^{2\phi_{u+1}}
R_{ijkl}$ and bringing out an expression $\nabla\phi_{u+1}$.

\par Then, by definition:

\begin{equation}
\label{heidegger}
\begin{split}
&Image^{1,+}_{\phi_{u+1}}[L_{g}(\Omega_1,\dots
,\Omega_p,\phi_1,\dots ,\phi_u)]= CurvTrans[L_{g}(\Omega_1,\dots
,\Omega_p,\phi_1,\dots ,\phi_u)]
\\& +LC[L_{g}(\Omega_1,\dots ,\Omega_p,\phi_1,\dots ,
\phi_u)]+W[L_{g}(\Omega_1,\dots ,\Omega_p,\phi_1,\dots , \phi_u)].
\end{split}
\end{equation}

 Much of this section and of the next ones consists of
 understanding the three sublinear combinations above and of
 using our inductive assumption on Corollary 1 in \cite{alexakis1} in
 order to derive our three Lemmas \ref{zetajones}, \ref{pool2}, \ref{pskovb}.

\subsection{Preliminary Work.}
\label{prelwork}

We will be generically denoting all the tensor fields that appear
in (\ref{hypothese2}) by $C^{l,i_1\dots i_a}_{g}(\Omega_1,\dots
,\Omega_p,\phi_1, \dots ,\phi_u)$. Also, $C^j_{g}(\Omega_1,\dots
,\Omega_p, \phi_1,\dots ,\phi_u)$ will stand for a generic
complete contraction with a weak character
$Weak(\vec{\kappa}_{simp})$, where $C^j_g$ is simply subsequent to
$\vec{\kappa}_{simp}$.

\par We need a definition in order to formulate our claim:

\begin{definition}
\label{jhc} We define $LC_\Phi[Xdiv_{i_1}\dots Xdiv_{i_a}
C^{l,i_1\dots i_a}_{g}(\Omega_1,\dots,\Omega_p,\phi_1, \dots
,\phi_u)]$ and $LC_\Phi[C^j_{g}(\Omega_1,\dots,\Omega_p,
\phi_1,\dots ,\phi_u)]$ to stand for the sublinear combinations in
$LC[Xdiv_{i_1}\dots Xdiv_{i_a} C^{l,i_1\dots
i_a}_{g}(\Omega_1,\dots,\Omega_p,\phi_1, \dots ,\phi_u)]$,
$LC[C^j_{g}(\Omega_1,\dots,\Omega_p, \phi_1,\dots ,\phi_u)]$,\footnote{Recall
that by definition these sublinear 
combinations consist of complete contractions of length $\sigma +u+1$
with weak character $Weak(\vec{\kappa}^{+}_{simp})$.} that
arise when we bring out a factor $\nabla\phi_{u+1}$  by applying
the transformation law (\ref{levicivita})  to a pair
 of indices $(\nabla_A, {}_B)$
  where $\nabla_A$ denotes
 a derivative index, and where either $\nabla_A$ or ${}_B$ is
 contracting against a factor $\nabla\phi_h$ ($1\le h\le u$).
\end{definition}

 We observe that a complete contraction in
\\$LC_\Phi [Xdiv_{i_1}\dots Xdiv_{i_a} C^{l,i_1\dots
i_a}_{g}(\Omega_1,\dots,\Omega_p,\phi_1, \dots ,\phi_u)]$,
\\$LC_\Phi[C^j_{g}(\Omega_1,\dots,\Omega_p, \phi_1,\dots ,\phi_u)]$
can only arise by applying the transformation law
(\ref{levicivita}) to a pair of indices $(\nabla_A,{}_B)$ in one
of the following two ways:

\par Firstly, if the index $\nabla_A$ is contracting against a
 factor $\nabla\phi_h$ and the index ${}_B$ is contracting
 against the selected factor and we bring out the third
 summand in (\ref{levicivita}).  Alternatively, if ${}_B$ is
 contracting against a factor $\nabla\phi_h$ and $\nabla_A$ is
 contracting against the selected factor and we bring out the
second summand in (\ref{levicivita}).
\newline

\par The aim of this subsection is to show the following:

\begin{lemma}
\label{bjorn} Consider (\ref{hypothese2}). Then, in the notation
of definition \ref{jhc}, we claim:

\begin{equation}
\label{trikala}
\begin{split}
&\Sum_{l\in L} a_l LC_\Phi[Xdiv_{i_1}\dots Xdiv_{i_a}
C^{l,i_1\dots i_a}_{g}(\Omega_1,\dots,\Omega_p,\phi_1, \dots
,\phi_u)]+
\\&\Sum_{j\in J} a_j
LC_\Phi[C^j_{g}(\Omega_1,\dots,\Omega_p, \phi_1,\dots ,\phi_u)]=
\\&\Sum_{l\in L'} a_l Xdiv_{i_1}\dots Xdiv_{i_a}C^{l,i_1\dots
i_a,i_{*}}_{g}(\Omega_1,\dots,\Omega_p,\phi_1, \dots ,
\phi_u)\nabla_{i_{*}}\phi_{u+1}+
\\&\Sum_{j\in J} a_j
C^j_{g}(\Omega_1,\dots,\Omega_p, \phi_1,\dots,\phi_{u+1}),
\end{split}
\end{equation}
where $\Sum_{l\in L'} a_l C^{l,i_1\dots
i_a,i_{*}}_{g}(\Omega_1,\dots,\Omega_p,\phi_1, \dots ,
\phi_u)\nabla_{i_{*}}\phi_{u+1}$, in the setting of Lemmas
\ref{zetajones} and \ref{pool2}, stands for a generic linear
 combination of acceptable tensor fields of 
length $\sigma +u+1$, with $a\ge \mu$ and with a
$(u+1)$-simple character $\vec{\kappa}^{+}_{simp}$. On the other
hand, in the setting of Lemma \ref{pskovb} it stands for a generic
linear combination in the form
$$\Sum_{t\in T_1\bigcup T_2} a_t C^{t,i_1\dots i_{z_t}}_{g}
(\Omega_1,\dots ,\Omega_p,\phi_1,\dots
,\phi_{u+1})\nabla_{i_1}\phi_{u+1}.$$ 
 $\Sum_{j\in J} a_j C^j_{g}(\Omega_1,\dots,\Omega_p,
\phi_1,\dots, \phi_{u+1})$ stands for a generic linear combination
of complete contractions that are simply subsequent to
$\vec{\kappa}^{+}_{simp}$ (in the
 setting of all three Lemmas).
\end{lemma}

{\bf Proof of Lemma \ref{bjorn}:}
\newline

\par In order to show this Lemma, we will again introduce some
preliminary definitions. We recall that in equation
(\ref{hypothesegen}) all the complete contractions of length
$\sigma +u$ have the same weak character.\footnote{See \cite{alexakis4} for 
a precise definition of this notion.} For each partial
contraction with a simple character $\vec{\kappa}_{simp}$, we let
$\{F_1,\dots , F_X\}$ stand for the set of {\it non-selected}
factors that are contracting against some factor $\nabla\phi_f$.

Then, for each $1\le h\le X$, we define an operation
$Hit^h_{\nabla\tau}$ that formally acts on each tensor field
$C^{l,i_1\dots i_a}_{g}(\Omega_1,\dots ,\Omega_p,\phi_1,\dots
,\phi_u)$ and each complete contraction
 $C^j_{g}(\Omega_1,\dots,\Omega_p,\phi_1,\dots ,\phi_u)$ in
 (\ref{hypothese2})
by hitting the factor $F_h$ by a derivative index $\nabla_u$
and then  contracting ${}_u$ against a factor $\nabla_u\tau$. We also define:

\begin{equation}
\label{tsira}
\begin{split}
&Hit_{\nabla\tau}[Xdiv_{i_1}\dots Xdiv_{i_a} C^{l,i_1\dots
i_a}_{g}(\Omega_1,\dots ,\Omega_p,\phi_1,\dots ,\phi_u)]=
\\&\Sum_{h=1}^X
Hit^h_{\nabla\tau}[Xdiv_{i_1}\dots Xdiv_{i_a} C^{l,i_1\dots
i_a}_{g}(\Omega_1,\dots ,\Omega_p,\phi_1,\dots ,\phi_u)],
\end{split}
\end{equation}
and also:

\begin{equation}
\label{tsira2} Hit_{\nabla\tau}[C^j_{g}(\Omega_1,\dots
,\Omega_p,\phi_1,\dots ,\phi_u)]=\Sum_{h=1}^X
Hit^h_{\nabla\tau}[C^j_{g}(\Omega_1,\dots ,\Omega_p,\phi_1,\dots
,\phi_u)].
\end{equation}

\par Since (\ref{hypothesegen}) holds formally, we derive
 that:

\begin{equation}
\label{hithypo}
\begin{split}
&\Sum_{l\in L} a_l Hit_{\nabla\tau}[Xdiv_{i_1}\dots
Xdiv_{i_\alpha}C^{l,i_1\dots i_a}_{g}(\Omega_1,\dots
,\Omega_p,\phi_1,\dots ,\phi_u)]+
\\&\Sum_{j\in J} a_j
Hit_{\nabla\tau}[C^j_{g}(\Omega_1,\dots,\Omega_p,\phi_1,\dots ,
\phi_u)]=0,
\end{split}
\end{equation}
modulo complete contractions of length $\ge\sigma +u+2$.
\newline

{\it Terminology:} For any complete contraction
 $C_{g}$ in the form (\ref{form2}) or (\ref{form1}) with a weak character
 $Weak(\vec{\kappa}_{simp})$, we have
  assigned symbols $F_1,\dots,F_X$ to their factors
$\nabla^{(m)}R_{ijkl},S_{*}\nabla^{(\nu)}R_{ijkl},\nabla^{(p)}\Omega_h$
in $C_g$ which are contracting against some $\nabla\phi_h$.
 Now, regarding the complete contractions in
 $Image^1_{\phi_{u+1}}[C_{g}]$, we impose the following
 conventions:

 \par For each contraction of length $\sigma +u+1$ in
$Image^{1,+}_{\phi_{u+1}}[C_{g}]$ we will still speak of the
 factors $F_1,\dots ,F_X$. (They can be identified, since all the
 complete contractions in $Image^{1,+}_{\phi_{u+1}}[C_g]$ still
 have a $u$-simple character $\vec{\kappa}_{simp}$).
 On the other hand, for each
contraction of length $\sigma+u$ in
$Image^{1,+}_{\phi_{u+1}}[C_{g}]$ (with a factor
$\nabla^{(v+2)}\phi_{u+1}$ that has arisen from some
$F_a=\nabla^{(m)}R_{ijkl}$ or $F_a=S_{*}\nabla^{(\nu)}R_{ijkl}$),
we will speak of the factors $F_1,\dots ,F_X$, only now $F_a$ will
be the factor $\nabla^{(v+2)}\phi_{u+1}$.

\par We now separately consider the
sublinear combinations in
\\ $Image^1_{\phi_{u+1}}[Xdiv_{i_1}\dots Xdiv_{i_a}
C^{l,i_1\dots i_a}_{g} (\Omega_1,\dots,\Omega_p,\phi_1, \dots
,\phi_u)]$ and \\ $Image^{1}_{\phi_{u+1}} [ C^j_{g}
(\Omega_1,\dots,\Omega_p,\phi_1,\dots ,\phi_u)]$
 that belong to $Image^{1,+}_{\phi_{u+1}}[L_{g}]$. We denote
 these sublinear combinations by
 \\$Image^{1,+}_{\phi_{u+1}}
[Xdiv_{i_1}\dots Xdiv_{i_a} C^{l,i_1\dots i_a}_{g}
(\Omega_1,\dots,\Omega_p,\phi_1,\dots ,\phi_u)]$,
\\$Image^{1,+}_{\phi_{u+1}} [ C^j_{g}
(\Omega_1,\dots,\Omega_p,\phi_1,\dots ,\phi_u)]$.
Then, (\ref{tsabes}) can be re-expressed as:

\begin{equation}
\label{testacoli}
\begin{split}
&\Sum_{l\in L} a_l Image^{1,+}_{\phi_{u+1}}[Xdiv_{i_1}\dots
Xdiv_{i_a} C^{l,i_1\dots i_a}_{g}
(\Omega_1,\dots,\Omega_p,\phi_1,\dots ,\phi_u)] +
\\& \Sum_{j\in J}
a_j Image^{1,+}_{\phi_{u+1}}[C^j_{g}
(\Omega_1,\dots,\Omega_p,\phi_1,\dots ,\phi_u)]= \Sum_{z\in Z} a_z
C^z_{g}(\Omega_1,\dots,\Omega_p,\phi_1, \dots ,\phi_u),
\end{split}
\end{equation}
modulo complete contractions of length $\ge\sigma +u+2$.

\par Thus, we again define the operation
$Hit_{\nabla\tau}$ on the linear combinations
$$Image^{1,+}_{\phi_{u+1}} [Xdiv_{i_1}\dots Xdiv_{i_a}C^{l,i_1\dots
i_a}_{g} (\Omega_1,\dots,\Omega_p,\phi_1,\dots ,\phi_u)],$$
$Image^{1,+}_{\phi_{u+1}} [C^{j}_{g}
(\Omega_1,\dots,\Omega_p,\phi_1,\dots ,\phi_u)]$ as in
(\ref{tsira}), (\ref{tsira2}).

\par On the other hand, we can also consider
\\$Image^{1}_{\phi_{u+1}}\{ Hit_{\nabla\tau}[Xdiv_{i_1} \dots
Xdiv_{i_a}C^{l,i_1\dots i_a}_{g}
(\Omega_1,\dots,\Omega_p,\phi_1,\dots ,\phi_u)]\}$ and define
$Image^{1,+}_{\phi_{u+1}}\{ Hit_{\nabla\tau}[Xdiv_{i_1}\dots
Xdiv_{i_a}C^{l,i_1\dots i_a}_{g}
(\Omega_1,\dots,\Omega_p,\phi_1,\dots ,\phi_u)]\}$ in the same way
as above.

\par This leads us to consider the quantity:

\begin{equation}
\label{nigeria}
\begin{split}
&Difference[Xdiv_{i_1}\dots Xdiv_{i_a}C^{l,i_1\dots i_a}_{g}
(\Omega_1,\dots,\Omega_p,\phi_1,\dots ,\phi_u)]=
\\&Hit_{\nabla\tau}\{Image^{1,+}_{\phi_{u+1}}
[Xdiv_{i_1}\dots Xdiv_{i_a}C^{l,i_1\dots i_a}_{g}
(\Omega_1,\dots,\Omega_p,\phi_1,\dots ,\phi_u)]\}-
\\& Image^{1,+}_{\phi_{u+1}}\{Hit_{\nabla\tau}
[Xdiv_{i_1}\dots Xdiv_{i_a}C^{l,i_1\dots i_a}_{g}
(\Omega_1,\dots,\Omega_p,\phi_1,\dots ,\phi_u)]\}.
\end{split}
\end{equation}

\par Analogously, we define a quantity:

$$Difference[C^j_{g}(\Omega_1,\dots,\Omega_p,\phi_1,\dots ,\phi_u)],$$
for each $j\in J$.
\newline

\par Then, since (\ref{hypothesegen}) holds formally,
and just by virtue of the transformation laws (\ref{curvtrans}),
(\ref{levicivita}) and our definitions above we have that:

\begin{equation}
\label{megalosympe}
\begin{split}
&\Sum_{l\in L} a_l Difference[Xdiv_{i_1}\dots Xdiv_{i_a}
C^{l,i_1\dots i_a}_{g}(\Omega_1,\dots,\Omega_p,\phi_1, \dots
,\phi_u)]+
\\&\Sum_{j\in J} a_j Difference[C^j_{g}(\Omega_1,\dots,\Omega_p,
\phi_1,\dots ,\phi_u)]=0
\end{split}
\end{equation}
modulo complete contractions of length $\ge\sigma +u+3$. This
equation, suitably analyzed, will imply Lemma \ref{bjorn}.
\newline

\par Easily, we observe that both \\$Difference[Xdiv_{i_1}\dots
Xdiv_{i_a}C^{l,i_1\dots i_a}_{g}
(\Omega_1,\dots,\Omega_p,\phi_1,\dots ,\phi_u)]$ and
\\$Difference[C^j_{g}(\Omega_1,\dots,\Omega_p,\phi_1,\dots
,\phi_u)]$ consist of complete contractions with length
$\ge\sigma+u+2$, each with one factor $\nabla\phi_{u+1}$ and one
factor $\nabla\tau$. This is because the contractions with length
$\sigma+u+1$ in (\ref{nigeria}) will cancel perfectly (without
correction terms with more factors).
\newline

\par Now, we write out:

\begin{equation}
\label{martini1}
\begin{split}
&\Sum_{l\in L} a_l Difference[Xdiv_{i_1}\dots Xdiv_{i_a}
C^{l,i_1\dots i_a}_{g}(\Omega_1,\dots,\Omega_p,\phi_1, \dots
,\phi_u)]=
\\&\Sum_{t\in T} a_t C^t_{g}(\Omega_1,\dots,\Omega_p,
\phi_1,\dots ,\phi_u,\tau),
\end{split}
\end{equation}

\begin{equation}
\label{martini2}
\begin{split}
&\Sum_{j\in J} a_j
Difference[C^j_{g}(\Omega_1,\dots,\Omega_p,\phi_1, \dots
,\phi_u)]=
\\&\Sum_{y\in Y} a_y C^y_{g}(\Omega_1,\dots,\Omega_p,
\phi_1,\dots ,\phi_u,\tau).
\end{split}
\end{equation}

We then consider the sets of complete contractions $C^t$ and
 $C^y$ for which the factor
$\nabla\tau$ is contracting against a given factor $F_h$, $1\le
h\le X$.\footnote{Recall that $\{F_h\}_{1\le h\le X}$ is the set of
all the non-selected factors that are contracting against some
factor $\nabla\phi_v, v\le u$.} We denote the index
 sets of those complete contractions by $Y^h, T^h$.  Of course,
 since (\ref{megalosympe}) holds formally, it follows that for each
$h\le X$:

\begin{equation}
\label{megalosympe2}
\begin{split}
&\Sum_{t\in T^h} a_t C^t_{g}(\Omega_1,\dots,\Omega_p, \phi_1,\dots
,\phi_u,\tau)+
\\&\Sum_{y\in Y^h} a_y
C^y_{g}(\Omega_1,\dots,\Omega_p, \phi_1,\dots ,\phi_u,\tau)=0,
\end{split}
\end{equation}
modulo complete contractions of length $\ge\sigma +u+2$.

\par Now, for each $h\in \{1,\dots, x\}$ we define $Set(h)=\{a_1,\dots ,
a_{r_h}\}$ as follows:
$\rho\in Set(h)$ if and only if $\nabla\phi_\rho$
 is contracting against $F_h$ in $\vec{\kappa}_{simp}$.
 For each $\rho\in Set(h)$, we define
$G^{\rho,\tau}[C^t_{g}(\Omega_1,\dots,\Omega_p, \phi_1,\dots
,\phi_u,\tau)]$, $G^{\rho,\tau}[C^y_{g}(\Omega_1,\dots,
\Omega_p,\phi_1,\dots ,\phi_u,\tau)]$ to stand for the complete
 contraction that arises
 by replacing the factors $\nabla_i\tau$, $\nabla_j\phi_\rho$ by a
 factor $g_{ij}$ (we thus obtain a complete
 contraction with length $\sigma+u$ and one internal contraction). We also define:

 \begin{equation}
 \label{safe}
 \begin{split}
&G\{\Sum_{l\in L} a_l Difference[Xdiv_{i_1}\dots Xdiv_{i_a}
C^{l,i_1\dots i_a}_{g}(\Omega_1,\dots,\Omega_p,\phi_1, \dots
,\phi_u)]+
\\&\Sum_{j\in J} a_j Difference[C^j_{g}(\Omega_1,\dots,\Omega_p,
\phi_1,\dots ,\phi_u)]\}=
\\&\Sum_{h=1}^X\Sum_{t\in T^h} a_t
\Sum_{\rho\in Set(h)}G^{\rho,\tau}\{
C^t_{g}(\Omega_1,\dots,\Omega_p,\phi_1,\dots ,\phi_u,\tau)\}+
\\&\Sum_{h=1}^X\Sum_{y\in Y^h} a_y
\Sum_{\rho\in Set(h)}G^{\rho,\tau}\{C^y_{g}(\Omega_1,
\dots,\Omega_p,\phi_1,\dots ,\phi_u,\tau)\}.
\end{split}
\end{equation}

\par Then, since (\ref{megalosympe}) holds formally, we derive
 that modulo complete contractions of length
$\ge\sigma +u+1$:

\begin{equation}
\label{kyriospontikos}
\begin{split}
&\Sum_{h=1}^X\Sum_{t\in T^h} a_t \Sum_{\rho\in
Set(h)}G^{\rho,\tau}\{
C^t_{g}(\Omega_1,\dots,\Omega_p,\phi_1,\dots,\phi_u,\tau)\}+
\\&\Sum_{h=1}^X\Sum_{y\in Y^h} a_y
\Sum_{\rho\in Set(h)}G^{\rho,\tau}\{C^y_{g}(\Omega_1,
\dots,\Omega_p,\phi_1,\dots ,\phi_u,\tau)\}=0.
\end{split}
\end{equation}

\par Now, we wish to explicitly understand how the left hand side of
(\ref{kyriospontikos}) arises from $\Sum_{l\in L} a_l
Xdiv_{i_1}\dots Xdiv_{i_a} C^{l,i_1\dots
i_a}_{g}(\Omega_1,\dots,\Omega_p,\phi_1, \dots ,\phi_u)$,
\\$\Sum_{j\in J} a_j C^j_{g}(\Omega_1,\dots,\Omega_p, \phi_1,\dots
,\phi_u)$. We will find that the left hand sides of these
equations are closely related to the LHS of the equation in Lemma
\ref{bjorn} (which we are trying to prove).

\par We introduce some notation in order to accomplish this.
\newline

\begin{definition}
\label{dsharp} We  define an operation $G^\sharp$ that formally
acts on the complete contractions in the linear combinations
$$LC_\Phi[Xdiv_{i_1}\dots Xdiv_{i_a} C^{l,i_1\dots
i_a}_{g}(\Omega_1,\dots,\Omega_p,\phi_1, \dots ,\phi_u)],LC_\Phi[C^f_{g}(\Omega_1,\dots,\Omega_p, \phi_1,\dots ,\phi_u)]:$$

\par Firstly, we define an operation $G^{\sharp, h, a_b}$ for
 each $h, 1\le h\le X$ and each $a_b\in Set(h)=
\{a_1,\dots ,a_{r_h}\}$. $G^{\sharp, h, a_b}$ formally acts by
erasing the factor $\nabla\phi_{a_b}$ that contracts against the
factor $F_h$ (say against an index ${}_\chi$) and then by hitting
the factor $F_h$ by a derivative $\nabla^\chi$ (so we again obtain
a complete contraction with an internal
 contraction and with length $\sigma+u$). Then, we define:

\begin{equation}
\label{tromos}
\begin{split}
&G^\sharp \{ LC_\Phi[Xdiv_{i_1}\dots Xdiv_{i_a} C^{l,i_1\dots
i_a}_{g}(\Omega_1,\dots,\Omega_p,\phi_1, \dots ,\phi_u)]\}
\\&=\Sum_{h=1}^X\Sum_{a_b\in Set(h)=\{a_1,\dots
,a_{r_h}\}}G^{\sharp, h,a_b}\{LC_\Phi[Xdiv_{i_1}\dots Xdiv_{i_a}
\\&C^{l,i_1\dots i_a}_{g}(\Omega_1,\dots,\Omega_p,\phi_1, \dots
,\phi_u)]\},
\end{split}
\end{equation}

\begin{equation}
\label{tromos2}
\begin{split}
& G^\sharp \{LC_\Phi [C^j_{g}(\Omega_1,\dots,\Omega_p,\phi_1,
\dots ,\phi_u)]\}
\\&=\Sum_{h=1}^X\Sum_{a_b\in Set(h)=\{a_1,\dots
,a_{r_h}\}}G^{\sharp,
h,a_b}\{LC_\Phi[C^j_{g}(\Omega_1,\dots,\Omega_p, \phi_1,
\dots,\phi_u)]\}.
\end{split}
\end{equation}
\end{definition}

\par Now, for our next definition we will treat each
factor $S_{*}\nabla^{(\nu)}R_{ijkl}$ as a sum of factors in the
form $\nabla^{(\nu)}R_{ijkl}$. We do this for the purpose of
picking out the factors $\nabla\phi_h$ that are contracting
against internal indices in each of the summands
$\nabla^{(\nu)}R_{ijkl}$ of the un-symmetrization of
$S_{*}\nabla^{(\nu)}R_{ijkl}$. Thus, we are now considering tensor
fields and
 complete contractions in the form (\ref{form1}) in paper \cite{alexakis3}.
 For each tensor field $C^{l,i_1\dots i_a}_{g}(\Omega_1,\dots,\Omega_p,\phi_1,
\dots ,\phi_u)$ and also for each complete contraction
$C^j_{g}(\Omega_1,\dots, \Omega_p,\phi_1, \dots ,\phi_u)$ we
consider the factors $\nabla^{(m)}R_{ijkl}$ which have a factor
$\nabla\phi_h$ contracting against an internal index.
 For any such factor $F_h$, we denote by $\Pi_l(h)$ (or $\Pi_j(h)$) the set of numbers
 $\pi$ for which $\nabla\phi_\pi$
is contracting against an internal index in the factor $F_h$. We
will denote by $\Pi$ (or $\Pi_l, \Pi_f$ if we wish to be more
precise) the set of all numbers $\pi$ for which $\nabla\phi_\pi$
is contracting against an internal index in some curvature factor
$\nabla^{(m)}R_{ijkl}$. For the purposes of the next definition we
assume that if $\pi\in \Pi$ then the factor $\nabla\phi_\pi$ is
contracting against the index ${}_i$ of a factor
$\nabla^{(m)}R_{ijkl}$. This assumption can be made with no loss
of generality for the purposes of this proof. If it were
contracting
 against one of the indices ${}_j,{}_k,{}_l$, we can apply standard
identities to make it contract against the index ${}_i$.
\newline

 Now, for each $l\in L$ and each $\pi\in \Pi_l(h)$ and each
given complete contraction $C_{g}(\Omega_1, \dots,\Omega_p,\phi_1,
\dots ,\phi_u)$ in $Xdiv_{i_1}\dots Xdiv_{i_a}C^{l,i_1\dots
i_a}_{g}(\Omega_1, \dots,\Omega_p,\phi_1, \dots ,\phi_u)$, 
we define an operation $Oper^{\pi,k,h}_{\phi_{u+1}}[C_g ]$
that acts as follows: If the index ${}_k$ in $F_h$ is contracting
against the selected factor,
$Oper^{\pi,k,h}_{\phi_{u+1}} [C_{g}(\Omega_1,
\dots,\Omega_p,\phi_1, \dots ,\phi_u)]$ stands for the
 complete contraction that arises from
$C_{g}(\Omega_1,\dots,\Omega_p,\phi_1, \dots , \phi_u)$ by
replacing the expression $\nabla^{(m)}_{r_1\dots
r_m}R_{ijkl}\nabla^i {\phi}_\pi$ by 
$-\nabla_k\phi_{u+1}\nabla^{(m)}_{r_1\dots r_m}R_{ijal}
\nabla^i{\phi}_\pi\nabla^a\tau$. If ${}_k$ is not contracting
against the selected factor, we let
\\$Oper^{\pi,k,h}_{\phi_{u+1}}[C_{g}(\Omega_1,
\dots,\Omega_p,\phi_1, \dots ,\phi_u)]=0$. Similarly, we define
the operation $Oper^{\pi,l,h}_{\phi_{u+1}}[\dots]$.

\par Then, for each $j\in J$ and each $\pi\in \Pi_j(h)$,
we define an operation
\\$Oper^{\pi,k,h}_{\phi_{u+1}}[C^j_{g}
(\Omega_1,\dots,\Omega_p,\phi_1, \dots ,\phi_u)]$ that acts as
follows: If ${}_k$ (in $F_h$) is contracting against the selected
factor, we  define $Oper^{\pi,k,h}_{\phi_{u+1}}
[C^f_{g}(\Omega_1, \dots,\Omega_p,\phi_1, \dots ,\phi_u)]$, to
stand for the
 complete contraction that arises from
$C^j_{g}(\Omega_1,\dots,\Omega_p,\phi_1, \dots , \phi_u)$ by
replacing the expression $\nabla^{(m)}_{r_1\dots
r_m}R_{ijkl}\nabla^i\phi_f$ by an expression
\\$-\nabla_k\phi_{u+1}\nabla^{(m)}_{r_1\dots \dots r_m}R_{ijal}
\nabla^i\phi_\pi\nabla^a\tau$. If ${}_k$ is not contracting
against the selected factor, we define
\\$Oper^{\pi,k,h}_{\phi_{u+1}}[C^j_{g}(\Omega_1,
\dots,\Omega_p,\phi_1, \dots ,\phi_u)]=0$. Similarly, we define
the operation \\$Oper^{\pi,l,h}_{\phi_{u+1}}[
C^j_{g}(\Omega_1,\dots,\Omega_p,\phi_1, \dots ,\phi_u)]$. This
operation extends to linear combinations.

\par We then define, for each $l\in L$:

\begin{equation}
\label{lebanon}
\begin{split}
&Special_{\phi_{u+1}}[Xdiv_{i_1}\dots Xdiv_{i_a} C^{l,i_1\dots
i_a}_{g}(\Omega_1, \dots,\Omega_p,\phi_1, \dots ,\phi_u)]=
\\&\Sum_{\pi\in \Pi(l)}
G^{\tau, \pi}\{ Oper^{\pi,k,h}_{\phi_{u+1}}[Xdiv_{i_1}\dots
Xdiv_{i_a}C^{l,i_1\dots i_a}_{g}(\Omega_1, \dots,\Omega_p,\phi_1,
\dots ,\phi_u)]+
\\&Oper^{\pi,l,h}_{\phi_{u+1}}[Xdiv_{i_1}\dots
Xdiv_{i_a}C^{l,i_1\dots i_a}_{g} (\Omega_1,\dots,\Omega_p,\phi_1,
\dots ,\phi_u)]\},
\end{split}
\end{equation}

and also for each $j\in J$:

\begin{equation}
\label{lebanon2}
\begin{split}
&Special_{\phi_{u+1}}[C^j_{g}(\Omega_1, \dots,\Omega_p,\phi_1,
\dots ,\phi_u)]=
\\&\Sum_{\pi\in \Pi(f)}
G^{\tau, \pi}\{ Oper^{\pi,k,h}_{\phi_{u+1}}[C^j_{g}(\Omega_1,
\dots,\Omega_p,\phi_1, \dots ,\phi_u)]+
\\&Oper^{\pi,l,h}_{\phi_{u+1}}[C^f_{g}
(\Omega_1,\dots,\Omega_p,\phi_1, \dots ,\phi_u)]\}.
\end{split}
\end{equation}

\par Thus, by construction each tensor field or complete contraction in
\\ $Special_{\phi_{u+1}}[Xdiv_{i_1}\dots Xdiv_{i_a}C^{l,i_1\dots
i_a}_{g}(\Omega_1, \dots,\Omega_p,\phi_1, \dots ,\phi_u)]$ or
\\$Special_{\phi_{u+1}} [C^j_{g}(\Omega_1,\dots,\Omega_p,\phi_1,
\dots ,\phi_u)]$ has length $\sigma+u$ and an internal contraction
in a factor $\nabla^{(p)}Ric_{ik}$. Moreover, we observe that each
complete contraction in $Special_{\phi_{u+1}}
[C^j_{g}(\Omega_1,\dots,\Omega_p, \phi_1,\dots ,\phi_u)]$ will
have at least one factor $\nabla\phi_b$, $b\in
Def(\vec{\kappa}_{simp})$\footnote{Recall that $Def(\vec{\kappa}_{simp})$
stands for the the set of numbers $h$ for which some factor
$\nabla\tilde{\phi}_h$ is contracting against
 the index ${}_i$ in some $S_{*}\nabla^{(\nu)} R_{ijkl}$.}
contracting against a derivative index of a factor
$\nabla^{(m)}R_{ijkl}$ or  $\nabla^{(p)}Ric$.
\newline

\par We then claim:

\begin{lemma}
\label{sigkapourh}

\par In the notation above, we claim that for each
$h, 1\le h\le X$:
\begin{equation}
\label{armenia}
\begin{split}
&(0=)\Sum_{t\in T^h} a_t \Sum_{\rho\in Set(h)}G^{\rho,\tau}\{
C^t_{g}(\Omega_1,\dots,\Omega_p,\phi_1,\dots ,\phi_u,\tau)\}+
\\&\Sum_{y\in Y^h} a_y
\Sum_{\rho\in Set(h)}G^{\rho,\tau}\{
C^y_{g}(\Omega_1,\dots,\Omega_p, \phi_1,\dots ,\phi_u,\tau)\}=
\\&\Sum_{l\in L} a_l G^\sharp \{LC_\Phi[Xdiv_{i_1}\dots Xdiv_{i_a}
C^{l,i_1\dots i_a}_{g}(\Omega_1,\dots,\Omega_p,\phi_1, \dots ,
\phi_u)]\}
\\& +Special_{\phi_{u+1}}[Xdiv_{i_1}\dots
Xdiv_{i_a}C^{l,i_1\dots i_a}_{g} (\Omega_1,\dots,\Omega_p,\phi_1,
\dots ,\phi_u)]\}+
\\& \Sum_{j\in J} a_j G^\sharp \{LC_\Phi[ C^j_{g}(\Omega_1,\dots,
\Omega_p,\phi_1, \dots ,\phi_u)]\}
\\&+Special_{\phi_{u+1}}[C^j_{g}(\Omega_1,\dots,\Omega_p,\phi_1, \dots
,\phi_u)].
\end{split}
\end{equation}
\end{lemma}

\par Thus, in view of equations (\ref{martini1}), (\ref{martini2})
the above Lemma in some sense provides us with information on the
linear combination

\begin{equation}
\label{pontik}
\begin{split}
&\Sum_{l\in L} a_l LC_\Phi[Xdiv_{i_1}\dots Xdiv_{i_a}
C^{l,i_1\dots i_a}_{g}(\Omega_1,\dots,\Omega_p,\phi_1, \dots
,\phi_u)]+
\\& \Sum_{j\in J} a_j LC_\Phi[ C^j_{g}(\Omega_1,\dots,
\Omega_p,\phi_1, \dots ,\phi_u)].
\end{split}
\end{equation}

{\it Proof of Lemma \ref{sigkapourh}:} We see our claim by the
definition of the operation $Difference[\dots]$, by book-keeping
and also from the definitions of the operations $G$ and
$G^\sharp$.

\par (Brief discussion:) Notice that for any complete contraction $C_g$ 
in
\\ $Xdiv_{i_1},\dots Xdiv_{i_a} C^{l,i_1\dots i_a}_g$ or $C^j_g$,
the complete contractions appearing in \\$Difference[C_g]$ arise
exclusively by applying the third term in the transformation law (\ref{levicivita})
 to two indices $(\nabla_s,{}_z)$\footnote{(where $\nabla_s$ is contracting
  against the factor $\nabla^s\tau$ that we have introduced--i.e.~provided that ${}_z$
 contracts against the selected factor in $C_g$)} in some factor $F_k$ in $C_g$
 (provided the factor $\nabla\phi_{u+1}$ that we then
 introduce contracts against the selected factor).
We then obtain a complete contraction $C'_g(\nabla\phi_{u+1})$
 that arises from $C_g$ by replacing the
expression $(\nabla_s,{}_z)\nabla^s\tau$ by an
expression $-\nabla_z\phi_{u+1}({}_s,\nabla^s\tau)$ (here the
 factor $\nabla^s\tau$ is now contracting against
  the ``position'' that the index ${}_z$ occupied).
Furthermore,  consider any factor $\nabla\phi_\rho$ that is
contracting against the factor $F_k$ and denote 
by ${}_b$ the index against which it contracts.

\par Now consider the case where one of the indices
 ${}_z,{}_b$ described in the previous paragraph
are derivative indices. We then observe that if we replace two
factors $\nabla_\alpha\phi_\rho,\nabla_\beta\tau$
 in $C'_g(\nabla\phi_{u+1})$ by $g_{\alpha\beta}$,
 we obtain precisely the complete contraction that arises
 in $LC_\Phi[C_g]$ when we apply the second or third term of the transformation law
 (\ref{levicivita}) to the
 indices ${}_z,{}_b$ and then replace the factor $\nabla\phi_\rho$
 by an internal contraction (as described in the definition
  of the operation $G^\sharp$). It is easy to see that this
gives a one-to-one correspondence between the LHS of (\ref{armenia}) and
 the terms $G^\sharp[\dots]$ in the RHS of (\ref{armenia}),
{\it except} for the terms arising in the LHS when both the indices
 ${}_b,{}_z$ discussed above are internal indices.
The same ``book-keeping'' reasoning then shows that 
the terms we obtain from the LHS of (\ref{armenia}) correspond
 to the terms $Special_{\phi_{u+1}}[\dots]$ in the RHS of (\ref{armenia}). $\Box$
\newline

\par The goal of our next Lemma will be to ``get rid'' of the
 sublinear combination

\begin{equation}
\label{megustas}
\begin{split}
&\Sum_{l\in L} a_l Special_{\phi_{u+1}}[Xdiv_{i_1}\dots
Xdiv_{i_a}C^{l,i_1\dots i_a}_{g} (\Omega_1,\dots,\Omega_p,\phi_1,
\dots ,\phi_u)]
\\&+\Sum_{j\in J} a_j
Special_{\phi_{u+1}}[C^j_{g}(\Omega_1,\dots,\Omega_p,\phi_1, \dots
,\phi_u)]
\end{split}
\end{equation}
in (\ref{armenia}),  and to replace it by a new linear combination
which will have one internal contraction involving a derivative
index (rather than an internal contraction in a  factor
$\nabla^{(p)}Ric$). To state our next Lemma, we will be using
tensor fields and
 complete contractions with an internal contraction in a factor
  $\nabla^{(m)}R_{ijkl}$, where that internal contraction
 will involve a derivative index. It will be useful to recall
 the operation $Sub_\omega$ from the Appendix in \cite{alexakis1} and the
 discussion directly above it.

In particular, we will show the following:

\begin{lemma}
\label{toxeri} In the notation of equation (\ref{armenia}),
 we claim that we can write:

\begin{equation}
\label{michalo}
\begin{split}
&\Sum_{l\in L} a_l Special_{\phi_{u+1}} [Xdiv_{i_1}\dots
Xdiv_{i_a} C^{l,i_1\dots i_a}_{g}(\Omega_1, \dots,\Omega_p,\phi_1,
\dots ,\phi_u)]
\\&+\Sum_{j\in J} a_j
Special_{\phi_{u+1}} [C^j_{g}(\Omega_1, \dots,\Omega_p,\phi_1,
\dots ,\phi_u)]=
\\&\Sum_{t\in T} a_t Xdiv_{i_1}\dots Xdiv_{i_a}
C^{t,i_1\dots i_a}_{g}(\Omega_1, \dots, \Omega_p,\phi_1,
\dots,\phi_{u+1})\\&+\Sum_{j\in J} a_j
C^j_{g}(\Omega_1,\dots,\Omega_p,\phi_1, \dots ,\phi_{u+1}),
\end{split}
\end{equation}
where the complete contractions and tensor fields on the right
hand side have an internal contraction in exactly one factor
$\nabla^{(m)}R_{ijkl}$ and that internal contraction involves a
 derivative index. Furthermore, in each tensor field and each complete
  contraction in the above equation precisely one of the factors
$\nabla\phi_x, x=1\dots ,u$ is missing. We
 accordingly denote by $T^x, J^x$ the corresponding index
sets of the complete contractions and tensor fields where
$\nabla\phi_x$ is missing.

\par We furthermore claim that in (\ref{michalo})
 each of the tensor fields
$Sub_{\phi_x}\{ C^{t,i_1\dots i_a}_{g}\}$ has the following
properties:
 In the case of Lemma \ref{zetajones} and \ref{pool2} it is acceptable and
moreover has a $(u+1)$-simple character $\vec{\kappa}^{+}_{simp}$.
In the case of Lemma \ref{pskovb},
$$\Sum_{t\in T^x} a_t Sub_{\phi_x}\{
Xdiv_{i_1}\dots Xdiv_{i_a} C^{t,i_1\dots i_a}_{g}(\Omega_1,
\dots,\Omega_p,\phi_1, \dots,\phi_{u+1})\}$$
 will be a generic linear combination of complete contractions like the ones indexed in
$T_1\bigcup T_2$ (in the conclusion of that Lemma).

\par On the other hand, for each $j\in J^x$,
 $Sub_{\phi_x}[C^j_{g}(\Omega_1,
\dots,\Omega_p,\phi_1, \dots,\phi_{u+1})]$ is simply subsequent to
$\vec{\kappa}^{+}_{simp}$.
\end{lemma}

\par Let us show how proving Lemma \ref{toxeri}, in conjunction
 with equation (\ref{armenia}) (which we have already proven),
 would imply our Lemma \ref{bjorn}.
\newline

{\it Proof that Lemma \ref{toxeri} implies Lemma \ref{bjorn}:}
\newline

\par By virtue of the Lemma \ref{toxeri} we may refer to
(\ref{armenia}) and replace:

$$\Sum_{l\in L} a_l Special_{\phi_{u+1}}[Xdiv_{i_1}\dots Xdiv_{i_a}
C^{l,i_1\dots i_a}_{g}(\Omega_1, \dots,\Omega_p,\phi_1, \dots
,\phi_u)]$$ as in the conclusion of Lemma \ref{toxeri}.
\newline

\par We then obtain an equation:

\begin{equation}
\label{joseph}
\begin{split}
&0=\Sum_{l\in L} a_l G^\sharp \{LC_\Phi[Xdiv_{i_1}\dots Xdiv_{i_a}
C^{l,i_1\dots i_a}_{g}(\Omega_1,\dots,\Omega_p,\phi_1, \dots ,
\phi_u)]\}+
\\& \Sum_{j\in J} a_j G^\sharp \{LC_\Phi[ C^j_{g}
(\Omega_1,\dots,\Omega_p,\phi_1, \dots ,\phi_u)]\}+
\\& \Sum_{t\in T} a_t Xdiv_{i_1}\dots
Xdiv_{i_a} C^{t,i_1\dots i_a}_{g}(\Omega_1, \dots,\Omega_p,\phi_1,
\dots ,\phi_{u+1})
\\& +\Sum_{j\in J} a_j
C^j_{g}(\Omega_1,\dots,\Omega_p,\phi_1, \dots ,\phi_{u+1}),
\end{split}
\end{equation}
where the linear combinations indexed in $T,J$ are as in
 (\ref{michalo}).

We can then derive Lemma \ref{bjorn} straightforwardly: We break
the above equation into sublinear combinations according to the
factor $\nabla\phi_x$ that is missing. Each of these sublinear
combinations must vanish separately, since the above holds
formally. We denote the respective sublinar combinations in each
$G^\sharp\{\dots\}$ by $G^{\sharp, x}\{\dots\}$. Thus we have that
for each $x$:

\begin{equation}
\label{josephx}
\begin{split}
&0=\Sum_{l\in L} a_l G^{\sharp,x} \{LC_\Phi[Xdiv_{i_1}\dots
Xdiv_{i_a} C^{l,i_1\dots i_a}_{g}(\Omega_1,\dots,\Omega_p,\phi_1,
\dots , \phi_u)]\}+
\\& \Sum_{j\in J} a_j G^{\sharp,x} \{LC_\Phi[ C^j_{g}
(\Omega_1,\dots,\Omega_p,\phi_1, \dots ,\phi_u)]\}+
\\& \Sum_{t\in T^x} a_t Xdiv_{i_1}\dots
Xdiv_{i_a} C^{t,i_1\dots i_a}_{g}(\Omega_1, \dots,\Omega_p,\phi_1,
\dots ,\phi_{u+1})
\\& +\Sum_{j\in J^x} a_j
C^j_{g}(\Omega_1,\dots,\Omega_p,\phi_1, \dots ,\phi_{u+1}).
\end{split}
\end{equation}

 We then use the operation $Sub_{\phi_x}$
(which acts on each complete
 contraction $C_{g}$ in the right hand side of
(\ref{joseph}) by picking out the one factor
 \\$\nabla^{r_a}\nabla^{(m)}_{r_1\dots r_m}
R_{r_{m+1}\dots r_{m+4}}$ or $\nabla^{r_a}
\nabla^{(p)}_{r_1\dots r_p}\Omega_h$ with the
 internal contraction,\footnote{The
 indices $({\nabla}^{r_a},{}_{r_a})$ are contracting against
each other.} and replacing it by an expression
$\nabla^{(m)}_{r_1\dots r_m}R_{r_{m+1}\dots r_{m+4}}
\nabla^{r_a}\phi_x$, $\nabla^{(p)}_{r_1\dots r_p}\Omega_h
\nabla^{r_a}\phi_x$).

\par Now, by just keeping track of the definitions above we find that:

\begin{equation}
\label{sthseira}
\begin{split}
&\sum_{x=1}^u Sub_{\phi_x}\{G^{\sharp,x} \{LC_\Phi[Xdiv_{i_1}\dots
Xdiv_{i_a} C^{l,i_1\dots i_a}_{g}(\Omega_1,\dots,\Omega_p,\phi_1,
\dots ,\phi_u)]\}\}=
\\& LC_\Phi[Xdiv_{i_1}\dots Xdiv_{i_a}
C^{l,i_1\dots i_a}_{g}(\Omega_1,\dots,\Omega_p,\phi_1, \dots
,\phi_u)],
\end{split}
\end{equation}
and also, for every $j\in J$:

\begin{equation}
\label{sthseira}
\begin{split}
&\sum_{x=1}^u Sub_{\phi_x}\{G^{\sharp,x}
\{LC_\Phi[C^j_{g}(\Omega_1,\dots,\Omega_p, \phi_1, \dots
,\phi_u)]\}=
\\&LC_\Phi[C^j_{g}(\Omega_1,\dots,\Omega_p,\phi_1,
\dots ,\phi_u)].
\end{split}
\end{equation}

 Hence, acting by $\Sum_x Sub_{\phi_x}$ on (\ref{joseph}), we
 deduce our Lemma \ref{bjorn}. $\Box$
\newline

\par Therefore, if we can show Lemma \ref{toxeri},  we can then
derive Lemma \ref{bjorn}.
\newline

{\it  Proof of Lemma \ref{toxeri}:}
\newline

\par Firstly we make a few observations: The complete contractions and tensor
fields in the sublinear combinations $Special_{\phi_{u+1}}[\dots]$
are each in the form:

\begin{equation}
\label{aoustria}
\begin{split}
&pcontr(\nabla^{(m_1)}R_{ijkl}\otimes\dots\otimes\nabla^{(m_{\sigma_1})}
R_{ijkl}\otimes \nabla^{(b)} Ric_{ij}\otimes
S_{*}\nabla^{(\nu_1)}R_{ijkl}\otimes
\\&\dots\otimes
S_{*}\nabla^{(\nu_t)} R_{ijkl}
 \otimes\nabla^{(b_1)}\Omega_1\otimes\dots\otimes \nabla^{(b_p)}\Omega_p\otimes
\\& \nabla\phi_{z_1}\dots \otimes\nabla\phi_{z_u}\otimes\nabla
\phi'_{z_{u+1}}\otimes
\dots\otimes\nabla\phi'_{z_{u+d}}\otimes\dots \otimes
\nabla\tilde{\phi}_{z_{u+d+1}}\otimes\dots\otimes\nabla
\tilde{\phi}_{z_{u+d+y}}),
\end{split}
\end{equation}
where for each of the complete contractions and tensor fields in
the above form one of the factors
$\nabla\phi_1,\dots,\nabla\phi_u$ is missing (the factor
$\nabla\phi_x$ in the notation of Lemma \ref{toxeri}).
Accordingly, we will re-express the LHS of our Lemma hypothesis as:

\begin{equation}
\label{theform}
\begin{split}
&\Sum_{l\in L} a_l Special_{\phi_{u+1}} [Xdiv_{i_1}\dots
Xdiv_{i_a} C^{l,i_1\dots i_a}_{g}(\Omega_1, \dots,\Omega_p,\phi_1,
\dots ,\phi_u)]
\\&+\Sum_{j\in J} a_j
Special_{\phi_{u+1}} [C^j_{g}(\Omega_1, \dots,\Omega_p,\phi_1,
\dots ,\phi_u)]=
\\&\Sum_{x=1}^X \{\Sum_{l\in L'_x} a_l Xdiv_{i_1}\dots Xdiv_{i_c}
C^{l,i_1\dots i_c,i_b}_{g}(\Omega_1,\dots,\Omega_{p},\phi_1, \dots
, \phi_u)\nabla_{i_b}\phi_{u+1}+
\\&\Sum_{j\in J'_x} a_j C^j_{g}(\Omega_1,
\dots,\Omega_{p},\phi_1, \dots ,\phi_{u+1})\},
\end{split}
\end{equation}
where complete contractions and tensor fields indexed in
$L_x',J'_x$ have the factor $\nabla\phi_x$ missing. (We will
explain momentarily how $c$ is related to $a$).

\par We will then prove the assertion of Lemma \ref{toxeri} for
each of the sublinear combinations indexed in $L'_x,J'_x$
separately. Clearly, just adding all those equations will then
show our whole claim.

\par Recall (from the Appendix in \cite{alexakis1}) the operation
$Ricto\Omega$ (also denoted by $UnRic$) which replaces the factor
$\nabla^{(p)}_{r_1\dots r_p}Ric_{ik}$ by a factor
$\nabla^{(p+2)}_{r_1\dots r_p ik}\Omega_{p+1}$. In view of
(\ref{armenia})  we derive that for every $x\le u$:

\begin{equation}
\label{istoria2}
\begin{split}
&\Sum_{l\in L'_x} a_l Xdiv_{i_1}\dots Xdiv_{i_c} C^{l,i_1\dots
i_c,i_b}_{g}(\Omega_1,\dots,\Omega_{p+1},\phi_1, \dots ,
\phi_u)\nabla_{i_b}\phi_{u+1}+
\\&\Sum_{j\in J'_x} a_j C^j_{g}(\Omega_1,
\dots,\Omega_{p+1},\phi_1, \dots ,\phi_{u+1})=0,
\end{split}
\end{equation}
modulo complete contractions of length $\ge\sigma +u+2$. This
holds because (\ref{theform}) holds formally.

\par We will be using the equations (\ref{istoria2}) to derive our
claim. We have two different proofs based on $x$: Either $x\in
Def(\vec{\kappa}_{simp})$ or $x\notin
Def(\vec{\kappa}_{simp})$.\footnote{Recall that
$Def(\vec{\kappa}_{simp})$ stands for the set of numbers $o$ for
which some $\nabla\tilde{\phi}_o$ is contracting against the index
${}_i$ in some factor $S_{*}\nabla^{(\nu)}R_{ijkl}$.} We start
with the case where $x\in Def(\vec{\kappa}_{simp})$.
\newline

{\it Proof of Lemma \ref{toxeri} in the case $x\in
Def(\vec{\kappa}_{simp})$:}
\newline

\par In this case
we will derive our claim in two steps. Firstly, (in the notation
of (\ref{theform})) we claim that we can write:

\begin{equation}
\label{mousikh}
\begin{split}
&\Sum_{l\in L'_x} a_l Xdiv_{i_1}\dots Xdiv_{i_c} C^{l,i_1\dots
i_c,i_b}_{g}(\Omega_1,\dots,\Omega_{p},\phi_1, \dots
,\phi_u)\nabla_{i_b}\phi_{u+1}=
\\&\Sum_{t\in T} a_t Xdiv_{i_1}\dots Xdiv_{i_a}C^{t,i_1\dots i_a}_{g}
(\Omega_1,\dots,\Omega_{p},\phi_1, \dots ,\phi_{u+1})+
\\&\Sum_{j\in J} a_jC^j_{g} (\Omega_1,\dots,\Omega_{p},\phi_1,
\dots ,\phi_{u+1})+\Sum_{j\in J'} a_jC^j_{g}
(\Omega_1,\dots,\Omega_{p},\phi_1, \dots ,\phi_{u+1}),
\end{split}
\end{equation}
where the tensor fields and complete contractions indexed in $T,J$
are as in the conclusion of  Lemma \ref{toxeri}, while the
contractions indexed in $J'$ have exactly one factor
$\nabla^{(p)}Ric$ but also have one of the factors $\nabla\phi_y,
y\in Def(\vec{\kappa}_{simp})$\footnote{Recall that
$Def(\vec{\kappa}_{simp})$ stands for the set of numbers $\rho$
for which some factor $\nabla\phi_\rho$ is contracting against the
index ${}_i$ in a factor $S_{*}\nabla^{(\nu)}R_{ijkl}$.}
contracting against a derivative index in a factor
$\nabla^{(m)}R_{ijkl}$ {\it or, in the case of Lemmas}
\ref{zetajones}, \ref{pool2}, the factor $\nabla\phi_{u+1}$ is
contracting against a {\it derivative index} of the (one of the)
selected factor(s) $\nabla^{(m)}R_{ijkl}$.

\par If we can show (\ref{mousikh}),
 then by replacing the above into (\ref{armenia}) and using the notation
 of (\ref{theform}), we obtain a new equation:

\begin{equation}
\label{armenia2}
\begin{split}
&\Sum_{l\in L} a_l G^{\sharp,x} \{LC_\Phi[Xdiv_{i_1}\dots
Xdiv_{i_a} C^{l,i_1\dots i_a}_{g}(\Omega_1,\dots,\Omega_p,\phi_1,
\dots , \phi_u)]\}+
\\& \Sum_{j\in J} a_j G^{\sharp,x} \{LC_\Phi[ C^j_{g}(\Omega_1,\dots,
\Omega_p,\phi_1, \dots ,\phi_u)]\}+
\\&\Sum_{t\in T} a_t Xdiv_{i_1}\dots Xdiv_{i_a}C^{t,i_1\dots i_a}_{g}
(\Omega_1,\dots,\Omega_{p},\phi_1, \dots ,\phi_{u+1})+
\\&\Sum_{j\in J} a_jC^j_{g} (\Omega_1,\dots,\Omega_{p},\phi_1,
\dots ,\phi_{u+1})+
\\&\Sum_{j\in J'} a_jC^j_{g}
(\Omega_1,\dots,\Omega_{p},\phi_1, \dots ,\phi_{u+1})+ \Sum_{j\in
J'_x} a_j C^j_{g}(\Omega_1,\dots,\Omega_p,\phi_1,\dots
,\phi_{u+1})=0;
\end{split}
\end{equation}
here the  linear combination indexed in $J'$ is not generic
 notation, it is precisely the linear combination appearing in
 the right hand side of (\ref{mousikh}). Furthermore, we observe that
  the complete contractions belonging to the sublinear combination
$$\Sum_{j\in J'_x} a_j C^j_{g}(\Omega_1,\dots,\Omega_p,\phi_1,\dots ,\phi_u)$$
will have at least one of the factors $\nabla\phi_x, x\in
Def(\vec{\kappa}_{simp})$ contracting against a derivative
 index of some factor $\nabla^{(m)}R_{ijkl}$. This follows by the
definition of the operation $Special_{\phi_{u+1}}[\dots ]$.

\par Now, applying the operation $Ricto\Omega_{p+1}$ to the above we
derive an equation:

\begin{equation}
\label{armenia3}
\begin{split}
&\Sum_{j\in J'} a_jC^j_{g} (\Omega_1,\dots,\Omega_{p+1},\phi_1,
\dots ,\phi_{u+1})+ \Sum_{j\in J'_x} a_j
C^j_{g}(\Omega_1,\dots,\Omega_{p+1}, \phi_1,\dots ,\phi_{u+1})=0,
\end{split}
\end{equation}
where $C^j(\Omega_1,\dots,\Omega_{p+1},\phi_1, \dots ,\phi_{u+1})$
arises from $C^j_g(\Omega_1,\dots,\Omega_{p},\phi_1, \dots
,\phi_{u+1})$ by replacing the factor $\nabla^{(y)}Ric$ by
$\nabla^{(y+2)}\Omega_{p+1}$.

\par Since the above holds formally, we may repeat the permutations
 by which we make it vanish formally to the linear combination

\begin{equation}
\label{armenia4}
\begin{split}
&\Sum_{j\in J'} a_jC^j_{g} (\Omega_1,\dots,\Omega_{p},\phi_1,
\dots ,\phi_{u+1})+ \Sum_{j\in J'_x} a_j
C^j_{g}(\Omega_1,\dots,\Omega_{p}, \phi_1,\dots \\&,\phi_{u+1})=0,
\end{split}
\end{equation}
modulo introducing correction terms of length $\sigma +u+1$ by
virtue of the formula $\nabla_aRic_{bc}-\nabla_bRic_{ac}=
\nabla^lR_{ablc}$, and also correction terms of length
$\sigma+u+2$ (which we do not care about). Thus, we derive that we
can write:

\begin{equation}
\label{sion}
\begin{split}
& \Sum_{j\in J'} a_jC^j_{g} (\Omega_1,\dots,\Omega_{p},\phi_1,
\dots ,\phi_{u+1})+ \Sum_{j\in J'_x} a_j
C^j_{g}(\Omega_1,\dots,\Omega_{p}, \phi_1,\dots ,\phi_{u+1})
\\&=\Sum_{j\in J} a_jC^j_{g} (\Omega_1,\dots,\Omega_{p},
\phi_1, \dots,\phi_{u+1}).
\end{split}
\end{equation}
Therefore, in the case $x\in Def(\vec{\kappa}_{simp})$, matters
are reduced to showing (\ref{mousikh}).
\newline

{\it Proof of (\ref{mousikh}):}
\newline

\par Now, we see that since all $C^{l,i_1,\dots ,i_a}_{g}$,
$l\in L$ in (\ref{armenia}) have a given simple character
$\vec{\kappa}_{simp}$, it follows that all tensor fields
$C^{l,i_1\dots i_c,i_b}(\Omega_1,\dots,\Omega_{p+1},\phi_1, \dots
, \phi_u)$, $l\in L'_x$ in (\ref{mousikh}) will have the same
$(u-1)$-simple character (the
 one defined by the factors $\nabla\phi_1,\dots
,\hat{\nabla}\phi_x,\dots ,\nabla\phi_u$).\footnote{See the
notation in (\ref{theform}).} We denote this $(u-1)$-simple
character by $\vec{\kappa}_{-1}^x$. Moreover,
 each $C^j_{g}(\Omega_1,\dots,
\Omega_{p+1},\phi_1,\dots ,\phi_u)$, $j\in J'_x$ must be
$(u-1)$-subsequent to $\vec{\kappa}_{-1}^x$.

\par Now, an observation: In the language of the introduction in \cite{alexakis4}, 
if (\ref{hypothese2}) falls under case II or III (i.e.~if we are proving
Lemma \ref{pool2} or Lemma \ref{pskovb}) then for each $l\in
L_\mu$ we have $Def^{*}(l)=0$ (i.e. in the tensor fields
$C^{l,i_1\dots i_\mu}_{g}$ in (\ref{hypothese2}) have no special
free indices in any factor $S_{*}\nabla^{(\nu)}R_{ijkl}$).
 Thus, we see by construction that
in (\ref{theform}) each $c$ is $\ge\mu$. On the other hand, in the
setting of Lemma \ref{zetajones}, we have that for any $F_h$ in
the form $S_{*}\nabla^{(\nu)}R_{ijkl}$ for which ${}_k$ is a free
index then we may write out:

\begin{equation}
\label{sircar}
\begin{split}
&G^{\tau,\pi}\{Oper^{k,h}_{\phi_{u+1}}[Xdiv_{i_1} \dots
Xdiv_{i_a}C^{l,i_1\dots i_a}_{g} (\Omega_1,\dots,\Omega_p,\phi_1,
\dots ,\phi_{u})]\}\}=
\\&\Sum_{r\in R} a_r Xdiv_{i_1} \dots Xdiv_{i_{a-1}}C^{r,i_1\dots
i_{a-1}}_{g}(\Omega_1,\dots,\Omega_p,\phi_1, \dots ,\phi_{u+1}).
\end{split}
\end{equation}
Note that if $l\in L_\mu$ (i.e. if $a=\mu$) then the tensor fields
on the right hand side will have the factor $\nabla\phi_{u+1}$
contracting against a derivative index in the selected factor.
Thus, in the notation we introduced:

\begin{equation}
\begin{split} 
&\Sum_{r\in R} a_r Xdiv_{i_1}\dots Xdiv_{i_{a-1}}
C^{r,i_1\dots i_{a-1}}_{g}(\Omega_1,\dots,\Omega_p,\phi_1, \dots
,\phi_{u+1})\\&=\Sum_{j\in J'} a_j
C^j_{g}(\Omega_1,\dots,\Omega_p,\phi_1, \dots ,\phi_{u+1}).
\end{split}
\end{equation}

\par Therefore, after this observation, we may assume that all
the tensor fields indexed in each $L'_x$ in the equation
(\ref{theform}) have $c\ge \mu$. Now, we will prove our Lemma by
an induction:
 We will use the notation
$$\Sum_{t\in T} a_t Xdiv_{i_1}\dots Xdiv_{i_a}
C^{t,i_1\dots i_a}_{g}(\Omega_1, \dots,\Omega_p,\phi_1,
\dots,\phi_{u+1}), $$

$$\Sum_{j\in J\bigcup J'} a_j C^j_{g}(\Omega_1,\dots,\Omega_p,\phi_1,
\dots ,\phi_{u+1})$$
 to denote generic linear combinations as
described above.
 We assume that for some $m\ge \mu$:

\begin{equation}
\label{vaticanisti}
\begin{split}
&\Sum_{l\in L'_x} a_l Xdiv_{i_1}\dots Xdiv_{i_c} C^{l,i_1\dots
i_c,i_b}_{g}(\Omega_1, \dots,\Omega_p,\phi_1, \dots
,\phi_u)\nabla_{i_b}\phi_{u+1}=
\\&\Sum_{d\in D^{m}} c_d Xdiv_{i_1}\dots Xdiv_{i_{a_d}}
C^{d,i_1\dots i_{a_d}}_{g}(\Omega_1, \dots,\Omega_p,\phi_1, \dots
,\phi_u,\phi_{u+1})+
\\& \Sum_{t\in T} a_t Xdiv_{i_1}\dots
Xdiv_{i_a} C^{t,i_1\dots i_a}_{g}(\Omega_1, \dots,\Omega_p,\phi_1,
\dots ,\phi_{u+1})\\&+\Sum_{j\in J\bigcup J'} a_j
C^j_{g}(\Omega_1,\dots,\Omega_p,\phi_1, \dots ,\phi_{u+1}),
\end{split}
\end{equation}
where here the vector fields $C^{d,i_1\dots i_{a_d}}_{g}(\Omega_1,
\dots,\Omega_p,\phi_1, \dots ,\phi_u,\phi_{u+1})$ are in the form
(\ref{aoustria}) with a factor $\nabla^{(p)}Ric$ and each has
$a_d\ge m$. Moreover, they each have a $(u-1)$-simple character
$\vec{\kappa}^x_{-1}$. We
 then claim that we can write:

\begin{equation}
\label{vaticanisti2}
\begin{split}
&\Sum_{l\in L'_x} a_l Xdiv_{i_1}\dots Xdiv_{i_a} C^{l,i_1\dots
i_a,i_b}_{g}(\Omega_1, \dots,\Omega_p,\phi_1, \dots
,\phi_u)\nabla_{i_b}\phi_{u+1}=
\\&\Sum_{d\in D^{m+1}} a_d Xdiv_{i_1}\dots Xdiv_{i_{a_d}}
C^{d,i_1\dots i_{a_d}}_{g}(\Omega_1, \dots,\Omega_p,\phi_1, \dots
,\phi_u,\phi_{u+1})+
\\& \Sum_{t\in T} a_t Xdiv_{i_1}\dots
Xdiv_{i_a} C^{t,i_1\dots i_a}_{g}(\Omega_1, \dots,\Omega_p,\phi_1,
\dots ,\phi_{u+1})\\&+\Sum_{j\in J\bigcup J'} a_j
C^j_{g}(\Omega_1,\dots,\Omega_p,\phi_1, \dots ,\phi_{u+1}),
\end{split}
\end{equation}
with the same notational conventions as above.

\par Clearly, if we can show the above inductive statement then by
iterative repetition we will derive (\ref{mousikh}).
\newline

{\it Proof that (\ref{vaticanisti}) implies (\ref{vaticanisti2}):}
We observe that all the tensor fields in (\ref{vaticanisti})
 have the same $(u-1)$-simple character, which we have denoted by
$\vec{\kappa}_{-1}^{x}$.

\par Now, by applying $Ricto\Omega_{p+1}$ (see the Appendix in \cite{alexakis1}) 
to (\ref{vaticanisti}) and using (\ref{istoria2})
 we derive:

\begin{equation}
\label{vaticanisti3}
\begin{split}
&\Sum_{d\in D^{m}} c_d Xdiv_{i_1}\dots Xdiv_{i_{a_d}}
Ricto\Omega_{p+1}[C^{d,i_1\dots i_{a_d}}_{g}(\Omega_1,
\dots,\Omega_p,\phi_1, \dots ,\phi_u,\phi_{u+1})]+
\\&\Sum_{j\in J'_x} a_j C^j_{g}(\Omega_1,
\dots,\Omega_{p+1},\phi_1, \dots ,\phi_u,\phi_{u+1})\\&+ \Sum_{j\in
J'} a_j Ricto\Omega_{p+1}[C^j_{g}(\Omega_1,\dots,\Omega_p,\phi_1,
\dots , \phi_u)]=0,
\end{split}
\end{equation}
modulo complete contractions of length $\ge\sigma +u+1$.

\par We pick out the index set $D^{m,*}\subset D^{m}$ of those
tensor fields for which $a_d=m$. Then (apart from certain 
``forbidden cases'' which
we discuss in the ``Digression'' below), we apply the first claim of
Lemma 4.10 in \cite{alexakis1} to (\ref{vaticanisti3}),\footnote{Note that (\ref{vaticanisti3}) 
formally falls under the inductive assumption of 
this Lemma, because we replaced a curvature term
by a factor $\nabla^{(y)}\Omega_{p+1}$, hence we are in a case in
which case Corollary 1 in \cite{alexakis4} already holds, by our
inductive assumption. Observe that if $m=\mu$ there is no danger of falling 
under a ``forbidden case'' by weight considerations, since we are assuming that the equation in 
our Lemma assumption does not contain ``forbidden terms''. 
The possibility of  ``forbidden cases'' when $m>\mu$ will 
be treated in the ``Digression'' below.} we
 deduce that for some linear combination of $(m+1)$-tensor
 fields
with $(u-1)$-simple character $\vec{\kappa}_{-1}^{x}$ (say
$\Sum_{s\in S} a_s C^{s,i_1\dots i_{m+1}}_{g}(\Omega_1,
\dots,\Omega_{p+1},\phi_1, \dots ,\phi_u,\phi_{u+1})$), we
 will have that:

\begin{equation}
\label{ainte}
\begin{split}
&\Sum_{d\in D^{m,*}} a_d Ricto\Omega_{p+1}[C^{d,i_1\dots
i_m}_{g}(\Omega_1, \dots,\Omega_p,\phi_1, \dots
,\phi_u,\phi_{u+1})]\nabla_{i_1}\upsilon\dots\nabla_{i_m}\upsilon
\\&-Xdiv_{i_{m+1}}\Sum_{s\in S} a_s
C^{s,i_1\dots i_{m+1}}_{g}(\Omega_1, \dots,\Omega_p,\phi_1, \dots
,\phi_u,\phi_{u+1})\nabla_{i_1}\upsilon\dots\nabla_{i_m}\upsilon
\\&=\Sum_{j\in J'} a_j
Ricto\Omega_{p+1}[C^{j,i_1\dots
i_m}_{g}(\Omega_1,\dots,\Omega_p,\phi_1, \dots
,\phi_u)\nabla_{i_1}\upsilon\dots\nabla_{i_m}\upsilon];
\end{split}
\end{equation}
(we are using the same generic notational conventions as above-the
$m$-tensor fields $Ricto\Omega_{p+1}[C^{j,i_1\dots
i_m}_{g}(\Omega_1,\dots,\Omega_p,\phi_1, \dots ,\phi_u)]$ again
have a $(u-1)$-simple character that is simply subsequent to
$\vec{\kappa}_{-1}^{x}$). Therefore, since the above must hold
formally, we conclude that:

\begin{equation}
\label{ainte}
\begin{split}
&\Sum_{d\in D^{m,*}} a_d Xdiv_{i_1}\dots Xdiv_{i_m}C^{d,i_1\dots
i_m}_{g}(\Omega_1, \dots,\Omega_p,\phi_1, \dots
,\phi_u,\phi_{u+1})-
\\&Xdiv_{i_1}\dots Xdiv_{i_m}Xdiv_{i_{m+1}}\Sum_{s\in S} a_s
C^{s,i_1\dots i_{m+1}}_{g}(\Omega_1, \dots,\Omega_p,\phi_1, \dots
,\phi_u,\phi_{u+1})=
\\& \Sum_{t\in T} a_t Xdiv_{i_1}\dots
Xdiv_{i_a} C^{t,i_1\dots i_a}_{g}(\Omega_1, \dots,\Omega_p,\phi_1,
\dots ,\phi_{u+1})+
\\&\Sum_{j\in J\bigcup J'} a_j
C^j_{g}(\Omega_1,\dots,\Omega_p,\phi_1, \dots ,\phi_{u+1}),
\end{split}
\end{equation}
where $C^{s,i_1\dots i_{m+1}}_{g}(\Omega_1, \dots,\Omega_p,
\phi_1,\dots ,\phi_u,\phi_{u+1})$ arises from \\$C^{s,i_1\dots
i_{m+1}}_{g}(\Omega_1, \dots, \Omega_{p+1},\phi_1,\dots
,\phi_u,\phi_{u+1})$ by formally replacing
 the factor \\$\nabla^{(p)}_{r_1\dots r_p}\Omega_{p+1}$ by a factor
$\nabla^{(p-2)}_{r_1\dots r_{p-2}}Ric_{r_{p-1}r_p}$. Furthermore,
$\sum_{t\in T}\dots$ is a generic linear combination as described
after (\ref{mousikh}).

This is precisely our desired inductive step. Therefore we have
derived our claim in the case where $x\in
Def(\vec{\kappa}_{simp})$, except for the ``forbidden cases''
which we now discuss:
\newline

{\it Digression: The ``forbidden cases''.} As noted above, 
the only case where Lemma 4.10 in \cite{alexakis4} 
cannot be applied to (\ref{vaticanisti3}) (because it 
falls under a ``forbidden case'' of that Lemma) 
is when $m>\mu$. 
 So, in that case we derive from
(\ref{vaticanisti}) that $D^{m,*}=D^m$, and then applying the ``weak substitute'' of 
the fundamental Proposition 2.1 in \cite{alexakis4} we 
derive:\footnote{See the Appendix of \cite{alexakis4}.} 

\begin{equation}
\label{kinaelaion}
\begin{split}
&\Sum_{d\in D^{m,*}} a_d Xdiv_{i_1}\dots
Xdiv_{i_m}C^{d,i_1\dots i_m}_{g}(\dots,\phi_{u+1})=
\\& \Sum_{t\in T} a_t Xdiv_{i_1}\dots Xdiv_{i_{m-1}} 
C^{t,i_1\dots i_{m-1}}_{g}(\dots,\phi_{u+1})+
\\&\Sum_{j\in J} a_j C^j_{g}(\dots,\phi_{u+1})+
\Sum_{f\in F} a_f Xdiv_{i_1}\dots Xdiv_{i_{m-1}}C^{f,i_1\dots i_{m-1}}_{g}(\dots,\phi_{u+1}),
\end{split}
\end{equation}
where the terms indexed in $J$ are simply subsequent to
$\vec{\kappa}_{simp}$ and have a factor $Ric$, while the
terms in $F$ have the $\nabla\phi_{u+1}$
contracting against a non-special index (and both terms above
have a factor $Ric$).

 Then, replacing the
above into (\ref{vaticanisti}) and applying Lemma 4.10 in \cite{alexakis4}
(notice it can now be applied, since the factor $\nabla\phi_{u+1}$ is
contracting against a non-special index), we derive that:
$$\Sum_{f\in F} a_f C^{f,i_1\dots i_{m-1}}_{g}(\dots,\Omega_{p+1},\phi_{u+1})\nabla_{i_1}
\upsilon\dots \nabla_{i_{m-1}}\upsilon=0.$$
(Here $C^{f,i_1\dots i_{m-1}}_{g}(\dots,\Omega_{p+1},\phi_{u+1})$ arises from
$C^{f,i_1\dots i_{m-1}}_{g}(\dots,\phi_{u+1})$ 
by applying $Ricto\Omega_{p+1}$).

Thus, we may erase the terms $\sum_{f\in F}\dots$ in
(\ref{kinaelaion}); with that new feature, (\ref{kinaelaion}) is
precisely our desired equation (\ref{ainte}). $\Box$
\newline

{\it Proof of Lemma \ref{toxeri} in the case $x\notin
Def(\vec{\kappa}_{simp})$:}
\newline

We recall that the factor $\nabla\phi_x$ ($x\notin
Def(\vec{\kappa}_{simp})$) is contracting against a factor
$T^{*}(x)=\nabla^{(m)}R_{ijkl}$. We then distinguish two cases:
Either in $\vec{\kappa}_{simp}$ there is some other $h'\ne x$ with
$h'\in Def(\vec{\kappa}_{simp})$ so that $\nabla\phi_{h'}$ is
contracting against $T^{*}(x)$ in $\vec{\kappa}_{simp}$, or there
is no such factor. Another way of describing these two cases is
that in the first case the factor $T^{*}(x)=\nabla^{(m)}R_{ijkl}$
has arisen from the de-symmetrization of some factor
$S_{*}\nabla^{(\nu)} R_{ijkl}$ (for which the factor
$\nabla\phi_{h'}$ was contracting against the index ${}_i$), while
in the second case $\nabla^{(m)}R_{ijkl}$ corresponds to a factor
$\nabla^{(m)} R_{ijkl}$ in $\vec{\kappa}_{simp}$. The second
subcase is easier, so we will start with that one.
\newline

{\it Second subcase:} We observe that in this setting we must have
$L'_x=\emptyset$ (refer to (\ref{istoria2})). 
Moreover, for each $C^j_g$, $j\in J'_h$, we must have at least one
of the factors $\nabla\phi_w$, $w\in Def(\vec{\kappa}_{simp})$ 
contracting against a derivative index of some factor $F_b$ in
$C^j_{g}(\Omega_1, \dots,\Omega_p, \phi_1,\dots
,\hat{\phi}_x,\dots,\phi_u,\phi_{u+1})$, where in addition $F_b\ne
T^{*}(x)$. In view of these observations, it is enough to show that
in this second subcase:

\begin{equation}
\label{openheart}
\begin{split} &\Sum_{j\in J'_x} a_j C^j_{g}(\Omega_1,
\dots,\Omega_p, \phi_1,\dots
,\hat{\phi}_x,\dots,\phi_u,\phi_{u+1})\\&=\Sum_{j\in J} a_j
C^j_{g}(\Omega_1, \dots,\Omega_p,
\phi_1,\dots,\hat{\phi}_x,\dots,\phi_u,\phi_{u+1}),
\end{split}\end{equation}
modulo complete contractions of greater length. Clearly, that will
imply our claim for this second subcase since $L'_x=\emptyset$ in
(\ref{istoria2}). We derive this equation rather easily: By
(\ref{istoria2}) we have:

\begin{equation}
\label{perhaps}
 \Sum_{j\in J'_x} a_j C^j_{g}(\Omega_1,
\dots,\Omega_{p+1}, \phi_1,\dots
,\hat{\phi}_h,\dots,\phi_u,\phi_{u+1})=0,
\end{equation}
where we recall that for each $C^j_g$, $j\in J'_x$, one of the
factors $\nabla\phi_c$, $c\in Def(\vec{\kappa}_{simp})$ is
contracting against a derivative index of some curvature factor.

 Now, since (\ref{perhaps}) holds formally (at the linearized level), we
may repeat the permutations by which we make the left hand
 side vanish formally to the linear combination

$$\Sum_{j\in J'_x} a_j C^j_{g}(\Omega_1,
\dots,\Omega_{p}, \phi_1,\dots
,\hat{\phi}_h,\dots,\phi_u,\phi_{u+1}),$$ and derive the right hand
side in (\ref{openheart}) as correction terms.
\newline

{\it  First subcase:} Recall that we are
assuming $x\notin Def(\vec{\kappa}_{simp})$, and moreover the
factor $\nabla\phi_x$ is contracting against some factor
$T^{*}(x)$ in $\vec{\kappa}_{simp}$ in the form
$S_{*}\nabla^{(\nu)}R_{ijkl}$. Recall also that we are now
assuming that $\nabla\phi_x$ is {\it not} the factor that
contracts against the index ${}_i$ in $T^{*}(x)$ in $\vec{\kappa}_{simp}$. That factor is
$\nabla\phi_{h'}$.

\par In this first subcase we define $\Sum_{j\in J''} a_j \dots$
 to stand for a generic linear
combination of complete contractions in the form (\ref{aoustria})
 with the
factor $\nabla\phi_{h'}$ contracting against a {\it derivative}
index in a factor $\nabla^{(p)}Ric$. We also denote by
$$\Sum_{x\in X} a_x C^{x,i_1\dots i_a}_{g}(\Omega_1,
\dots,\Omega_{p}, \phi_1,\dots
,\hat{\phi}_h,\dots,\phi_u,\phi_{u+1})$$ a generic linear
combination of tensor fields in the form (\ref{aoustria}) with a
factor $Ric_{ij}$ (with no derivatives) where $\nabla\phi_{h'}$ is
contracting against the index ${}_i$ in that factor $Ric_{ij}$, and
where $a\ge\mu$. Finally, we define $\Sum_{j\in J'}a_j\dots$ to
stand for the same generic linear combination as in the previous
case.

\par It then follows that in this first subcase we can re-express the terms
 in the equation (\ref{theform}) as:

\begin{equation}
\label{fascinating} \begin{split} &\Sum_{l\in L'_x} a_l
Xdiv_{i_1}\dots Xdiv_{i_a} C^{l,i_1\dots i_a,i_b}_{g}(\Omega_1,
\dots,\Omega_p,\phi_1, \dots ,\phi_u)\nabla_{i_b}\phi_{u+1}+
\\&\Sum_{j\in J'_x} a_j C^j_{g}(\Omega_1,
\dots,\Omega_{p+1},\phi_1, \dots ,\phi_{u+1})=
\\& \Sum_{t\in T} a_t Xdiv_{i_1}\dots
Xdiv_{i_a} C^{t,i_1\dots i_a}_{g}(\Omega_1, \dots,\Omega_p,\phi_1,
\dots ,\phi_{u+1})+
\\&Xdiv_{i_1}\dots Xdiv_{i_a}\Sum_{x\in X} a_x C^{x,i_1\dots
i_a}_{g}(\Omega_1, \dots,\Omega_{p}, \phi_1,\dots
,\hat{\phi}_h,\dots,\phi_u,\phi_{u+1})+
\\&\Sum_{j\in J\bigcup J'\bigcup J''} a_j C^j_{g}(\Omega_1,
\dots,\Omega_p,\phi_1, \dots ,\phi_{u+1}).
\end{split}
\end{equation}
This just follows by the definitions and by applying the second
Bianchi identity to the factor $\nabla^{(p)}Ric$ if necessary (this can be done, because $p>0$).

\par We then prove Lemma \ref{toxeri} in this setting via
an inductive statement: Let us suppose that the minimum rank among
the tensor fields indexed in $X$ in (\ref{fascinating}) is $m\ge
\mu$ and the corresponding tensor fields are indexed in
$X^{m}\subset X$. We then claim that we can write:

\begin{equation}
\label{lightbrigade}
\begin{split}
&Xdiv_{i_1}\dots Xdiv_{i_m}\Sum_{x\in X^{m}} a_x C^{x,i_1\dots
i_m}_{g}(\Omega_1, \dots,\Omega_{p}, \phi_1,\dots
,\hat{\phi}_h,\dots,\phi_u,\phi_{u+1})=
\\&Xdiv_{i_1}\dots Xdiv_{i_{m+1}}\Sum_{x\in X^{m+1}} a_x C^{x,i_1\dots
i_{m+1}}_{g}(\Omega_1, \dots,\Omega_{p}, \phi_1,\dots
,\hat{\phi}_h,\dots,\phi_u,\phi_{u+1})+
\\& \Sum_{t\in T} a_t Xdiv_{i_1}\dots
Xdiv_{i_a} C^{t,i_1\dots i_a}_{g}(\Omega_1, \dots,\Omega_p,\phi_1,
\dots ,\phi_{u+1})+
\\&\Sum_{j\in J\bigcup J'\bigcup J''} a_j C^j_{g}(\Omega_1,
\dots,\Omega_p,\phi_1, \dots ,\phi_{u+1}),
\end{split}
\end{equation}
where for all the linear combinations on the right hand side
 we are using generic notation.

 \par We will show (\ref{lightbrigade}) momentarily. For the time being, we
 observe that if we can show (\ref{lightbrigade}) then by iterative repetition
 we are reduced to proving our claim for the first subcase above under the
extra assumption that $X=\emptyset$.

 \par Under this extra assumption we claim that the sum in
 $J'\bigcup J''$ in the above satisfies:

\begin{equation}
\label{biaggio} \begin{split} &\Sum_{j\in J'\bigcup J''} a_j
C^j_{g}(\Omega_1, \dots,\Omega_p,\phi_1, \dots ,\phi_{u+1})=
\\&\Sum_{j\in J} a_j C^j_{g}(\Omega_1, \dots,\Omega_p,\phi_1,
\dots ,\phi_{u+1}).
\end{split}
\end{equation}
Notice that proving the above two equations will complete the
proof of Lemma \ref{toxeri} in this subcase. We first prove
(\ref{biaggio}) (assuming we have shown (\ref{lightbrigade})) by
the usual argument:

Plug (\ref{lightbrigade}) into (\ref{fascinating}) and then apply
$Ricto\Omega_{p+1}$ to derive:
 $$\Sum_{j\in J'\bigcup J''} a_j
C^j_{g}(\Omega_1, \dots,\Omega_{p+1},\phi_1, \dots
,\phi_{u+1})=0.$$ Now, we use
 the fact that the resulting equation holds formally:
  We may arrange that the factor $\nabla\phi_{h'}$ is
 contracting against the first index in the factor
 $\nabla^{(p)}_{r_1\dots r_p}Ric_{ij}$, hence in the permutations by
 which we make the LHS of the above
 formally zero the first index is not moved. Thus we
see that the correction terms arising when we repeat those
permutations for
 $$\Sum_{j\in J'\bigcup J''} a_j
C^j_{g}(\Omega_1, \dots,\Omega_{p},\phi_1, \dots ,\phi_{u+1})$$
are indeed as in the right hand side of (\ref{biaggio}).
\newline

{\it Proof of (\ref{lightbrigade}):} First of all,
observe that if the factor $S_{*}R_{ijkl}\nabla^i\tilde{\phi}_h$ in
$\vec{\kappa}_{simp}$ is contracting against some factor $\nabla\phi_{h''}$
(in addition to the factors $\nabla\phi_h, \nabla\phi_{h'}$), then
(\ref{lightbrigade}) is obvious since then by definition $X^m=X=\emptyset$.
Thus, we may now assume that only the factors $\nabla\phi_x$ and
$\nabla\phi_{h'}$ are contracting against $S_{*}R_{ijkl}\nabla^i\tilde{\phi}_h$ in
$\vec{\kappa}_{simp}$.

\par In that setting, we firstly show (\ref{lightbrigade})
for $m=\mu$: We only have to refer
 equation (\ref{ainte}) when $m=\mu$,
and replace $\phi_{h'}$ by $\phi_x$: We derive an equation:

\begin{equation}
 \label{luchino}
\begin{split}
&Xdiv_{i_1}\dots Xdiv_{i_\mu}\Sum_{x\in X^{\mu}} a_x C^{x,i_1\dots
i_\mu}_{g}(\Omega_1, \dots,\Omega_{p}, \phi_1,\dots
,\hat{\phi}_h,\dots,\phi_u,\phi_{u+1})=
\\&Xdiv_{i_1}\dots Xdiv_{i_{\mu+1}}\Sum_{x\in \overline{X}^{\mu+1}} a_x C^{x,i_1\dots
i_{m+1}}_{g}(\Omega_1, \dots,\Omega_{p}, \phi_1,\dots
,\hat{\phi}_h,\dots,\phi_u,\phi_{u+1})+
\\& \Sum_{t\in T} a_t Xdiv_{i_1}\dots
Xdiv_{i_\mu} C^{t,i_1\dots i_\mu}_{g}(\Omega_1, \dots,\Omega_p,\phi_1,
\dots ,\phi_{u+1})+
\\&\Sum_{j\in J\bigcup J'\bigcup J''} a_j C^j_{g}(\Omega_1,
\dots,\Omega_p,\phi_1, \dots ,\phi_{u+1}),
\end{split}
\end{equation}
(using generic notation of the tensor fields in the RHS).

\par Now, we
consider the tensor fields indexed in $\overline{X}^{\mu+1}$ which have a free
index in the factor $Ric_{ab}\nabla^a\phi_{h'}$\footnote{In
other words, the index ${}_b$ is free.} and we ``forget'' the $Xdiv$ structure of $Xdiv_b$.
Therefore, we are reduced to proving (\ref{lightbrigade}) with
with two additionnal features if $m=\mu$: Firstly that the tensor fields  with rank $m=\mu$ {\it do not}
have a free index in the expression $Ric_{ab}\nabla^a\tilde{\phi}_{h'}$ and also
that for those tensor fields the index ${}_b$ in that expression
is contracting either against a non-special index in
some curvature factor or against some index in a factor
$\nabla^{(A)}\Omega_f$ with $A\ge 3$.\footnote{(The
last property follows since $\nabla^b$ has arisen by
``forgetting'' an $Xdiv$).}
\newline

\par Armed with this additionnal hypothesis for the case $m=\mu$,
 we will now show (\ref{lightbrigade}) for any $m\ge \mu$:

\par We apply the
operation $Ricto\Omega$ to the Lemma hypothesis (using the notation of
(\ref{fascinating})), and we pick out the sublinear combination
with an expression
$\nabla^{(2)}_{ij}\Omega_{p+1}\nabla^i\phi_{h'}$. It follows
 that this expression (which we denote by $E_{g}$)
vanishes separately. Thus, we derive an equation:

\begin{equation}
\label{banc}
\begin{split}
&E_{g}=X_{*}div_{i_1}\dots X_{*}div_{i_a}\Sum_{x\in X^m} a_x
C^{x,i_1\dots i_a}_{g}(\Omega_1, \dots,\Omega_{p+1}, \phi_1,\dots
,\hat{\phi}_x,\dots,\phi_u,\phi_{u+1})
\\& +\Sum_{j\in J'} a_j C^j_{g}(\Omega_1, \dots,
\Omega_{p+1},
\phi_1,\dots ,\hat{\phi}_h,\dots,\phi_u,\phi_{u+1})=0,
\end{split}
\end{equation}
where $X_{*}div_i$ stands for the sublinear combination in $Xdiv_i$
where $\nabla_i$ is in addition not allowed to hit the expression
 $\nabla^{(2)}_{ij}\Omega_{p+1}\nabla^i\phi_{h'}$.
 Now, we formally replace the expression
$\nabla^{(2)}_{ij}\Omega_{p+1}\nabla^i\phi_{h'}$ by a factor
$\nabla_jY$. We denote the resulting tensor fields and complete
contractions by:

$$C^{x,i_1\dots i_a}_{g}(\Omega_1, \dots,\Omega_{p},Y,
\phi_1,\dots ,\hat{\phi}_h,\dots,\phi_u,\phi_{u+1}),$$
$$C^j_{g}(\Omega_1, \dots,\Omega_{p},Y,
\phi_1,\dots ,\hat{\phi}_h,\dots,\phi_u,\phi_{u+1}).$$
Then, since the above equation holds formally we derive that:

\begin{equation}
\label{banc2}
\begin{split}
&X_{*}div_{i_1}\dots X_{*}div_{i_a}\Sum_{x\in X} a_x C^{x,i_1\dots
i_a}_{g}(\Omega_1, \dots,\Omega_{p},Y, \phi_1,\dots
,\hat{\phi}_h,\dots,\phi_u,\phi_{u+1})+
\\& \Sum_{j\in J'} a_j C^j_{g}(\Omega_1, \dots,
\Omega_{p},Y, \phi_1,\dots ,\hat{\phi}_h,\dots,\phi_u,
\phi_{u+1})=0.
\end{split}
\end{equation}
$X_{*}div_i$ in this setting stands  for the sublinear combination
in $Xdiv_i$ where $\nabla_i$ is not allowed to hit the factor
$\nabla Y$.

\par But then (subject to certain exceptions which we explain below)
applying 4.6 in \cite{alexakis4}),\footnote{In particular, the exceptions are when
there are tensor fields of minimum rank in (\ref{banc2}) that fall under
one of the ``forbidden cases'' of that Lemma. The derivation of (\ref{lightbrigade})
in that case will be discussed below.} (or Lemma 4.7 in \cite{alexakis4}
if $\sigma=3$)\footnote{Notice that by the
conventions above this Lemma can be applied since the tensor
fields of minimum rank do not have a free index in $\nabla Y$.}
to the above,\footnote{Equation
(\ref{banc2}) formally falls under the inductive assumptions of
these Lemmas, since we have reduced the weight.} we derive that
there is a linear combination of $(m+1)$-tensor fields,
$$\Sum_{q\in Q} a_q
C^{q,i_1\dots i_{m+1}}_{g}(\Omega_1, \dots,\Omega_{p},Y,
\phi_1,\dots ,\hat{\phi}_x,\dots,\phi_u,\phi_{u+1}),$$ just like
the ones indexed in $X^{m}$ only with another free index, so that:

\begin{equation}
\label{banc3}
\begin{split}
&\Sum_{x\in X} a_x C^{x,i_1\dots i_m}_{g}(\Omega_1,
\dots,\Omega_{p},Y, \phi_1,\dots
,\hat{\phi}_x,\dots,\phi_u,\phi_{u+1})
\nabla_{i_1}\upsilon\dots\nabla_{i_m}\upsilon-\Sum_{q\in Q} a_q
\\& 
X_{*}div_{i_{m+1}} C^{q,i_1\dots i_{m+1}}_{g}(\Omega_1,
\dots,\Omega_{p},Y, \phi_1,\dots
,\hat{\phi}_x,\dots,\phi_u,\phi_{u+1})
\nabla_{i_1}\upsilon\dots\nabla_{i_m}\upsilon
\\&+ \Sum_{j\in J'} a_j C^{j,i_1\dots i_m}_{g}
(\Omega_1, \dots,\Omega_{p},Y, \phi_1,\dots ,\hat{\phi}_x,
\dots,\phi_u,\phi_{u+1})\nabla_{i_1}\upsilon\dots
\nabla_{i_m}\upsilon=0.
\end{split}
\end{equation}

{\it The Exceptions:} In the exceptional cases, we
apply Lemma 4.10 in \cite{alexakis4} to (\ref{banc2}) to
derive (\ref{lightbrigade}) directly. (Notice that this Lemma can be applied since 
$m>\mu$ in this case; this is because of the additionnal hypothesis in the case $m=\mu$ which ensures that 
we do not fall under the forbidden case when $m=\mu$).
\newline

{\it Derivation of (\ref{lightbrigade}) from
(\ref{banc3}):\footnote{In the non-exceptional cases.}} Now,
formally replace $\nabla_aY$ by an expression
$Ric_{ia}\nabla^a\phi_{h'}$ in the above; we thus again obtain an
 equation:

\begin{equation}
\label{banc3'}
\begin{split}
&\Sum_{x\in X} a_x C^{x,i_1\dots i_m}_{g}(\Omega_1,
\dots,\Omega_{p}, \phi_1,\dots
,\hat{\phi}_x,\dots,\phi_u,\phi_{u+1})
\nabla_{i_1}\upsilon\dots\nabla_{i_m}\upsilon-
\\& \Sum_{q\in Q} a_q
X_{*}div_{i_{m+1}} C^{q,i_1\dots i_{m+1}}_{g}(\Omega_1,
\dots,\Omega_{p}, \phi_1,\dots
,\hat{\phi}_x,\dots,\phi_u,\phi_{u+1})
\nabla_{i_1}\upsilon\dots\nabla_{i_m}\upsilon
\\& +\Sum_{j\in J'} a_j C^{j,i_1\dots i_m}_{g}
(\Omega_1, \dots,\Omega_{p}, \phi_1,\dots ,\hat{\phi}_x,
\dots,\phi_u,\phi_{u+1})\nabla_{i_1}\upsilon\dots
\nabla_{i_m}\upsilon=0.
\end{split}
\end{equation}
 Observe that making
the $X_{*}div$ into and $Xdiv$ introduces tensor fields with a
factor $\nabla Ric_{ij}\nabla^i\phi_{h'}$; where we may then apply
the second Bianchi identity to this expression and make the factor
$\nabla\phi_{h'}$ contract against the derivative index in $\nabla
Ric$. We obtain correction terms that are in the form:

$$\Sum_{t\in T} a_t Xdiv_{i_1}\dots Xdiv_{i_a} C^{t,i_1\dots
i_a}_{g}(\Omega_1,\dots ,\Omega_p,\phi_1,\dots ,\phi_{u+1}).$$

\par Thus, replacing the $\nabla\upsilon$s by $Xdiv$s (see the 
last Lemma in the Appendix in \cite{alexakis1}),
 we obtain
 our desired (\ref{lightbrigade}). $\Box$

\section{An analysis of the  sublinear combination 
\\$CurvTrans[L_g]$.} 
\label{theanalysis}
\subsection{Brief outline of this section: 
How to ``get rid'' of the terms with $\sigma+u$ factors
in (\ref{heidegger}), modulo correction terms we can control.}
\label{anal3sub}

\par Let us recapitulate to recall our main achievements so far and to outline
how our argument will proceed: We have set out to prove Lemmas
\ref{zetajones}, \ref{pool2} (and, eventually Lemma \ref{pskovb}, under the inductive
assumption of Proposition \ref{giade}, along with all the
Corollaries and Lemmas that the inductive assumption of
Proposition \ref{giade} implies. The main assumption for
all these Lemmas is equation (\ref{hypothesegen}), whose
left hand side we denote by $L_g(\Omega_1,\dots,\Omega_p,\phi_1,\dots ,\phi_u)$
 (or just $L_g$, for short).

 \par In Lemma \ref{lemtsabes} we showed that the sublinear
 combination $Image^{1,+}_{\phi_{u+1}}[L_g]$ must vanish
 separately (modulo complete contractions that we may ignore).
 The equation (\ref{tsabes}) is the main assumption for this section. In equation
 (\ref{heidegger}) we broke up $Image^{1,+}_{\phi_{u+1}}[L_g]$
 into three sublinear combinations $CurvTrans[L_g]$, $LC[L_g]$, $W[L_g]$
 which we will study separately in
 the next few subsections (the reader may wish
 to recall these three sublinear combinations now).

Finally, in Lemma \ref{bjorn} we showed that the sublinear
combination $LC_{\Phi}[L_g]$ (in $LC[L_g]$) can be replaced by the
right hand side of the equation in (\ref{bjorn}). Since that right
hand side consists of generic terms that are allowed in the
conclusions of the Lemmas \ref{zetajones}, \ref{pool2} and
Lemma \ref{pskovb} in case A, we may interpret this result as saying that the
sublinear combination $LC_{\Phi}[L_g]$ can be {\it ignored} when we
study $LC[L_g]$ further down.
\newline

{\it Synopsis of subsections  \ref{ct2}, \ref{otherplaces}:} We
commence this section with a study of the sublinear combination
$CurvTrans[L_g]$. The next two subsections
(\ref{otherplaces} and \ref{ct2}) are devoted to that
goal. Our analysis will
proceed as follows: We will firstly seek to understand the
sublinear combination $CurvTrans[L_g]$ as it arises from the
application of the formula (\ref{curvtrans}). We will observe that
the terms we obtain in $CurvTrans[L_g]$ can be grouped up into a
few sublinear combinations, defined by certain algebraic
properties.
 After we do this grouping,  we will apply
  the curvature identity, 
\begin{equation}
\label{curvature}
\nabla^{(2)}_{ab}X_c-\nabla^{(2)}_{ba}X_c=R_{abdc}X^d,
\end{equation}  
 which will introduce
  corrections terms of length $\sigma+u+1$, some of which will be
   important in deriving our Lemmas \ref{zetajones}, \ref{pool2}, 
\ref{pskovb} (in particular the
     sublinear combinations in $Leftover[\dots]$ will be the important ones), and many will
 be generic terms (i.e. generic terms allowed in the conclusions of our three Lemmas).
 Finally, after this analysis and the algebraic manipulations,
 we will still be left with sublinear combinations in
 $CurvTrans[L_g]$ of length $\sigma+u$. Roughly speaking, these
 sublinear combinations will either be linear combinations of
 iterated  $Xdiv$'s with high enough rank, or they will be terms that are
 simply subsequent to $pre\vec{\kappa}^{+}_{simp}$. We will then show
 that these sublinear combinations can be {\it re-written} as
 linear combinations of $Xdiv$'s of tensors fields
 {\it with $\sigma+u+1$ factors}, of the general type that
  is allowed in the conclusions of our Lemmas.
\newline

{\it Notational conventions:} Now, abusing notation, we will
denote $Image^{1,+}_{\phi_{u+1}}[L_g]$ by
$Image^{1}_{\phi_{u+1}}[L_g]$ for the rest of this section.
 Furthermore, we recall that in the setting of Lemmas
\ref{zetajones} and \ref{pool2} the {\it selected factor}
discussed in the definition of $Image^{1,+}_{\phi_{u+1}}[L_g]$ is
{\it always} the crucial factor, defined, in the statements of
Lemmas \ref{zetajones} and \ref{pool2}. On the other hand, in the
setting of \ref{pskovb}, we have declared that the selected
factor is some factor (or set of factors) that we pick once and
for all; it does not have to be the crucial factor. Recall that 
$Image^{1,+}_{\phi_{u+1}}[L_g]$ has been defined in definition
\ref{plus}--this sublinear combination depends on the choice of
selected factor.
 Therefore, in the next subsections, we will
sometimes be making distinctions when we discuss the sublinear
combination $CurvTrans[L_g]$; these distinctions will depend on
which of the Lemmas \ref{zetajones}, \ref{pool2} or \ref{pskovb}
we are proving.
\newline

\subsection{A study of the sublinear combination
$CurvTrans[L_{g}]$ in the setting of Lemmas \ref{pool2} and
 \ref{pskovb}  (when the selected factor is in the form
$\nabla^{(m)}R_{ijkl}$).}
\label{ct2}

\par Let us firstly recall that in the setting of Lemma \ref{pool2}
the notions of ``selected''and ``crucial'' factor coincide. In the setting of Lemma
 \ref{pskovb} they need not coincide. Furthermore, in the
setting of Lemma \ref{pool2}, we will be denoting by $Free(Max)$
 the number of free indices in the crucial factor in the
tensor fields $C^{l,i_1\dots i_\mu}_g, l\in \bigcup_{z\in
Z'_{Max}}L^z$ (see the statement of Lemma
\ref{pool2}).\footnote{Recall that by definition, if the $\mu$-tensorf 
fields of maximal refined double character in (\ref{hypothese2}) 
have $s$ special free indices in the crucial factor $\nabla^{(m)}R_{ijkl}$ ($s=1$ or $s=2$) then 
 all other $\mu$-tensor fields in the
 assumption of Lemma \ref{pool2} will have at most
$Free(Max)$ free indices in any factor $\nabla^{(m)}R_{ijkl}$ that
contains $s$ special free indices.} Whenever we make a claim regarding 
\ref{pskovb} in this subsection, we will be assuming 
that the selected factor is in the from $\nabla^{(m)}R_{ijkl}$.

In this case, we recall that $$CurvTrans[Xdiv_{i_1}\dots
Xdiv_{i_a}C^{l,i_1\dots i_a}_{g} (\Omega_1,\dots
,\Omega_p,\phi_1,\dots ,\phi_u)]$$ stands for the sublinear
combination in $$Image^1_{\phi_{u+1}}[Xdiv_{i_1}\dots
Xdiv_{i_a}C^{l,i_1\dots i_a}_{g} (\Omega_1,\dots
,\Omega_p,\phi_1,\dots ,\phi_u)]$$\footnote{Recall that we are
denoting $Image^{1,+}_{\phi_{u+1}}[\dots]$ by
$Image^1_{\phi_{u+1}}[\dots]$, abusing notation.} that consists of
complete contractions with length $\sigma +u$, with no internal
contractions that arise by replacing the (one of the) selected
factor(s) $\nabla^{(m)}R_{ijkl}$  by one of the four linear
expressions $\nabla^{(m+2)}_{r_1\dots r_mil}\phi_{u+1}g_{jk}$ etc
on the right hand side of (\ref{curvtrans}), provided no internal
contraction arises in that way. We will be treating the function
$\nabla^{(A)}\phi_{u+1}$ as a function 
$\nabla^{(A)}\Omega_{p+1}$ in this subsection.

  In this setting, a complete contraction in
$$CurvTrans
[L_{g} (\Omega_1,\dots ,\Omega_p,\phi_1,\dots ,\phi_u)]$$ will be
called {\it extra acceptable} if it is acceptable (see the discussion
after (\ref{form2}) in \cite{alexakis4}) {\it and} in addition it has all the
 factors $\nabla\phi_{g_i}$ contracting against the factor
$\nabla^{(m+2)}\phi_{u+1}$,\footnote{Here $\{\nabla\phi_{g_i}\}$
stands for the set of terms $\nabla\phi_h$ that are contracting
against the selected factor $\nabla^{(m)}R_{ijkl}$ in 
$\vec{\kappa}_{simp}$.} {\it and}
the two rightmost indices in $\nabla^{(A)}\phi_{u+1}$ are not
contracting against any factor $\nabla\phi_h$. We
straightforwardly observe that all the complete
 contractions in each sublinear combination

\begin{equation}
\label{putaname} CurvTrans [Xdiv_{i_1}\dots
Xdiv_{i_a}C^{l,i_1\dots i_a}_{g} (\Omega_1,\dots
,\Omega_p,\phi_1,\dots ,\phi_u)]
\end{equation}
 are extra acceptable.
\newline

\par As before, we will be using the equation:

\begin{equation}
\label{humanrights}
\begin{split}
& \Sum_{l\in L} a_l CurvTrans[Xdiv_{i_1}\dots Xdiv_{i_a}
C^{l,i_1\dots i_a}_{g}(\Omega_1,\dots ,\Omega_p,\phi_1,\dots
,\phi_u)]+
\\& \Sum_{j\in J} a_j
CurvTrans[C^j_{g}(\Omega_1,\dots ,\Omega_p,\phi_1,\dots
,\phi_u)]=0,
\end{split}
\end{equation}
which holds modulo complete contractions of length $\ge\sigma +u+1$.
\newline

\par Now, we separately
study the sublinear combinations in the left hand side of the above.

\par We start with the sublinear combinations
$CurvTrans[C^j_{g} (\Omega_1,\dots ,\Omega_p,\phi_1,\dots
,\phi_u)]$ for each $j\in J$.

\par Let us introduce some notation. We will denote by
$$\Sum_{j\in J} a_j C^j_{g} (\Omega_1,\dots
,\Omega_p,\phi_{u+1}, \phi_1,\dots ,\phi_u)$$ a generic linear
combination of complete contractions in the form (\ref{plusou})
with a weak character $Weak(pre\vec{\kappa}^{+}_{simp})$ and with
at least one factor $\nabla\phi_f$, $f\in
Def(\vec{\kappa}_{simp})$ contracting against a derivative index
in some factor $\nabla^{(p)}R_{ijkl}$. We then straightforwardly
observe that:

\begin{equation}
\label{suda} CurvTrans[\Sum_{j\in J} a_j C^j_{g}(\Omega_1,\dots
,\Omega_p,\phi_1,\dots ,\phi_u)]= \Sum_{j\in J} a_j C^j_{g}
(\Omega_1,\dots ,\Omega_p,\phi_{u+1}, \phi_1,\dots ,\phi_u).
\end{equation}

\par Next, we proceed to carefully study the sublinear
combinations (\ref{putaname}) in (\ref{humanrights}). We will need to introduce
 some further notational conventions.

\begin{definition}
\label{boh9hse}
\par For this entire subsection, we denote by
$$\Sum_{p\in P} a_p C^{p,i_1\dots
i_a,i_{*}}_{g} (\Omega_1,\dots ,\Omega_p,\phi_1,\dots
,\phi_u)\nabla_{i_{*}}\phi_{u+1}$$ a generic linear combination of
acceptable tensor fields with length $\sigma +u+1$, and $a\ge \mu$
and with a $u$-simple character $\vec{\kappa}_{simp}$ and
a weak $(u+1)$-character $Weak(\vec{\kappa}_{simp}^{+})$.

\par Furthermore, we denote by
$$\Sum_{u\in U} a_u Xdiv_{i_1}\dots Xdiv_{i_a}C^{u,i_1\dots
i_a}_{g} (\Omega_1,\dots ,\Omega_p,\phi_{u+1},\phi_1,\dots
,\phi_u)$$
 a generic linear combination of $Xdiv$s of extra acceptable $a$-tensor
fields $(a\ge \mu$) with the following features: The tensor fields 
$C^{u,i_1\dots i_a}_{g} (\Omega_1,\dots ,\Omega_p,\phi_{u+1},\phi_1,\dots
,\phi_u)$ must have of length $\sigma +u$, be in the form
(\ref{plusou}) and have a $u$-simple character
$pre\vec{\kappa}^{+}_{simp}$, where all the factors
$\nabla\phi_{g_1},\dots ,\nabla\phi_{g_z}$ are precisely those
$\nabla\phi$'s that are contracting against the first $z$ indices
in the factor $\nabla^{(A)}\phi_{u+1}$ with $A\ge z+2$.
Furthermore, we require that if $a=\mu$  then either at
least one of the free indices ${}_{i_1},\dots,{}_{i_\mu}$ is a
derivative index (and moreover if it belongs to a factor 
$\nabla^{(B)}\Omega_h$ then $B\ge 3$), or none of these indices is a 
special index in a factor $S_*\nabla^{(\nu)}R_{ijkl}$.
\newline

\par Our next definitions will be only for the setting of Lemma \ref{pool2}.

\par We will  denote by
$$\Sum_{t\in T} a_t C^{t,i_1\dots
i_{\mu}}_{g} (\Omega_1,\dots ,\Omega_p,\phi_1,\dots
,\phi_u)\nabla_{i_1}\phi_{u+1}$$ a generic linear combination of
acceptable $(\mu-1)$-tensor fields of length $\sigma+u+1$ with
$(u+1)$-simple character $\vec{\kappa}^{+}_{simp}$, for which
either the selected (=crucial=critical) factor
 contains fewer than $Free(Max)-1$ free indices, or it contains exactly
$Free(Max)-1$ free indices but its refined double character is doubly 
subsequent to each $\vec{L^z}', z\in Z'_{Max}$.\footnote{In other words,
in the statements of Lemma \ref{pool2} this corresponds to a
generic linear combination $\Sum_{\nu\in N} a_\nu C^{\nu,i_1\dots
,i_{\mu-1}i_\mu}_{g}(\Omega_1,\dots ,\Omega_p,\phi_1,\dots
,\phi_u)\nabla_{i_\mu}\phi_{u+1}$.}

\par In addition (again only in the setting of Lemma \ref{pool2})
 we denote by $$\Sum_{u\in
U_1} a_u C^{u,i_1\dots i_{\mu -1}}_{g} (\Omega_1,\dots
,\Omega_p,\phi_{u+1},\phi_1,\dots ,\phi_u)$$ a generic linear
combination of extra acceptable $(\mu -1)$ tensor fields in
 the form (\ref{plusou}), of length $\sigma+u$
with a simple character $pre\vec{\kappa}_{simp}^{+}$, with the
extra property that the factor $\nabla^{(A)}\phi_{u+1}$ has fewer
than $Free(Max)-1$ free indices. 
We also require that none of the free indices are special indices in a factor $S_{*}\nabla^{(\nu)}R_{ijkl}$.

\par Moreover, in the case where the maximal refined double characters
$\vec{L^z}', z\in Z'_{Max}$ have two internal free indices in the
crucial factor, we denote by $$\Sum_{u\in U_2} a_u C^{u,i_1\dots
i_{\mu -1}}_{g} (\Omega_1,\dots ,\Omega_p,\phi_{u+1},\phi_1,\dots
,\phi_u)$$  a generic linear combination of extra acceptable $(\mu
-1)$ tensor fields in the form (\ref{plusou}) with a simple
character $pre\vec{\kappa}^{+}_{simp}$ and with the feature that
it has $Free(Max)-1$ free indices in the factor
$\nabla^{(A)}\phi_{u+1}$ and with the additional property that one
of the free indices that belong to the factor
$\nabla^{(A)}_{r_1\dots r_A}\phi_{u+1}$ is the last index
${}_{r_A}$. We also require that none of the free indices are special indices in a factor $S_{*}\nabla^{(\nu)}R_{ijkl}$.
\end{definition}

\par Armed with this definition, we may now study the sublinear
 combination (\ref{putaname}) in detail.

\par We recall a few notational conventions we have made in the setting of Lemma \ref{pool2}:

\par Recall that in the setting of Lemma \ref{pool2},
the crucial factor(s) in each $C^{l,i_1\dots i_\mu}$, $l\in L^z,
z\in Z'_{Max}$ will all have either two, one or no internal free
indices.

If the maximal refined double characters
$\vec{L^z}, z\in Z'_{Max}$ have one or two internal free indices
belonging to the (a) crucial factor, we recall that we have denoted
 by $I_{*,l}\subset I_l$ the set of all internal free indices
 that belong to a crucial factor. Furthermore, 
if there are two such free indices,  we have declared that in the (each)
crucial factor of the form $T=\nabla^{(m)}R_{ijkl}$, the internal free indices
will be the indices ${}_i,{}_k$. In that setting, we may then assume wlog that the
 indices ${}_{i_1},{}_{i_3},\dots ,{}_{i_{2k_l+1}}\in I_{*,l}$ are the
 indices ${}_i$ in the crucial factors. Also in this case we may
 assume wlog that if $k$ is odd then ${}_k$ and ${}_{k+1}$ belong
 to the same crucial factor.  We then claim:

\begin{lemma}
\label{burk}
In the setting of Lemma \ref{pool2}, if the tensor
fields $C^{l,i_1\dots i_\mu}_g, l\in L^z, z\in Z'_{Max}$
have two internal free indices in the crucial factor(s) then:

\begin{equation}
\label{ankourasteis}
\begin{split}
&\Sum_{l\in L_\mu} a_l CurvTrans[Xdiv_{i_1}\dots
Xdiv_{i_\mu}C^{l,i_1\dots i_\mu}_{g} (\Omega_1,\dots
,\Omega_p,\phi_1,\dots ,\phi_u)]=
\\& \Sum_{z\in Z'_{Max}}\Sum_{l\in L^z} a_l
\Sum_{h=0}^{k_l}Xdiv_{i_1}\dots \hat{Xdiv}_{i_{2h+1}}\dots
 Xdiv_{i_\mu}\\&C^{l,i_1\dots i_\mu}_{g} (\Omega_1,\dots ,
\Omega_p,\phi_1,\dots,\phi_u)\nabla_{i_1}\phi_{u+1}
\\&+\Sum_{u\in U} a_u Xdiv_{i_1}\dots Xdiv_{i_a}C^{u,i_1\dots
i_a}_{g} (\Omega_1,\dots ,\Omega_p,\phi_{u+1},\phi_1,\dots
,\phi_u)+
\\&\Sum_{u\in U_1} a_u Xdiv_{i_1}\dots Xdiv_{i_{\mu-1}}C^{u,i_1\dots
i_{\mu-1}}_{g} (\Omega_1,\dots ,\Omega_p,\phi_{u+1},\phi_1,\dots
,\phi_u)+
\\&\Sum_{u\in U_2} a_u Xdiv_{i_1}\dots Xdiv_{i_{\mu-1}}C^{u,i_1\dots
i_{\mu-1}}_{g} (\Omega_1,\dots ,\Omega_p,\phi_{u+1},\phi_1,\dots
,\phi_u)+
\\&\Sum_{t\in T} a_t Xdiv_{i_2}\dots Xdiv_{i_\mu}C^{t,i_1\dots
i_\mu}_{g} (\Omega_1,\dots ,\Omega_p,\phi_1,\dots
,\phi_u)\nabla_{i_1}\phi_{u+1}.
\end{split}
\end{equation}

\par Next claim: Again in the setting of Lemma \ref{pool2},
if the tensor fields $C^{l,i_1\dots i_\mu}_g, l\in L^z, z\in Z'_{Max}$
have one internal free index in the crucial factor(s)
then:

\begin{equation}
\label{ankourasteis2}
\begin{split}
&\Sum_{l\in L_\mu} a_l CurvTrans[Xdiv_{i_1}\dots
Xdiv_{i_\mu}C^{l,i_1\dots i_\mu}_{g} (\Omega_1,\dots
,\Omega_p,\phi_1,\dots ,\phi_u)]=
\\&\Sum_{z\in Z'_{Max}}\Sum_{l\in L^z} a_l
\Sum_{i_h\in I_{*,l}} Xdiv_{i_1}\dots \hat{Xdiv}_{i_h}\dots
Xdiv_{i_\mu}C^{l,i_1\dots i_\mu}_{g} (\Omega_1,\dots
,\Omega_p,\phi_1, \dots,\phi_u)\\&\nabla_{i_h}\phi_{u+1}
+\Sum_{u\in U} a_u Xdiv_{i_1}\dots Xdiv_{i_\mu}C^{u,i_1\dots
i_{\mu}}_{g} (\Omega_1,\dots ,\Omega_p,\phi_{u+1},\phi_1,\dots
,\phi_u)+
\\&\Sum_{t\in T} a_t Xdiv_{i_2}\dots Xdiv_{i_{\mu}}C^{t,i_1\dots
i_{\mu}}_{g} (\Omega_1,\dots ,\Omega_p,\phi_1,\dots
,\phi_u)\nabla_{i_1}\phi_{u+1}.
\end{split}
\end{equation}

\par In the setting of Lemma \ref{pskovb} and also
in the setting of Lemma \ref{pool2} if the tensor fields
 $C^{l,i_1\dots i_\mu}_g, l\in L^z, z\in Z'_{Max}$
have no internal free index in the crucial factor(s) then:

\begin{equation}
\label{ankourasteis3}
\begin{split}
&\Sum_{l\in L_\mu} a_l CurvTrans[Xdiv_{i_1}\dots
Xdiv_{i_\mu}C^{l,i_1\dots i_\mu}_{g} (\Omega_1,\dots
,\Omega_p,\phi_1,\dots ,\phi_u)]
\\&=\Sum_{u\in U} a_u Xdiv_{i_1}\dots Xdiv_{i_\mu}C^{u,i_1\dots
i_{\mu}}_{g} (\Omega_1,\dots ,\Omega_p,\phi_{u+1},\phi_1,\dots
,\phi_u).
\end{split}
\end{equation}

\par Moreover, for both Lemmas \ref{pool2} and \ref{pskovb}:

\begin{equation}
\label{truth}
\begin{split}
&\Sum_{l\in (L\setminus L_\mu)} a_l CurvTrans[Xdiv_{i_1}\dots
Xdiv_{i_a}C^{l,i_1\dots i_a}_{g} (\Omega_1,\dots
,\Omega_p,\phi_1,\dots ,\phi_u)]
\\&=\Sum_{p\in P} a_p Xdiv_{i_1}\dots Xdiv_{i_a}
C^{p,i_1\dots i_a,i_{*}}_{g} (\Omega_1,\dots
,\Omega_p,\phi_1,\dots ,\phi_u)\nabla_{i_{*}}\phi_{u+1}+
\\&\Sum_{u\in U} a_u Xdiv_{i_1}\dots Xdiv_{i_a}C^{u,i_1\dots
i_a}_{g} (\Omega_1,\dots ,\Omega_p,\phi_{u+1},\phi_1,\dots
,\phi_u).
\end{split}
\end{equation}
\end{lemma}

{\it Proof of Lemma \ref{burk}:} The proof just follows by
 applying the transformation law (\ref{curvtrans}). For the first two equations in
the Lemma, we ``complete the divergence''  to get
 the terms on the first two lines of the right hand sides.\footnote{We
 explain the notion of ``completing the divergence'': We observe
that we obtain terms in $CurvTrans[Xdiv_{i_1}\dots
Xdiv_{i_a}C^{l,i_1\dots i_a}_{g}]$ which are in the form
$X_{*}div_{i_1}Xdiv_{i_2}\dots Xdiv_{i_\mu}C^{*,i_1\dots
i_\mu}_g$, where the free index ${}_{i_1}$ does not belong to the
factor $\nabla^{(B)}\phi_{u+1}$, and $X_{*}div_{i_1}$ means that
$\nabla^{i_1}$ is not allowed to hit the factor
$\nabla^{(B)}\phi_{u+1}$. Then adding and subtracting a term
$Hitdiv_{i_1}Xdiv_{i_2}\dots Xdiv_{i_\mu}C^{*,i_1\dots i_\mu}_g$
($Hitdiv_{i_1}$ means that we force $\nabla^{i_1}$ to hit the
factor $\nabla^{(B)}\phi_{u+1}$), we obtain the terms of length
$\sigma+u+1$ in (\ref{ankourasteis}), (\ref{ankourasteis2}) (when
we subtract the term in question), by also applying the curvature
identity (\ref{curvature}).}
  Also, for the first two equations,
the proof of our claim also relies on  the definition of {\it maximal}
refined double characters. {\it Note:} For (\ref{ankourasteis2})
we also use the fact that (\ref{hypothese2}) does not fall under the 
``special cases'' outlined at the very end of the introduction. $\Box$
\newline

In conclusion, we have shown that in the setting of Lemma
\ref{pool2} when there are two internal free indices in the
crucial factor in the tensor fields $C^{l,i_1\dots i_\mu}_{g}$,
$L\in L^z, z\in Z'_{Max}$ we will have:

\begin{equation}
\label{akinola}
\begin{split}
&CurvTrans[L_{g}(\Omega_1,\dots ,\Omega_p,\phi_1,\dots ,\phi_u)]=
\\& \Sum_{z\in Z'_{Max}}\Sum_{l\in L^z} a_l
\Sum_{h=0}^{k_l}Xdiv_{i_1}\dots \hat{Xdiv}_{i_{2h+1}}\dots
 Xdiv_{i_\mu}C^{l,i_1\dots i_\mu}_{g} (\Omega_1,\dots ,
\Omega_p,\phi_1,\dots,\phi_u)
\\&+\nabla_{i_1}\phi_{u+1}\Sum_{u\in U} a_u Xdiv_{i_1}\dots Xdiv_{i_a}C^{u,i_1\dots
i_a}_{g} (\Omega_1,\dots ,\Omega_p,\phi_{u+1},\phi_1,\dots
,\phi_u)+
\\&\Sum_{u\in U_1} a_u Xdiv_{i_1}\dots Xdiv_{i_{\mu-1}}
C^{u,i_1\dots i_{\mu-1}}_{g} (\Omega_1,\dots
,\Omega_p,\phi_{u+1},\phi_1, \dots ,\phi_u)
\\&+\Sum_{u\in U_2} a_u Xdiv_{i_1}\dots Xdiv_{i_{\mu-1}}C^{u,i_1\dots
i_{\mu-1}}_{g} (\Omega_1,\dots ,\Omega_p,\phi_{u+1},\phi_1,
\dots,\phi_u)
\\& +\Sum_{j\in J} a_j C^j_{g}(\Omega_1,\dots
,\Omega_p,\phi_{u+1},\phi_1,\dots ,\phi_u)+
\\&\Sum_{p\in P} a_p Xdiv_{i_1}\dots Xdiv_{i_a}
C^{p,i_1\dots i_a,i_{*}}_{g} (\Omega_1,\dots
,\Omega_p,\phi_1,\dots ,\phi_u)\nabla_{i_{*}}\phi_{u+1}+
\\&\Sum_{t\in T} a_t Xdiv_{i_1}\dots Xdiv_{i_{\mu-1}}
C^{t,i_1\dots i_\mu}_{g} (\Omega_1,\dots ,\Omega_p,\phi_1,\dots
,\phi_u)\nabla_{i_\mu}\phi_{u+1}.
\end{split}
\end{equation}

Also, in the setting of Lemma \ref{pool2} when there is one
internal free index in the crucial factor in the tensor fields
$C^{l,i_1\dots i_\mu}_{g}$, $l\in L^z, z\in Z'_{Max}$ we will
have:

\begin{equation}
\label{akinola2}
\begin{split}
&CurvTrans[L_{g}(\Omega_1,\dots ,\Omega_p,\phi_1,\dots ,\phi_u)]=
\\&\Sum_{z\in Z'_{Max}}\Sum_{l\in L^z} a_l
\Sum_{i_h\in I_{*,l}} Xdiv_{i_1}\dots \hat{Xdiv}_{i_h}\dots
Xdiv_{i_\mu}C^{l,i_1\dots i_\mu}_{g} (\Omega_1,\dots
,\Omega_p,\phi_1, \dots,\phi_u)
\\&+\nabla_{i_h}\phi_{u+1}\Sum_{u\in U} a_u Xdiv_{i_1}\dots Xdiv_{i_a}C^{u,i_1\dots
i_a}_{g} (\Omega_1,\dots ,\Omega_p,\phi_{u+1},\phi_1,\dots
,\phi_u)
\\&+\Sum_{u\in U_1} a_u Xdiv_{i_1}\dots Xdiv_{i_{\mu-1}}
C^{u,i_1\dots i_{\mu-1}}_{g} (\Omega_1,\dots
,\Omega_p,\phi_{u+1},\phi_1, \dots ,\phi_u)+
\\& \Sum_{j\in J} a_j C^j_{g}(\Omega_1,\dots
,\Omega_p,\phi_{u+1},\phi_1,\dots ,\phi_u)+
\\&\Sum_{p\in P} a_p Xdiv_{i_1}\dots Xdiv_{i_a}
C^{p,i_1\dots i_a,i_{*}}_{g} (\Omega_1,\dots
,\Omega_p,\phi_1,\dots ,\phi_u)\nabla_{i_{*}}\phi_{u+1}+
\\&\Sum_{t\in T} a_t Xdiv_{i_1}\dots Xdiv_{i_{\mu-1}}
C^{t,i_1\dots i_\mu}_{g} (\Omega_1,\dots ,\Omega_p,\phi_1,\dots
,\phi_u)\nabla_{i_\mu}\phi_{u+1},
\end{split}
\end{equation}
and lastly in the setting of Lemma \ref{pskovb}
 we will have:

\begin{equation}
\label{akinola3}
\begin{split}
&CurvTrans[L_{g}(\Omega_1,\dots ,\Omega_p,\phi_1,\dots ,\phi_u)]=
\\&\Sum_{u\in U} a_u Xdiv_{i_1}\dots Xdiv_{i_a}C^{u,i_1\dots
i_a}_{g} (\Omega_1,\dots ,\Omega_p,\phi_{u+1},\phi_1,\dots
,\phi_u)+
\\& +\Sum_{j\in J} a_j C^j_{g}(\Omega_1,\dots
,\Omega_p,\phi_{u+1},\phi_1,\dots ,\phi_u)+
\\&\Sum_{p\in P} a_p Xdiv_{i_1}\dots Xdiv_{i_a}
C^{p,i_1\dots i_a,i_{*}}_{g} (\Omega_1,\dots
,\Omega_p,\phi_1,\dots ,\phi_u)\nabla_{i_{*}}\phi_{u+1}.
\end{split}
\end{equation}

\par We then make three claims, for each of the three subcases
above. We are interested in ``getting rid of'' the sublinear combinations
 of length $\sigma +u$ that we have been left with in
 $CurvTrans[L_{g}(\Omega_1,\dots ,\Omega_p,
\phi_1,\dots,\phi_u)]$.
We first consider  the setting of
 Lemma \ref{pool2} and there are tensor fields indexed in
$L_\mu\subset L$ with two internal free indices in the crucial
factor $\nabla^{(m)}R_{ijkl}$; call this the first subcase. We claim:

\begin{lemma}
\label{en9ousiasmos} Consider the setting of
 Lemma \ref{pool2} when there are tensor fields indexed in
$L_\mu\subset L$ with two free indices in the crucial factor
$\nabla^{(m)}R_{ijkl}$. Then, refer to (\ref{akinola}). By virtue
of our inductive assumption on Proposition \ref{giade}, we claim
that the sublinear combination
 of length $\sigma +u$ in (\ref{akinola}) will be equal to:

\begin{equation}
\label{enduo}
\begin{split}
&\Sum_{u\in U} a_u Xdiv_{i_1}\dots Xdiv_{i_a}C^{u,i_1\dots
i_a}_{g} (\Omega_1,\dots ,\Omega_p,\phi_{u+1},\phi_1,\dots
,\phi_u)+
\\&\Sum_{u\in U_1} a_u Xdiv_{i_1}\dots Xdiv_{i_{\mu-1}}
C^{u,i_1\dots i_{\mu-1}}_{g} (\Omega_1,\dots
,\Omega_p,\phi_{u+1},\phi_1, \dots ,\phi_u)+
\\&\Sum_{u\in U_2} a_u Xdiv_{i_1}\dots Xdiv_{i_{\mu-1}}
C^{u,i_1\dots i_{\mu-1}}_{g} (\Omega_1,\dots
,\Omega_p,\phi_{u+1},\phi_1, \dots ,\phi_u)
\\& +\Sum_{j\in J} a_j C^j_{g}(\Omega_1,\dots
,\Omega_p,\phi_{u+1},\phi_1,\dots ,\phi_u)=
\\&\Sum_{r\in R_\alpha} a_r Xdiv_{i_1}\dots Xdiv_{i_{\mu-1}}
C^{r,i_1\dots i_{\mu -1}}_{g} (\Omega_1,\dots
,\Omega_p,\phi_1,\dots ,\phi_{u+1})
\\&+\Sum_{r\in
R_\beta} a_r Xdiv_{i_1}\dots Xdiv_{i_{\mu-1}}C^{r,i_1\dots i_{\mu
-1}}_{g} (\Omega_1,\dots ,\Omega_p,\phi_1,\dots ,\phi_{u+1})
\\&+\Sum_{r\in R_\gamma} a_r Xdiv_{i_1}\dots Xdiv_{i_a}C^{r,i_1\dots i_a}_{g}
(\Omega_1,\dots ,\Omega_p,\phi_1,\dots ,\phi_u,\phi_{u+1})+
\\&\Sum_{j\in J} a_j C^j_{g}(\Omega_1,\dots
,\Omega_p,\phi_1,\dots ,\phi_{u+1});
\end{split}
\end{equation}
here the $(\mu -1)$-tensor fields indexed in $R_\alpha$ have
 a $(u+1)$-simple character $\vec{\kappa}^{+}_{simp}$, and
 also  have $Free(Max)-1$ free indices belonging to the
 crucial factor $\nabla^{(m)}_{r_1\dots r_m}R_{ijkl}$ and moreover all of them
  are of the form ${}_{r_1},\dots ,{}_{r_m},{}_j$ and furthermore
$\nabla\phi_{u+1}$ is contracting against the index ${}_i$ (so in
particular they are doubly subsequent to all $\vec{L^z}', z\in
Z'_{Max}$).

 Also, the $(\mu-1)$-tensor fields
that are indexed in $R_\beta$ have a refined double character that
is doubly subsequent to each $\vec{L^z}', z\in Z'_{Max}$. (In
particular they have fewer than $Free(Max)-1$ free indices in the
crucial factor). Finally, each $a$-tensor field ($a\ge \mu$)
indexed in $R_\gamma$ has a $u$-simple character
$\vec{\kappa}_{simp}$ and a weak $(u+1)$-character
$Weak(\vec{\kappa}^{+}_{simp})$. Lastly, each $C^j$ has length
$\sigma +u+1$ and is simply subsequent to
$\vec{\kappa}^{+}_{simp}$.
\end{lemma}

\par Next, we consider the  setting of Lemma
\ref{pool2} where there is one internal free index in the tensor
fields $C^{l,i_1\dots i_\mu}_g, l\in L^z, z\in Z'_{Max}$ (call
this the {\it second subcase}). We also consider the setting of
Lemma \ref{pskovb} (where there are no internal free indices in
any factor $\nabla^{(m)}R_{ijkl}$ in the $\mu$-tensor fields we
are considering; call this the {\it third subcase}).

\begin{lemma}
\label{en9ousiasmos2} By virtue of our inductive assumption on
Proposition \ref{giade}, in both (second and third) subcases above
we claim that the sublinear combination of length $\sigma +u$ that
we have been left with in $CurvTrans[L_{g}]$ (see (\ref{akinola2})
and (\ref{akinola3})) can be written as:

\begin{equation}
\label{enduo2}
\begin{split}
&\Sum_{u\in U} a_u Xdiv_{i_1}\dots Xdiv_{i_a}C^{u,i_1\dots
i_a}_{g} (\Omega_1,\dots ,\Omega_p,\phi_{u+1},\phi_1,\dots
,\phi_u)+
\\&\Sum_{j\in J} a_j C^j_{g}(\Omega_1,\dots
,\Omega_p,\phi_{u+1},\phi_1,\dots ,\phi_u)=
\\&(\Sum_{r\in R_\beta} a_r Xdiv_{i_1}\dots Xdiv_{i_{\mu-1}}
C^{r,i_1\dots i_{\mu-1}}_{g} (\Omega_1,\dots ,\Omega_p,\phi_{u+1},
\phi_1,\dots,\phi_{u+1}))+
\\&\Sum_{r\in R_\gamma} a_r Xdiv_{i_1}\dots Xdiv_{i_a}C^{r,i_1\dots i_a}_{g}
(\Omega_1,\dots ,\Omega_p,\phi_1,\dots ,\phi_u,\phi_{u+1})+
\\& \Sum_{j\in J} a_j C^j_{g}(\Omega_1,\dots
,\Omega_p,\phi_1,\dots ,\phi_u,\phi_{u+1});
\end{split}
\end{equation}
 The $(\mu-1)$-tensor fields that are indexed in
$R_\beta$ arise only in the second subcase and have a (refined) double
character that is subsequent to each $\vec{L^z}'$, $z\in
Z'_{Max}$. The tensor fields in $R_\beta,R_\gamma$ are as above.
\end{lemma}

{\it Proof of Lemmas \ref{en9ousiasmos} and
\ref{en9ousiasmos2}:}
\newline

\par We will prove Lemma \ref{en9ousiasmos}. We will then indicate how
Lemma \ref{en9ousiasmos2} follows by the same argument.

\par By virtue of the equation (\ref{heidegger}) we have that:

\begin{equation}
\label{saito}
\begin{split}
&\Sum_{u\in U} a_u Xdiv_{i_1}\dots Xdiv_{i_a}C^{u,i_1\dots
i_a}_{g} (\Omega_1,\dots ,\Omega_p,\phi_{u+1},\phi_1,\dots
,\phi_u)+
\\&\Sum_{u\in U_1} a_u Xdiv_{i_1}\dots Xdiv_{i_{\mu-1}}C^{u,i_1\dots
i_{\mu-1}}_{g} (\Omega_1,\dots ,\Omega_p,\phi_{u+1},\phi_1,\dots
,\phi_u)+
\\&\Sum_{u\in U_2} a_u Xdiv_{i_1}\dots Xdiv_{i_{\mu-1}}C^{u,i_1\dots
i_{\mu-1}}_{g} (\Omega_1,\dots
,\Omega_p,\phi_{u+1},\phi_1,\dots,\phi_u)
\\&+\Sum_{j\in J} a_j C^j_{g}(\Omega_1,\dots
,\Omega_p,\phi_{u+1},\phi_1,\dots ,\phi_u)=0,
\end{split}
\end{equation}
modulo contractions of length $\ge\sigma +u+1$.
\newline

\par We focus on the left hand side linear combination in
(\ref{saito}), where we treat the function $\phi_{u+1}$ as a
function $\Omega_{p+1}$, and we apply the eraser to the factors
$\nabla\phi_g$ that are contracting against
$\nabla^{(A)}\phi_{u+1}$. (We will be applying this operation in
order to apply the inductive assumption of Corollary
1 in \cite{alexakis4} on various occasions below; after we have applied
our Corollary, we will then re-introduce the factors
$\nabla\phi_g$ that we erased.\footnote{By abuse of notation we
will sometimes use $\overline{\vec{\kappa}}_{simp}$ to also denote
the simple character with the factors $\nabla\phi_g$ put back
in.}) We observe that the tensor fields of length $\sigma +u$
that we obtain via this operation
 are acceptable, and will all have the same
 simple character which we denote by $\overline{\vec{\kappa}}_{simp}$. Furthermore, the
complete contractions $C^j_{g}(\Omega_1,\dots
,\Omega_p,\phi_{u+1},\phi_1,\dots ,\phi_u)$ that arise after
applying the eraser will be  subsequent to the simple character
$\overline{\vec{\kappa}}_{simp}$.

\par Thus, we can apply our inductive assumption on Corollary
1 in \cite{alexakis4}\footnote{Since the tensor fields indexed
in $U_1$ have no special free indices in factors $S_{*}\nabla^{(\nu)}R_{ijkl}$, 
there is no danger of falling under a ``forbidden case''.}
 to the left hand side of (\ref{saito}) (to which we
have applied the eraser). We conclude that there is a linear
combination of acceptable $\mu$-tensor fields (indexed in $T$
below) with a simple character $\overline{\vec{\kappa}}_{simp}$ so
that:

\begin{equation}
\label{mooo}
\begin{split}
&\Sum_{u\in U_1} a_u C^{u,i_1\dots i_{\mu-1}}_{g} (\Omega_1,\dots
,\Omega_p,\phi_{u+1},\phi_1,\dots
,\phi_u)\nabla_{i_1}\upsilon\dots \nabla_{i_{\mu-1}}\upsilon+
\\&\Sum_{u\in U_2} a_u C^{u,i_1\dots
i_{\mu-1}}_{g} (\Omega_1,\dots ,\Omega_p,\phi_{u+1},\phi_1,\dots
,\phi_u)\nabla_{i_1}\upsilon\dots \nabla_{i_{\mu-1}}\upsilon -
\\&\Sum_{t\in T} a_u Xdiv_{i_\mu}C^{t,i_1\dots
i_{\mu-1}i_\mu}_{g} (\Omega_1,\dots
,\Omega_p,\phi_{u+1},\phi_1,\dots
,\phi_u)\nabla_{i_1}\upsilon\dots \nabla_{i_{\mu-1}}\upsilon=
\\& \Sum_{f\in F} a_f C^{f,i_1\dots i_{\mu-1}}_{g}
(\Omega_1,\dots,\Omega_p,\phi_{u+1},\phi_1,\dots ,\phi_u)
\nabla_{i_1}\upsilon\dots \nabla_{i_{\mu-1}}\upsilon,
\end{split}
\end{equation}
modulo complete contractions of length $\ge\sigma +u+\mu$. Here
each $C^{f,i_1\dots i_{\mu-1}}_{g}$ is simply subsequent to $\overline{\vec{\kappa}}_{simp}$.
 Now, we index in $T_2\subset T$
the tensor fields with precisely $Free(Max)-1$ factors
$\nabla\upsilon$ contracting against $\nabla^{(A)}\phi_{u+1}$ and
in $T_1\subset T$ the tensor fields with fewer than $Free(Max)-1$
factors $\nabla\upsilon$ contracting against
$\nabla^{(A)}\phi_{u+1}$. Since the above holds formally, we then
derive two equations: Firstly:

\begin{equation}
\label{internetia}
\begin{split}
&\Sum_{u\in U_2} a_u C^{u,i_1\dots i_{\mu-1}}_{g} (\Omega_1,\dots
,\Omega_p,\phi_{u+1},\phi_1,\dots
,\phi_u)\nabla_{i_1}\upsilon\dots \nabla_{i_{\mu-1}}\upsilon-
\\&\Sum_{t\in T_2} a_u Xdiv_{i_\mu}C^{t,i_1\dots
i_{\mu-1}i_\mu}_{g} (\Omega_1,\dots
,\Omega_p,\phi_{u+1},\phi_1,\dots
,\phi_u)\nabla_{i_1}\upsilon\dots \nabla_{i_{\mu-1}}\upsilon=
\\& \Sum_{f\in F} a_f C^{f,i_1\dots i_{\mu-1}}_{g}
(\Omega_1,\dots,\Omega_p,\phi_{u+1},\phi_1,\dots ,\phi_u)
\nabla_{i_1}\upsilon\dots \nabla_{i_{\mu-1}}\upsilon,
\end{split}
\end{equation}
modulo complete contractions of length $\ge\sigma +u+\mu$.
 Moreover, we may assume with no
loss of generality that for each of the tensor fields
\\$C^{t,i_1\dots i_{\mu-1}i_\mu}_{g} (\Omega_1,\dots
,\Omega_p,\phi_{u+1},\phi_1,\dots ,\phi_u)$ above, one of the free
indices ${}_{i_1},\dots ,{}_{i_{\mu-1}}$ that belongs to
 the factor $\nabla^{(A)}\phi_{u+1}=\nabla^{(A)}\Omega_{p+1}$ is the last index
 ${}_{r_A}$ in that factor (as is  also the case of the contractions indexed in $U_2$).
\newline

\par Secondly, we derive:

 \begin{equation}
\label{mooo2}
\begin{split}
&\Sum_{u\in U_1} a_u C^{u,i_1\dots i_{\mu-1}}_{g} (\Omega_1,\dots
,\Omega_p,\phi_{u+1},\phi_1,\dots
,\phi_u)\nabla_{i_1}\upsilon\dots \nabla_{i_{\mu-1}}\upsilon-
\\&\Sum_{t\in T_1} a_u Xdiv_{i_\mu}C^{t,i_1\dots
i_{\mu-1}i_\mu}_{g} (\Omega_1,\dots
,\Omega_p,\phi_{u+1},\phi_1,\dots
,\phi_u)\nabla_{i_1}\upsilon\dots \nabla_{i_{\mu-1}}\upsilon=
\\&\Sum_{f\in F} a_f C^{f,i_1\dots i_{\mu-1}}_{g}
(\Omega_1,\dots,\Omega_p,\phi_{u+1},\phi_1,\dots ,\phi_u)
\nabla_{i_1}\upsilon\dots \nabla_{i_{\mu-1}}\upsilon,
\end{split}
\end{equation}
modulo complete contractions of length $\ge\sigma +u+\mu$.

\par Now, we seek to use the fact that the above equations hold formally
 to derive information about the correction terms of length
$\sigma +\mu +u$ in (\ref{internetia}) and (\ref{mooo2}).
 In (\ref{internetia}) we
use the fact that the index ${}_{i_1}$ is the last index in the
factor $\nabla^{(A)}\phi_{u+1}$. We may then assume (using the
eraser) that in the permutations by which we make the left hand
side of (\ref{internetia}) formally zero, the last index in the
factor $\nabla^{(A)}\phi_{u+1}$ (which is contracting against the
factor $\nabla\upsilon$) is not permuted. We conclude
 that we can write out:

\begin{equation}
\label{froma}
\begin{split}
&\Sum_{u\in U_2} a_u C^{u,i_1\dots i_{\mu-1}}_{g} (\Omega_1,\dots
,\Omega_p,\phi_{u+1},\phi_1,\dots
,\phi_u)\nabla_{i_1}\upsilon\dots \nabla_{i_{\mu-1}}\upsilon-
\\&\Sum_{t\in T_2} a_u Xdiv_{i_\mu}C^{t,i_1\dots
i_{\mu-1}i_\mu}_{g} (\Omega_1,\dots
,\Omega_p,\phi_{u+1},\phi_1,\dots
,\phi_u)\nabla_{i_1}\upsilon\dots \nabla_{i_{\mu-1}}\upsilon=
\\&\Sum_{r\in R_\alpha} a_r
C^{r,i_1\dots i_{\mu -1}}_{g} (\Omega_1,\dots
,\Omega_p,\phi_1,\dots ,\phi_{u+1})\nabla_{i_1}\upsilon \dots
\nabla_{i_{\mu-1}}\upsilon
\\& +\Sum_{f\in F} a_f C^{f,i_1\dots i_{\mu-1}}_{g}
(\Omega_1,\dots,\Omega_p,\phi_{u+1},\phi_1,\dots ,\phi_u)
\nabla_{i_1}\upsilon\dots \nabla_{i_{\mu-1}}\upsilon
\\&+\Sum_{z\in Z} a_z C^{z,i_1\dots
i_{\mu-1}}_{g} (\Omega_1,\dots ,\Omega_p,\phi_1,\dots
,\phi_u,\phi_{u+1})\nabla_{i_1}\upsilon\dots
\nabla_{i_{\mu-1}}\upsilon,
\end{split}
\end{equation}
where each tensor field
 indexed in $Z$ has length $\sigma +u+1$ and a
 factor $\nabla^{(b)}\phi_{u+1}$, $b\ge 2$. The above holds modulo
 correction terms of length greater than the $RHS$.

\par By a similar argument, we derive that  we can write:

\begin{equation}
\label{froma2}
\begin{split}
&\Sum_{u\in U_1} a_u C^{u,i_1\dots i_{\mu-1}}_{g} (\Omega_1,\dots
,\Omega_p,\phi_{u+1},\phi_1,\dots
,\phi_u)\nabla_{i_1}\upsilon\dots \nabla_{i_{\mu-1}}\upsilon-
\\&\Sum_{t\in T_1} a_u Xdiv_{i_\mu}C^{t,i_1\dots
i_{\mu-1}i_\mu}_{g} (\Omega_1,\dots
,\Omega_p,\phi_{u+1},\phi_1,\dots
,\phi_u)\nabla_{i_1}\upsilon\dots \nabla_{i_{\mu-1}}\upsilon=
\\&\Sum_{f\in F} a_f C^{f,i_1\dots i_{\mu-1}}_{g}
(\Omega_1,\dots,\Omega_p,\phi_{u+1},\phi_1,\dots ,\phi_u)
\nabla_{i_1}\upsilon\dots \nabla_{i_{\mu-1}}\upsilon+
\\&\Sum_{r\in R_\beta} a_r
C^{r,i_1\dots i_{\mu -1}}_{g} (\Omega_1,\dots
,\Omega_p,\phi_1,\dots ,\phi_{u+1})\nabla_{i_1}\upsilon \dots
\nabla_{i_{\mu-1}}\upsilon
\\&+\Sum_{z\in Z} a_z C^{z,i_1\dots
i_{\mu-1}}_{g} (\Omega_1,\dots ,\Omega_p,\phi_1,\dots
,\phi_u,\phi_{u+1})\nabla_{i_1}\upsilon\dots
\nabla_{i_{\mu-1}}\upsilon.
\end{split}
\end{equation}
The above holds modulo
 correction terms of length greater than the $RHS$.

\par Therefore, we make the $\nabla\upsilon$'s into $Xdiv$'s (by applying the 
last Lemma in the Appendix of \cite{alexakis1})
 in the above two
equations we deduce that:

\begin{equation}
\label{fromab}
\begin{split}
&\Sum_{u\in U_2} a_u Xdiv_{i_1}\dots Xdiv_{i_{\mu-1}}C^{u,i_1\dots
i_{\mu-1}}_{g} (\Omega_1,\dots ,\Omega_p,\phi_{u+1},\phi_1,\dots
,\phi_u)-
\\&\Sum_{t\in T_2} a_u Xdiv_{i_1}\dots Xdiv_{i_{\mu-1}}Xdiv_{i_\mu}
C^{t,i_1\dots i_{\mu-1}i_\mu}_{g} (\Omega_1,\dots ,\Omega_p,
\phi_{u+1},\phi_1,\dots,\phi_u)
\\&=\Sum_{j\in J} a_j C^{j}_{g}
(\Omega_1,\dots,\Omega_p,\phi_{u+1},\phi_1,\dots ,\phi_u)+
\\&\Sum_{r\in R_\alpha} a_r Xdiv_{i_1}\dots Xdiv_{i_{\mu-1}}
C^{r,i_1\dots i_{\mu -1}}_{g} (\Omega_1,\dots
,\Omega_p,\phi_1,\dots ,\phi_{u+1})
\\& +\Sum_{z\in Z} a_z C^{z,i_1\dots
i_{\mu-1}}_{g} (\Omega_1,\dots ,\Omega_p,\phi_1,\dots
,\phi_u,\phi_{u+1}),
\end{split}
\end{equation}
(modulo complete contractions longer than the RHS) and also that:

\begin{equation}
\label{fromab2}
\begin{split}
&\Sum_{u\in U_1} a_u Xdiv_{i_1}\dots Xdiv_{i_{\mu-1}}C^{u,i_1\dots
i_{\mu-1}}_{g} (\Omega_1,\dots ,\Omega_p,\phi_{u+1},\phi_1,\dots
,\phi_u)-
\\&\Sum_{t\in T_1} a_u Xdiv_{i_1}\dots Xdiv_{i_{\mu-1}}Xdiv_{i_\mu}
C^{t,i_1\dots i_{\mu-1}i_\mu}_{g} (\Omega_1,\dots ,\Omega_p,
\phi_{u+1},\phi_1,\dots,\phi_u)
\\&= \Sum_{j\in J} a_j C^{j}_{g}
(\Omega_1,\dots,\Omega_p,\phi_{u+1},\phi_1,\dots ,\phi_u)+
\\&\Sum_{r\in R_\beta} a_r Xdiv_{i_1}\dots Xdiv_{i_{\mu-1}}
C^{r,i_1\dots i_{\mu -1}}_{g} (\Omega_1,\dots
,\Omega_p,\phi_1,\dots ,\phi_{u+1})
\\&+\Sum_{z\in Z} a_z C^{z,i_1\dots
i_{\mu-1}}_{g} (\Omega_1,\dots ,\Omega_p,\phi_1,\dots
,\phi_u,\phi_{u+1}),
\end{split}
\end{equation}
(modulo complete contractions longer than the RHS) where here we
have added the ``missing'' factors $\nabla\phi_g, g\in
(\vec{\kappa})_1$ (recall $(\vec{\kappa}_1)$ stands for
 the set of numbers $g$ for which $\nabla\phi_g$ is
 contracting against $\nabla^{(A)}\phi_{u+1}$ in $pre\vec{\kappa}^{+}_{simp}$)
onto the factor $\nabla^{(m)}\phi_{u+1}$ (which
has $M\ge 2$ for all the tensor fields) for all the tensor fields
and complete contractions above.

\par In view of these equations, we deduce that we may assume
$U_1\bigcup U_2=\emptyset$ in (\ref{enduo}). Then, in order to
show our Lemma \ref{en9ousiasmos}, we only have to show that we
can write:

\begin{equation}
\label{enduoweak}
\begin{split}
&\Sum_{u\in U} a_u Xdiv_{i_1}\dots Xdiv_{i_a}C^{u,i_1\dots
i_a}_{g} (\Omega_1,\dots ,\Omega_p,\phi_{u+1},\phi_1,\dots
,\phi_u)+
\\&\Sum_{j\in J} a_j C^j_{g}(\Omega_1,\dots
,\Omega_p,\phi_{u+1},\phi_1,\dots ,\phi_u)=
\\&\Sum_{r\in R_\gamma} a_r Xdiv_{i_1}\dots
Xdiv_{i_a}C^{r,i_1\dots i_a}_{g}
(\Omega_1,\dots ,\Omega_p,\phi_1,\dots ,\phi_u,\phi_{u+1})+
\\&\Sum_{j\in J} a_j C^j_{g}(\Omega_1,\dots
,\Omega_p,\phi_1,\dots ,\phi_{u+1})\\&+ \Sum_{z\in Z} a_z
C^{z,i_1\dots i_{\mu-1}}_{g} (\Omega_1,\dots
,\Omega_p,\phi_1,\dots ,\phi_u,\phi_{u+1});
\end{split}
\end{equation}
(recall that by definition $a\ge \mu$). 

\par In order to see this, we firstly
 again apply the eraser to the factors $\nabla\phi_g, g\in (\vec{\kappa}_1)$ in (\ref{enduo}).
Then, we pick out the tensor fields $C^{u,i_1\dots i_a}_{g}$ with
the smallest
 number $a=\delta$ of free indices,
where $\delta\ge \mu$.  We suppose they are indexed in
$U_\delta\subset U$. We will then apply Corollary 1 in \cite{alexakis1}
to the above, but first we will make a small note regarding the
 potential appearance of terms in
one of the ``forbidden forms''. By  the definition of $\sum_{u\in U}\dots$,
the only 
way that terms in (\ref{enduoweak}) can be ``forbidden'' for
 Corollary 1 in \cite{alexakis4} is if they have $\delta>\mu$. 
Thus in that case, we apply the 
  ``weaker version'' of the Proposition \ref{giade}, from the Appendix in \cite{alexakis4}; the
  correction terms that we obtain are of the form we require.

\par Now, in the remaining cases where no tensor field appearing
in $U_\delta$ is ``forbidden'', we use our inductive assumption of
Corollary 1 in \cite{alexakis4} and deduce that there exists some linear
combination of acceptable $(\delta+1)$-tensor fields (indexed
 in $T_\delta$ below with simple
character $\overline{\vec{\kappa}}$, so that:

\begin{equation}
\label{froma3}
\begin{split}
&\Sum_{u\in U_\delta} a_u C^{u,i_1\dots i_\delta}_{g}
(\Omega_1,\dots ,\Omega_p,\phi_{u+1},\phi_1,\dots
,\phi_u)\nabla_{i_1}\upsilon\dots \nabla_{i_\delta}\upsilon -
\\& Xdiv_{i_{\delta+1}}\Sum_{t\in T_\delta} a_t C^{t,i_1\dots
i_\delta,i_{\delta+1}}_{g} (\Omega_1,\dots
,\Omega_p,\phi_{u+1},\phi_1,\dots
,\phi_u)\nabla_{i_1}\upsilon\dots \nabla_{i_\delta}\upsilon+
\\&\Sum_{f\in F} a_f C^{f,i_1\dots i_{\mu-1}}_{g}
(\Omega_1,\dots,\Omega_p,\phi_{u+1},\phi_1,\dots ,\phi_u)
\nabla_{i_1}\upsilon\dots \nabla_{i_{\mu-1}}\upsilon =0,
\end{split}
\end{equation}
modulo complete contractions of length $\ge\sigma +u+1$.
Thus, since the above must hold formally, we deduce that:

\begin{equation}
\label{froma3b}
\begin{split}
&\Sum_{u\in U_\delta} a_u C^{u,i_1\dots i_\delta}_{g}
(\Omega_1,\dots ,\Omega_p,\phi_{u+1},\phi_1,\dots
,\phi_u)\nabla_{i_1}\upsilon\dots \nabla_{i_\delta}\upsilon -
\\& Xdiv_{i_{\delta+1}}\Sum_{t\in T_\delta} a_tC^{t,i_1\dots
i_\delta,i_{\delta+1}}_{g} (\Omega_1,\dots
,\Omega_p,\phi_{u+1},\phi_1,\dots
,\phi_u)\nabla_{i_1}\upsilon\dots \nabla_{i_\delta}\upsilon =
\\&\Sum_{r\in R_\gamma} a_r C^{r,i_1\dots i_\delta}_{g}
(\Omega_1,\dots ,\Omega_p,\phi_{u+1},\phi_1,\dots ,\phi_u)
\nabla_{i_1}\upsilon\dots \nabla_{i_\delta}\upsilon+
\\&\Sum_{z\in Z} a_z C^{z,i_1\dots
i_{\delta}}_{g} (\Omega_1,\dots ,\Omega_p,\phi_1,\dots
,\phi_u,\phi_{u+1})\nabla_{i_1}\upsilon\dots
\nabla_{i_{\delta}}\upsilon,
\end{split}
\end{equation}
modulo longer complete contractions.
Hence, as before, we add the missing factors $\nabla\phi_g$
 onto the factor $\nabla^{(A)}\phi_{u+1}$, 
$A\ge 2$ and make the $\nabla\upsilon$s into $Xdiv$s 
(applying the last Lemma in the Appendix of \cite{alexakis1})
to deduce that:

\begin{equation}
\label{froma3b'}
\begin{split}
&\Sum_{u\in U_\delta} a_u Xdiv_{i_1}\dots Xdiv_{i_\delta}
C^{u,i_1\dots i_\delta}_{g} (\Omega_1,\dots ,\Omega_p,\phi_{u+1},
\phi_1,\dots,\phi_u) -
\\& Xdiv_{i_1}\dots Xdiv_{i_\delta} Xdiv_{i_{\delta+1}}
\Sum_{t\in T_\delta} a_tC^{t,i_1\dots i_\delta,i_{\delta+1}}_{g}
(\Omega_1,\dots ,\Omega_p,\phi_{u+1},\phi_1,\dots ,\phi_u) 
\\&=\Sum_{f\in F} a_f C^{f,i_1\dots i_{\mu-1}}_{g}
(\Omega_1,\dots,\Omega_p,\phi_{u+1},\phi_1,\dots ,\phi_u)
\nabla_{i_1}\upsilon\dots \nabla_{i_{\mu-1}}\upsilon+
\\&\Sum_{r\in R_\gamma} a_r Xdiv_{i_1}\dots Xdiv_{i_\delta}
C^{r,i_1\dots i_\delta}_{g} (\Omega_1,\dots
,\Omega_p,\phi_{u+1},\phi_1,\dots ,\phi_u)+
\\&\Sum_{z\in Z} a_z C^{z,i_1\dots
i_{\mu-1}}_{g} (\Omega_1,\dots ,\Omega_p,\phi_1,\dots
,\phi_u,\phi_{u+1}).
\end{split}
\end{equation}

\par Therefore, by iteratively repeating this step we may assume
that $U=\emptyset$ and we are
 reduced to showing that:

\begin{equation}
\label{resolution}
\begin{split}
&\Sum_{j\in J} a_j C^j_{g}(\Omega_1,\dots
,\Omega_p,\phi_{u+1},\phi_1,\dots ,\phi_u)=
\Sum_{j\in J} a_j C^j_{g}(\Omega_1,\dots
,\Omega_p,\phi_1,\dots ,\phi_{u+1})\\&+ \Sum_{z\in Z} a_z
C^{z,i_1\dots i_{\mu-1}}_{g} (\Omega_1,\dots
,\Omega_p,\phi_1,\dots ,\phi_u,\phi_{u+1}).\end{split}
\end{equation}

\par But this follows easily: First of all, we pick out each
factor $\nabla^{(A)}_{r_1\dots r_A}\Omega_h$ (including
$\nabla^{(A)}\phi_{u+1}$) and we pull to the left all the indices
that are contracting against a factor $\nabla\phi_f$. We can do
this modulo introducing correction terms as in the right hand side
of (\ref{resolution}). So now  for each $C^j_{g}$, we have that
the factors $\nabla^{(A)}_{r_1\dots r_A}\Omega_h$ have the
property that their indices that are contracting against factors
$\nabla\phi_f$ are all pulled out to the left. Moreover, since we
are dealing with complete contractions with the same weak
character, we may speak of {\it the set} of numbers $A(h)=\{a_1,\dots
,a_{b_h}\}$, for which the factors $\nabla\phi_{a_1},\dots
,\nabla\phi_{a_{b_h}}$ are contracting against
$\nabla^{(A)}\Omega_h$, for each $h,1\le h\le p+1$.

\par Now, we arbitrarily pick out an ordering for each set
$A(h)$. Modulo introducing a linear combination
$$\Sum_{j\in J} a_j C^j_{g}(\Omega_1,\dots
,\Omega_p,\phi_1,\dots ,\phi_{u+1})+ \Sum_{z\in Z} a_z C^z_{g}
(\Omega_1,\dots ,\Omega_p,\phi_1,\dots ,\phi_{u+1}),$$ we may
assume that the factors $\nabla\phi_{a_1},\dots
,\nabla\phi_{a_{b_h}}$ are contracting against the left $b_h$
indices of $\nabla^{(A)}\Omega_h$ {\it in the order that we have
picked}.

\par Now, we define $Sym[C^j_{g}(\Omega_1,\dots
,\Omega_p,\phi_{u+1},\phi_1,\dots ,\phi_u)]$ to stand for the complete
 contraction that is obtained from
$C^j_{g}(\Omega_1,\dots,\Omega_p,\phi_{u+1},\phi_1,\dots ,
\phi_{u})$ by symmetrizing over the indices in each factor
$\nabla^{(A)}\Omega_h$ {\it that are not contracting against a
factor $\nabla\phi_h$}.

\par By just applying the equation (\ref{curvature}) we then deduce that:

\begin{equation}
\label{resolution2}
\begin{split}
&\Sum_{j\in J} a_j C^j_{g}(\Omega_1,\dots
,\Omega_p,\phi_{u+1},\phi_1,\dots ,\phi_u)\\&= \Sum_{j\in J} a_j
Sym[C^f_{g}(\Omega_1,\dots ,\Omega_p,\phi_{u+1},\phi_1,\dots
,\phi_u)]
\\&+\Sum_{j\in J} a_j C^j_{g}(\Omega_1,\dots
,\Omega_p,\phi_1,\dots ,\phi_{u+1})+ \Sum_{z\in Z} a_z C^z_{g}
(\Omega_1,\dots ,\Omega_p,\phi_1,\dots ,\phi_{u+1}).
\end{split}
\end{equation}

\par Finally, we use the eraser to deduce that:

\begin{equation}
\label{resolution3}
\begin{split}
&\Sum_{j\in J} a_j Sym[C^j_{g}(\Omega_1,\dots
,\Omega_p,\phi_{u+1},\phi_1,\dots ,\phi_u)]=
\\&\Sum_{z\in Z} a_z C^z_{g} (\Omega_1,\dots ,\Omega_p,
\phi_1,\dots,\phi_u,\phi_{u+1}).
\end{split}
\end{equation}
$\Box$

\par The proof of Lemma \ref{en9ousiasmos2} is entirely identical.
We just do not have the sublinear combinations indexed in
$U_1\bigcup U_2$, hence we just iteratively apply (\ref{froma3b'})
and then (\ref{resolution}). $\Box$
\newline

\par In conclusion, in view of equations (\ref{akinola}),
(\ref{akinola2}), (\ref{akinola3}) and  Lemmas  \ref{en9ousiasmos}
and \ref{en9ousiasmos2}, we have shown the following:
\newline

{\it Conclusions:} In the setting of Lemma \ref{pool2} when there are two
internal free indices in the selected factor
$\nabla^{(m)}R_{ijkl}$ in $\vec{L^z}'$, $z\in Z'_{Max}$, we have derived
an equation:

\begin{equation}
\label{akinola'}
\begin{split}
&CurvTrans[L_{g}(\Omega_1,\dots ,\Omega_p,\phi_1,\dots ,\phi_u)]=
\\& \Sum_{z\in Z'_{Max}}\Sum_{l\in L^z} a_l
\Sum_{h=1}^{2k_l+1}Xdiv_{i_1}\dots \hat{Xdiv}_{i_h}\dots
 Xdiv_{i_\mu}C^{l,i_1\dots i_\mu}_{g} (\Omega_1,\dots ,
\Omega_p,\phi_1,\dots,\phi_u)
\\&\nabla_{i_1}\phi_{u+1}+\Sum_{t\in T} a_t Xdiv_{i_1}\dots Xdiv_{i_{\mu-1}}
C^{t,i_1\dots i_\mu}_{g} (\Omega_1,\dots ,\Omega_p,\phi_1,\dots
,\phi_u)\nabla_{i_\mu}\phi_{u+1}
\\&+\Sum_{r\in R_\alpha} a_r Xdiv_{i_1}\dots Xdiv_{i_{\mu-1}}
C^{r,i_1\dots i_{\mu -1}}_{g} (\Omega_1,\dots
,\Omega_p,\phi_1,\dots ,\phi_{u+1})
\\&+\Sum_{r\in
R_\beta} a_r Xdiv_{i_1}\dots Xdiv_{i_{\mu-1}}C^{r,i_1\dots i_{\mu
-1}}_{g} (\Omega_1,\dots ,\Omega_p,\phi_1,\dots ,\phi_{u+1})
\\&\Sum_{r\in R_\gamma} a_r Xdiv_{i_1}\dots Xdiv_{i_a}C^{r,i_1\dots i_a}_{g}
(\Omega_1,\dots ,\Omega_p,\phi_1,\dots ,\phi_u,\phi_{u+1})+
\\&\Sum_{j\in J} a_j C^j_{g}(\Omega_1,\dots
,\Omega_p,\phi_1,\dots ,\phi_{u+1})+ \sum_{z\in Z} a_z
C^z_g(\Omega_1,\dots,\Omega_p,\phi_1,\dots,\phi_{u+1}),
\end{split}
\end{equation}
(using the notational conventions of equation (\ref{akinola}) and
Lemma \ref{en9ousiasmos}).

\par In the setting of Lemma \ref{pool2} when there is one
internal free index in the selected factor $\nabla^{(m)}R_{ijkl}$
in $\vec{L^z}',z\in Z'_{Max}$, we have derived an equation:

\begin{equation}
\label{akinola2'}
\begin{split}
&CurvTrans[L_{g}(\Omega_1,\dots ,\Omega_p,\phi_1,\dots ,\phi_u)]=
\\&\Sum_{z\in Z'_{Max}}\Sum_{l\in L^z} a_l
\Sum_{i_h\in I_{*,l}} Xdiv_{i_1}\dots \hat{Xdiv}_{i_h}\dots
Xdiv_{i_\mu}C^{l,i_1\dots i_\mu}_{g} (\Omega_1,\dots
,\Omega_p,\phi_1, \dots,\phi_u)
\\&\nabla_{i_h}\phi_{u+1}+\Sum_{t\in T} a_t Xdiv_{i_1}\dots Xdiv_{i_{\mu-1}}
C^{t,i_1\dots i_\mu}_{g} (\Omega_1,\dots ,\Omega_p,\phi_1,\dots
,\phi_u)\nabla_{i_\mu}\phi_{u+1}
\\&+\Sum_{r\in
R_\beta} a_r Xdiv_{i_1}\dots Xdiv_{i_{\mu-1}}C^{r,i_1\dots i_{\mu
-1}}_{g} (\Omega_1,\dots ,\Omega_p,\phi_1,\dots ,\phi_{u+1})+
\\&\Sum_{r\in R_\gamma} a_r Xdiv_{i_1}\dots Xdiv_{i_a}C^{r,i_1\dots i_a}_{g}
(\Omega_1,\dots ,\Omega_p,\phi_1,\dots ,\phi_u,\phi_{u+1})+
\\&\Sum_{j\in J} a_j C^j_{g}(\Omega_1,\dots
,\Omega_p,\phi_1,\dots ,\phi_{u+1})+ \sum_{z\in Z} a_z
C^z_g(\Omega_1,\dots,\Omega_p,\phi_1,\dots,\phi_{u+1}).
\end{split}
\end{equation}

\par In
the setting of Lemma \ref{pskovb} with a selected factor
$\nabla^{(m)}R_{ijkl}$, we have derived an equation:

\begin{equation}
\label{akinola3'}
\begin{split}
&CurvTrans[L_{g}(\Omega_1,\dots ,\Omega_p,\phi_1,\dots ,\phi_u)]=
\\&\Sum_{r\in R_\gamma} a_r Xdiv_{i_1}\dots Xdiv_{i_a}C^{r,i_1\dots i_a}_{g}
(\Omega_1,\dots ,\Omega_p,\phi_1,\dots ,\phi_u,\phi_{u+1})+
\\& \Sum_{j\in J} a_j C^j_{g}(\Omega_1,\dots
,\Omega_p,\phi_1,\dots ,\phi_u,\phi_{u+1})+
\\&\sum_{z\in Z} a_z C^z_g(\Omega_1,\dots,\Omega_p,\phi_1,\dots,\phi_{u+1}),
\end{split}
\end{equation}
(using the notational conventions spelled out in Lemma
\ref{en9ousiasmos2}).

\subsection{A study of the sublinear combination $CurvTrans[L_{g}]$ in the
context of Lemmas \ref{zetajones} and \ref{pskovb} (when the
selected factor is of the form $S_{*}\nabla^{(\nu)}R_{ijkl}$).}
\label{otherplaces}

In our study of $CurvTrans[L_{g}]$ it will be important to recall
our inductive assumptions on Proposition \ref{giade}.
 Recall that in the settings where we are inductively assuming
Proposition \ref{giade}, we may also apply Corollary 1 from \cite{alexakis1} and
 Lemma 4.6 from \cite{alexakis4}.

\par Recall the $(u+1)$-simple character
$\vec{\kappa}^{+}_{simp}$ (which corresponds to contractions with
$\sigma+u+1$ factors), and also the $(u+1)$-simple character
$pre\vec{\kappa}^{+}_{simp}$ that corresponds to contractions with
$\sigma+u$ factors (see the paper \cite{alexakis3} for a precise definition 
of simple character, and the Definition \ref{laborat} for the 
definition of $pre\vec{\kappa}^{+}_{simp}$).

 We then have denoted by:

 $$\Sum_{j\in J} a_j C^j_{g}(\Omega_1,\dots
 ,\Omega_p,\phi_1,\dots ,\phi_{u+1})$$
a generic linear combination of complete contractions of length
$\ge\sigma +u+1$ which is simply subsequent to
$\vec{\kappa}^{+}_{simp}$.

\par The aim of this subsection is to derive the equation (\ref{keepinminda}) below. 
However, in order not to burden the reader with many new definitions 
from the outset, we will commence with some simple 
calculations and then introduce the necessary 
 notation needed  along the way. 
Thus, rather than stating the objective of this subsection from 
the outset, we will reach it at the end of this section as a consequence of some 
(seemingly unmotivated) calculations.
\newline

\par Now, in this subsection the selected factor is of the form
$T=S_{*}\nabla^{(\nu)}R_{ijkl}$. We denote by $\nabla\phi_{Min}$
the factor that is contracting against the index ${}_i$ of $T$ in
$\vec{\kappa}_{simp}$. We recall that the $CurvTrans[L_g]$ stands
for the sublinear combination in $Image^{1,+}_{\phi_{u+1}}[L_g]$
of complete contractions with length $\sigma+u$ and a weak
character $Weak(pre\vec{\kappa}^{+}_{simp})$. Therefore, the factor
$\nabla\phi_{Min}$ {\it must} be contracting against the factor
$\nabla^{(A)}\phi_{u+1}$, for each complete contraction in
$CurvTrans[L_g]$.

\par  We observe that the sublinear combination in
any term \\$Image^{1,+}_{\phi_{u+1}}[\Sum_{l\in L} a_l Xdiv_{i_1}\dots
Xdiv_{i_a} C^{l,i_1\dots i_a}_{g}(\Omega_1,\dots
,\Omega_p,\phi_1,\dots ,\phi_u)]$ for which $\nabla\phi_{Min}$
{\it is}
 contracting against the factor $\nabla^{(A)}\phi_{u+1}$ can only arise by
replacing the crucial factor $S_{*}\nabla^{(\nu)}_{r_1\dots r_n}
R_{ijkl}$ by either $S_{*}\nabla^{(\nu+2)}_{r_1\dots r_\nu ik}
\phi_{u+1}g_{jl}$ or $-S_{*}\nabla^{(\nu+2)}_{r_1\dots r_\nu
il}\phi_{u+1}g_{jk}$ (here $S_{*}$ again stands for symmetrization
over the indices ${}_{r_1},\dots ,{}_{r_\nu},{}_j$). Accordingly,
we denote by
\\$CurvTrans^I [\Sum_{l\in L} a_l Xdiv_{i_1}\dots Xdiv_{i_a}
C^{l,i_1\dots i_a}_{g}(\Omega_1,\dots ,\Omega_p,\phi_1,\dots
,\phi_u)]$,
$$CurvTrans^{II}[\Sum_{l\in L} a_l
Xdiv_{i_1}\dots Xdiv_{i_a}
C^{l,i_1\dots i_a}(\Omega_1,\dots ,\Omega_p,\phi_1,\dots ,\phi_u)]$$ the
sublinear combinations that arise by making these substitutions.
\newline

\par Now, we carefully study the linear combinations
\\$CurvTrans^I [\Sum_{l\in L} a_l Xdiv_{i_1}\dots Xdiv_{i_a}
C^{l,i_1\dots i_a}(\Omega_1,\dots ,\Omega_p,\phi_1,\dots
,\phi_u)]$, \\$CurvTrans^{II} [\Sum_{l\in L} a_l Xdiv_{i_1}\dots
Xdiv_{i_a} C^{l,i_1\dots i_a}_{g}(\Omega_1,\dots
,\Omega_p,\phi_1,\dots ,\phi_u)]$
 and \\$CurvTrans[C^j_g(\Omega_1,\dots
,\Omega_p,\phi_1,\dots ,\phi_u)]$. In
order to state our next claim, we introduce some notational
conventions. Let the total number of factors $\nabla\phi_h$
contracting against the selected factor in $\vec{\kappa}_{simp}$ be
$\pi$, and in particular let those factors be
$\nabla\tilde{\phi}_{Min},\nabla\phi'_{e_1},\dots
,\nabla\phi'_{e_{\pi-1}}$. Also, {\it in the setting of Lemma
\ref{zetajones}} we consider the total number of free indices in
the selected factor, for each $C^{l,i_1\dots i_\mu}_g$, $l\in
\bigcup_{z\in Z'_{Max}} L^z$. We denote that number by $Free(Max)$.
\newline

\begin{definition}
\label{maciunas} Recall the simple character
$pre\vec{\kappa}^{+}_{simp}$.

We denote by
 $$\Sum_{j\in J} a_j C^j_{g}(\Omega_1,\dots
 ,\Omega_p,\phi_{u+1},\phi_1,\dots ,\phi_{u})$$
a generic linear combination of complete contractions of length
$\sigma +u$ which are simply subsequent to
$pre\vec{\kappa}^{+}_{simp}$. (Thus $\phi_{u+1}$ is regarded here
as a factor $\Omega_h$).

Moreover, we denote by $\Sum_{d\in D^\sharp} a_d C^d_{g}
(\Omega_1,\dots ,\Omega_p,\phi_{u+1},\phi_1,\dots ,\phi_u)$ a
generic linear combination of complete contractions with length
$\sigma+u$, a weak $u$-simple character
$Weak(pre\vec{\kappa}^{+}_{simp})$ and with $\nabla\phi_{Min}$
contracting against the first index in the factor
$\nabla^{(P)}\phi_{u+1}$, where we additionally require that $P\ge
3$.

\par We will also let $\sum_{p\in P} a_p C^{p,i_1\dots i_a,i_{*}}_g
\nabla_{i_{*}}\phi_{u+1}$ be a generic linear combination of
$a$-tensor fields ($a\ge \mu$) with length $\sigma+u+1$ and 
the following additional properties: In the setting of Lemma \ref{zetajones} 
they must have a
$(u+1)$-simple character $\vec{\kappa}^{+}_{simp}$; 
in the setting of Lemma \ref{pskovb}
they must have a $u$-simple character $\vec{\kappa}_{simp}$
and a weak $(u+1)$-character $Weak(\vec{\kappa}_{simp}^{+})$.

\par Now, a few definitions that are only applicable
when $\pi=1$.\footnote{See the  notation above.}
We denote by $\Sum_{d\in D} a_d C^{d,i_1,\dots ,i_a}_{g}
(\Omega_1,\dots ,\Omega_p,\phi_{u+1},\phi_1,\dots ,\phi_u)$ a
generic linear combination of acceptable $a$-tensor fields with
 length $\sigma +u$, $a\ge \mu$, with simple character
$pre\vec{\kappa}^{+}_{simp}$, and with an expression
 $\nabla^{(2)}_{ij}\phi_{u+1}\nabla^i\tilde{\phi}_{Min}$. If
 $a=\mu$ then we additionnaly require that if we formally replace 
 the expression $\nabla^{(2)}_{ij}\phi_{u+1}\nabla^i\tilde{\phi}_{Min}$ by a factor 
$\nabla_j Y$ then the resulting tensor field is not forbidden in the sense 
of Lemma 4.6 in \cite{alexakis4}.

\par Finally, only in the setting of Lemma \ref{zetajones},
we will denote by
\\$\sum_{d\in D_{nc}} a_d C^{d,i_1\dots i_{\mu-1}}_g
(\Omega_1,\dots,\Omega_p,\phi_{u+1},\phi_1,\dots)$ a generic linear
combination of tensor fields with length $\sigma +u$, a factor
$\nabla^{(2)}_{ij}\phi_{u+1}\nabla^i\phi_{Min}$
 (${}_j$ is not a free index) and simple character
$pre\vec{\kappa}^{+}_{simp}$, and with one of the free indices
${}_{i_1},\dots,{}_{\mu-1}$ being a 
derivative index. If this index belongs to a factor
$\nabla^{(B)}\Omega_h$ tehn $B\ge 3$.
\end{definition}

\par Armed with all the above notational conventions, we refer
back to (\ref{tsabes}), and we set out
to understand the form of the sublinear combinations:

$$CurvTrans[C^j_{g}(\Omega_1,\dots ,\Omega_p,\phi_1,\dots ,\phi_u)]$$
and $CurvTrans[Xdiv_{i_1}\dots Xdiv_{i_a} C^{l,i_1\dots
i_a}_{g}(\Omega_1,\dots ,\Omega_p,\phi_1,\dots ,\phi_u)]$.
\newline

\begin{lemma}

\begin{equation}
\label{bammena2}
\begin{split}
&CurvTrans[C^j_{g}(\Omega_1,\dots ,\Omega_p, \phi_1,\dots
,\phi_u)]=\Sum_{d\in D^\sharp} a_d C^d_{g}(\Omega_1,\dots ,
\Omega_p,\phi_1,\dots ,\phi_{u+1})+
\\&\Sum_{j\in J} a_j C^j_{g}(\Omega_1,\dots ,\Omega_p,
\phi_1,\dots ,\phi_{u+1})+\Sum_{z\in Z} a_z C^z_{g}(\Omega_1,\dots
,\Omega_p, \phi_1,\dots ,\phi_{u+1})
\\& +\Sum_{j\in J} a_j C^j_{g}
(\Omega_1,\dots ,\Omega_p,\phi_{u+1},\phi_1,\dots ,\phi_u).
\end{split}
\end{equation}
\end{lemma}

\begin{lemma}
\label{morena} Consider any tensor field $C^{l,i_1\dots i_a}_{g}
(\Omega_1,\dots ,\Omega_p,\phi_1,\dots , \phi_u)$, $l\in L$  where none of
the indices ${}_k,{}_l$ in the crucial factor
 are free indices. Consider the special factor
 $S_{*}\nabla^{(\nu)}_{r_1\dots r_\nu}R_{ijkl}$ and
 denote by $\rho$ the number of the indices
${}_{r_1},\dots ,{}_{r_\nu},{}_j$ that are free in $C^{l,i_1\dots i_a}_{g}$, and
we denote them by ${}_{i_1},\dots ,{}_{i_\rho}$, for convenience. We claim that:

\begin{equation}
\label{thangave}
\begin{split}
&CurvTrans^I\{ Xdiv_{i_1}\dots Xdiv_{i_a} C^{l,i_1\dots i_a}_{g}
(\Omega_1,\dots ,\Omega_p,\phi_1,\dots ,\phi_u)\}+
\\&CurvTrans^{II}\{ Xdiv_{i_1}\dots
Xdiv_{i_a} C^{l,i_1\dots i_a}_{g}(\Omega_1,\dots
,\Omega_p,\phi_1,\dots , \phi_u)\}=
\\& -\frac{1}{\nu+1}\Sum_{y=1}^\rho Xdiv_{i_1}\dots \hat{Xdiv}_{i_y}\dots Xdiv_{i_a}
 C^{l,i_1\dots i_a}_{g} (\Omega_1,\dots
,\Omega_p,\phi_1,\dots ,\phi_u) \nabla_{i_y}\phi_{u+1}
\\&+\Sum_{p\in P} a_p Xdiv_{i_1}\dots Xdiv_{i_a}
C^{p,i_1,\dots i_a,i_{*}}_{g} (\Omega_1,\dots
,\Omega_p,\phi_1,\dots ,\phi_u)\nabla_{i_{*}} \phi_{u+1}
\\&+\Sum_{j\in J} a_j C^j_{g}(\Omega_1,\dots ,\Omega_p,
\phi_1,\dots ,\phi_{u+1})
\\&+\Sum_{d\in D} a_d Xdiv_{i_1}\dots Xdiv_{i_a}C^{d,i_1\dots i_a}_{g}(\Omega_1,\dots
,\Omega_p,\phi_{u+1}, \phi_1,\dots ,\phi_u)
\\& +\Sum_{d\in D^\sharp} a_d C^d_{g}
(\Omega_1,\dots ,\Omega_p,\phi_{u+1},\phi_1,\dots ,\phi_u)+
\Sum_{z\in Z} a_z C^z_{g}
(\Omega_1,\dots ,\Omega_p,\phi_1,\dots ,\phi_{u+1}).
\end{split}
\end{equation}
\end{lemma}

{\it Note:} Let us observe that in the setting of Lemma
\ref{zetajones}, the tensor fields $C^{l,i_1\dots i_a}_{g}
(\Omega_1,\dots ,\Omega_p,\phi_1,\dots ,\phi_u)
\nabla_{i_y}\phi_{u+1}$ are $(u+1)$-simply subsequent to
$\vec{\kappa}^{+}_{simp}$. We will denote the sublinear
combination:

$$-\Sum_{y=1}^\rho Xdiv_{i_1}\dots \hat{Xdiv}_{i_y}\dots Xdiv_{i_a}
 C^{l,i_1\dots i_a}_{g} (\Omega_1,\dots
,\Omega_p,\phi_1,\dots ,\phi_u) \nabla_{i_y}\phi_{u+1}$$ by
$Leftover_{\phi_{u+1}}[Xdiv_{i_1}\dots Xdiv_{i_a}
 C^{l,i_1\dots i_a}_{g}
(\Omega_1,\dots ,\Omega_p,\phi_1,\dots ,\phi_u)]$.
\newline

{\it Proof of the two Lemmas above:} The proof follows straightforwardly by applying
the transformation laws (\ref{curvtrans}) and (\ref{curvature}).
$\Box$
\newline

\par Now, we focus on the case of the tensor fields
$C^{l,i_1\dots i_a}_{g} (\Omega_1,\dots ,\Omega_p,\phi_1,
\dots,\phi_u)$  where one of the indices ${}_k,{}_l$  in the selected
factor is a free index (with no loss of generality we assume that
${}_k$ is the free index ${}_{i_1}$ and ${}_l$ is
 not a free index). For convenience, we assume that the
  rest of the indices in $S_{*}\nabla^{(\nu)}_{r_1\dots r_\nu}
R_{ijkl}$ that are free are precisely ${}_{i_2},\dots
,{}_{i_{\rho+1}}$.
 For each such tensor field we will denote by $\epsilon$
 the number of indices ${}_{r_1},\dots ,{}_{r_\nu},{}_j$ in the crucial
 factor $S_{*}\nabla^{(\nu)}_{r_1\dots r_\nu}R_{ijkl}$
that are neither free nor contracting against a factor
$\nabla\phi_{e_1},\dots \nabla\phi_{e_{\pi-1}}$.

\begin{lemma}
\label{n>pi-1}
With the notational conventions above we claim:
\begin{equation}
\label{thangave5}
\begin{split}
&CurvTrans^I \{Xdiv_{i_1}\dots Xdiv_{i_a} C^{l,i_1\dots i_a}_{g}
(\Omega_1,\dots ,\Omega_p,\phi_1,\dots ,\phi_u)\}+
\\&CurvTrans^{II}\{Xdiv_{i_1}\dots
Xdiv_{i_a} C^{l,i_1\dots i_a}_{g} (\Omega_1,\dots
,\Omega_p,\phi_1,\dots ,\phi_u)\}=
\\&-\frac{1}{\nu+1}\Sum_{y=1}^\rho
Xdiv_{i_1}\dots \hat{Xdiv}_{i_{y+1}}\dots Xdiv_{i_a} C^{l,i_1\dots
i_a}_{g} (\Omega_1,\dots ,\Omega_p,\phi_1,\dots
,\phi_u)\nabla_{i_{y+1}}\phi_{u+1}
\\&+\Sum_{p\in P} a_p Xdiv_{i_1}\dots Xdiv_{i_a}
C^{p,i_1,\dots i_a,i_{*}}_{g} (\Omega_1,\dots
,\Omega_p,\phi_1,\dots ,\phi_u)\nabla_{i_{*}} \phi_{u+1}
\\&-\frac{\epsilon}{\nu+1}Xdiv_{i_2}\dots Xdiv_{i_a}C^{l,i_1\dots i_a}_{g} (\Omega_1,\dots
,\Omega_p,\phi_1,\dots ,\phi_u)\nabla_{i_1}\phi_{u+1}+
\\&\Sum_{d\in D^\sharp} a_d C^d_{g} (\Omega_1,\dots
,\Omega_p,\phi_{u+1},\phi_1,\dots ,\phi_u)+
\\&(\Sum_{d\in D_{nc}} a_d Xdiv_{i_1}\dots Xdiv_{i_{a-1}}
C^{d,i_1,\dots ,i_{a-1}}_{g} (\Omega_1,\dots
,\Omega_p,\phi_1,\dots ,\phi_u,\phi_{u+1}))+
\\&\Sum_{z\in Z} a_z C^z_{g} (\Omega_1,\dots
,\Omega_p,\phi_1,\dots ,\phi_u,\phi_{u+1});
\end{split}
\end{equation}
here here the sublinear combination $\sum_{d\in D_{nc}}a_d\dots$ arises
only in the setting of Lemma \ref{zetajones}, when $Free(Max)>1$.
\end{lemma}

\begin{definition}
\label{orizw}
\par We will denote the linear combination
\begin{equation}
\label{xristkpan} \begin{split} &-\frac{1}{\nu+1} \Sum_{y=1}^\rho
Xdiv_{i_1}\dots \hat{Xdiv}_{i_{y+1}}\dots
 Xdiv_{i_a} C^{l,i_1\dots i_a}_{g} (\Omega_1,\dots ,
\Omega_p,\phi_1,\dots,\phi_u)\nabla_{i_{y+1}}\phi_{u+1}
\\&-\frac{\epsilon}{\nu+1}Xdiv_{i_2}\dots Xdiv_{i_a}C^{l,i_1\dots
i_a}_{g} (\Omega_1,\dots ,\Omega_p,\phi_1,\dots
,\phi_u)\nabla_{i_1}\phi_{u+1}
\end{split}
\end{equation}
 by
$$Leftover_{\phi_{u+1}}[Xdiv_{i_1}\dots Xdiv_{i_a}
C^{l,i_1\dots i_a}_{g} (\Omega_1,\dots ,\Omega_p,\phi_1,
\dots,\phi_u)].$$
\end{definition}

{\it Proof of Lemma \ref{n>pi-1}:} All the above claims follow by
the definitions. Only for
 (\ref{thangave5}) we must
 also use the equation (\ref{curvature}) also, for (\ref{thangave5}) we
 use the first Bianchi identity. $\Box$
\newline

\par In conclusion, we have shown that in the setting of Lemmas \ref{zetajones}
and \ref{pskovb} (when 
the selected factor is of the from $S_{*}\nabla^{(\nu)}R_{ijkl}$),
$CurvTrans[L_{g}]$ can now be expressed as:

\begin{equation}
\label{keepinmind}
\begin{split}
&\Sum_{l\in L} a_l CurvTrans[Xdiv_{i_1}\dots
Xdiv_{i_a}C^{l,i_1\dots i_a}_{g}(\Omega_1\dots
,\Omega_p,\phi_1,\dots ,\phi_u)]+
\\& \Sum_{j\in J} a_j CurvTrans[C^j_{g}
(\Omega_1\dots ,\Omega_p,\phi_1,\dots ,\phi_u)]=
\\&\Sum_{l\in L} a_l Leftover_{\phi_{u+1}}
[Xdiv_{i_1}\dots Xdiv_{i_a}C^{l,i_1\dots i_a}_{g} (\Omega_1\dots
,\Omega_p,\phi_1,\dots ,\phi_u)]+
\\& \Sum_{d\in D^\sharp} a_d C^d_{g}(\Omega_1,\dots,
\Omega_p,\phi_{u+1},\phi_1,\dots ,\phi_u)+ \Sum_{j\in J} a_j
C^j_{g}(\Omega_1,\dots, \Omega_p,\phi_{u+1},\phi_1,\dots ,\phi_u)
\\&+(\Sum_{d\in D_{nc}} a_d Xdiv_{i_1}\dots Xdiv_{i_{a-1}}
C^{d,i_1,\dots ,i_{a-1}}_{g} (\Omega_1,\dots
,\Omega_p,\phi_1,\dots ,\phi_u,\phi_{u+1}))+
\\&\Sum_{d\in D} a_d Xdiv_{i_1}\dots Xdiv_{i_a}C^{d,i_1\dots i_a}_{g}(\Omega_1,\dots
,\Omega_p,\phi_{u+1}, \phi_1,\dots ,\phi_u)+
\\&\Sum_{p\in P} a_p Xdiv_{i_1}\dots Xdiv_{i_a}
C^{p,i_1,\dots i_a,i_{*}}_{g} (\Omega_1,\dots
,\Omega_p,\phi_1,\dots ,\phi_u)\nabla_{i_{*}} \phi_{u+1}+
\\&\Sum_{z\in Z} a_z C^z_{g} (\Omega_1,\dots
,\Omega_p,\phi_1,\dots ,\phi_u,\phi_{u+1}).
\end{split}
\end{equation}

\par Our aim is now to ``get rid'' of all the sublinear combinations
 indexed in $D,D_{nc},D^\sharp$, modulo introducing correction terms
with $\sigma+u+1$ factors that are allowed in the conclusions of 
 Lemma \ref{zetajones} and \ref{pskovb}. 
The rest of this subsection is devoted to that goal.
\newline

In view of (\ref{tsabes}) and (\ref{keepinmind}), we derive an equation:

\begin{equation}
\label{tonkyrio=0}
\begin{split}
&0=\big{(}\Sum_{d\in D_{nc}} a_d Xdiv_{i_1}\dots Xdiv_{i_{\mu-1}}
C^{d,i_1,\dots ,i_{\mu-1}}_{g} (\Omega_1,\dots
,\Omega_p,\phi_{u+1},\phi_1,\dots ,\phi_u)\big{)}
\\&+\Sum_{d\in D} a_d Xdiv_{i_1}\dots Xdiv_{i_a}
C^{d,i_1\dots i_a}_{g}(\Omega_1,\dots,
\Omega_p,\phi_{u+1},\phi_1,\dots ,\phi_u)+
\\&\Sum_{d\in D^\sharp} a_d
C^d_{g}(\Omega_1,\dots, \Omega_p,\phi_{u+1},\phi_1,\dots ,\phi_u)+
\Sum_{j\in J} a_j C^j_{g}(\Omega_1,\dots,
\Omega_p,\phi_{u+1},\phi_1,\dots ,\phi_u),\end{split}
\end{equation}
modulo complete contractions of length $\ge\sigma +u+1$.
 We then claim:

\begin{lemma}
\label{bosnios} Refer to (\ref{tonkyrio=0}). We claim that we can write:

\begin{equation}
\label{einkid}
\begin{split}
&(\Sum_{d\in D_{nc}} a_d Xdiv_{i_1}\dots Xdiv_{i_{\mu-1}}
C^{d,i_1,\dots ,i_{\mu-1}}_{g} (\Omega_1,\dots
,\Omega_p,\phi_{u+1},\phi_1,\dots ,\phi_u))+
\\&+\Sum_{d\in D} a_d Xdiv_{i_1}\dots Xdiv_{i_a}
C^{d,i_1\dots i_a}_{g}(\Omega_1,\dots,
\Omega_p,\phi_{u+1},\phi_1,\dots ,\phi_u)+
\\&\Sum_{d\in D^\sharp} a_d
C^d_{g}(\Omega_1,\dots, \Omega_p,\phi_{u+1},\phi_1,\dots ,\phi_u)
+\Sum_{j\in J} a_j C^j_{g}(\Omega_1,\dots,
\Omega_p,\phi_{u+1},\phi_1,\dots ,\phi_u)
\\&=(\Sum_{p\in P'} a_p Xdiv_{i_1}\dots Xdiv_{i_{\mu-1}}
C^{p,i_1,\dots i_{\mu-1},i_{*}}_{g} (\Omega_1,\dots
,\Omega_p,\phi_1,\dots ,\phi_u)\nabla_{i_{*}} \phi_{u+1})+
\\&\Sum_{a>\mu-1}\Sum_{p\in P} a_p Xdiv_{i_1}\dots Xdiv_{i_a}
C^{p,i_1,\dots i_a,i_{*}}_{g} (\Omega_1,\dots
,\Omega_p,\phi_{u+1},\phi_1,\dots ,\phi_u)\nabla_{i_{*}}
\phi_{u+1}
\\&+\Sum_{j\in J} a_j C^j_{g}(\Omega_1,\dots,
\Omega_p,\phi_1,\dots ,\phi_{u+1})+ \Sum_{z\in Z} a_z C^z_{g}
(\Omega_1,\dots ,\Omega_p,\phi_1,\dots ,\phi_{u+1}).
\end{split}
\end{equation}
Here the sublinear combination $\sum_{p\in P'}\dots$ arises only
in the setting of Lemma \ref{zetajones}. In that case, it stands
for a generic linear combination of acceptable tensor fields with a $(u+1)$-simple
character $\vec{\kappa}^{+}_{simp}$ and with fewer than
$Free(Max)$ free indices in the selected (crucial) factor.
Equation (\ref{einkid}) holds modulo terms of length
$\ge\sigma+u+2$.
\end{lemma}

{\it Proof of Lemma \ref{bosnios}:}

\par We show the above via an induction. However the base case of our
induction depends on which setting we are in. In the setting of
Lemma \ref{pskovb} the linear combination $\sum_{d\in
D_{nc}}\dots$ is not present. Also, in the case of Lemma
\ref{zetajones}, if $Free(Max)=1$ then $\sum_{d\in
D_{nc}}\dots$ is not present. In those cases we may skip to after
equation (\ref{pesax2}). In the setting of Lemma \ref{zetajones}
with $Free(Max)\ge 2$ we must first ``get rid'' of the
sublinear combination $\sum_{d\in D_{nc}}\dots$.

\par So, we now assume that $D_{nc}\ne\emptyset$.
We refer to (\ref{tonkyrio=0}), and we recall that all the tensor fields
involved have a {\it fixed} simple character, which we have
denoted by $pre\vec{\kappa}^{+}_{simp}$.
Picking out the sublinear combination in (\ref{tonkyrio=0})
which consists of complete contractions with
 a factor $\nabla^{(2)}\phi_{u+1}$ we derive
 a new equation:

\begin{equation}
\label{paparides}
\begin{split}
&\Sum_{d\in D_{nc}} a_d X_{*}div_{i_1}\dots X_{*}div_{i_{\mu-1}}
C^{d,i_1,\dots ,i_{\mu-1}}_{g} (\Omega_1,\dots
,\Omega_p,\phi_{u+1},\phi_1,\dots ,\phi_u)+
\\&+\Sum_{d\in D} a_d X_{*}div_{i_1}\dots X_{*}div_{i_a}
C^{d,i_1\dots i_a}_{g}(\Omega_1,\dots,
\Omega_p,\phi_{u+1},\phi_1,\dots ,\phi_u)+
\\& \Sum_{j\in J} a_j C^j_{g}(\Omega_1,\dots,
\Omega_p,\phi_{u+1},\phi_1,\dots ,\phi_u)=0,\end{split}
\end{equation}
which holds modulo complete contractions of length $\ge\sigma +u+1$. Here
$X_{*}div_i$ stands for the sublinear combination in $Xdiv_i$
where $\nabla_i$ is not allowed to hit the factor
$\nabla^{(2)}\phi_{u+1}$. Therefore, applying the eraser
 to the factor $\nabla\phi_{Min}$, in the  above equation
and then the Lemma 4.10 from \cite{alexakis4} (which we are now inductively
assuming because we have lowered the absolute value of the
weight)\footnote{By the definition of the terms 
$\sum_{d\in D_{nc}}\dots$ there is no danger of falling 
under a ``forbidden case'' of that Lemma.} we derive that there is a linear combination of
acceptable $\mu$-tensor fields, with a simple character
$pre\vec{\kappa}^{+}_{simp}$, say

$$\Sum_{h\in H} a_h C^{h,i_1\dots ,i_\mu}_{g}(\Omega_1,\dots,
\Omega_p,\phi_{u+1},\phi_1,\dots ,\phi_u),$$ each with a factor
$\nabla^{(2)}_{ij}\phi_{u+1}\nabla^i\tilde{\phi}_{Min}$
so that:

\begin{equation}
\label{pesax}
\begin{split}
&\Sum_{d\in D_{nc}} a_d C^{d,i_1,\dots ,i_{\mu-1}}_{g}
(\Omega_1,\dots ,\Omega_p,\phi_{u+1},\phi_1,\dots
,\phi_u)\nabla_{i_1}\upsilon\dots
\nabla_{i_{\mu-1}}\upsilon\\&-X_{*}div_{i_\mu}\Sum_{h\in H} a_h
C^{h,i_1\dots ,i_\mu}_{g}(\Omega_1,\dots,
\Omega_p,\phi_{u+1},\phi_1,\dots ,\phi_u)\nabla_{i_1}\upsilon\dots
\nabla_{i_{\mu-1}}\upsilon
\\&=\Sum_{j\in J} a_j C^{j,i_1,\dots ,i_{\mu-1}}_{g}
(\Omega_1,\dots ,\Omega_p,\phi_{u+1},\phi_1,\dots
,\phi_u)\nabla_{i_1}\upsilon\dots \nabla_{i_{\mu-1}}\upsilon,
\end{split}
\end{equation}
where each $C^{j,i_1,\dots ,i_{\mu-1}}_{g}$ is simply subsequent
to $pre\vec{\kappa}^{+}_{simp}$ (the above holds modulo terms of
length $\ge\sigma+u+2$).

\par Now, two observations: Firstly, in the generic notation we
have introduced, we have:

\begin{equation}
\begin{split} 
&\Sum_{h\in H} a_h C^{h,i_1\dots ,i_\mu}_{g}(\Omega_1,\dots,
\Omega_p,\phi_{u+1},\phi_1,\dots ,\phi_u)
\\&=\Sum_{d\in D} a_d
C^{d,i_1\dots ,i_\mu}_{g}(\Omega_1,\dots,
\Omega_p,\phi_{u+1},\phi_1,\dots ,\phi_u).
\end{split}
\end{equation}
 Secondly, since
(\ref{pesax}) holds formally,  by making the
$\nabla\upsilon$'s into $Xdiv$'s (using the last Lemma 
in the Appendix of \cite{alexakis1}), we derive:

\begin{equation}
\label{pesax2}
\begin{split}
&\Sum_{d\in D_{nc}} a_d Xdiv_{i_1}\dots
Xdiv_{i_{\mu-1}}C^{d,i_1,\dots ,i_{\mu-1}}_{g} (\Omega_1,\dots
,\Omega_p,\phi_{u+1},\phi_1,\dots ,\phi_u) \\&-Xdiv_{i_1}\dots
Xdiv_{i_\mu}\Sum_{h\in H} a_h C^{h,i_1\dots
,i_\mu}_{g}(\Omega_1,\dots, \Omega_p,\phi_{u+1},\phi_1,\dots
,\phi_u)=
\\&\Sum_{d\in D^\sharp} a_d C^d_{g}
(\Omega_1,\dots ,\Omega_p,\phi_{u+1},\phi_1,\dots,\phi_u)+
\\&\Sum_{p\in P'} a_p Xdiv_{i_1}\dots
Xdiv_{i_{\mu-1}} C^{p,i_1,\dots ,i_{\mu-1},i_{*}}_{g}
(\Omega_1,\dots ,\Omega_p,\phi_1,\dots
,\phi_u)\nabla_{i_{*}}\phi_{u+1}+
\\&\Sum_{j\in J} a_j C^{j}_{g}
(\Omega_1,\dots ,\Omega_p,\phi_{u+1},\phi_1,\dots ,\phi_u)+
\Sum_{z\in Z} a_z C^z_{g} (\Omega_1,\dots ,\Omega_p,\phi_1,\dots
,\phi_{u+1}),
\end{split}
\end{equation}
(modulo length $\ge\sigma+u+2$). Thus, by virtue of the two above
equations we are reduced to showing our claim under the extra
assumption that $D_{nc}=\emptyset$.
\newline

 Next, we refer back to (\ref{tonkyrio=0}) and we claim
that we can write:

\begin{equation}
\label{phlon}
\begin{split}
&\Sum_{d\in D} a_d Xdiv_{i_1}\dots Xdiv_{i_a} C^{d,i_1,\dots
,i_a}_{g} (\Omega_1,\dots ,\Omega_p,\phi_{u+1},\phi_1,\dots
,\phi_u)=
\\& \Sum_{d\in D'} a_d Xdiv_{i_1}\dots
Xdiv_{i_a} C^{d,i_1,\dots ,i_a}_{g} (\Omega_1,\dots
,\Omega_p,\phi_{u+1},\phi_1,\dots ,\phi_u)+
\\&\Sum_{j\in J} a_j
C^j_{g} (\Omega_1,\dots ,\Omega_p,\phi_{u+1},\phi_1,\dots
,\phi_u)+\Sum_{z\in Z} a_z C^z_{g}
(\Omega_1,\dots ,\Omega_p,\phi_1,\dots ,\phi_{u+1}),
\end{split}
\end{equation}
modulo length $\ge\sigma+u+2$.

\par Here $$\Sum_{d\in D'} a_d C^{d,i_1,\dots ,i_a}_{g} (\Omega_1,\dots
,\Omega_p,\phi_{u+1},\phi_1,\dots ,\phi_u)$$ stands for a generic
linear combination of tensor fields in the general form
(\ref{plusou}) with a factor $\nabla^{(m)}\phi_{u+1}$ that contracts according to the pattern
\\$\nabla^{(3)}_{sr_1z}\phi_{u+1}\nabla^{r_1}\phi_{Min}$, where neither
of the indices ${}_s,{}_z$ is free.

\par The above equation can be proven as follows: Firstly, apply the
 eraser (in the equation (\ref{tonkyrio=0}))
 to the factor $\nabla\phi_{Min}$ that is contracting
against the factor $\nabla^{(2)}\phi_{u+1}$ (and thus obtain a new
true equation which we denote by (\ref{tonkyrio=0})'), and then
pick out the sublinear
 combination that contains a factor $\nabla\phi_{u+1}$
 with only one derivative (and thus obtain another
true equation which we denote by (\ref{tonkyrio=0})'')--for the
next construction we re-name the function $\phi_{u+1}$
$Y$.\footnote{(This is after we have applied the eraser to the
factor $\nabla\phi_{Min}$ that contracted against
$\nabla^{(B)}\phi_{u+1}$).} We are then in a 
position to apply Corollary 2 from \cite{alexakis4} (if $\sigma>3$)
 or Lemma 4.7 from \cite{alexakis4} (if $\sigma=3$)
  to the new true equation (\ref{tonkyrio=0})'',
and derive (\ref{phlon}): We start with the case $\sigma=3$:
Observe that (\ref{tonkyrio=0})'' satisfies the requirements of
Lemma 4.7 by weight considerations since we are assuming that
the assumption of Lemma \ref{zetajones} does 
not contain ``forbidden'' tensor fields. Thus, we apply Lemma
4.7 in \cite{alexakis4} and in the end we replace the function $Y$ by a
function $\nabla_s\phi_{u+1}\nabla^s\phi_{Min}$. Then, picking out
the sublinear combination with the function $\nabla\phi_{Min}$
differentiated only once,\footnote{Observe that this sublinear
combination will vanish separately.} we derive (\ref{phlon}).
Now, the case $\sigma>3$ follows by the exact same argument, only
instead of Lemma 4.7 we apply Corollary 2 from \cite{alexakis4}.
 Corollary 2 can be  applied to
(\ref{tonkyrio=0})'' since by definition there is no 
danger of falling under a ``forbidden case'' of that Lemma.

\par But then, we refer to (\ref{phlon}) and we
 observe that we can write:

\begin{equation}
\label{matchup} \begin{split} & \Sum_{d\in D'} a_d Xdiv_{i_1}\dots
Xdiv_{i_a} C^{d,i_1,\dots ,i_a}_{g} (\Omega_1,\dots
,\Omega_p,\phi_{u+1},\phi_1,\dots ,\phi_u)=
\\&\Sum_{d\in D^\sharp} a_d C^{d}_{g} (\Omega_1,\dots
,\Omega_p,\phi_{u+1},\phi_1,\dots ,\phi_u)+ \\&\Sum_{p\in P} a_p
Xdiv_{i_1}\dots Xdiv_{i_a} C^{d,i_1,\dots i_a,i_{*}}_{g}
(\Omega_1,\dots ,\Omega_p,\phi_1,\dots ,\phi_u)\nabla_{i_{*}}
\phi_{u+1}+
\\&\Sum_{z\in Z} a_z C^z_{g}
(\Omega_1,\dots ,\Omega_p,\phi_1,\dots ,\phi_{u+1}).
\end{split}
\end{equation}

This just follows from the identity:

\begin{equation}
\label{tory}
\begin{split}
(\nabla_a\nabla^{(2)}_{r_1 r_2}\phi_{u+1})\nabla^{r_1}\phi_{Min}=
\nabla^{(3)}_{r_1a r_2}\phi_{u+1}\nabla^{r_1}\phi_{Min}+
R_{ar_{1}kr_2}\nabla^k\phi_{u+1}.
\end{split}
\end{equation}

\par Thus, replacing (\ref{phlon}) and (\ref{matchup}) into 
(\ref{tonkyrio=0})  we are now reduced to showing Lemma \ref{bosnios}
 in the case $D_{nc}=\emptyset, D=\emptyset$.

\par Now, one more Lemma:

\begin{lemma}
\label{karnagio} Assume (\ref{tonkyrio=0}) with $D_{nc}=D=\emptyset$. We claim:

\begin{equation}
\label{tail}
\begin{split}
&\Sum_{d\in D^\sharp} a_d C^d_{g}(\Omega_1,\dots,
\Omega_p,\phi_{u+1},\phi_1,\dots ,\phi_u)= \Sum_{j\in J} a_j
C^j_{g}(\Omega_1,\dots, \Omega_p,\phi_1,\dots ,\phi_u,\phi_{u+1})
\\&+\Sum_{z\in Z} a_z C^z_{g}
(\Omega_1,\dots ,\Omega_p,\phi_1,\dots ,\phi_{u+1}),
\end{split}
\end{equation}
(modulo length $\ge\sigma+u+2$), where here each $C^j$ on the
right hand side has length $\sigma +u+1$ and a weak character
$Weak(\vec{\kappa}^{+}_{simp})$ but also has the factor
$\nabla\phi_{Min}$ contracting against a derivative index of the
factor $\nabla^{(m)}R_{ijkl}$ (and thus is simply  subsequent to
$\vec{\kappa}^{+}_{simp}$).

\par Moreover, we claim that we can write:
\begin{equation}
\label{tail2}
\begin{split}
&\Sum_{j\in J} a_j C^j_{g}(\Omega_1,\dots,
\Omega_p,\phi_{u+1},\phi_1,\dots ,\phi_u) =
\\&\Sum_{j\in J} a_j C^j_{g}(\Omega_1,\dots,
\Omega_p,\phi_1,\dots ,\phi_u,\phi_{u+1}) +\Sum_{z\in Z} a_z
C^z_{g} (\Omega_1,\dots ,\Omega_p,\phi_1,\dots ,\phi_{u+1}),
\end{split}
\end{equation}
where each $C^j_g$ on the right hand side is simply subsequent to
$\vec{\kappa}^{+}_{simp}$.
\end{lemma}

{\it Proof:} In order to introduce a strict dichotomy between the
linear combinations indexed in $D^\sharp, J$, we permute the
indices in the factor $\nabla^{(B)}\phi_{u+1}$ in each
$C^j_{g}(\Omega_1,\dots, \Omega_p,\phi_{u+1},\phi_1,\dots
,\phi_u)$ to make the factors
$\nabla\tilde{\phi}_{Min},\nabla\phi'_{e_1},\dots
\nabla\phi'_{e_{\pi-1}}$ contract against the first $\pi$ indices,
in that order. We can clearly do this modulo introducing
correction terms in the general form:

$$\Sum_{j\in J} a_j C^j_{g}(\Omega_1,\dots,
\Omega_p,\phi_1,\dots ,\phi_u,\phi_{u+1}) +\Sum_{z\in Z} a_z
C^z_{g} (\Omega_1,\dots ,\Omega_p,\phi_1,\dots ,\phi_{u+1}).$$
Then, we inquire on the number of derivatives on the factor
$\nabla^{(B)}\phi_{u+1}$ in each $C^j$. If $B>\pi+1$ we just
re-name $C^j$ into $C^d$ and index it in $D^\sharp$. We are
reduced to showing our claim under the hypothesis that all $C^j$
have $B\le \pi+1$ derivatives on the factor
$\nabla^{(B)}\phi_{u+1}$. We proceed under that assumption.

 Trivially, since (\ref{tonkyrio=0}) holds formally and since
 we are assuming $D_{nc}=D=\emptyset$, we derive that:

$$\Sum_{d\in D^\sharp} a_d linC^d_{g}(\Omega_1,\dots,
\Omega_p,\phi_{u+1},\phi_1,\dots ,\phi_u)=0,$$

$$\Sum_{j\in J} a_j linC^j_{g}(\Omega_1,\dots,
\Omega_p,\phi_{u+1},\phi_1,\dots ,\phi_u)=0.$$

 Therefore our claim
(\ref{tail2}) follows by just repeating the permutations by which
we make the left hand sides of the last two equations formally
zero, whereas (\ref{tail}) follows by the same fact, and also by
using the fact that the first $\pi$ indices in the factor
$\nabla^{(B)}\phi_{u+1}$ are not permuted (which can be proven as
usual using the eraser). $\Box$
\newline

\par Thus, in the setting of Lemma \ref{zetajones}
and of Lemma\ref{pskovb} (if the selected factor is in the
 form $S_{*}\nabla^{(\nu)}R_{ijkl}$) we have shown that:

\begin{equation}
\label{keepinminda}
\begin{split}
&\Sum_{l\in L} a_l CurvTrans[Xdiv_{i_1}\dots
Xdiv_{i_a}C^{l,i_1\dots i_a}_{g}(\Omega_1\dots
,\Omega_p,\phi_1,\dots ,\phi_u)]+
\\& \Sum_{j\in J} a_j CurvTrans[C^j_{g}
(\Omega_1\dots ,\Omega_p,\phi_1,\dots ,\phi_u)]=
\\&\Sum_{l\in L} a_l Leftover_{\phi_{u+1}}
[Xdiv_{i_1}\dots Xdiv_{i_a}C^{l,i_1\dots i_a}_{g} (\Omega_1\dots
,\Omega_p,\phi_1,\dots ,\phi_u)]+
\\&(\Sum_{p\in P'} a_p Xdiv_{i_1}\dots Xdiv_{i_{\mu-1}}
C^{p,i_1,\dots i_{\mu-1},i_{*}}_{g} (\Omega_1,\dots
,\Omega_p,\phi_1,\dots ,\phi_u)\nabla_{i_{*}} \phi_{u+1})+
\\&\Sum_{p\in P} a_p Xdiv_{i_1}\dots Xdiv_{i_a}
C^{p,i_1,\dots i_a,i_{*}}_{g} (\Omega_1,\dots
,\Omega_p,\phi_1,\dots ,\phi_u)\nabla_{i_{*}} \phi_{u+1}+
\\& \Sum_{j\in J} a_j C^j_{g}
(\Omega_1,\dots ,\Omega_p,\phi_1,\dots
,\phi_u,\phi_{u+1})+\Sum_{z\in Z} a_z C^z_{g} (\Omega_1,\dots
,\Omega_p,\phi_1,\dots ,\phi_{u+1});
\end{split}
\end{equation}
(the linear combination $\sum_{p\in P'}\dots$ arises only in the
setting of Lemma \ref{zetajones}).

\section{A study of the sublinear combinations $LC[L_{g}]$ and $W[L_{g}]$
in (\ref{heidegger}). Computations and cancellations.}
\label{easyproof}

\subsection{General discussion of ideas:}

 The main conclusions we retain from the previous two
subsections are equations (\ref{keepinminda}) and
(\ref{akinola'}), (\ref{akinola2'}), (\ref{akinola3'}).

\par We will denote by $CurvTrans^{study}[L_g]+\sum_{z\in Z}\dots $ the right hand sides
of those equations.\footnote{$\sum_{z\in Z}\dots$ stands for
the sublinear combination of complete contractions indexed in $Z$
in the right hand sides of (\ref{keepinminda}) and
(\ref{akinola'}), (\ref{akinola2'}), (\ref{akinola3'}),
 and $CurvTrans^{study}[L_g]$ stands for the terms
 with a factor $\nabla\phi_{u+1}$ in the RHSs in
(\ref{keepinminda}) and (\ref{akinola'}), (\ref{akinola2'}),
(\ref{akinola3'}).} Thus, we have shown that the sublinear
combinations $CurvTrans[L_g]$ can be {\it re-expressed} as linear
combinations $CurvTrans^{study}[L_g]+\sum_{z\in Z}\dots$.

\par We then {\it replace} $CurvTrans[L_g]$ by $CurvTrans^{study}[L_g]+\sum_{z\in Z}\dots$
in (\ref{heidegger}), obtaining a new equation:

\begin{equation}
\label{prheideger2}
CurvTrans^{study}[L_g]+LC[L_g]+W[L_g]+\sum_{z\in Z}\dots=0,
\end{equation}
which holds modulo complete contractions of length $\ge\sigma+u+2$ (notice all
 complete contractions in the LHS of the above have length $\sigma+u+1$).
 Thus, picking out the sublinear combination of complete contractions
 with length $\sigma+u+1$ and with a factor $\nabla\phi_{u+1}$ (with only one
derivative)\footnote{This sublinear combination must vanish separately.} we obtain a new equation:

\begin{equation}
\label{heidegger2} CurvTrans^{study}[L_g]+LC[L_g]+W[L_g]=0,
\end{equation}
which holds modulo complete contractions of length $\ge\sigma+u+2$.

\par For the rest of this paper and the next one in this series, 
we will try to understand the sublinear combinations
$LC_{\phi_{u+1}}[L_g]+W_{\phi_{u+1}}[L_g]$ in the above equation.
\newline

\par We now focus on the sublinear combinations
$LC[L_{g}]$ and $W[L_{g}]$ in (\ref{heidegger2}) which by
definition consist of terms with length $\sigma +u+1$ {\it and} have a weak
$(u+1)$-character $Weak(\vec{\kappa}^{+}_{simp})$.  We recall
that $LC[L_{g}]$ and $W[L_{g}]$ were defined to be specific
sublinear combinations of $Image^{1,+}_{\phi_{u+1}}
 [L_{g}]$: $LC[L_g]$ stands for the sublinear combination that
arises {\it in
 $Image^{1,+}_{\phi_{u+1}} [L_{g}]$} by applying the formula
 (\ref{levicivita}), and $W[L_g]$ stands for the sublinear
 combination of terms that arises by virtue of the transformation
 $R_{ijkl}(e^{2\phi_{u+1}}g)\rightarrow
 e^{2\phi_{u+1}}R_{ijkl}(g)$.

\par Furthermore, recall that we have broken up
$$LC[Xdiv_{i_1}\dots Xdiv_{i_a}
C^{l,i_1\dots i_a}_{g}(\Omega_1,\dots
,\Omega_p,\phi_{u+1},\phi_1,\dots ,\phi_{u+1})]$$ into two
sublinear combinations:
$$LC_\Phi[Xdiv_{i_1}\dots Xdiv_{i_a}
C^{l,i_1\dots i_a}_{g}(\Omega_1,\dots
,\Omega_p,\phi_{u+1},\phi_1,\dots ,\phi_{u+1})],$$

$$LC^{No\Phi}[Xdiv_{i_1}\dots Xdiv_{i_a}
C^{l,i_1\dots i_a}_{g}(\Omega_1,\dots
,\Omega_p,\phi_{u+1},\phi_1,\dots ,\phi_{u+1})].$$ Recall that
$LC_\Phi[Xdiv_{i_1}\dots Xdiv_{i_a}
C^{l,i_1\dots i_a}_{g}(\Omega_1,\dots
,\Omega_p,\phi_{u+1},\phi_1,\dots ,\phi_{u+1})]$ is the sublinear
combination in $LC[Xdiv_{i_1}\dots
Xdiv_{i_a}C^{l,i_1\dots i_a}_{g}(\Omega_1,\dots
,\Omega_p,\phi_1,\dots ,\phi_u)]$ that arises by applying
 (\ref{levicivita}) to
 two indices, at least one of which is contracting against
a factor $\nabla\phi_y$. Recall (\ref{trikala}).

Recall that $LC^{No\Phi}[Xdiv_{i_1}\dots Xdiv_{i_a}
C^{l,i_1\dots i_a}_{g}(\Omega_1,\dots ,\Omega_p,\phi_1,\dots
,\phi_u)]$ stands for the sublinear combination that arises in
\\ $LC[Xdiv_{i_1}\dots
Xdiv_{i_a} C^{l,i_1\dots i_a}_{g}(\Omega_1,\dots
,\Omega_p,\phi_1,\dots ,\phi_u)]$ by
 applying (\ref{levicivita}) to two indices that are not
 contracting against a factor $\nabla\phi_y$.
 We have analogously defined
$LC^{No\Phi}[ C^j_{g}(\Omega_1,\dots ,\Omega_p,\phi_1,\dots
,\phi_u)]$.

\par Our aim for this section is to understand the sublinear
combinations:

\begin{equation}
\label{asdf1} \begin{split} & LC^{No\Phi}[Xdiv_{i_1}\dots
Xdiv_{i_a} C^{l,i_1\dots i_a}_{g}(\Omega_1,\dots
,\Omega_p,\phi_1,\dots ,\phi_{u})]+\\&Xdiv_{i_1}\dots Xdiv_{i_a}
C^{l,i_1\dots i_a}_{g}(\Omega_1\cdot\phi_{u+1},\dots
,\Omega_p,\phi_1,\dots ,\phi_{u})+\dots +
\\&Xdiv_{i_1}\dots
Xdiv_{i_a} C^{l,i_1\dots
i_a}_{g}(\Omega_1,\dots,\Omega_X\cdot\phi_{u+1},\dots
,\Omega_p,\phi_1,\dots ,\phi_{u}),
\end{split}
\end{equation}

\begin{equation}
\label{asdf2}
\begin{split}
&LC^{No\Phi}[ C^j_{g}(\Omega_1,\dots ,\Omega_p,\phi_1,\dots
,\phi_{u})]+
C^j_{g}(\Omega_1\cdot\phi_{u+1},\dots ,\Omega_p,\phi_1,\dots
,\phi_{u})\\&+\dots +C^j_{g}(\Omega_1,\dots ,\Omega_X\cdot
\phi_{u+1},\dots ,\Omega_p,\phi_1,\dots ,\phi_{u}),
\end{split}
\end{equation}

\begin{equation}
\label{asdf3} W[Xdiv_{i_1}\dots Xdiv_{i_a} C^{l,i_1\dots
i_a}_{g}(\Omega_1,\dots ,\Omega_p,\phi_1,\dots ,\phi_{u})],
\end{equation}

\begin{equation}
\label{asdf4}W[C^j_{g}(\Omega_1,\dots ,\Omega_p,\phi_1,\dots
,\phi_{u})].\end{equation} (Recall that $\Omega_1,\dots ,\Omega_X$
are the factors $\Omega_h$ that are not contracting against a
factor $\nabla\phi_h$ in $\vec{\kappa}_{simp}$).

\par A few notes: In the setting of Lemmas \ref{zetajones} and
\ref{pool2}, we will be able to {\it discard} sublinear
combinations in the above four linear combinations that consist of
complete contractions with length $\sigma+u+1$ and where the
factor $\nabla\phi_{u+1}$ is contracting against a derivative
index in the crucial factor when it is in the form
$\nabla^{(m)}R_{ijkl}$,\footnote{i.e. in the setting of Lemma
\ref{pool2}} or contracting against any index ${}_{r_1},\dots
,{}_{r_\nu},{}_j$ in the crucial factor if it is in the form
$S_{*}\nabla^{(\nu)} R_{ijkl}$.\footnote{i.e. if we are in the
setting of Lemma \ref{zetajones}.} We denote generic such
linear combinations by $\Sum_{q\in Q} a_q C^q(\Omega_1,\dots
,\Omega_p,\phi_1,\dots ,\phi_{u+1})$. The contractions $C^q_{g}$
do not have to be acceptable. In particular, they
 might have a factor $\nabla\Omega_h$. Observe that such generic
 linear combinations are {\it allowed} in the right hand sides of
 Lemmas \ref{zetajones} and \ref{pool2}: They are special cases of
 the linear combinations $\sum_{t\in T} \dots$ in the right hand
 sides of those Lemmas.
\newline

\par We introduce another notational convention we will be using
 throughout this section:

\begin{definition}
\label{contributeur} We  denote by $\sum_{h\in H} a_h
C^{i_1\dots i_a,i_{*}}_g(\Omega_1,\dots ,\Omega_p,\phi_1,\dots
,\phi_u)\nabla_{i_{*}}\phi_{u+1}$ a  generic linear combination of
the forms $\sum_{p\in P}a_p\dots$ as in the statements of Lemmas
\ref{zetajones} and \ref{pool2} (if we are in the setting of those
Lemmas), or a generic linear combination of the form $\sum_{t\in
T_1\bigcup T_2\bigcup T_3\bigcup T_4}\dots$, in the notation of
Lemma \ref{pskovb} (if we are in the setting of case A of that
Lemma).

We will be calling the tensor fields in those linear combinations
``contributors''.
\end{definition}

\par Finally, we also recall that $\sum_{j\in J} a_j C^j_g(\Omega_1,\dots,\Omega_p,
\phi_1,\dots ,\phi_{u+1})$ stands for a generic linear combination
of complete contractions of length $\sigma+u+1$, with a weak
$(u+1)$-character $Weak(\vec{\kappa}^{+}_{simp})$ and which are
 $u$-simply subsequent to $\vec{\kappa}_{simp})$.

\par Now, we proceed to study the four expressions (\ref{asdf1}),
(\ref{asdf2}), (\ref{asdf3}), (\ref{asdf4}).

\par The easiest to study are (\ref{asdf2}) and (\ref{asdf4}).
We straightforwardly observe that for each $j\in J$ we must have:

\begin{equation}
\label{dwyer}
\begin{split}
&LC^{No\Phi}[ C^j_{g}(\Omega_1,\dots ,\Omega_p,\phi_1,\dots
,\phi_{u})]+
C^j_{g}(\Omega_1\cdot\phi_{u+1},\dots ,\Omega_p,\phi_1,\dots
,\phi_{u})\\&+\dots +C^j_{g}(\Omega_1,\dots ,\Omega_X\cdot
\phi_{u+1},\dots ,\Omega_p,\phi_1,\dots ,\phi_{u})  \\&
=\Sum_{j\in J} a_j C^j_{g}(\Omega_1,\dots ,\Omega_p,\phi_1,\dots
,\phi_{u+1}),
\end{split}
\end{equation}

\begin{equation}
\label{dwyer'}
\begin{split}
&W[ C^j_{g}(\Omega_1,\dots ,\Omega_p,\phi_1,\dots ,\phi_{u})]=
 \Sum_{j\in J} a_j C^j_{g}(\Omega_1,\dots
,\Omega_p,\phi_1,\dots ,\phi_{u+1}).
\end{split}
\end{equation}

The harder challenge is to understand the linear
combinations (\ref{asdf1}), (\ref{asdf3}). In order to understand these
 two sublinear combinations, we will break them up into further sublinear combinations:

\begin{definition}
\label{polypoly}
 We define $LC^{No\Phi}_{\phi_{u+1}}[C^{l,i_1\dots
i_a}_{g}(\Omega_1,\dots ,\Omega_p,\phi_1,\dots ,\phi_u)]$ to stand
for the sublinear combination that arises in

$$Image^1_{\phi_{u+1}}[C^{l,i_1\dots i_a}_{g}(\Omega_1,\dots
,\Omega_p,\phi_1,\dots ,\phi_u)]$$ when we apply the
transformation law (\ref{levicivita}) to any two indices in the
tensor field $C^{l,i_1\dots i_a}_{g}$ (we will call these original
indices), that are not contracting against a factor $\nabla\phi_y,
y\le u$ and bring out a factor $\nabla_{i_{*}}\phi_{u+1}$ for
which ${}_{i_{*}}$ is either contracting against the crucial
factor or is a free index.

\par We also define

$$LC^{No\Phi, div}_{\phi_{u+1}}[Xdiv_{i_1}\dots Xdiv_{i_a}
C^{l,i_1\dots i_a}_{g}(\Omega_1,\dots ,\Omega_p,\phi_1,\dots
,\phi_u)]$$ to stand for the sublinear combination that arises in
$$Image^1_{\phi_{u+1}}[Xdiv_{i_1}\dots Xdiv_{i_a}
C^{l,i_1\dots i_a}_{g}(\Omega_1,\dots ,\Omega_p,\phi_1,\dots
,\phi_u)]$$ when we apply (\ref{levicivita}) to two indices in the
same factor, and at least one of those indices is of the form
$\nabla^{i_h}$ ($1\le h\le a$) (we call such indices divergence
indices), and the other is not contracting against a factor
$\nabla\phi_y$.

\par Secondly, we denote by $W[C^{l,i_1\dots i_a}_{g}(\Omega_1,\dots
,\Omega_p,\phi_1,\dots ,\phi_u)]$ the sublinear combination in

\begin{equation}
\label{strata} Image^1_{\phi_{u+1}}[C^{l,i_1\dots
i_a}_{g}(\Omega_1,\dots ,\Omega_p,\phi_1,\dots ,\phi_u)]
\end{equation}
 that arises when we replace a factor $\nabla^{(m)}_{r_1\dots r_m}R_{ijkl}$ by a factor
 $\nabla^{(m)}_{r_1\dots r_m}[e^{2\phi_{u+1}}R_{ijkl}]$ (by virtue of (\ref{curvtrans})) and
 then bring out an expression $e^{2\phi_{u+1}}\nabla^{(m-1)}R_{ijkl}\nabla\phi_{u+1}$
  by hitting the factor $e^{2\phi_{u+1}}$ by one of the derivatives
  $\nabla_{r_1},\dots ,\nabla_{r_m}$.

\par Furthermore, we define $W^{div}[Xdiv_{i_1}\dots Xdiv_{i_a}C^{l,i_1\dots
i_a}_{g}(\Omega_1,\dots ,\Omega_p,\phi_1,\dots ,\phi_u)]$ to stand
for the sublinear combination in \\$W[Xdiv_{i_1}\dots
Xdiv_{i_a}C^{l,i_1\dots i_a}_{g} (\Omega_1,\dots
,\Omega_p,\phi_1,\dots ,\phi_u)]$ that arises by picking a factor
$\nabla^{i_y\dots i_z}\nabla^{(m)}_{r_1\dots r_m}R_{ijkl}$ in some
summand in \\$Xdiv_{i_1}\dots Xdiv_{i_a}C^{l,i_1\dots
i_a}_{g}(\Omega_1,\dots ,\Omega_p,\phi_1,\dots ,\phi_u)$
(${}^{i_y},\dots ,{}^{i_z}$ are divergence indices) and replacing
it by $\nabla^{i_y\dots i_z}[e^{2\phi_{u+1}}\nabla^{(m)}_{r_1\dots
r_m}R_{ijkl}]$ and then bringing out an expression
$e^{2\phi_{u+1}}\nabla^{(m'-1)}R_{ijkl}\nabla\phi_{u+1}$ by
hitting  the factor $e^{2\phi_{u+1}}$ by one of the
divergence indices ${\nabla}^{i_y},\dots ,\nabla^{i_z}$.
\end{definition}

\par It then follows directly from the above definition that:

\begin{equation}
\label{ela}
\begin{split}
&LC^{No\Phi}_{\phi_{u+1}}[Xdiv_{i_1}\dots Xdiv_{i_a} C^{l,i_1\dots
i_a}_{g}(\Omega_1,\dots ,\Omega_p,\phi_1,\dots ,\phi_u)]
\\& +W[Xdiv_{i_1}\dots Xdiv_{i_a} C^{l,i_1\dots
i_a}_{g}(\Omega_1,\dots ,\Omega_p,\phi_1,\dots ,\phi_u)]=
\\&LC^{No\Phi, div}_{\phi_{u+1}}[Xdiv_{i_1}\dots Xdiv_{i_a}
C^{l,i_1\dots i_a}_{g}(\Omega_1,\dots ,\Omega_p,\phi_1,\dots
,\phi_u)]+
\\&Xdiv_{i_1}\dots Xdiv_{i_a}LC^{No\Phi}_{\phi_{u+1}}
[C^{l,i_1\dots i_a}_{g}(\Omega_1,\dots ,\Omega_p,\phi_1,\dots
,\phi_u)]+
\\&W^{div}[Xdiv_{i_1}\dots Xdiv_{i_a}C^{l,i_1\dots
i_a}_{g}(\Omega_1,\dots ,\Omega_p,\phi_1,\dots ,\phi_u)]+
\\& Xdiv_{i_1}\dots Xdiv_{i_a}W[C^{l,i_1\dots
i_a}_{g}(\Omega_1,\dots ,\Omega_p,\phi_1,\dots ,\phi_u)],
\end{split}
\end{equation}
subject to a small clarification regarding the notion
of $Xdiv$ in the linear combinations $Xdiv_{i_1}\dots
Xdiv_{i_a}LC^{No\Phi}_{\phi_{u+1}} [C^{l,i_1\dots
i_a}_{g}(\Omega_1,\dots ,\Omega_p,\phi_1,\dots ,\phi_u)]$:

\begin{definition}
\label{inyourheart}
We have defined $Xdiv_{i_y}$ to stand for the
sublinear combination in $div_{i_y}$ where $\nabla_{i_y}$ is not
allowed to hit the factor to which ${}_{i_y}$ belongs, nor any
 factor $\nabla\phi_h$.

\par Now, for each free index ${}_{i_y}$ (that belongs to a factor $T$
 in the form $\nabla^{(p)}\Omega_h$ or $\nabla^{(m)}R_{ijkl}$),
 and each tensor field $C^{*,i_1\dots i_a}_{g}
(\Omega_1,\dots ,\Omega_p,\phi_1,\dots ,\phi_{u+1})$ in
\\ $LC^{No\Phi}_{\phi_{u+1}} [C^{l,i_1\dots
i_a}_{g}(\Omega_1,\dots ,\Omega_p,\phi_1,\dots ,\phi_u)]$, 
either ${}_{i_y}$ still belongs to the factor
$\nabla^{(p)}\Omega_h$ or $\nabla^{(m)}R_{ijkl}$ in $C^{*,i_1\dots
i_a}_{g}(\Omega_1,\dots ,\Omega_p,\phi_1, \dots ,\phi_{u+1})$, or
${}_{i_y}$ belongs to an un-contracted metric tensor $g$, or
${}_{i_y}$ belongs to a factor $\nabla_{i_y}\phi_{u+1}$. In
the first case, we define $Xdiv_{i_y}C^{*,i_1\dots
i_a}_{g}(\Omega_1,\dots ,\Omega_p,\phi_1,\dots ,\phi_{u+1})$ as in
the above paragraph. In the second case, we see that the
un-contracted metric tensor in \\$C^{*,i_1\dots i_a}_{g}
(\Omega_1,\dots ,\Omega_p,\phi_1,\dots ,\phi_{u+1})$ must have
arisen by applying the third summand on the right hand side on
(\ref{levicivita}) to a pair of indices $(\nabla_a, {}_b)$
 in a factor $T$ of the form $\nabla^{(p)}\Omega_h$ or
$\nabla^{(m)}R_{ijkl}$, where either $\nabla_a$ or ${}_b$ is the
free index ${}_{i_y}$. In that case, we define
\\$Xdiv_{i_y}C^{*,i_1\dots i_a}_{g}(\Omega_1,\dots
,\Omega_p,\phi_1,\dots ,\phi_{u+1})$ to stand for the sublinear
combination in $div_{i_y}C^{*,i_1\dots i_a}_{g}(\Omega_1,\dots
,\Omega_p,\phi_1,\dots ,\phi_{u+1})$ where $\nabla_{i_y}$ is not
allowed to hit the factor $T$ nor the uncontracted metric tensor
nor the factor $\nabla_{i_y}\phi_{u+1}$ nor any $\nabla\phi_h,
h\le u$. In the third case, the expression $\nabla_{i_y}\phi$ must
have arisen by applying one of the first two terms in the
transformation law (\ref{levicivita}) to the factor T. In that
case, we define $Xdiv_{i_y}$ to stand for the sublinear
combination in $div_{i_y}$ where $\nabla_{i_y}$ is not allowed to
hit the factor $T$ nor the factor $\nabla_{i_y}\phi_{u+1}$ nor any
$\nabla\phi_h, h\le u$.
\end{definition}

\par With this clarification equation (\ref{ela}) just follows by
the definition of \\$LC^{No\Phi}_{\phi_{u+1}}[\dots]$ and $W[\dots]$
and the transformation law (\ref{levicivita}). Now, we will
 subdivide the right hand side of (\ref{ela}) into further
sublinear combinations:
\newline

{\bf A study of the right hand side of (\ref{ela}):}
\newline

\par We will firstly study the last two lines in (\ref{ela}).

We define $$W^{targ}[C^{l,i_1\dots i_a}_{g}(\Omega_1,\dots
,\Omega_p,\phi_1,\dots ,\phi_u)]$$ to stand for the sublinear
combination in
$$W[C^{l,i_1\dots i_a}_{g}(\Omega_1,\dots
,\Omega_p,\phi_1,\dots ,\phi_u)]$$ for which
$\nabla\phi_{u+1}$ is contracting against the crucial factor.
\newline

 We define
$$W^{free}[C^{l,i_1\dots
i_a}_{g}(\Omega_1,\dots ,\Omega_p,\phi_1,\dots ,\phi_u)]$$ to
stand for the sublinear combination in $W[C^{l,i_1\dots
i_a}_{g}(\Omega_1,\dots ,\Omega_p,\phi_1,\dots ,\phi_u)]$ for
which the index ${}_\alpha$ in $\nabla_\alpha\phi_{u+1}$ is a free index.

\par We also define
$$W^{targ,div}[Xdiv_{i_1}\dots Xdiv_{i_a}
C^{l,i_1\dots i_a}_{g}(\Omega_1,\dots ,\Omega_p,\phi_1,\dots
,\phi_u)]$$ to stand for the sublinear combination in
$W^{div}[\dots]$ where in addition the factor $\nabla\phi_{u+1}$
that we bring out is contracting against the (a) crucial factor.

\par Thus it follows that:

\begin{equation}
\label{harun1}
\begin{split}
&W[\Sum_{l\in L} a_l Xdiv_{i_1}\dots Xdiv_{i_a} C^{l,i_1\dots
i_a}_{g}(\Omega_1,\dots ,\Omega_p,\phi_1,\dots ,\phi_u)]=
\\& W^{targ,div}
[\Sum_{l\in L} a_l Xdiv_{i_1}\dots Xdiv_{i_a} C^{l,i_1\dots
i_a}_{g}(\Omega_1,\dots ,\Omega_p,\phi_1,\dots ,\phi_u)]+
\\& \Sum_{l\in L} a_l Xdiv_{i_1}\dots Xdiv_{i_a}
W^{targ}[C^{l,i_1\dots i_a}_{g}(\Omega_1,\dots
,\Omega_p,\phi_1,\dots ,\phi_u)]+
\\& \Sum_{l\in L} a_l Xdiv_{i_1}\dots Xdiv_{i_a}
W^{free}[C^{l,i_1\dots i_a}_{g}(\Omega_1,\dots
,\Omega_p,\phi_1,\dots ,\phi_u)].
\end{split}
\end{equation}

\par In view of the above, we study the three
sublinear combinations in the right hand side of (\ref{harun1})
separately. We derive by definition:

\begin{equation}
\label{eudokia2}
\begin{split}
&Xdiv_{i_1}\dots Xdiv_{i_a}W^{free}[C^{l,i_1\dots
i_a}_{g}(\Omega_1,\dots ,\Omega_p,\phi_1,\dots ,\phi_u)]=
\\&\Sum_{q\in Q} a_q C^q_{g}(\Omega_1,\dots ,\Omega_p,\phi_1,\dots
,\phi_{u+1}).
\end{split}
\end{equation}

\par Next, we will study the remaining two sublinear combinations
in the right hand side of (\ref{harun1}) together:
 Recall that the total number of factors
$\nabla^{(m)}R_{ijkl}$ or $S_{*}\nabla^{(\nu)}R_{ijkl}$ is
$s=\sigma_1+\sigma_2$. Now,  in the setting of Lemma
\ref{zetajones}, if $l\in L_\mu$ and
$C^{l,i_1\dots i_a}_{g}$ has one free index (say ${}_{i_1}$ with no
loss of generality) being the index ${}_k$ in the crucial factor, we
have:

\begin{equation}
\label{evraz1}
\begin{split}
&Xdiv_{i_1}\dots Xdiv_{i_\mu}W^{targ}[C^{l,i_1\dots
i_\mu}_{g}(\Omega_1,\dots ,\Omega_p,\phi_1,\dots ,\phi_u)]+
\\&W^{targ,div}[Xdiv_{i_1}\dots Xdiv_{i_\mu}C^{l,i_1\dots
i_\mu}_{g}(\Omega_1,\dots ,\Omega_p,\phi_1,\dots ,\phi_u)]=
\\&2(s-1)Xdiv_{i_2}\dots Xdiv_{i_\mu}C^{l,i_1\dots
i_\mu}_{g}(\Omega_1,\dots ,\Omega_p,\phi_1,\dots
,\phi_u)\nabla_{i_1}\phi_{u+1}+
\\&\Sum_{h\in H} a_h Xdiv_{i_1}\dots Xdiv_{i_\mu}C^{h,i_1\dots
i_\mu,i_{*}}_{g}(\Omega_1,\dots ,\Omega_p,\phi_1,\dots
,\phi_u)\nabla_{i_{*}}\phi_{u+1}
\\&+\Sum_{q\in Q} a_q C^q_{g}(\Omega_1,\dots
,\Omega_p, \phi_1,\dots ,\phi_{u+1}),
\end{split}
\end{equation}
where
$$\Sum_{h\in H} a_h Xdiv_{i_1}\dots Xdiv_{i_\mu}C^{h,i_1\dots
i_\mu,i_{*}}_{g}(\Omega_1,\dots ,\Omega_p,\phi_1,\dots
,\phi_u)\nabla_{i_{*}}\phi_{u+1}$$ stands for a generic linear
combination of {\it acceptable} contributors.\footnote{See the
definition \ref{contributeur}.}
 If the tensor field
$C^{l,i_1\dots i_a}_{g}$ does not have a free index ${}_{i_h}$
that is an index ${}_k$ or ${}_l$ in the crucial factor, then we
calculate:

\begin{equation}
\label{evraz1b}
\begin{split}
&Xdiv_{i_1}\dots Xdiv_{i_\mu}W^{targ}[C^{l,i_1\dots
i_\mu}_{g}(\Omega_1,\dots ,\Omega_p,\phi_1,\dots ,\phi_u)]+
\\&W^{targ,div}[Xdiv_{i_1}\dots Xdiv_{i_\mu}C^{l,i_1\dots
i_\mu}_{g}(\Omega_1,\dots ,\Omega_p,\phi_1,\dots ,\phi_u)]=
\\&\Sum_{h\in H} a_h Xdiv_{i_1}\dots Xdiv_{i_\mu}C^{h,i_1\dots
i_\mu,i_{*}}_{g}(\Omega_1,\dots ,\Omega_p,\phi_1,\dots
,\phi_u)\nabla_{i_{*}}\phi_{u+1}
\\&+\Sum_{q\in Q} a_q C^q_{g}(\Omega_1,\dots
,\Omega_p, \phi_1,\dots ,\phi_{u+1}).
\end{split}
\end{equation}

If $l\in L\setminus L_\mu$, (hence $a>\mu$):

\begin{equation}
\label{evraz2}
\begin{split}
&Xdiv_{i_1}\dots Xdiv_{i_a}W^{targ}[C^{l,i_1\dots
i_a}_{g}(\Omega_1,\dots ,\Omega_p,\phi_1,\dots ,\phi_u)]+
\\&W^{targ,div}[Xdiv_{i_1}\dots Xdiv_{i_a}C^{l,i_1\dots
i_a}_{g}(\Omega_1,\dots ,\Omega_p,\phi_1,\dots ,\phi_u)]=
\\&\Sum_{h\in H} a_h Xdiv_{i_1}\dots Xdiv_{i_{a-1}}C^{h,i_1\dots
i_{a-1},i_{*}}_{g}(\Omega_1,\dots ,\Omega_p,\phi_1,\dots
,\phi_u)\nabla_{i_{*}}\phi_{u+1}
\\&+\Sum_{q\in Q} a_q C^q_{g}(\Omega_1,\dots
,\Omega_p, \phi_1,\dots ,\phi_{u+1}).
\end{split}
\end{equation}

\par On the other hand, in the setting of Lemma \ref{pool2},
recall that $I_{l,*}$ stands for the index set of special free
indices in factors $\nabla^{(m)}R_{ijkl}$ in $C^{l,i_1\dots
i_\mu}_g$; we then  compute:

\begin{equation}
\label{evraz3}
\begin{split}
&Xdiv_{i_1}\dots Xdiv_{i_\mu}W^{targ}[C^{l,i_1\dots
i_\mu}_{g}(\Omega_1,\dots ,\Omega_p,\phi_1,\dots ,\phi_u)]+
\\&W^{targ,div}[Xdiv_{i_1}\dots Xdiv_{i_\mu}C^{l,i_1\dots
i_\mu}_{g}(\Omega_1,\dots ,\Omega_p,\phi_1,\dots ,\phi_u)]=
\\&2(s-1)\sum_{i_h\in I_{*,l}}Xdiv_{i_1}\dots \hat{Xdiv}_{i_h}\dots Xdiv_{i_\mu}C^{l,i_1\dots
i_\mu}_{g}(\Omega_1,\dots ,\Omega_p,\phi_1,\dots
,\phi_u)\\&\nabla_{i_1}\phi_{u+1}
+\Sum_{h\in H} a_h Xdiv_{i_1}\dots Xdiv_{i_\mu}C^{h,i_1\dots
i_\mu,i_{*}}_{g}(\Omega_1,\dots ,\Omega_p,\phi_1,\dots
,\phi_u)\nabla_{i_{*}}\phi_{u+1}
\\&+\Sum_{z\in Z} a_z C^z_{g}(\Omega_1,\dots ,\Omega_p,
\phi_1,\dots ,\phi_{u+1})+\Sum_{q\in Q} a_q C^q_{g}(\Omega_1,\dots
,\Omega_p, \phi_1,\dots ,\phi_{u+1}).
\end{split}
\end{equation}

\par Moreover, in the setting of Lemma \ref{pool2} we again have the equation
(\ref{evraz2}).
\newline

{\it A study of  $LC^{No\Phi,
div}_{\phi_{u+1}}[Xdiv_{i_1}\dots Xdiv_{i_a} C^{l,i_1\dots
i_a}_{g}(\Omega_1,\dots ,\Omega_p,\phi_1,\dots ,\phi_u)]$,
\\$Xdiv_{i_1}\dots Xdiv_{i_a}LC^{No\Phi}_{\phi_{u+1}} [C^{l,i_1\dots
i_a}_{g}(\Omega_1,\dots ,\Omega_p,\phi_1,\dots ,\phi_u)]$ in the
RHS of (\ref{ela}):}
\newline

\par Next, we will study the first two lines in the right hand
side of (\ref{ela}).

We define
$$LC^{No\Phi,targ}_{\phi_{u+1}}
[C^{l,i_1\dots i_a}_{g}(\Omega_1,\dots ,\Omega_p,\phi_1,\dots
,\phi_u)]$$ to stand for the sublinear combination in
$$LC^{No\Phi}_{\phi_{u+1}}
[C^{l,i_1\dots i_a}_{g}(\Omega_1,\dots ,\Omega_p,\phi_1,\dots
,\phi_u)]$$ for which $\nabla\phi_{u+1}$ is contracting against
the (a) crucial factor.
\newline

Furthermore, we define
$$LC^{No\Phi,free}_{\phi_{u+1}}
[C^{l,i_1\dots i_a}_{g}(\Omega_1,\dots ,\Omega_p,\phi_1,\dots
,\phi_u)]$$ to stand for the sublinear combination in
$$LC^{No\Phi}_{\phi_{u+1}} [C^{l,i_1\dots
i_a}_{g}(\Omega_1,\dots ,\Omega_p,\phi_1,\dots ,\phi_u)]$$ for
which the index in $\nabla\phi_{u+1}$ is a free index. It follows
(by just applying the notational conventions of definition
\ref{polypoly}) that:

\begin{equation}
\label{kavaleria} \begin{split} & LC^{No\Phi}_{\phi_{u+1}}
[C^{l,i_1\dots i_a}_{g}(\Omega_1,\dots ,\Omega_p,\phi_1,\dots
,\phi_u)]=LC^{No\Phi,targ}_{\phi_{u+1}} [C^{l,i_1\dots
i_a}_{g}(\Omega_1,\dots ,\Omega_p,\phi_1,\dots ,\phi_u)]
\\&+LC^{No\Phi,free}_{\phi_{u+1}} [C^{l,i_1\dots
i_a}_{g}(\Omega_1,\dots ,\Omega_p,\phi_1,\dots ,\phi_u)].
\end{split}
\end{equation}

\par We easily observe that:

\begin{equation}
\label{eudokia}
\begin{split}
&\\&Xdiv_{i_1}\dots Xdiv_{i_a}LC^{No\Phi,free}_{\phi_{u+1}}
[C^{l,i_1\dots i_a}_{g}(\Omega_1,\dots ,\Omega_p,\phi_1,\dots
,\phi_u)]=
\\&\Sum_{q\in Q} a_q C^q_{g}(\Omega_1,\dots ,\Omega_p,\phi_1,\dots
,\phi_{u+1}).
\end{split}
\end{equation}
These above facts were straightforward. We now explain a more delicate
fact. We consider the set $I_l=\{{}_{i_1},\dots ,{}_{i_a}\}$ of free
indices. We break it up into two subsets. We say ${}_i\in I_1\subset
I$ if ${}_i$ belongs to the crucial factor. We say ${}_i\in I_2$ if it
does not belong to the crucial factor. We then denote by
$$LC^{No\Phi, div,I_2}_{\phi_{u+1}}[Xdiv_{i_1}\dots Xdiv_{i_a}
C^{l,i_1\dots i_a}_{g}(\Omega_1,\dots ,\Omega_p,\phi_1,\dots
,\phi_u)]$$ the sublinear combination in
$$LC^{No\Phi, div}_{\phi_{u+1}}[Xdiv_{i_1}\dots Xdiv_{i_a}
C^{l,i_1\dots i_a}_{g}(\Omega_1,\dots ,\Omega_p,\phi_1,\dots
,\phi_u)]$$ that arises when we apply the transformation law
(\ref{levicivita}) to a pair of indices $(\nabla_{i_y}, {}_a)$
where ${}_{i_y}\in I_2$ and ${}_a$ is an index in the tensor
 field $$C^{l,i_1\dots i_a}_{g}(\Omega_1,\dots ,\Omega_p,
\phi_1,\dots,\phi_u)$$, {\it or} ${}_a$
is another derivative index $\nabla_{i_x},{}_{i_x}\in I_2$
 {\it or} ${}_a$ is a derivative index $\nabla_{i_x}$ with
${}_{i_x}\in I_1$. We also denote by

$$LC^{No\Phi, div,I_1}_{\phi_{u+1}}[Xdiv_{i_1}\dots Xdiv_{i_a}
C^{l,i_1\dots i_a}_{g}(\Omega_1,\dots ,\Omega_p,\phi_1,\dots
,\phi_u)]$$ the sublinear combination in
$$LC^{No\Phi, div}_{\phi_{u+1}}[Xdiv_{i_1}\dots Xdiv_{i_a}
C^{l,i_1\dots i_a}_{g}(\Omega_1,\dots ,\Omega_p,\phi_1,\dots
,\phi_u)]$$ that arises when we apply the transformation law
(\ref{levicivita}) to a pair of indices $(\nabla_{i_y}, {}_b)$
where ${}_{i_y}\in I_1$ and ${}_b$ is an index in the tensor
 field \\$C^{l,i_1\dots i_a}_{g}(\Omega_1,\dots ,\Omega_p,
\phi_1,\dots,\phi_u)$ (same as before) {\it or} ${}_b$ is another
derivative index $\nabla_{i_x}$, ${}_{i_x}\in I_1$.
\newline

\par We derive by definition that:

\begin{equation}
\label{harun2}
\begin{split}
&LC^{No\phi}_{\phi_{u+1}} [\Sum_{l\in L} a_l Xdiv_{i_1}\dots
Xdiv_{i_a} C^{l,i_1\dots i_a}_{g}(\Omega_1,\dots
,\Omega_p,\phi_1,\dots ,\phi_u)]=
\\&LC^{No\phi,div,I_1}_{\phi_{u+1}}
[\Sum_{l\in L} a_l Xdiv_{i_1}\dots Xdiv_{i_a} C^{l,i_1\dots
i_a}_{g}(\Omega_1,\dots ,\Omega_p,\phi_1,\dots ,\phi_u)]+
\\&LC^{No\phi,div,I_2}_{\phi_{u+1}}
[\Sum_{l\in L} a_l Xdiv_{i_1}\dots Xdiv_{i_a} C^{l,i_1\dots
i_a}_{g}(\Omega_1,\dots ,\Omega_p,\phi_1,\dots ,\phi_u)]+
\\& \Sum_{l\in L} a_l Xdiv_{i_1}\dots Xdiv_{i_a}
LC^{No\phi,targ}_{\phi_{u+1}}[C^{l,i_1\dots
i_a}_{g}(\Omega_1,\dots ,\Omega_p,\phi_1,\dots ,\phi_u)]+
\\& \Sum_{l\in L} a_l Xdiv_{i_1}\dots Xdiv_{i_a}
LC^{No\phi,free}_{\phi_{u+1}} [C^{l,i_1\dots
i_a}_{g}(\Omega_1,\dots ,\Omega_p,\phi_1,\dots ,\phi_u)].
\end{split}
\end{equation}

 We then observe  that:

\begin{equation}
\label{greatexp}
\begin{split}
&LC^{No\Phi, div,I_2}_{\phi_{u+1}}[Xdiv_{i_1}\dots Xdiv_{i_a}
C^{l,i_1\dots i_a}_{g}(\Omega_1,\dots ,\Omega_p,\phi_1,\dots
,\phi_u)]=
\\&\Sum_{q\in Q} a_q
C^q_{g}(\Omega_1,\dots,\Omega_p,\phi_1,\dots ,\phi_{u+1}).
\end{split}
\end{equation}

\par Now, let us denote by
$X_{*}div_{i_1}\dots X_{*}div_{i_a}C^{l,i_1\dots
i_a}_{g}(\Omega_1\cdot\phi_{u+1},\dots ,\Omega_p,\phi_1,\dots
,\phi_u)$, $\dots$, $X_{*}div_{i_1}\dots
X_{*}div_{i_a}C^{l,i_1\dots i_a}_{g}(\Omega_1,\dots
,\Omega_X\cdot\phi_{u+1},\dots ,\Omega_p,\phi_1,\dots ,\phi_u)$
the sublinear combination in $Xdiv_{i_1}\dots
Xdiv_{i_a}C^{l,i_1\dots i_a}_{g}(\Omega_1\cdot\phi_{u+1},\dots
,\Omega_p,\phi_1,\dots ,\phi_u)$, $\dots$, $Xdiv_{i_1}\dots
Xdiv_{i_a}C^{l,i_1\dots i_a}_{g}(\Omega_1,\dots
,\Omega_X\cdot\phi_{u+1},\dots,\Omega_p,\phi_1,\dots ,\phi_u)$
that arises when each $\nabla_i, i\in I_1$ is {\it} not allowed to
hit the factor $\phi_{u+1}$.

\par By definition, it follows that for each $h, 1\le h\le X$:

\begin{equation}
\label{sickboy}
\begin{split}
&Xdiv_{i_1}\dots Xdiv_{i_a}C^{l,i_1\dots
i_a}_{g}(\Omega_1,\dots,\Omega_h\cdot \phi_{u+1},\dots
,\Omega_p,\phi_1,\dots ,\phi_u)=
\\& \Sum_{i_y\in
I_1}Xdiv_{i_1}\dots \hat{Xdiv}_{i_y}\dots
Xdiv_{i_\mu}C^{l,i_1\dots i_a}_{g}(\Omega_1,\dots
,\Omega_p,\phi_1,\dots ,\phi_u)\nabla_{i_y}\phi_{u+1}+
\\&X_{*}div_{i_1}\dots X_{*}div_{i_a}C^{l,i_1\dots
i_a}_{g}(\Omega_1,\dots ,\Omega_h\cdot \phi_{u+1},\dots
\Omega_p,\phi_1,\dots ,\phi_u).
\end{split}
\end{equation}

\begin{definition}
\label{definition} We consider each linear combination
$$LC^{No\Phi,targ}_{\phi_{u+1}} [C^{l,i_1\dots
i_a}_{g}(\Omega_1,\dots ,\Omega_p,\phi_1,\dots ,\phi_u)]$$ and we
break it into two sublinear combinations:

$$LC^{No\Phi,targ,A}_{\phi_{u+1}}
[C^{l,i_1\dots i_a}_{g}(\Omega_1,\dots ,\Omega_p,\phi_1,\dots
,\phi_u)]$$ will stand for the sublinear combination that arises
when we apply the transformation law (\ref{levicivita}) to any
 factor other than the selected one.

$$LC^{No\Phi,targ,B}_{\phi_{u+1}}
[C^{l,i_1\dots i_a}_{g}(\Omega_1,\dots ,\Omega_p,\phi_1,\dots
,\phi_u)]$$ will stand for the sublinear combination that arises
when we apply the transformation law (\ref{levicivita}) to the (a)
selected factor.
\end{definition}

\par We then compute another delicate cancellation:

\begin{equation}
\label{greatexp2}
\begin{split}
&Xdiv_{i_1}\dots Xdiv_{i_a}
LC^{No\Phi,targ,A}_{\phi_{u+1}}[C^{l,i_1\dots
i_a}_{g}(\Omega_1,\dots ,\Omega_p,\phi_1,\dots ,\phi_u)] +
\\&X_{*}div_{i_1}\dots X_{*}div_{i_a}C^{l,i_1\dots i_a}_{g}(\Omega_1\cdot\phi_{u+1},\dots
,\Omega_p,\phi_1,\dots ,\phi_u)+\dots
\\& + X_{*}div_{i_1}\dots
X_{*}div_{i_a}C^{l,i_1\dots i_a}_{g}(\Omega_1,\dots,
\Omega_X\cdot\phi_{u+1},\dots ,\Omega_p,\phi_1,\dots ,\phi_u)=
\\&\Sum_{h\in H} a_h Xdiv_{i_1}\dots Xdiv_{i_a} C^{h,i_1\dots
i_a,i_{*}}_{g}(\Omega_1,\dots ,\Omega_p,\phi_1,\dots
,\phi_u)\nabla_{i_{*}}\phi_{u+1}+
\\&\Sum_{q\in Q} a_q C^q_{g}(\Omega_1,\dots ,\Omega_p,\phi_1,\dots
,\phi_{u+1})+\Sum_{z\in Z} a_z C^z_{g}(\Omega_1,\dots
,\Omega_p,\phi_1,\dots ,\phi_{u+1});
\end{split}
\end{equation}
here each $C^{h,i_1\dots i_a,i_{*}}_{g}(\Omega_1,\dots
,\Omega_p,\phi_1,\dots ,\phi_u)\nabla_{i_{*}}\phi_{u+1}$ is an
{\it acceptable} contributor.\footnote{See definition \ref{contributeur}.}
\newline

\par Next, we seek to understand the hardest linear combination in
(\ref{harun2}):

\begin{equation}
\label{zaoush}
\begin{split} &Xdiv_{i_1}\dots Xdiv_{i_a}
LC^{No\Phi,targ,B}_{\phi_{u+1}}[C^{l,i_1\dots
i_a}_{g}(\Omega_1,\dots ,\Omega_p,\phi_1,\dots ,\phi_u)]+
\\&LC^{No\Phi, div,I_1}_{\phi_{u+1}}[Xdiv_{i_1}\dots
Xdiv_{i_a}C^{l,i_1\dots i_a}_{g}(\Omega_1,\dots
,\Omega_p,\phi_1,\dots ,\phi_u)],
\end{split}
\end{equation}
where we recall that we have denoted by
$$LC^{No\Phi, div,I_1}_{\phi_{u+1}}[Xdiv_{i_1}\dots Xdiv_{i_a}
C^{l,i_1\dots i_a}_{g}(\Omega_1,\dots ,\Omega_p,\phi_1,\dots
,\phi_u)]$$ the sublinear combination in
$$LC^{No\Phi, div}_{\phi_{u+1}}[Xdiv_{i_1}\dots Xdiv_{i_a}
C^{l,i_1\dots i_a}_{g}(\Omega_1,\dots ,\Omega_p,\phi_1,\dots
,\phi_u)]$$ that arises when we apply the 
transformation law (\ref{levicivita}) to a
pair of indices $(\nabla_{i_y}, {}_a)$ where ${}_{i_y}\in I_1$ and
${}_a$ is {\it not} a divergence index $\nabla_{i_y}, i_y\in I_2$
(we have already counted those pairs).

\par Firstly, we present our claim in the setting of Lemma \ref{zetajones}.
Recall that in this setting the selected(=crucial) factor is unique.

\par We introduce some language conventions in order to formulate our claim:

\begin{definition}
\label{trelos} We consider each tensor field $C^{l,i_1\dots
i_a}_{g}, l\in L$, and we denote by $\gamma_l$ the total number of
indices (free and non-free) that {\it do not} belong to the
crucial factor and {\it are not} contracting against a factor
$\nabla\phi_y$. We also recall that the number of free indices
that belong to the crucial factor is $|I_1|$, and we let  $\nu_l$
stand for the number of derivatives on the crucial factor
$S_{*}\nabla^{(\nu)}_{r_1\dots r_\nu} R_{ijkl}$. We denote by
$\epsilon_l$ the number of indices in the form ${}_{r_1},\dots
,{}_{r_\nu}, {}_j$ in the crucial factor that are not free and are
not contracting against a factor $\nabla\phi_y$.
\end{definition}
{\it Note:} By abuse of notation, we will write $\gamma,\nu,\epsilon$ 
instead of $\gamma_l,\nu_l,\epsilon_l$. 

\begin{lemma}
\label{travail}
Consider the setting of Lemma \ref{zetajones}, and consider any tensor field
$C^{l,i_1\dots i_a}_{g}$, $l\in L$,\footnote{Recall that by hypothesis $a\ge \mu$.}
 which has  a special free index
${}_{i_1}={}_k$ in the crucial factor. Then:

\begin{equation}
\label{zaoush2}
\begin{split} &Xdiv_{i_1}\dots Xdiv_{i_a}
LC^{No\Phi,targ,B}_{\phi_{u+1}}[C^{l,i_1\dots
i_a}_{g}(\Omega_1,\dots ,\Omega_p,\phi_1,\dots ,\phi_u)]+
\\&LC^{No\Phi, div,I_1}_{\phi_{u+1}}[Xdiv_{i_1}\dots
Xdiv_{i_a}C^{l,i_1\dots i_a}_{g}(\Omega_1,\dots
,\Omega_p,\phi_1,\dots ,\phi_u)]=
\\&-(\gamma+\frac{\nu\cdot\epsilon}{\nu+1}-1)
Xdiv_{i_2}\dots Xdiv_{i_a} C^{l,i_1\dots i_a}_{g}(\Omega_1,\dots
,\Omega_p,\phi_1,\dots ,\phi_u)\nabla_{i_1}\phi_{u+1}+
\\&\Sum_{t\in T} a_t Xdiv_{i_1}\dots Xdiv_{i_{a-1}}
C^{t,i_1\dots i_{a-1},i_{*}}_{g} (\Omega_1,\dots
,\Omega_p,\phi_1,\dots ,\phi_u)\nabla_{i_{*}}\phi_{u+1}+
\\&\Sum_{h\in H} a_h Xdiv_{i_1}\dots Xdiv_{i_\mu}C^{h,i_1\dots
i_\mu,i_{*}}_{g}(\Omega_1,\dots ,\Omega_p,\phi_1,\dots
,\phi_u)\nabla_{i_{*}}\phi_{u+1}+
\\&\Sum_{q\in Q} a_q
C^q_{g}(\Omega_1,\dots ,\Omega_p, \phi_1,\dots ,\phi_{u+1});
\end{split}
\end{equation}
here
$$\Sum_{t\in T} a_t C^{t,i_1\dots i_{a-1},i_{*}}_{g}
(\Omega_1,\dots ,\Omega_p,\phi_1,\dots
,\phi_u)\nabla_{i_{*}}\phi_{u+1}$$ stands for a generic linear
combination of acceptable $(a-1)$-tensor fields for which ${}_{i_{*}}$
is the index ${}_k$ in the crucial factor, but we have fewer than
$|I_1|-1$ free
 indices in the crucial factor. In particular, if $a=\mu$,
then the $(\mu-1)$-tensor field will have a refined double
 character that is doubly subsequent to each $\vec{L^z}', z\in Z'_{Max}$.

If $C^{l,i_1\dots i_a}_{g}$ does not have a free index in the
position ${}_k$ or ${}_l$ in the selected factor
$S_{*}\nabla^{(\nu)}R_{ijkl}$, then:

\begin{equation}
\label{zaoush2b}
\begin{split} &Xdiv_{i_1}\dots Xdiv_{i_a}
LC^{No\Phi,targ,B}_{\phi_{u+1}}[C^{l,i_1\dots
i_a}_{g}(\Omega_1,\dots ,\Omega_p,\phi_1,\dots ,\phi_u)]+
\\&LC^{No\Phi, div,I_1}_{\phi_{u+1}}[Xdiv_{i_1}\dots
Xdiv_{i_a}C^{l,i_1\dots i_a}_{g}(\Omega_1,\dots
,\Omega_p,\phi_1,\dots ,\phi_u)]=
\\&\Sum_{h\in H} a_h Xdiv_{i_1}\dots Xdiv_{i_a}
C^{h,i_1\dots i_a,i_{*}}_{g}(\Omega_1,\dots ,\Omega_p,\phi_1,\dots
,\phi_u)\nabla_{i_{*}}\phi_{u+1}+
\\&\Sum_{q\in Q} a_q
C^q_{g}(\Omega_1,\dots ,\Omega_p, \phi_1,\dots ,\phi_{u+1}).
\end{split}
\end{equation}
\end{lemma}

{\it Proof of Lemma \ref{travail}:}
\newline

\par We start with the first claim, which is the hardest.
We will show the above by breaking the left hand side into
numerous sublinear combinations. Recall that we are assuming that
the free index ${}_{i_1}$ is the index ${}_k$ in the crucial factor
$S_{*}\nabla^{(\nu)} R_{ijkl}$, while the other free indices that
belong to the crucial factor are ${}_{i_2},\dots ,{}_{i_{|I_1|}}$. Firstly,
let us analyze the sublinear combination

$$LC^{No\Phi, div,I_1}_{\phi_{u+1}}[Xdiv_{i_1}\dots
Xdiv_{i_a}C^{l,i_1\dots i_a}_{g}(\Omega_1,\dots
,\Omega_p,\phi_1,\dots ,\phi_u)].$$ We break this sum into four
sublinear combinations: First, we consider the sublinear combination
that arises when we apply the transformation law
(\ref{levicivita}) to a pair of indices $(\nabla_{i_1},{}_b)$ where
${}_b$ is an original index in $C^{l,i_1\dots i_a}_{g}$. Secondly,
we consider the sublinear combination that arises when we apply the
transformation law (\ref{levicivita}) to a pair of divergence
indices, $(\nabla_{i_1},\nabla_{i_k})$, $2\le k\le |I_1|$.
Thirdly, we consider the sublinear combination that arises when we
apply the transformation law (\ref{levicivita}) to a pair of
divergence indices $({\nabla}_{i_k},{\nabla}_{i_l})$, $2\le k,l\le |I_1|$. Fourthly, we
consider the sublinear combination that arises when we apply the
transformation law (\ref{levicivita}) to a pair of indices
$({}_{i_k},{}_b)$, $2\le k\le |I_1|$ and ${}_b$ being an original index
in $C^{l,i_1\dots i_a}_{g}$. We respectively denote those
sublinear combinations by $LC^{No\Phi,
div,I_1,\alpha}_{\phi_{u+1}}[\dots]$, $LC^{No\Phi,
div,I_1,\beta}_{\phi_{u+1}}[\dots]$, $LC^{No\Phi,
div,I_1,\gamma}_{\phi_{u+1}}[\dots]$, $LC^{No\Phi,
div,I_1,\delta}_{\phi_{u+1}}[\dots]$.

\par We then observe that:

\begin{equation}
\label{hef1} \begin{split} &LC^{No\Phi,
div,I_1,\alpha}_{\phi_{u+1}}[Xdiv_{i_1}\dots
Xdiv_{i_a}C^{l,i_1\dots i_a}_{g}(\Omega_1,\dots
,\Omega_p,\phi_1,\dots ,\phi_u)]=
\\&-(\gamma-1) C^{l,i_1\dots
i_a}_{g}(\Omega_1,\dots ,\Omega_p,\phi_1,\dots
,\phi_u)\nabla_{i_1}\phi_{u+1}\\&+\Sum_{q\in Q} a_q
C^q_{g}(\Omega_1,\dots ,\Omega_p,\phi_1,\dots ,\phi_{u+1}).
\end{split}
\end{equation}

The second sublinear combination is a little more complicated.

\begin{equation}
\label{hef2} \begin{split} &LC^{No\Phi,
div,I_1,\beta}_{\phi_{u+1}}[Xdiv_{i_1}\dots
Xdiv_{i_a}C^{l,i_1\dots i_a}_{g}(\Omega_1,\dots
,\Omega_p,\phi_1,\dots ,\phi_u)]=
\\&-(|I_1|-1) Xdiv_{i_2}\dots Xdiv_{i_\alpha}C^{l,i_1\dots
i_a}_{g}(\Omega_1,\dots ,\Omega_p,\phi_1,\dots
,\phi_u)\nabla_{i_1}\phi_{u+1}+
\\& \Sum_{q\in Q} a_q
C^q_{g}(\Omega_1,\dots ,\Omega_p,\phi_1,\dots ,\phi_{u+1}).
\end{split}
\end{equation}

\par On the other hand, we also see that:

\begin{equation}
\label{hef3} \begin{split} &LC^{No\Phi,
div,I_1,\gamma}_{\phi_{u+1}}[Xdiv_{i_1}\dots
Xdiv_{i_a}C^{l,i_1\dots i_a}_{g}(\Omega_1,\dots
,\Omega_p,\phi_1,\dots ,\phi_u)]=
\\&\Sum_{q\in Q} a_q
C^q_{g}(\Omega_1,\dots ,\Omega_p,\phi_1,\dots ,\phi_{u+1}).
\end{split}
\end{equation}

\par Lastly, to describe the fourth sublinear combination, we
introduce some notation: For each $k, 2\le k\le |I_1|$ we define
$\hat{C}^{l,i_1\dots i_a}_{g}$ to stand for the sublinear
combination which arises from $C^{l,i_1\dots i_a}_{g}$ by
performing a cyclic permutation of the indices
${}_{i_k},{}_{i_1},{}_l$ (${}_{i_k}$ is picked out arbitrarily
among ${}_{i_2},\dots ,{}_{i_{|I_1|}}$) in the crucial factor
$S_{*}\nabla^{(\nu)}_{r_1\dots r_\nu} R_{ijkl}$. We then conclude:

\begin{equation}
\label{hef4} \begin{split} &LC^{No\Phi,
div,I_1,\delta}_{\phi_{u+1}}[Xdiv_{i_1}\dots
Xdiv_{i_a}C^{l,i_1\dots i_a}_{g}(\Omega_1,\dots
,\Omega_p,\phi_1,\dots ,\phi_u)]=
\\&-(|I_1|-1) Xdiv_{i_2}\dots Xdiv_{i_\alpha}\hat{C}^{l,i_1\dots
i_a}_{g}(\Omega_1,\dots ,\Omega_p,\phi_1,\dots
,\phi_u)\nabla_{i_1}\phi_{u+1}
\\&+ \Sum_{q\in Q} a_q
C^q_{g}(\Omega_1,\dots ,\Omega_p,\phi_1,\dots ,\phi_{u+1}).
\end{split}
\end{equation}

\par Now, just by the first and second Bianchi identity we
 observe that:

\begin{equation}
\label{hef5} \begin{split} & -(|I_1|-1) Xdiv_{i_2}\dots
Xdiv_{i_\alpha}C^{l,i_1\dots i_a}_{g}(\Omega_1,\dots
,\Omega_p,\phi_1,\dots ,\phi_u)\nabla_{i_1}\phi_{u+1}
\\&-(|I_1|-1) Xdiv_{i_2}\dots Xdiv_{i_\alpha}\hat{C}^{l,i_1\dots
i_a}_{g}(\Omega_1,\dots ,\Omega_p,\phi_1,\dots
,\phi_u)\nabla_{i_1}\phi_{u+1}
\\&= \Sum_{q\in Q} a_q
C^q_{g}(\Omega_1,\dots ,\Omega_p,\phi_1,\dots ,\phi_{u+1}).
\end{split}
\end{equation}

\par Next, we seek to analyze the sublinear combination:

$$Xdiv_{i_1}\dots Xdiv_{i_a}
LC^{No\Phi,targ,B}_{\phi_{u+1}}[C^{l,i_1\dots
i_a}_{g}(\Omega_1,\dots ,\Omega_p,\phi_1,\dots ,\phi_u)].$$ We
have by definition that this sublinear combination can only arise
by applying the last term in (\ref{levicivita}) to a pair of
indices in the selected factor. We now break it into five sublinear
combinations: We define the first sublinear combination to be the
one that arises when we apply the fourth summand in
(\ref{levicivita}) to a pair of indices $({}_{i_1},{}_b)$, where
${}_b$ is an original non-free index in the selected factor in
$C^{l,i_1\dots i_a}_{g}$. We define the second sublinear
combination to be one that arises
 by applying the last term in (\ref{levicivita}) to two free
 indices
$({}_{i_1},{}_{i_k})$, $2\le k\le |I_1|$ in the selected factor. We
define the third to be the one that arises by applying the last
term in (\ref{levicivita}) to a pair $({}_{i_k},{}_b)$ where $k\ge
2$ and the index ${}_b$ is an original non-free index in the
crucial factor. The fourth sublinear combination arises  when we
apply the last term in (\ref{levicivita}) to a pair
$({}_{i_k},{}_{i_l})$ of free indices in the crucial factor, $2\le
k,l\le |I_1|$. Lastly, the fifth sublinear combination is the one
that arises by applying the last term in (\ref{levicivita}) to a
pair of non-free original indices in $C^{l,i_1\dots i_a}_{g}$. We
denote these sublinear combinations by
$LC^{No\Phi,targ,B,\alpha}_{\phi_{u+1}}[\dots]$,
$LC^{No\Phi,targ,B,\beta}_{\phi_{u+1}}[\dots]$,
$LC^{No\Phi,targ,B,\gamma}_{\phi_{u+1}}[\dots]$,
$LC^{No\Phi,targ,B,\delta}_{\phi_{u+1}}[\dots]$,
$LC^{No\Phi,targ,B,\varepsilon}_{\phi_{u+1}}[\dots]$. It follows
that:

\begin{equation}
\label{cheeg1}
\begin{split}
&Xdiv_{i_1}\dots Xdiv_{i_a}
LC^{No\Phi,targ,B,\alpha}_{\phi_{u+1}}[C^{l,i_1\dots
i_a}_{g}(\Omega_1,\dots ,\Omega_p,\phi_1,\dots ,\phi_u)]=
\\& -\frac{\nu\cdot\epsilon}{\nu+1}Xdiv_{i_2}\dots
Xdiv_{i_\alpha}C^{l,i_1\dots i_a}_{g}(\Omega_1,\dots
,\Omega_p,\phi_1,\dots ,\phi_u)\nabla_{i_1}\phi_{u+1}+
\\&\Sum_{h\in H} a_h Xdiv_{i_1}\dots Xdiv_{i_a}C^{h,i_1\dots
i_\mu,i_{*}}_{g}(\Omega_1,\dots ,\Omega_p,\phi_1,\dots
,\phi_u)\nabla_{i_{*}}\phi_{u+1},
\end{split}
\end{equation}

\begin{equation}
\label{cheeg2}
\begin{split}
&Xdiv_{i_1}\dots Xdiv_{i_a}
LC^{No\Phi,targ,B,\beta}_{\phi_{u+1}}[C^{l,i_1\dots
i_a}_{g}(\Omega_1,\dots ,\Omega_p,\phi_1,\dots ,\phi_u)]=
\\& -\frac{\nu\cdot(|I_1|-1)}{\nu+1}Xdiv_{i_2}\dots
Xdiv_{i_\alpha}C^{l,i_1\dots i_a}_{g}(\Omega_1,\dots
,\Omega_p,\phi_1,\dots ,\phi_u)\nabla_{i_1}\phi_{u+1}+
\\&\Sum_{t\in T} a_t Xdiv_{i_1}\dots Xdiv_{i_\mu}C^{h,i_1\dots
i_{a-1},i_{*}}_{g}(\Omega_1,\dots ,\Omega_p,\phi_1,\dots
,\phi_u)\nabla_{i_{*}}\phi_{u+1}+
\\&\Sum_{q\in Q} a_q
C^q_{g}(\Omega_1,\dots ,\Omega_p,\phi_1,\dots ,\phi_{u+1}),
\end{split}
\end{equation}

\begin{equation}
\label{cheeg3}
\begin{split}
&Xdiv_{i_1}\dots Xdiv_{i_a}
LC^{No\Phi,targ,B,\gamma}_{\phi_{u+1}}[C^{l,i_1\dots
i_a}_{g}(\Omega_1,\dots ,\Omega_p,\phi_1,\dots ,\phi_u)]=
\\& -\frac{\nu\cdot(|I_1|-1)}{\nu+1}Xdiv_{i_2}\dots
Xdiv_{i_\alpha}\hat{C}^{l,i_1\dots i_a}_{g}(\Omega_1,\dots
,\Omega_p,\phi_1,\dots ,\phi_u)\nabla_{i_1}\phi_{u+1}+
\\&\Sum_{h\in H} a_h Xdiv_{i_1}\dots Xdiv_{i_a}C^{h,i_1\dots
i_\mu,i_{*}}_{g}(\Omega_1,\dots ,\Omega_p,\phi_1,\dots
,\phi_u)\nabla_{i_{*}}\phi_{u+1} +
\\&\Sum_{q\in Q} a_q
C^q_{g}(\Omega_1,\dots ,\Omega_p,\phi_1,\dots ,\phi_{u+1}),
\end{split}
\end{equation}

\begin{equation}
\label{cheeg4}
\begin{split}
&Xdiv_{i_1}\dots Xdiv_{i_a}
LC^{No\Phi,targ,B,\delta}_{\phi_{u+1}}[C^{l,i_1\dots
i_a}_{g}(\Omega_1,\dots ,\Omega_p,\phi_1,\dots ,\phi_u)]=
\\&\Sum_{q\in Q} a_q C^q_{g}(\Omega_1,\dots
,\Omega_p,\phi_1,\dots ,\phi_{u+1}),
\end{split}
\end{equation}

\begin{equation}
\label{cheeg5}
\begin{split}
&Xdiv_{i_1}\dots Xdiv_{i_a}
LC^{No\Phi,targ,B,\epsilon}_{\phi_{u+1}}[C^{l,i_1\dots
i_a}_{g}(\Omega_1,\dots ,\Omega_p,\phi_1,\dots ,\phi_u)]=
\\&\Sum_{h\in H} a_h Xdiv_{i_1}\dots Xdiv_{i_a}C^{h,i_1\dots
i_\mu,i_{*}}_{g}(\Omega_1,\dots ,\Omega_p,\phi_1,\dots
,\phi_u)\nabla_{i_{*}}\phi_{u+1} +
\\&\Sum_{q\in Q} a_q
C^q_{g}(\Omega_1,\dots ,\Omega_p,\phi_1,\dots ,\phi_{u+1}).
\end{split}
\end{equation}

\par Adding all the above we derive the first claim of our Lemma,
(where the selected factor has a special free index ${}_{i_1}={}_k$).

\par The second claim of our Lemma (where there is no special index in the crucial factor)
 follows more easily. We now have that ${}_{i_1}$ is
 not a special free index, so we will now consider all the
 sublinear combinations above where ${}_{i_1}$ is not mentioned,
 and also whenever we mentioned above one of the free indices
  ${}_{i_2},\dots ,{}_{i_{|I_1|}}$ we will now read ``one of the free indices
   ${}_{i_1},\dots ,{}_{i_{|I_1|}}$'' (since the index ${}_{i_1}$ is not
  special now). We then find that all the relevant equations
  will hold, with the exception of (\ref{hef4}),
 (\ref{cheeg3}), which now  become:

\begin{equation}
\label{hef4'} \begin{split} &LC^{No\Phi,
div,I_1,\delta}_{\phi_{u+1}}[Xdiv_{i_1}\dots
Xdiv_{i_a}C^{l,i_1\dots i_a}_{g}(\Omega_1,\dots
,\Omega_p,\phi_1,\dots ,\phi_u)]=
\\& \Sum_{q\in Q} a_q
C^q_{g}(\Omega_1,\dots ,\Omega_p,\phi_1,\dots ,\phi_{u+1}),
\end{split}
\end{equation}

\begin{equation}
\label{cheeg3'}
\begin{split}
&Xdiv_{i_1}\dots Xdiv_{i_a}
LC^{No\Phi,targ,B,\gamma}_{\phi_{u+1}}[C^{l,i_1\dots
i_a}_{g}(\Omega_1,\dots ,\Omega_p,\phi_1,\dots ,\phi_u)]=
\\&\Sum_{h\in H} a_h Xdiv_{i_1}\dots Xdiv_{i_a}C^{h,i_1\dots
i_\mu,i_{*}}_{g}(\Omega_1,\dots ,\Omega_p,\phi_1,\dots
,\phi_u)\nabla_{i_{*}}\phi_{u+1} +
\\&\Sum_{q\in Q} a_q
C^q_{g}(\Omega_1,\dots ,\Omega_p,\phi_1,\dots ,\phi_{u+1}),
\end{split}
\end{equation}
by application of the first and second Bianchi identity.
This concludes the proof of our Lemma. $\Box$
\newline

\par Now, we study the sublinear combination (\ref{zaoush})
 in the setting of Lemma \ref{pool2}.
We recall the discussion from the introduction in \cite{alexakis4} on the
{\it crucial factor}.\footnote{Recall that in this setting the 
``crucial'' and ``selected'' factors coincide.} We recall that the crucial factor is defined
{\it in terms of the $u$-simple character $\vec{\kappa}_{simp}$}
by examining the tensor fields
 $C^{l,i_1\dots i_\mu}_{g}$, $l\in L^z, z\in Z'_{Max}$ (for a precise definition
 see the discussion above the statement of Lemma \ref{zetajones}).
Once it has been defined, we may speak of the crucial factor(s)
for {\it any} tensor field with the $u$-simple character
$\vec{\kappa}_{simp}$ (in fact, even for any complete contractions
with a weak character $Weak(\vec{\kappa}_{simp})$). In each tensor
field $C^{l,i_1\dots i_a}, l\in L$, we will denote by $\{
T_1,\dots ,T_M\}$ the set of crucial factors.

We will now separately consider the sublinear combinations in
\begin{equation}
\label{tamas}
\begin{split}
&Xdiv_{i_1}\dots Xdiv_{i_a}
LC^{No\Phi,targ,B}_{\phi_{u+1}}[C^{l,i_1\dots
i_a}_{g}(\Omega_1,\dots ,\Omega_p,\phi_1,\dots ,\phi_u)]+
\\&LC^{No\Phi, div,I_1}_{\phi_{u+1}}[Xdiv_{i_1}\dots
Xdiv_{i_a}C^{l,i_1\dots i_a}_{g}(\Omega_1,\dots
,\Omega_p,\phi_1,\dots ,\phi_u)]
\end{split}
\end{equation}
that have a factor $\nabla\phi_{u+1}$ contracting against $T_1,
\dots ,T_M$. We use the symbols $LC^{No\Phi,targ,B,T_i}$
and $LC^{No\Phi, div,I_1,T_i}$ to illustrate that we are
considering those sublinear combinations. Again, we will denote by

$$\Sum_{t\in T} a_t C^{t,i_1\dots i_{\mu-1},i_{*}}_{g}
(\Omega_1,\dots ,\Omega_p,\phi_1,\dots
,\phi_u)\nabla_{i_{*}}\phi_{u+1}$$ a generic linear combination of
{\it acceptable} $(\mu-1)$-tensor fields which have a simple
character $\vec{\kappa}^{+}_{simp}$ but are doubly subsequent to
each $\vec{L^z}, z\in Z'_{Max}$.

\begin{definition}
\label{epsilongamma2} 
Here,
$\epsilon^l_i$ will stand for the number of derivative indices
in the crucial factor $T_i=\nabla^{(m)}R_{ijkl}$
 that are not free and are not contracting against a factor
 $\nabla\phi_h$. $\gamma^l_i$ will stand for the number of indices
  in the other factors in $C^{l,i_1\dots i_a}_{g}$
 that are not contracting against a factor $\nabla\phi_h$.
\end{definition}
{\it Note:} By abuse of notation, we will be writing $\epsilon_i,\gamma_i$ 
instead of $\epsilon^l_i,\gamma^l_i$ from now on). We claim:

\begin{lemma}
\label{travail2}
 If $C^{l,i_1\dots i_a}_{g}(\Omega_1,\dots ,\Omega_p,\phi_1,\dots ,\phi_u)$ has one
internal free index (say ${}_{i_1}$) in the crucial factor
$T_i=\nabla^{(m)}R_{i_1jkl}$, then:

\begin{equation}
\label{zaoush3}
\begin{split} &Xdiv_{i_1}\dots Xdiv_{i_a}
LC^{No\Phi,targ,B,T_i}_{\phi_{u+1}}[C^{l,i_1\dots
i_a}_{g}(\Omega_1,\dots ,\Omega_p,\phi_1,\dots ,\phi_u)]+
\\&LC^{No\Phi, div,I_1,T_i}_{\phi_{u+1}}[Xdiv_{i_1}\dots
Xdiv_{i_a}C^{l,i_1\dots i_a}_{g}(\Omega_1,\dots
,\Omega_p,\phi_1,\dots ,\phi_u)]=
\\&-(\gamma_i+\epsilon_i)Xdiv_{i_1}\dots Xdiv_{i_a}
C^{l,i_1\dots i_a}_{g}(\Omega_1,\dots ,\Omega_p,\phi_1,\dots
,\phi_u)\nabla_{i_1}\phi_{u+1}+
\\&\Sum_{t\in T} a_t Xdiv_{i_1}\dots Xdiv_{i_{a-1}}
C^{t,i_1\dots i_{a-1},i_{*}}_{g} (\Omega_1,\dots
,\Omega_p,\phi_1,\dots ,\phi_u)\nabla_{i_{*}}\phi_{u+1}+
\\&\Sum_{h\in H} a_h Xdiv_{i_1}\dots Xdiv_{i_\mu}
C^{h,i_1\dots i_\mu,i_{*}}_{g}(\Omega_1,\dots
,\Omega_p,\phi_1,\dots ,\phi_u)\nabla_{i_{*}}\phi_{u+1}+
\\&\Sum_{q\in Q} a_q
C^q_{g}(\Omega_1,\dots ,\Omega_p, \phi_1,\dots ,\phi_{u+1}).
\end{split}
\end{equation}

\par Moreover, in the case of Lemma \ref{pool2} and if
$C^{l,i_1\dots i_a}_{g}(\Omega_1,\dots ,\Omega_p,\phi_1,\dots
,\phi_u)$ has two internal free indices (say ${}_{i_1}$ and
${}_{i_2}$) in the crucial factor $T_i=\nabla^{(m)}R_{i_1ji_2l}$,
then:

\begin{equation}
\label{zaoush3b}
\begin{split} &Xdiv_{i_1}\dots Xdiv_{i_a}
LC^{No\Phi,targ,B,T_i}_{\phi_{u+1}}[C^{l,i_1\dots
i_a}_{g}(\Omega_1,\dots ,\Omega_p,\phi_1,\dots ,\phi_u)]+
\\&LC^{No\Phi, div,I_1,T_i}_{\phi_{u+1}}[Xdiv_{i_1}\dots
Xdiv_{i_a}C^{l,i_1\dots i_a}_{g}(\Omega_1,\dots
,\Omega_p,\phi_1,\dots ,\phi_u)]=
\\&-(\gamma_i+\epsilon_i)Xdiv_{i_2}\dots Xdiv_{i_a}
C^{l,i_1\dots i_a}_{g}(\Omega_1,\dots ,\Omega_p,\phi_1,\dots
,\phi_u)\nabla_{i_1}\phi_{u+1}
\\&-(\gamma_i+\epsilon_i)Xdiv_{i_1}Xdiv_{i_3}\dots Xdiv_{i_a}
C^{l,i_1\dots i_a}_{g}(\Omega_1,\dots ,\Omega_p,\phi_1,\dots
,\phi_u)\nabla_{i_2}\phi_{u+1}+
\\&\Sum_{t\in T} a_t Xdiv_{i_1}\dots Xdiv_{i_{a-1}}
C^{t,i_1\dots i_{a-1},i_{*}}_{g} (\Omega_1,\dots
,\Omega_p,\phi_1,\dots ,\phi_u)\nabla_{i_{*}}\phi_{u+1}+
\\&\Sum_{h\in H} a_h Xdiv_{i_1}\dots Xdiv_{i_\mu}C^{h,i_1\dots
i_\mu,i_{*}}_{g}(\Omega_1,\dots ,\Omega_p,\phi_1,\dots
,\phi_u)\nabla_{i_{*}}\phi_{u+1}+
\\&\Sum_{q\in Q} a_q
C^q_{g}(\Omega_1,\dots ,\Omega_p, \phi_1,\dots ,\phi_{u+1}).
\end{split}
\end{equation}
Finally, if the crucial factor $T_i$ has no internal free
 indices then:
\begin{equation}
\label{zaoush3c}
\begin{split} &Xdiv_{i_1}\dots Xdiv_{i_a}
LC^{No\Phi,targ,B,T_i}_{\phi_{u+1}}[C^{l,i_1\dots
i_a}_{g}(\Omega_1,\dots ,\Omega_p,\phi_1,\dots ,\phi_u)]+
\\&LC^{No\Phi, div,I_1,T_i}_{\phi_{u+1}}[Xdiv_{i_1}\dots
Xdiv_{i_a}C^{l,i_1\dots i_a}_{g}(\Omega_1,\dots
,\Omega_p,\phi_1,\dots ,\phi_u)]=
\\&\Sum_{h\in H} a_h Xdiv_{i_1}\dots Xdiv_{i_\mu}
C^{h,i_1\dots i_\mu,i_{*}}_{g}(\Omega_1,\dots
,\Omega_p,\phi_1,\dots ,\phi_u)\nabla_{i_{*}}\phi_{u+1}+
\\&\Sum_{q\in Q} a_q
C^q_{g}(\Omega_1,\dots ,\Omega_p, \phi_1,\dots ,\phi_{u+1}).
\end{split}
\end{equation}
\end{lemma}

{\it Proof of Lemma \ref{travail2}:}
\newline

\par The proof of this claim is similar  to the
previous one. We start with the case where the crucial factor $T_i$ has one internal
free index.

\par We again denote by ${}_{i_1}$ the one internal free index
in the crucial factor $T_i$  and we denote by ${}_{i_2},\dots
,{}_{i_{|I_1|}}$ the other free indices. We divide the sublinear
combination $LC^{No\Phi,div,I_1}_{\phi_{u+1}}[\dots]$ into further
sublinear combinations (indexed by ${}_\alpha,\dots ,{}_\epsilon$)
as in the previous case.
 We calculate:

\begin{equation}
\label{ghef1} \begin{split} &LC^{No\Phi,
div,I_1,\alpha}_{\phi_{u+1}}[Xdiv_{i_1}\dots
Xdiv_{i_a}C^{l,i_1\dots i_a}_{g}(\Omega_1,\dots
,\Omega_p,\phi_1,\dots ,\phi_u)]=
\\&-\gamma C^{l,i_1\dots
i_a}_{g}(\Omega_1,\dots ,\Omega_p,\phi_1,\dots
,\phi_u)\nabla_{i_1}\phi_{u+1}+\Sum_{q\in Q} a_q
C^q_{g}(\Omega_1,\dots ,\Omega_p,\phi_1,\dots ,\phi_{u+1});
\end{split}
\end{equation}
(we have used the first Bianchi identity here).

The second sublinear combination is a little more complicated.

\begin{equation}
\label{ghef2} \begin{split} &LC^{No\Phi,
div,I_1,\beta}_{\phi_{u+1}}[Xdiv_{i_1}\dots
Xdiv_{i_a}C^{l,i_1\dots i_a}_{g}(\Omega_1,\dots
,\Omega_p,\phi_1,\dots ,\phi_u)]=
\\&-(|I_1|-1) Xdiv_{i_2}\dots Xdiv_{i_\alpha}C^{l,i_1\dots
i_a}_{g}(\Omega_1,\dots ,\Omega_p,\phi_1,\dots
,\phi_u)\nabla_{i_1}\phi_{u+1}+
\\& \Sum_{q\in Q} a_q
C^q_{g}(\Omega_1,\dots ,\Omega_p,\phi_1,\dots ,\phi_{u+1}).
\end{split}
\end{equation}

\par On the other hand, we also see that:

\begin{equation}
\label{ghef3} \begin{split} &LC^{No\Phi,
div,I_1,\gamma}_{\phi_{u+1}}[Xdiv_{i_1}\dots
Xdiv_{i_a}C^{l,i_1\dots i_a}_{g}(\Omega_1,\dots
,\Omega_p,\phi_1,\dots ,\phi_u)]=
\\&\Sum_{q\in Q} a_q
C^q_{g}(\Omega_1,\dots ,\Omega_p,\phi_1,\dots ,\phi_{u+1}).
\end{split}
\end{equation}

\par Lastly, to describe the fourth sublinear combination, we
introduce some notation: For each $k, 2\le k\le |I_1|$ we define
$\hat{C}^{l,i_1\dots i_a}_{g}$ to stand for the sublinear
combination which arises from $C^{l,i_1\dots i_a}_{g}$ by
performing a cyclic permutation of the indices ${}_{i_k},{}_{i_1},{}_j$ in the
crucial factor $\nabla^{(m)} R_{ijkl}$. We then have that:

\begin{equation}
\label{ghef4} \begin{split} &LC^{No\Phi,
div,I_1,\delta}_{\phi_{u+1}}[Xdiv_{i_1}\dots
Xdiv_{i_a}C^{l,i_1\dots i_a}_{g}(\Omega_1,\dots
,\Omega_p,\phi_1,\dots ,\phi_u)]=
\\&-(|I_1|-1) Xdiv_{i_2}\dots Xdiv_{i_\alpha}\hat{C}^{l,i_1\dots
i_a}_{g}(\Omega_1,\dots ,\Omega_p,\phi_1,\dots
,\phi_u)\nabla_{i_1}\phi_{u+1}+
\\& \Sum_{q\in Q} a_q
C^q_{g}(\Omega_1,\dots ,\Omega_p,\phi_1,\dots ,\phi_{u+1});
\end{split}
\end{equation}
(we have used the second Bianchi).

\par Now, just by the first and second Bianchi identity we
 derive that:

\begin{equation}
\label{ghef5} \begin{split} & -(|I_1|-1) Xdiv_{i_2}\dots
Xdiv_{i_a}C^{l,i_1\dots i_a}_{g}(\Omega_1,\dots
,\Omega_p,\phi_1,\dots ,\phi_u)\nabla_{i_1}\phi_{u+1}
\\&-(|I_1|-1) Xdiv_{i_2}\dots Xdiv_{i_a}\hat{C}^{l,i_1\dots
i_a}_{g}(\Omega_1,\dots ,\Omega_p,\phi_1,\dots
,\phi_u)\nabla_{i_1}\phi_{u+1}=
\\& \Sum_{q\in Q} a_q
C^q_{g}(\Omega_1,\dots ,\Omega_p,\phi_1,\dots ,\phi_{u+1}).
\end{split}
\end{equation}

\par Now, we study the second sublinear combination:

$$Xdiv_{i_1}\dots Xdiv_{i_a}
LC^{No\Phi,targ,B,T_i}_{\phi_{u+1}}[C^{l,i_1\dots
i_a}_{g}(\Omega_1,\dots ,\Omega_p,\phi_1,\dots ,\phi_u)].$$
By definition and by the transformation law (\ref{levicivita})
 that this sublinear combination can only arise by applying
the last term in (\ref{levicivita}) to a pair of indices in the
crucial factor $T_i$. We now break it into five sublinear
combinations: We define the first sublinear combination to be the
one that arises when we apply the fourth summand in
(\ref{levicivita}) to a pair of indices $({}_{i_1},{}_b)$, where
${}_b$ is an original non-free index in the crucial factor in
$C^{l,i_1\dots i_a}_{g}$. We define the second sublinear
combination to be the one that arises
 by applying the last term in (\ref{levicivita}) to two free
 indices
$({}_{i_1},{}_{i_k})$, $2\le k\le |I_1|$ in the crucial factor
$T_i$. We define the third to be the one that arises by applying
the last term in (\ref{levicivita}) to a pair $({}_{i_k},{}_b)$
where $k\ge 2$ and the index ${}_b$ is an original non-free index
in the crucial factor. The fourth is when we apply the last term
in (\ref{levicivita}) to a pair $({}_{i_k},{}_{i_l})$ of free
indices in the crucial factor, $2\le k,l\le |I_1|$. Lastly, the
fifth sublinear combination is the one that arises by applying the
last term in (\ref{levicivita}) to a pair of non-free original
indices in $C^{l,i_1\dots i_a}_{g}$. We denote these sublinear
combinations by
$LC^{No\Phi,targ,B,T_i,\alpha}_{\phi_{u+1}}[\dots]$,
$LC^{No\Phi,targ,B,T_i,\beta}_{\phi_{u+1}}[\dots]$,
$LC^{No\Phi,targ,B,T_i,\gamma}_{\phi_{u+1}}[\dots]$,
$LC^{No\Phi,targ,B,T_i,\delta}_{\phi_{u+1}}[\dots]$,
$LC^{No\Phi,targ,B,T_i,\varepsilon}_{\phi_{u+1}}[\dots]$. It
follows that:

\begin{equation}
\label{bcheeg1}
\begin{split}
&Xdiv_{i_1}\dots Xdiv_{i_a}
LC^{No\Phi,targ,B,T_i,\alpha}_{\phi_{u+1}}[C^{l,i_1\dots
i_a}_{g}(\Omega_1,\dots ,\Omega_p,\phi_1,\dots ,\phi_u)]=
\\& -\epsilon_iXdiv_{i_2}\dots
Xdiv_{i_\alpha}C^{l,i_1\dots i_a}_{g}(\Omega_1,\dots
,\Omega_p,\phi_1,\dots ,\phi_u)\nabla_{i_1}\phi_{u+1}+
\\&\Sum_{h\in H} a_h Xdiv_{i_1}\dots Xdiv_{i_a}C^{h,i_1\dots
i_\mu,i_{*}}_{g}(\Omega_1,\dots ,\Omega_p,\phi_1,\dots
,\phi_u)\nabla_{i_{*}}\phi_{u+1},
\end{split}
\end{equation}

\begin{equation}
\label{bcheeg2}
\begin{split}
&Xdiv_{i_1}\dots Xdiv_{i_a}
LC^{No\Phi,targ,B,T_i,\beta}_{\phi_{u+1}}[C^{l,i_1\dots
i_a}_{g}(\Omega_1,\dots ,\Omega_p,\phi_1,\dots ,\phi_u)]=
\\& -(|I_1|-1)Xdiv_{i_2}\dots
Xdiv_{i_\alpha}C^{l,i_1\dots i_a}_{g}(\Omega_1,\dots
,\Omega_p,\phi_1,\dots ,\phi_u)\nabla_{i_1}\phi_{u+1}+
\\&\Sum_{t\in T} a_t Xdiv_{i_1}\dots Xdiv_{i_\mu}C^{h,i_1\dots
i_{a-1},i_{*}}_{g}(\Omega_1,\dots ,\Omega_p,\phi_1,\dots
,\phi_u)\nabla_{i_{*}}\phi_{u+1}+
\\&\Sum_{q\in Q} a_q
C^q_{g}(\Omega_1,\dots ,\Omega_p,\phi_1,\dots ,\phi_{u+1}),
\end{split}
\end{equation}

\begin{equation}
\label{bcheeg3}
\begin{split}
&Xdiv_{i_1}\dots Xdiv_{i_a}
LC^{No\Phi,targ,B,T_i,\gamma}_{\phi_{u+1}}[C^{l,i_1\dots
i_a}_{g}(\Omega_1,\dots ,\Omega_p,\phi_1,\dots ,\phi_u)]=
\\& -(|I_1|-1)Xdiv_{i_2}\dots
Xdiv_{i_\alpha}\hat{C}^{l,i_1\dots i_a}_{g}(\Omega_1,\dots
,\Omega_p,\phi_1,\dots ,\phi_u)\nabla_{i_1}\phi_{u+1}+
\\&\Sum_{h\in H} a_h Xdiv_{i_1}\dots Xdiv_{i_a}C^{h,i_1\dots
i_\mu,i_{*}}_{g}(\Omega_1,\dots ,\Omega_p,\phi_1,\dots
,\phi_u)\nabla_{i_{*}}\phi_{u+1} +
\\&\Sum_{q\in Q} a_q
C^q_{g}(\Omega_1,\dots ,\Omega_p,\phi_1,\dots ,\phi_{u+1}),
\end{split}
\end{equation}

\begin{equation}
\label{bcheeg4}
\begin{split}
&Xdiv_{i_1}\dots Xdiv_{i_a}
LC^{No\Phi,targ,B,T_i,\delta}_{\phi_{u+1}}[C^{l,i_1\dots
i_a}_{g}(\Omega_1,\dots ,\Omega_p,\phi_1,\dots ,\phi_u)]=
\\&\Sum_{q\in Q} a_q C^q_{g}(\Omega_1,\dots
,\Omega_p,\phi_1,\dots ,\phi_{u+1}),
\end{split}
\end{equation}

\begin{equation}
\label{bcheeg5}
\begin{split}
&Xdiv_{i_1}\dots Xdiv_{i_a}
LC^{No\Phi,targ,B,T_i,\varepsilon}_{\phi_{u+1}}[C^{l,i_1\dots
i_a}_{g}(\Omega_1,\dots ,\Omega_p,\phi_1,\dots ,\phi_u)]=
\\&\Sum_{h\in H} a_h Xdiv_{i_1}\dots Xdiv_{i_a}C^{h,i_1\dots
i_\mu,i_{*}}_{g}(\Omega_1,\dots ,\Omega_p,\phi_1,\dots
,\phi_u)\nabla_{i_{*}}\phi_{u+1} +
\\&\Sum_{q\in Q} a_q
C^q_{g}(\Omega_1,\dots ,\Omega_p,\phi_1,\dots ,\phi_{u+1}).
\end{split}
\end{equation}

\par Adding all the above we obtain our conclusion in the case
 where the crucial factor $T_i$ has precisely one special free index.

\par We now consider the case where the are two special free
 indices in $T_i$. We have assumed that these two special free indices
 are ${}_{i_1}, {}_{i_2}$ in the crucial factor $T_i=\nabla^{(m)}R_{i_1ji_2l}$.
 We will moreover
slightly alter our notational conventions:
 Now, we will still speak of the non-special free
 indices in the crucial factor, but they will in fact be
${}_{i_3},\dots ,{}_{i_{|I_1|}}$. Moreover, in the above sublinear
 combinations when we referred to the indices ${}_{i_1}$ or
$\nabla^{i_1}$ we will now read ``one of the indices
${}_{i_1},{}_{i_2}$ or $\nabla^{i_1},\nabla^{i_2}$''. Lastly, we
define $LC^{No\Phi,div,I_1,T_i,\varepsilon}$
 to stand for the sublinear
 combination that arises by applying the transformation law
 (\ref{levicivita}) to the pair of divergence indices
$\nabla^{i_1},\nabla^{i_2}$ when they have hit the same
 factor.

We now want to describe our first sublinear combination.
 To do so, we need just a little more notation. We denote by
$\tilde{C}^{l,i_1\dots i_a}_{g}$ the tensor field that arises by
switching the indices ${}_{i_1}={}_{i},{}_l$ and by $\tilde{C'}^{l,i_1\dots
i_a}_{g}$ the tensor field that
 arises by switching the indices ${}_{i_2}={}_k,{}_j$. We calculate:

\begin{equation}
\label{vhef1} \begin{split} &LC^{No\Phi,
div,I_1,T_i,\alpha}_{\phi_{u+1}}[Xdiv_{i_1}\dots
Xdiv_{i_a}C^{l,i_1\dots i_a}_{g}(\Omega_1,\dots
,\Omega_p,\phi_1,\dots ,\phi_u)]=
\\&-(\gamma-1) Xdiv_{i_2}\dots
Xdiv_{i_a}C^{l,i_1\dots i_a}_{g}(\Omega_1,\dots
,\Omega_p,\phi_1,\dots ,\phi_u)\nabla_{i_1}\phi_{u+1}
\\&-(\gamma-1) Xdiv_{i_1}Xdiv_{i_3}
\dots Xdiv_{i_a}C^{l,i_1\dots i_a}_{g} (\Omega_1,\dots
,\Omega_p,\phi_1,\dots ,\phi_u)\nabla_{i_2}\phi_{u+1}
\\&-Xdiv_{i_2}\dots
Xdiv_{i_a}\tilde{C}^{l,i_1\dots i_a}_{g}(\Omega_1,\dots
,\Omega_p,\phi_1,\dots ,\phi_u)\nabla_{i_1}\phi_{u+1}
\\&-Xdiv_{i_1}Xdiv_{i_3}\dots Xdiv_{i_a}
\tilde{C'}^{l,i_1\dots i_a}_{g} (\Omega_1,\dots
,\Omega_p,\phi_1,\dots ,\phi_u)\nabla_{i_2}\phi_{u+1}+
\\&\Sum_{q\in Q} a_q
C^q_{g}(\Omega_1,\dots ,\Omega_p,\phi_1,\dots ,\phi_{u+1}).
\end{split}
\end{equation}
Here we note that if $a=\mu$ then the $(\mu-1)$-tensor fields
 $\tilde{C}^{l,i_1\dots i_a}_{g}$ and
$\tilde{C'}^{l,i_1\dots i_a}_{g}$ will be {\it doubly
 subsequent} to the maximal refined double characters
$\vec{L'^z}, z\in Z'_{Max}$.

The second sublinear combination is a little more complicated:

\begin{equation}
\label{vhef2} \begin{split} &LC^{No\Phi,
div,I_1,\beta}_{\phi_{u+1}}[Xdiv_{i_1}\dots
Xdiv_{i_a}C^{l,i_1\dots i_a}_{g}(\Omega_1,\dots
,\Omega_p,\phi_1,\dots ,\phi_u)]=
\\&-(|I_1|-2) Xdiv_{i_2}\dots Xdiv_{i_a}C^{l,i_1\dots
i_a}_{g}(\Omega_1,\dots ,\Omega_p,\phi_1,\dots
,\phi_u)\nabla_{i_1}\phi_{u+1}
\\&-(|I_1|-2) Xdiv_{i_1}Xdiv_{i_3}\dots
Xdiv_{i_a}C^{l,i_1\dots i_a}_{g}(\Omega_1,\dots
,\Omega_p,\phi_1,\dots ,\phi_u)\nabla_{i_2}\phi_{u+1}
\\& +\Sum_{q\in Q} a_q
C^q_{g}(\Omega_1,\dots ,\Omega_p,\phi_1,\dots ,\phi_{u+1}).
\end{split}
\end{equation}

\par On the other hand, we also see that:

\begin{equation}
\label{vhef3} \begin{split} &LC^{No\Phi,
div,I_1,T_i,\gamma}_{\phi_{u+1}}[Xdiv_{i_1}\dots
Xdiv_{i_a}C^{l,i_1\dots i_a}_{g}(\Omega_1,\dots
,\Omega_p,\phi_1,\dots ,\phi_u)]=
\\&\Sum_{q\in Q} a_q
C^q_{g}(\Omega_1,\dots ,\Omega_p,\phi_1,\dots ,\phi_{u+1}).
\end{split}
\end{equation}

\par To describe the fourth sublinear combination, we
introduce some notation: For each $k, 2\le k\le |I_1|$ we define
$\hat{C}^{l,i_1\dots i_a}_{g}$ to stand for the sublinear
combination which arises from $C^{l,i_1\dots i_a}_{g}$ by
performing a cyclic permutation of the indices
${}_{i_k},{}_{i_1},{}_j$ in the crucial factor
$\nabla^{(m)}_{r_1\dots r_m} R_{ijkl}$. We also define
$\hat{C'}^{l,i_1\dots i_a}_{g}$ to stand for the sublinear
combination which arises from $C^{l,i_1\dots i_a}_{g}$ by
performing a cyclic permutation of the indices
${}_{i_k},{}_{i_2},{}_l$ in the crucial factor
$\nabla^{(m)}_{r_1\dots r_m} R_{ijkl}$. We then derive:

\begin{equation}
\label{vhef4} \begin{split} &LC^{No\Phi,
div,I_1,T_i,\delta}_{\phi_{u+1}}[Xdiv_{i_1}\dots
Xdiv_{i_a}C^{l,i_1\dots i_a}_{g}(\Omega_1,\dots
,\Omega_p,\phi_1,\dots ,\phi_u)]=
\\&-(|I_1|-2) Xdiv_{i_2}\dots Xdiv_{i_\alpha}
\hat{C}^{l,i_1\dots i_a}_{g}(\Omega_1,\dots ,\Omega_p,
\phi_1,\dots,\phi_u)\nabla_{i_1}\phi_{u+1}
\\& -(|I_1|-2) Xdiv_{i_1}Xdiv_{i_3}\dots
Xdiv_{i_\alpha}\hat{C'}^{l,i_1\dots i_a}_{g}(\Omega_1,\dots
,\Omega_p,\phi_1,\dots ,\phi_u)\nabla_{i_2}\phi_{u+1}
\\&+\Sum_{h\in H} a_h Xdiv_{i_1}\dots Xdiv_{i_a}C^{h,i_1\dots
i_\mu,i_{*}}_{g}(\Omega_1,\dots ,\Omega_p,\phi_1,\dots
,\phi_u)\nabla_{i_{*}}\phi_{u+1}
\\& +\Sum_{q\in Q} a_q
C^q_{g}(\Omega_1,\dots ,\Omega_p,\phi_1,\dots ,\phi_{u+1}).
\end{split}
\end{equation}

\par Lastly, in this case we  compute:

\begin{equation}
\label{vhef5} \begin{split} &LC^{No\Phi,
div,I_1,T_i,\varepsilon}_{\phi_{u+1}}[Xdiv_{i_1}\dots
Xdiv_{i_a}C^{l,i_1\dots i_a}_{g}(\Omega_1,\dots
,\Omega_p,\phi_1,\dots ,\phi_u)]=
\\&- Xdiv_{i_2}\dots Xdiv_{i_a}
C^{l,i_1\dots i_a}_{g}(\Omega_1,\dots ,\Omega_p,
\phi_1,\dots,\phi_u)\nabla_{i_1}\phi_{u+1}-
\\& - Xdiv_{i_1}Xdiv_{i_3}\dots
Xdiv_{i_\alpha}C^{l,i_1\dots i_a}_{g}(\Omega_1,\dots
,\Omega_p,\phi_1,\dots ,\phi_u)\nabla_{i_2}\phi_{u+1}.
\end{split}
\end{equation}

\par Now, it is quite straightforward in this case to understand the sublinear combinations
in
$$Xdiv_{i_1}\dots Xdiv_{i_a}
LC^{No\Phi,targ,B,T_i}_{\phi_{u+1}}[C^{l,i_1\dots
i_a}_{g}(\Omega_1,\dots ,\Omega_p,\phi_1,\dots ,\phi_u)].$$
 We calculate:

\begin{equation}
\label{bcheeg1}
\begin{split}
&Xdiv_{i_1}\dots Xdiv_{i_a}
LC^{No\Phi,targ,B,T_i,\alpha}_{\phi_{u+1}}[C^{l,i_1\dots
i_a}_{g}(\Omega_1,\dots ,\Omega_p,\phi_1,\dots ,\phi_u)]=
\\& -\epsilon_iXdiv_{i_2}\dots
Xdiv_{i_\alpha}C^{l,i_1\dots i_a}_{g}(\Omega_1,\dots
,\Omega_p,\phi_1,\dots ,\phi_u)\nabla_{i_1}\phi_{u+1}+
\\&-\epsilon_iXdiv_{i_1}Xdiv_{i_3}\dots
Xdiv_{i_\alpha}C^{l,i_1\dots i_a}_{g}(\Omega_1,\dots
,\Omega_p,\phi_1,\dots ,\phi_u)\nabla_{i_2}\phi_{u+1}+
\\&\Sum_{h\in H} a_h Xdiv_{i_1}\dots Xdiv_{i_a}C^{h,i_1\dots
i_\mu,i_{*}}_{g}(\Omega_1,\dots ,\Omega_p,\phi_1,\dots
,\phi_u)\nabla_{i_{*}}\phi_{u+1},
\end{split}
\end{equation}

\begin{equation}
\label{bcheeg2}
\begin{split}
&Xdiv_{i_1}\dots Xdiv_{i_a}
LC^{No\Phi,targ,B,T_i,\beta}_{\phi_{u+1}}[C^{l,i_1\dots
i_a}_{g}(\Omega_1,\dots ,\Omega_p,\phi_1,\dots ,\phi_u)]=
\\& -(|I_1|-2)Xdiv_{i_2}\dots
Xdiv_{i_\alpha}C^{l,i_1\dots i_a}_{g}(\Omega_1,\dots
,\Omega_p,\phi_1,\dots ,\phi_u)\nabla_{i_1}\phi_{u+1}
\\& -(|I_1|-2)Xdiv_{i_1}Xdiv_{i_3}\dots
Xdiv_{i_\alpha}C^{l,i_1\dots i_a}_{g}(\Omega_1,\dots
,\Omega_p,\phi_1,\dots ,\phi_u)\nabla_{i_2}\phi_{u+1}+
\\&\Sum_{t\in T} a_t Xdiv_{i_1}\dots Xdiv_{i_\mu}C^{h,i_1\dots
i_{a-1},i_{*}}_{g}(\Omega_1,\dots ,\Omega_p,\phi_1,\dots
,\phi_u)\nabla_{i_{*}}\phi_{u+1}+
\\&\Sum_{q\in Q} a_q
C^q_{g}(\Omega_1,\dots ,\Omega_p,\phi_1,\dots ,\phi_{u+1}),
\end{split}
\end{equation}

\begin{equation}
\label{bcheeg3}
\begin{split}
&Xdiv_{i_1}\dots Xdiv_{i_a}
LC^{No\Phi,targ,B,T_i,\gamma}_{\phi_{u+1}}[C^{l,i_1\dots
i_a}_{g}(\Omega_1,\dots ,\Omega_p,\phi_1,\dots ,\phi_u)]=
\\& -(|I_1|-2)Xdiv_{i_2}\dots
Xdiv_{i_\alpha}\hat{C}^{l,i_1\dots i_a}_{g}(\Omega_1,\dots
,\Omega_p,\phi_1,\dots ,\phi_u)\nabla_{i_1}\phi_{u+1}
\\& -(|I_1|-2)Xdiv_{i_1}Xdiv_{i_3}\dots
Xdiv_{i_\alpha}\hat{C'}^{l,i_1\dots i_a}_{g}(\Omega_1,\dots
,\Omega_p,\phi_1,\dots ,\phi_u)\nabla_{i_2}\phi_{u+1}+
\\&\Sum_{h\in H} a_h Xdiv_{i_1}\dots Xdiv_{i_a}C^{h,i_1\dots
i_\mu,i_{*}}_{g}(\Omega_1,\dots ,\Omega_p,\phi_1,\dots
,\phi_u)\nabla_{i_{*}}\phi_{u+1} +
\\&\Sum_{q\in Q} a_q
C^q_{g}(\Omega_1,\dots ,\Omega_p,\phi_1,\dots ,\phi_{u+1}),
\end{split}
\end{equation}

\begin{equation}
\label{bcheeg4}
\begin{split}
&Xdiv_{i_1}\dots Xdiv_{i_a}
LC^{No\Phi,targ,B,T_i,\delta}_{\phi_{u+1}}[C^{l,i_1\dots
i_a}_{g}(\Omega_1,\dots ,\Omega_p,\phi_1,\dots ,\phi_u)]=
\\&\Sum_{q\in Q} a_q C^q_{g}(\Omega_1,\dots
,\Omega_p,\phi_1,\dots ,\phi_{u+1}),
\end{split}
\end{equation}

\begin{equation}
\label{bcheeg5}
\begin{split}
&Xdiv_{i_1}\dots Xdiv_{i_a}
LC^{No\Phi,targ,B,\epsilon}_{\phi_{u+1}}[C^{l,i_1\dots
i_a}_{g}(\Omega_1,\dots ,\Omega_p,\phi_1,\dots ,\phi_u)]=
\\&\Sum_{h\in H} a_h Xdiv_{i_1}\dots Xdiv_{i_a}
C^{h,i_1\dots i_\mu,i_{*}}_{g}(\Omega_1,\dots ,\Omega_p,
\phi_1,\dots ,\phi_u)\nabla_{i_{*}}\phi_{u+1} +
\\&\Sum_{q\in Q} a_q
C^q_{g}(\Omega_1,\dots ,\Omega_p,\phi_1,\dots ,\phi_{u+1}).
\end{split}
\end{equation}

We thus derive the claim of our Lemma in the case where the
crucial factor has two internal free indices, by adding all the
above equations.
\newline

The last case, where the crucial factor has no internal free
indices follows more easily. We now have that ${}_{i_1}$ is
 not a special free index, so we will now consider all the
 sublinear combinations above where ${}_{i_1}$ is not mentioned,
  an also whenever we mentioned above one of the free indices
  ${}_{i_2},\dots ,{}_{i_{|I_1|}}$ we will now read ``one of the free indices
   ${}_{i_1},\dots ,{}_{i_{|I_1|}}$'' (since the index ${}_{i_1}$ is not
  special now). We then have that all the relevant equations
  will hold, with the exception of (\ref{hef4}),
 (\ref{cheeg3}), which now  become:

\begin{equation}
\label{hef4'} \begin{split} &LC^{No\Phi,
div,I_1,\delta}_{\phi_{u+1}}[Xdiv_{i_1}\dots
Xdiv_{i_a}C^{l,i_1\dots i_a}_{g}(\Omega_1,\dots
,\Omega_p,\phi_1,\dots ,\phi_u)]=
\\& \Sum_{q\in Q} a_q
C^q_{g}(\Omega_1,\dots ,\Omega_p,\phi_1,\dots ,\phi_{u+1}),
\end{split}
\end{equation}

\begin{equation}
\label{cheeg3'}
\begin{split}
&Xdiv_{i_1}\dots Xdiv_{i_a}
LC^{No\Phi,targ,B,\gamma}_{\phi_{u+1}}[C^{l,i_1\dots
i_a}_{g}(\Omega_1,\dots ,\Omega_p,\phi_1,\dots ,\phi_u)]=
\\&\Sum_{h\in H} a_h Xdiv_{i_1}\dots Xdiv_{i_a}C^{h,i_1\dots
i_\mu,i_{*}}_{g}(\Omega_1,\dots ,\Omega_p,\phi_1,\dots
,\phi_u)\nabla_{i_{*}}\phi_{u+1} +
\\&\Sum_{q\in Q} a_q
C^q_{g}(\Omega_1,\dots ,\Omega_p,\phi_1,\dots ,\phi_{u+1}),
\end{split}
\end{equation}
by application of the first or second Bianchi identity. $\Box$
\newline

\par We are now in a position to plug in all the equations from
 this section into (\ref{heidegger2}) and derive Lemma \ref{zetajones}:

{\bf Proof of Lemma \ref{zetajones}:} We use the equations
(\ref{heidegger2}), (\ref{harun1}), (\ref{harun2}) and plug in all
the other equations from this section, and also use
(\ref{keepinminda}). We write out our conclusion concisely:

\begin{equation}
\label{provezetajones}
\begin{split}
&-\Sum_{z\in Z'_{Max}}\Sum_{l\in L^z} a_l
(\gamma+\epsilon-1-2(s-1)-X) Xdiv_{i_2}\dots Xdiv_{i_a}
\\& C^{l,i_1\dots i_a}_{g}(\Omega_1,\dots ,\Omega_p,\phi_1,
\dots,\phi_u)\nabla_{i_1}\phi_{u+1}+
\\&\Sum_{t\in T} a_t
Xdiv_{i_2}\dots Xdiv_{i_\mu} C^{t,i_1\dots
i_\mu}_{g}(\Omega_1,\dots ,\Omega_p,\phi_1,\dots
,\phi_u)\nabla_{i_1}\phi_{u+1}+
\\&\Sum_{h\in H} a_h Xdiv_{i_1}\dots Xdiv_{i_\mu}
C^{h,i_1\dots i_\mu,i_{*}}_{g}(\Omega_1,\dots
,\Omega_p,\phi_1,\dots ,\phi_u)\nabla_{i_{*}}\phi_{u+1}+
\\&\Sum_{q\in Q} a_q
C^q_{g}(\Omega_1,\dots ,\Omega_p, \phi_1,\dots ,\phi_{u+1})=0,
\end{split}
\end{equation}
modulo complete contractions of length $\ge\sigma +u+2$. Here
$$\Sum_{t\in T} a_t
Xdiv_{i_2}\dots Xdiv_{i_\mu} C^{t,i_1\dots
i_\mu}_{g}(\Omega_1,\dots ,\Omega_p,\phi_1,\dots
,\phi_u)\nabla_{i_1}\phi_{u+1}$$ stands for a generic linear
combination of $(\mu-1)$-tensor fields with a $(u+1)$-simple
character $\vec{\kappa}^{+}_{simp}$ but that are doubly subsequent
to each $\vec{L^z}', z\in Z'_{Max}$.

\par Moreover, by the definition of weight we derive
the following elementary identity:

\begin{equation}
\label{elementary}
\begin{split}
&\text{(Total number of indices in each complete contraction in}
\\&Xdiv_{i_1}\dots Xdiv_{i_a} C^{l,i_1\dots
i_a}_{g}(\Omega_1,\dots ,\Omega_p,\phi_1,\dots ,\phi_u))
\\&=n+2s=\gamma+\epsilon +2|\Phi|+2|I_1| +|I_2|+1.
\end{split}
\end{equation}
Thich shows us that the quantity $(\gamma+\epsilon -1-2(s-1)-X)$
is {\it fixed} for each $l\in L^z, z\in Z'_{Max}$; (i.e. the same
in all the terms in the first line in (\ref{provezetajones})).

A counting argument shows that $(\gamma+\epsilon -1-2(s-1)-X)=0$
if and only if 
$\sigma_1=0$ (i.e. there are no factors 
$\nabla^{(m)}R_{ijkl}$ in $\vec{\kappa}_{simp}$), 
and the tensor fields $C^{l,i_1\dots i_\mu}_g$, $l\in L^z$ 
have exactly one ``exceptionnal index''; we define an index ${}_b$ 
in $C^{l,i_1\dots i_g}_g$ (in the form (\ref{form2})
to be exceptional if it belongs to a factor $S*\nabla^{(\nu)}R_{ijkl}$ 
or $\nabla^{(B)}\Omega_h$, is non-special, and moreover:
 if it belongs to the crucial factor $S_{*}\nabla^{(\nu)}R_{ijkl}$ then it must be  
non-free and not contracting against a factor $\nabla\phi_h$;
if it belongs to a non-crucial factor then it must not 
contracting against a factor $\nabla\phi_h$. 
Notice that by weight considerations, if one of the 
tensor fields $C^{l,i_1\dots i_\mu}_g$, $l\in L^z,z\in Z'_{Max}$ 
has exactly one exceptional index, then all of them do.

 We will check that Lemma \ref{zetajones}
indeed holds in this very special case (we will call it the {\it
unfortunate case}) in a Mini-Appendix at the end of this paper.
In all the remaining cases, we derive Lemma \ref{zetajones} by just dividing 
(\ref{provezetajones}) by  $(\gamma+\epsilon -1-2(s-1)-X)$. 
\newline

{\bf Proof of Lemma \ref{pool2}:} We again use equations
(\ref{heidegger}) and (\ref{harun1})
 and replace according to all the equations of this
 subsection, also using (\ref{akinola'}), (\ref{akinola2'}). Then,
we first consider the case where the maximal refined double
characters $\vec{L^z}, z\in Z'_{Max}$ have 2 free indices
${}_i={}_{i_1},{}_k={}_{i_2}$
 in the crucial factor(s), and we denote by
 $M$ the number of crucial factors and by ${}_{i_1},{}_{i_2}, {}_{i_3},{}_{i_4},\dots ,
{}_{i_{2M-1}},{}_{i_{2M}}$ the special free indices that belong
 to those factors (${}_{i_1},{}_{i_2}$ are the indices ${}_i,{}_k$
  in the first crucial factor etc.) we deduce:

\begin{equation}
\label{provepool2a}
\begin{split}
&-\Sum_{z\in Z'_{Max}}\Sum_{l\in L^z} \sum_{i_h\in
I_{l,*}}a_l(\gamma_i+\epsilon_i -1
-2(s-1)-X)Xdiv_{i_1}\dots\hat{Xdiv}_{i_h}\dots Xdiv_{i_\mu}
\\& \tilde{C}^{l,i_1\dots i_a}_{g}(\Omega_1,\dots ,\Omega_p,\phi_1,
\dots,\phi_u)\nabla_{i_1}\phi_{u+1}
\\&+\Sum_{t\in T} a_t
Xdiv_{i_2}\dots Xdiv_{i_\mu} C^{t,i_1\dots
i_\mu}_{g}(\Omega_1,\dots ,\Omega_p,\phi_1,\dots
,\phi_u)\nabla_{i_1}\phi_{u+1}+
\\&\Sum_{h\in H} a_h Xdiv_{i_1}\dots Xdiv_{i_\mu}
C^{h,i_1\dots i_\mu,i_{*}}_{g}(\Omega_1,\dots
,\Omega_p,\phi_1,\dots ,\phi_u)\nabla_{i_{*}}\phi_{u+1}+
\\&\Sum_{q\in Q} a_q
C^q_{g}(\Omega_1,\dots ,\Omega_p, \phi_1,\dots ,\phi_{u+1})=0,
\end{split}
\end{equation}
modulo complete contractions of length $\ge\sigma +u+2$. Here
again $$\Sum_{t\in T} a_t C^{t,i_1\dots i_\mu}_{g}(\Omega_1,\dots
,\Omega_p,\phi_1, \dots,\phi_u)\nabla_{i_1}\phi_{u+1}$$ is a
generic linear combination of $(\mu+1)$-tensor fields with
$(u+1)$-simple character $\vec{\kappa}^{+}_{simp}$ and a refined
double
 character that is subsequent to each $\vec{L^z}',
 z\in Z'_{Max}$.
Moreover, the same elementary observation as above continues
 to hold, hence we deduce that $(\gamma_i+\epsilon_i
-2(s-1)-X)$ is independent of $l\in L^z$ (and of the choice of
crucial factor), and also that it is non-zero 
(by a counting argument again--we use the fact that 
$\gamma\ge 2$). Thus we derive Lemma \ref{pool2} by dividing by  this constant.

\par Finally, consider the case where the maximal refined double
characters $\vec{L^z}, z\in  Z'_{Max}$ have one special free index
in the crucial factor(s). Recall that $I_{l,*}$ stands for the
index set of the special free indices that belong to the (one of
the) crucial factor(s). We deduce:

\begin{equation}
\label{provepool2b}
\begin{split}
&-\Sum_{z\in Z'_{Max}} \Sum_{l\in L^z}\sum_{i_h\in I_{l,*}} a_l
(\gamma+\epsilon-2(s-1)-X)Xdiv_{i_1}\dots \hat{Xdiv}_{i_h}\dots
Xdiv_{i_\mu}
\\&\tilde{C}^{l,i_1\dots i_\mu}_{g}(\Omega_1,\dots ,\Omega_p,\phi_1,\dots
,\phi_u)\nabla_{i_1}\phi_{u+1}+
\\&\Sum_{t\in T} a_t
Xdiv_{i_2}\dots Xdiv_{i_a} C^{t,i_1\dots i_a}_{g}(\Omega_1,\dots
,\Omega_p,\phi_1,\dots ,\phi_u)\nabla_{i_1}\phi_{u+1}+
\\&\Sum_{h\in H} a_h Xdiv_{i_1}\dots Xdiv_{i_\mu}
C^{h,i_1\dots i_\mu,i_{*}}_{g}(\Omega_1,\dots
,\Omega_p,\phi_1,\dots ,\phi_u)\nabla_{i_{*}}\phi_{u+1}+
\\&\Sum_{q\in Q} a_q
C^q_{g}(\Omega_1,\dots ,\Omega_p, \phi_1,\dots ,\phi_{u+1})+
\Sum_{z\in Z} a_z C^z_{g}(\Omega_1,\dots ,\Omega_p, \phi_1,\dots
,\phi_{u+1})=0,
\end{split}
\end{equation}
modulo complete contractions of length $\ge\sigma +u+2$. Moreover,
the same elementary observation as above continue to hold, hence
we deduce that $(\gamma_i+\epsilon_i + -2(s-1)-X)$ is
independent of $l$ and also that it is non-zero (since in this
case $\gamma\ge 3$). This shows Lemma \ref{pool2} in this case also. $\Box$

\subsection{Mini-Appendix: Proof of Lemma
\ref{zetajones} in the unfortunate
case.}

\par In order to show Lemma \ref{zetajones} in this setting, we
 recall that here $\sigma_1=0$ and
$X=0$ (hence all
factors $\nabla^{(B)}\Omega_x$ must be contracting against some
factor $\nabla\phi'_h$). 
We also observe that by weight considerations all tensor fields
in (\ref{hypothese2}) {\it other than} 
the ones indexed in $L^z$, for a given $z\in Z'_{Max}$  
 can have at most $M-1$ free indices belonging to the 
crucial factor $S_{*}\nabla^{(\nu)}R_{ijkl}\nabla^i\tilde{\phi}_1$. 

We make a further observation: If for some $C^{l,i_1\dots i_\mu}_g$ 
the exceptionnal index belongs to a factor $\nabla^{(2)}\Omega_h$ 
(with exactly two derivatives), or to a simple factor $S_{*}R_{ijkl}$,
then all tensor fields $C^{l,i_1\dots i_\mu}_g$ 
must have that property; this follows by weight considerations. We call this subcase A. 
The other case we call subcase B. 

Let us prove Lemma \ref{zetajones} in subcase B, which is the
 hardest: We observe that by explicitly constructing 
divergences of $(\mu+1)$-tensor fields (indexed in $H$ below), 
as allowed in Lemma \ref{zetajones},
 we can write:

\begin{equation}
\label{cotton.mather} 
\begin{split}
&\sum_{l\in L^z} a_l C^{l,i_1\dots i_\mu}_g=
Xdiv_{i_{\mu+1}}\sum_{h\in H} a_h C^{h,i_1\dots i_{\mu+1}}_g+
\\&\sum_{l\in \tilde{L}^z} a_l C^{l,i_1\dots i_\mu}_g+
\sum_{j\in J} a_j C^{j,i_1\dots i_\mu}_g.
\end{split}
\end{equation}
Here the terms indexed in $J$ are simply subsequent 
to $\vec{\kappa}_{simp}$. The terms indexed in $\tilde{L}^z$ have all 
the features of the terms in $L^z$ in the LHS (in 
particular the have the same, maximal, refined double character), and in addition:
\begin{enumerate}
 \item The exceptional index does not belong to the crucial factor.

\item If the index ${}_l$ in the crucial factor contracts against  
 a special index 
in a factor $T=S_*\nabla^{(\nu)}_{r_1\dots r_\nu}R_{ijkl'}$,
then it contracts against the index ${}_k$; moreover, if the exceptional index 
belongs to the factor $T$ (say it is the index ${}_{r_\nu}$,  then 
the indices ${}_l,{}_{r_\nu}$ are symmetrized over. 

\item If the index ${}^l$ against which the index ${}_l$ in 
the crucial factor contracts is a non-special index (and hence we may assume 
wlog that it is a derivative index--denote it by $\nabla^l$) then 
it {\it does not} belong to some specified factor $T'$ in $\vec{\kappa}_{simp}$.
\end{enumerate}

In view of the above, we may assume that all the tensor fields indexed in $L^z$ 
in (\ref{hypothese2}) have the features described above. Let us then 
break the index set $L^z$ into two subsets: We say that $l\in L^z_1$ if and only if 
the index ${}_l$ in the crucial factor contracts against 
a special index. We set $L^z_2=L^z\setminus L^z_1$.

We will then prove:

\begin{equation}
\label{paristanw1}
\sum_{l\in L^z_1} a_l C^{l,(i_1\dots i_\mu)}_g=0,
\end{equation}

\begin{equation}
\label{paristanw2}
\sum_{l\in L^z_2} a_l C^{l,(i_1\dots i_\mu)}_g=0.
\end{equation}
Clearly, if we can proe the above then Lemma \ref{zetajones} will follow in subcase B of 
the unfortunate case.
\newline

{\it Proof of (\ref{paristanw1}):} For future reference, let us
 denote by $L^z_{1,b}\subset L^z_1$ the index set of tensor
  fields for which the index ${}_l$ in the crucial factor contracts against the factor 
 $S_*\nabla^{(\nu')}R_{i'j'k'l'}\nabla^{i'}\tilde{\phi}_b$.
 
The main tool we will use in this proof
will be used for the other claims in this subsection. 
We consider $Image^1_Y[L_g]=0$ (the first conformal variation 
of (\ref{hypothese2}), and we pick out the sublinear combination 
$Image^{1,*}_Y[L_g]$ of terms with length $\sigma+u$, 
with the crucial factor $S_{*}\nabla^{(\nu)}R_{ijkl}$ 
being replaced by a factor $\nabla^{(\nu+2)}Y$, and the factor $\nabla\phi_1$ 
contracting against  some other factor. We derive that:
$$Image^{1,*}_Y[L_g]=0,$$
modulo complete contractions of length $\sigma+u+1$. 
We further break up $Image^{1,*}_Y[L_g]$, into sublinear combinations
$Image^{1,*,b}_Y[L_g]$, where a complete contraction 
belongs to $Image^{1,*,b}_Y[L_g]$ if and only if 
$\nabla\phi_i$ contracts against the
 factor $S_{*}\nabla^{(\nu')}R_{ijkl}\nabla^i\tilde{\phi}_b$. 
Clearly, we also have:
$$Image^{1,*}_Y[L_g]=0,$$
modulo complete contractions of length $\sigma+u+1$. 
Now, in the above equation we formally replace the two factors $\nabla_a\phi_1,\nabla_c\phi_b$ 
by a factor $g_{ac}$. This clearly produces a new true equation. 
We then act on the new true equation by $Ricto\Omega_{p+1}$,\footnote{See 
the Appendix in \cite{alexakis1}.} thus deriving a new true equation:

\begin{equation}
\label{deutsche}
Z_g(\Omega_1,\dots,\Omega_p,\Omega_{p+1},\phi_1,\dots,\phi_u)=0,
\end{equation}
which holds modulo terms of length $\sigma+u+1$. 

This equation can in fact be re-expressed in a  more useful form:
Let us denote by $L_*\subset L_\mu\bigcup L_{>\mu}$ the index set of terms for which 
one special index (say the index ${}_l$ in the 
crucial factor contracts against a special index (say ${}_{l'}$)
in the factor $S_*\nabla^{(\nu')}R_{i'j'k'l'}$.\footnote{Notice that $L^z\subset L_*$.} We denote 
by $\overline{C}^{l,i_1\dots i_a}_g(\dots ,Y,\Omega_{p+1})$
 the tensor field that arises from $C^{l,i_1\dots i_\mu}$
 by replacing the expression 
$S_*\nabla^{(\nu)}_{r_1\dots r_\nu}R_{ijkl}\otimes S_*\nabla^{(\nu')}_{t_1\dots t_{\nu'}}{R_{i'j'k'}}^l\otimes\nabla^i\tilde{\phi}_1\otimes\nabla^{i'}\tilde{\phi}_b$
by $\nabla^{(\nu+2)}_{r_1\dots r_\nu jk}Y\otimes \nabla^{(\nu'+2)}_{t_1\dots t_{\nu'} j'k'}\Omega_{p+1}$.
Denote by  $Cut[\vec{\kappa}_{simp}]$ the simple character of
 the resulting tensor field. Then (\ref{deutsche})  can be re-expressed as:

\begin{equation}
\label{kupriakon} 
\sum_{l\in L_*} a_l Xdiv_{i_1}\dots Xdiv_{i_a} 
C^{l,i_1\dots i_\mu}_g(\dots,Y,\Omega_{p+1})+\sum_{j\in J} a_j C^j=0.
\end{equation}
modulo complete contractions of length $\sigma+u+1$. The terms indexed 
in $J$ are simply subsequent to $Cut[\vec{\kappa}_{simp}]$. We now apply the 
Eraser to all factors $\nabla\phi_h$ that contract 
against the factor $\nabla^{(B)}Y$; by abuse of notation, we still denote the 
resulting tensor fields, complete contractions etc by $C^{l,i_1\dots i_\mu}_g,C^j_g$.

Now, we apply the inverse integration 
by part to the above equation,\footnote{This has been defined at
in section 3 of \cite{alexakis5}.} deriving an integral equation:

\begin{equation}
\label{kupriakon'} 
\int_{M^n}\sum_{l\in L_*} a_l  
\hat{C}^{l}_g(\dots,Y,\Omega_{p+1})+\sum_{j\in J} a_j C^jdV_g=0.
\end{equation}
Here the complete contractions $\hat{C}^{l}_g(\dots,Y,\Omega_{p+1})$ 
arise from the tensor fields 
$C^{l,i_1\dots i_\mu}_g(\dots,Y,\Omega_{p+1})$ by 
making the free indices into internal contractions.\footnote{For future reference, let us write
$\sum_{l\in L^z} a_l \hat{C}^l_g=\Delta^{M}Y\sum_{l\in L^z} a_l C'^l_g$.} 
We then consider the silly divergence formula for the above, 
which arises by integrating by parts wrt. $\nabla^{(B)}Y$; denote the resulting 
equation by $$silly_Y[\sum_{l\in L_*} a_l  
\hat{C}^{l}_g(\dots,Y,\Omega_{p+1})+\sum_{j\in J} a_j C^j].$$ We consider
 the  sublinear combination $silly_{Y,*}[\dots]$ which consists of 
terms with length $\sigma+u$, and $\mu-M$ internal contractions and all factors $\nabla\phi_h$
differentiated only once.
Clearly, $silly_{Y,*}[\dots]=0$. This equation can be re-expressed in the form:  

$$\sum_{l\in L^z_{1,b}} a_l Spread^M_{\nabla,\nabla}[C'^l_g]+\sum_{j\in J} a_j C^j_g=0. $$
Here $Spread_{\nabla,\nabla}$ stands for an operation that hits a pair of 
different factors (not in the form $\nabla\phi$) by derivatives 
$\nabla^a,\nabla_a$ that contract against each 
other and then adding over the resulting terms. From the above, we derive that:
$$\sum_{l\in L^z_{1,b}} a_l C'^l_g+\sum_{j\in J} a_j C^j_g=0.$$
Now, we formally replace all $\mu-M$ internal contractions in the above
by factors $\nabla\upsilon$ (this gives rise to a new true equation), and then 
formally replacing the expression $Y\otimes S_{*}\nabla^{(\nu')}_{t_1\dots t_{\nu'}j'k'}\Omega_{p+1}$ 
by $S_*\nabla^{(\nu)}_{r_1\dots r_\nu}R_{ij(free)l}\otimes
\nabla^i\tilde{\phi}_1\otimes S_*\nabla^{(\nu')}_{t_1\dots t_{\nu'}}
{R_{i'j'k'}}^l\nabla^{r_1}\phi_{x_1}\dots \nabla^{x_z}\phi_z$, where
the indices ${}_{r_1},\dots ,{}_{r_\nu},{}_j$  that are not contracting 
 against a factor $\nabla\phi_x$ are free.\footnote{We add factors $\nabla\phi_x$ 
and free indices to the first factor according to the form of the crucial 
factor in the tensor fields $C^{l,i_1\dots i_\mu}_g$, 
$l\in L^z$ in (\ref{hypothese2}).} 
The resulting true equation is our claim of 
Lemma \ref{zetajones} for the index set  $L^z_{1,b}$. Thus we have derived (\ref{paristanw1}).
\newline

The proof of (\ref{paristanw2}) 
follows by an adaptation of the argument above. We consider the equation 
$Image^1_Y[L_g]=0$ and pick out the sublinear combination $Image^{1,*}_Y[L_g]$ where $\nabla\phi_1$ contracts 
 against the factor $T'$. Clearly 
$Image^{1,*}_Y[L_g]=0$. (Among the tensor fields indexed in $L^z$, 
such terms can only arise from the terms in $L_g=0$ where the $Xdiv_k$ 
for the index ${}_k$ in the crucial factor is forced to hit the factor $T'$;
then the factor $S_{*}\nabla^{(\nu)}_{r_1\dots r_\nu}R_{ijkl}\nabla^i\tilde{\phi}_1$ must be replaced by
 $\nabla^{(\nu+2)}_{r_1\dots r_\nu jl}Y\nabla_k\phi_1$). Now, from $Image^{1,*}_Y[L_g]=0$,
  we again repeat the argument with the inverse integration 
by parts and the silly divergence formula, picking
 out in that equation the sublinear combination 
 with $\mu-M$ internal contractions and all functions $\phi_h$ 
 differentiated only once, with the factor $\nabla\phi_1$ contracting
  against the factor $T'$ and with $M$ particular contractions between 
 the factor $T'$ and another factor $T''$. 
We derive that this sublinear combination, say $silly_+[\dots]$ must vanish 
separately; $silly_+[\dots]=0$. It is also in one-to-one correspondence
with the sublinear combination $\sum_{l\in L^z_2} a_l C^{l,(i_1\dots i_\mu)}_g$, 
by the same reasoning as above. 
We can then reproduce the formal operations as above, and derive 
the claim of Lemma \ref{zetajones} for the index set $L^z_2$. 
\newline

{\it Proof of Lemma \ref{zetajones} in subcase A of the unfortunate case:}  
We again consider the equation 
$Image^1_Y[L_g]=0$ and pick out the sublinear 
combination $Image^{1,*}_Y[L_g]=0$ where $\nabla\phi_1$ contracts 
 against the factor $T'$. This sublinear combination must vanish separately, 
 thus we derive  a new true equation, which we again denote by $Image^{1,*}_Y[L_g]=0$.
Now, we again apply the inverse integration by parts (replacing
 the $Xdiv$'s by internal contractions), deriving an integral equation.
  From this new integral equation, we derive a silly divergence
   formula, by integrating by parts with respect to the factor $\nabla^{(B)}Y$.
We pick out the sublinear combination of terms with 
length $\sigma+u$, $\mu-M$ internal contractions, all functions $\phi_h$
differentiated once, the factor $\nabla\phi_1$ contracting 
against $T'$, where $T'$ contains no  exceptional indices,
 and where there are $M$ exceptional indices on some factor $T''$, 
 contracting against $M$ exceptional indices  in some other factor $T''\ne T'$. Again, this sublinear combination
 vanishes separately, and is in 
one-to-one correspondence
with the sublinear combination $\sum_{l\in L^z} a_l C^{l,(i_1\dots i_\mu)}_g$, 
by the same reasoning as above. 
We can then reproduce the formal operations as above, and derive 
the claim of Lemma \ref{zetajones} for the index set $L^z_2$. $\Box$

\subsection{Mini-Appendix: A postponed claim.}
\label{minappV}

\par We directly derive
Proposition \ref{giade} when the tensor fields of maximal refined
double character in (\ref{hypothese2}) are in the special forms
described at the end of the introduction.

{\it Proof of the postponed claim:} Observe that by weight 
considerations, all tensor fields in (\ref{hypothese2}) will have rank $\mu$. 
Also, none will have special free inices in a factor $S_* R_{ijkl}$. 

 For each
 $l\in L_\mu$ we denote by $C^{l,i_1|A}_g$
the vector field that arises by replacing the factor
 $S_{*}R_{ijkl}\nabla^i\tilde{\phi}_1$
 by $\nabla^{(2)}_{jk}Y\nabla_l\phi_1$; we also denote by
$C^{l,i_1|A}_g$ the vector field that arises by replacing the factor
 $S_{*}R_{ijkl}\nabla^i\tilde{\phi}_1$
 by $-\nabla^{(2)}_{jl}Y\nabla_k\phi_1$. Denote the $(u-1)$-simple 
character of the resulting tensor fields by $\vec{\kappa}'_{simp}$.  
We then derive an equation:

$$\sum_{l\in L_\mu} a_l Xdiv_{i_1}\dots Xdiv_{i_\mu} [C^{l,i_1\dots i_\mu|A}_g+C^{l,i_1\dots i_\mu|B}_g]=0.$$
Now, we apply Lemma 4.10 in \cite{alexakis4} to the above\footnote{There 
are no special free indices in factors $S_*\nabla^{(\nu)}R_{ijkl}$ 
here, hence no danger of ``forbidden cases''.} and derive:

$$\sum_{l\in L_\mu} a_l [C^{l,i_1\dots i_\mu|A}_g+C^{l,i_1\dots i_\mu|B}_g]
\nabla_{i_1}\upsilon\dots \nabla_{i_\mu}\upsilon=0.$$
Finally, formally replacing the expression $\nabla^{(2)}_{(ab)}Y\nabla_c\phi_1$
by $S_{*}R_{i(ab)c}\nabla^1\tilde{\phi}_1$ we derive our claim. $\Box$

\section{Part B: A proof of Lemas 3.3, 3,4, 3.5 in \cite{alexakis4}}

 \par The rest of this paper is devoted to proving Lemma \ref{pskovb} and also 
Lemmas 3.3, 3.4 in \cite{alexakis4}.
  In this proof we will derive the strongest consequence of the local equation (\ref{hypothese2}),
 which is the common assumption of all three Lemmas above. This consequence,
 a {\it new, less complicated} local equation is called
 the ``grand conclusion''. In case A of Lemma \ref{hypothese2},
 the   ``grand conclusion'' coincides with the claim of Lemma \ref{pskovb}. 
 
  In order to derive Lemma \ref{pskovb} in case B, in 
section \ref{caseB}  we will apply the grand conclusion (and certain 
other equations derived below)  in less 
straightforward ways. 
In the next section we provide more technical details
regarding the derivation of the ``grand conclusion''. 

\subsection{A general discussion regarding the derivation of the ``grand conclusion'':}

 The derivation of the ``grand conclusion'' can be divided into two parts:
 In the first part, we repeat the analysis performed in \cite{alexakis6}: We consider the first conformal variation 
 $Image^1_{\phi_{u+1}}[L_g]=0$ of (\ref{hypothese2}).\footnote{Recall that is general, the first
  conformal variation of a true equation $L_g(\psi_1,\dots ,\psi_t)=0$
   that depends on an auxilliary Riemannian metric $g$ is defined through the formula:
  $Image^1_{\phi}[L_g(\psi_1,\dots, \psi_t)]:=\frac{d}{dt}|_{t=0}L_{e^{2t\phi}}(\psi_1,\dots,\psi_t)$; 
  of course $Image^1_{\phi}[L_g(\psi_1,\dots, \psi_t)]=0$. } 
Since the terms in the LHS of the equation $L_g=0$ can be written out as general
complete contractions of the form (\ref{form1}), 
the first conformal variation can be calculated 
by virtue of the transformation laws of the curvature
 tensor and the Levi-Civita connection under conformal 
changes of the metric, \ref{curvtrans}), (\ref{levicivita}).

 Now as part A, in $Image^{1}_{\phi_{u+1}}[L_g]=0$, we pick out the sublinear combination 
$Image^{1,+}_{\phi_{u+1}}[L_g]$ of terms which either {\it are} or {\it can give rise
 to}\footnote{After application of the curvature identity.} terms 
of the type that appear in the claim of 
Lemma \ref{pskovb}. As we have shown in \cite{alexakis6}, 
this sublinear combination will vanish separately, modulo junk terms which we may disregard:

\begin{equation}
\label{souha}
Image^{1,+}_{\phi_{u+1}}[L_g]+\sum_{z\in Z} a_z C^z_g(\Omega_1,\dots,\Omega_p,\phi_1,\dots,\phi_{u+1})=0.
\end{equation}
In section \ref{first.half} 
 we repeat the analysis of (\ref{souha}) that we performed in \cite{alexakis6}. 
(This analysis involves much calculation, but also the appropriate application of the
 inductive assumption of Proposition \ref{giade}).
The resulting equation, however, completely fails to give us the desired conclusion of 
Lemma \ref{pskovb}.\footnote{The defficiencies 
of the analysis of equation (\ref{souha}) are explained at the end of
 subsection \ref{a.study}.} 

In section \ref{uphrxeden}, we seek to analyze a {\it second} 
equation which arises in the conformal variation 
$Image^1_{\phi_{u+1}}[L_g]=0$ of (\ref{hypothese2}), and 
which we had previously discarded:
We pick out the sublinear combination $Image^{1,\beta}_{\phi_{u+1}}[L_g]$
which consists of terms with an {\it internal contraction} 
(see subsection \ref{first.steps} below for more details). Again, this sublinear combination 
vanishes separately, modulo junk terms which we may disregard:
\begin{equation}
\label{francesca}
Image^{1,\beta}_{\phi_{u+1}}[L_g]+\sum_{z\in Z} a_z C^z_g(\Omega_1,\dots,
\Omega_p,\phi_1,\dots,\phi_{u+1})=0.
\end{equation}
In principle,  one might think that this new equation should
 be of no interest, since the claim of Lemma \ref{pskovb} 
involves terms with no internal contractions. However, we are able to suitably analyze
 (\ref{francesca}), to derive  a new equation $Im^{1,\beta}_{\phi_{u+1}}[L_g]=0$
  below (see \ref{antrikos})), which, combined with the analysis of 
  (\ref{souha}) yields the ``grand conclusion'', as follows: 
 
 \par We define a simple formal  operation $Soph\{\dots\}$ which acts on the terms in 
  $Im^{1,\beta}_{\phi_{u+1}}[L_g]=0$ to produce terms of the type that appear in the claim of Lemma
\ref{pskovb}. We then {\it add} the resulting equation to our analysis 
of (\ref{souha}) and observe many miraculous cancellations of ``bad terms''.
The resulting local equations are collectively called the ``grand conclusion''; they are
the equations (\ref{olaxreiazontai1}), (\ref{olaxreiazontai2}), (\ref{olaxreiazontai3}) below. 
\newline

{\bf Recall Language Conventions:} Firstly, we recall that our hypothesis is equation
(\ref{hypothese2}). We recall that in that equation all the tensor
fields have a given simple character $\vec{\kappa}_{simp}$ and are
acceptable and all the complete contractions $C^j_{g}$ have a weak
character $Weak(\vec{\kappa}_{simp})$ (they are not assumed to be
acceptable).

\par We recall that a free index that is of the form
${}_i,{}_j,{}_k,{}_l$ in some factor $\nabla^{(m)}R_{ijkl}$ or ${}_k,{}_l$ in
some factor $S_{*}\nabla^{(\nu)} R_{ijkl}$ is called special.

\par Next, we recall the discussion regarding the ``selected factor''.
Let us firstly re-explain how the various {\it factors} in the various
 complete contractions and tensor fields in (\ref{hypothese2}) can be distinguished:
We recall that for all the tensor fields and complete
contractions appearing in (\ref{hypothese2}) and for each $\nabla\phi_f, 1\le
f\le u$, there will be a unique factor $\nabla^{(m)} R_{ijkl}$,
$S_{*}\nabla^{(\nu)} R_{ijkl}$ or $\nabla^{(B)}\Omega_h$ against which
$\nabla\phi_f$ is contracting. Therefore, for each $\nabla\phi_f,1\le f\le u$
we may unambiguously speak of {\it the} factor against which
$\nabla\phi_f$ is contracting for each of the tensor fields and
contractions in (\ref{hypothese2}). Furthermore, for each $h,1\le
h,\le p$ we may also unambiguously speak of the factor
$\nabla^{(B)}\Omega_h$ in each tensor field and contraction in
(\ref{hypothese2}).

\par On the other hand, we may have factors $\nabla^{(m)} R_{ijkl}$ in
the terms in 
(\ref{hypothese2}) that are not contracting against any factors
$\nabla\phi_h$. We notice that there is the same number of such
factors in each of the tensor fields and the complete contractions
in (\ref{hypothese2}) (since they all have the same weak character
$Weak(\vec{\kappa}_{simp})$). We sometimes refer to such
factors as ``generic factors of the form $\nabla^{(m)}
R_{ijkl}$''.

We recall the two cases that we have distinguished in the setting of Lemma \ref{pskovb}:
Recall that we have denoted by $M$ the number of free
 indices in the critical factor for
the tensor fields $C^{l,i_1\dots i_\mu}_g$, $l\in \bigcup_{z\in
Z'_{Max}} L^z$ (i.e. the tensor fields with the maximal refined 
double character $\vec{L}^z$, for a given $z\in Z'_{Max}$). 
We also denote by $M'(=\alpha)$ the number of free
indices in the second critical factor for the tensor fields
$C^{l,i_1\dots i_\mu}_g, l\in \bigcup_{z\in Z'_{Max}} L^z$. We
recall that case $A$ is when $M'\ge 2$, and case $B$ is when
$M'\le 1$.

{\it The ``special subcase'' of case B:} 
We introduce a ``special subcase'' of case B, in which case 
the derivation of ``grand conclusion'' will be somewhat different.
We say that (\ref{hypothese2}) falls under the ``special subcase'' 
when the tensor fields of maximal refined double character 
in (\ref{hypothese2}) are in the form:

\begin{equation}
\label{tokolpo}
\begin{split} 
&contr(\nabla_{(free)}R_{\sharp\sharp\sharp\sharp}\otimes R_{\sharp\sharp\sharp\sharp}
\otimes\dots\otimes R_{\sharp\sharp\sharp\sharp}\otimes 
\\&S_*R_{ix\sharp\sharp}\otimes\dots
\otimes S_*R_{ix\sharp\sharp}\otimes 
\nabla^{(2)}_{y\sharp}\Omega_1\otimes\dots\otimes\nabla^{(2)}_{y\sharp}\Omega_p\otimes\nabla\phi_1
\otimes\dots\otimes\nabla\phi_u). 
\end{split}
\end{equation}
(In the above, each index ${}_\sharp$ must contract against 
another index in the form ${}_\sharp,{}_x,{}_y$; the indices ${}_x$ are either contracting against 
another index ${}_\sharp,{}_x,{}_y$ or are free, and the indices ${}_y$
 are either contracting against indices ${}_\sharp,{}_x,{}_y$ or are 
free or contract against a factor $\nabla\phi_h$).
In some instances, our 
argument will be modified to treat those ``special subcases''. 
\newline

\par Now, we will be deriving three equations below, (\ref{olaxreiazontai1}), 
(\ref{olaxreiazontai2}), (\ref{olaxreiazontai3}), which will be
collectively  called the ``grand conclusion''. We next explicitly spell out the
hypotheses under which the grand conclusion will be derived in the
list below: The main assumption will be equation
(\ref{hypothese2}), and we also have the extra assumptions stated
in Lemma \ref{pskovb}; we re-iterate that list of assumptions
below. 
\newline

\begin{enumerate}
\item{In both cases A and B, no tensor field of rank $\mu$ can have an internal free
index in any factor $\nabla^{(m)} R_{ijkl}$, nor a free index of the
form ${}_k,{}_l$ in any factor $S_{*}\nabla^{(\nu)} R_{ijkl}$ (this is the main assumption of
Lemma \ref{pskovb}).}

\item{In both cases A and B, there are no $\mu$-tensor fields in (\ref{hypothese2})
with a  free index in the form ${}_j$ in some factor $S_*R_{ijkl}$, 
with no derivatives. THis is a re-statement 
of the assumption $L_\mu^+=\emptyset$.}

\item{In both cases A and B,  no
$(\mu+1)$-tensor field in (\ref{hypothese2}) contains a factor
$S_{*}R_{ijkl}$ with two internal free
indices (this is a re-statement of the assumption $L''_{+}=\emptyset$).}
\end{enumerate}

Now, in order to describe the terms that appear in
the RHSs of our equations below, 
we recal some notational conventions:

\par We introduce the notion of a ``contributor'', adapted to this setting:

\begin{definition}
\label{contributeur2} 

\par In the setting of Lemma \ref{pskovb}, $\sum_{h\in H} a_h C^{h,i_1\dots i_a}_g$ will stand
for a generic linear combination of $a$-tensor fields ($a\ge \mu$)
with a $u$-simple character $\vec{\kappa}_{simp}$, a weak
$(u+1)$-character
$Weak(\vec{\kappa}_{simp}^{+})$\footnote{$\vec{\kappa}_{simp}^{+}$
is some chosen $(u+1)$-simple character where $\nabla\phi_{u+1}$
is {\it not} contracting against a special index.} and
the following additionnal  features: Either the tensor fields above are acceptable, or
they are unacceptable with one unacceptable factor
$\nabla\Omega_x$, which either contracts against
$\nabla\phi_{u+1}$ or does not contract against any
$\nabla\phi_h$. Furthermore, if a tensor field $C^{h,i_1\dots
i_a,i_{*}}_g(\Omega_1,\dots ,\Omega_p,\phi_1,\dots
,\phi_u)\nabla_{i_{*}}\phi_{u+1}$ has $a=\mu$ and one unacceptable
factor $\nabla\Omega_x$ which does not contract against
$\nabla\phi_{u+1}$ then it must have a $(u+1)$-simple character
$\vec{\kappa}^{+}_{simp}$,\footnote{In other words
$\nabla\phi_{u+1}$ is not contracting against a special index.}
and moreover $\nabla\phi_{u+1}$ must be contracting against a
derivative index, and moreover if it is contracting against a
factor $\nabla^{(B)}\Omega_x$ then $B\ge 3$. Moreover, if $a=\mu$
and the factor $\nabla\phi_{u+1}$ is contracting against
$\nabla\Omega_h$ then we require that at least one of the $\mu$
free indices should be non-special, and there should be a
removable index in the tensor field \\$C^{i_1\dots
i_a,i_{*}}_g(\Omega_1,\dots ,\Omega_p,\phi_1,\dots
,\phi_u)\nabla_{i_{*}}\phi_{u+1}$. Finally, if a tensor field
\\$C^{i_1\dots i_a,i_{*}}_g(\Omega_1,\dots ,\Omega_p,\phi_1,\dots
,\phi_u)\nabla_{i_{*}}\phi_{u+1}$ has $a=\mu$ but the factor
$\nabla\phi_{u+1}$ is contracting against a special index then all
$\mu$ free indices must be non-special.

We will be calling the tensor fields in those linear combinations
``contributors''.
\end{definition}

We recall that 
$\sum_{j\in J} a_j C^j_g(\Omega_1,\dots,\Omega_p,\phi_1,\dots,\phi_{u+1})$ stands for
 a generic linear combination of complete contractions
 of length $\sigma+u+1$ in the form (\ref{form1}) which are simple subsequent 
 to the $u$-simple character $\vec{\kappa}_{simp}$.\footnote{Recall that this 
 means that one of the factors $\nabla\phi_h$ which are 
 supposed to contract against an index ${}_i$ in a 
 factor $S_{*}\nabla^{(\nu)}R_{ijkl}$ in the simple character $\vec{\kappa}_{simp}$
  is now contracting against a derivative index in a factor $\nabla^{(m)}R_{ijkl}$.}
\newline

\par Let us introduce one new piece of notation that will be useful further
down. We will define a generic linear combination of $(\mu-1)$-tensor fields
which appears in the grand conclusion only in case B, in the special subcase.

\begin{definition}
\label{generalni}
We denote by
$$\Sum_{b\in B'} a_b C^{b,i_1\dots i_{\mu-1}i_\mu}_{g}
(\Omega_1,\dots ,\Omega_p,\phi_1,\dots
,\phi_u)\nabla_{i_\mu}\phi_{u+1}$$ a generic linear combination of
acceptable $(\mu-1)$-tensor fields with length $\sigma+u+1$, with
a $u$-simple character $\vec{\kappa}_{simp}$, and with the factor
$\nabla\phi_{u+1}$ contracting against a non-special index in the
selected factor and moreover each of the $(\mu-1)$ free indices 
belong to a different factor. 

Moreover, in the special subcase when in addition $\mu=1$, we additionaly 
require that the index ${}_{i_\mu}$ should be 
a derivative index, and moreover if it belongs to  
a factor $\nabla^{(B)}\Omega_h$ then $B\ge 3$.  Furthermore,
we require that if we change the selected factor from $T_a$ to $T_b$, 
then $\Sum_{b\in B'} a_b C^{b,i_1}_{g}
(\Omega_1,\dots ,\Omega_p,\phi_1,\dots
,\phi_u)$ changes by erasing the index $\nabla_{i_1}$ 
from $T_a$ and adding a  derivative $\nabla_{i_1}$ onto $T_b$. 
\end{definition}

\section{The first half of the ``grand conclusion''.}
\label{first.half}

 As explained, our first step 
in deriving the ``grand conclusion'' is 
 to repeat the analysis of the equation (\ref{souha}) from \cite{alexakis6}. We refer the reader 
 to  subsection 2.2 in \cite{alexakis6} for the strict definition 
of the sublinear combination $Image^{1,+}_{\phi_{u+1}}[L_g]$.
 {\it We recall that the sublinear
 combination $Image^{1,+}_{\phi_{u+1}}[L_g]$ in $Image^{1}_{\phi_{u+1}}[L_g]$
 is defined once we have picked a (set of) selected factor(s)
in $\vec{\kappa}_{simp}$.} We also recall that the equation  (\ref{souha})
 has been proven to hold (modulo 
complete contractions of length $\ge\sigma+u+2$).\footnote{We recall 
that the ``junk terms'' $\sum_{z\in A} a_z C^z_g$ have length $\sigma+u+1$ 
and a factor $\nabla^{(A)}\phi_{u+1}$ with $A\ge 2$.}
 
 We also recall (from the end of subsection 2.2 in \cite{alexakis6}) 
 the natural break-up of the sublinear combination $Image^{1,+}_{\phi_{u+1}}[L_g]$
 into three ``pieces'' according to which {\it rule} of conformal variation a given term has arisen from:
\begin{equation}
\label{heidegger}
\begin{split}
&Image^{1,+}_{\phi_{u+1}}[L_{g}(\Omega_1,\dots
,\Omega_p,\phi_1,\dots ,\phi_u)]= CurvTrans[L_{g}(\Omega_1,\dots
,\Omega_p,\phi_1,\dots ,\phi_u)]
\\& +LC[L_{g}(\Omega_1,\dots ,\Omega_p,\phi_1,\dots ,
\phi_u)]+W[L_{g}(\Omega_1,\dots ,\Omega_p,\phi_1,\dots , \phi_u)].
\end{split}
\end{equation}
(The sublinear combination $CurvTrans[L_g\dots]$ consists of terms with length $\sigma+u$ 
and a factor $\nabla^{(A)}\phi_{u+1}, A\ge 2$, while the sublinear combinations
$LC[L_g], W[L_g]$ consist of terms with length $\sigma+u+1$ and a factor $\nabla\phi_{u+1}$).
We then recall the study of the sublinear combination $CurvTrans[L_g\dots]$ 
that we performed in \cite{alexakis6}:\footnote{As noted there, the analysis
 also applies in the setting of Lemma \ref{pskovb}.}

We recall that the sublinear combination $CurvTrans[L_g]$
will be {\it zero, by definition} if the selected factor is
of the form $\nabla^{(A)}\Omega_h$. We also recall that
if the selected factor is of the form $\nabla^{(m)}R_{ijkl}$ or
 $S_{*}\nabla^{(\nu)}R_{ijkl}$, then we have proven that the sublinear combination
$CurvTrans[L_g]$ in (\ref{heidegger}) can be expressed as in equations
(3.31) and  (3.47)   in \cite{alexakis6}, 
 respectively;\footnote{We re-express the
conclusions of equations (3.31) and (3.47)
in new notation in the next subsection.}
 {\it and moreover we proved in \cite{alexakis6} that these equations also hold
in the setting of Lemma \ref{pskovb}}.  We also recall that
 the terms inside parentheses in (3.31) and  (3.47) are defined to be
{\it zero} in the setting of Lemma \ref{pskovb}. 
The result of this analysis in \cite{alexakis6}, as a new local 
equation, namely equation (4.1) in that paper, which we reproduce here:

\begin{equation}
\label{heidegger2} CurvTrans^{study}[L_g]+LC[L_g]+W[L_g]=0;
\end{equation}
this holds modulo complete contractions of length $\ge\sigma+u+2$.

The purpose of the first half of this paper will be to study the
sublinear combinations $LC_{\phi_{u+1}}[\dots],W_{\phi_{u+1}}[\dots]$ in
equation (\ref{heidegger2}).
\newline

\par Now,  our aim in this section  will be to study
 the sublinear combinations $LC[L_g]$ and $W[L_g]$
in the context of Lemma \ref{pskovb} from \cite{alexakis6}. We recall the
Lemma 2.2 in \cite{alexakis6}, which  has also been proven in the setting of Lemma \ref{pskovb}.
In view of this, we only need to study the sublinear
combinations $LC^{No\Phi}[L_g]$ and $W[L_g]$.
\newline

Our aim here is to  recall our description of the term  $CurvTrans^{study}[L_g]$ appearing in
 (\ref{heidegger2}),\footnote{We recall that this analysis was performed in \cite{alexakis6}.} and then 
 the understand the sublinear combinations \\$LC^{No\Phi}[L_g], W[L_g]$ in (\ref{heidegger}), in
 the setting of Lemma \ref{pskovb}. 

\subsection{Proof of Lemma \ref{pskovb}: A description of
 the sublinear
 combination $CurvTrans[L_{g}]$ in this setting.}

\par In this subsection we just reproduce the 
equations (\ref{akinola3'}) and (\ref{keepinminda}) from part A, 
for the reader's convenience. These equations describe the sublinear 
combination $CurvTrans^{study}[L_g]$ appearing in (\ref{heidegger2}).
Recall that (3.31) and  (3.47)
in \cite{alexakis6}) correspond to the cases where the selected
factor(s) is (are) in the form $S_{*}\nabla^{(\nu)}R_{ijkl}$,
$\nabla^{(m)}R_{ijkl}$.  
We recall that if the selected factor is
in the form $\nabla^{(A)}\Omega_h$ {\it then by definition
$CurvTrans[L_g]=CurvTrans^{study}[L_g]=0$}.

{\it Notation:} If the selected factor is of the form
$S_{*}\nabla^{(\nu)}R_{ijkl}$, then for each $l\in L$ we denote by
$I_1^l$ the index set of free indices that belong to the selected
factor $T_l$, where $T_l =S_{*}\nabla^{(\nu_l)}_{r_1\dots
r_\nu}R_{ijkl}$ (note there are $\nu_l$ derivatives on the
selected factor). By our Lemma hypothesis (the assumption 1 in the
list of the previous subsection),  each of the
free indices ${}_i\in I_1^l$ will be one of the indices
${}_{r_1},\dots, ,{}_{r_\nu},{}_j$ in the factor $T_l$ and $\nu_l>0$.
Then, equation (3.47) in \cite{alexakis6} can be re-expressed as:

\begin{equation}
\label{keepinminda'} \begin{split}
 &\Sum_{l\in L} a_l
CurvTrans[Xdiv_{i_1}\dots Xdiv_{i_a}C^{l,i_1\dots
i_a}_{g}(\Omega_1\dots ,\Omega_p,\phi_1,\dots ,\phi_u)]+
\\& \Sum_{j\in J} a_j CurvTrans[C^j_{g}
(\Omega_1\dots ,\Omega_p,\phi_1,\dots ,\phi_u)]=
\\&\Sum_{l\in L_\mu} a_l \frac{1}{\nu_l}\sum_{i_s\in I^l_1}
Xdiv_{i_1}\dots \hat{Xdiv}_{i_s}\dots Xdiv_{i_\mu}C^{l,i_1\dots
i_\mu}_{g} (\Omega_1\dots ,\Omega_p,\phi_1,\dots
,\phi_u)\nabla_{i_s}\phi_{u+1}
\\&+\Sum_{h\in H} a_h Xdiv_{i_1}\dots Xdiv_{i_a}
C^{h,i_1,\dots i_a,i_{*}}_{g} (\Omega_1,\dots
,\Omega_p,\phi_1,\dots ,\phi_u)\nabla_{i_{*}} \phi_{u+1}+
\\& \Sum_{j\in J} a_j C^j_{g}
(\Omega_1,\dots ,\Omega_p,\phi_1,\dots
,\phi_u,\phi_{u+1})+\Sum_{z\in Z} a_z C^z_{g} (\Omega_1,\dots
,\Omega_p,\phi_1,\dots ,\phi_{u+1}).
\end{split}
\end{equation}
(In fact, in this case the tensor fields indexed in $H$ are also
acceptable--but we will not be using this fact).

\par On the other hand, when the selected factor is of the form
$\nabla^{(m)}R_{ijkl}$, then (3.31) in \cite{alexakis6} can be re-expressed
in the form:

\begin{equation}
\label{akinola3''} \begin{split}
 &\Sum_{l\in L} a_l
CurvTrans[Xdiv_{i_1}\dots Xdiv_{i_a}C^{l,i_1\dots
i_a}_{g}(\Omega_1\dots ,\Omega_p,\phi_1,\dots ,\phi_u)]+
\\& \Sum_{j\in J} a_j CurvTrans[C^j_{g}
(\Omega_1\dots ,\Omega_p,\phi_1,\dots ,\phi_u)]=
\\&\Sum_{h\in H} a_h Xdiv_{i_1}\dots Xdiv_{i_a}
C^{h,i_1,\dots i_a,i_{*}}_{g} (\Omega_1,\dots
,\Omega_p,\phi_1,\dots ,\phi_u)\nabla_{i_{*}} \phi_{u+1}+
\\& \Sum_{j\in J} a_j C^j_{g}
(\Omega_1,\dots ,\Omega_p,\phi_1,\dots
,\phi_u,\phi_{u+1})+\Sum_{z\in Z} a_z C^z_{g} (\Omega_1,\dots
,\Omega_p,\phi_1,\dots ,\phi_{u+1}).
\end{split}
\end{equation}
(In this case also the tensor fields indexed in $H$ are
acceptable, but this will not matter).

\subsection{A study of the sublinear combinations
$LC^{No\Phi}[L_{g}]$, $W[L_{g}]$ in the setting of Lemma
\ref{pskovb}.}
\label{a.study}

\par As explained in the beginning of this section, our aim is to find analogues of the
 equations concerning the sublinear combinations\footnote{(which appear in (\ref{heidegger2}))}
$LC^{No\Phi}[L_{g}]$, $W[L_{g}]$  in \cite{alexakis6}. 

\par It is again immediate by the definitions that for each $j\in J$ in
 (\ref{hypothese2}) we must have:

\begin{equation}
\label{daffy1} W[C^j_g(\Omega_1,\dots ,\Omega_p, \phi_1,\dots
,\phi_{u})]=\sum_{j\in J} a_j C^j_g(\Omega_1,\dots ,\Omega_p,
\phi_1,\dots ,\phi_{u+1}),
\end{equation}

\begin{equation}
\label{daffy2} LC^{No\Phi}[C^j_g(\Omega_1,\dots ,\Omega_p,
\phi_1,\dots ,\phi_{u})]=\sum_{j\in J} a_j C^j_g(\Omega_1,\dots
,\Omega_p, \phi_1,\dots ,\phi_{u+1});
\end{equation}
(using generic notation on the right hand side).

\par So the challenge is again
to understand the two sublinear combinations:

$$LC^{No\phi}[Xdiv_{i_1}\dots Xdiv_{i_a}C^{l,i_1\dots i_a}_g(\Omega_1,\dots ,\Omega_p,
\phi_1,\dots ,\phi_{u})],$$

$$W[Xdiv_{i_1}\dots Xdiv_{i_a}C^{l,i_1\dots i_a}_g(\Omega_1,\dots ,\Omega_p,
\phi_1,\dots ,\phi_{u})].$$

\par We again divide the two linear combinations above into the
same linear combinations as in section 
\ref{theanalysis} in part A. Thus,
our aim here is to find analogues of the equations
(\ref{eudokia2}), (\ref{evraz1}), (\ref{avraz1b}), (\ref{evraz2}), (\ref{evraz3}),  
(\ref{eudokia}), 
 (\ref{greatexp}), (\ref{greatexp}) in that paper to find analogues of
 Lemmas \ref{travail} and \ref{travail2} and to find an
 analogue of those Lemmas when the
 selected factor is of the form $\nabla^{(A)}\Omega_s$.

{\bf Important Point:} In the cases where the selected factor is a
curvature term (in the form $\nabla^{(m)}R_{ijkl}$ or
$S_{*}\nabla^{(\nu)}R_{ijkl}$), the major difference with the
setting of Lemmas 3.1 and 3.2 in \cite{alexakis4} (which were proven in \cite{alexakis6}) is the role of
the contractions that belong to the linear combinations
$\sum_{q\in Q}\dots$. Recall that in \cite{alexakis4}
$\sum_{q\in Q} \dots$ stood for a {\it generic} linear combination
of complete contractions with length $\sigma+u+1$ for which the
factor $\nabla\phi_{u+1}$ was contracting against a {\it
non-special}\footnote{Recall that a special
 index is an index of the form ${}_i,{}_j,{}_k,{}_l$ in a factor $\nabla^{(m)}R_{ijkl}$,
or an index of the form ${}_k,{}_l$ in a factor $S_{*}\nabla^{(\nu)}R_{ijkl}$.} index in the selected
factor (which was in one of the forms $\nabla^{(m)}R_{ijkl}$, $S_{*}\nabla^{(\nu)}R_{ijkl}$).

 In the settings of Lemmas \ref{zetajones}, \ref{pool2},
such generic complete contractions were {\it simply subsequent} to the
$(u+1)$-simple character $\vec{\kappa}^{+}_{simp}$ that we were
interested in, and hence they could be {\it disregarded} 
(since they were {\it allowed} in the conclusion of our Lemma).
 In this setting of Lemma \ref{pskovb}, the complete contractions indexed in
$\sum_{q\in Q}\dots$ will have precisely a $(u+1)$-simple character
$\vec{\kappa}^{+}_{simp}$ that we are interested in. Hence they {\it can not} be
disregarded, as they were in part A of this paper.

\par In this sense, our aim for this subsection will be to obtain
 a precise
understanding of the terms belonging to the sublinear combinations
$\sum_{q\in Q} \dots$ in the right hand sides of equations
(4.12), (4.13), (4.14), (4.15), (4.16),  (4.18), 
 (4.20), (4.22) in \cite{alexakis6}
 and in the conclusions of Lemmas 4.1 and 4.2 there.
\newline

{\bf A study of the sublinear combination $W[Xdiv_{i_1}\dots Xdiv_{i_a} C^{l,i_1\dots i_a}_g]$
 in this context: The analogues of equations (\ref{evraz1}), (\ref{evraz1b}), 
(\ref{evraz2}),
(\ref{evraz3}), (\ref{eudokia}), (\ref{greatexp}), (\ref{greatexp}) 
 in the setting of Lemma \ref{pskovb}.}
\newline

\par For each $C^{l,i_1\dots i_a}_{g}$ we will
denote by $\{ T_1,\dots ,T_{b_l}\}$ the set of selected
factors. Also, for each selected factor $T_i$ we will denote
 by $I_1^{T_i}$ the set of free indices that belong to $T_i$
and by $I_2^{T_i}$ the set of free indices that {\it do not} belong
to $T_i$. We recall that by our
 Lemma hypotheses, if a factor $\nabla^{(m)}R_{ijkl}$ in
$C^{l,i_1\dots
 i_\mu}_{g}$, $l\in L_\mu$ has free indices then those free indices {\it must}
 be derivative indices. Also, if a factor $S_{*}\nabla^{(\nu)}_{r_1\dots r_\nu}R_{ijkl}$
 has free indices, they must be of the form ${}_{r_1},\dots ,
{}_{r_\nu},{}_j$. We then consider any $C^{l,i_1\dots i_\mu}_g$, $l\in L_\mu$ 
with a free index in the selected
factor $T_i$; we then calculate the analogue of (\ref{evraz1}), (\ref{evraz1b}),
(\ref{evraz3}):

 \begin{equation}
 \label{Wsec}
\begin{split}
&Xdiv_{i_1}\dots Xdiv_{i_\mu}W^{targ,T_i}[C^{l,i_1\dots
i_\mu}_{g}(\Omega_1,\dots ,\Omega_p,\phi_1,\dots ,\phi_u)]+
\\&W^{targ,T_i,div}[Xdiv_{i_1}\dots Xdiv_{i_\mu}C^{l,i_1\dots
i_\mu}_{g}(\Omega_1,\dots ,\Omega_p,\phi_1,\dots ,\phi_u)]=
\\&2\sigma^{*}\Sum_{i_h\in I^{T_i}_1} Xdiv_{i_1}\dots
\hat{Xdiv}_{i_h}\dots Xdiv_{i_\mu}C^{l,i_1\dots
i_\mu}_{g}(\Omega_1,\dots ,\Omega_p,\phi_1,\dots
,\phi_u)\nabla_{i_h}\phi_{u+1}
\\&+\Sum_{h\in H} a_h Xdiv_{i_1}\dots Xdiv_{i_\mu}C^{h,i_1\dots
i_\mu,i_{*}}_{g}(\Omega_1,\dots ,\Omega_p,\phi_1,\dots
,\phi_u)\nabla_{i_{*}}\phi_{u+1}
\\&+\Sum_{z\in Z} a_z C^z_{g}(\Omega_1,\dots ,\Omega_p,
\phi_1,\dots ,\phi_{u+1});
 \end{split}
 \end{equation}
 here $\sigma^{*}$ stands for $\sigma_1+\sigma_2$ if the crucial
 factor is of the form $\nabla^{(p)}\Omega_h$ and for
 $(\sigma_1+\sigma_2-1)$ if it is in any other form.
 
\par On the other hand,  if $C^{l,i_1\dots i_a}_{g}$ has $a=\mu$ and
does not have free indices in the selected factor $T_i$, or
 if $a>\mu$, we derive the analogue of (\ref{evraz2}):

 \begin{equation}
 \label{Wsec'}
\begin{split}
&Xdiv_{i_1}\dots Xdiv_{i_a}W^{targ,T_i}[C^{l,i_1\dots
i_a}_{g}(\Omega_1,\dots ,\Omega_p,\phi_1,\dots ,\phi_u)]+
\\&W^{targ,T_i,div}[Xdiv_{i_1}\dots Xdiv_{i_a}C^{l,i_1\dots
i_a}_{g}(\Omega_1,\dots ,\Omega_p,\phi_1,\dots ,\phi_u)]=
\\&\Sum_{h\in H} a_h Xdiv_{i_1}\dots Xdiv_{i_\mu}C^{h,i_1\dots
i_\mu,i_{*}}_{g}(\Omega_1,\dots ,\Omega_p,\phi_1,\dots
,\phi_u)\nabla_{i_{*}}\phi_{u+1}
\\&+\Sum_{z\in Z} a_z C^z_{g}(\Omega_1,\dots ,\Omega_p,
\phi_1,\dots ,\phi_{u+1}).
 \end{split}
 \end{equation}

\par The first hard case is to find the  analogue of equation (\ref{greatexp}).
Some notation will prove useful to do this: 
\begin{definition}
\label{pairs}
We define
$(I_2^{T_i})^{2,dif}$ to stand for the set of pairs of indices
$({}_{i_k},{}_{i_l})$, where ${}_{i_k},{}_{i_l}\in I^{T_i}_2$, for
which ${}_{i_k}$ and ${}_{i_l}$ {\it do not} belong to the same
factor. Then, for each pair $({}_{i_k},{}_{i_l})\in (I^{T_i}_2)^{2,dif}$, we
define
$$[C^{l,i_1\dots i_a,i_{*}|T_i}_{g}
(\Omega_1,\dots ,\Omega_p,\phi_1,\dots ,\phi_u) g^{i_k
i_l}]\nabla_{i_{*}}\phi_{u+1}$$ to stand for the $(a-2)$-tensor
field that arises from $C^{l,i_1\dots i_a}_{g}(\Omega_1,\dots
,\Omega_p,\phi_1,\dots ,\phi_u)$ by contracting the indices
${}_{i_k},{}_{i_l}$ against each other and adding a derivative
index $\nabla_{i_{*}}$ on the selected
 factor $T_i$, and then contracting it against a factor
$\nabla_{i_{*}}\phi_{u+1}$ (and also if the selected factor is
 of the form $S_{*}\nabla^{(\nu)} R_{ijkl}$ performing an extra
 $S_{*}$-symmetrization).
\end{definition}

\par It then follows that the sublinear combination
\\$\Sum_{q\in Q} a_q C^q_{g}(\Omega_1,\dots ,\Omega_p, \phi_1,\dots
,\phi_{u+1})$ in (\ref{greatexp}) is precisely of the form:

\begin{equation}
\label{steven}
\begin{split}
& -\Sum_{i=1}^{b_l} \Sum_{(i_k,i_l)\in I^{T_i,2,def}_2}
Xdiv_{i_1}\dots \hat{Xdiv}_{i_k}\dots \hat{Xdiv}_{i_l}\dots
Xdiv_{i_a}
\\&[C^{l,i_1\dots i_a,i_{*}|T_i}_{g}(\Omega_1,\dots
,\Omega_p,\phi_1,\dots ,\phi_u)g^{i_ki_l}
\nabla_{i_{*}}\phi_{u+1}].
\end{split}
\end{equation}

\par Next, we seek to understand the term
$\Sum_{q\in Q} a_q C^q_{g}(\Omega_1,\dots ,\Omega_p, \phi_1,\dots
,\phi_{u+1})$ in $Xdiv_{i_1}\dots Xdiv_{i_a}\{
LC^{No\Phi,targ,A}_{\phi_{u+1}}[C^{l,i_1\dots
i_a}_{g}(\Omega_1,\dots ,\Omega_p,\phi_1,\dots ,\phi_u)]\}$, 
see (\ref{greatexp2}). With our new notation, we easily observe that:

\begin{equation}
\label{toxazo}
\begin{split}
&Xdiv_{i_1}\dots Xdiv_{i_a}\{
LC^{No\Phi,targ,A}_{\phi_{u+1}}[C^{l,i_1\dots
i_a}_{g}(\Omega_1,\dots ,\Omega_p,\phi_1,\dots ,\phi_u)]\}=
\\&\Sum_{h\in H} a_h Xdiv_{i_1}\dots Xdiv_{i_a} C^{h,i_1\dots
i_a,i_{*}}_{g}(\Omega_1,\dots ,\Omega_p,\phi_1,\dots
,\phi_u)\nabla_{i_{*}}\phi_{u+1}.
\end{split}
\end{equation}

\par A harder task is to understand the sublinear combination
\\$\Sum_{q\in Q} a_q C^q_{g}(\Omega_1,\dots ,\Omega_p, \phi_1,\dots
,\phi_{u+1})$ in (\ref{eudokia}) and (\ref{eudokia2}).
\newline

{\it Analogues of (\ref{eudokia}) and (\ref{eudokia2}):}
\newline

{\it Analogue of (\ref{eudokia}):} We recall
 that (by definition)
 \\$LC^{No\Phi, free}_{\phi_{u+1}}[C^{l,i_1\dots
i_a}_{g}(\Omega_1,\dots ,\Omega_p, \phi_1,\dots ,\phi_u)]$ is a
linear combination that consists of tensor fields in the form:

\begin{equation}
\label{gora} C^{t,i_1\dots \hat{i}_f\dots i_a}_{g}(\Omega_1,\dots
,\Omega_p, \phi_1,\dots ,\phi_u)\nabla_{i_f}\phi_{u+1},
\end{equation}
 where the index ${}_{i_f}$
(which is a free index) belongs
 to the factor $\nabla\phi_{u+1}$.
 We observe (by virtue of the transformation 
law (\ref{levicivita})) that the tensor fields
\\$C^{t,i_1\dots \hat{i}_f\dots i_a}_{g}(\Omega_1,\dots ,\Omega_p,
\phi_1,\dots ,\phi_u)$
 will then be acceptable except possibly for one unacceptable
factor $\nabla\Omega_h$.
For each $l\in L$, we will write out that linear combination
as:

\begin{equation}
\label{tyess} \Sum_{t\in T^l} a_t C^{t,i_1\dots\hat{i}_f\dots
i_a}_{g}(\Omega_1,\dots ,\Omega_p, \phi_1,\dots
,\phi_u)\nabla_{i_{f_t}}\phi_{u+1}.
\end{equation}

\par We then observe that (\ref{eudokia}) again holds,
where the sublinear combination $\Sum_{q\in Q} a_q
C^q_{g}(\Omega_1,\dots ,\Omega_p, \phi_1,\dots ,\phi_{u+1})$ is in
fact of the form:

\begin{equation}
\label{den9elw}\begin{split}& \Sum_{t\in T^l} a_t Xdiv_{i_1}\dots
\hat{Xdiv}_{i_{f_t}}\dots Xdiv_{i_a}\{
\nabla^{i_{f_t}}_{sel}[C^{t,i_1\dots \hat{i}_{f_t}\dots
i_a}_{g}\\&(\Omega_1,\dots ,\Omega_p, \phi_1,\dots
,\phi_u)\nabla_{i_{f_t}}\phi_{u+1}]\};
\end{split}
\end{equation}
here
$$\{\nabla^{i_{f_t}}_{sel}[C^{t,i_1\dots \hat{i}_{f_t}\dots
i_a}_{g}(\Omega_1,\dots ,\Omega_p, \phi_1,\dots
,\phi_u)\nabla_{i_{f_t}}\phi_{u+1}]\}$$ stands for the
$(a-1)$-tensor field that arises from \\$C^{t,i_1\dots
\hat{i}_{f_t}\dots i_a}_{g}(\Omega_1,\dots ,\Omega_p, \phi_1,\dots
,\phi_u)\nabla_{i_{f_t}}\phi_{u+1}$ by hitting the selected factor
(if it is unique) with a derivative index $\nabla^{i_{f_t}}$, or if
there are multiple selected factors $T_i,i=1,\dots ,b_l$  then:

\begin{equation}
\label{tatla}
\begin{split}
&\{\nabla^{i_{f_t}}_{sel}[C^{t,i_1\dots \hat{i}_{f_t}\dots
i_a}_{g}(\Omega_1,\dots ,\Omega_p, \phi_1,\dots
,\phi_u)\nabla_{i_{f_t}}\phi_{u+1}]\}=
\\&\Sum_{i=1}^{b_l}\{\nabla^{i_{f_t}}_{T_i}
[C^{t,i_1\dots \hat{i}_{f_t}\dots i_a}_{g}(\Omega_1,\dots
,\Omega_p, \phi_1,\dots ,\phi_u)\nabla_{i_{f_t}}\phi_{u+1}]\}
\end{split}
\end{equation}
(where $\nabla^{i_{f_t}}_{T_i}$ means that $\nabla^{i_{f_t}}$ is forced to hit the factor $T_i$).

\par We observe,  that if $l\in L\setminus L_\mu$ then
the right hand side in (\ref{den9elw}) is a
 contributor.\footnote{See Definition \ref{contributeur2}} If $l\in L_\mu$, we
 must understand it in more detail:

\par We recall that if an index ${}_{i_y}\in I^{T_i}_2$ belongs to
 a factor $\nabla^{(m)}R_{ijkl}$, then by the hypothesis of Lemma
 \ref{pskovb},
 ${}_{i_y}$ must be a derivative index. If it belongs to a factor
 $S_{*}\nabla^{(\nu)} R_{ijkl}$, ${}_{i_y}$
must be one of the indices ${}_{r_1},\dots ,{}_{r_\nu},{}_j$. For each
${}_{i_y}\in I_2$ that belongs to a factor $\nabla^{(m)}R_{ijkl}$, we denote by
$\sigma(i_y)$ the number of indices in that factor
 that are {\it not} contracting against a factor
$\nabla\phi_w$, {\it minus one}. We also define $C^{l,i_1\dots
\hat{i}_y\dots i_\mu}_{g}(\Omega_1,\dots ,\Omega_p, \phi_1,\dots
,\phi_u)$ to stand for the vector field that arises from
$C^{l,i_1\dots  i_\mu}_{g}$ by formally erasing the derivative index ${}_{i_y}$.

 For each ${}_{i_y}\in I^{T_i}_2$ that
belongs to a factor $\nabla^{(A)}\Omega_h$ we denote by
$\sigma(i_y)$
the number of indices in that factor
 that are {\it not} contracting against a factor
$\nabla\phi_y$, {\it minus one}. We again define $C^{l,i_1\dots
\hat{i}_y\dots i_\mu}_{g}(\Omega_1,\dots ,\Omega_p, \phi_1,\dots
,\phi_u)$ to stand for the tensor field that arises from
$C^{l,i_1\dots i_\mu}_{g}$ by erasing the derivative index ${}_{i_y}$.

\par For each ${}_{i_y}\in I^{T_i}_2$ that belongs to a factor
$T=S_{*}\nabla^{(\nu)} R_{ijkl}$,
 we denote by $\sigma(i_y)$ the number
$[(\epsilon_T-1)+\frac{2\nu}{\nu+1}]$ (recall that
$\epsilon_T$ stands
for the number of indices ${}_{r_1},\dots ,{}_{r_\nu},{}_j$ in that factor that
{\it are not} contracting against a factor $\nabla\phi_f$) notice
that if $\nu=0$ this number is zero). Now, we
 also define $C^{l,i_1\dots
\hat{i}_y\dots i_\mu}_{g}(\Omega_1,\dots ,\Omega_p,
\phi_1,\dots,\phi_u)$ as follows: If $\nu \ge 1$, we assume for
convenience that ${}_{i_y}$ is a derivative index (wlog because of the $S_*$-symmetrization). We then define
$C^{l,i_1\dots \hat{i}_y\dots i_\mu}_{g}(\Omega_1,\dots ,\Omega_p,
\phi_1, \dots ,\phi_u)$ to stand for the $(\mu -1)$-tensor field
that arises from $C^{l,i_1\dots i_\mu}_{g}$ by replacing the
factor $S_{*}\nabla^{(\nu)}_{i_yr_2\dots r_\nu}R_{ijkl}$ by
$S_{*}\nabla^{(\nu-1)}_{r_2\dots r_\nu}R_{ijkl}$. If $\nu =0$, we
define $C^{l,i_1\dots \hat{i}_y\dots i_\mu}_{g}(\Omega_1,
\dots,\Omega_p, \phi_1, \dots ,\phi_u)$ to be zero.

\par Applying (\ref{levicivita}), it then follows that the analogue of
(\ref{eudokia2}) in this setting will be:

\begin{equation}
\label{neweudokia}
\begin{split}
&Xdiv_{i_1}\dots Xdiv_{i_a}LC^{No\Phi,free}_{\phi_{u+1}}
[C^{l,i_1\dots i_a}_{g}(\Omega_1,\dots ,\Omega_p,\phi_1,\dots
,\phi_u)]=
\\&-\Sum_{i=1}^{b_l}\Sum_{i_y\in I^{T_i}_2} \sigma(i_y)Xdiv_{i_1}\dots
\hat{Xdiv}_{i_y}\dots Xdiv_{i_a} [\nabla^{i_y}_{T_i}
C^{l,i_1\dots\hat{i}_y\dots i_a}_{g}\\&(\Omega_1,\dots , \Omega_p,
\phi_1,\dots ,\phi_u)\nabla_{i_y}\phi_{u+1}]
+\Sum_{z\in Z}
a_z C^z_{g}(\Omega_1,\dots ,\Omega_p,\phi_1,\dots ,\phi_{u+1});
\end{split}
\end{equation}
(as noted above, if $a>\mu$ then the second line is a
contributor).
\newline

{\it Analogue of (\ref{eudokia2}):} In order to determine the
sublinear combination \\$\Sum_{q\in Q} a_q C^q_{g}(\Omega_1,\dots
,\Omega_p, \phi_1,\dots,\phi_u)$ in (\ref{eudokia2}), we will
start with the case $l\in L_\mu$. We define a number $\tau(i_y)$
for each ${}_{i_y}\in I^{T_i}_2$: Firstly, if ${}_{i_y}$ belongs to a factor
$\nabla^{(A)}\Omega_h$, we define $\tau(i_y)=0$. If ${}_{i_y}$ belongs
to a factor $\nabla^{(m)}R_{ijkl}$ (and by the hypothesis of Lemma
 \ref{pskovb} it must be a derivative index), we define $\tau(i_y)=2$. Finally, if ${}_{i_y}$
belongs to a factor $S_{*}\nabla^{(\nu)} R_{ijkl}$ (and by our
Lemma hypothesis it must be one of the indices ${}_{r_1},\dots
,{}_{r_\nu},{}_j$), then we define $\tau(i_y)=\frac{2\nu}{\nu+1}$. We then
observe that if $l\in L_\mu$, the sublinear combination
$\Sum_{q\in Q} a_q C^q_{g}(\Omega_1,\dots ,\Omega_p,
\phi_1,\dots,\phi_u)$ in (4.12) will be of the form:

\begin{equation}
\label{analyse664}
\begin{split}
&\Sum_{i=1}^{b_l} \Sum_{i_y\in I_2} \tau(i_y)Xdiv_{i_1}\dots
\hat{Xdiv}_{i_y}\dots Xdiv_{i_\mu} [\nabla^{i_y}_{T_i}
C^{l,i_1\dots\hat{i}_y\dots i_\mu}_{g}(\Omega_1,\dots , \Omega_p,
\\&\phi_1,\dots ,\phi_u)\nabla_{i_y}\phi_{u+1}],
\end{split}
\end{equation}
while if $l\in L\setminus L_\mu$, it will simply be of the form:
$$\Sum_{h\in H} a_h C^{h,i_1\dots i_a,i_{*}}_{g}(\Omega_1,\dots ,\Omega_p,
\phi_1,\dots,\phi_u)\nabla_{i_{*}}\phi_{u+1}.$$
\newline

{\bf Analogue of Lemmas \ref{travail}, \ref{travail2} in \cite{alexakis6}:}
\newline

Finally, we want to find a version of Lemmas
4.1 and 4.2 in this context, and also to find a version
 of these two Lemmas in the case where the selected factor is of the form $\nabla^{(A)}\Omega_h$.
As above, this boils down to understanding the sublinear combination 
 $\Sum_{q\in Q} a_q
C^q_{g}(\Omega_1,\dots , \Omega_p, \phi_1,\dots ,\phi_u)$ 
 that appear in the statements of those Lemmas. 
In order to understand this sublinear combination, several
more pieces of notation are needed:

\begin{definition}
\label{polonos} Consider any tensor field  $C^{l,i_1\dots i_a}_g$
in the form (\ref{form2}) with a simple character $\vec{\kappa}_{simp}$.

 Firstly, for each selected factor $T_i$ and for each pair of free indices $({}_{i_k},{}_{i_l})$ with
${}_{i_k}\in I^{T_i}_1, {}_{i_l}\in I^{T_i}_2$ we denote by
$\nabla^{i_{*}}_{T_i}[C^{l,i_1\dots i_a}_{g}(\Omega_1, \dots
,\Omega_p, \phi_1,\dots ,\phi_u)g^{i_ki_l}]
\nabla_{i_{*}}\phi_{u+1}$ the $(a-2)$-tensor field that formally arises
from $C^{l,i_1\dots i_a}_{g}$
 by contracting the index ${}_{i_k}$ against ${}_{i_l}$ and then adding a derivative index $\nabla_{i_{*}}$
 on the selected factor $T_i$ and contracting it against a
 factor $\nabla_{i_{*}}\phi_{u+1}$.

Furthermore, we denote by $F_1,\dots ,F_{\sigma -1}$ the set of
real factors other than $T_i$,\footnote{Recall that the ``real
factors'' are in one of the forms
$\nabla^{(m)}R_{ijkl},S_{*}\nabla^{(\nu)}R_{ijkl},\nabla^{(B)}\Omega_x$.}
and other than any $\nabla\phi_w$. If $|I_1^{T_i}|\ge 2$, we
define $I^{T_i,2,*}_1$ to stand for the set of pairs $({}_{i_k},{}_{i_l})$
for which ${}_{i_k},{}_{i_l}\in I_1^{T_i}$, ${}_{i_k}\ne {}_{i_l}$ and at
 least one of the two indices (say ${}_{i_k}$ with no loss of
 generality) is a derivative index.
 If $|I^{T_i}_1|\le 1$ we define $I^{T_i,2,*}_1=\emptyset$.

 \par For each $({}_{i_k},{}_{i_l})\in I_1^{T_i,2,*}$ and each $S,1\le S
\le \sigma-1$, we then additionally define:

\begin{equation}
\label{malaki} \tilde{C}^{l,i_1\dots ,\hat{i}_k\dots
,i_a,i_z|S}_{g}(\Omega_1,\dots , \Omega_p, \phi_1,\dots
,\phi_u)\nabla_{i_l}\phi_{u+1}
\end{equation}
to stand for the $(a -1)$-tensor field that arises from
$C^{l,i_1\dots ,i_a}_{g}(\Omega_1,\dots , \Omega_p, \phi_1,\dots
,\phi_u)$ by first erasing the free index ${}_{i_k}$ (which is a derivative
index), then contracting ${}_{i_l}$ against a factor
$\nabla\phi_{u+1}$, and finally hitting the $S^{th}$ real factor
 by a derivative free index $\nabla_{i_z}$.
\end{definition}

\par Now, to understand the sublinear combinations:

\begin{equation}
\label{apalakia}
\begin{split}
&Xdiv_{i_1}\dots Xdiv_{i_\mu}
LC^{No\Phi,targ,T_i,B}_{\phi_{u+1}} [C^{l,i_1\dots
i_\mu}_{g}(\Omega_1, \dots ,\Omega_p, \phi_1,\dots ,\phi_u)]+
\\&LC^{No\Phi,div,T_i,I_1}_{\phi_{u+1}}Xdiv_{i_1}\dots
Xdiv_{i_a} C^{l,i_1\dots i_a}_{g}(\Omega_1, \dots ,\Omega_p,
\phi_1,\dots ,\phi_u)]
\end{split}
\end{equation}
 in this context, we
distinguish three subcases: Either the selected factor(s) is (are)
of the form $\nabla^{(m)}R_{ijkl}$, or of the form
$S_{*}\nabla^{(\nu)} R_{ijkl}$, or of the from
$\nabla^{(A)}\Omega_h$. We start with the first subcase, where the
selected factor(s) is (are) of the form $\nabla^{(m)}R_{ijkl}$.

 One more piece of notation: 
\begin{definition}
\label{msharp} 
 For each selected factor $T_i$ we denote by
$m^\sharp_i$ the number of derivative indices in the selected
factor $T_i=\nabla^{(m)}R_{ijkl}$ that are {\it not} contracting
against a factor $\nabla\phi_f$ and {\it are not} free. Also,
recall (from part A) that $\gamma_i$ stands for the number of indices in
$C^{l,i_1\dots i_\mu}_{g}$ that do not belong to the selected
factor $T_i$ and are not contracting against factors
$\nabla\phi_h$.
\end{definition}

We then compute that for each $l\in L_\mu$, where the selected factor
$T_i$ is of the form $T_i=\nabla^{(m_i)}R_{ijkl}$:

\begin{equation}
\label{xontros}
\begin{split}
&Xdiv_{i_1}\dots Xdiv_{i_\mu} LC^{No\Phi,targ,T_i,B}_{\phi_{u+1}}
[C^{l,i_1\dots i_\mu}_{g}(\Omega_1, \dots ,\Omega_p, \phi_1,\dots
,\phi_u)]+
\\&LC^{No\Phi,div,T_i,I_1}_{\phi_{u+1}}Xdiv_{i_1}\dots
Xdiv_{i_a} C^{l,i_1\dots i_a}_{g}(\Omega_1, \dots ,\Omega_p,
\phi_1,\dots ,\phi_u)]= -\Sum_{i_h\in I^{T_i}_1}
\\& [2(m^\sharp_i+2)+\gamma_i]
Xdiv_{i_1}\dots\hat{Xdiv}_{i_h}\dots Xdiv_{i_\mu} 
C^{l,i_1\dots
i_\mu}_{g}(\Omega_1, \dots ,\Omega_p, \phi_1,\dots
,\phi_u)\nabla_{i_h}\phi_{u+1}-
\\&3\Sum_{(i_k,i_l)\in I^{T_i,2,*}_1}
Xdiv_{i_1}\dots \hat{Xdiv}_{i_l}\dots Xdiv_{i_\mu} C^{l,i_1\dots
i_\mu}_{g}(\Omega_1, \dots ,\Omega_p, \phi_1,\dots
,\phi_u)\nabla_{i_l}\phi_{u+1}
\\&+\Sum_{(i_k,i_l)\in I^{T_i,2,*}_1}
\Sum_{S=1}^{\sigma -1}Xdiv_{i_1}\dots \hat{Xdiv}_{i_k}\dots
\hat{Xdiv}_{i_l}\dots Xdiv_{i_\mu}Xdiv_{i_z}
\\&\tilde{C}^{l,i_1\dots
,\hat{i}_k\dots ,i_\mu,i_z|S}_{g}(\Omega_1,\dots , \Omega_p,
\phi_1,\dots ,\phi_u)\nabla_{i_l}\phi_{u+1}
\\&-\Sum_{i_k\in I^{T_i}_1, i_l\in I^{T_i}_2}
Xdiv_{i_1}\dots\hat{Xdiv}_{i_k}\dots \hat{Xdiv}_{i_l}\dots
Xdiv_{i_\mu}\\&\nabla^{i_{*}}_{sel}[C^{l,i_1\dots i_\mu}_{g}
(\Omega_1, \dots,\Omega_p, \phi_1,\dots ,\phi_u)g^{i_ki_l}]
\nabla_{i_{*}}\phi_{u+1}
\\&+\Sum_{h\in H} a_h
Xdiv_{i_1}\dots Xdiv_{i_a} C^{h,i_1\dots
i_a,i_{*}}_{g}(\Omega_1,\dots ,\Omega_p,
 \phi_1,\dots ,\phi_u)\nabla_{i_{*}}\phi_{u+1}.
\end{split}
\end{equation}

{\it Explanation of the calculations 
that bring out (\ref{xontros}):} The expression multiplied by $-2(m^\sharp_i
+2)$ arises in two ways: Firstly, by applying the third summand in
(\ref{levicivita}) to  pairs $(\nabla_{i_k},{}_b)$ where ${}_b$ is
an original non-free index in $C^{l,i_1\dots ,i_a}_{g}$ that is
contracting against the
 crucial factor, and also (when ${}_b={}_i,{}_j,{}_k,{}_l$)
applying the second Bianchi
 identity twice. The second way is by applying the
last summand in (\ref{levicivita}) to two indices
$({}_{i_k},{}_b)$ where both ${}_{i_k}, {}_b$ belong to the
selected
 factor and then ``completing the divergence'' that we have
 created.  With this observation, the rest of the
 sublinear combinations on the right hand side can be
 checked by book-keeping.

\par By the same analysis, if $l\in L\setminus L_\mu$, we calculate:

\begin{equation}
\label{xontrosII}
\begin{split}
&Xdiv_{i_1}\dots Xdiv_{i_a}LC^{No\Phi,targ,T_i,B}_{\phi_{u+1}}
[C^{l,i_1\dots i_a}_{g}(\Omega_1, \dots ,\Omega_p, \phi_1,\dots
,\phi_u)]+
\\&LC^{No\Phi,div,I_1}_{\phi_{u+1}}Xdiv_{i_1}\dots Xdiv_{i_a}
C^{l,i_1\dots i_a}_{g}(\Omega_1, \dots ,\Omega_p, \phi_1,\dots
,\phi_u)]=
\\&\Sum_{h\in H} a_h
C^{h,i_1\dots i_{a-1},i_{*}}_{g}(\Omega_1,\dots ,\Omega_p,
 \phi_1,\dots ,\phi_u)\nabla_{i_{*}}\phi_{u+1}-
\\&\Sum_{i_k\in I^{T_i}_1, i_l\in I^{T_i}_2}
Xdiv_{i_1}\dots \hat{Xdiv}_{i_k}\dots \hat{Xdiv}_{i_l}\dots
 Xdiv_{i_a}\\&\nabla^{i_{*}}_{T_i}[C^{l,i_1\dots i_a}_{g}
(\Omega_1, \dots,\Omega_p, \phi_1,\dots ,\phi_u)g^{i_ki_l}]
\nabla_{i_{*}}\phi_{u+1}.
\end{split}
\end{equation}

\par The analogues of these two equations in the case where the selected factor
is of the form $S_{*}\nabla^{(\nu)} R_{ijkl}$ or
$\nabla^{(A)}\Omega_h$ are straightforward (recall that in this
setting there is a
 unique selected factor so we will write $I_1,I_2$ rather than
$I^{T_i}_1,I_2^{T_i}$). We only have to recall the trivial
 formula:

\begin{equation}
\label{lia} \nabla^{(m)}_{r_1\dots r_m}R_{ir_{m+1}kl}=
S_{*}\nabla^{(m)}_{r_1\dots r_m}R_{ir_{m+1}kl}+ \sum
\nabla^{(m)}_{ir_1\dots r_{m-1}}R_{r_{m-1}r_mkl} +\sum Q(R),
\end{equation}
where $\sum Q(R)$ stands for a generic linear combination of
quadratic expressions in curvature, while $\sum
\nabla^{(m)}_{ir_1\dots r_{m-1}}R_{r_{m}r_{m+1}kl}$ stands for a
generic linear combination of tensors $\nabla^{(m)} R_{abcd}$
where the index ${}_i$ is a derivative index.
\newline

{\it Analogue of Lemmas \ref{travail}, \ref{travail2} when the
selected factor is of the form $S_{*}\nabla^{(\nu)}R_{ijkl}$ 
or $\nabla^{(A)}\Omega_h$:}
We first consider the case where the selected
factor is of the form $S_{*}\nabla^{(\nu)} R_{ijkl}$ (in this case
it will be unique) and we use (\ref{lia}).
 For the selected factor $T=\nabla^{(\nu)}_{r_1\dots r_\nu}R_{ijkl}$ we will
 denote by $\nu^\sharp$ the number of indices
${}_{r_1},\dots ,{}_{r_\nu},{}_j$ that are not free and not
contracting against a factor $\nabla\phi_h$.\footnote{We should write $\nu_l,\nu^\sharp_l$ to stress
that these numbers depend on the tensor field $C^{l,i_1\dots i_\mu}_g$, $l\in L_\mu$.
However, for simplicity of notation, we will not do so.} Then, for each $l\in
L_\mu$, we calculate:

\begin{equation}
\label{xontros2ndcase}
\begin{split}
&Xdiv_{i_1}\dots Xdiv_{i_a}LC^{No\Phi,targ,B}_{\phi_{u+1}}
[C^{l,i_1\dots i_\mu}_{g}(\Omega_1, \dots ,\Omega_p, \phi_1,
\dots,\phi_u)]+
\\&LC^{No\Phi,div,I_1}_{\phi_{u+1}}[Xdiv_{i_1}\dots
Xdiv_{i_\mu} C^{l,i_1\dots i_\mu}_{g}(\Omega_1, \dots ,\Omega_p,
\phi_1, \dots,\phi_u)]=
\\& -\Sum_{i_h\in I_1}[\gamma
+(\nu^\sharp +1)+(\nu^\sharp +\frac{\nu}{\nu +1})]
Xdiv_{i_1}\dots\hat{Xdiv}_{i_h}\dots Xdiv_{i_\mu}
\\&C^{l,i_1\dots i_\mu}_{g}(\Omega_1, \dots ,\Omega_p,
\phi_1,\dots,\phi_u)\nabla_{i_h}\phi_{u+1}]
\\& -3\Sum_{(i_k,i_l)\in I^{2,*}_1}
Xdiv_{i_1}\dots\hat{Xdiv}_{i_l}\dots Xdiv_{i_\mu} C^{l,i_1\dots
i_\mu}_{g}(\Omega_1, \dots ,\Omega_p, \phi_1,\dots
,\phi_u)\nabla_{i_l}\phi_{u+1}
\\&+\Sum_{(i_k,i_l)\in I^{2,*}_1}\Sum_{S=1}^{\sigma -1}Xdiv_{i_1}\dots \hat{Xdiv}_{i_k}\dots
\hat{Xdiv}_{i_l}\dots Xdiv_{i_\mu}Xdiv_{i_z}
\\&\tilde{C}^{l,i_1\dots
,\hat{i}_k\dots ,i_\mu,i_z|S}_{g}(\Omega_1,\dots , \Omega_p,
\phi_1,\dots ,\phi_u)
\\&\nabla_{i_l}\phi_{u+1}-\Sum_{i_k\in I_1, i_l\in I_2}Xdiv_{i_1}\dots
\hat{Xdiv}_{i_k}\dots \hat{Xdiv}_{i_l}\dots Xdiv_{i_\mu}
\\&\nabla^{i_{*}}_{sel}[C^{l,i_1\dots i_\mu}_{g}(\Omega_1, \dots
,\Omega_p, \phi_1,\dots ,\phi_u)g^{i_ki_l}]
\nabla_{i_{*}}\phi_{u+1}
\\&+\Sum_{h\in H} a_h Xdiv_{i_1}\dots Xdiv_{i_a}
C^{h,i_1\dots i_a,i_{*}}_{g}(\Omega_1,\dots ,\Omega_p,
 \phi_1,\dots ,\phi_u)\nabla_{i_{*}}\phi_{u+1},
\end{split}
\end{equation}
while for each $l\in L\setminus L_\mu$, we again have
(\ref{xontrosII}). A small extra explanation for this case: The
 sublinear combination multiplied by
 $-(\nu^\sharp +1+\frac{\nu}{\nu+1})$ arises by virtue of
 applying the last summand in (\ref{levicivita}) to a pair of
 indices in the selected factor, one of which is free and one
 of which is not. We see that this formula can only be applied
 to indices
 $({}_{i_k},{}_b)$ in the selected factor if at least one of
 these indices is a derivative index.

Finally, in the case where the selected factor is of the form
$\nabla^{(A)}\Omega_h$, we denote by $A^\sharp$ the number of
indices in $\nabla^{(A)}\Omega_h$ that are not free and not
contracting against a factor $\nabla\phi_f$. We then derive that for
each $l\in L$, ($a\ge \mu$):

\begin{equation}
\label{xontros''}
\begin{split}
&Xdiv_{i_1}\dots Xdiv_{i_a}LC^{No\Phi,targ,B}_{\phi_{u+1}}
[C^{l,i_1\dots i_a}_{g}(\Omega_1, \dots ,\Omega_p, \phi_1,\dots
,\phi_u)]+
\\&LC^{No\Phi,div,I_1}_{\phi_{u+1}}Xdiv_{i_1}\dots Xdiv_{i_a}
C^{l,i_1\dots i_a}_{g}(\Omega_1, \dots ,\Omega_p, \phi_1,
\dots,\phi_u)]=
\\& -\Sum_{i_h\in I_1}(\gamma  +2A^\sharp)
Xdiv_{i_1}\dots \hat{Xdiv}_{i_h}\dots Xdiv_{i_a} C^{l,i_1\dots
i_a}_{g}(\Omega_1, \dots ,\Omega_p, \phi_1,\dots
\\&,\phi_u)\nabla_{i_h}\phi_{u+1}
\\& -3\Sum_{(i_k,i_l)\in I^{2,*}_1}
Xdiv_{i_1}\dots \hat{Xdiv}_{i_l}\dots Xdiv_{i_a} C^{l,i_1\dots
i_a}_{g}(\Omega_1, \dots ,\Omega_p, \phi_1,\dots
\\&,\phi_u)\nabla_{i_k}\phi_{u+1}
\\&+\Sum_{(i_k,i_l)\in I^{2,*}_1} \Sum_{S=1}^{\sigma -1}Xdiv_{i_1}\dots \hat{Xdiv}_{i_k}\dots
\hat{Xdiv}_{i_l}\dots Xdiv_{i_\mu}Xdiv_{i_z}
\\&\tilde{C}^{l,i_1\dots
,\hat{i}_k\dots ,i_\mu,i_z|S}_{g}(\Omega_1,\dots , \Omega_p,
\phi_1,\dots ,\phi_u)\nabla_{i_l}\phi_{u+1}
\\&-\Sum_{i_k\in I_1, i_l\in I_2}Xdiv_{i_1}\dots
\hat{Xdiv}_{i_k}\dots \hat{Xdiv}_{i_l}\dots Xdiv_{i_a}
\\&\nabla^{i_{*}}_{crit}[C^{l,i_1\dots i_a}_{g}(\Omega_1, \dots
,\Omega_p, \phi_1,\dots ,\phi_u)g^{i_ki_l}]
\nabla_{i_{*}}\phi_{u+1}
\\&+\Sum_{h\in H} a_hXdiv_{i_1}\dots Xdiv_{i_a}
C^{h,i_1\dots i_a,i_{*}}_{g}(\Omega_1,\dots ,\Omega_p,
 \phi_1,\dots ,\phi_u)\nabla_{i_{*}}\phi_{u+1}.
\end{split}
\end{equation}

Finally, for each $l\in L\setminus L_\mu$, we again have
(\ref{xontrosII}).
\newline

\par Thus, having computed the sublinear
 combinations $LC^{No\Phi}[L_g], W[L_g]$ in the setting of Lemma \ref{pskovb}, we will
 now substitute them into (\ref{heidegger2}).

 We first consider the case where
the selected factor is of the form $S_{*}\nabla^{(\nu)} R_{ijkl}$.
(Recall that in this setting there is only one selected
 factor for each $C^{l,i_1\dots i_\mu}_{g},l\in L^z$, $z\in Z'_{Max}$).
We then derive a local equation:

\begin{equation}
\label{proolaxreiazontai1}
\begin{split}
&\Sum_{l\in L_\mu} a_l\{ -\Sum_{i_h\in I_1}(\gamma
+(2\nu^\sharp+1)+\frac{\nu}{\nu +1}-2(\sigma_1+\sigma_2-1)-X)
\\& Xdiv_{i_1}\dots\hat{Xdiv}_{i_h}\dots
Xdiv_{i_\mu}C^{l,i_1\dots i_\mu}_{g}(\Omega_1, \dots ,\Omega_p,
\phi_1,\dots ,\phi_u)\nabla_{i_h}\phi_{u+1}
\\&-3\Sum_{(i_k,i_l)\in I^{2,*}_1}
Xdiv_{i_1}\dots\hat{Xdiv}_{i_k}\dots Xdiv_{i_\mu} C^{l,i_1\dots
i_\mu}_{g}(\Omega_1, \dots ,\Omega_p, \phi_1,\dots
,\phi_u)\nabla_{i_k}\phi_{u+1}
\\&-\frac{|I_1|}{\nu+1}Xdiv_{i_2}\dots Xdiv_{i_\mu}C^{l,i_1
\dots i_\mu}_{g}(\Omega_1, \dots ,\Omega_p, \phi_1,\dots ,
\phi_u)\nabla_{i_1}\phi_{u+1}+
\\&\Sum_{i_y\in I_2} (-\sigma(i_y)+\tau(i_y))
Xdiv_{i_1}\dots \hat{Xdiv}_{i_y}\dots Xdiv_{i_\mu}
\\& [\nabla^{i_y}_{sel} C^{l,i_1\dots\hat{i}_y\dots
i_\mu}_{g}(\Omega_1,\dots , \Omega_p, \phi_1,\dots
,\phi_u)\nabla_{i_y}\phi_{u+1}]+
\\&\Sum_{(i_k,i_l)\in I^{2,*}_1}
\Sum_{S=1}^{\sigma -1}Xdiv_{i_1}\dots \hat{Xdiv}_{i_k}\dots
\hat{Xdiv}_{i_l}\dots Xdiv_{i_\mu}Xdiv_{i_z}
\\&\tilde{C}^{l,i_1\dots
,\hat{i}_k\dots ,i_\mu,i_z|S}_{g}(\Omega_1,\dots , \Omega_p,
\phi_1,\dots ,\phi_u)\nabla_{i_l}\phi_{u+1}\}
\\& -\Sum_{l\in L} a_l \Sum_{(i_k,i_l)\in
I^{2,def}_2} Xdiv_{i_1}\dots \hat{Xdiv}_{i_k}\dots
\hat{Xdiv}_{i_l}\dots Xdiv_{i_a}
\\&[C^{l,i_1\dots
i_a,i_{*}}_{g}(\Omega_1,\dots ,\Omega_p,\phi_1,\dots
,\phi_u)g^{i_ki_l}]\nabla_{i_{*}}\phi_{u+1}
\\&-\Sum_{l\in L} a_l \Sum_{i_k\in I_1, i_l\in I_2}Xdiv_{i_1}\dots
\hat{Xdiv}_{i_k}\dots \hat{Xdiv}_{i_l}\dots Xdiv_{i_a}
\\& \nabla^{i_{*}}_{sel}
[C^{l,i_1\dots i_a}_{g}(\Omega_1, \dots ,\Omega_p, \phi_1,\dots
,\phi_u)g^{i_ki_l}] \nabla_{i_{*}}\phi_{u+1}+
\\&\Sum_{j\in J} a_j C^j_{g}(\Omega_1,\dots ,\Omega_p,
 \phi_1,\dots ,\phi_u,\phi_{u+1})+
\\&\Sum_{h\in H} a_h Xdiv_{i_1}\dots Xdiv_{i_a}
C^{h,i_1\dots i_a,i_{*}}_{g}(\Omega_1,\dots ,\Omega_p,
 \phi_1,\dots ,\phi_u)\nabla_{i_{*}}\phi_{u+1}=0,
\end{split}
\end{equation}
modulo complete contractions of length $\ge\sigma +u+2$.

In the case where the selected factor(s) is (are) of the form
$\nabla^{(m)}R_{ijkl}$, we derive:

\begin{equation}
\label{proolaxreiazontai2}
\begin{split}
&\Sum_{l\in L_\mu} a_l\{-\Sum_{i=1}^{b_l} \Sum_{i_h\in
I^{T_i}_1}(\gamma_i +2(m^\sharp_i +2)-2(\sigma_1+\sigma_2-1)-X)
\\& Xdiv_{i_1}\dots\hat{Xdiv}_{i_h}\dots
Xdiv_{i_\mu}C^{l,i_1\dots i_\mu}_{g}(\Omega_1, \dots ,\Omega_p,
\phi_1,\dots ,\phi_u)\nabla_{i_h}\phi_{u+1}
\\& -3\Sum_{(i_k,i_l)\in I^{T_i,2,*}_1}
Xdiv_{i_1}\dots\hat{Xdiv}_{i_l}\dots Xdiv_{i_\mu} C^{l,i_1\dots
i_\mu}_{g}(\Omega_1, \dots ,\Omega_p, \phi_1,\dots
,\phi_u)\nabla_{i_l}\phi_{u+1}
\\&+\Sum_{i=1}^{b_l}
\Sum_{i_y\in I^{T_i}_2} (-\sigma(i_y)+\tau(i_y)) Xdiv_{i_1}\dots
\hat{Xdiv}_{i_y}\dots Xdiv_{i_\mu}
\\&\nabla^{i_y}_{T_i} C^{l,i_1\dots\hat{i}_y\dots
i_\mu}_{g}(\Omega_1,\dots , \Omega_p, \phi_1,\dots
,\phi_u)\nabla_{i_y}\phi_{u+1}
\\&+\Sum_{(i_k,i_l)\in I^{2,*}_1}\Sum_{S=1}^{\sigma
-1}Xdiv_{i_1}\dots \hat{Xdiv}_{i_k}\dots \hat{Xdiv}_{i_l}\dots
Xdiv_{i_\mu}Xdiv_{i_z}
\\& \tilde{C}^{l,i_1\dots ,\hat{i}_k\dots
,i_\mu,i_z|S}_{g}(\Omega_1,\dots , \Omega_p, \phi_1,\dots
,\phi_u)\nabla_{i_l}\phi_{u+1}\}
\\& -\Sum_{l\in L} a_l \Sum_{i=1}^{b_l}
\Sum_{(i_k,i_l)\in (I^{T_i}_2)^{2,def}} Xdiv_{i_1}\dots
\hat{Xdiv}_{i_k}\dots \hat{Xdiv}_{i_l}\dots Xdiv_{i_a}
\\& [C^{l,i_1\dots i_a,i_{*}|T_i}_{g}(\Omega_1,\dots ,\Omega_p,
\phi_1,\dots,\phi_u)g^{i_ki_l}]\nabla_{i_{*}}\phi_{u+1}
\\&-\Sum_{l\in L} a_l \Sum_{i=1}^{b_l}\Sum_{i_k\in
I^{T_i}_1, i_l\in I^{T_i}_2}Xdiv_{i_1}\dots
\hat{Xdiv}_{i_k}\dots \hat{Xdiv}_{i_l}\dots Xdiv_{i_a}
\\&\nabla^{i_{*}}_{T_i}[C^{l,i_1\dots i_a}_{g}(\Omega_1,
\dots,\Omega_p, \phi_1,\dots ,\phi_u)g^{i_ki_l}]
\nabla_{i_{*}}\phi_{u+1}+
\\&\Sum_{j\in J} a_j C^j_{g}(\Omega_1,\dots ,\Omega_p,
 \phi_1,\dots ,\phi_u,\phi_{u+1})+
\\&\Sum_{h\in H} a_h Xdiv_{i_1}\dots Xdiv_{i_a}
C^{h,i_1\dots i_a,i_{*}}_{g}(\Omega_1,\dots ,\Omega_p,
 \phi_1,\dots ,\phi_u)\nabla_{i_{*}}\phi_{u+1}=0,
\end{split}
\end{equation}
modulo complete contractions of length $\ge\sigma +u+2$.

Finally, in the case where the selected factor is of the form
$\nabla^{(A)}\Omega_h$ (in which case it is again unique),
we derive:

\begin{equation}
\label{proolaxreiazontai3}
\begin{split}
&\Sum_{l\in L_\mu} a_l \{\Sum_{i_h\in I_1} (\gamma
 +2A^\sharp -2\sigma_1-2\sigma_2)
\\& Xdiv_{i_1}\dots \hat{Xdiv}_{i_h}\dots
Xdiv_{i_\mu}C^{l,i_1\dots i_\mu}_{g}(\Omega_1, \dots ,\Omega_p,
\phi_1,\dots ,\phi_u)\nabla_{i_h}\phi_{u+1}
\\& -3\Sum_{(i_k,i_l)\in I^{2,*}_1}
Xdiv_{i_1}\dots\hat{Xdiv}_{i_l}\dots Xdiv_{i_\mu} C^{l,i_1\dots
i_a}_{g}(\Omega_1, \dots ,\Omega_p, \phi_1,\dots
,\phi_u)\nabla_{i_l}\phi_{u+1}
\\&+\Sum_{i_y\in I_2} (-\sigma(i_y)+\tau(i_y))
Xdiv_{i_1}\dots \hat{Xdiv}_{i_y}\dots Xdiv_{i_\mu}
\\&[\nabla^{i_y}_{sel} C^{l,i_1\dots\hat{i}_y\dots
i_\mu}_{g}(\Omega_1,\dots , \Omega_p, \phi_1,\dots
,\phi_u)\nabla_{i_y}\phi_{u+1}]+
\\&\Sum_{(i_k,i_l)\in I^{2,*}_1}\Sum_{S=1}^{\sigma
-1}Xdiv_{i_1}\dots \hat{Xdiv}_{i_k}\dots \hat{Xdiv}_{i_l}\dots
Xdiv_{i_\mu}Xdiv_{i_z}
\\&\tilde{C}^{l,i_1\dots ,\hat{i}_k\dots
,i_\mu,i_z|S}_{g}(\Omega_1,\dots , \Omega_p, \phi_1,\dots
,\phi_u)\nabla_{i_l}\phi_{u+1}\}
\\& -\Sum_{l\in L} a_l \Sum_{(i_k,i_l)\in
I^{2,def}_2} Xdiv_{i_1}\dots \hat{Xdiv}_{i_k}\dots
\hat{Xdiv}_{i_l}\dots Xdiv_{i_a}
\\&[C^{l,i_1\dots
i_a,i_{*}}_{g}(\Omega_1,\dots ,\Omega_p,\phi_1,\dots
,\phi_u)g^{i_ki_l}]\nabla_{i_{*}}\phi_{u+1}
\\&-\Sum_{l\in L} a_l \Sum_{i_k\in I_1, i_l\in I_2}Xdiv_{i_1}\dots
\hat{Xdiv}_{i_k}\dots \hat{Xdiv}_{i_l}\dots Xdiv_{i_a}
\\&\nabla^{i_{*}}_{sel}[C^{l,i_1\dots i_a}_{g}(\Omega_1, \dots
,\Omega_p, \phi_1,\dots ,\phi_u)g^{i_ki_l}]
\nabla_{i_{*}}\phi_{u+1}+
\\&\Sum_{j\in J} a_j C^j_{g}(\Omega_1,\dots ,\Omega_p,
 \phi_1,\dots ,\phi_u,\phi_{u+1})+
\\&\Sum_{h\in H} a_h Xdiv_{i_1}\dots Xdiv_{i_a}
C^{h,i_1\dots i_a,i_{*}}_{g}(\Omega_1,\dots ,\Omega_p,
 \phi_1,\dots ,\phi_u)\nabla_{i_{*}}\phi_{u+1}=0,
\end{split}
\end{equation}
modulo complete contractions of length $\ge\sigma +u+2$.

\par Now, in order to motivate our next section, let us briefly illustrate how the above three equations {\it
are not} in the form required by Lemma \ref{pskovb} (in case A).
There are two important defficiencies of the three equations above
that we wish to highlight: Firstly, observe that in all equations
above the sublinear combinations indexed in $\Sum_{(i_k,i_l)\in
I^{2,def}_2}$ have rank $(\mu-2)$. The desired conclusion of Lemma
\ref{pskovb} requires the minimum rank of the tensor fields to be
$\mu-1$. Secondly, the coefficient $(\gamma
+(2\nu^\sharp+1)+\frac{\nu}{\nu +1}-2(\sigma_1+\sigma_2-1))$ (and also the
coefficients in the fist lines of (\ref{proolaxreiazontai2}), (\ref{proolaxreiazontai3}))
{\it depends} on the number $\nu$ of derivatives on the selected factor
in each individual tensor field in the second line in
(\ref{proolaxreiazontai1}) (and similarly for the other two
equations). This is not the case in the conclusion of Lemma
\ref{pskovb} where all tensor fields of rank $\mu-1$ must be
multiplied by a {\it universal} constant (either 1 or $\alpha\choose{2}$).

\par Thus, we observe that we {\it can not} derive our Lemma
\ref{pskovb} merely by repeating the  analysis of equation
$Image^{1,+}_{\phi_{u+1}}[L_g]$ that we performed to
 derive the Lemmas 3.1 and 3.2 in \cite{alexakis6}. The second ingredient in the
proof of Lemma \ref{pskovb} will be an analysis of the equation
$Image^{1,\beta}_{\phi_{u+1}}[L_g]=0$ (see Definition
2.3 and equation (2.9) in \cite{alexakis6}),
coupled with a formal operation $Soph\{\dots \}$ which will
then {\it convert} the terms in that equation into the type of terms that are required
in Lemma \ref{pskovb}. 

\section{The second part of the ``grand conclusion'': 
\\A study of  $Image^{1,\beta}_{\phi_{u+1}}[L_g]=0$.}
 \label{uphrxeden}

\subsection{Basic calculations in the equation $Image^{1,\beta}_{\phi_{u+1}}[L_g]=0$, and 
 first steps in its analysis.}
\label{first.steps}

\par We recall from definition \ref{newdefinition} that
$Image^{1,\beta}_{\phi_{u+1}}[L_g]$ consists of the complete
contractions in $Image^{1}_{\phi_{u+1}}[L_g]$ that have one
internal contraction, and also either have length $\sigma+u$ and
a factor $\nabla^{(A)}\phi_{u+1}$, $A\ge 2$, or have length
$\sigma+u+1$ and a factor $\nabla\phi_{u+1}$. {\it Recall that the
selected factor(s), the crucial factor(s) etc. are completely
irrelevant in the context of $Image^{1,\beta}_{\phi_{u+1}}[L_g]$}. We
recall that we have denoted by
$Image^{1,\beta,\sigma+u}_{\phi_{u+1}}[L_g]$ the sublinear
combination of contractions with $\sigma+u$ factors in
$Image^{1,\beta}_{\phi_{u+1}}[L_g]$, and by
$Image^{1,\beta,\sigma+u+1}_{\phi_{u+1}}[L_g]$ the sublinear
combination of contractions with $\sigma+u+1$ factors in
$Image^{1,\beta}_{\phi_{u+1}}[L_g]$.

\par We recall that in (2.9) we showed that modulo complete
 contractions with at least $\sigma+u+2$ factors:

\begin{equation}
\label{jgallo} Image^{1,\beta,\sigma+u}_{\phi_{u+1}}[L_g]+
Image^{1,\beta,\sigma+u+1}_{\phi_{u+1}}[L_g]+\sum_{z\in Z} a_z
C^z_{g}(\Omega_1,\dots,\Omega_p,\phi_1,\dots,\phi_{u+1})=0.
\end{equation}
(We recall that $\sum_{z\in Z}\dots$ stands for a {\it generic}
linear combination of contractions with $\sigma+u+1$, one of which
is in the from $\nabla^{(B)}\phi_{u+1}$, $B\ge 2$).

\par Now, in the rest of this section we will study the LHS of the equation
(\ref{jgallo}) and repeatedly apply the inductive assumption of Corollary 
1 in \cite{alexakis4} to 
it in many different guises, in order to derive a new local equation which will help us in
 deriving Lemma \ref{pskovb}. In particular, the new local
  equation combined with (\ref{proolaxreiazontai1}), 
(\ref{proolaxreiazontai2}), (\ref{proolaxreiazontai3}) will give us our ``grand conclusion''. 
 
 \par Firstly, we focus of the
complete contractions of length $\sigma +u$ in
$Image^{1,\beta}_{\phi_{u+1}}[L_{g}]$.

\par We initially seek to understand {\it how} the complete
 contractions in $Image^{1,\beta}_{\phi_{u+1}}[L_{g}]$ can
arise. For the purposes of this subsection, when we study the
transformation law of each factor of the form $S_{*}\nabla^{(\nu)}
R_{ijkl}$, we will be treating each such factor as a sum of
 factors in the form $\nabla^{(m)}R_{ijkl}$.
With this convention we immediately see that the complete
contractions of
 length $\sigma +u$ in $Image^{1,\beta}_{\phi_{u+1}}
[L_{g}]$ can only arise by replacing a factor
$\nabla^{(m)}R_{ijkl}$ by an expression
$\nabla^{(m)}[\nabla^{(2)}_{cd}\phi_{u+1}\otimes g_{ab}]$, {\it
provided the indices ${}_a,{}_b$ both contract against the same
factor}.

\par We then recall the operation $Sub_{\omega}$
(from the Appendix in \cite{alexakis1})
 with which we will act on the complete
 contractions in $Image^{1,\beta}_{\phi_{u+1}}[L_{g}]$:

\begin{definition}
\label{crown} If a complete contraction $C_{g}(\Omega_1,\dots
,\Omega_p,\phi_1,\dots ,\phi_u)$ has its internal contraction
 in a factor of the form $\nabla^{(p)}_{r_1\dots r_p}Ric_{ij}$, then 
 $Sub_\omega[C_{g}]$  stands for the complete contraction
that arises from $C_{g}$ by replacing the factor
$\nabla^{(p)}_{r_1\dots r_p} Ric_{ij}$ by a factor
$-\nabla^{(p+2)}_{r_1\dots r_pij}\omega$.
 If the internal contraction involves at least one
 derivative index (hence it is in the form
$(\nabla^a,{}_a)$),  then $Sub_\omega[C_{g}(\Omega_1,\dots
,\Omega_p,\phi_1,\dots ,\phi_u)]$
  stands for the complete contraction that arises from
$C_{g}(\Omega_1,\dots ,\Omega_p,\phi_1,\dots ,\phi_u)$ by
replacing the expression $(\nabla^a,{}_a)$ by an
 expression $(\nabla^a\omega, {}_a)$.\footnote{In other words the 
derivative index $\nabla^a$ is erased, and a new factor $\nabla\omega$ 
is introduced, which is then contracted against the index ${}_a$ (${}_a$ is the index
 that originally contracted against $\nabla^a$).} 
 The above extends to linear
 combinations of complete contractions and also to tensor
 fields.
\end{definition}
Hence, by applying  the last Lemma in the Appendix of \cite{alexakis1} to the
equation (\ref{jgallo}), we derive that:

\begin{equation}
\label{baldouin} Sub_\omega \{
Image^{1,\beta}_{\phi_{u+1}}[L_{g}]\}+ \Sum_{t\in T} a_t
C^t_{g}(\Omega_1,\dots , \Omega_p,\phi_1,\dots
,\phi_{u+1},\omega)=0,
\end{equation}
modulo complete contractions of length $\ge\sigma +u+3$. Here
\\$\Sum_{t\in T} a_t C^t_{g}(\Omega_1,\dots , \Omega_p,\phi_1,\dots
,\phi_u,\omega)$ stands for a generic linear combination of
complete contractions with either length $\sigma +u+2$ and a
factor $\nabla^{(A)}\omega, A\ge 2$, or a factor
$\nabla^{(A)}\phi_{u+1}$, $A\ge 2$, {\it or} with length
$\sigma+u+1$ and a factor $\nabla^{(A)}\phi_{u+1}, A\ge 2$ {\it
and} a factor $\nabla^{(B)}\omega, B\ge 2$.
 This equation follows from the last Lemma
 in the Appendix of \cite{alexakis1}, and the definition of the various terms in (\ref{jgallo}).
\newline

{\it A brief analysis of equation (\ref{baldouin}):} We observe
that by definition \\$Sub_\omega
\{Image^{1,\beta}_{\phi_{u+1}}[L_{g}]\}$ consists of three main
sublinear combinations,\footnote{Plus junk terms of greater length.}
 depending on the total number of factors
that a given complete contraction in $Sub_\omega
\{Image^{1,\beta}_{\phi_{u+1}}[L_{g}]\}$ contains: A given
complete contraction in \\$Sub_\omega
\{Image^{1,\beta}_{\phi_{u+1}}[L_{g}]\}$ may contain $\sigma+u,\sigma+u+1$ or
$\sigma+u+2$ factors in total. Accordingly, we denote the
corresponding sublinear combinations by $Sub_\omega^{\sigma+u}
\{Image^{1,\beta}_{\phi_{u+1}}[L_{g}]\}$, $Sub_\omega^{\sigma+u+1}
\{Image^{1,\beta}_{\phi_{u+1}}[L_{g}]\}$, $Sub_\omega^{\sigma+u+2}
\{Image^{1,\beta}_{\phi_{u+1}}[L_{g}]\}$. By definition, the
complete contractions in $Sub_\omega^{\sigma+u}
\{Image^{1,\beta}_{\phi_{u+1}}[L_{g}]\}$
 will have two factors $\nabla^{(A)}\phi_{u+1},\nabla^{(B)}\omega$, $A,B\ge 2$.
 The complete contractions in
\\$Sub_\omega^{\sigma+u+1} \{Image^{1,\beta}_{\phi_{u+1}}[L_{g}]\}$
 will have two factors, either in the form
 $\nabla^{(A)}\phi_{u+1},\nabla\omega$, $A\ge 2$ {\it or}
in the form $\nabla\phi_{u+1},\nabla^{(B)}\omega$, $B\ge 2$.
Finally, the complete contractions  in
\\$Sub_\omega^{\sigma+u+2} \{Image^{1,\beta}_{\phi_{u+1}}[L_{g}]\}$
 will have two factors in the form $\nabla\phi_{u+1},\nabla\omega$.

\par Thus, (\ref{baldouin}) can be re-expressed in the form:

\begin{equation}
\label{baldouin'}
\begin{split}
& Sub^{\sigma+u}_\omega \{
Image^{1,\beta}_{\phi_{u+1}}[L_{g}]\}+
Sub^{\sigma+u+1}_\omega \{
Image^{1,\beta}_{\phi_{u+1}}[L_{g}]\}+
\\& Sub^{\sigma+u+2}_\omega \{
Image^{1,\beta}_{\phi_{u+1}}[L_{g}]\}+
 \Sum_{t\in T} a_t
C^t_{g}(\Omega_1,\dots , \Omega_p,\phi_1,\dots
,\phi_{u+1},\omega)=0.
\end{split}
\end{equation}
(modulo complete contractions of length
$\ge\sigma+u+3$). For future reference, we will also divide the
index set $T$ into subsets $T^{\sigma+u+1},T^{\sigma+u+2}$,
according to the number of factors that a given complete contraction
$C^t_{g}(\Omega_1,\dots , \Omega_p,\phi_1,\dots
,\phi_{u+1},\omega)$ contains.
\newline

\par We  first seek to understand the sublinear
 combination $Sub^{\sigma
+u}_\omega \{ Image^{1,\beta}_{\phi_{u+1}} [L_{g}]\}$ in
$Sub_\omega \{ Image^{1,\beta}_{\phi_{u+1}} [L_{g}]\}$ which
consists of complete contractions of length $\sigma +u$.  We start by  studying how $Sub_\omega \{
Image^{1,\beta}_{\phi_{u+1}} [L_{g}]\}$ arises from the equation $L_g=0$ and we will then
prove two important equations regarding this sublinear
combination, namely equations (\ref{mikul}) and (\ref{tzieria}).
{\it Important Note:} For future reference, we note here that all
the discussion until (\ref{mikul}) and (\ref{tzieria}), and also
the proofs of these equations also hold {\it without} the
assumptions that $L^{*}_\mu\bigcup L^{+}_\mu\bigcup
L''_{+}=\emptyset$. This will be used in the proof 
 of Lemmas 3.3, 3.4 from \cite{alexakis4}. Furthermore, the
assumptions $L''_{+}=\emptyset$ will not be used until
after equation (\ref{gounzales1}).
\newline

{\bf How does the sublinear combination $Sub^{\sigma +u}_\omega \{
Image^{1,\beta}_{\phi_{u+1}} [L_{g}]\}$ arise?}

\par For each tensor field $C^{l,i_1\dots i_a}_{g}$ and each
complete contraction $C^j_{g}$ in $L_{g}$, we look for the pairs
of non-anti-symmetric internal indices $({}_a,{}_b), ({}_c,{}_d)$ in two
different factors $\nabla^{(m)}R_{ijkl}$ (recall that we are
treating the factors $S_{*}\nabla^{(\nu)} R_{ijkl}$
 as sums of tensors in the form $\nabla^{(m)}R_{ijkl}$) for which
  ${}_a,{}_c$ and ${}_b,{}_d$ are contracting against each other.

In each tensor field $C^{l,i_1\dots i_a}_{g}$ or complete
contraction $C^j_{g}$, we denote the
 set of those pairs by $INT^2_l, INT^2_j$, respectively. For each such pair
$[({}_a{}_,{}_b), ({}_c,{}_d)]\in INT^2_l$ or $[({}_a,{}_b), ({}_c,{}_d)]\in INT^2_j$, we
denote by $Rep^1[C^{l,i_1,\dots i_a}_{g}]$, $Rep^1[C^j_{g}]$ the
complete
 contraction or tensor field that arises by replacing the
 first factor $\nabla^{(m)}R_{ijkl}$ by the linear expression
$(+-)\nabla^{(m+2)}\phi\otimes g_{ab}$ on the right hand side
 of (\ref{curvtrans}) and then applying $Sub_\omega$
 (here, of course, the $+/-$ sign comes from (\ref{curvtrans})).
We also denote by $Rep^2[C^j_{g}]$, $Rep^2[C^{l,i_1,\dots
,i_a}_{g}]$ the tensor field (or complete contraction) that arises
by replacing the second factor $\nabla^{(m')}R_{ijkl}$ by the
linear expression $(+-)\nabla^{(m'+2)}\phi_{u+1}\otimes g_{cd}$ on
the right hand side of (\ref{curvtrans}) and then applying
$Sub_\omega$.

\par We define:

\begin{equation}
\label{reno} \begin{split}&Rep[C^{l,i_1,\dots i_a}_{g}]=\Sum_{[(a,b), (c,d)]\in
INT^2_l} \{Rep^1_{[(a,b), (c,d)]} [C^{l,i_1,\dots
,i_a}_{g}]\\&+Rep^2_{[(a,b), (c,d)]}[C^{l,i_1,\dots ,i_a}_{g}]\}
\end{split}\end{equation}
and
\begin{equation}
\label{reno2} Rep[C^j_{g}]=\Sum_{[(a,b), (c,d)]\in INT^2_l}\{
Rep^1_{[(a,b), (c,d)]}[C^j_{g}]+Rep^2_{[(a,b), (c,d)]} [C^j_{g}]\}.
\end{equation}

\par We straightforwardly observe that:

\begin{equation}
\label{mikulski}
\begin{split}
&Sub^{\sigma +u}_\omega \{ Image^{1,\beta}_{\phi_{u+1}}
[L_{g}]\}=\Sum_{j\in J} a_j
Rep[C^j_{g}(\Omega_1,\dots ,\Omega_p,\phi_1,\dots ,\phi_u)]
\\&\Sum_{l\in L} a_l Xdiv_{i_1}\dots Xdiv_{i_a}
Rep[C^{l,i_1,\dots ,i_a}_{g}(\Omega_1,\dots ,\Omega_p,\phi_1,\dots
,\phi_u)].
\end{split}
\end{equation}

\par Moreover, by (\ref{baldouin}) we have that:

\begin{equation}
\label{reno3} Sub^{\sigma +u}_\omega \{
Image^{1,\beta}_{\phi_{u+1}} [L_{g}]\}=0,
\end{equation}
modulo complete contractions of length $\ge\sigma +u+1$.

\par Before we move on to examine the other sublinear combinations of
\\$Sub^{\sigma +u}_\omega \{ Image^{1,\beta}_{\phi_{u+1}}
[L_{g}]\}$, we want to somehow apply Corollary 1 in \cite{alexakis4} to
(\ref{reno3}). We introduce some notational conventions:

\begin{definition}
\label{leaves}

\par We will denote by
$$\Sum_{j\in J^{\sigma +u}} a_j C^j_{g}
(\Omega_1,\dots ,\Omega_p,\phi_1,\dots ,\phi_u,\phi_{u+1},
\omega)$$ a generic linear combination of complete contractions
with length $\sigma +u$ and two factors $\nabla^{(A)}\phi_{u+1},
\nabla^{(B)}\omega$ and  at least one factor $\nabla\phi_f, f\in
Def(\vec{\kappa}_{simp})$ contracting against a derivative
 index in a factor $\nabla^{(m)}R_{ijkl}$ {\it or} contracting against
 one of the first $A-2$ indices in a factor
$\nabla^{(A)}\phi_{u+1}$ or $\nabla^{(A)}\omega$.

We denote by
$$\Sum_{u\in U_1} a_u
C^{u,i_1\dots ,i_a}_{g}(\Omega_1,\dots ,\Omega_p,\phi_1, \dots
,\phi_u,\phi_{u+1},\omega)$$ ($a\ge \mu$) a generic linear combination of
acceptable tensor fields with length $\sigma +u+1$, a factor
$\nabla\phi_{u+1}$ (which we treat as a factor
$\nabla\phi$)\footnote{In particular this factor does not contain a free index.}
 and a factor $\nabla^{(A)}\omega$, $A\ge 2$
(which we treat as a factor $\nabla^{(A)}\Omega_{p+1}$), with the
additional feature that the $u$-simple character of $C^{u,i_1\dots
,i_a}_{g}$ arises from $\vec{\kappa}_{simp}$ by replacing one
factor $T=\nabla^{(m)}R_{ijkl}$ or $T=S_{*}\nabla^{(\nu)}R_{ijkl}$
by a factor $T'=\nabla^{(\nu+2)}\Omega_{p+1}$, where all the
factors $\nabla\phi_h$ that contracted against $T$ now contract
against $T'$; moreover if $a=\mu$ and the factor
$\nabla^{(A)}\omega$ has $A=2$ and is contracting against a factor $\nabla\phi_h$
 then we additionnaly require that it does not contain free indices.

 We denote by
$$\Sum_{u\in U_2} a_u
C^{u,i_1\dots i_a}_{g}(\Omega_1,\dots ,\Omega_p,\phi_1, \dots
,\phi_u,\phi_{u+1},\omega)$$ a generic linear combination of
acceptable tensor fields with length $\sigma +u+1$, a factor
$\nabla\omega$ (which we treat as a factor
$\nabla\phi$)\footnote{In particular this factor does not contain a free index.}
 and a factor $\nabla^{(A)}\phi_{u+1}$, $A\ge 2$
(which we treat as a factor $\nabla^{(A)}\Omega_{p+1}$), with the
 restriction that the $u$-simple character of $C^{u,i_1\dots
i_a}_{g}$ arises from $\vec{\kappa}_{simp}$ by replacing one
factor $T=\nabla^{(m)}R_{ijkl}$ or $T=S_{*}\nabla^{(\nu)}R_{ijkl}$
by a factor $T'=\nabla^{(A)}\Omega_{p+1}$ where all the factors
$\nabla\phi_h$ that contracted against $T$ now contract against
$T'$; moreover if $a=\mu$ and the factor
$\nabla^{(A)}\phi_{u+1}$ has $A=2$ and is contracting against a factor $\nabla\phi_h$
 then we additionnaly require that it does not contain free indices.

\par In addition, we will denote by
$$\Sum_{u\in U^\sharp_1} a_u
C^{u,i_1\dots i_a}_{g}(\Omega_1,\dots ,\Omega_p,\phi_1, \dots
,\phi_u,\phi_{u+1},\omega),$$

$$\Sum_{u\in U^\sharp_2} a_u
C^{u,i_1\dots i_a}_{g}(\Omega_1,\dots ,\Omega_p,\phi_1, \dots
,\phi_u,\phi_{u+1},\omega)$$ generic linear combinations as above
{\it with} one un-acceptable factor of
 the form $\nabla\Omega_h$, $h\le p$, (with only one derivative) contracting against a factor
$\nabla\phi_{u+1}$ or $\nabla\omega$, respectively.

We will denote by
$$\Sum_{j\in J^{\sigma +u+1}} a_j C^j_{g}
(\Omega_1,\dots ,\Omega_p,\phi_1,\dots ,\phi_u,\phi_{u+1},
\omega)$$ a generic linear combination of complete contractions
with length $\sigma +u+1$ and two factors either
$\nabla^{(A)}\phi_{u+1},\nabla\omega$ or
$\nabla^{(A)}\omega,\nabla\phi_{u+1}$
 ($A\ge 2$) and  at least one factor $\nabla\phi_f,
f\in Def(\vec{\kappa}_{simp})$\footnote{Recall that
$Def(\vec{\kappa}_{simp})$ stands for the set of numbers $o$ for
which $\nabla\tilde{\phi}_o$ is contracting against the index
${}_i$ in some factor $S_{*}\nabla^{(\nu)}R_{ijkl}$.} contracting
against a
 derivative
 index in a factor $\nabla^{(m)}R_{ijkl}$ {\it or} contracting against
 one of the first $A-2$ indices in a factor
$\nabla^{(A)}\phi_{u+1}$ or $\nabla^{(A)}\omega$.
\newline

\par For future reference, we will also put down some definitions
regarding complete contractions of length $\sigma +u+2$ in
$Sub_\omega \{ Image^{1,\beta}_{\phi_{u+1}} [L_{g}]\}$ (with a
factor $\nabla\phi_{u+1}$ and  a factor $\nabla\omega$): We will
denote by

$$\Sum_{m\in M} a_m C^{m,i_1\dots i_{a+2}}_{g}
(\Omega_1,\dots ,\Omega_p,\phi_1, \dots
,\phi_u)\nabla_{i_{a+1}}\phi_{u+1}\nabla_{i_{a+2}}\omega$$ (where
$a\ge \mu$) a generic linear combination of {\it acceptable}
tensor fields of length $\sigma +u+2$ with a $u$-simple character
$\vec{\kappa}_{simp}$ (this $u$-simple character only encodes
information on the factors $\nabla\phi_h, h\le u$).

\par We also denote by $$\Sum_{m\in M^\sharp} a_m C^{m,i_1\dots i_{a+2}}_{g}
(\Omega_1,\dots ,\Omega_p,\phi_1, \dots
,\phi_u)\nabla_{i_{a+1}}\phi_{u+1}\nabla_{i_{a+2}}\omega$$ a
generic linear combination of tensor fields as above, with {\it one}
un-acceptable factor of the form $\nabla_{i_b}\Omega_f$ where
${}_{i_b}$ is a free index and in fact $b=a+1$ {\it or} $b=a+2$.

We will denote by
$$\Sum_{j\in J^{\sigma +u+2}} a_j C^j_{g}
(\Omega_1,\dots ,\Omega_p,\phi_1,\dots ,\phi_u,\phi_{u+1},
\omega)$$ a generic linear combination of complete contractions
with length $\sigma +u+2$ and two factors
$\nabla\phi_{u+1},\nabla\omega$ and  at least one factor
$\nabla\phi_f, f\in Def(\vec{\kappa}_{simp})$ contracting
 against a derivative
 index in a factor $\nabla^{(m)}R_{ijkl}$.
\newline

\par Finally, we will denote by
$$\Sum_{w\in W} a_w
C^w_{g}(\Omega_1,\dots ,\Omega_p,\phi_1, \dots
,\phi_u,\phi_{u+1},\omega)$$ a generic linear combination of
complete contractions that have length $\ge\sigma +u+1$ and
 two factors $\nabla^{(A)}\phi_{u+1},\nabla^{(B)}\omega$ with
$A,B\ge 2$.\footnote{Notice that for length $\sigma+u+1$,
the sublinear combination $\sum_{w\in W}\dots$
corresponds exactly to the generic linear combination $\sum_{t\in T_1}\dots$
in (\ref{baldouin'}).}
\end{definition}

Now, return to (\ref{mikulski}). We wish to analyze the terms in
this equation. We firstly observe that:

\begin{equation}
\label{middle} \Sum_{j\in J} a_j
Rep[C^j_{g}(\Omega_1,\dots,\Omega_p,\phi_1,\dots ,\phi_u)]=
\Sum_{j\in J^{\sigma+u}} a_j
C^j_{g}(\Omega_1,\dots,\Omega_p,\phi_1,\dots
,\phi_u,\phi_{u+1},\omega).
\end{equation}

Now, we claim that we can write out:

\begin{equation}
\label{mikul}
\begin{split}
&\Sum_{l\in L} a_l Xdiv_{i_1}\dots Xdiv_{i_a} Rep[C^{l,i_1,\dots
,i_a}_{g}(\Omega_1,\dots ,\Omega_p,\phi_1,\dots ,\phi_u)]=
\\& \Sum_{u\in U_1} a_u Xdiv_{i_1}\dots Xdiv_{i_a}
C^{u,i_1\dots ,i_a}_{g}(\Omega_1,\dots ,\Omega_p,\phi_1, \dots
,\phi_u,\phi_{u+1},\omega)+
\\&\Sum_{u\in U_2} a_u Xdiv_{i_1}\dots Xdiv_{i_a}
C^{u,i_1\dots ,i_a}_{g}(\Omega_1,\dots ,\Omega_p,\phi_1, \dots
,\phi_u,\phi_{u+1},\omega)+
\\&\Sum_{j\in J^{\sigma +u}} a_j
C^j_{g}(\Omega_1,\dots ,\Omega_p,\phi_1, \dots
,\phi_u,\phi_{u+1},\omega)+
\\& \Sum_{w\in W} a_w
C^w_{g}(\Omega_1,\dots ,\Omega_p,\phi_1, \dots
,\phi_u,\phi_{u+1},\omega);
\end{split}
\end{equation}
(this holds perfectly, {\it not} modulo complete contractions of
greater length).

\par If we can prove this, we can then replace into (\ref{baldouin})  to
derive an equation:

\begin{equation}
\label{gonnahappen}
\begin{split}
&Sub^{\sigma+u}_\omega[L_{g}]=\Sum_{u\in U_1} a_u Xdiv_{i_1}\dots
Xdiv_{i_a} C^{u,i_1\dots ,i_a}_{g}(\Omega_1,\dots
,\Omega_p,\phi_1, \dots ,\phi_u,\phi_{u+1},\omega)+
\\&\Sum_{u\in U_2} a_u Xdiv_{i_1}\dots Xdiv_{i_a}
C^{u,i_1\dots ,i_a}_{g}(\Omega_1,\dots ,\Omega_p,\phi_1, \dots
,\phi_u,\phi_{u+1},\omega)+
\\&\Sum_{j\in J^{\sigma +u}_1} a_j
C^j_{g}(\Omega_1,\dots ,\Omega_p,\phi_1, \dots
,\phi_u,\phi_{u+1},\omega)+
\\& \Sum_{w\in W} a_w
C^w_{g}(\Omega_1,\dots ,\Omega_p,\phi_1, \dots
,\phi_u,\phi_{u+1},\omega),
\end{split}
\end{equation}
(the above holds perfectly), where the tensor fields indexed in $J^{\sigma+u}_1$ are of
 the same general form as the tensor fields we generically
index in $J^{\sigma+u}$. In other words, provided we can show
(\ref{mikul}) we derive that modulo complete contractions of
length $\ge\sigma +u+1$:

\begin{equation}
\label{suffering} \Sum_{j\in J^{\sigma +u}_1} a_j
C^j_{g}(\Omega_1,\dots ,\Omega_p,\phi_1, \dots
,\phi_u,\phi_{u+1},\omega)=0,
\end{equation}
(modulo length $\ge\sigma+u+1$).

\par Then, using the above we will show that we can write:

\begin{equation}
\label{tzieria}
\begin{split}
&\Sum_{j\in J^{\sigma +u}_1} a_j C^j_{g}(\Omega_1,\dots
,\Omega_p,\phi_1, \dots ,\phi_u,\phi_{u+1},\omega)=
\\&\Sum_{j\in J^{\sigma +u+1}} a_j
C^j_{g}(\Omega_1,\dots ,\Omega_p,\phi_1, \dots
,\phi_u,\phi_{u+1},\omega)\\&+ \Sum_{w\in W} a_w
C^w_{g}(\Omega_1,\dots ,\Omega_p,\phi_1, \dots
,\phi_u,\phi_{u+1},\omega).
\end{split}
\end{equation}
(The above holds exactly, not modulo longer complete
 contractions). Moreover, the left hand side is not generic
 notation, but stands for the same sublinear combination in
(\ref{mikul}). The right hand side consists of generic 
linear combinations as defined in definition \ref{leaves}.

\subsection{ Further steps in the analysis of $Image^{1,\beta}_{\phi_{u+1}}[L_g]=0$: 
Proof of (\ref{mikul}) and (\ref{tzieria}):}

\par We first show (\ref{mikul}).
We will prove this equation using the inductive assumption on
Proposition \ref{giade} and the usual inductive argument: We 
will  break up
(\ref{reno3}) into sublinear combinations that have the same
$u$-weak character (here the weak character also 
takes into account the two new factors 
 $\Omega_{p+1}=\phi_{u+1}, \Omega_{p+2}=\omega$),
and then inductively apply Corollary 1 from \cite{alexakis4} and
at each step we convert the $\nabla\upsilon$'s in $Xdiv$'s.

\par Specifically: For each $C^{l,i_1\dots
i_a}_{g}$, $C^j_{g}$ we divide their curvature factors
\\$(\nabla^{(m)}R_{ijkl},S_{*}\nabla^{(\nu)} R_{ijkl}$) into
{\it categories}: $K_1,\dots,K_d$. We decide that two curvature factors
$T,T'$ in two complete contractions or tensor fields in
(\ref{mikul}) belong to the same category if they are contracting
against the same factors $\nabla\phi_h$. We also decide that the
curvature factors that do not contract against any factors
$\nabla\phi_h$ belong to the last category $K_d$.

\par We then recall that each tensor field and each complete
contraction in (\ref{mikul}) has arisen by replacing one curvature
factor $\nabla^{(v)}R_{ijkl}$ by
$\nabla^{(v)}[\nabla^{(2)}_{ab}\phi_{u+1}]$ and one other
curvature factor $\nabla^{(c)}R_{i'j'k'l'}$ by
$\nabla^{(c)}[\nabla^{(2)}_{a'b'}\omega]$. We then index the
tensor fields and complete contractions in (\ref{mikul}) in the
sets $L^{\alpha,\beta}, J^{\alpha,\beta}$ ($1\le\alpha,\beta\le
d$) according to the rule that a tensor field or complete
contraction in (\ref{mikul}) belongs to $L^{\alpha,\beta},
J^{\alpha,\beta}$ if and only if it has arisen by replacing a
curvature factor that belongs to the category $K_\alpha$ by
$\nabla^{(c+2)}_{ij}\phi_{u+1}$ and a curvature factor that belongs
to $K_\beta$ by $\nabla^{(v+2)}_{ij}\omega$.

\par We then see that (\ref{mikul}) can be re-written in the form:

\begin{equation}
\label{mpourzoua} \Sum_{1\le\alpha,\beta\le d}
\{\Sum_{L^{\alpha,\beta}} a_l Xdiv_{i_1}\dots
Xdiv_{i_a}C^{l,i_1\dots i_a}_{g}+\Sum_{j\in J^{\alpha,\beta}} a_j
C^j_{g}\}=0,
\end{equation}
modulo complete contractions of length $\ge\sigma +u+1$. Then,
since the above must hold formally and hence sublinear
combinations with different weak characters must vanish
separately, we derive that for each pair $1\le \alpha,\beta\le d$:

\begin{equation}
\label{mpourzoua2} \Sum_{l\in L^{\alpha,\beta}} a_l
Xdiv_{i_1}\dots Xdiv_{i_a}C^{l,i_1\dots i_a}_{g}+\Sum_{j\in
J^{\alpha,\beta}} a_j C^j_{g}=0,
\end{equation}
(modulo complete contractions of length $\ge\sigma +u+1$).

\par We will then show our claim (\ref{mikul}) for each of the
index sets $L^{\alpha,\beta}$ separately. We distinguish three
cases: Either both the categories $K_a,K_b$ represent a
``genuine'' factor $\nabla^{(m)}R_{ijkl}$ {\it or} one of them
represents such a ``genuine'' factor and the  other represents a factor
$S_{*}\nabla^{(\nu)} R_{ijkl}$ {\it or} both represent factors
$S_{*}\nabla^{(\nu)} R_{ijkl}$. We denote these three cases by
$(i),(ii),(iii)$ respectively.
\newline

{\it Proof of (\ref{mikul}) in case (i):} We observe that all
 tensor fields indexed in $L^{\alpha,\beta}$ have the {\it same}
$u$-simple character, say $\tilde{\kappa}_{simp}$ (we are treating
the factors $\nabla^{(u)}\phi_{u+1}, \nabla^{(y)}\omega$ as
functions $\Omega_{p+1},\Omega_{p+2}$). Moreover, each of the
tensor fields indexed in $L^{\alpha,\beta}$ has the property that any factors
$\nabla\phi_h$ that are contracting against the factors
$\nabla^{(u)}\phi_{u+1}, \nabla^{(y)}\omega$ will be contracting against
one of the left-most  $u-2$ or $y-2$ indices. This follows by
the definition above. We also note that all the complete
contractions indexed in $J^{\alpha,\beta}$ are simply subsequent
to $\tilde{\kappa}_{simp}$. We denote by $\tau$ the minimum rank
of the tensor fields appearing in $L^{\alpha,\beta}$ (by
hypothesis $\tau\ge \mu$). We index those tensor fields in the set
$L^{\alpha,\beta|\tau}$. Therefore, (after using the Eraser, defined in the 
Appendix of \cite{alexakis1}, if
necessary) we
 may apply the inductive assumption of Corollary
1 in \cite{alexakis4}\footnote{Observe that (\ref{mpourzoua2}) falls
under the inductive assumption of Corollary 1 in \cite{alexakis4},
because the tensor fields there have length $\sigma+u$ and $p+2$
factors $\nabla^{(y)}\Omega_h$. Moreover observe that here is no danger of falling under
a ``forbidden'' case of Corollary 1 in \cite{alexakis4}, since all the $\mu$-tensor fields in
(\ref{hypothese2}) have no special free indices, thus there will be no special
free indices among the tensor fields of minimum rank in (\ref{mpourzoua2}).} to derive that there exists a
linear combination of acceptable tensor fields, $\Sum_{h\in
H^{\alpha,\beta}} a_h C^{h,i_1\dots i_{\tau+1}}_{g}$
 of $u$-simple character $\tilde{\kappa}_{simp}$ and with
the additional restriction that  any factor $\nabla\phi_h$ that is
contracting against $\nabla^{(u)}\phi_{u+1}$ or
$ \nabla^{(y)}\omega$ is contracting against one of the
  $u-2$ or $y-2$ leftmost indices and moreover with
rank $\tau+1$ so that:

\begin{equation}
\label{mpourzoua3} \begin{split} & \Sum_{l\in
L^{\alpha,\beta|\tau}} a_l C^{l,i_1\dots
i_\tau}_{g}\nabla_{i_1}\upsilon\dots \nabla_{i_\tau}\upsilon-
\Sum_{h\in H^{\alpha,\beta}} a_h Xdiv_{i_{\tau+1}}C^{h,i_1\dots
i_{\tau+1}}_{g}\nabla_{i_1}\upsilon\dots \nabla_{i_\tau}\upsilon=
\\&\Sum_{j\in J^{\sigma+u}} a_j C^{j,i_1\dots
i_\tau}_{g}\nabla_{i_1}\upsilon\dots \nabla_{i_\tau}\upsilon,
\end{split}
\end{equation}
modulo complete contractions of length $\ge\sigma +u+1+\tau$. Here
the tensor fields indexed in $J^{\sigma+u}$ have the same
properties as the contractions indexed in $J^{\sigma+u}$.

\par In fact since the above holds formally, by just paying attention to the correction terms
of greater length that arise in (\ref{mpourzoua3}) we derive:

\begin{equation}
\label{mpourzoua4} \begin{split} & \Sum_{l\in
L^{\alpha,\beta|\tau}} a_l C^{l,i_1\dots
i_\tau}_{g}\nabla_{i_1}\upsilon\dots \nabla_{i_\tau}\upsilon-
\Sum_{h\in H^{\alpha,\beta}} a_h Xdiv_{i_{\tau+1}}C^{h,i_1\dots
i_{\tau+1}}_{g}\nabla_{i_1}\upsilon\dots \nabla_{i_\tau}\upsilon
\\&= \Sum_{u\in U_1} a_u
C^{u,i_1\dots ,i_\tau}_{g}(\Omega_1,\dots ,\Omega_p,\phi_1, \dots
,\phi_u,\phi_{u+1},\omega)\nabla_{i_1}\upsilon\dots
\nabla_{i_\tau}\upsilon+
\\&\Sum_{u\in U_2} a_u
C^{u,i_1\dots ,i_\tau}_{g}(\Omega_1,\dots ,\Omega_p,\phi_1, \dots
,\phi_u,\phi_{u+1},\omega)\nabla_{i_1}\upsilon\dots
\nabla_{i_\tau}\upsilon+
\\&\Sum_{j\in J^{\sigma +u}} a_j
C^{j,i_1\dots i_\tau}_{g}(\Omega_1,\dots ,\Omega_p,\phi_1, \dots
,\phi_u,\phi_{u+1},\omega)\nabla_{i_1}\upsilon\dots
\nabla_{i_\tau}\upsilon+
\\& \Sum_{w\in W} a_w
C^{w,i_1\dots i_\tau}_{g}(\Omega_1,\dots ,\Omega_p,\phi_1, \dots
,\phi_u,\phi_{u+1},\omega)\nabla_{i_1}\upsilon\dots
\nabla_{i_\tau}\upsilon,
\end{split}\end{equation}
modulo complete contractions of length $\ge\sigma +u+2+\tau$. Here
the tensor fields indexed in $W$ have the same properties as the
complete contractions indexed in $W$: They have length
$\sigma+u+1$ but also have two factors $\nabla^{(A)}\phi_{u+1},
\nabla^{(B)}\omega$, $A,B\ge 2$.

\par Therefore, making the $\nabla\upsilon$'s 
into $Xdiv$'s\footnote{(We are using the last Lemma from the Appendix in \cite{alexakis1}).} we
derive an equation:

\begin{equation}
\label{mpourzoua5} \begin{split} & Xdiv_{i_1}\dots
Xdiv_{i_\tau}\Sum_{l\in L^{\alpha,\beta|\tau}} a_l C^{l,i_1\dots
i_\tau}_{g}- \Sum_{h\in H^{\alpha,\beta}} a_h Xdiv_{i_1}\dots
Xdiv_{i_\tau}Xdiv_{i_{\tau+1}}
\\&C^{h,i_1\dots i_{\tau+1}}_{g}= Xdiv_{i_1}\dots
Xdiv_{i_\tau}\Sum_{u\in U_1} a_u C^{u,i_1\dots
,i_\tau}_{g}(\Omega_1,\dots ,\Omega_p,\phi_1, \dots
,\phi_u,\phi_{u+1},\omega)+
\\&\Sum_{u\in U_2} a_u
Xdiv_{i_1}\dots Xdiv_{i_\tau}C^{u,i_1\dots
,i_\tau}_{g}(\Omega_1,\dots ,\Omega_p,\phi_1, \dots
,\phi_u,\phi_{u+1},\omega)+
\\&\Sum_{j\in J^{\sigma +u}} a_j
C^j_{g}(\Omega_1,\dots ,\Omega_p,\phi_1, \dots
,\phi_u,\phi_{u+1},\omega)+
\\& \Sum_{w\in W} a_w
C^{w}_{g}(\Omega_1,\dots ,\Omega_p,\phi_1, \dots
,\phi_u,\phi_{u+1},\omega).
\end{split}\end{equation}
Thus, substituting into (\ref{mpourzoua2}) and inductively
 repeating this step,\footnote{At the very last step of this inductive argument,
 we may have to apply the ``weak substitute'' of Proposition \ref{giade}, 
 from the Appendix of \cite{alexakis4}. In that case our result will follow
 since in that case the  minimum rank among those terms will be $>\mu$.} we derive our claim in this
first case.
\newline

{\it Proof of (\ref{mikul}) in the case (ii):} Now, in the second
case we assume with no loss of generality that $K_\alpha$
corresponds to a factor $\nabla^{(m)}R_{ijkl}$ in
$\vec{\kappa}_{simp}$ and $K_\beta$ corresponds to a factor
$S_{*}\nabla^{(\nu)} R_{ijkl}\nabla^i\tilde{\phi}_f$ in
$\vec{\kappa}_{simp}$.
 We again observe that all the tensor fields
$C^{l,i_1\dots i_a}_{g}(\Omega_1,\dots ,\Omega_p,\phi_1,\dots
,\phi_{u+1},\omega)$ have the same $u$-simple character
$\vec{\kappa}_{simp}'$ and any factor $\nabla\phi_h$ that is
contracting against the factor $\nabla^{(p)}\phi_{u+1}$ in any
tensor field $C^{l,i_1\dots i_a}_{g}$ in (\ref{mpourzoua2}) must
be contracting against
 one of the first $p-2$ indices there.

 Then, for each vector field
$C^{l,i_1\dots i_a}_{g}, l\in L^{\alpha,\beta}$ we inquire whether
the factor $\nabla^{(p)}_{r_1\dots r_p}\omega$ (for which 
 ${}_{r_{p-1}}$ is contracting against the factor
$\nabla\tilde{\phi}_f$) has $p>2$ or $p=2$. In the first case, we
straightforwardly observe that we can write:

\begin{equation}
\label{grafw}\begin{split} & Xdiv_{i_1}\dots
Xdiv_{i_a}C^{l,i_1\dots i_a}_{g}(\Omega_1,\dots ,\Omega_p,\phi_1,
\dots ,\phi_u,\phi_{u+1},\omega)=
\\&\Sum_{j\in
J^{\sigma+u}} a_j C^j_{g}(\Omega_1,\dots ,\Omega_p,\phi_1, \dots
,\phi_u,\phi_{u+1},\omega)+
\\&\Sum_{u\in U_2} a_u
Xdiv_{i_1}\dots Xdiv_{i_a}C^{u,i_1\dots ,i_a}_{g}(\Omega_1,\dots
,\Omega_p,\phi_1, \dots ,\phi_u,\phi_{u+1},\omega)
\\& +\Sum_{w\in W} a_w C^w_{g}(\Omega_1,\dots ,\Omega_p,
\phi_1, \dots,\phi_u,\phi_{u+1},\omega);
\end{split}
\end{equation}
(the above holds perfectly, not modulo terms of greater length).
This just follows by applying the curvature identity 
to ${}_{r_{p-2}},{}_{r_{p-1}}$ in $\nabla^{(p)}_{r_1\dots
r_p}\omega$.
Thus, we may assume with no loss of generality that
all the tensor fields indexed in $L^{\alpha,\beta}$ have a factor
$\nabla^{(2)}\omega$.

\par By the same argument, we may also assume that all complete contractions indexed in
$J^{\alpha,\beta}$ either have a factor $\nabla^{(p)}\omega, p\ge
3$ or a factor $\nabla^{(2)}\omega$ but one of the factors $\nabla
\phi_h$, $h\in Def(\vec{\kappa}_{simp})$ is contracting against a
derivative index in some factor $\nabla^{(m)}R_{ijkl}$.
Accordingly, we break up $J^{\alpha,\beta}$ into
$J^{\alpha,\beta,I}$, $J^{\alpha,\beta,II}$.

\par Therefore, picking out the sublinear combination in
(\ref{mikulski}) with a factor $\nabla^{(2)}\omega_1$, we derive
 the equation:

\begin{equation}
\label{grafw2}
\begin{split}
& \Sum_{l\in L^{\alpha,\beta}} a_l X_{*}div_{i_1}\dots
X_{*}div_{i_a}C^{l,i_1\dots i_a}_{g}(\Omega_1,\dots
,\Omega_p,\phi_1, \dots ,\phi_u,\phi_{u+1},\omega)+
\\&\Sum_{j\in J^{\alpha,\beta,II}} a_j
C^j_{g}(\Omega_1,\dots ,\Omega_p,\phi_1, \dots
,\phi_u,\phi_{u+1},\omega)=0,
\end{split}
\end{equation}
where $X_{*}div_i$ stands for the sublinear combination in
$Xdiv_i$ where $\nabla_i$ is not allowed to hit the factor
$\nabla^{(2)}\omega$. Furthermore, we may assume that for each of
the tensor fields above, the factor $\nabla^{(2)}\omega$ does not
contain a free index, yet the rank of all the tensor fields is
$\ge \mu$: This assumption can be made with no loss of generality
since if a tensor field in (\ref{grafw2}) has $a=\mu$ then it will
not contain a free index in $\nabla^{(2)}\omega$ by definition, while if $a>\mu$, we may
just {\it neglect} the $X_{*}div_i$, where ${}_i$ is the (unique)
free index belonging to $\nabla^{(2)}\omega_1$--the resulting
iterated $Xdiv$ will still have rank $\ge\mu$.

\par Then, by using the eraser onto the factor $\nabla\phi_h$ that contracts
against
$\nabla^{(2)}\omega$ and applying our inductive assumption of
Corollary 2 from \cite{alexakis4}\footnote{Observe that since $\nabla\omega$ does
not contain free indices we do not have to worry about the ``forbidden cases''.}
 (or Lemma 4.7 if $\sigma=3$),
 we derive that we can write:

\begin{equation}
\label{grafw3} \begin{split} &\Sum_{l\in L^{\alpha,\beta}} a_l
 Xdiv_{i_1}\dots
Xdiv_{i_a}C^{l,i_1\dots i_a}_{g}(\Omega_1,\dots ,\Omega_p,\phi_1,
\dots ,\phi_u,\phi_{u+1},\omega)=
\\&\Sum_{l\in L'^{\alpha,\beta}} a_l
 Xdiv_{i_1}\dots
Xdiv_{i_a}C^{l,i_1\dots i_a}_{g}(\Omega_1,\dots ,\Omega_p,\phi_1,
\dots ,\phi_u,\phi_{u+1},\omega)+
\\&\Sum_{j\in J^{\alpha,\beta,II}} a_j
C^j_{g}(\Omega_1,\dots ,\Omega_p,\phi_1, \dots
,\phi_u,\phi_{u+1},\omega),
\end{split}
\end{equation}
where the tensor fields indexed in $L'^{\alpha,\beta}$ have a
factor $\nabla^{(3)}_{ijk}\omega\nabla^j\tilde{\phi}_f$. Also
\\$\Sum_{j\in J^{\alpha,\beta,II}} a_j \dots$ on the RHS is generic
 notation. Thus, using the identity \\$\nabla^{(3)}_{ijk}\omega\nabla^j\tilde{\phi}_f=
 \nabla^{(3)}_{jik}\omega\nabla^j\tilde{\phi}_f+R_{ijlk}\nabla^l\omega_1\nabla^k\tilde{\phi}_f$
 we derive our claim.
\newline

{\it Proof of (\ref{mikul}) in case $(iii)$:} Finally, the last
case, where both $K_\alpha,K_\beta$ are in the form
$S_{*}\nabla^{(\nu)}R_{ijkl}$:  As before, we consider any tensor
field $C^{l,i_1\dots i_a}_{g}$ in (\ref{mikulski}) and we note
that it must contain two expressions $\nabla^{(y)}_{r_1\dots
r_y}\phi_{u+1}\nabla^{r_{y-1}}\tilde{\phi}_b$ and
$\nabla^{(p)}_{t_1\dots t_p}\omega\nabla^{t_{p-1}}\tilde{\phi}_c$
with $y,p\ge 2$. Moreover all tensor fields have the same
$u$-simple character which we denote by $\tilde{\kappa}_{simp}$.
Furthermore the last index in both these factors is not
contracting against a factor $\nabla\phi_h$.  If either of the
numbers $y,p$ is $>2$, we then apply the curvature identity to derive
that we can write:

\begin{equation}
\label{groove}
\begin{split} & Xdiv_{i_1}\dots
Xdiv_{i_a}C^{l,i_1\dots i_a}_{g}(\Omega_1,\dots ,\Omega_p,\phi_1,
\dots ,\phi_u,\phi_{u+1},\omega)=
\\&\Sum_{j\in
J^{\sigma+u}} a_j C^j_{g}(\Omega_1,\dots ,\Omega_p,\phi_1, \dots
,\phi_u,\phi_{u+1},\omega)+
\\&\Sum_{u\in U_1} a_u
Xdiv_{i_1}\dots Xdiv_{i_a}C^{u,i_1\dots ,i_a}_{g}(\Omega_1,\dots
,\Omega_p,\phi_1, \dots ,\phi_u,\phi_{u+1},\omega)+
\\&\Sum_{u\in U_2} a_u
Xdiv_{i_1}\dots Xdiv_{i_a}C^{u,i_1\dots ,i_a}_{g}(\Omega_1,\dots
,\Omega_p,\phi_1, \dots ,\phi_u,\phi_{u+1},\omega)
\\& +\Sum_{w\in W}
a_w C^w_{g}(\Omega_1,\dots ,\Omega_p,\phi_1, \dots
,\phi_u,\phi_{u+1},\omega).
\end{split}
\end{equation}

\par Therefore, we may now prove our claim under the assumption
that all tensor fields $C^{l,i_1\dots ,i_a}_{g}$, $l\in
L^{\alpha,\beta}$ have two factors $\nabla^{(2)}\phi_{u+1},
\nabla^{(2)}\omega$, with the indices before last contracting
against factors $\nabla\phi_h$, $h\in Def(\vec{\kappa}_{simp})$.

\par Analogously, we again divide $J^{\alpha,\beta}$ into two subsets.
 We say $j\in J^{\alpha,\beta,I}$ if at least
one of the factors $\nabla^{(u)}\phi_{u+1}$ or
$\nabla^{(y)}\omega$ has $u\ge 3$ or $y\ge 3$. We say $j\in
J^{\alpha,\beta,II}$ if we have two factors
$\nabla^{(2)}\phi_{u+1}$ and $\nabla^{(2)}\omega$. In this second
case we see that by definition we must have at least one factor
$\nabla\phi_h$,
 $h\in Def(\vec{\kappa}_{simp})$ contracting against a
 derivative index in some factor $\nabla^{(m)}R_{ijkl}$.

\par Now, picking out the tensor fields in (\ref{mikulski}) with
two factors $\nabla^{(2)}\phi_{u+1},\nabla^{(2)}\omega$ we derive
that:

\begin{equation}
\label{grafw3}
\begin{split}
& \Sum_{l\in L^{\alpha,\beta}} a_l X_{*}div_{i_1}\dots
X_{*}div_{i_a}C^{l,i_1\dots i_a}_{g}(\Omega_1,\dots
,\Omega_p,\phi_1, \dots ,\phi_u,\phi_{u+1},\omega)+
\\&\Sum_{j\in J^{\alpha,\beta,II}} a_j C^j_{g}
(\Omega_1,\dots ,\Omega_p,\phi_1, \dots
,\phi_u,\phi_{u+1},\omega)=0,
\end{split}
\end{equation}
where here $X_{*}div_i$ means $\nabla_i$ is not allowed to hit the
factor $\nabla^{(2)}\phi_{u+1}$ nor the factor
$\nabla^{(2)}\omega$.

\par We then claim that we can write:

\begin{equation}
\label{grafw4} \begin{split}& \Sum_{l\in L^{\alpha,\beta}} a_l
Xdiv_{i_1}\dots Xdiv_{i_a}C^{l,i_1\dots i_a}_{g} (\Omega_1,\dots
,\Omega_p,\phi_1, \dots ,\phi_u,\phi_{u+1},\omega)=
\\&\Sum_{l\in L'^{\alpha,\beta}} a_l
Xdiv_{i_1}\dots Xdiv_{i_a}C^{l,i_1\dots i_a}_{g}(\Omega_1,\dots
,\Omega_p,\phi_1, \dots ,\phi_u,\phi_{u+1},\omega)+
\\& \Sum_{j\in
J^{\sigma+u}} a_j C^j_{g}(\Omega_1,\dots ,\Omega_p,\phi_1, \dots
,\phi_u,\phi_{u+1},\omega),
\end{split}
\end{equation}
where each tensor field indexed in $L'^{\alpha,\beta}$ has a
simple character $\tilde{\kappa}_{simp}$ and either an expression
$\nabla^{(3)}_{ijk}\phi_{u+1}\nabla^j\tilde{\phi}_h$ or an
expression $\nabla^{(3)}_{ijk}\omega\nabla^j\tilde{\phi}_h$ and
has $a\ge \mu$. If we can do this, then repeating the curvature
identity as above we can derive our claim (\ref{mikul}) in case
$(iii)$.
\newline

\par The proof of equation (\ref{grafw4}) 
is rather technical and the methods used there are  
 not relevant to our further study of the 
sublinear combination \\$Image^{1,\beta}_{\phi_{u+1}}[L_g]=0$.
Thus, in order not to distract the reader 
from the main points of the argument, we 
will present the proof of (\ref{grafw4}) 
in the Mini-Appendix \ref{pfgrafw4} below. 
\newline

{\bf Proof of (\ref{tzieria}):}  (\ref{tzieria}) follows from
(\ref{suffering}) by the usual argument where we make the
linearized complete contractions hold formally and then repeat the
permutations to the non-linearized setting: We have that

$$\Sum_{j\in J^{\sigma +u}_1} a_j
C^j_{g}(\Omega_1,\dots ,\Omega_p,\phi_1, \dots
,\phi_u,\phi_{u+1},\omega)=0,$$ modulo longer complete
contractions, so the above holds
 formally at the linearized level, so repeating the
 permutations to the non-linearized level we get correction
 terms of the desired form.
 $\Box$
\newline

\par In conclusion, using (\ref{mikul}) and (\ref{tzieria}), we can replace the
first 
 sublinear combination $Sub_\omega[\dots]$
in (\ref{baldouin'}) to obtain a new equation:
\begin{equation}
\label{baldouin''}
\begin{split}
& 0=\Sum_{u\in U_1} a_u Xdiv_{i_1}\dots Xdiv_{i_a}
C^{u,i_1\dots ,i_a}_{g}(\Omega_1,\dots ,\Omega_p,\phi_1, \dots
,\phi_u,\phi_{u+1},\omega)+
\\&\Sum_{u\in U_2} a_u Xdiv_{i_1}\dots Xdiv_{i_a}
C^{u,i_1\dots ,i_a}_{g}(\Omega_1,\dots ,\Omega_p,\phi_1, \dots
,\phi_u,\phi_{u+1},\omega)+
\\&\Sum_{j\in J^{\sigma +u}} a_j
C^j_{g}(\Omega_1,\dots ,\Omega_p,\phi_1, \dots
,\phi_u,\phi_{u+1},\omega)+
\\& \Sum_{w\in W} a_w
C^w_{g}(\Omega_1,\dots ,\Omega_p,\phi_1, \dots
,\phi_u,\phi_{u+1},\omega)+
Sub^{\sigma+u+1}_\omega \{
Image^{1,\beta}_{\phi_{u+1}}[L_{g}]\}+
\\& Sub^{\sigma+u+2}_\omega \{
Image^{1,\beta}_{\phi_{u+1}}[L_{g}]\}+
 \Sum_{t\in T^{\sigma+u+1}\bigcup T^{\sigma+u+2}} a_t
C^t_{g}(\Omega_1,\dots , \Omega_p,\phi_1,\dots
,\phi_{u+1},\omega).
\end{split}
\end{equation}
\par In particular, since the minimum length of the complete
 contractions in the above is
$\sigma+u+1$, if we denote by $F_g$
 the sublinear combination of terms with $\sigma+u+1$ factors with two factors
$\nabla^{(A)}\phi_{u+1},\nabla^{(B)}\omega$, then $lin\{ F_g\}=0$
formally. Now, notice that the sublinear combination $F_g$ is in
fact:
\begin{equation}
 \begin{split}
&F_g=\Sum_{w\in W} a_w
C^w_{g}(\Omega_1,\dots ,\Omega_p,\phi_1, \dots
,\phi_u,\phi_{u+1},\omega)+\\&
\Sum_{t\in T^{\sigma+u+1}} a_t C^t_{g}(\Omega_1,\dots ,
\Omega_p,\phi_1,\dots ,\phi_{u+1},\omega).
\end{split}\end{equation}

\par Then, using the fact that $lin\{F_g\}=0$ formally,
we repeat the permutations that make the LHS of this
 equation zero to the {\it non-linear} setting, to derive:

\begin{equation}
\label{reconquise} \begin{split} &F_g=
\Sum_{t\in T^{\sigma+u+2}} a_t
C^t_{g}(\Omega_1,\dots , \Omega_p,\phi_1,\dots ,\phi_{u+1},\omega).
\end{split}
\end{equation}
(This equation holds perfectly, not ``modulo longer terms''),
using generic notation on the RHS. Therefore, plugging the above
into (\ref{baldouin''}), we derive a new equation:

\begin{equation}
\label{baldouin'''}
\begin{split}
& \Sum_{u\in U_1} a_u Xdiv_{i_1}\dots Xdiv_{i_a}
C^{u,i_1\dots ,i_a}_{g}(\Omega_1,\dots ,\Omega_p,\phi_1, \dots
,\phi_u,\phi_{u+1},\omega)+
\\&\Sum_{u\in U_2} a_u Xdiv_{i_1}\dots Xdiv_{i_a}
C^{u,i_1\dots ,i_a}_{g}(\Omega_1,\dots ,\Omega_p,\phi_1, \dots
,\phi_u,\phi_{u+1},\omega)+
\\&\Sum_{j\in J^{\sigma +u}} a_j
C^j_{g}(\Omega_1,\dots ,\Omega_p,\phi_1, \dots
,\phi_u,\phi_{u+1},\omega)+
Sub^{\sigma+u+1}_\omega \{
Image^{1,\beta}_{\phi_{u+1}}[L_{g}]\}+
\\& Sub^{\sigma+u+2}_\omega \{
Image^{1,\beta}_{\phi_{u+1}}[L_{g}]\}+
 \Sum_{t\in T^{\sigma+u+2}} a_t
C^t_{g}(\Omega_1,\dots , \Omega_p,\phi_1,\dots
,\phi_{u+1},\omega)=0.
\end{split}
\end{equation}

\par We will keep the equation (\ref{baldouin'''}) in mind.
 We
now set out to study the sublinear combination
$Sub^{\sigma +u+1}_\omega \{
Image^{1,\beta}_{\phi_{u+1}} [L_{g}]\}$ in the above, which (we recall)
 consists of the complete contractions in
 $Sub_\omega \{ Image^{1,\beta}_{\phi_{u+1}}
[L_{g}]\}$ which have $\sigma+u+1$ factors.
\newline

\subsection{A study of the sublinear combination \\$Sub^{\sigma +u+1}_\omega \{
Image^{1,\beta}_{\phi_{u+1}} [L_{g}]\}$.}

\par In order to understand {\it how} the sublinear combination $Sub^{\sigma +u+1}_\omega \{
Image^{1,\beta}_{\phi_{u+1}} [L_{g}]\}$ arises from $L_g=0$,
 we must study certain {\it special patterns} of 
particular contractions among the different complete contractions in $L_g$:

 We think of $L_{g}$ as a linear combination of
complete contractions (i.e.~we momentarily forget its
structure-that it contains a linear combination of
$X$-divergences), and we also break the $S_{*}$-symmetrization in
the factors $S_{*}\nabla^{(\nu)}R_{ijkl}$--i.e. we treat those
terms as sums of factors $\nabla^{(\nu)}R_{ijkl}$. Thus we write out:

\begin{equation}
\label{panagiota} L_{g}(\Omega_1,\dots ,\Omega_p,\phi_1, \dots
,\phi_u)=\Sum_{v\in V} a_v C^v_{g}(\Omega_1,\dots
,\Omega_p,\phi_1, \dots ,\phi_u).
\end{equation}
 We then separately examine each complete contraction $C^v_{g}$ in
$L_{g}$ and we consider all the {\it  pairs of pairs} of indices,
 $[({}_a,{}_b), ({}_c,{}_d)]$ where ${}_a,{}_b$ are two
non-antisymmetric internal indices in a factor
$\nabla^{(m)}R_{ijkl}$ and the two indices ${}_c,{}_d$ belong to
some other factor and ${}_c$ is a derivative index, so that
${}_a$ is contracting against ${}_c$ and ${}_b$ against ${}_d$.

\par We denote the set of such pairs in each complete contraction
$C^v_{g}$ in $L_{g}$ by $INT^2_v$ (or just $INT^2$ for short).
Then, for each $[({}_a,{}_b), ({}_c,{}_d)]\in INT^2$, we denote by
$Rep^1_{[(a,b),(c,d)]}[C^v_{g}]$  the complete contraction that
 formally arises from $C^v_g$ by replacing the first factor $\nabla^{(m)}R_{ijkl}$ by an
expression $(+-)\nabla^{(m+2)}\phi_{u+1} g_{ab}$ and then
replacing the two indices $({}_c,{}_d)$ in the second factor 
(which now contract against each other, 
by virtue of the term $g_{ab}$ that we brought out) by an
 expression $(\nabla^d\omega,{}_d)$.\footnote{Recall that ${}_c$ is a derivative index, so
this formal operation is well-defined--in other words we may {\it formally} erase
the index ${}_c$ and bring in a new factor $\nabla\omega$ 
which will then contract against the index ${}_d$.} We also denote by
$Rep^2[C_{g}]$ the complete contraction
 that formally arises from $C^v_g$  by replacing the first factor
$\nabla^{(m)}R_{ijkl}$ by an expression $(+-)\nabla^{(m+2)}\omega
g_{ab}$,\footnote{See (\ref{curvtrans}).} and then replacing
 the two  indices $({}_c,{}_d)$ in the second
 factor by an expression $(\nabla^t\phi_{u+1},{}_t)$. We then define
(slightly abusing notation):

\begin{equation}
\label{xalia} Rep^1[C^v_{g}]=\Sum_{[(a,b), (c,d)]\in INT^2}
Rep^1_{[(a,b),(c,d)]}[C^v_{g}],
\end{equation}
and we similarly define $Rep^2[C^v_{g}]$.

It follows by definition that:

 \begin{equation}
 \label{tota}
Sub^{\sigma +u+1}_\omega \{ Image^{1,\beta}_{\phi_{u+1}}
[L_{g}]\}=\Sum_{v\in V} a_v (Rep^1+Rep^2)[C^v_{g}(\Omega_1,\dots
,\Omega_p,\phi_1, \dots ,\phi_u)].
 \end{equation}

\par Thus, we derive that:

\begin{equation}
\label{food}
\begin{split}
&Sub^{\sigma +u+1}_{\phi_{u+1}}\{
Image^{1,\beta}_{\phi_{u+1}}[L_{g}(\Omega_1,\dots ,\Omega_p,
\phi_1, \dots,\phi_u)]\}=
\\&\Sum_{l\in L} a_l
(Rep^1+Rep^2)[Xdiv_{i_1}\dots Xdiv_{i_a} C^{l,i_1\dots
i_a}_{g}(\Omega_1,\dots ,\Omega_p,\phi_1, \dots,\phi_u)]
\\&+\Sum_{j\in J} a_j (Rep^1+Rep^2)[ C^j_{g}(\Omega_1,\dots
,\Omega_p,\phi_1, \dots,\phi_u)].
\end{split}
\end{equation}

\par Now, we will be separately
studying each of the sublinear combinations
$$Sub^{\sigma +u+1}_\omega \{
Image^{1,\beta}_{\phi_{u+1}} [Xdiv_{i_1}\dots
Xdiv_{i_a}C^{l,i_1\dots i_a}_{g}]\},$$ $Sub^{\sigma +u+1}_\omega
\{ Image^{1,\beta}_{\phi_{u+1}} [C^j_{g}]\}$. By the definition 
of the formal operation $Rep$  
we derive:

\begin{equation}
\label{original} \Sum_{j\in J} a_j Sub^{\sigma +u+1}_\omega \{
Image^{1,\beta}_{\phi_{u+1}} [C^j_{g}]\}= \Sum_{j\in J^{\sigma
+u+1}} a_j C^j_{g}(\Omega_1,\dots , \Omega_p,\phi_1,
\dots,\phi_u,\phi_{u+1},\omega).
\end{equation}

\par On the other hand, in order to understand each sublinear
combination $Sub^{\sigma +u+1}_\omega \{
Image^{1,\beta}_{\phi_{u+1}} [Xdiv_{i_1}\dots
Xdiv_{i_a}C^{l,i_1\dots i_a}_{g}]\}$, we firstly seek to write out
that linear combination in a {\it normalized form}: We impose the
restriction that all the factors $\nabla\phi_1,\dots
,\nabla\phi_{u+1}$ or $\nabla\omega$ that are contracting
 against a factor $\nabla^{(A)}\Omega_f$ ($A\ge 2$)
 in the tensor fields appearing in the RHSs of all equations until 
(\ref{gounzales2}) must be contracting against
the leftmost indices. If this condition does not hold for
 some contraction that appears in $Sub^{\sigma +u+1}_{\phi_{u+1}}\{
Image^{1,\beta}_{\phi_{u+1}}[L_{g}]\}$, we apply the curvature identity
enough times to make it hold. (We will refer to this below as the
{\it shifting operation}).

We then distinguish three cases regarding the two internal
 indices $({}_a,{}_b)$ discussed in the  definition above:
Either in our tensor field $C^{l,i_1\dots i_a}_g$
neither of the indices ${}_a$, ${}_b$
 is a free index, or precisely one of them is a free index,
 or that both of them are free indices. We will accordingly denote
those sublinear combinations by

$$Sub^{\sigma +u+1,I}_\omega \{
Image^{1,\beta}_{\phi_{u+1}} [Xdiv_{i_1}\dots
Xdiv_{i_a}C^{l,i_1\dots i_a}_{g}]\},$$

$$Sub^{\sigma +u+1,II}_\omega \{
Image^{1,\beta}_{\phi_{u+1}} [Xdiv_{i_1}\dots
Xdiv_{i_a}C^{l,i_1\dots i_a}_{g}]\},$$

$$Sub^{\sigma +u+1,III}_\omega \{
Image^{1,\beta}_{\phi_{u+1}} [Xdiv_{i_1}\dots
Xdiv_{i_a}C^{l,i_1\dots i_a}_{g}]\}.$$

\par Then, for each $l\in L$, by the above
 discussion, we calculate:

\begin{equation}
\label{paulo} \begin{split} & Sub^{\sigma +u+1,I}_\omega \{
Image^{1,\beta}_{\phi_{u+1}} [Xdiv_{i_1}\dots
Xdiv_{i_a}C^{l,i_1\dots i_a}_{g}]\}=
\\& \Sum_{u\in U_1\bigcup U^\sharp_1} a_u Xdiv_{i_1}\dots Xdiv_{i_a}
C^{u,i_1\dots ,i_a}_{g}(\Omega_1,\dots ,\Omega_p,\phi_1, \dots
,\phi_u,\phi_{u+1},\omega)+
\\&\Sum_{u\in U_2\bigcup U^\sharp_2} a_u Xdiv_{i_1}\dots Xdiv_{i_a}
C^{u,i_1\dots ,i_a}_{g}(\Omega_1,\dots ,\Omega_p,\phi_1, \dots
,\phi_u,\phi_{u+1},\omega)+ \\& \Sum_{m\in M\bigcup M^\sharp} a_m
C^{m,i_1\dots i_{a+2}}_{g} (\Omega_1,\dots ,\Omega_p,\phi_1, \dots
,\phi_u)\nabla_{i_{a+1}}\phi_{u+1}\nabla_{i_{a+2}} \omega+
\\&\Sum_{j\in J^{\sigma+u+1}} a_j C^j_{g}(\Omega_1,\dots ,
\Omega_p,\phi_1,\dots ,\phi_{u+1},\omega).
\end{split}
\end{equation}

\par Now, we seek to study each $Sub^{\sigma +u+1,II}_\omega \{
Image^{1,\beta}_{\phi_{u+1}} [Xdiv_{i_1}\dots
Xdiv_{i_a}C^{l,i_1\dots i_a}_{g}]\}$. We will draw a different
conclusion, in the cases where $l\in L_\mu$ and where $l\in
L\setminus L_\mu$. In general, for any $l\in L$, we make special
note of the
 factors $S_{*}\nabla^{(\nu)} R_{ijkl}$ in $C^{l,i_1\dots i_a}_{g}$
that contain free indices.
 By the hypothesis of Lemma \ref{pskovb} $l\in L_\mu$, those free indices will
 be of the form ${}_{r_1},\dots ,{}_{r_\nu},{}_j$.

\begin{definition}
\label{eurovision}
 We denote the set of free indices that belong to  factors
 $S_{*}\nabla^{(\nu)}R_{ijkl}$ with $\nu>0$ by
 $I^\sharp$. For future reference we denote by $I^\sharp_{*}$
 the set of free indices that belong to a
 factor $S_{*}\nabla^{(\nu)} R_{ijkl}$ with $\nu=0$.
\end{definition}
 Recall that for tensor
 fields $C^{l,i_1\dots i_\mu}_{g}, l\in L_\mu$ we will have
$I^\sharp_{*}=\emptyset$, by virtue of
 the assumptions in the beginning of this paper.

\par Now, another definition:

\begin{definition}
\label{eurovision2} For each ${}_{i_h}\in I^\sharp$, we denote by
$C^{l,i_1\dots i_a|f(i_h)}_{g}$ the tensor field that arises from
$C^{l,i_1\dots i_a}_{g}$ by replacing the factor
$T({}_{i_h})=S_{*}\nabla^{(\nu)}_{i_hr_2\dots r_\nu}R_{ijkl}$ to which
${}_{i_h}$ belongs (assume with no loss of generality that
${}_{i_h}={}_{r_1}$) by a factor
$\frac{1}{\nu+1}\nabla^{(\nu)}_{r_2\dots r_{\nu -1}j}R_{ii_hkl}$.

Then, for each ${}_{i_h}\in I^\sharp\bigcup I^\sharp_{*}$, we denote by
$$Rep^{i_h,1,\phi_{u+1},\omega}[C^{l,i_1\dots i_a|f(i_h)}_{g}
(\Omega_1,\dots ,\Omega_p,\phi_1,\dots ,\phi_u)]$$ the
$(a-1)$-tensor field  that arises from $C^{l,i_1\dots
i_a|f(i_h)}_{g}$
 by replacing the factor $\frac{1}{\nu+1}\nabla^{(\nu)}_{r_2\dots r_{\nu -1}j} R_{ii_hkl}$ by a factor
$\frac{1}{\nu+1}\nabla^{(\nu+2)}_{r_2\dots r_{\nu -1}jil}\omega$
and then making the index ${}^k$ contract against a factor
${\nabla}_k\phi_{u+1}$. We also denote by
$$Rep^{i_h,1,\omega,\phi_{u+1}}[C^{l,i_1\dots i_a|f(i_h)}_{g}
(\Omega_1,\dots ,\Omega_p,\phi_1,\dots ,\phi_u)]$$
the $(a-1)$-tensor field that arises from
$$Rep^{i_h,1,\phi_{u+1},\omega}[C^{l,i_1\dots i_a|
f(i_h)}_{g}(\Omega_1,\dots ,\Omega_p,\phi_1,\dots , \phi_u)]$$ by
replacing $\phi_{u+1}$ by $\omega$ and $\omega$ by $\phi_{u+1}$.

\par Analogously, we denote by
$Rep^{i_h,2,\phi_{u+1},\omega}[C^{l,i_1\dots i_a|f(i_h)}_{g}
(\Omega_1,\dots ,\Omega_p,\phi_1,\dots ,\phi_u)]$ the
$(a-1)$-tensor field  that arises from $C^{l,i_1\dots
i_a|f(i_h)}_{g}$ by replacing the factor
$\frac{1}{\nu+1}\nabla^{(\nu)}_{r_2\dots r_{\nu -1}j}R_{ii_hkl}$ by a
factor $-\frac{1}{\nu+1}\nabla^{(\nu+2)}_{r_2\dots r_{\nu -1}jik}
\omega$ and then making the index ${}^l$ contract against a factor
${\nabla}_l\phi_{u+1}$. We also denote by
\\ $Rep^{i_h,2,\omega,\phi_{u+1}}
[C^{l,i_1\dots i_a|f(i_h)}_{g} (\Omega_1,\dots
,\Omega_p,\phi_1,\dots ,\phi_u)]$ the $(a-1)$-tensor field that
arises from
$$Rep^{i_h,2,\phi_{u+1},\omega}[C^{l,i_1\dots i_a|
f(i_h)}_{g}(\Omega_1,\dots ,\Omega_p,\phi_1,\dots , \phi_u)]$$ by
replacing $\phi_{u+1}$ by $\omega$ and $\omega$ by $\phi_{u+1}$.

\par Finally, for each ${}_{i_h}\in I^\sharp$ we denote by
$$FRep^{i_h,1,\phi_{u+1},\omega} [C^{l,i_1\dots i_a|f(i_h)}_{g}
(\Omega_1,\dots ,\Omega_p,\phi_1,\dots ,\phi_u)]$$ the $(a-1)$
tensor fields that
 arises from $$Rep^{i_h,1,\phi_{u+1},\omega}
[C^{l,i_1\dots i_a|f(i_h)}_{g} (\Omega_1,\dots
,\Omega_p,\phi_1,\dots ,\phi_u)]$$ by replacing the expression
$\frac{1}{\nu+1}\nabla^{(\nu+2)}_{r_2\dots r_{\nu -1}jik}
\omega\nabla^i\phi_f,\dots $ by an expression \\$\frac{1}{\nu+1}
S_{*}\nabla^{(\nu-1)}_{r_2\dots r_{\nu -1}}R_{ijlk}
\nabla^l\omega\nabla^i\phi_f$. (This is well-defined since
${}_{i_h}\in I^\sharp$, therefore $\nu>0$). Accordingly we define
 $FRep^{i_h,2,\phi_{u+1},\omega}$, $FRep^{i_h,1,\omega,
\phi_{u+1}}$, $FRep^{i_h,2,\omega,\phi_{u+1}}$.
\end{definition}

\par Then, in view of the hypotheses of Lemma \ref{pskovb},
and by applying the shifting operation to
$\nabla^{(\nu+2)}\omega$,
we derive that for each $l\in L_\mu$:

\begin{equation}
\label{mumu}
\begin{split}
&Sub^{\sigma +u+1,II}_\omega \{ Image^{1,\beta}_{\phi_{u+1}}
[Xdiv_{i_1}\dots Xdiv_{i_\mu}C^{l,i_1\dots
i_\mu}_{g}(\Omega_1,\dots , \Omega_p,\phi_1,\dots ,\phi_u)]\}
\\&=\Sum_{i_h\in I^\sharp}Xdiv_{i_1}\dots \hat{Xdiv}_{i_h}
\dots Xdiv_{i_\mu}\\&\{ FRep^{i_h,1,\omega,\phi_{u+1}} [C^{l,i_1\dots
i_\mu|f(i_h)}_{g} (\Omega_1,\dots ,\Omega_p,\phi_1,\dots ,\phi_u)]
\\& +FRep^{i_h,1,\phi_{u+1},\omega}[C^{l,i_1\dots i_\mu|
f(i_h)}_{g}(\Omega_1,\dots ,\Omega_p,\phi_1,\dots ,\phi_u)]
\\& +FRep^{i_h,2,\omega,\phi_{u+1}}[C^{l,i_1\dots i_\mu
|f(i_h)}_{g}(\Omega_1,\dots ,\Omega_p,\phi_1,\dots ,\phi_u)] +\\&
FRep^{i_h,2,\phi_{u+1},\omega}[C^{l,i_1\dots i_\mu
|f(i_h)}_{g}(\Omega_1,\dots ,\Omega_p,\phi_1,\dots , \phi_u)]\}+
\\&\Sum_{j\in J^{\sigma+u+1}} a_j
C^j_g(\Omega_1,\dots ,\Omega_p,\phi_1,\dots ,\phi_u,\phi_{u+1},
\omega).
\end{split}
\end{equation}

\par On the other hand, by virtue of our  assumptions
on the tensor fields in (\ref{hypothese2})
and by the analysis above, we easily see that for each
$l\in L\setminus L_\mu$:

\begin{equation}
\label{gounzales}
\begin{split}
&Sub^{\sigma +u+1,II}_\omega \{ Image^{1,\beta}_{\phi_{u+1}}
[Xdiv_{i_1}\dots Xdiv_{i_a}C^{l,i_1\dots i_a}_{g}(\Omega_1,\dots
,\Omega_p, \phi_1, \dots,\phi_u)]\}=
\\& \Sum_{u\in U_1} a_u Xdiv_{i_1}\dots Xdiv_{i_a}
C^{u,i_1\dots ,i_a}_{g}(\Omega_1,\dots ,\Omega_p,\phi_1, \dots
,\phi_u,\phi_{u+1},\omega)+
\\&\Sum_{u\in U_2} a_u Xdiv_{i_1}\dots Xdiv_{i_a}
C^{u,i_1\dots ,i_a}_{g}(\Omega_1,\dots ,\Omega_p, \phi_1,
\dots,\phi_u,\phi_{u+1},\omega)+
\\& \Sum_{m\in M} a_m C^{m,i_1\dots i_{a+2}}_{g}
(\Omega_1,\dots ,\Omega_p,\phi_1,\dots ,\phi_u)
\nabla_{i_{a+1}}\phi_{u+1}\nabla_{i_{a+2}}\omega+
\\&\Sum_{j\in J^{\sigma+u+1}} a_j
C^j_{g}(\Omega_1,\dots ,\Omega_p,\phi_1, \dots
,\phi_u,\phi_{u+1},\omega).
\end{split}
\end{equation}

\par Finally, we study the linear combinations
$$Sub^{\sigma +u+1,III}_\omega \{
Image^{1,\beta}_{\phi_{u+1}} [Xdiv_{i_1}\dots
Xdiv_{i_a}C^{l,i_1\dots i_a}_{g}]\}.$$ Clearly, by the hypothesis
of Lemma \ref{pskovb} that no $\mu$-tensor fields have special
free indices, we derive that for each $l\in L_\mu$:

$$Sub^{\sigma +u+1,III}_\omega \{
Image^{1,\beta}_{\phi_{u+1}} [Xdiv_{i_1}\dots
Xdiv_{i_a}C^{l,i_1\dots i_a}_{g}]\}=0.$$

\par In addition, for each $l\in L_k, k\ge \mu+2$ we
straightforwardly  obatin:

\begin{equation}
\label{gounzales1}
\begin{split}
&Sub^{\sigma +u+1,III}_\omega \{ Image^{1,\beta}_{\phi_{u+1}}
[Xdiv_{i_1}\dots Xdiv_{i_a}C^{l,i_1\dots i_a}_{g}(\Omega_1,\dots
,\Omega_p, \phi_1, \dots,\phi_u)]\}
\\&= \Sum_{u\in U_1} a_u Xdiv_{i_1}\dots Xdiv_{i_a}
C^{u,i_1\dots ,i_a}_{g}(\Omega_1,\dots ,\Omega_p,\phi_1, \dots
,\phi_u,\phi_{u+1},\omega)+
\\&\Sum_{u\in U_2} a_u Xdiv_{i_1}\dots Xdiv_{i_a}
C^{u,i_1\dots ,i_a}_{g}(\Omega_1,\dots ,\Omega_p,\phi_1,
 \dots,\phi_u,\phi_{u+1},\omega)+
\\& \Sum_{m\in M} a_m C^{m,i_1\dots i_{a+2}}_{g}
 (\Omega_1,\dots ,\Omega_p,\phi_1,\dots ,\phi_u)
\nabla_{i_{a+1}}\phi_{u+1}\nabla_{i_{a+2}}\omega+
\\&\Sum_{j\in J^{\sigma+u+1}} a_j
C^j_{g}(\Omega_1,\dots ,\Omega_p,\phi_1, \dots
,\phi_u,\phi_{u+1},\omega).
\end{split}
\end{equation}

\par So, we next set out to study the linear combination
$$Sub^{\sigma +u+1,III}_\omega \{
Image^{1,\beta}_{\phi_{u+1}} [Xdiv_{i_1}\dots
Xdiv_{i_a}C^{l,i_1\dots i_a}_{g}(\Omega_1,\dots ,\Omega_p, \phi_1,
\dots,\phi_u)]\}$$ for each $l\in L_{\mu+1}$. 
We just introduce one more piece of notation:

\begin{definition}
\label{matia}
We denote by $\sum_{b\in B} a_b C^{b,i_1\dots i_{\mu-1}}_g
(\Omega_1,\dots, \Omega_p,\phi_{u+1},\phi_1,\dots,\phi_u)$ a generic linear combination 
of acceptable $(\mu-1)$-tensor fields in the form (\ref{form2}), with length 
$\sigma+u$, weight $-n+\mu+3$, and a $u$-simple
 character that arises from $\vec{\kappa}_{simp}$ by
  formally replacing a factor $\nabla^{(m)}R_{ijkl}$ by $\nabla^{(m+2)}Y$. 
We denote the other
factors other than $\nabla\phi$'s in $C^{z,i_1\dots i_{\mu-1}}_{g}$
 by $F_1,\dots ,F_{\sigma-1}$. 

\par We then define
$Hit^K_\omega[C^{b,i_1\dots i_{\mu -1}}_{g}(\Omega_1,\dots
,\Omega_p, \phi_{u+1},\phi_1, \dots,\phi_u)]$ to stand for the
$(\mu-1)$-tensor field that arises from $C^{b,i_1\dots i_{\mu
-1}}_{g}(\Omega_1,\dots ,\Omega_p, \phi_{u+1},\phi_1,
\dots,\phi_u)$ by hitting the factor $F^K$  by a
derivative $\nabla_{i_{*}}$ and then contracting ${}_{i_{*}}$
against a factor $\nabla\omega$.
\newline

\par We denote by
$$\Sum_{\zeta\in Z} a_\zeta C^{\zeta,i_1\dots i_{\mu -1}}_{g}(\Omega_1,\dots ,\Omega_p,
\phi_1, \dots,\phi_u,\phi_{u+1})$$
 a generic linear combination of
$(\mu-1)$-acceptable tensor fields with weight
 $-n+\mu+3$, length $\sigma+u+1$ and 
a $u$-simple character $\vec{\kappa}_{simp}$
 and with one factor $\nabla\phi_{u+1}$
contracting against a factor
$\nabla^{(\nu)} R_{ijkl}$. We denote this factor
$\nabla^{(\nu)} R_{ijkl}$ by $F$ and we denote the other
factors other than $\nabla\phi$'s in $C^{\zeta,i_1\dots i_{\mu
-1}}_{g}$ by $F_1,\dots ,F_{\sigma-1}$. 

\par We then define
$Hit^K_\omega[C^{\zeta,i_1\dots i_{\mu -1}}_{g}(\Omega_1,\dots
,\Omega_p, \phi_1, \dots,\phi_u,\phi_{u+1})]$ to stand for the
$(\mu-1)$-tensor field that arises from $C^{\zeta,i_1\dots i_{\mu
-1}}_{g}(\Omega_1,\dots ,\Omega_p, \phi_1,
\dots,\phi_u,\phi_{u+1})$ by hitting the factor $F^K$  by a
derivative $\nabla_{i_{*}}$ and then contracting ${}_{i_{*}}$
against a factor $\nabla\omega$ (or $\nabla\phi_{u+1}$).
\end{definition}

 Two types of linear combinations that we will be
encountering are linear combinations in the forms:

$$\Sum_{b\in B} a_b \Sum_{K=1}^{\sigma-1}Hit^K_\omega
[C^{b,i_1\dots i_{\mu -1}}_{g}(\Omega_1,\dots ,
\Omega_p,\phi_{u+1}, \phi_1,
\dots,\phi_u)],$$
$$\Sum_{\zeta\in Z} a_\zeta \Sum_{K=1}^{\sigma-1}Hit^K_\omega
[C^{\zeta,i_1\dots i_{\mu -1}}_{g}(\Omega_1,\dots ,\Omega_p, \phi_1,
\dots,\phi_u,\phi_{u+1})].$$ 
For complete contractions as above, we denote by 
$$Switch\{ Hit^K_\omega
[C^{b,i_1\dots i_{\mu -1}}_{g}(\Omega_1,\dots ,
\Omega_p,\phi_{u+1}, \phi_1,
\dots,\phi_u)]\}$$ the tensor fields that arise by
 interchanging the functions $\phi_{u+1}, \omega$.

\par Now, by the same argument as for equation (\ref{mumu}),
we derive that for each $l\in L_{\mu +1}$:

\begin{equation}
\label{gounzales2}
\begin{split}
&Sub^{\sigma +u+1,III}_\omega \{ Image^{1,\beta}_{\phi_{u+1}}
[Xdiv_{i_1}\dots Xdiv_{i_a}C^{l,i_1\dots i_a}_{g}(\Omega_1,\dots
,\Omega_p, \phi_1, \dots,\phi_u)]\}=
\\&\Sum_{b\in B} a_b Xdiv_{i_1}\dots Xdiv_{i_{\mu-1}}\Sum_{K=1}^{\sigma-1}\{Hit^K_\omega
[C^{b,i_1\dots i_{\mu -1}}_{g}(\Omega_1,\dots ,
\Omega_p,\phi_{u+1}, \phi_1,
\dots,\phi_u)]
\\&Switch\{Hit^K_\omega
[C^{b,i_1\dots i_{\mu -1}}_{g}(\Omega_1,\dots ,
\Omega_p,\phi_{u+1}, \phi_1,
\dots,\phi_u)]\}\}+
\\& \Sum_{m\in M} a_m Xdiv_{i_1}\dots Xdiv_{i_a}C^{m,i_1\dots i_{a+2}}_{g}
(\Omega_1,\dots ,\Omega_p,\phi_1,\dots ,\phi_u)
\nabla_{i_{a+1}}\phi_{u+1}\nabla_{i_{a+2}}\omega
\\&+\Sum_{j\in J^{\sigma+u+1}} a_j
C^j_{g}(\Omega_1,\dots ,\Omega_p,\phi_1, \dots
,\phi_u,\phi_{u+1},\omega).
\end{split}
\end{equation}

\par In conclusion, combining equations (\ref{reno3}),
(\ref{original}), (\ref{paulo}),
(\ref{mumu}), (\ref{gounzales}), (\ref{gounzales1}),
 (\ref{gounzales2}), and replacing them into (\ref{baldouin'''}),
we have shown that:

\begin{equation}
\label{uptohere}
\begin{split}
&\Sum_{l\in L} a_l Sub^{\sigma +u+2}_\omega
\{Image^{1,\beta}_{\phi_{u+1}} [Xdiv_{i_1}\dots Xdiv_{i_a}
C^{l,i_1\dots ,i_a}_{g}(\Omega_1,\dots ,\Omega_p, \phi_1,
\dots,\phi_u)]\}
\\&+\Sum_{j\in J} a_j Sub^{\sigma +u+2}_\omega
\{Image^{1,\beta}_{\phi_{u+1}} [
C^j_{g}(\Omega_1,\dots ,\Omega_p, \phi_1,
\dots,\phi_u)]\}+
\\&\Sum_{l\in L_\mu} a_l \Sum_{i_h\in I^\sharp}Xdiv_{i_1}\dots \hat{Xdiv}_{i_h}
\dots Xdiv_{i_\mu}
\\&\{FRep^{i_h,1,\omega,\phi_{u+1}} [C^{l,i_1\dots i_\mu|f(i_h)}_{g}
(\Omega_1,\dots ,\Omega_p,\phi_1,\dots ,\phi_u)]
\\& +\Sum_{l\in L_\mu} a_l FRep^{i_h,1,\phi_{u+1},\omega}
[C^{l,i_1\dots i_\mu| f(i_h)}_{g}(\Omega_1,\dots
,\Omega_p,\phi_1,\dots ,\phi_u)]
\\& +\Sum_{l\in L_\mu} a_l FRep^{i_h,2,\omega,\phi_{u+1}}[C^{l,i_1\dots i_\mu
|f(i_h)}_{g}(\Omega_1,\dots ,\Omega_p,\phi_1,\dots ,\phi_u)]
+\\&\Sum_{l\in L_\mu} a_l
FRep^{i_h,2,\phi_{u+1},\omega}[C^{l,i_1\dots i_\mu
|f(i_h)}_{g}(\Omega_1,\dots ,\Omega_p,\phi_1,\dots , \phi_u)]\}+
\\&\Sum_{b\in B} a_b Xdiv_{i_1}\dots Xdiv_{i_{\mu-1}}\Sum_{K=1}^{\sigma-1}\{Hit^K_\omega
[C^{b,i_1\dots i_{\mu -1}}_{g}(\Omega_1,\dots ,
\Omega_p,\phi_{u+1}, \phi_1,
\dots,\phi_u)]+
\\&Switch\{Hit^K_\omega
[C^{b,i_1\dots i_{\mu -1}}_{g}(\Omega_1,\dots ,
\Omega_p,\phi_{u+1}, \phi_1,
\dots,\phi_u)]\}\}+
\\& \Sum_{u\in U_1\bigcup U^\sharp_1} a_u Xdiv_{i_1}\dots Xdiv_{i_a}
C^{u,i_1\dots ,i_a}_{g}(\Omega_1,\dots ,\Omega_p,\phi_1, \dots
,\phi_u,\phi_{u+1},\omega)+
\\&\Sum_{u\in U_2\bigcup U^\sharp_2} a_u Xdiv_{i_1}\dots Xdiv_{i_a}
C^{u,i_1\dots ,i_a}_{g}(\Omega_1,\dots ,\Omega_p, \phi_1,
\dots,\phi_u,\phi_{u+1},\omega)+
\\& \Sum_{m\in M} a_m Xdiv_{i_1}\dots Xdiv_{i_a} C^{m,i_1\dots i_{a+2}}_{g}
(\Omega_1,\dots ,\Omega_p,\phi_1,\dots ,\phi_u)
\nabla_{i_{a+1}}\phi_{u+1}\nabla_{i_{a+2}}\omega
\\&+\Sum_{j\in J^{\sigma+u+1}} a_j
C^j_{g}(\Omega_1,\dots ,\Omega_p,\phi_1, \dots
,\phi_u,\phi_{u+1},\omega) +
\\&\Sum_{\zeta\in Z} a_\zeta Xdiv_{i_1}\dots Xdiv_{\mu-1}
Hit^K_\omega[C^{\zeta,i_1\dots i_{\mu -1}}_{g}(\Omega_1,\dots
,\Omega_p, \phi_1, \dots,\phi_u,\phi_{u+1})]
\\&= \Sum_{t\in T^{\sigma+u+2}} a_t
C^t_{g}(\Omega_1,\dots ,\Omega_p, \phi_1,
\dots,\phi_u,\phi_{u+1},\omega),
\end{split}
 \end{equation}
modulo complete contractions of length $\ge\sigma +u+3$. The
sublinear combination in the RHS is a generic 
sublinear combination as defined below (\ref{baldouin}).
 Notice that the
 minimum length among the complete contractions
 above is $\sigma+u+1$. The complete contractions (and tensor fields)
with $\sigma+u+1$ factors are indexed in 
$U_1,U^\sharp_1,U_2,U^\sharp_2, B, J^{\sigma+u+1}$.

\par Therefore the above equation implies:

\begin{equation}
\label{toki}
\begin{split}
&\Sum_{b\in B} a_b Xdiv_{i_1}\dots 
Xdiv_{i_{\mu-1}}\Sum_{K=1}^{\sigma-1}\{Hit^K_\omega
[C^{b,i_1\dots i_{\mu -1}}_{g}(\Omega_1,\dots ,
\Omega_p,\phi_{u+1}, \phi_1,
\dots,\phi_u)]
\\&+Switch\{Hit^K_\omega
[C^{b,i_1\dots i_{\mu -1}}_{g}(\Omega_1,\dots ,
\Omega_p,\phi_{u+1}, \phi_1,
\dots,\phi_u)]\}\}+
\\& \Sum_{u\in U_1\bigcup U^\sharp_1} a_u Xdiv_{i_1}\dots Xdiv_{i_a}
C^{u,i_1\dots ,i_a}_{g}(\Omega_1,\dots ,\Omega_p,\phi_1, \dots
,\phi_u,\phi_{u+1},\omega)+
\\&\Sum_{u\in U_2\bigcup U^\sharp_2} a_u Xdiv_{i_1}\dots Xdiv_{i_a}
C^{u,i_1\dots ,i_a}_{g}(\Omega_1,\dots ,\Omega_p, \phi_1,
\dots,\phi_u,\phi_{u+1},\omega)+
\\&\Sum_{j\in J^{\sigma+u+1}} a_j
C^j_{g}(\Omega_1,\dots ,\Omega_p,\phi_1, \dots
,\phi_u,\phi_{u+1},\omega)=0,
\end{split}
\end{equation}
modulo complete contractions of length $\ge\sigma+u+2$.
We claim that in the {\it non-special cases}:

\begin{equation}
\label{tokicor}
\begin{split}
&\Sum_{b\in B} a_b Xdiv_{i_1}\dots Xdiv_{i_{\mu-1}}\Sum_{K=1}^{\sigma-1}\{Hit^K_\omega
[C^{b,i_1\dots i_{\mu -1}}_{g}(\Omega_1,\dots ,
\Omega_p,\phi_{u+1}, \phi_1,
\dots,\phi_u)]
\\&+Switch\{Hit^K_\omega
[C^{b,i_1\dots i_{\mu -1}}_{g}(\Omega_1,\dots ,
\Omega_p,\phi_{u+1}, \phi_1,
\dots,\phi_u)]\}\}=
\\& \Sum_{u\in U_1\bigcup U^\sharp_1} a_u Xdiv_{i_1}\dots Xdiv_{i_a}
C^{u,i_1\dots ,i_a}_{g}(\Omega_1,\dots ,\Omega_p,\phi_1, \dots
,\phi_u,\phi_{u+1},\omega)+
\\&\Sum_{u\in U_2\bigcup U^\sharp_2} a_u Xdiv_{i_1}\dots Xdiv_{i_a}
C^{u,i_1\dots ,i_a}_{g}(\Omega_1,\dots ,\Omega_p, \phi_1,
\dots,\phi_u,\phi_{u+1},\omega)+
\\&\Sum_{j\in J^{\sigma+u+1}} a_j
C^j_{g}(\Omega_1,\dots ,\Omega_p,\phi_1, \dots
,\phi_u,\phi_{u+1},\omega)+
\\&\Sum_{\zeta\in Z} a_\zeta  Xdiv_{i_1}\dots Xdiv_{i_{\mu-1}}\{\Sum_{K=1}^{\sigma-1}Hit^K_\omega
[C^{\zeta,i_1\dots i_{\mu -1}}_{g}(\Omega_1,\dots ,\Omega_p, \phi_1,
\dots,\phi_u,\phi_{u+1})]
\\&+Switch[\Sum_{K=1}^{\sigma-1}Hit^K_\omega
[C^{\zeta,i_1\dots i_{\mu -1}}_{g}(\Omega_1,\dots ,\Omega_p, \phi_1,
\dots,\phi_u,\phi_{u+1})]]\}+
\\& \Sum_{m\in M} a_m C^{m,i_1\dots i_{a+2}}_{g}
(\Omega_1,\dots ,\Omega_p,\phi_1,\dots ,\phi_u)
\nabla_{i_{a+1}}\phi_{u+1}\nabla_{i_{a+2}}\omega+
\\&\Sum_{j\in
J^{\sigma+u+2}} a_j C^j_{g}(\Omega_1,\dots ,\Omega_p,\phi_1, \dots
,\phi_u,\phi_{u+1},\omega) \\&= \Sum_{t\in T} a_t
C^t_{g}(\Omega_1,\dots ,\Omega_p, \phi_1,
\dots,\phi_u,\phi_{u+1},\omega),
\end{split}
\end{equation}
{\it modulo complete contractions of length $\ge\sigma+u+3$}. 
Here the terms indexed in $U_1, U_2$ in the RHS are
 {\it generic} linear combination in the forms described in Definition \ref{leaves}.
 In the special cases, we claim: 
 
\begin{equation}
\label{tokicor'}
\begin{split}
&\Sum_{b\in B} a_b Xdiv_{i_1}\dots Xdiv_{i_{\mu-1}}
\Sum_{K=1}^{\sigma-1}\{Hit^K_\omega
[C^{b,i_1\dots i_{\mu -1}}_{g}(\Omega_1,\dots ,
\Omega_p,\phi_{u+1}, \phi_1,
\dots,\phi_u)]
\\&+Switch\{Hit^K_\omega
[C^{b,i_1\dots i_{\mu -1}}_{g}(\Omega_1,\dots ,
\Omega_p,\phi_{u+1}, \phi_1,
\dots,\phi_u)]\}\}+
\\& \Sum_{u\in U_1\bigcup U^\sharp_1} a_u Xdiv_{i_1}\dots Xdiv_{i_a}
C^{u,i_1\dots ,i_a}_{g}(\Omega_1,\dots ,\Omega_p,\phi_1, \dots
,\phi_u,\phi_{u+1},\omega)+
\\&\Sum_{u\in U_2\bigcup U^\sharp_2} a_u Xdiv_{i_1}\dots Xdiv_{i_a}
C^{u,i_1\dots ,i_a}_{g}(\Omega_1,\dots ,\Omega_p, \phi_1,
\dots,\phi_u,\phi_{u+1},\omega)+
\\&\Sum_{j\in J^{\sigma+u+1}\bigcup J^{\sigma+u+2}} a_j
C^j_{g}(\Omega_1,\dots ,\Omega_p,\phi_1, \dots
,\phi_u,\phi_{u+1},\omega)
\\&= \Sum_{m\in M_{\mu-1}} a_m
C^{m,i_1\dots i_{mu-1}}_{g}(\Omega_1,\dots ,\Omega_p, \phi_1,
\dots,\phi_u,\phi_{u+1},\omega).
\end{split}
\end{equation}
Here the tensor fields indexed in $M_{\mu-1}$ are acceptable,
 have length $\sigma+u+2$, $u$-simple
  character $\vec{\kappa}_{simp}$ and moreover 
each of the $\mu-1$ free indices belongs to a different factor.

The harder challenge is to prove (\ref{tokicor}), so we start with that equation.
\newline

{\it Proof of (\ref{tokicor}):} Let us pick out the sublinear combination in
(\ref{toki}) with a factor
 $\nabla\omega$ contracting against a given factor
  $F_1$.\footnote{We assume for convenience that $F_1$ is a 
  well-defined factor in $\vec{\kappa}_{simp}$. If it were not, we
 just pick out the sublinear combination where
$\nabla\phi_{u+1}$ contracts against any generic factor
 $\nabla^{(m)}R_{ijkl}$ and the same argument applies.} 
 Since this sublinear combination 
 must vanish separately, we derive that:

\begin{equation}
\label{toki'}
\begin{split}
&\Sum_{b\in B} a_b Xdiv_{i_1}\dots Xdiv_{i_{\mu-1}}Hit^1_\omega
[C^{b,i_1\dots i_{\mu -1}}_{g}(\Omega_1,\dots ,
\Omega_p,\phi_{u+1}, \phi_1,
\dots,\phi_u)]+
\\& \Sum_{u\in \overline{U}_1\bigcup \overline{U}^\sharp_1} a_u Xdiv_{i_1}\dots Xdiv_{i_a}
C^{u,i_1\dots ,i_a}_{g}(\Omega_1,\dots ,\Omega_p,\phi_1, \dots
,\phi_u,\phi_{u+1},\omega)+
\\&\Sum_{j\in J^{\sigma+u+1}} a_j
C^j_{g}(\Omega_1,\dots ,\Omega_p,\phi_1, \dots
,\phi_u,\phi_{u+1},\omega)=0,
\end{split}
\end{equation}
modulo complete contractions of length $\ge\sigma+u+2$. Here the index sets 
$\overline{U}_1\subset U_1, \overline{U}^\sharp_1\subset U^\sharp_1$,
 are the index sets of terms with
 a factor $\nabla\omega$ contracting against the factor $F_1$. 
 
 Our aim is to 
 derive an equation like the above, only with 
 the factor $\nabla\omega$ contracting against 
a derivative index in the factor $F_1$, and moreover, if $F_1$
 is of the form $\nabla^{(B)}\Omega_h$, 
then we additionaly require that $B\ge 3$.  Call this the $*$-property.
  Now, if $F_1$ is a curvature factor, we apply the inductive assumption of Lemma 
  4.10 in \cite{alexakis4} to ensure  that 
  in all terms in $\overline{U}_1$ the factor $\nabla\phi_{u+1}$ 
  is not contracting against a 
  special index. Now, if $F_1$ is a factor $\nabla^{(B)}\Omega_x$, 
  we apply the inductive assumption of 
  Lemma 4.1 in \cite{alexakis4} if necessary 
  to assume wlog that $\overline{U}_1^\sharp=\emptyset$. 
  Finally, if needed, 
 we apply the inductive assumptions onf Corollaries 2 or 3 in
\cite{alexakis4} (if $F_1$ is a simple factor in the
 form $S_{*}\nabla^{(\nu)}R_{ijkl}$, or a
simple factor in the form $\nabla^{(A)}\Omega_h$, respectively) to ensure that
for each $u\in \overline{U}_1$  $\nabla\omega$ is 
not contracting against a factor $S_{*}R_{ijkl}$
 or $\nabla^{(2)}\Omega_h$.\footnote{In all the above 
 applications of Lemmas and Corollaries from \cite{alexakis4}, 
 we observe that by weight considerations, the fact that (\ref{hypothese2})
 does not fall under the special cases ensures that there is no danger of 
 falling under a forbidden case of those Lemmas/Corollaries.} Therefore, we may assume
  wlog that the $*$-property holds in (\ref{toki'}).
 
  We then apply the 
 Eraser to the factor $\nabla\phi_{u+1}$ (see the Appendix in \cite{alexakis1})
  and derive a new equation:
 
\begin{equation}
\label{toki''}
\begin{split}
&\Sum_{b\in B} a_b Xdiv_{i_1}\dots Xdiv_{i_{\mu-1}}
 C^{b,i_1\dots i_{\mu -1}}_{g}(\Omega_1,\dots ,
\Omega_p,\phi_{u+1}, \phi_1,
\dots,\phi_u)+
\\& \Sum_{u\in \overline{U}_1\bigcup \overline{U}^\sharp_1} a_u Xdiv_{i_1}\dots Xdiv_{i_a}
Erase_{\omega}[C^{u,i_1\dots ,i_a}_{g}(\Omega_1,\dots ,\Omega_p,\phi_1, \dots
,\phi_u,\phi_{u+1},\omega)]+
\\&\Sum_{j\in J^{\sigma+u+1}} a_j
C^j_{g}(\Omega_1,\dots ,\Omega_p,\phi_1, \dots
,\phi_u,\phi_{u+1})=0,
\end{split}
\end{equation}
modulo complete contractions of length $\ge\sigma+u+1$. 
 
 Now, apply the inductive assumption of Corollary 1 in \cite{alexakis4} to the above. 
 We derive that there exists a linear combination of acceptable 
 $\mu$-tensor fields with a simple character 
 $Pre(\vec{\kappa}_{simp})$ (indexed in $P$ below)
  such that:

\begin{equation}
\label{toki''}
\begin{split}
&\Sum_{b\in B} a_b 
 C^{b,i_1\dots i_{\mu -1}}_{g}(\Omega_1,\dots ,
\Omega_p,\phi_{u+1}, \phi_1,
\dots,\phi_u)\nabla_{i_1}\upsilon\dots\nabla_{i_{\mu-1}}\upsilon+
\\&\sum_{p\in P} a_p Xdiv_{i_\mu} C^{p,i_1\dots i_\mu}_g(\Omega_1,\dots ,
\Omega_p,\phi_{u+1}, \phi_1,
\dots,\phi_u)\nabla_{i_1}\upsilon\dots\nabla_{i_{\mu-1}}\upsilon=
\\&\Sum_{j\in J^{\sigma+u+1}} a_j
C^j_{g}(\Omega_1,\dots ,\Omega_p,\phi_1, \dots
,\phi_u,\phi_{u+1}, \upsilon^{\mu-1})=0,
\end{split}
\end{equation}
modulo complete contractions of length $\ge\sigma+u+\mu$ (the terms in 
the above have lengh $\sigma+u+\mu-1$). 
 Then, keeping track of the greater length correction terms 
 that arise in the above, we derive  a new equation:
    
\begin{equation}
\label{toki''}
\begin{split}
&\Sum_{b\in B} a_b 
 C^{b,i_1\dots i_{\mu -1}}_{g}(\Omega_1,\dots ,
\Omega_p,\phi_{u+1}, \phi_1,
\dots,\phi_u)\nabla_{i_1}\upsilon\dots\nabla_{i_{\mu-1}}\upsilon+
\\&\sum_{p\in P} a_p Xdiv_{i_\mu} C^{p,i_1\dots i_\mu}_g(\Omega_1,\dots ,
\Omega_p,\phi_{u+1}, \phi_1,
\dots,\phi_u)\nabla_{i_1}\upsilon\dots\nabla_{i_{\mu-1}}\upsilon+
\\&\Sum_{j\in J^{\sigma+u+1}} a_j
C^j_{g}(\Omega_1,\dots ,\Omega_p,\phi_1, \dots
,\phi_u,\phi_{u+1}, \upsilon^{\mu-1})=
\\&\Sum_{\zeta\in Z} a_\zeta
 C^{\zeta,i_1\dots i_{\mu -1}}_{g}(\Omega_1,\dots ,
\Omega_p, \phi_1,
\dots,\phi_u,\phi_{u+1})\nabla_{i_1}\upsilon\dots\nabla_{i_{\mu-1}}\upsilon+
\\&\sum_{j\in J} a_j C^{j,i_1\dots i_{\mu -1}}_{g}(\Omega_1,\dots ,
\Omega_p, \phi_1,
\dots,\phi_u,\phi_{u+1})\nabla_{i_1}\upsilon\dots\nabla_{i_{\mu-1}}\upsilon.
\end{split}
\end{equation}
 This equation holds {\it perfectly}--not modulo 
 longer terms. The terms indexed in $J$ have length
  $\sigma+u+\mu$, a factor $\nabla\omega$  and a $u$-simple 
  character $\vec{\kappa}_{simp}$; 
  the terms indexed in $T$, a factor $\nabla^{(B)}\phi_{u+1}$. 
  Thus, invoking the last Lemma in 
  the Appendix of \cite{alexakis1}, we derive: 
  
  \begin{equation}
\label{toki''}
\begin{split}
&\Sum_{b\in B} a_b Xdiv_{i_1}\dots Xdiv_{i_{\mu-1}}
 C^{b,i_1\dots i_{\mu -1}}_{g}(\Omega_1,\dots ,
\Omega_p,\phi_{u+1}, \phi_1,
\dots,\phi_u)+
\\&\sum_{p\in P} a_p Xdiv_{i_1}\dots Xdiv_{i_{\mu-1}}
Xdiv_{i_\mu} C^{p,i_1\dots i_\mu}_g(\Omega_1,\dots ,
\Omega_p,\phi_{u+1}, \phi_1,
\dots,\phi_u)+
\\&\Sum_{j\in J^{\sigma+u+1}} a_j
C^j_{g}(\Omega_1,\dots ,\Omega_p,\phi_1, \dots
,\phi_u,\phi_{u+1})=
\\&\Sum_{\zeta\in Z} a_\zeta
Xdiv_{i_1}\dots Xdiv_{i_{\mu-1}} C^{\zeta,i_1\dots i_{\mu -1}}_{g}(\Omega_1,\dots ,
\Omega_p, \phi_1,
\dots,\phi_u,\phi_{u+1})+
\\&\sum_{j\in J} a_j C^{j,i_1\dots i_{\mu -1}}_{g}(\Omega_1,\dots ,
\Omega_p, \phi_1,
\dots,\phi_u,\phi_{u+1}).
\end{split}
\end{equation}
 Thus, by operating on the above by the operation $\sum_{K=1}^{\sigma-1}Hit^K[\dots]$
 (this clearly produces a true equation), and then
interchanging the two functions $\phi_{u+1},\omega$ 
(this also produces a new true equation), we derive (\ref{tokicor}). 
\newline

{\it Proof of (\ref{tokicor'}):} We just neglect
 the algebraic structure of $\Sum_{K=1}^{\sigma-1}Hit^K_\omega$ in 
 the LHS and apply the Lemma 4.10 in \cite{alexakis4} to (\ref{toki}. 
We   use the fact that 
the LHS of the resulting equation vanishes formally at the linearized level, 
and then repeat the formal permutations of indices
to the non-linearized level, and finaly replace the
 $\mu-1$ factors $\nabla\upsilon$ by $Xdiv$'s (using the
  last Lemma in the Appendix of \cite{alexakis1}). $\Box$
\newline

\par Therefore, replacing (\ref{tokicor}) (or (\ref{tokicor'}))
 into (\ref{uptohere}) we derive:

\begin{equation}
\label{uptohere2}
\begin{split}
&\Sum_{l\in L} a_l Sub^{\sigma +u+2}_\omega \{Image
^{1,\beta}_{\phi_{u+1}} [Xdiv_{i_1}\dots Xdiv_{i_a}C^{l,i_1\dots
,i_a}_{g}(\Omega_1,\dots ,\Omega_p, \phi_1, \dots,\phi_u)]\}+
\\&Sub^{\sigma +u+2}_\omega \{
Image^{1,\beta}_{\phi_{u+1}} [C^j_{g}(\Omega_1,\dots ,
\Omega_p,\phi_1, \dots,\phi_u)]\}
\\&+\Sum_{l\in L_\mu} a_l \Sum_{i_h\in I^\sharp}Xdiv_{i_1}\dots \hat{Xdiv}_{i_h}
\dots Xdiv_{i_\mu}
\\&\{FRep^{i_h,1,\omega,\phi_{u+1}}[C^{l,i_1\dots i_a|f(i_h)}_{g}
(\Omega_1,\dots ,\Omega_p,\phi_1,\dots ,\phi_u)]
\\& + \{FRep^{i_h,1,\phi_{u+1},\omega}[C^{l,i_1\dots i_a|
f(i_h)}_{g}(\Omega_1,\dots ,\Omega_p,\phi_1,\dots ,\phi_u)]
\\& + FRep^{i_h,2,\omega,\phi_{u+1}}[C^{l,i_1\dots i_a
|f(i_h)}_{g}(\Omega_1,\dots ,\Omega_p,\phi_1,\dots ,\phi_u)] +\\&
 FRep^{i_h,2,\phi_{u+1},\omega}
[C^{l,i_1\dots i_a|f(i_h)}_{g}(\Omega_1,\dots ,\Omega_p,
\phi_1,\dots ,\phi_u)]\}\}+
\\&(\Sum_{m\in M_{\mu-1}} a_m Xdiv_{i_1}\dots
Xdiv_{i_{\mu-1}}C^{m,i_1\dots i_{\mu+1}}_{g} (\Omega_1,\dots
,\Omega_p,\phi_1,\dots ,\phi_u)
\nabla_{i_\mu}\phi_{u+1}\nabla_{i_{\mu+1}}\omega)
\\& +\Sum_{m\in M} a_m Xdiv_{i_1}\dots Xdiv_{i_a} C^{m,i_1\dots i_{a+2}}_{g}
(\Omega_1,\dots ,\Omega_p,\phi_1,\dots ,\phi_u)
\nabla_{i_{a+1}}\phi_{u+1}\nabla_{i_{a+2}}\omega+
\\& \Sum_{\zeta\in Z} a_\zeta Xdiv_{i_1}\dots Xdiv_{\mu-1}
Hit^K_\omega[C^{\zeta,i_1\dots i_{\mu -1}}_{g}(\Omega_1,\dots
,\Omega_p, \phi_1, \dots,\phi_u,\phi_{u+1})]+
\\&\Sum_{j\in
J^{\sigma+u+2}} a_j C^j_{g}(\Omega_1,\dots ,\Omega_p,\phi_1, \dots
,\phi_u,\phi_{u+1},\omega)
 \\&= \Sum_{t\in T} a_t
C^t_{g}(\Omega_1,\dots ,\Omega_p, \phi_1,
\dots,\phi_u,\phi_{u+1},\omega),
\end{split}
 \end{equation}
modulo complete contractions of length $\ge\sigma +u+3$.
\newline

\par Hence, we are reduced to studying the sublinear combinations
\begin{equation}
\label{montazeri} Sub^{\sigma +u+2}_\omega
\{Image^{1,\beta}_{\phi_{u+1}} [Xdiv_{i_1}\dots Xdiv_{i_a}
C^{l,i_1\dots ,i_a}_{g}(\Omega_1,\dots ,\Omega_p, \phi_1,
\dots,\phi_u)]\}
\end{equation}
and
$$Sub^{\sigma +u+2}_\omega \{
Image^{1,\beta}_{\phi_{u+1}} [C^j_{g}(\Omega_1,\dots ,
\Omega_p,\phi_1, \dots,\phi_u)]\}.$$

\par As before, we straightforwardly derive:

\begin{equation}
\label{arxidhs}
\begin{split}& Sub^{\sigma +u+2}_\omega \{
Image^{1,\beta}_{\phi_{u+1}} [C^j_{g}(\Omega_1,\dots ,
\Omega_p,\phi_1, \dots,\phi_u)]\}\\&=\Sum_{j\in J^{\sigma +u+2}} a_j
C^j_{g}(\Omega_1,\dots , \Omega_p,\phi_1,
\dots,\phi_u,\phi_{u+1},\omega).
\end{split}\end{equation}

\par To analyze the sublinear combination
(\ref{montazeri}) we firstly seek to understand how it arizes:
\newline

{\bf A study of the sublinear combination \\ $Sub^{\sigma
+u+2}_\omega \{Image^{1,\beta}_{\phi_{u+1}} [Xdiv_{i_1}\dots
Xdiv_{i_a} C^{l,i_1\dots ,i_a}_{g}(\Omega_1,\dots ,\Omega_p,
\phi_1, \dots,\phi_u)]\}$:}
\newline

As before, we write out
$Xdiv_{i_1}\dots Xdiv_{i_a} C^{l,i_1\dots
,i_a}_{g}(\Omega_1,\dots ,\Omega_p, \phi_1, \dots,\phi_u)$ as a
linear combination of complete contractions, say

\begin{equation}
\label{poonen2} 
Xdiv_{i_1}\dots Xdiv_{i_a}
C^{l,i_1\dots ,i_a}_{g}=\sum_{x\in X} a_x C^x_g.
\end{equation}
Then, for each $C^x_g$ we identify the {\it (ordered)
sets of pairs of pairs of indices}, $[({}_a,{}_b), ({}_c,{}_d)]$ where ${}_a,{}_b$
 belong to the same factor and ${}_c,{}_d$ belong to the same factor, and either 
 ${}_a$ is contracting against ${}_c$ and ${}_b$ against
 ${}_d$ on vice versa, {\it and both ${}_a,{}_c$ are derivative indices}.
Denote this set of ordered pairs by $Z^x$. Then, for
each $[({}_a,{}_b), ({}_c,{}_d)]\in Z^x$ we let $B_{[(a,b), (c,d)]}\{C^x_g\}$
stand for the complete contraction that formally arises from
 $C^x_g$ by applying the last summand in (\ref{levicivita}) to
the indices $({}_a,{}_b)$ (recall that one of them
 is a derivative index, so this is a well-defined operation),
 thus making the indices ${}_c,{}_d$ contract against each other. Then 
apply $Sub_\omega$ to this complete contraction we have created.
This replaces the internal contraction between ${}_c,{}_d$ by a factor $\nabla\omega$
(since  ${}_c$ is a derivative index). Denote the complete
contraction we thus obtain by $\overline{B}_{[(a,b), (c,d)]}\{C^x_g\}$. It follows by 
the definition of $Sub_\omega^{\sigma+u+2}\{\dots\}$ that:

\begin{equation}
\label{omorfoula}
\begin{split}
&Sub^{\sigma +u+2}_\omega \{Image^{1,\beta}_{\phi_{u+1}}
[Xdiv_{i_1}\dots Xdiv_{i_a} C^{l,i_1\dots ,i_a}_{g}(\Omega_1,\dots
,\Omega_p, \phi_1, \dots,\phi_u)]\}=
\\&\sum_{x\in X} a_x \sum_{[(a,b), (c,d)]\in Z^x}\overline{B}_{[(a,b), (c,d)]}\{C^x_g\}.
\end{split}
\end{equation}

\par Having obtained an understanding of
 how $$Sub^{\sigma +u+2}_\omega \{Image^{1,\beta}_{\phi_{u+1}}
[Xdiv_{i_1}\dots Xdiv_{i_a} C^{l,i_1\dots ,i_a}_{g}(\Omega_1,\dots
,\Omega_p, \phi_1, \dots,\phi_u)]\}$$ arises, we now proceed to
express it in a more useful form:

\par We distinguish cases
 depending on the {\it form} of the indices
 $({}_a,{}_b), ({}_c,{}_d)$, for each complete contraction
$C_{g}(\Omega_1,\dots ,\Omega_p,\phi_1, \dots,\phi_u)$ appearing in
$$Xdiv_{i_1}\dots Xdiv_{i_a}
C^{l,i_1\dots ,i_a}_{g}(\Omega_1,\dots ,\Omega_p,\phi_1,
\dots,\phi_u).$$
 Recall that we have defined an index in $C_{g}$ to be a {\it divergence index} if
it is an index $\nabla^{i_f}$ which has arisen by taking an $Xdiv$ operation,  $Xdiv_{i_f}$,
with respect to some free index ${}_{i_f}$; we
have also defined an index in $C_g$ to be
 an {\it original index} in $C^{l,i_1\dots ,i_a}_{g}$ if the
index appears in the tensor
 field $C^{l,i_1\dots ,i_a}_{g}$ (before we take any $Xdiv$'s).

\par Now, we place each complete contraction in
$$Sub^{\sigma +u+2}_\omega \{Image^{1,\beta}_{\phi_{u+1}}
[Xdiv_{i_1}\dots Xdiv_{i_a} C^{l,i_1\dots ,i_a}_{g}(\Omega_1,\dots
,\Omega_p, \phi_1, \dots,\phi_u)]\}$$
 into one of the sublinear combinations
$$Sub^{\sigma +u+2,K}_\omega \{Image^{1,\beta}_{\phi_{u+1}}
[Xdiv_{i_1}\dots Xdiv_{i_a} C^{l,i_1\dots ,i_a}_{g}(\Omega_1,\dots
,\Omega_p, \phi_1, \dots,\phi_u)]\}$$
($K=\alpha,\beta,\gamma,\delta$) based on the pair $(a,b), (c,d)$
from which it arose.

Specifically:

\begin{definition}
\label{breakglass}
Refer to (\ref{omorfoula}) and pick out a 
term in the RHS. For any given index 
 ${}_a,{}_b,{}_c,{}_d$ (recall that we are now assuming that ${}_a,{}_c$ are
 derivative indices), we inquire whether it is an original
 index or a divergence  index $\nabla^{i_h},h=1,\dots ,a$.
 Accordingly, we place the term $\overline{B}_{[(a,b), (c,d)]}\{C^x_g\}$ 
into one of the four sublinear combinations $Sub_\omega^{\sigma+u+2,\alpha}$, 
$Sub_\omega^{\sigma+u+2,\beta}$, $Sub_\omega^{\sigma+u+2,\gamma}$, 
$Sub_\omega^{\sigma+u+2,\delta}$ according to the following rule: 

We declare that $C_{g}$
 belongs to $Sub^{\sigma +u+2,\alpha}_\omega \{
 Image^{1,\beta}_{\phi_{u+1}} [Xdiv_{i_1}\dots Xdiv_{i_a}
C^{l,i_1\dots ,i_a}_{g}]\}$
 if and only if ${}_a$ and ${}_c$ are divergence indices.

\par We declare that $C_{g}$
 belongs to $Sub^{\sigma +u+2,\beta}_\omega
 \{Image^{1,\beta}_{\phi_{u+1}} [Xdiv_{i_1}\dots Xdiv_{i_a}
C^{l,i_1\dots ,i_a}_{g}]\}$ if and only if only one of the indices
${}_a,{}_b,{}_c,{}_d$ is a divergence index (say ${}_a$ with no loss of
generality). We declare that $C_{g}$
 belongs to \\$Sub^{\sigma +u+2,\gamma}_\omega
 \{Image^{1,\beta}_{\phi_{u+1}} [Xdiv_{i_1}\dots Xdiv_{i_a}
C^{l,i_1\dots ,i_a}_{g}]\}$  if either ${}_a,{}_b$ or ${}_c,{}_d$ are
 both divergence indices. Finally, we declare that $C_{g}$
 belongs to \\$Sub^{\sigma +u+2,\delta}_\omega
 \{Image^{1,\beta}_{\phi_{u+1}} [Xdiv_{i_1}\dots Xdiv_{i_a}
C^{l,i_1\dots ,i_a}_{g}]\}$  if all four indices ${}_a,{}_b,{}_c,{}_d$
 are original indices.
\end{definition}

\par Now, another piece of notation: We denote by
$$ \Sum_{m\in M^{\sharp\sharp}} a_m
C^{m,i_1\dots i_{a+2}}_{g} (\Omega_1,\dots ,\Omega_p, \phi_1,\dots
,\phi_u) \nabla_{i_{a+1}}\phi_{u+1}\nabla_{i_{a+2}}\omega$$
a generic linear combination of tensor fields of length
$\sigma+u+2$ with {\it two} un-normalized factors
$\nabla\Omega_h,\nabla\Omega_{h'}$, that are contracting against
factors $\nabla\phi_{u+1}$, $\nabla\omega$ respectively.
We also require that if $a=\mu$ then all free indices must be non-special.
 Now, the
first thing we easily notice is that for each $l\in L$:

\begin{equation}
\label{maaa}
\begin{split}
&Sub^{\sigma +u+2,\delta}_\omega \{ Image^{1,\beta}_{\phi_{u+1}}
[Xdiv_{i_1}\dots Xdiv_{i_a} C^{l,i_1\dots ,i_a}_{g}]\}=\Sum_{m\in M\bigcup M^\sharp\bigcup M^{\sharp\sharp}}
\\&Xdiv_{i_1}\dots Xdiv_{i_a}
a_m C^{m,i_1\dots i_{a+2}}_{g} (\Omega_1,\dots ,\Omega_p,
\phi_1,\dots ,\phi_u)
\nabla_{i_{a+1}}\phi_{u+1}\nabla_{i_{a+2}}\omega.
\end{split}
\end{equation}

\par In order to describe $Sub^{\sigma +u+2,\alpha}_\omega
\{ Image^{1,\beta}_{\phi_{u+1}} [Xdiv_{i_1}\dots Xdiv_{i_a}
C^{l,i_1\dots ,i_a}_{g}]\}$, we define $(I\times I)^\sharp$ to
stand for the subset of $(I\times I)$ that consists of
 all ordered pairs of free indices that belong to different
 factors. For each $l\in L$ we then compute:

\begin{equation}
\label{maaa2}
\begin{split}
&Sub^{\sigma +u+2,\alpha}_\omega \{ Image^{1,\beta}_{\phi_{u+1}}
[Xdiv_{i_1}\dots Xdiv_{i_a} C^{l,i_1\dots ,i_a}_{g}]\}=
\\&\Sum_{(i_k,i_l)\in (I\times I)^\sharp}
Xdiv_{i_1}\dots \hat{Xdiv}_{i_k}\dots\hat{Xdiv}_{i_l}\dots
Xdiv_{i_a}C^{l,i_1\dots ,i_a}_{g}
\nabla_{i_k}\phi_{u+1}\nabla_{i_l}\omega.
\end{split}
\end{equation}

\par In particular, we observe that if
$l\in L_K, K\ge \mu +2$, the right hand side of the above is a
 generic linear combination of the form:

\begin{equation}
\label{mhnmoupeis} \begin{split} &\Sum_{u\in U_1} a_u
Xdiv_{i_1}\dots Xdiv_{i_a} C^{u,i_1\dots ,i_a}_{g}(\Omega_1,\dots
,\Omega_p,\phi_1, \dots,\phi_u,\phi_{u+1},\omega)+
\\&\Sum_{u\in U_2} a_u
Xdiv_{i_1}\dots Xdiv_{i_a} C^{u,i_1\dots ,i_a}_{g}(\Omega_1,\dots
,\Omega_p,\phi_1, \dots,\phi_u,\phi_{u+1},\omega).
\end{split}
\end{equation}

\par Now, to describe each $Sub^{\sigma +u+2,\beta}_\omega
\{ Image^{1,\beta}_{\phi_{u+1}} [Xdiv_{i_1}\dots Xdiv_{i_a}
C^{l,i_1\dots ,i_a}_{g}]\}$, we introduce more notation:
For each $l\in L_\mu$ and each
${}_{i_h}\in I_l$,\footnote{Recall that $I_l$ stands for
the set of free indices in the tensor field $C^{l,i_1\dots i_a}_g$.}
 we define $T(i_h)$ to
stand for the factor to which ${}_{i_h}$ belongs.
 We observe that if $l\in L_\mu$ then for each factor $T(i_h)$
of the form $T(i_h)=\nabla^{(m)}R_{ijkl}$ or
$T(i_h)=\nabla^{(A)}\Omega_f$, ${}_{i_h}$ must be a derivative index,
since in the setting of Lemma \ref{pskovb} no $\mu$-tensor field contains special free indices.
 If $T(i_h)$ is of the form $S_{*}\nabla^{(\nu)} R_{ijkl}$, we
 then have that ${}_{i_h}$ must be one of the indices
${}_{r_1},\dots ,{}_{r_\nu},{}_j$ (by the first assumption in the
introduction). In that case, we write out $T(i_h)$ as a sum of
tensors of the form $\nabla^{(\nu)} R_{ijkl}$.

\begin{definition}
\label{intoclub}
\par With the above convention, for each $l\in L$,\footnote{Recall
 that $L=L_\mu\bigcup L_{>\mu}$ is the index set of the tensor
fields $C^{l,i_1\dots i_a}_g$ in our Lemma hypothesis
(\ref{hypothese2}).}
 we denote by $I^{*}_l\subset I_l$\footnote{Recall that $I_l$ stands for
the set of free indices in $C^{l,i_1\dots i_a}$.} the
set of free indices that are derivative indices. We denote by
$I^{+}_l=I_l\setminus I^{*}_l$.\footnote{Observe that the indices that belong to
$I^{+}_l$ will be the index ${}_j$ in some factor $\nabla^{(\nu)}
R_{ijkl}$ that has arisen from  a de-symmetrization as above.} For
each ${}_i\in I^{*}_l$, we denote by $Set(T(i))$\footnote{Recall that $T(i)$ stands
 for the factor in $C^{l,i_1\dots i_a}_g$  to which the index ${}_i$ belongs.}
 to be the set of all
the indices in $T(i)$ that are not free and not contracting
against a factor $\nabla\phi_h$. For each ${}_i\in I^{+}_l$, we
denote by $Set(T(i))$ the set of derivative indices in the factor
$\nabla^{(\nu)} R_{ijkl}$ that are not free and not contracting
against a factor $\nabla\phi_h$.

\par Then, for each ${}_i\in I^{*}_l$ and each $t\in Set(T(i))$,
 let
$Repla^{i,t}_{\phi_{u+1},\omega}[C^{l,i_1\dots ,i_a}_{g}]$ be the
$(a-1)$ tensor field that formally arises from $C^{l,i_1\dots ,i_a}_{g}$ by
erasing the index ${}_i$ and making the index ${}_t$ contract
against a factor $\nabla\phi_{u+1}$ and also making the index
${}^t$ contract against a factor $\nabla\omega$. We denote by
$Repla^{i,t}_{\omega,\phi_{u+1}}[C^{l,i_1\dots ,i_a}_{g}]$ the
$(a-1)$-tensor field that arises from
$Repla^{i,t}_{\phi_{u+1},\omega}[C^{l,i_1\dots ,i_a}_{g}]$ by
switching $\phi_{u+1}$ and $\omega$.

\par For each $i_h\in I^{+}_l$ and each $t\in Set(T(i))$,
we denote by $Repla^{i,t}_{\phi_{u+1},\omega} [C^{l,i_1\dots
,i_a}_{g}]$ the $(a-1)$ tensor field that
 arises from $C^{l,i_1\dots i_a}_{g}$
by erasing the index ${}_t$ and making the index ${}_i$ contract
against a factor $\nabla\phi_{u+1}$. We also make the index ${}^t$
contract against a factor $\nabla\omega$. We again denote by
$Repla^{i,t}_{\omega,\phi_{u+1}}[C^{l,i_1\dots i_a}_{g}]$ the
$(a-1)$-tensor field that arises from
$Repla^{i,t}_{\phi_{u+1},\omega} [C^{l,i_1\dots i_a}_{g}]$ by
switching $\phi_{u+1}$ and $\omega$.
\end{definition}

\par We then calculate that for each $l\in L_\mu$:

\begin{equation}
\label{beta}
\begin{split}
&Sub^{\sigma +u+2,\beta}_\omega \{ Image^{1,\beta}_{\phi_{u+1}}
[Xdiv_{i_1}\dots Xdiv_{i_a} C^{l,i_1\dots ,i_a}_{g}]\}=\Sum_{i_h\in I^{*}_l\bigcup I^{+}_l}\Sum_{t\in Set[T(i)]}
\\&Xdiv_{i_1}\dots \hat{Xdiv}_{i_h}\dots Xdiv_{i_a}\{
Repla^{i,t}_{\phi_{u+1},\omega}[C^{l,i_1\dots ,i_a}_{g}]
+Repla^{t,i}_{\phi_{u+1},\omega}[C^{l,i_1\dots ,i_a}_{g}]\}.
\end{split}
\end{equation}

{\it Convention:} For future reference, we will further subdivide the
index sets $I^{*}_l,I^{+}_l$: If there is a unique selected factor
we define $I^{*,1}_l=I^{*}_l \bigcap I_1$ and analogously
$I^{*,2}_l=I^{*}_l\bigcap I_2$ and similarly for $I^{+}_l$.
(Recall that $I_1$ (or $I_2$) stand for the sets of free indices
that belong (do not belong) to the selected factor, when the
selected factor is unique). If there are multiple selected factors
$\{T_i\}_{i=1}^{b_l}$, we define $I^{*,1,T_i}_l=I^{*}_l \bigcap
I_1^{T_i}$ and analogously $I^{*,2,T_i}_l=I^{*}_l\bigcap
I^{T_i}_2$. (Recall that
 $I^{T_i}_1$ is the set of free indices that belong to the
 selected factor $T_i$ and $I^{T_i}_2$ is the set of free
indices that do not belong to the selected factor $T_i$).
\newline

\par Analogously, we deduce that for each $l\in L\setminus
 L_\mu$:

\begin{equation}
\label{beta2}
\begin{split}
&Sub^{\sigma +u+2,\beta}_\omega \{ Image^{1,\beta}_{\phi_{u+1}}
[Xdiv_{i_1}\dots Xdiv_{i_a} C^{l,i_1\dots ,i_a}_{g}]\}=\Sum_{m\in M\bigcup M^\sharp} a_m
\\&Xdiv_{i_1}\dots Xdiv_{i_a} C^{u,i_1\dots ,i_a, i_{a+1},
i_{a+2}}_{g}(\Omega_1,\dots ,\Omega_p,\phi_1,
\dots,\phi_u)\nabla_{i_{a+1}}\phi_{u+1}\nabla_{i_{a+2}}\omega.
\end{split}
\end{equation}

\par Finally, we seek to understand
$Sub^{\sigma +u+2,\gamma}_\omega \{ Image^{1,\beta}_{\phi_{u+1}}
[Xdiv_{i_1}\dots Xdiv_{i_a} C^{l,i_1\dots ,i_a}_{g}]\}$.

\begin{definition}
\label{definition2}
 For each $C^{l,i_1\dots i_a}_{g}$, we denote by $I^d_l$
the set
 of pairs of free indices that belong to the same factor,
such that at least one of them is a derivative index.
\end{definition}

\par Now, for each $l\in L$ and each $({}_{i_k},{}_{i_l})\in I^d_l$,
we assume with no loss of generality that ${}_{i_k}$ is a derivative
index. We also denote by $\{ F_1,\dots ,F_{\sigma-1}\}$ the set of
real factors (i.e. factors that are not in the form
$\nabla\phi_h$) in $C^{l,i_1\dots ,i_a}_{g}$ {\it other
 than} the factor to which ${}_{i_k},{}_{i_l}$ belong. We then denote by
$Re^{K,\phi_{u+1},\omega}_{i_k,i_l}[C^{l,i_1\dots ,i_a}_{g}]$ the
$(a-1)$-tensor field that arises from $C^{l,i_1\dots ,i_a}_{g}$ by
erasing ${}_{i_k}$, contracting ${}_{i_l}$ against a factor
$\nabla\phi_{u+1}$ and then hitting  the factor $F_K$ by a
derivative $\nabla_z$ and contracting ${}_z$ against a factor
$\nabla^z\omega$. We denote by
$Re^{K,\omega,\phi_{u+1}}_{i_k,i_l}[C^{l,i_1\dots ,i_a}_{g}]$ the
$(a-1)$-tensor field that arises from
$Re^{K,\omega,\phi_{u+1}}_{i_k,i_l}[C^{l,i_1\dots ,i_a}_{g}]$ by
switching $\phi_{u+1}$ and $\omega$.
 We then calculate that for each $l\in L$:

\begin{equation}
\label{gamma1}
\begin{split}
&Sub^{\sigma +u+2,\gamma}_\omega \{ Image^{1,\beta}_{\phi_{u+1}}
[Xdiv_{i_1}\dots Xdiv_{i_a} C^{l,i_1\dots ,i_a}_{g}]\}=
\\&\Sum_{(i_k,i_l)\in I^d_l} \Sum_{K=1}^{\sigma
-1}Xdiv_{i_1}\dots\hat{Xdiv}_{i_k}\dots \hat{Xdiv}_{i_l}\dots
Xdiv_{i_a}
\\&\{ Re^{K,\omega,\phi_{u+1}}_{i_k,i_l}[C^{l,i_1\dots
,i_a}_{g}]+Re^{K,\phi_{u+1},\omega}_{i_k,i_l}[C^{l,i_1\dots
,i_a}_{g}]\}.
\end{split}
\end{equation}

\par In conclusion, we have shown that:

\begin{equation}
\label{finalconclusion1}
\begin{split}
&Sub^{\sigma +u+2}_\omega \{ Image^{1,\beta}_{\phi_{u+1}}
L_{g}(\Omega_1,\dots ,\Omega_p,\phi_1, \dots,\phi_u)\}=
\\&\Sum_{l\in L} a_l \Sum_{(i_k,i_l)\in (I\times I)^\sharp}
Xdiv_{i_1}\dots \hat{Xdiv}_{i_k}\dots\hat{Xdiv}_{i_l}\dots
Xdiv_{i_a}C^{l,i_1\dots ,i_a}_{g}
\nabla_{i_k}\phi_{u+1}\nabla_{i_l}\omega+
\\&\Sum_{l\in L_\mu} a_l
\Sum_{i_h\in I^{*}_l\bigcup I^{+}_l}\Sum_{t\in T(i_h)}
Xdiv_{i_1}\dots \hat{Xdiv}_{i_h}\dots Xdiv_{i_a}\\&\{
Repla^{i,t}_{\phi_{u+1},\omega}[C^{l,i_1\dots ,i_a}_{g}]
+Repla^{t,i}_{\phi_{u+1},\omega}[C^{l,i_1\dots ,i_a}_{g}]\}
\\&+\Sum_{m\in M\bigcup M^\sharp\bigcup M^{\sharp\sharp}}
a_m C^{m,i_1\dots i_{a+2}}_{g} (\Omega_1,\dots ,\Omega_p,
\phi_1,\dots ,\phi_u)
\nabla_{i_{a+1}}\phi_{u+1}\nabla_{i_{a+2}}\omega+
\\&\Sum_{l\in L} a_l \Sum_{(i_k,i_l)\in I^d_l} \Sum_{K=1}^{\sigma
-1}Xdiv_{i_1}\dots\hat{Xdiv}_{i_k}\dots \hat{Xdiv}_{i_l}\dots
Xdiv_{i_a}\\&\{ Re^{K,\omega,\phi_{u+1}}_{i_k,i_l}[C^{l,i_1\dots
,i_a}_{g}]+Re^{K,\phi_{u+1},\omega}_{i_k,i_l}[C^{l,i_1\dots
,i_a}_{g}]\}
\\&+\Sum_{j\in J^{\sigma +u+2}} a_j C^j_{g}(\Omega_1,\dots ,
\Omega_p,\phi_1,\dots ,\phi_u,\phi_{u+1},\omega).
\end{split}
\end{equation}

\par Replacing the above into (\ref{uptohere2}) we obtain a
new equation, after all this extensive analysis of the equation
$Image^{1,\beta}_{\phi_{u+1}}[L_{g}]=0$:

\begin{equation}
\label{profinalconclusion2}
\begin{split}
&\Sum_{l\in L} a_l \Sum_{(i_k,i_l)\in (I\times I)^\sharp}
Xdiv_{i_1}\dots \hat{Xdiv}_{i_k}\dots\hat{Xdiv}_{i_l}\dots
Xdiv_{i_a}C^{l,i_1\dots ,i_a}_{g}
\nabla_{i_k}\phi_{u+1}\nabla_{i_l}\omega
\\&+\Sum_{l\in L_\mu} a_l
\Sum_{i_h\in I^{*}_l\bigcup I^{+}_l}\Sum_{t\in T(i_h)}
Xdiv_{i_1}\dots \hat{Xdiv}_{i_h}\dots Xdiv_{i_a}
\\&\{ Repla^{i,t}_{\phi_{u+1},\omega}[C^{l,i_1\dots ,i_a}_{g}]
+Repla^{t,i}_{\phi_{u+1},\omega}[C^{l,i_1\dots ,i_a}_{g}]\}
\\&+\Sum_{l\in L_\mu} a_l \Sum_{i_f\in I^\sharp}
Xdiv_{i_1}\dots \hat{Xdiv}_{i_f}\dots Xdiv_{i_\mu} 
\\&\{FRep^{i_h,1,\omega,\phi_{u+1}} [C^{l,i_1\dots
i_a|f(i_h)}_{g}(\Omega_1,\dots ,\Omega_p, \phi_1,\dots ,\phi_u)]
\\& +FRep^{i_h,1,\phi_{u+1},\omega}[C^{l,i_1\dots i_a|
f(i_h)}_{g}(\Omega_1,\dots ,\Omega_p,\phi_1,\dots ,\phi_u)]
\\& + FRep^{i_h,2,\omega,\phi_{u+1}}[C^{l,i_1\dots i_a
|f(i_h)}_{g}(\Omega_1,\dots ,\Omega_p,\phi_1,\dots ,\phi_u)] +\\&
FRep^{i_h,2,\phi_{u+1},\omega}[C^{l,i_1\dots i_a
|f(i_h)}_{g}(\Omega_1,\dots ,\Omega_p,\phi_1,\dots , \phi_u)]\}+
\\&(\Sum_{m\in M_{\mu-1}} a_m Xdiv_{i_1}\dots
Xdiv_{i_{\mu-1}}C^{m,i_1\dots i_{\mu-1}}_{g} (\Omega_1,\dots
,\Omega_p,\phi_1,\dots , \phi_u)\nabla_{i_\mu}
\phi_{u+1}\nabla_{i_{\mu+1}}\omega)
\\&+\Sum_{m\in M\bigcup M^\sharp\bigcup M^{\sharp\sharp}}
a_m C^{m,i_1\dots i_{a+2}}_{g} (\Omega_1,\dots ,\Omega_p,
\phi_1,\dots ,\phi_u)
\nabla_{i_{a+1}}\phi_{u+1}\nabla_{i_{a+2}}\omega
\\&+\Sum_{l\in L} a_l \Sum_{(i_k,i_l)\in I^d_l} \Sum_{K=1}^{\sigma
-1}Xdiv_{i_1}\dots\hat{Xdiv}_{i_k}\dots \hat{Xdiv}_{i_l}\dots
Xdiv_{i_a}\\&\{ Re^{K,\omega,\phi_{u+1}}_{i_k,i_l}[C^{l,i_1\dots
,i_a}_{g}]+Re^{K,\phi_{u+1},\omega}_{i_k,i_l}[C^{l,i_1\dots
,i_a}_{g}]\}+
\\&\Sum_{\zeta\in Z} a_\zeta Xdiv_{i_1}\dots Xdiv_{\mu-1}
Hit^K_\omega[C^{\zeta,i_1\dots i_{\mu -1}}_{g}(\Omega_1,\dots
,\Omega_p, \phi_1, \dots,\phi_u,\phi_{u+1})]+
\\&+\Sum_{j\in J^{\sigma +u+2}} a_j C^j_{g}
(\Omega_1,\dots ,
\Omega_p,\phi_1,\dots ,\phi_u,\phi_{u+1},\omega)=
\\&\Sum_{t\in T^{\sigma+u+2}} a_t C^t_{g}(\Omega_1,\dots ,
\Omega_p,\phi_1,\dots ,\phi_u,\phi_{u+1},\omega),
\end{split}
\end{equation}
modulo complete contractions of length $\ge\sigma+u+3$. The sublinear 
combination $\sum_{m\in M_{\mu-1}}\dots$ appears only 
in the special subcase of case B. The linear
combination on the RHS stands for generic notation
 (see the notational convention introduced after (\ref{baldouin'})).

\par In fact, we observe that the minimum length of the complete
contractions above is $\sigma+u+2$, and that all terms on the LHS have
 two factors $\nabla\phi_{u+1},\nabla\omega$, while each term on the RHS has at least one
term $\nabla^{(A)}\phi_{u+1}$ or $\nabla^{(A)}\omega$, with $A\ge 2$.

Therefore, since the above holds formally, we derive:

\begin{equation}
\label{finalconclusion2}
\begin{split}
&\Sum_{l\in L} a_l \Sum_{(i_k,i_l)\in (I\times I)^\sharp}
Xdiv_{i_1}\dots \hat{Xdiv}_{i_k}\dots\hat{Xdiv}_{i_l}\dots
Xdiv_{i_a}C^{l,i_1\dots ,i_a}_{g}
\nabla_{i_k}\phi_{u+1}\nabla_{i_l}\omega+
\\&\Sum_{l\in L_\mu} a_l
\Sum_{i_h\in I^{*}_l\bigcup I^{+}_l}\Sum_{t\in T(i_h)}
Xdiv_{i_1}\dots \hat{Xdiv}_{i_h}\dots Xdiv_{i_a}
\\&\{Repla^{i,t}_{\phi_{u+1},\omega}[C^{l,i_1\dots ,i_a}_{g}]
+Repla^{t,i}_{\phi_{u+1},\omega}[C^{l,i_1\dots ,i_a}_{g}]\}
\\&+\Sum_{l\in L_\mu} a_l \Sum_{i_f\in I^\sharp}
Xdiv_{i_1}\dots \hat{Xdiv}_{i_f}\dots Xdiv_{i_\mu} \\&\{
FRep^{i_h,1,\omega,\phi_{u+1}} [C^{l,i_1\dots
i_a|f(i_h)}_{g}(\Omega_1,\dots ,\Omega_p, \phi_1,\dots ,\phi_u)]
\\& +FRep^{i_h,1,\phi_{u+1},\omega}[C^{l,i_1\dots i_a|
f(i_h)}_{g}(\Omega_1,\dots ,\Omega_p,\phi_1,\dots ,\phi_u)]
\\& + FRep^{i_h,2,\omega,\phi_{u+1}}[C^{l,i_1\dots i_a
|f(i_h)}_{g}(\Omega_1,\dots ,\Omega_p,\phi_1,\dots ,\phi_u)] +\\&
FRep^{i_h,2,\phi_{u+1},\omega}[C^{l,i_1\dots i_a
|f(i_h)}_{g}(\Omega_1,\dots ,\Omega_p,\phi_1,\dots , \phi_u)]\}+
\\&(\Sum_{m\in M_{\mu-1}} a_m Xdiv_{i_1}\dots
Xdiv_{i_{\mu-1}}C^{m,i_1\dots i_{\mu-1}}_{g} (\Omega_1,\dots
,\Omega_p,\phi_1,\dots , \phi_u)\nabla_{i_\mu}
\phi_{u+1}\nabla_{i_{\mu+1}}\omega)
\\&+\Sum_{m\in M\bigcup M^\sharp\bigcup M^{\sharp\sharp}}
a_m C^{m,i_1\dots i_{a+2}}_{g} (\Omega_1,\dots ,\Omega_p,
\phi_1,\dots ,\phi_u)
\nabla_{i_{a+1}}\phi_{u+1}\nabla_{i_{a+2}}\omega
\\&+\Sum_{l\in L} a_l \Sum_{(i_k,i_l)\in I^d_l} \Sum_{K=1}^{\sigma
-1}Xdiv_{i_1}\dots\hat{Xdiv}_{i_k}\dots \hat{Xdiv}_{i_l}\dots
Xdiv_{i_a}\\&\{ Re^{K,\omega,\phi_{u+1}}_{i_k,i_l}[C^{l,i_1\dots
,i_a}_{g}]+Re^{K,\phi_{u+1},\omega}_{i_k,i_l}[C^{l,i_1\dots
,i_a}_{g}]\}+
\\&\Sum_{\zeta\in Z} a_\zeta Xdiv_{i_1}\dots Xdiv_{\mu-1}
Hit^K_\omega[C^{\zeta,i_1\dots i_{\mu -1}}_{g}(\Omega_1,\dots
,\Omega_p, \phi_1, \dots,\phi_u,\phi_{u+1})]+
\\&+\Sum_{j\in J^{\sigma +u+2}} a_j C^j_{g}
(\Omega_1,\dots ,
\Omega_p,\phi_1,\dots ,\phi_u,\phi_{u+1},\omega)=0,
\end{split}
\end{equation}
modulo complete contractions of length $\ge\sigma+u+3$.

\par We denote the above equation by:

\begin{equation}
\label{antrikos} Im^{1,\beta}_{\phi_{u+1}}[L_{g}]=0,
\end{equation}
for short. We repeat that the contractions appearing  in the above
equation all have length $\sigma +u+2$, and the equation holds
modulo complete
 contractions of length $\ge\sigma+u+3$.
\newline

{\bf The operation $Soph$:} We now define a formal operation $Soph$ that acts on the complete
contractions above: For each 
$C_{g}(\Omega_1,\dots , \Omega_p,\phi_1,\dots
,\phi_u,\phi_{u+1},\omega)$ (with factors $\nabla\phi_{u+1}$,
$\nabla\omega$ which are necessarily contracting against different factors),
 we first replace the two factors $\nabla_i\phi_{u+1},
\nabla_j\omega$ by a factor $g_{ij}$. Then, we add a derivative
index $\nabla_u$ onto the selected factor and contract it against
a factor $\nabla^u\phi_{u+1}$ (if there are
 multiple selected
 factors we perform the same operation for each of them and then add).
 Finally, we multiply the complete
contraction by a factor $\frac{1}{2}$. 
This definition extends to tensor fields and linear combinations.

\par Observe that when this operation acts on the complete contractions in 
$Im^{1,\beta}_{\phi_{u+1}}[L_{g}]$, it produces
 complete contractions of length $\sigma +u+1$
with a factor $\nabla\phi_{u+1}$ and with a weak character
$Weak(\vec{\kappa}^{+}_{simp})$.
\newline

\par Observe that
since (\ref{antrikos}) holds formally, it follows that:

\begin{equation}
\label{antrikos2} Soph\{ Im^{1,\beta}[L_{g}]\}=0,
\end{equation}
modulo complete contractions of length $\ge\sigma +u+2$.

\subsection{Preparation for the grand conclusion.}
Schematically, our goal for the rest of this section will be to {\it add} 
(\ref{antrikos2}) to the equation (\ref{proolaxreiazontai1}) (or (\ref{proolaxreiazontai2}),
(\ref{proolaxreiazontai3}), depending on the form of the selected factor), thus 
deriving a new true equation which we denote by:

\begin{equation}
\label{grandgoal} Image^{1,+}_{\phi_{u+1}}[L_{g}]+Soph\{
Im^{1,\beta}[L_{g}]\} =0;
\end{equation}
this holds modulo complete contractions of length $\ge\sigma
+u+2$. {\it This} new true equation is the 
``grand conclusion'', which is the main aim of our present
paper. The ``grand conclusion'' will almost directly imply Lemma 
\ref{pskovb} in case A.
 It will also be the main tool in deriving Lemma \ref{pskovb} in
  case B--this will be done in section \ref{caseB}.

\par A few easy observations:

\begin{equation}
\label{crucifi}
\begin{split}
&\Sum_{j\in J^{\sigma +u+2}} a_j Soph[C^j_{g}(\Omega_1,\dots ,
\Omega_p,\phi_1,\dots ,\phi_u,\phi_{u+1},\omega)]=
\\&\Sum_{j\in J} a_j C^j_{g}(\Omega_1,\dots ,
\Omega_p,\phi_1,\dots ,\phi_u,\phi_{u+1}).
\end{split}
\end{equation}

\par We also observe that for $m\in M$, 
\begin{equation}\label{yegor} 
 Soph\{ Xdiv_{i_1}\dots Xdiv_{i_a} C^{m,i_1\dots i_{a+2}}_{g}
(\Omega_1,\dots ,\Omega_p, \phi_1,\dots
,\phi_u)\nabla_{i_{a+1}}\phi_{u+1}\nabla_{i_{a+2}}\omega\}
\end{equation}
 is an acceptable contributor (see Definition \ref{contributeur2}),
 while if $m\in M^\sharp$ (\ref{yegor} 
 is an un-acceptable contributor with one un-acceptable
 factor, and for $m\in M^{\sharp\sharp}$ (\ref{yegor}) 
 is a linear combination of terms 
with all the properties of contributors, {\it but} there will 
be {\it two} un-acceptable factors $\nabla\Omega,\nabla\Omega_{h'}$
 that are contracting against each other
 (and if $a=\mu$ then all free indices are non-special). We have denoted by
$$\Sum_{h\in H} a_h Xdiv_{i_1}\dots Xdiv_{i_a}
C^{l,i_1\dots i_a}_{g}(\Omega_1,\dots ,\Omega_p, \phi_1,
\dots,\phi_{u+1})$$ a generic linear combination of contributors
(acceptable or with one unacceptable factor $\nabla\Omega_h$ as
in the conclusion of Lemma \ref{pskovb}); we also denote by
$$\Sum_{h\in H^{\S\S}} a_h Xdiv_{i_1}\dots Xdiv_{i_a}
C^{l,i_1\dots i_a}_{g}(\Omega_1,\dots ,\Omega_p, \phi_1,
\dots,\phi_{u+1})$$ generic linear combinations of terms line the ones indexed in $M^{\sharp\sharp}$.
\newline

\par By definition, we observe that
\begin{equation}
\label{enakaiduo} Soph\{ \Sum_{\zeta\in Z} a_\zeta Xdiv_{i_1}\dots
Xdiv_{a_{\mu-1}} \Sum_{K=1}^{\sigma-1}Hit^K_\omega[C^{\zeta,i_1\dots
i_{\mu -1}}_{g}(\Omega_1,\dots ,\Omega_p, \phi_1,
\dots,\phi_u,\phi_{u+1})]\}
\end{equation}
is a contributor,\footnote{See Definition \ref{contributeur2}.}
 because acting by $Soph\{\dots\}$ on the operation
$\sum_{K=1}^{\sigma-1}Hit_\omega^K$ gives rise to
another $Xdiv$ (see the Definition \ref{matia} and the discussion under it).

\par Now, we proceed to derive some delicate cancellations occurring in (\ref{grandgoal}).

\par Observe that:

\begin{equation}
\label{xusia}
\begin{split}
&Soph\{ \Sum_{l\in L} a_l \Sum_{(i_k,i_l)\in (I\times I)^\sharp}
Xdiv_{i_1}\dots \hat{Xdiv}_{i_k}\dots\hat{Xdiv}_{i_l}\dots
Xdiv_{i_a}C^{l,i_1\dots ,i_a}_{g}
\nabla_{i_k}\phi_{u+1}\nabla_{i_l}\omega\}
\\&=\Sum_{i=1}^{b_l}\{ \Sum_{(i_k,i_l)\in
(I^{T_i}_2)^{2,dif}} Xdiv_{i_1}\dots \hat{Xdiv}_{i_k}\dots
\hat{Xdiv}_{i_l}\dots Xdiv_{i_a} [C^{l,i_1\dots i_a,i_{*}|T_i}_{g}
\\&(\Omega_1,\dots
,\Omega_p,\phi_1,\dots
,\phi_u)g^{i_ki_l}]\nabla_{i_{*}}\phi_{u+1}+ \Sum_{i_k\in I^{T_i}_1, i_l\in I^{T_i}_2}
Xdiv_{i_1}\dots \hat{Xdiv}_{i_k}\dots \hat{Xdiv}_{i_l}\\&\dots
Xdiv_{i_a} \nabla^{i_{*}}_{T_i}[C^{l,i_1\dots i_a}_{g}(\Omega_1,
\dots ,\Omega_p, \phi_1,\dots ,\phi_u)g^{i_ki_l}]
\nabla_{i_{*}}\phi_{u+1}\}.
\end{split}
\end{equation}

\par For our next observation, we will look at each
$l\in L$ and pick out each pair of indices
$({}_{i_k},{}_{i_l})\in I^d_l$. We assume with no loss of generality
 that ${}_{i_k}$ is a derivative index (recall that $I^d_l$ stands
  for the set of pairs of indices that belong to the same factor in
  $C^{l,i_1\dots i_a}_g$ and at least one of which is a derivative index). 
Then, for each such pair
$({}_{i_k},{}_{i_l})$,  we denote by

$$\dot{C}^{l,i_1\dots \hat{i}_k\dots i_\mu,i_{*}|T_i}_{g}
(\Omega_1,\dots,\Omega_p, \phi_1, \dots,\phi_u)\nabla_{i_{*}}
\phi_{u+1}$$
the tensor field that arises from $C^{l,i_1\dots i_\mu}$ by
erasing the index ${}_{i_k}$ and adding a
 derivative index $\nabla_{i_{*}}$ onto the selected factor $T_i$
and then contracting ${\nabla}_{i_{*}}$ against a factor
$\nabla^{i_{*}}\phi_{u+1}$. We see that for each
$l\in L$:

\begin{equation}
\label{youaretheone2} \begin{split} & Soph\{ \Sum_{l\in L} a_l
\Sum_{(i_k,i_l)\in I^d_l} \Sum_{K=1}^{\sigma
-1}Xdiv_{i_1}\dots\hat{Xdiv}_{i_k}\dots \hat{Xdiv}_{i_l}\dots
Xdiv_{i_a}\\&\{ Re^{K,\omega,\phi_{u+1}}_{i_k,i_l}[C^{l,i_1\dots
,i_a}_{g}] +Re^{K,\phi_{u+1},\omega}_{i_k,i_l}[C^{l,i_1\dots
,i_a}_{g}]\}\}= \Sum_{i=1}^{b_l}\Sum_{(i_k,i_l)\in I^d_l}
\\&Xdiv_{i_1}\dots \hat{Xdiv}_{i_k}\dots Xdiv_{i_a}
\dot{C}^{l,i_1\dots \hat{i}_k\dots i_\mu,i_{*}|T_i}_{g}
(\Omega_1,\dots,\Omega_p, \phi_1, \dots,\phi_u)\nabla_{i_{*}}
\phi_{u+1}.
\end{split}
\end{equation}

\par Furthermore, for each selected factor $T_i$,
let us denote by $I^{d,non-T_i}_l\subset I^d_l$ the subset of
$I^d_l$  that consists of pairs of free indices that do not belong
to the selected factor $T_i$.  We observe:

\begin{equation}
\label{youaretheone1} \begin{split} & \Sum_{i=1}^{b_l}
\Sum_{i_y\in I^{T_i}_2} \sigma(i_y)Xdiv_{i_1}\dots
\hat{Xdiv}_{i_y}\dots Xdiv_{i_\mu} [\nabla^{i_y}_{T_i}
C^{l,i_1\dots\hat{i}_y\dots i_\mu}_{g}(\Omega_1,\dots , \Omega_p,
\phi_1,\dots ,\phi_u)\\&\nabla_{i_y}\phi_{u+1}]
 =\Sum_{i=1}^{b_l}
Soph\{ \Sum_{i_h\in I^{*,2|T_i}_l\bigcup I^{+,2|T_i}_l} \Sum_{t\in
Set[T(i)]} Xdiv_{i_1}\dots \hat{Xdiv}_{i_h}\dots Xdiv_{i_a}\\&\{
Repla^{i,t}_{\phi_{u+1},\omega}[C^{l,i_1\dots ,i_a}_{g}]
+Repla^{t,i}_{\phi_{u+1},\omega} [C^{l,i_1\dots ,i_a}_{g}]\}+
 2\Sum_{i=1}^{b_l}\Sum_{(i_k,i_l)\in I^{d,non-T_i}_l}
\\&Xdiv_{i_1}\dots \hat{Xdiv}_{i_k}\dots Xdiv_{i_a}
\dot{C}^{l,i_1\dots \hat{i}_k\dots i_\mu,i_{*}|T_i}_{g}
(\Omega_1,\dots,\Omega_p, \phi_1, \dots,\phi_u)\nabla_{i_{*}}
\phi_{u+1}\}.
\end{split}
\end{equation}

\par Now, a few more delicate observations. For each
$l\in L_\mu$ and each selected factor $T_i$, we denote by
$I^{\sharp,T_i}_2\subset I^\sharp$\footnote{Recall that $I^\sharp$
stands for the set of free indices in $C^{l,i_1\dots i_a}_g$ that
belong to a factor $S_{*}\nabla^{(\nu)}R_{ijkl}$.} the
 index set of the indices that belong to a factor
$S_{*}\nabla^{(\nu)} R_{ijkl}\ne T_i$. We also denote by
$I^{\sharp,T_i}_1= I^\sharp\setminus I^{\sharp,T_i}_2$.
Furthermore, for each tensor field $C^{l,i_1\dots ,i_\mu}_{g}$ and
each free index ${}_{i_h}$ in that tensor field, we will set
$2_{i_h}=2$ if ${}_{i_h}$ belongs to a factor $\nabla^{(m)}R_{ijkl}$ or
$S_{*}\nabla^{(\nu)} R_{ijkl}$ and $2_{i_h}=0$ if it belongs to a
factor $\nabla^{(A)}\Omega_h$.
 Then, comparing
the discussion above (\ref{analyse664}) and (\ref{mumu}), we
 derive that:

\begin{equation}
\label{xilary}
\begin{split}
&\Sum_{i=1}^{b_l}\Sum_{i_y\in I^{T_i}_2} \tau(i_y)Xdiv_{i_1}\dots
\hat{Xdiv}_{i_y}\dots Xdiv_{i_\mu}[\nabla^{i_y}_{T_i}C^{l,i_1\dots
\hat{i_y}\dots i_\mu}_{g}(\Omega_1,\dots
,\Omega_p)\nabla_{i_y}\phi_{u+1}]
\\&+Soph\{ \Sum_{i_h\in I^{\sharp,T_i}_2}Xdiv_{i_1}\dots \hat{Xdiv}_{i_h}
\dots Xdiv_{i_a}\\&\{ FRep^{i_h,1,\omega,\phi_{u+1}}[C^{l,i_1\dots
i_a|f(i_h)}_{g} (\Omega_1,\dots ,\Omega_p,\phi_1,\dots ,\phi_u)]
\\& +FRep^{i_h,1,\phi_{u+1},\omega}[C^{l,i_1\dots i_a|
f(i_h)}_{g}(\Omega_1,\dots ,\Omega_p,\phi_1,\dots ,\phi_u)]
\\& +FRep^{i_h,2,\omega,\phi_{u+1}}[C^{l,i_1\dots i_a
|f(i_h)}_{g}(\Omega_1,\dots ,\Omega_p,\phi_1,\dots ,\phi_u)] 
\\&+ FRep^{i_h,2,\phi_{u+1},\omega}[C^{l,i_1\dots i_a
|f(i_h)}_{g}(\Omega_1,\dots ,\Omega_p,\phi_1,\dots , \phi_u)]\}\}=
\\& \Sum_{i=1}^{b_l}\Sum_{i_y\in I^{T_i}_2}
2_{i_y}Xdiv_{i_1}\dots \hat{Xdiv}_{i_y}\dots
Xdiv_{i_\mu}[\nabla^{i_y}_{T_i} C^{l,i_1\dots \hat{i_y}\dots
i_\mu}_{g}(\Omega_1,\dots
,\Omega_p)\nabla_{i_y}\phi_{u+1}].
\end{split}
\end{equation}

\par On the other hand, we observe that:

\begin{equation}
\label{kokoko}
\begin{split}
&Soph\{\Sum_{m\in M_{\mu-1}} a_m Xdiv_{i_1}\dots Xdiv_{i_{\mu-1}}
C^{m,i_1\dots ,i_{\mu-1}}_{g}(\Omega_1,\dots,\Omega_p,
\phi_1,\dots ,\phi_u,\phi_{u+1},\omega) \}
\\&=\Sum_{b\in B'} a_b Xdiv_{i_1}\dots Xdiv_{i_{\mu-1}}
C^{b,i_1\dots ,i_{\mu-1}}_{g}(\Omega_1,\dots ,\Omega_p,
\phi_1,\dots ,\phi_{u+1}).
\end{split}
\end{equation}
(Recall that the sublinear combination indexed in $M_{\mu-1}$ 
appears only in the special subcase of case B).
The linear combination indexed in $B'$ is a generic linear combination 
defined in Definition \ref{generalni}.

\par Next, we note some further cancellations,
 for each $C^{l,i_1\dots ,i_\mu}_{g}$ with
at least one free index in the selected factor. In the case where
the selected factor is of the form $\nabla^{(A)}\Omega_h$ (in
which case it is unique) we must have:

\begin{equation}
\label{xilary2a} \begin{split} & Soph\{ \Sum_{i_h\in
I^{*,1}_l}\Sum_{t\in Set(T(i_h))} Xdiv_{i_1}\dots
\hat{Xdiv}_{i_h}\dots Xdiv_{i_a}\\&\{
Repla^{i_h,t}_{\phi_{u+1},\omega}[C^{l,i_1\dots ,i_a}_{g}]
+Repla^{t,i_h}_{\phi_{u+1},\omega}[C^{l,i_1\dots ,i_a}_{g}]\}\}
\\&=\Sum_{i_h\in I_1} A^\sharp Xdiv_{i_1}\dots
\hat{Xdiv}_{i_h}\dots Xdiv_{i_\mu}C^{l,i_1\dots
i_\mu}_{g}(\Omega_1,\dots ,\Omega_p,\phi_1,\dots , \phi_u)
\nabla_{i_h}\phi_{u+1};
\end{split}
\end{equation}
(Recall that $A^\sharp$ stands for the number of indices in the factor 
$\nabla^{(A)}\Omega_h$ that are not free and not contracting 
against a factor $\nabla\phi_h$). On the other hand, if the selected factor(s) is
(are) of the form $\nabla^{(m)}R_{ijkl}$ we will have:

\begin{equation}
\label{xilary2b} \begin{split} & Soph\{ \Sum_{i_h\in
I^{*,1}_l}\Sum_{t\in Set(T(i_h))} Xdiv_{i_1}\dots
\hat{Xdiv}_{i_h}\dots Xdiv_{i_a}\\&\{
Repla^{i_h,t}_{\phi_{u+1},\omega}[C^{l,i_1\dots ,i_a}_{g}]
+Repla^{t,i_h}_{\phi_{u+1},\omega}[C^{l,i_1\dots ,i_a}_{g}]\}\}
=\Sum_{i=1}^{b_l}\Sum_{i_h\in I^{T_i}_1}
(m_{i}^\sharp+4) \\&Xdiv_{i_1}\dots \hat{Xdiv}_{i_h}\dots
Xdiv_{i_\mu}C^{l,i_1\dots i_\mu}_{g}(\Omega_1,\dots ,
\Omega_p,\phi_1,\dots , \phi_u)\nabla_{i_h}\phi_{u+1};
\end{split}
\end{equation}
(recall that $m^\sharp_{i}$ stands for the number of
 derivative indices in the selected factor
$T_i=\nabla^{(m)}R_{ijkl}$ that are not free and not contracting
 against a factor $\nabla\phi_h$).

Finally, in the case where the selected factor is of the form
$S_{*}\nabla^{(\nu)} R_{ijkl}$ (in which case it is again
 unique),  for each $l\in L_\mu$
with at least one free index in the selected factor, we find:

\begin{equation}
\label{xilary2c}
\begin{split}
&Soph\{ \Sum_{i_h\in I^{*,1}_l\bigcup I^{+,1}_1}\Sum_{t\in
Set(T(i_h))} Xdiv_{i_1}\dots \hat{Xdiv}_{i_h}\dots Xdiv_{i_a} \\&\{
Repla^{i_h,t}_{\phi_{u+1},\omega} [C^{l,i_1\dots ,i_a}_{g}]+
Repla^{t,i_h}_{\phi_{u+1},\omega} [C^{l,i_1\dots
,i_a}_{g}]\}\}+
 Soph\{\Sum_{i_f\in I^\sharp_1}Xdiv_{i_1}\dots
\\&\hat{Xdiv}_{i_h} \dots Xdiv_{i_a}\{
FRep^{i_h,1,\omega,\phi_{u+1}}[C^{l,i_1\dots i_a|f(i_h)}_{g}
(\Omega_1,\dots ,\Omega_p,\phi_1,\dots ,\phi_u)]
\\& +FRep^{i_h,1,\phi_{u+1},\omega}[C^{l,i_1\dots i_a|
f(i_h)}_{g}(\Omega_1,\dots ,\Omega_p,\phi_1,\dots ,\phi_u)]
\\& +FRep^{i_h,2,\omega,\phi_{u+1}}[C^{l,i_1\dots i_a
|f(i_h)}_{g}(\Omega_1,\dots ,\Omega_p,\phi_1,\dots ,\phi_u)] +\\&
FRep^{i_h,2,\phi_{u+1},\omega}[C^{l,i_1\dots i_a
|f(i_h)}_{g}(\Omega_1,\dots ,\Omega_p,\phi_1,\dots , \phi_u)]\}\}=
\\& \Sum_{i_h\in I_1} (\nu^\sharp +2) Xdiv_{i_1}\dots
\hat{Xdiv}_{i_h}\dots Xdiv_{i_\mu} C^{l,i_1\dots
i_\mu}_{g}(\Omega_1,\dots ,\Omega_p,\phi_1,\dots ,
\phi_u)\nabla_{i_h}\phi_{u+1};
\end{split}
\end{equation}
(recall that $\nu^\sharp$ stands for the number of indices
 ${}_{r_1},\dots ,{}_{r_\nu},{}_j$ in the selected factor
$S_{*}\nabla^{(\nu)}_{r_1\dots r_\nu}R_{ijkl}$  that are not free
 and not contracting against a factor $\nabla\phi_h$).

\section{The grand conclusion, and the proof of Lemma \ref{pskovb}.}

\subsection{The grand conclusion.}

\par Now, we combine all the cancellations we have noted in the previous subsection
 to derive the ``grand conclusion''.
When the selected factor(s) is (are) of the form $S_{*}\nabla^{(\nu)}R_{ijkl}$
 or $\nabla^{(m)}R_{ijkl}$, the grand conclusion will be the equation:

\begin{equation}
\label{rampling}
\begin{split}
&Image^{1,+}_{\phi_{u+1}}[L_{g}]+Soph\{Im^{1,\beta}_{\phi_{u+1}}[L_{g}]\}
+\{L_g(\Omega_1\cdot\phi_{u+1},\dots,\Omega_p,\phi_1,\dots,\phi_u)\\&+\dots+
L_g(\Omega_1,\dots,\Omega_X\cdot\phi_{u+1},\dots,\Omega_p,\phi_1,\dots,\phi_u)\}=0.
\end{split}
\end{equation}
(Recall that $\Omega_1,\dots,\Omega_X$
are the factors in $\vec{\kappa}_{simp}$ that are not contracting
against any factor $\nabla\phi_h$. The terms in $\{\dots\}$ appear
{\it only} when the selected factor(s) is (are) curvature terms).

When the selected factor is of the form $\nabla^{(B)}\Omega_x$, the grand conclusion will
 be the equation:

\begin{equation}
\label{rampling'}
\begin{split}
&Image^{1,+}_{\phi_{u+1}}[L_{g}]+Soph\{Im^{1,\beta}_{\phi_{u+1}}[L_{g}]\}=0.
\end{split}
\end{equation}

\par For future reference, we put down a few facts before we
write out the ``grand conclusion'':

{\it Recall notation}: Recall that
$s$ stands for the total number of factors
$\nabla^{(m)}R_{ijkl}$, $S_{*}\nabla^{(\nu)} R_{ijkl}$ in the
simple character $\vec{\kappa}_{simp}$ (all the tensor fields in (\ref{hypothese2})
have this given simple character--see the introduction 
of the present paper for a simplified discussion of this notion).
 Recall that for each $C^{l,i_1\dots i_\mu}_g$, $l\in L_\mu$:
$\gamma$ (or $\gamma_i$ if there are multiple selected factors $T_i$)
stands for the number of indices in $C^{l,i_1\dots i_\mu}_g$
that do not belong to the selected factor and are not
contracting against a factor $\nabla\phi_h$. We also recall that $I_1$ (or $I_1^{T_i}$
 if there are multiple selected factors) stands for the set of
 free indices that belong to the selected factor, and
$I_2$ (or $I_2^{T_i}$
 if there are multiple selected factors) stands for the set of
 free indices that {\it do not} belong to the selected factor. We also recall that for each
$l\in L_{\mu}$ and each free index
${}_{i_h}\in I_2$ (or ${}_{i_h}\in I^{T_i}_2$) which belongs to $C^{l,i_1\dots i_\mu}_g$,
$2_{i_h}$ stands for the number $2$ if the free index ${}_{i_h}$ belongs to
a factor of the from $\nabla^{(m)}R_{ijkl}$ or $S_{*}\nabla^{(\nu)}R_{ijkl}$, and it will be zero if it
belongs to a factor of the form $\nabla^{(B)}\Omega_x$. Now, we define
$\overline{2}_{i_h}$ to equal number $2$ if the free index ${}_{i_h}$ belongs to
a factor of the from $\nabla^{(m)}R_{ijkl}$ or
 $S_{*}\nabla^{(\nu)}R_{ijkl}$, and to equal 1 if it
belongs to a factor of the form $\nabla^{(B)}\Omega_x$.

Finally, we recall: When the selected factor is of the form $S_{*}\nabla^{(\nu)}R_{ijkl}$
then  (for each $\mu$-tensor field $C^{l,i_1\dots i_\mu}_g$, $l\in L_\mu$)
$\nu^\sharp$ stands for the number of indices in the selected factor $S_{*}\nabla^{(\nu)}R_{ijkl}$
that are not free and not contacting against a factor $\nabla\phi_h$. When the
 selected factor(s) is (are) of the form $\nabla^{(m)}R_{ijkl}$, then
(for each $\mu$-tensor field $C^{l,i_1\dots i_\mu}_g$, $l\in L_\mu$ and)
for each selected factor $\nabla^{(m_i)}R_{ijkl}$,
$m_i^\sharp$ stands for the number of derivative indices that are
 not free and not contracting against a factor $\nabla\phi_h$.  Lastly, when
the selected factor is of the form $\nabla^{(A)}\Omega_k$ then
(for each $\mu$-tensor field $C^{l,i_1\dots i_\mu}_g$, $l\in
L_\mu$) 
$A^\sharp$ stands for the number of indices in
$\nabla^{(A)}\Omega_k$ that are not free and not contracting
against a factor $\nabla\phi_h$.
\newline

Then, if the selected factor is of the form
$S_{*}\nabla^{(\nu)} R_{ijkl}$, the ``grand conclusion'' may be written in the form:

\begin{equation}
\label{olaxreiazontai1}
\begin{split}
&\Sum_{l\in L_\mu} a_l \{-\Sum_{i_h\in I_1}(\gamma +(|I_1|-1)
+\nu^\sharp -2(s-1)-X)  Xdiv_{i_1}\dots\hat{Xdiv}_{i_h}
Xdiv_{i_\mu}\\& C^{l,i_1\dots i_\mu}_{g}(\Omega_1, \dots
,\Omega_p, \phi_1,\dots ,\phi_u)\nabla_{i_h}\phi_{u+1}+\Sum_{i_h\in I_2}\overline{2}_{i_h}
\\&  Xdiv_{i_1}\dots \hat{Xdiv}_{i_h}\dots
Xdiv_{i_\mu}\nabla^{i_{*}}_{sel}[C^{l,i_1\dots \hat{i}_h\dots
i_\mu,i_{*}}_{g}(\Omega_1,\dots ,\Omega_p,\phi_1,\dots ,
\phi_u)]\nabla_{i_{*}}\phi_{u+1}
\\&-  \{\Sum_{(i_k,i_l)\in
I^{d,non-sel}_l} Xdiv_{i_1}\dots \hat{Xdiv}_{i_k}\dots
Xdiv_{i_\mu} \dot{C}^{l,i_1\dots \hat{i}_k\dots i_\mu,i_{*}}_{g}
(\Omega_1,\dots,\Omega_p, \phi_1, \dots,\phi_u)\\&\nabla_{i_{*}}
\phi_{u+1}
+\Sum_{(i_k,i_l)\in I^{d,sel}_l}
\Sum_{S=1}^{\sigma -1}Xdiv_{i_1}\dots \hat{Xdiv}_{i_k}
\dots\hat{Xdiv}_{i_l}\dots Xdiv_{i_\mu}Xdiv_{i_z}
\\&\tilde{C}^{l,i_1\dots,\hat{i}_k\dots ,i_\mu,i_z|S}_{g}
(\Omega_1,\dots , \Omega_p, \phi_1,\dots
,\phi_u)\nabla_{i_l}\phi_{u+1}\}+
\\&\Sum_{j\in J^{\sigma +u+2}} a_j C^j_{g}(\Omega_1,\dots ,\Omega_p,
 \phi_1,\dots ,\phi_u,\phi_{u+1})+
\\&(\Sum_{b\in B'} a_b Xdiv_{i_1}\dots Xdiv_{i_{\mu-1}}
C^{b,i_1\dots ,i_{\mu-1}}_{g}(\Omega_1,\dots ,\Omega_p,
\phi_1,\dots ,\phi_{u+1}))+
\\&\Sum_{h\in H\bigcup H^{\S\S}} a_h Xdiv_{i_1}\dots Xdiv_{i_a}
C^{h,i_1\dots i_a,i_{*}}_{g}(\Omega_1,\dots ,\Omega_p,
 \phi_1,\dots ,\phi_u)\nabla_{i_{*}}\phi_{u+1}=0,
\end{split}
\end{equation}
modulo complete contractions of length $\ge\sigma +u+2$. 
 $\sum_{b\in B'}\dots$ (defined in Definition \ref{generalni}) 
appears {\it only} in the special subcase of case B.

\par If the selected factor(s) is (are) of the from
$\nabla^{(m)}R_{ijkl}$, the ``grand conclusion'' we obtain is very analogous:

\begin{equation}
\label{olaxreiazontai2}
\begin{split}
&\Sum_{l\in L_\mu} a_l\Sum_{i=1}^{b_l}\{ -\Sum_{i_h\in
I^{T_i}_1}(\gamma_i +(|I^{T_i}_1|-1) +m^\sharp_i -2(s-1)-X)
\\&Xdiv_{i_1}\dots \hat{Xdiv}_{i_h}\dots Xdiv_{i_\mu}C^{l,i_1\dots
i_\mu}_{g}(\Omega_1, \dots ,\Omega_p, \phi_1,\dots
,\phi_u)\nabla_{i_h}\phi_{u+1}+
\\& \Sum_{i_h\in I^{T_i}_2}\overline{2}_{i_h}\nabla^{i_{*}}_{T_i}[C^{l,i_1\dots \hat{i}_h\dots
i_\mu,i_{*}}_{g}(\Omega_1,\dots ,\Omega_p,\phi_1,\dots ,
\phi_u)]\nabla_{i_{*}}\phi_{u+1} -\Sum_{(i_k,i_l)\in I^{d,non-T_i}_l}
\\&Xdiv_{i_1}\dots \hat{Xdiv}_{i_k}\dots Xdiv_{i_a}
\dot{C}^{l,i_1\dots \hat{i}_k\dots i_\mu,i_{*}}_{g}
(\Omega_1,\dots,\Omega_p, \phi_1, \dots,\phi_u)\nabla_{i_{*}}
\phi_{u+1}
\\&+\Sum_{(i_k,i_l)\in I^{d,T_i}_l}
\Sum_{S=1}^{\sigma -1}Xdiv_{i_1}\dots \hat{Xdiv}_{i_k}
\dots\hat{Xdiv}_{i_l}\dots Xdiv_{i_\mu}Xdiv_{i_z}
\\&\tilde{C}^{l,i_1\dots,\hat{i}_k\dots ,i_\mu,i_z|S}_{g}
(\Omega_1,\dots , \Omega_p, \phi_1,\dots
,\phi_u)\nabla_{i_l}\phi_{u+1}\}
\\&+\Sum_{j\in J^{\sigma +u+2}} a_j C^j_{g}(\Omega_1,\dots ,\Omega_p,
 \phi_1,\dots ,\phi_u,\phi_{u+1})+
 \\&(\Sum_{b\in B'} a_b Xdiv_{i_1}\dots Xdiv_{i_{\mu-1}}
C^{b,i_1\dots ,i_{\mu-1}}_{g}(\Omega_1,\dots ,\Omega_p,
\phi_1,\dots ,\phi_{u+1}))+
\\&\Sum_{h\in H\bigcup H^{\S\S}} a_h Xdiv_{i_1}\dots Xdiv_{i_a}
C^{h,i_1\dots i_a,i_{*}}_{g}(\Omega_1,\dots ,\Omega_p,
 \phi_1,\dots ,\phi_u)\nabla_{i_{*}}\phi_{u+1}=0,
\end{split}
\end{equation}
modulo complete contractions of length $\ge\sigma
+u+2$.  $\sum_{b\in B'}\dots$ (defined in Definition \ref{generalni}) 
appears {\it only} in the special subcase of case B.

Finally, in the case where the selected factor is of the form
$\nabla^{(A)}\Omega_1$, (\ref{rampling'}) gives us the ``grand conclusion'':

\begin{equation}
\label{olaxreiazontai3}
\begin{split}&\Sum_{l\in L_\mu} a_l\{ -\Sum_{i_h\in I_1}(\gamma +(|I_1|-1)
+A^\sharp -2s)  Xdiv_{i_1}\dots \hat{Xdiv}_{i_h}\dots
Xdiv_{i_\mu}C^{l,i_1\dots i_\mu}_{g}
\\&(\Omega_1, \dots ,\Omega_p,
\phi_1,\dots ,\phi_u)\nabla_{i_h}\phi_{u+1}+
\\& \Sum_{i_h\in I_2}2_{i_h}\nabla^{i_{*}}_{sel}[C^{l,i_1\dots \hat{i}_h\dots
i_\mu,i_{*}}_{g}(\Omega_1,\dots ,\Omega_p,\phi_1,\dots ,
\phi_u)]\nabla_{i_{*}}\phi_{u+1}-\Sum_{(i_k,i_l)\in I^{d,non-sel}_l} \\&\{ Xdiv_{i_1}\dots
\hat{Xdiv}_{i_k}\dots Xdiv_{i_\mu} \dot{C}^{l,i_1\dots
\hat{i}_k\dots i_\mu,i_{*}}_{g} (\Omega_1,\dots,\Omega_p, \phi_1,
\dots,\phi_u)\nabla_{i_{*}} \phi_{u+1}
\\&+\Sum_{(i_k,i_l)\in I^{d,sel}_l}\Sum_{S=1}^{\sigma -1}Xdiv_{i_1}\dots \hat{Xdiv}_{i_k}
\dots\hat{Xdiv}_{i_l}\dots Xdiv_{i_\mu}Xdiv_{i_z}
\\&\tilde{C}^{l,i_1\dots,\hat{i}_k\dots ,i_\mu,i_z|S}_{g}
(\Omega_1,\dots , \Omega_p, \phi_1,\dots
,\phi_u)\nabla_{i_l}\phi_{u+1}\}+
\\&\Sum_{j\in J^{\sigma +u+2}} a_j C^j_{g}(\Omega_1,\dots ,\Omega_p,
 \phi_1,\dots ,\phi_u,\phi_{u+1})+
\\&(\Sum_{b\in B'} a_b Xdiv_{i_1}\dots Xdiv_{i_{\mu-1}}
C^{b,i_1\dots ,i_{\mu-1}}_{g}(\Omega_1,\dots ,\Omega_p,
\phi_1,\dots ,\phi_{u+1}))+
\\&\Sum_{h\in H\bigcup H^{\S\S}} a_h Xdiv_{i_1}\dots Xdiv_{i_a}
C^{h,i_1\dots i_a,i_{*}}_{g}(\Omega_1,\dots ,\Omega_p,
 \phi_1,\dots ,\phi_u)\nabla_{i_{*}}\phi_{u+1}=0,
\end{split}
\end{equation}
modulo complete contractions of length
$\ge\sigma+u+2$.  $\sum_{b\in B'}\dots$ (defined in Definition \ref{generalni}) 
appears {\it only} in the special subcase of case B.

\par We observe that because of our Lemma assumption that
$L^{*}_\mu=\emptyset$,\footnote{See the notation in the {\it statement} of Lemma \ref{pskovb}.}
 it follows that all the $(\mu-1)$-tensor
fields above are acceptable. Moreover, by construction they each have a $(u+1)$-simple character
$\vec{\kappa}_{simp}^{+}$.

\par Now, we will show in a ``Mini-Appendix'' below
 that using the above, we may write:

\begin{equation}
\label{xiaodongcao}
\begin{split}
&\Sum_{h\in H^{\S\S}} a_h Xdiv_{i_1}\dots Xdiv_{i_a} C^{h,i_1\dots
i_a,i_{*}}_{g}(\Omega_1,\dots ,\Omega_p,
 \phi_1,\dots ,\phi_u)\nabla_{i_{*}}\phi_{u+1}=
\\&\Sum_{h\in H} a_h Xdiv_{i_1}\dots Xdiv_{i_a}
C^{h,i_1\dots i_a,i_{*}}_{g}(\Omega_1,\dots ,\Omega_p,
 \phi_1,\dots ,\phi_u)\nabla_{i_{*}}\phi_{u+1}
\\&+\Sum_{j\in J^{\sigma +u+2}} a_j C^j_{g}(\Omega_1,\dots ,\Omega_p,
 \phi_1,\dots ,\phi_u,\phi_{u+1}),
\end{split}
\end{equation}
{\it using generic notation in the right hand side--the sublinear
 combination in the left hand side is exactly
the one appearing in (\ref{olaxreiazontai1}),
(\ref{olaxreiazontai2}), (\ref{olaxreiazontai3})}.

\par In view of (\ref{xiaodongcao}) (which we prove below in the appendix), we may
assume that $H^{\S\S}=\emptyset$ in (\ref{olaxreiazontai1}),
(\ref{olaxreiazontai2}), (\ref{olaxreiazontai3}),
whenever we refer to these equations.
\newline

\subsection{Proof of Lemma \ref{pskovb} in case A.} We pick the selected factor(s) to
be the second critical factor(s) (see the statement of Lemma \ref{pskovb}). Recall that in this case A
 the second critical factor has at least two
 free indices.

 For convenience, in each sublinear combination
$$\big{(}\sum_{(i_k,i_l)\in I^{d,T_i}_l}\big{)}\Sum_{S=1}^{\sigma -1}
\tilde{C}^{l,i_1\dots,\hat{i}_k\dots ,i_\mu,i_z|S}_{g}
(\Omega_1,\dots , \Omega_p, \phi_1,\dots
,\phi_u)\nabla_{i_l}\phi_{u+1},$$ 
we will assume that among all the factors
$F_1,\dots ,F_{\sigma -1}$, the first critical factor(s) is
 (are) $F_1$ (or $F_1,\dots ,F_a$).\footnote{The expression
$\sum_{(i_k,i_l)\in I^{d,T_i}_l}$ will only
 be present when the selected (second critical) factor
 is generic in the form $\nabla^{(m)}R_{ijkl}$.}

\par We then claim that among all the $(\mu-1)$-tensor
fields in (\ref{olaxreiazontai1}), (\ref{olaxreiazontai2}),
(\ref{olaxreiazontai3})
(all of which have a $(u+1)$-simple character
$\vec{\kappa}^{+}_{simp})$, the sublinear combination of {\it
maximal} refined double character will be precisely:

\begin{equation}
\label{mona3ia} \Sum_{z\in Z'_{Max}}\Sum_{l\in L^z}
a_l\big{(}\sum_{(i_k,i_l)\in I^{d,T_i}_l}\big{)}
\tilde{C}^{l,i_1\dots\hat{i}_k\dots ,i_\mu,i_z|1}_{g}
(\Omega_1,\dots ,\Omega_p,\phi_1,\dots ,\phi_u)
\nabla_{i_l}\phi_{u+1},
\end{equation}
in the case where there is only one critical factor, and:

\begin{equation}
\label{mona3ia2} \Sum_{z\in Z'_{Max}} \Sum_{l\in L^z} a_l
\big{(}\sum_{(i_k,i_l)\in I^{d,T_i}_l}\big{)}\Sum_{S=1}^a
\tilde{C}^{l,i_1\dots\hat{i}_k\dots ,i_\mu,i_z|S}_{g}
(\Omega_1,\dots , \Omega_p, \phi_1,\dots
,\phi_u)\nabla_{i_l}\phi_{u+1},
\end{equation}
in the case where there are $a>1$ critical factors.

\par This fact essentially follows just by our definitions:
Firstly observe that the tensor fields in the above two equations
 have $M+1$ free indices in some factor.
Now, by definition of the {\it maximal} refined double character
we observe that for each $l\in L_\mu$, each factor in $C^{l,i_1\dots
i_\mu}_{g}$ can have at most $M$ free
 indices in one of its factors. Hence,  each
 tensor field of rank $\mu-1$ in the above three
equations {\it other than the tensor fields with a tilde sign,
$\tilde{C}$}, will again have at most $M$ 
free indices in any one of its factors (and this is double subsequent to the terms in
 (\ref{mona3ia}), (\ref{mona3ia2})).

\par Moreover, for each $l\in L_\mu\setminus
\bigcup_{z\in Z'_{Max}} L^z$, we observe
by definition that each tensor field in

\begin{equation}
\label{mona3ia3}
 \Sum_{S=1}^{\sigma-1}
\tilde{C}^{l,i_1\dots,\hat{i}_k\dots ,i_\mu,i_z|S}_{g}
(\Omega_1,\dots , \Omega_p, \phi_1,\dots
,\phi_u)\nabla_{i_l}\phi_{u+1}
\end{equation}
will either have at most $M$ free indices in any given
 factor or will have $M+1$ free indices in one
 factor but then its refined double character will be subsequent to
$\vec{L^z}, z\in Z'_{Max}$: This second claim just follows by the
 construction of the tensor fields above:
If $l\in L_\mu\setminus \bigcup_{z\in Z'_{Max}} L^z$ then the refined double character of
$C^{l,i_1\dots i_\mu}_g$ will be either doubly subsequent or ``equipolent'' 
to each refined double character $\vec{L}^z$, $z\in Z'_{Max}$ (which corresponds to the tensor fields
$C^{l,i_1\dots i_\mu}_g, l\in L^z, z\in Z'_{Max}$).\footnote{See the introduction for  
a discussion of these notions.} Now
$\tilde{C}^{l,i_1\dots,\hat{i}_k\dots ,i_\mu,i_z|S}_{g}
(\Omega_1,\dots , \Omega_p,\phi_1,\dots ,\phi_u)$ formally arises
 from $C^{l,i_1\dots i_\mu}$ by erasing the free
index ${}_{i_k}$ from the (selected) factor $T_i$ and adding a new
free index $\nabla_{i_z}$ onto another factor, with $M$ free
indices. Thus, our claim just follows by the definition of ordering among
 refined double characters.\footnote{See \cite{alexakis4} for a strict
  definition of this notion--see also the introduction of 
  the present paper for a simplified description of this notion.}

\par So we observe that the ``grand conclusion'' 
 proves Lemma \ref{pskovb} in the case A: The sublinear combination
(\ref{mona3ia}) in the grand conclusion is precisely the first
 line in (\ref{esmen}). All the other $(\mu-1)$-tensor
fields in the grand conclusion are in the general
 form $\sum_{\nu\in N} a_\nu\dots$ described in the claim of Lemma \ref{pskovb}.
Also, the tensor fields indexed in $H$ (with rank $\ge\mu$) are in the
 same general form as the tensor fields indexed in
 $T_1\bigcup T_2\bigcup T_3\bigcup T_4$ in (\ref{esmen}). $\Box$
\newline

{\bf Notes Regarding case B:} We will prove Lemma \ref{pskovb} in case B in section 
\ref{caseB} (and our proof there will heavily 
rely on the grand conclusion above). We only end this 
section with a remark, which will be essential in the proof of Lemma \ref{pskovb}
in case B:

{\it Important Remark:} The quantities in parentheses in the first lines of
(\ref{olaxreiazontai1}), (\ref{olaxreiazontai2}),
(\ref{olaxreiazontai3}) are {\it universal}, i.e.~they only
 depend on the simple character $\vec{\kappa}_{simp}$, and on the {\it form}
 of the selected factor $T_d$ (meaning whether it is of the form $S_{*}\nabla^{(\nu)}R_{ijkl}$,
 $\nabla^{(m)}R_{ijkl}$ or $\nabla^{(A)}\Omega_h$):
We denote those quantities (inside the parentheses)
by $q_z$. We observe that in the case of
(\ref{olaxreiazontai3}):

\begin{equation}
\label{ellhn1}
q_d=n-2u-\mu-1.
\end{equation}
In the case of (\ref{olaxreiazontai1}):

\begin{equation}
\label{ellhn2}
q_d=n-2u-\mu-1-X.
\end{equation}
Whereas in the case of (\ref{olaxreiazontai2}):

\begin{equation}
\label{ellhn3}
q_d=n-2u-\mu-3-X.
\end{equation}
(We will define $Q_d=|I_1|\cdot q_d$, for future reference).

\subsection{Mini-Appendix: Proof of (\ref{xiaodongcao}).}
\label{secxiao}

\par To prove this claim we will need to distinguish two cases: Either $\sigma=4$ or $\sigma>4$.
We will start with the case $\sigma>4$ which is the easiest.
\newline

{\it Proof of (\ref{xiaodongcao}) when $\sigma>4$:} In this setting,
we refer back to the grand conclusion. For each tensor field
$C^{h,i_1\dots i_a,i_{*}}_g\nabla_{i_{*}}\phi_{u+1}$, $h\in H^{\S\S}$ we
define
 $$X^\sharp div_{i_1}\dots X^\sharp div_{i_a}C^{h,i_1\dots i_a,i_{*}}_g\nabla_{i_{*}}\phi_{u+1}$$
to stand for the sublinear combination in
 $Xdiv_{i_1}\dots Xdiv_{i_a}C^{h,i_1\dots i_a,i_{*}}_g\nabla_{i_{*}}\phi_{u+1}$
where each $\nabla_{i_v}$ is not allowed to hit {\it either of the two}  factors
$\nabla\Omega_x,\nabla\Omega_{x'}$ (which are contracting against each other).

\par We may then straightforwardly use the fact that the grand conclusion
 holds formally to derive an equation:

\begin{equation}
\label{rixardos}
\begin{split}
&\sum_{h\in H\bigcup H^{\S\S}} a_h X^\sharp div_{i_1}\dots
X^\sharp div_{i_a} C^{h,i_1\dots i_a,i_{*}}_{g}(\Omega_1,\dots
,\Omega_p,
 \phi_1,\dots ,\phi_u)\nabla_{i_{*}}\phi_{u+1}+
\\&\sum_{j\in J} a_j C^{j}_{g}(\Omega_1,\dots ,\Omega_p,
 \phi_1,\dots ,\phi_{u+1})=0,
\end{split}
\end{equation}
where $\sum_{j\in J}\dots$ above stands for a generic
 linear combination of complete contractions of length
$\sigma+u+1$ with a weak character $Weak(\vec{\kappa}_{simp}^{+})$,
 with two factors $\nabla\Omega_x,\nabla\Omega_{x'}$ 
contracting against each other and which are 
{\it simply subsequent} to $\vec{\kappa}_{simp}$.

\par Now, we state a Lemma (which will be applied to other settings in the future), which
fits perfectly with the equation above:

\begin{lemma}
\label{yarusal}
Consider a linear combination of tensor fields,
\\$\sum_{\tau\in T} a_\tau C^{\tau,i_1\dots i_a}_g(\Omega_1,\dots,
\Omega_{p'},\phi_1,\dots ,\phi_{u'})$, each with a given simple character
 $\overline{\kappa}_{simp}$, and each with $a\ge V$ (for some given $V$).
 We assume that this simple character falls
 under the inductive assumption of Proposition \ref{giade}.

We consider the tensor fields $C^{\tau,i_1\dots i_a}_g(\Omega_1,\dots,
\Omega_{p'},\phi_1,\dots ,\phi_{u'})\nabla_i\Omega_x\nabla^i\phi_q$
which arise from the above by just multiplying by $\nabla_i\Omega_x\nabla^i\phi_q$.
We assume an equation:

\begin{equation}
\label{hmarton}
\begin{split}
&\sum_{\tau\in T} a_\tau X^\sharp div_{i_1}\dots X^\sharp div_{i_a}
[C^{\tau,i_1\dots i_a}_g(\Omega_1,\dots,\Omega_{p'},\phi_1,\dots ,\phi_{u'})
\nabla_i\Omega_x\nabla^i\phi_q]+
\\&\sum_{j\in J} a_j C^{j}_{g}(\Omega_1,\dots ,\Omega_{p'},
 \phi_1,\dots ,\phi_{u'})\nabla_i\Omega_x\nabla^i\phi_q=0,
\end{split}
\end{equation}
where $X^\sharp div_i$ stands for the sublinear combination in $Xdiv_i$ where $\nabla_i$ is in
addition not allowed to hit the expression $\nabla_i\Omega_x\nabla^i\phi_q$. $\sum_{j\in J}
a_j C^{j}_{g}(\Omega_1,\dots ,\Omega_{p'}, \phi_1,\dots ,\phi_{u'})$ stands for
 a generic linear combination of complete contractions with a weak character
$Weak(\overline{\kappa}_{simp})$ and simply subsequent to
$\overline{\kappa}_{simp}$. Furthermore, any terms of rank $\mu$ must
 have all $\mu$ free indices being non-special. 
\newline

\par Our conclusion is then that we can write:

\begin{equation}
\label{hmartoncor}
\begin{split}
&\sum_{\tau\in T} a_\tau Xdiv_{i_1}\dots Xdiv_{i_a}
[C^{\tau,i_1\dots i_a}_g(\Omega_1,\dots,\Omega_{p'},\phi_1,\dots ,\phi_{u'})
\nabla_i\Omega_x\nabla^i\phi_q]=
\\&\sum_{\tau\in T'} a_\tau Xdiv_{i_1}\dots Xdiv_{i_a}
C^{\tau,i_1\dots i_a}_g(\Omega_1,\dots,\Omega_{p'},\phi_1,\dots ,\phi_{u'}|\Omega_x,\phi_q)
\\&+\sum_{j\in J} a_j C^{j}_{g}(\Omega_1,\dots ,\Omega_{p'},
 \phi_1,\dots ,\phi_{u'})\nabla_i\Omega_x\nabla^i\phi_q,
\end{split}
\end{equation}
where the linear combination $\sum_{j\in J}\dots$ stands for a generic linear
 combination in the form described above. On the other hand, the linear combination
$\sum_{\tau\in T'}$ stands for a generic linear combination of tensor fields where we have two
 factors $\nabla^{(A)}\Omega_x,\nabla^{(B)}\phi_q$  with $A=2,B=1$ or $A=1,B=2$
respectively, and in each case the term with one derivative is contracting against
the other term (with two derivatives).
\end{lemma}

We claim that the above Lemma \ref{yarusal} (which we
will prove momentarily), when applied to (\ref{rixardos})
 implies our claim on the sublinear combination $\sum_{h\in H^{\S\S}}\dots$, in the case $\sigma>4$.
This follows immediately once we set $\phi_q=\Omega_{x'}$, and
once we observe that the tensor fields we obtain from
(\ref{rixardos}) by {\it erasing} the expression
$\nabla_i\Omega_x\nabla^i\Omega_{x'}$ have a simple character that
falls under the inductive assumption of Proposition \ref{giade}
(because we are increasing the weight).
\newline

{\it Proof of Lemma \ref{yarusal}:}

\par The proof follows the usual inductive scheme:

We will assume that for some $A\ge V$, we can write:

\begin{equation}
\label{hmartoncor'}
\begin{split}
&\sum_{\tau\in T} a_\tau Xdiv_{i_1}\dots Xdiv_{i_a}
[C^{\tau,i_1\dots i_a}_g(\Omega_1,\dots,\Omega_{p'},\phi_1,\dots ,\phi_{u'})
\nabla_i\Omega_x\nabla^i\phi_q]=
\\&\sum_{\tau\in T^A} a_\tau Xdiv_{i_1}\dots Xdiv_{i_a}
[C^{\tau,i_1\dots i_a}_g(\Omega_1,\dots,\Omega_{p'},\phi_1,\dots ,\phi_{u'})
\nabla_i\Omega_x\nabla^i\phi_q]+
\\&\sum_{\tau\in T'} a_\tau Xdiv_{i_1}\dots Xdiv_{i_a}
C^{\tau,i_1\dots i_a}_g(\Omega_1,\dots,\Omega_{p'},\phi_1,\dots ,\phi_{u'}|\Omega_x,\phi_q)
\\&+\sum_{j\in J} a_j C^{j}_{g}(\Omega_1,\dots ,\Omega_{p'},
 \phi_1,\dots ,\phi_{u'})\nabla_i\Omega_x\nabla^i\phi_q,
\end{split}
\end{equation}
where $\sum_{\tau\in T^A}$ on the RHS stands for a linear combination of tensor fields which are in the
general form of the tensor fields in $\sum_{\tau\in T}$, only with rank $\ge A$.
 We can use the above to replace $\sum_{\tau\in T}\dots$ in our
 Lemma hypothesis by $\sum_{\tau\in T^A}\dots$.
We will then show that we can write:

\begin{equation}
\label{hmartoncor''}
\begin{split}
&\sum_{\tau\in T^{A}} a_\tau Xdiv_{i_1}\dots Xdiv_{i_a}
[C^{\tau,i_1\dots i_a}_g(\Omega_1,\dots,\Omega_{p'},\phi_1,\dots ,\phi_{u'})
\nabla_i\Omega_x\nabla^i\phi_q]=
\\&\sum_{\tau\in T^{A+1}} a_\tau Xdiv_{i_1}\dots Xdiv_{i_a}
[C^{\tau,i_1\dots i_a}_g(\Omega_1,\dots,\Omega_{p'},\phi_1,\dots ,\phi_{u'})
\nabla_i\Omega_x\nabla^i\phi_q]+
\\&\sum_{\tau\in T'} a_\tau Xdiv_{i_1}\dots Xdiv_{i_a}
C^{\tau,i_1\dots i_a}_g(\Omega_1,\dots,\Omega_{p'},\phi_1,\dots ,\phi_{u'}|\Omega_x,\phi_q)
\\&+\sum_{j\in J} a_j C^{j}_{g}(\Omega_1,\dots ,\Omega_{p'},
 \phi_1,\dots ,\phi_{u'})\nabla_i\Omega_x\nabla^i\phi_q,
\end{split}
\end{equation}
with the same notational conventions. Clearly, if we can show that (\ref{hmartoncor'})
 implies (\ref{hmartoncor''}) then by iterating this step, we will derive our Lemma
\ref{yarusal}.
\newline

{\it Proof that (\ref{hmartoncor'}) implies (\ref{hmartoncor''}):}
As explained, we may assume that \\$\sum_{\tau\in T}\dots=\sum_{\tau\in T^A}$. We denote by $T_{*}^A\subset T^A$
the index set of tensor fields with rank exactly $A$. Then, applying the eraser to the expression
$\nabla_i\Omega_x\nabla^i\phi_q$, we derive an equation:

\begin{equation}
\label{hmartoncor'2}
\begin{split}
&\sum_{\tau\in T^A} a_\tau Xdiv_{i_1}\dots Xdiv_{i_a}
[C^{\tau,i_1\dots i_a}_g(\Omega_1,\dots,\Omega_{p'},\phi_1,\dots ,\phi_{u'})]+
\\&+\sum_{j\in J} a_j C^{j}_{g}(\Omega_1,\dots ,\Omega_{p'},
 \phi_1,\dots ,\phi_{u'})=0.
\end{split}
\end{equation}

\par Then, with certain exceptions,\footnote{These  exceptions are when there
are tensor fields in $T_{*}^A$ which are ``forbidden tensor fields
of Corollary 1 in \cite{alexakis4} with rank $m\ge \mu+1$. (It follows that
the forbidden tensor fields of rank $\mu$ cannot arise here, since all
$\mu$ free indices of the tensor fields in $H^{\S\S}$ must have
all their free indices being non-special).} we may apply Corollary
1 in \cite{alexakis4} to the above, and derive that there is some linear
combination of acceptable $(A+1)$-tensor fields with a $u'$-simple
character $\vec{\kappa}_{simp}'$ (indexed in $H$ below) so that:

\begin{equation}
\label{coldfeet}
\begin{split}
&\sum_{\tau\in T^A_{*}} a_\tau
C^{\tau,i_1\dots i_A}_g(\Omega_1,\dots,\Omega_{p'},\phi_1,\dots ,\phi_{u'})
\nabla_{i_1}\upsilon\dots \nabla_{i_A}\upsilon-
\\&\sum_{\tau\in T^A_{*}} a_\tau
Xdiv_{i_{A+1}}C^{\tau,i_1\dots i_{A+1}}_g(\Omega_1,\dots,\Omega_{p'},\phi_1,\dots ,\phi_{u'})
\nabla_{i_1}\upsilon\dots \nabla_{i_A}\upsilon=
\\&\sum_{j\in J} a_j C^{j,i_1\dots i_A}_{g}(\Omega_1,\dots ,\Omega_{p'},
 \phi_1,\dots ,\phi_{u'})\nabla_{i_1}\upsilon\dots\nabla_{i_A}\upsilon
\end{split}
\end{equation}
(the tensor fields indexed in $J$ are simply subsequent to
$\vec{\kappa}'_{simp}$).\footnote{In the exceptional cases above, our
claim (\ref{hmartoncor''}) follows from from the ``weak version'' of Proposition 
\ref{giade} presented in \cite{alexakis4}, with $\Phi=\nabla_s\Omega_x\nabla^s\phi_q$.}

\par Now, if we multiply the above by $\nabla_i\Omega_x\nabla^i\phi_q$
and make the $\nabla\upsilon$s into $Xdiv$'s
(which {\it are} allowed to hit either of the
 factors $\nabla\Omega_x,\nabla\phi_q$), then we derive
 (\ref{hmartoncor''}) and hence also  Lemma \ref{yarusal}. $\Box$
\newline

{\it The proof of (\ref{xiaodongcao}) in the case $\sigma=4$:}
\newline

\par In this case, we will prove (\ref{xiaodongcao}) directly,
immitating the ideas in \cite{alexakis3}:

{\it Proof of (\ref{xiaodongcao}) in the case where $\sigma=4$:}
We have that in this case, the tensor fields indexed in $H^{\S\S}$  will have two
factors $\nabla\Omega_x,\nabla\Omega_{x'}$ (contracting against each other)
 {\it and two other factors}, which we denote by $T_1,T_2$.\footnote{In the case where both
  $T_1,T_2$ are generic terms in the form $\nabla^{(m)}R_{ijkl}$, the labelling
   $T_1,T_2$ is arbitrary; in all other cases, we will have a
   well-defined factor $T_1$ and a well-defined $T_2$.} We also recall
 that all tensor fields indexed in $H^{\S\S}$ have rank $\ge \mu$ ($\ge 1$),
 and if they do have rank $\mu$ they will also have a removable index, by construction.
We distinguish the following cases regarding the {\it form} of the
factors $T_1$: Either both these factors are of the form
$\nabla^{(p)}\Omega_j$, or one ($T_1$, say) is in the form
$\nabla^{(p)}\Omega_j$ and the other is a curvature term (either
 in the form $\nabla^{(m)}R_{ijkl}$ or $S_{*}\nabla^{(\nu)}R_{ijkl}$) or both $T_1,T_2$ are
curvature factors (either in the form $\nabla^{(m)}R_{ijkl}$ or $S_{*}\nabla^{(\nu)}R_{ijkl}$).
Label these cases A,B,C respectively.

\par In case $A$, we will assume with
no loss of generality (up to re-labelling factors) that
$T_1=\nabla^{(c)}\Omega_1$,  $T_2=\nabla^{(c')}\Omega_2$ (and also
$\Omega_x=\Omega_3,\Omega_{x'}=\Omega_4$). Then, by ``manually''
constructing divergences, we can show that:

\begin{equation}
\label{burritou}
\begin{split}
& \sum_{h\in H^{\S\S}} a_h Xdiv_{i_1}\dots
Xdiv_{i_a} [C^{h,i_1\dots
i_a}_g(\Omega_1,\Omega_2,\phi_1,\dots,\phi_{u+1})\nabla_j\Omega_3\nabla^j\Omega_{4}]=
\\&\sum_{h\in H} a_h Xdiv_{i_1}\dots Xdiv_{i_a} [C^{h,i_1\dots
i_a}_g(\Omega_1,\dots,\Omega_4,\Omega\phi_1,\dots,\phi_{u+1})]+
\\&(Const)_{*} Xdiv_{i_1}\dots Xdiv_{i_b}[C^{*,i_1\dots
i_b}_{g}(\Omega_1,\Omega_2,\phi_1,\dots,\phi_{u+1})\nabla_j\Omega_3\nabla^j\Omega_{4}],
\end{split}
\end{equation}
where the tensor field $C^{*,i_1\dots
i_b}_{g}(\Omega_1,\Omega_2,\phi_1,\dots,\phi_{u+1})\nabla_j\Omega_3\nabla^j\Omega_{4}$
is in the form:
$$pcontr(\nabla^{(A)}_{r_1\dots r_A}\Omega_1\otimes\nabla^{(B)}_{t_1\dots t_B}\Omega_2
\otimes\nabla\phi_1\otimes\dots\otimes\nabla\phi_{u+1}),$$ with the
following restrictions: All indices in both
$\nabla^{(A)}\Omega_1,\nabla^{(B)}\Omega_2$ are either free or
contracting against some $\nabla\phi_h$. Also, if we denote by
$\beta$ the number of factors $\nabla\phi_h$ that are contracting
against $\nabla^{(B)}\Omega_2$ (notice $\beta$ is encoded in
$\vec{\kappa}_{simp}$, then $B=2$ if $\beta\le 2$, while $B=\beta$
if $\beta>2$. The linear combination $\sum_{h\in H}\dots$ on the
right hand side of the above stands for a {\it generic} linear
combination of the form $\sum_{t\in T'}\dots$ allowed in the RHS of (\ref{xiaodongcao}).

\par Then, using the above, we derive that $(Const)=0$, and that
concludes the proof of (\ref{xiaodongcao}) in this case.

\par In case B, using the same technique of constructing
``explicit'' $X$-divergences, we derive an equation
(\ref{burritou}) only {\it without} the last term
$(Const)_{*}\dots$. That immediately implies (\ref{xiaodongcao}) in
this case.

\par Finally in case C, we distinguish subcases on the factors
$T_1,T_2$: In subcase $(i)$, both factors will be in the form
$\nabla^{(m)}R_{ijkl}$. In subcase $(ii)$, $T_1$ will be in the
form $\nabla^{(m)}R_{ijkl}$ and $T_2$ will be in the form
$S_{*}\nabla^{(\nu)}R_{ijkl}$. In subcase $(iii)$, both $T_1,T_2$
will be in the form $S_{*}\nabla^{(\nu)}R_{ijkl}$.

\par Now, in subcase $(i)$ we show (\ref{burritou}) by the same
argument, only now the tensor field $C^{*,i_1\dots
i_b}_{g}(\Omega_1,\Omega_2,\phi_1,\dots,\phi_u)\nabla_j\Omega_3\nabla^j\Omega_{4}$
is in the form:
$$pcontr(\nabla^{(A)}_{r_1\dots r_A}R_{ijkl}\otimes\nabla^{(B)}_{t_1\dots
t_B}{{R_{i'}}^{jk}}_{l'}
\otimes\nabla\phi_1\otimes\dots\otimes\nabla\phi_{u+1}),$$ with the
following restrictions: The indices ${}_i,{}_l,{}_{i'},{}_{l'}$
are free. All derivative indices in both curvature factors are
either free or contracting against some $\nabla\phi_h$. Also, if
we denote by $\beta$ the number of factors $\nabla\phi_h$ that are
contracting against $T_2$ (notice $\beta$ is encoded in
$\vec{\kappa}_{simp}$,  $B=2$ if $\beta\le 2$, while $B=\beta$
if $\beta>2$). We again derive that $(Const)=0$, which implies that
(\ref{burritou}) is our desired equation (\ref{xiaodongcao}).

\par In subcase $(ii)$ we use this technique to derive an equation:

\begin{equation}
\label{burritou'}
\begin{split}
& \sum_{h\in H^{\S\S}} a_h Xdiv_{i_1}\dots
Xdiv_{i_a} [C^{h,i_1\dots
i_a}_g(\Omega_1,\Omega_2,\phi_1,\dots,\phi_{u+1})\nabla_j\Omega_3\nabla^j\Omega_{4}]=
\\&\sum_{h\in H} a_h Xdiv_{i_1}\dots Xdiv_{i_a} [C^{h,i_1\dots
i_a}_g(\Omega_1,\dots,\Omega_4,\phi_1,\dots,\phi_{u+1})]
\\&+\sum_{j\in J} a_j
[C^j_g(\Omega_1,\Omega_2,\phi_1,\dots,\phi_{u+1})\nabla_j
\Omega_3\nabla^j\Omega_{4}],
\end{split}
\end{equation}
where the terms indexed in $J$ are simply subsequent to
$\vec{\kappa}_{simp}$.

\par Finally, in subcase $(iii)$ we explicitly write:

\begin{equation}
\label{burritou''} \begin{split}
& \sum_{h\in H^{\S\S}} a_h
Xdiv_{i_1}\dots Xdiv_{i_a} [C^{h,i_1\dots
i_a}_g(\Omega_1,\Omega_2,\phi_1,\dots,\phi_{u+1})\nabla_j\Omega_3\nabla^j\Omega_{4}]=
\\&\sum_{h\in H} a_h Xdiv_{i_1}\dots Xdiv_{i_a} [C^{h,i_1\dots
i_a}_g(\Omega_1,\dots,\Omega_4,\Omega\phi_1,\dots,\phi_{u+1})]+
\\&(Const)_{*} Xdiv_{i_1}\dots Xdiv_{i_b}[C^{*,i_1\dots
i_b}_{g}(\Omega_1,\Omega_2,\phi_1,\dots,\phi_{u+1})\nabla_j\Omega_3\nabla^j\Omega_{4}]
\\&+\sum_{j\i J} a_j
[C^j_g(\Omega_1,\Omega_2,\phi_1,\dots,\phi_{u+1})\nabla_j
\Omega_3\nabla^j\Omega_{4}],
\end{split}
\end{equation}
where the tensor field $C^{*,i_1\dots
i_b}_{g}(\Omega_1,\Omega_2,\phi_1,\dots,\phi_u)\nabla_j\Omega_3\nabla^j\Omega_{4}$
is in the form:
$$pcontr(S_{*}\nabla^{(A)}_{r_1\dots r_A}R_{ijkl}\otimes S_{*}\nabla^{(B)}_{t_1\dots
t_B}{{R_{i'j'}}^{k}}_{l'}
\otimes\nabla\phi_1\otimes\dots\otimes\nabla\phi_{u+1}),$$ with the
following restrictions: The indices ${}_l,{}_{l'}$ are free. All
 indices ${}_{r_1},\dots ,{}_{r_A},{}_j$, ${}_{t_1},\dots ,{}_{t_B},{}_{j'}$
are either free or contracting against some $\nabla\phi_h$. Also,
if we denote by $\beta$ the number of factors $\nabla\phi'_h$ that
are contracting against $T_2$ (notice $\beta$ is encoded in
$\vec{\kappa}_{simp}$, then $B=1$ if $1\le \beta\le 2$, $B=0$ if
$\beta=0$, while $B=\beta-1$ if $\beta>2$. We again derive that
$(Const)=0$, which implies that (\ref{burritou''}) is our desired
equation (\ref{xiaodongcao}). $\Box$

\subsection{Mini-Appendix: Proof of  (\ref{grafw4}).}
\label{pfgrafw4}

 Firstly, we observe that for any
tensor field $C^{l,i_1\dots i_a}_{g}(\Omega_1,\dots
,\Omega_p,\phi_1,\dots,\phi_u,\phi_{u+1},\omega)$ with $a=\mu$ we
will not have any of the free indices ${}_{i_i},\dots ,{}_{i_\mu}$
belonging to any of the two factors
$\nabla^{(2)}\phi_{u+1},\nabla^{(2)}\omega$ (this holds because we are in
the setting of Lemma \ref{pskovb}) hence no $\mu$-tensor fields in
(\ref{hypothese2}) have special free indices in any factor
$S_{*}\nabla^{(\nu)}R_{ijkl}$, while any free index in any of the
factors $\nabla^{(2)}\phi_{u+1}$, $\nabla^{(2)}\omega_1$ would
necessarily  have arisen from a special free index in some factor
$S_{*}R_{ijkl}$ by (\ref{curvtrans}).

\par Now, we apply the eraser to the factors $\nabla\phi_h$ that are contracting
against $\nabla^{(2)}\phi_{u+1}$, $\nabla^{(2)}\omega$. We are left
 with factors $\nabla\phi_{u+1}$, $\nabla\omega$ and denote the
  tensor fields and complete contractions we
   are left with  by $\overline{C}^{l,i_1\dots i_a}_{g}$,
$\overline{C}^j_{g}$. Thus we obtain an equation:

\begin{equation}
\label{ianmorgan}
\begin{split}
& \Sum_{l\in L^{\alpha,\beta}} a_l X_{*}div_{i_1}\dots
X_{*}div_{i_a}\overline{C}^{l,i_1\dots i_a}_{g}(\Omega_1,\dots
,\Omega_p,\phi_1, \dots ,\phi_u,\phi_{u+1},\omega)+
\\&\Sum_{j\in J^{\alpha,\beta,II}}a_j \overline{C}^j_{g}
(\Omega_1,\dots ,\Omega_p,\phi_1, \dots,\phi_u,\phi_{u+1},
\omega)=0
\end{split}
\end{equation}
(modulo complete contractions of length $\ge \sigma + u+1$).

\par We regard
  the factor $\nabla\omega$ as a factor $\nabla\phi_{u+2}$.
 We observe that the tensor fields
$\overline{C}^{l,i_1\dots i_a}_{g}$ all have the same
$(u-2)$-simple character (the one
  defined by the factors $\nabla\phi_h, h\le u$), say
$\overline{\kappa}_{simp}$ and each $\overline{C}^j_{g}$ is simply
subsequent to that $(u-2)$-simple character. The tensor fields in
the above either have rank either $a\ge\mu+1$ but may contain free
indices in the factors $\nabla^{(2)}\phi_{u+1}$,
$\nabla^{(2)}\omega$, {\it or} they have rank $\mu$ and in
addition no free indices belong to the factors
$\nabla^{(2)}\phi_{u+1},\nabla^{(2)}\omega$. We denote by
$L^{\alpha,\beta,\sharp}$ the index set of tensor fields of rank
exactly $\mu+1$ where both factors $\nabla^{(2)}\phi_{u+1}$,
$\nabla^{(2)}\omega$
 contain a free index (say the indices ${}_{i_1},{}_{i_2}$ wlog).
We will prove that we can write:

\begin{equation}
\label{odromos}
\begin{split}
&\Sum_{l\in L^{\alpha,\beta,\sharp}} a_l Xdiv_{i_1}\dots
Xdiv_{i_{\mu+1}}\overline{C}^{l,i_1\dots i_a}_{g}(\Omega_1,\dots
,\Omega_p,\phi_1, \dots ,\phi_u,\phi_{u+1},\omega)=
\\&\Sum_{l\in L'^{\alpha,\beta}} a_l Xdiv_{i_1}\dots
Xdiv_{i_{\mu+1}}\overline{C}^{l,i_1\dots i_a}_{g}(\Omega_1,\dots
,\Omega_p,\phi_1, \dots ,\phi_u,\phi_{u+1},\omega)+
\\&\Sum_{j\in J^{\alpha,\beta,II}}a_j \overline{C}^j_{g}
(\Omega_1,\dots ,\Omega_p,\phi_1, \dots,\phi_u,\phi_{u+1},
\omega)
\end{split}
\end{equation}
(where the tensor fields indexed in $L'^{\alpha,\beta}$ have all
the features of the tensor fields indexed in $L^{\alpha,\beta}$
and in addition have the factor $\nabla^{(2)}\omega$ {\it not}
containing a free index {\it and} with one index in the
 factor $\nabla^{(2)}\phi_{u+1}$ contracting against a non-special index.

\par Notice that if we can prove the above, then we are reduced
to showing our claim under the additionnal assumption that
$L^{\alpha,\beta,\sharp}=\emptyset$.
Let us check how our claim then follows
under this assumtpion. We will then show (\ref{odromos}) below.
\newline

{\it Proof of our claim assuming (\ref{odromos}):} We
 break (\ref{ianmorgan}) into sublinear combinations with the same
 $u$-weak character\footnote{The one defined by $\nabla\omega$ and
$\nabla\phi_h$, $1\le h\le u+1$,
 $h\ne\alpha,h\ne\beta$.} (suppose those sublinear
 combinations are indexed in the sets $L^{\alpha,\beta,v}$, $v\in V$); we
 derive an equation for each $v\in V$:

\begin{equation}
\label{newfound}
\begin{split}
& \Sum_{l\in L^{\alpha,\beta,v}} a_l X_{*}div_{i_1}\dots
X_{*}div_{i_a}\overline{C}^{l,i_1\dots i_a}_{g}(\Omega_1,\dots
,\Omega_p,\phi_1, \dots ,\phi_u,\phi_{u+1},\omega)+
\\&\Sum_{j\in J^{\alpha,\beta,II}}a_j C^j_{g}(\Omega_1,
\dots ,\Omega_p,\phi_1, \dots,\phi_u,\phi_{u+1},\omega)=0.
\end{split}
\end{equation}
Now, by virtue of the assumption $L^{\alpha,\beta,\sharp}=\emptyset$, we may assume that all
the tensor fields in the above
 equation have rank $a\ge\mu$ and also have no free indices in the
 factors $\nabla\phi_{u+1},\nabla\phi_{u+2}$. Then, applying
 Lemma 2.5 in \cite{alexakis5}\footnote{A note to show why (\ref{newfound})
{\it does not} fall under the ``forbidden cases'' of Lemma 2.5 in \cite{alexakis5}:
We observe that the tensor fields of minimum rank $\mu$ in (\ref{newfound})
with both factors $\nabla\phi_{u+1},\nabla\omega$
contracting against special indices can only arise from the
$\mu$-tensor fields in (\ref{hypothese2})--but those
tensor fields will have no special free indices, thus
(\ref{newfound}) does not fall under a forbidden case of that Lemma.}
we derive (\ref{grafw4}). $\Box$
\newline

{\it Proof of (\ref{odromos}):} We initially pick out the sublinear combination
in (\ref{ianmorgan}) where both $\nabla\omega,\nabla\phi_{u+1}$
contract against the same factor $T$. Clearly this sublinear combination
must vanish separately, and we will denote the new true
equation that we thus obtain by New[(\ref{ianmorgan})].
Thus, the sublinear combination of tensor fields indexed in
$L^{\alpha,\beta}$ which contained at most one free index among the
 factors $\nabla\phi_{u+1},\nabla\omega$ contributes a
 linear combination of iterated $Xdiv$'s of rank at least $\mu$ to
New[(\ref{ianmorgan})]. We denote the equation we have obtained by:

\begin{equation}
\label{proxwra}
\begin{split}
&\sum_{l\in L^{\alpha,\beta,\sharp}} a_l
Doubdiv_{i_1i_2}Xdiv_{i_3}\dots Xdiv_{i_{\mu+1}}
\overline{C}^{l,i_1\dots i_{\mu+1}}_g+
\\&\sum_{f\in F} a_f Xdiv_{i_1}\dots Xdiv_{i_\mu}
 C^{f,i_1\dots i_\mu}_g(\nabla\phi_{u+1},\nabla\omega)+
\sum_{j\in J} a_j C^j_g =0
\end{split}
\end{equation}
($Doubdiv_{i_1i_2}$ means that both derivatives $\nabla^{i_1},\nabla^{i_2}$
are forced to hit
the same factor). We then symmetrize the two factors
$\nabla\phi_{u+1},\nabla\omega$ and thus obtain a new true
equation, which we denote by:

\begin{equation}
\label{proxwra}
\begin{split}
&\sum_{l\in L^{\alpha,\beta,\sharp}} a_l
Doubdiv_{i_1i_2}Sym[\overline{C}^{l,i_1\dots i_{\mu+1}}_g]Xdiv_{i_3}\dots Xdiv_{i_\mu}+
\\&\sum_{f\in F} a_f Xdiv_{i_1}\dots Xdiv_{i_\mu} Sym[C]^{f,i_1\dots i_\mu}_g
(\nabla\phi_{u+1},\nabla\omega)+\sum_{j\in J} a_j Sym[C]^j_g =0
\end{split}
\end{equation}
(here $Sym[\dots]$ stands for the symmetrization over the two
factors $\nabla\phi_{u+1},\nabla\omega$).

 We denote by  $F^a\subset F$ the index set of tensor fields for which
  both the factors $\nabla\phi_{u+1},\nabla\omega$
are contracting against internal indices in some factor
$\nabla^{(m)}R_{ijkl}$. We denote by $F^b\subset F$ the index set
of tensor fields for which  one of the indices
$\nabla\phi_{u+1},\nabla\omega$ are contracting against a special
index in some factor
$S_{*}\nabla^{(\nu)}R_{ijkl}\nabla^i\tilde{\phi}_x$.\footnote{We
will assume it is the factor $\nabla\omega$, wlog.}
 By replacing the
two factors $\nabla_a\phi_{u+1}\nabla_b\omega$ by $g_{ab}$ in the
first case and the two factors
$\nabla^i\tilde{\phi}_x\nabla^k\omega$ by $g^{ik}$ in the second,
and then using the operation $Ricto\Omega$ and iteratively
applying Corollary 1 in \cite{alexakis4},\footnote{A note to illustrate
why the ``forbidden cases'' of Corollary 1 in \cite{alexakis4} do not
interfere with our argument: Observe that for the terms indexed in
$F^b$ there will be a removable index by construction,
 therefore the ``forbidden cases'' do not
obstruct our iteration; the terms indexed in $F^a$ with rank $\mu$
must necessarily have arisen from the
tensor fields or rank $\mu$ in (\ref{hypothese2}).
Therefore they will have only non-special free indices,
therefore at the first iteration, Corollary 1 in \cite{alexakis4}
can be applied. On the other hand, it is possible
that at a subsequent step in the iteration
we may obtain a ``forbidden'' tensor field of rank $>\mu$;
in that case we use the ``weak substitute'' of Corollary 1 in \cite{alexakis4}, 
presented in the Appendix in \cite{alexakis6}.}
 we derive that we can write:

\begin{equation}
\label{kelarid}
\begin{split}
&\sum_{f\in F^a\bigcup F^b} a_f Xdiv_{i_1}\dots Xdiv_{i_\mu} Sym[C]^{f,i_1\dots i_\mu}_g
(\nabla\phi_{u+1},\nabla\omega)=
\\&\sum_{f\in F^{OK}} a_f Xdiv_{i_1}\dots Xdiv_{i_a} Sym[C]^{f,i_1\dots i_a}_g
(\nabla\phi_{u+1},\nabla\omega)+\sum_{j\in J} a_j Sym[C]^j_g,
\end{split}
\end{equation}
where the terms indexed in $F^{OK}$ have all the properties of
the terms indexed in $F$ in (\ref{proxwra}) and in addition
have at most one/none of the factors $\nabla\phi_{u+1},\nabla\omega$ contracting against
special indices in factors of the form $\nabla^{(m)}R_{ijkl},S_{*}\nabla^{(\nu)}R_{ijkl}$.
Therefore, we may assume that $F^a=F^b=\emptyset$ in (\ref{proxwra}).

\par Now, we refer to (\ref{proxwra}) and replace the expression
$\nabla_a\phi_{u+1}\nabla_b\omega$ by $g_{ab}$; we denote the resulting equation by
(\ref{proxwra})'. We then apply $Sub_\omega$ to (\ref{proxwra}) (see the Appendix in \cite{alexakis1})
(obtaining a new true equation which we denote by $D_g=0$) and we
then apply our inductive assumption of Lemma 4.10 in \cite{alexakis4} to
$D_g=0$.\footnote{When we apply Lemma 4.10 in \cite{alexakis4} we treat the
$Xdiv_{i_1}[\dots \nabla_{i_1}\omega]$ as a linear combination of
$(\mu-1)$-tensor fields--i.e.~we ``forget'' the $Xdiv$ structure
with respect to the factor $\nabla\omega$, thus  the terms of minimum rank
the factor $\nabla\omega$ contains a removable index, thus our assumption
does not fall under a ``forbidden case'' of that Lemma.}
In order to describe the resulting equation, we just denote by
$Cut[\overline{C}]^{l,i_1\dots i_{\mu+1}}_{g}$ the tensor field
that arises from $\overline{C}^{l,i_1\dots i_{\mu+1}}_{g}$ by
erasing the factor $\nabla_{i_1}\phi_{u+1}$ (along with the free
index ${}_{i_1}$). We thus derive that there exists a linear
combination of acceptable
 $\mu$-tensor fields (indexed in $K$ below),
 with a simple character $Cut(\vec{\kappa}_{simp})$
 {\it and with the index ${}_{i_{\mu+2}}$
 belonging to a real factor} so that:

\begin{equation}
\label{clever}
\begin{split}
&\Sum_{l\in L^{\alpha,\beta,\sharp}} a_l Xdiv_{i_2}
Cut[\overline{C}]^{l,i_1\dots i_{\mu+1}}_{g}(\Omega_1,\dots
,\Omega_p,\phi_1, \dots ,\phi_u,\phi_{u+1},\omega)
\nabla_{i_3}\upsilon\dots\nabla_{i_{\mu+1}}\upsilon
\\&=\sum_{k\in K} a_k Xdiv_{i_{\mu+2}}C^{k,i_3\dots i_{\mu+2}}_g (\Omega_1,\dots
,\Omega_p,\phi_1, \dots ,\phi_u,\omega)
\nabla_{i_3}\upsilon\dots\nabla_{i_{\mu+1}}\upsilon.
\end{split}
\end{equation}
Now, just multiplying the above by an expression
$\nabla_s\phi_{u+1}\nabla^s\upsilon$ and then
replacing the $\nabla\upsilon$s by $Xdiv$'s we derive (\ref{odromos}). $\Box$
\newline

In section \ref{caseB}, we derive the other half of \ref{pskovb}.

\section{Proof of  Lemma \ref{pskovb} in case B.}
\label{caseB}

\subsection{Introduction: A sketch of the strategy.} 
 
 In order to derive Lemma \ref{pskovb} (which corresponds
  to case B of Lemma 3.5 in \cite{alexakis4}) 
  we will use all the tools that were developed in thi paper.
 Most importantly the ``grand conclusion''  
 but also the two separate equations that were added
  in order to derive it the``grand conclusion''.
  \newline

{\it Main Strategy:} For each $\mu$-tensor field $C^{l,i_1\dots i_\mu}_g$, 
 $l\in L^z, z\in Z'_{Max}$ in (\ref{hypothese2}), we 
 will {\it canonicaly} pick out  a prescribed free index ${}_{i_1}$.
   We then consider the $(\mu-1)$-tensor fields 
$C^{l,i_1\dots i_{\mu}}_g\nabla_{i_1}\phi_{u+1}$, 
$l\in \bigcup_{z\in Z'_{Max}} L^z$.\footnote{These 
 $(\mu-1)$-tensor fields arise from $C^{l,i_1\dots i_\mu}_g$ by just 
 contracting the free index ${}_{i_1}$ against a new factor $\nabla\phi_{u+1}$.}
 We then prove (schematicaly) that there will exist a linear combination of $(\mu+1)$-tensor 
 fields, $\sum_{h\in H} a_h C^{h,i_1\dots i_{\mu+1}}_g\nabla_{i_1}\phi_{u+1}$,  
 each $C^{h,i_1\dots i_{\mu+1}}_g$ a partial contraction in the form (\ref{form2}), 
with the same $u$-simple character $\vec{\kappa}_{simp}$, such that:

\begin{equation}
\label{schemclaimpskovbB} 
\begin{split}
&\Sum_{l\in L^z} a_l 
Xdiv_{i_2}\dots Xdiv_{i_\mu}C^{l,i_1\dots i_{\mu}}_{g} (\Omega_1,\dots
,\Omega_p,\phi_1,\dots ,\phi_u)\nabla_{i_1}\phi_{u+1}=
\\&\sum_{h\in H} a_h Xdiv_{i_2}\dots Xdiv_{i_{\mu+1}}
C^{h,i_1\dots i_{\mu+1}}_g\nabla_{i_1}\phi_{u+1}\\&+
\sum_{j\in J} a_j C^{j,i_1}_g (\Omega_1,\dots
,\Omega_p,\phi_1,\dots ,\phi_u)\nabla_{i_1}\phi_{u+1};
\end{split}
\end{equation}
here the terms indexed in $J$ are ``junk terms''; they have length $\sigma+u$ (like the 
tensor felds indexed in $L_1$ and $H$)
and are in the general form (\ref{form1}). They are ``junk terms'' because 
one of the factors $\nabla\phi_h, 1\le h\le u$) which 
are supposed to contract against the
 index ${}_i$ in some factor $S_*\nabla^{(\nu)}R_{ijkl}$
  for all the tensor fields indexed in $L_\mu$
  now contracts against a derivative index of some factor $\nabla^{(m)}R_{ijkl}$.\footnote{In 
  the formal language introduced in \cite{alexakis4}, in this second scenario
  we would say that $C^{j,i_1}_g(\Omega_1,\dots
,\Omega_p,\phi_1,\dots ,\phi_u)$ is ``simply subequent'' to
 the simple character $\vec{\kappa}_{simp}$. }
 \newline
 
 After we have derived an equation  of the form 
 (\ref{schemclaimpskovbB}), 
 Lemma \ref{pskovb} follows by just applying the inductive
  claim of Proposition \ref{giade} to the above. 
\newline

In order to derive (\ref{schemclaimpskovbB}), we 
 will subdivide case B of Lemma \ref{pskovb} into subcases and treat them separately. 
 In certain cases we must derive {\it systems} of equations combining the ``grand conclusion'' 
 with other equations that we derived above. In certain very special subcases 
 (such as when $\mu=1$ in (\ref{hypothese2})),
  we will resort to ad hoc methods to derive (\ref{schemclaimpskovbB}).

 We wish to stress again that when proving of Lemma \ref{pskovb} 
 we are stll making all the inductive asumptions (on the parameters, 
 $n, \sigma, \Phi, \sigma_1+\sigma_2$) regarding the validity of Proposition \ref{giade} 
 and also all of its consequences. Hence we are allowed to apply our inductive assumption of 
 Proposition \ref{giade} or Lemmas 4.6, 4.8 etc from \cite{alexakis4}
 \newline
 
 {\bf The subcases of Lemma \ref{pskovb}:} We
   distinguish subcases for Lemma \ref{pskovb} according to the maximal refined 
   double character among the $\mu$-tensor fields in (\ref{hypothese2}). 
   We  refer the reader toto the introduction 
  for a loose discussion of the notion of maximal refined 
  double characters.\footnote{The proper definition appears 
in \cite{alexakis4}} In particular, we recall that 
the maximal refined double character contains a decreasing list of numbers, 
$\vec{R\lambda}_{Max}$, which corresponds to the distribution of free
indices among the different factors in the 
$\mu$-tensor fields in $\bigcup_{z\in Z'_{Max}} L^z$. 
We have denoted $\vec{R\lambda}_{Max}=(M,B_1,\dots, B_\pi)$.
  The subcases are then as follows:

\begin{enumerate}

\item{$M=1$, $\pi> 0$.}

 \item{$M\ge 2$ and
$B_1=\dots B_\pi=1$, $\pi>1$.}

\item{$M=\mu\ge 3$.}

\item{$M=\mu-1\ge 2$.}

\item{$M=\mu=2$.}

\item{$M=\mu=1$.}
\end{enumerate}

 {\bf Technical remarks regarding the ``grand conclusion'':}  
The ``grand conclusion'' is a new local equation, which is a consequence of 
(\ref{hypothese2}); it applies in the setting of Lemma
 \ref{pskovb}. It will be one of the main tools 
in deriving Lemma \ref{pskovb} in the present 
paper. For the reader's convenience, we recall 
a certaine conventions which we hav introduced:
\newline

{\bf Recall conventions}: Recall firstly that the ``grand conclusion'' 
is derived once we specify a particular factor/set of factors in 
$\vec{\kappa}_{simp}$.\footnote{In particular, ``specifying  one factor'' means that 
we pick out a factor in (\ref{form2}) which is either in the form, 
$\nabla^{(B)}\Omega_x$, for some given $x$, or in the 
form $S_{*}\nabla^{(\nu)}R_{ijkl}\nabla^i\tilde{\phi}_h$ 
for some given $h$, or in the form 
$\nabla^{(m)}_{r_1\dots r_m}R_{ijkl}\nabla^{r_a}\phi_{h'}$ for some given $h'$. Specifying  a 
``set of factors'' means that we pick out the set of factors 
$\nabla^{(m)}R_{ijkl}$ in $\vec{\kappa}_{simp}$ which are 
not contracting against any factor $\nabla\phi_h$.} This is called 
the ``selected factor''/``selected set of factors''.  
Given such a choice of factor/set of factors, we construct 
a new $(u+1)$-simple character $\vec{\kappa}^+_{simp}$ by 
contracting the chosen factor/one of the set of 
chosen factors against a new factor $\nabla\phi_{u+1}$; the 
new factor $\nabla\phi_{u+1}$ must {\it not} 
contract against a special index. Given such a choice
 of chosen factor/set of factors (and thus a new $(u+1)$-simple character), 
the grand conclusion will involve tensor fields  in the form (\ref{form2})
 with a $(u+1)$-simple character $\vec{\kappa}^+_{simp}$ and 
complete contractions in the form (\ref{form1}) 
with a weak $(u+1)$-character $Weak(\vec{\kappa}_{simp})$. When 
the selected factor is in the form $S_*\nabla^{(\nu)}R_{ijkl}$, the grand conclusion 
is the equation (\ref{olaxreiazontai1}). When 
the selected factor is in the form $\nabla^{(m)}R_{ijkl}$
it is the equation (\ref{olaxreiazontai2}), while when 
it is in the form $\nabla^{(A)}\Omega_h$ it is 
the equation (\ref{olaxreiazontai3}).

\subsection{A useful technical Lemma.}

The next Lemma will allow us to assume wlog that all the tensor fileds indexed in $H$ 
in the grand conclusion (i.e.~all ``contributors'' there) are acceptable, and have 
the factor $\nabla\phi_{u+1}$ {\it not} contracting against a special index. 
\newline

 The ``Technical Lemma'' below has certain
``forbidden cases'', which we spell out here
for reference purposes. A tensor field $C^{l,i_1\dots i_\mu}_g$ in (\ref{hypothese2})  
is ``forbidden'' (for the purposes of the next Lemma) if
it has $\sigma_2>0$, each of the $\mu$ free 
indices belonging to a different factor, 
all factors $\nabla^{(\nu)}R_{ijkl}$/$\nabla^{(p)}\Omega_h$ 
must contract against none/at most one factor $\nabla\phi_y$,
 and  {\it either}
there are no removable free  indices,\footnote{(See definition 4.1 in \cite{alexakis4}). 
In this setting, all factors must be 
in the form $R_{ijkl},S_*R_{ijkl}$ without derivatives, 
{\it or} in the form $\nabla^{(2)}\Omega_h$.} {\it or} there is exaclty 
one removable free index.\footnote{One way to think of this is that 
for such a $\mu$-tensor field one free index belongs to a factor 
$\nabla_{(free)}R_{\sharp\sharp\sharp\sharp}$ 
or $\nabla^{(3)}_{(free)\sharp\sharp}\Omega_h$,
and {\it all} its other factors are in the form  
$R_{ijkl},S_*R_{ijkl}$ or $\nabla^{(2)}\Omega_h$.}  

{\it Remark:} There ``forbidden cases'' will only 
force us to give a special proof of Lemma \ref{pskovb} 
in the subcase $\mu=1$, when the terms in (\ref{hypothese2}) are ``forbidden'' as 
defined above. Those cases will be treated in the 
Mini-Appendix at the end of this paper.

\begin{lemma}
\label{funny} Refer to the grand conclusion. We denote by $H_{Bad,1}\subset H$ the index
set of tensor fields in $H$ which have an unacceptable factor
$\nabla\Omega_h$. We denote by $H_{Bad,2}\subset H$ the index set
of tensor fields in $H$ which have the factor $\nabla\phi_{u+1}$
contracting against a special index,\footnote{Therefore, these
tensor fields will be acceptable by Definition
\ref{contributeur}.} if $\sigma\ge 4$ (if $\sigma=3$
we just set $H_{Bad,2}=\emptyset$).

\par Then, ({\it unless} the tensor fields of maximal refined double character
in (\ref{hypothese2}) are in one of the  ``forbidden forms'' above)
we claim that we can write:

\begin{equation}
\label{fnnyeq}
\begin{split}
& \sum_{h\in H_{Bad,1}\bigcup H_{Bad,2}} a_h
 Xdiv_{i_1}\dots Xdiv_{i_a}
C^{h,i_1\dots i_a,i_{*}}_{g}(\Omega_1,\dots ,\Omega_p,
 \phi_1,\dots ,\phi_u)\nabla_{i_{*}}\phi_{u+1}
 \\&=\sum_{h\in H_{OK}} a_h
 Xdiv_{i_1}\dots Xdiv_{i_a}
C^{h,i_1\dots i_a,i_{*}}_{g}(\Omega_1,\dots ,\Omega_p,
 \phi_1,\dots ,\phi_u)\nabla_{i_{*}}\phi_{u+1}+
 \\&\sum_{j\in J} a_j C^{j,i_{*}}_{g}(\Omega_1,\dots ,\Omega_p,
 \phi_1,\dots ,\phi_u)\nabla_{i_{*}}\phi_{u+1}.
\end{split}
\end{equation}
The terms indexed in $H_{OK}$ in the RHS stand for a generic linear combination of
acceptable contributors\footnote{See definition
\ref{contributeur}.} with a $(u+1)$-simple character
$\vec{\kappa}^{+}_{simp}$. The terms indexed in $J$ are $u$-simply
subsequent to $\vec{\kappa}_{simp}$.
\end{lemma}

\par We observe that if we can show the above then we can assume
wlog that all tensor fields in the grand conclusion are acceptable
and have a $(u+1)$-simple character $\vec{\kappa}_{simp}^{+}$.
\newline

{\it Proof of Lemma \ref{funny}:} We divide the index set
$H_{Bad,1}$ into subsets $H_{Bad,1}^\alpha$,
$H_{Bad,1}^\beta$
according to whether the factor $\nabla\Omega_h$ is contracting
against a factor $\nabla\phi_{u+1}$ or not, respectively.

We firstly pick out the sublinear combination in the grand
conclusion with a factor $\nabla\Omega_h$ contracting against the
factor $\nabla\phi_{u+1}$. We thus derive a new equation:

\begin{equation}
\label{funnyeqn} \begin{split} & \sum_{h\in H_{Bad,1}^\alpha} a_h
X_{*}div_{i_1}\dots X_{*}div_{i_a} C^{h,i_1\dots
i_a,i_{*}}_{g}(\Omega_1,\dots ,\Omega_p,
 \phi_1,\dots ,\phi_u)\nabla_{i_{*}}\phi_{u+1}
 \\&+\sum_{j\in J} a_j C^{j,i_{*}}_{g}(\Omega_1,\dots ,\Omega_p,
 \phi_1,\dots ,\phi_u)\nabla_{i_{*}}\phi_{u+1}=0.
\end{split}
\end{equation}
Then applying Lemma 4.1 from \cite{alexakis4}\footnote{The fact that 
we have excluded the forbidden cases ensures that 
 the terms of minimum rank 
in (\ref{funnyeqn}) does
not fall under the ``forbidden case'' of that Lemma.}
or 4.2 from \cite{alexakis4}\footnote{The fact that 
we have excluded the forbidden cases ensures that 
 the terms of minimum rank 
in (\ref{funnyeqn}) does
not fall under the ``forbidden case'' of that Lemma.} we
derive that we can write:

\begin{equation}
\label{fnnyeqb}
\begin{split}
& \sum_{h\in H^\alpha_{Bad,1}} a_h
 Xdiv_{i_1}\dots Xdiv_{i_a}
C^{h,i_1\dots i_a,i_{*}}_{g}(\Omega_1,\dots ,\Omega_p,
 \phi_1,\dots ,\phi_u)\nabla_{i_{*}}\phi_{u+1}
 \\&=\sum_{h\in H_{OK}} a_h
 Xdiv_{i_1}\dots Xdiv_{i_a}
C^{h,i_1\dots i_a,i_{*}}_{g}(\Omega_1,\dots ,\Omega_p,
 \phi_1,\dots ,\phi_u)\nabla_{i_{*}}\phi_{u+1}
 \\&+\sum_{j\in J} a_j C^{j,i_{*}}_{g}(\Omega_1,\dots ,\Omega_p,
 \phi_1,\dots ,\phi_u)\nabla_{i_{*}}\phi_{u+1}.
\end{split}
\end{equation}
Thus, we are reduced to the case $H^\alpha_{Bad,1}=\emptyset$.
Now, we pick out the sublinear combination in the grand conclusion
with a factor $\nabla\Omega_h$ not contracting against
$\nabla\phi_{u+1}$. We thus derive an equation:
\begin{equation}
\label{funnyeqn2} \begin{split}& \sum_{h\in H_{Bad,1}^\beta} a_h
X_{*}div_{i_1}\dots X_{*}div_{i_a} C^{h,i_1\dots
i_a,i_{*}}_{g}(\Omega_1,\dots ,\Omega_p,
 \phi_1,\dots ,\phi_u)\nabla_{i_{*}}\phi_{u+1}
 \\&+\sum_{j\in J} a_j C^{j,i_{*}}_{g}(\Omega_1,\dots ,\Omega_p,
 \phi_1,\dots ,\phi_u)\nabla_{i_{*}}\phi_{u+1}=0.
\end{split}
\end{equation}
Now applying Corollary 4.6 in \cite{alexakis4}\footnote{The fact that 
we have excluded the forbidden cases ensures that 
 the terms of minimum rank 
in (\ref{funnyeqn}) does
not fall under the ``forbidden case'' of that Lemma.}
 (if $\sigma\ge 4$) or Lemma
4.7 in \cite{alexakis4} (if $\sigma=3$)\footnote{The fact that 
we have excluded the forbidden cases ensures that 
 the terms of minimum rank 
in (\ref{funnyeqn}) does
not fall under the ``forbidden case'' of that Lemma.} to the above we derive that we
can write:
\begin{equation}
\label{funnyeqn22} \begin{split} & \sum_{h\in H^\beta_{Bad,1}} a_h
Xdiv_{i_1}\dots Xdiv_{i_a} C^{h,i_1\dots
i_a,i_{*}}_{g}(\Omega_1,\dots ,\Omega_p,
 \phi_1,\dots ,\phi_u)\nabla_{i_{*}}\phi_{u+1}=
 \\&\sum_{h\in H_{OK}} a_h
 Xdiv_{i_1}\dots Xdiv_{i_a}
C^{h,i_1\dots i_a,i_{*}}_{g}(\Omega_1,\dots ,\Omega_p,
 \phi_1,\dots ,\phi_u)\nabla_{i_{*}}\phi_{u+1}+
 \\&\sum_{j\in J} a_j C^{j,i_{*}}_{g}(\Omega_1,\dots ,\Omega_p,
 \phi_1,\dots ,\phi_u)\nabla_{i_{*}}\phi_{u+1}.
\end{split}
\end{equation}
Thus, we may additionally assume that $H^\beta_{Bad,1}=\emptyset$.
Finally, applying Lemma 4.10 in \cite{alexakis4}\footnote{The fact that 
we have excluded the forbidden cases ensures that 
 the terms of minimum rank 
in (\ref{funnyeqn}) does
not fall under the ``forbidden case'' of that Lemma.} to the above we
derive that we can write:

\begin{equation}
\label{funnyeqn22} \begin{split} & \sum_{h\in H_{Bad,2}} a_h
Xdiv_{i_1}\dots Xdiv_{i_a} C^{h,i_1\dots
i_a,i_{*}}_{g}(\Omega_1,\dots ,\Omega_p,
 \phi_1,\dots ,\phi_u)\nabla_{i_{*}}\phi_{u+1}
 \\&=\sum_{h\in H_{OK}} a_h
 Xdiv_{i_1}\dots Xdiv_{i_a}
C^{h,i_1\dots i_a,i_{*}}_{g}(\Omega_1,\dots ,\Omega_p,
 \phi_1,\dots ,\phi_u)\nabla_{i_{*}}\phi_{u+1}
 \\&+\sum_{j\in J} a_j C^{j,i_{*}}_{g}(\Omega_1,\dots ,\Omega_p,
 \phi_1,\dots ,\phi_u)\nabla_{i_{*}}\phi_{u+1}.
\end{split}
\end{equation}
This concludes the proof of our Lemma. $\Box$

\subsection{Proof of Lemma \ref{pskovb} in the subcase $M=1,\pi>0$:}

\par In this case, it follows by the definition of
the maximal refined double character that {\it all} $\mu$-tensor
fields $C^{l,i_1\dots i_\mu}_{g}$ in $L_\mu$ must have
$M=1,\pi=\mu-1>0$.

\par In this setting, we claim:

\begin{equation}
\label{apery}
\begin{split}
&\sum_{l\in L_\mu} a_l \sum_{k=1}^\mu
Xdiv_{i_1}\dots\hat{Xdiv}_{i_k}\dots Xdiv_{i_\mu} C^{l,i_1\dots
i_\mu}_{g}(\Omega_1,\dots ,\Omega_p,\phi_1,\dots
,\phi_u)\nabla_{i_k}\phi_{u+1}
\\&+\sum_{h\in H\bigcup H'} a_h Xdiv_{i_1}\dots Xdiv_{i_a} C^{h,i_1\dots i_ai_{a+1}}_{g}
(\Omega_1,\dots ,\Omega_p,\phi_1,\dots
,\phi_u)\nabla_{i_{a+1}}\phi_{u+1}
\\&+\sum_{j\in J} a_j
 C^{j,i_1}_{g}(\Omega_1,\dots ,\Omega_p,\phi_1,\dots ,\phi_u)\nabla_{i_1}\phi_{u+1},
\end{split}
\end{equation}
where the tensor fields indexed in $H$ have the property that they
are acceptable
 with a $u$-simple character $\vec{\kappa}_{simp}$ and they satisfy $a\ge \mu$.
  The tensor fields indexed in $H'$
have $a\ge \mu$, and are contributors, see 
Definition \ref{contributeur} above; in particular they 
have a $u$-simple character $\vec{\kappa}_{simp}$
but they also  have one unacceptable factor $\nabla\Omega_h$ (with only
one derivative).
The complete contractions in $J$ are $u$-simply subsequent to
$\vec{\kappa}_{simp}$.
\newline

\par We will now show how Lemma \ref{pskovb} can
be derived from (\ref{apery}) in this subcase.
\newline

{\it Lemma \ref{pskovb} follows from (\ref{apery}) in this subcase:}
Firstly, we observe that by breaking (\ref{apery}) into sublinear
 combinations with the same weak $(u+1)$-character, we obtain a new set of
 true equations (since (\ref{apery}) holds formally and the weak character
is invariant under the permutations that make the LHS of
(\ref{apery}) vanish formally). So, for each $z\in Z'_{Max}$ 
 we pick out the sublinear combination in (\ref{apery}) with a 
 $(u+1)$-weak character
$Weak(\vec{\kappa}^{z}_{ref-doub})$. We assume with no loss of
generality (just by re-labelling free indices) that the sublinear combination of contractions in the
first line of (\ref{apery}) with weak character
$Weak(\vec{\kappa}^{z}_{ref-doub})$ are the summands $k=1,\dots
,V_z$; we also denote by $H_z,H'_z,J_z$ the index sets of
contractions with a  weak character
$Weak(\vec{\kappa}^{z}_{ref-doub})$. We denote $V_z=V$ for brevity
and thus derive an equation:

\begin{equation}
\label{apery'}
\begin{split}
&\sum_{l\in L_\mu} a_l \sum_{k=1}^V
Xdiv_{i_1}\dots\hat{Xdiv}_{i_k}\dots Xdiv_{i_\mu} C^{l,i_1\dots
i_\mu}_{g}(\Omega_1,\dots ,\Omega_p,\phi_1,\dots
,\phi_u)\nabla_{i_k}\phi_{u+1}
\\&+\sum_{h\in H_z\bigcup H_z'} a_h Xdiv_{i_1}\dots Xdiv_{i_a} C^{h,i_1\dots i_ai_{a+1}}_{g}
(\Omega_1,\dots ,\Omega_p,\phi_1,\dots
,\phi_u)\nabla_{i_{a+1}}\phi_{u+1}
\\&+\sum_{j\in J_z} a_j
 C^{j,i_1}_{g}(\Omega_1,\dots ,\Omega_p,\phi_1,\dots ,\phi_u)\nabla_{i_1}\phi_{u+1}.
\end{split}
\end{equation}

\par Now, our aim is to apply the inductive
assumption of Corollary 1 in \cite{alexakis4} to (\ref{apery'}). In order
to do this, we first apply Lemma \ref{funny} to
(\ref{apery'}) to derive a new equation where $H=\emptyset$ 
(thus all tensor fields are
acceptable and have the same $(u+1)$-simple character).

Thus, we apply the inductive assumption of Corollary 
1 in \cite{alexakis4} to (\ref{apery'}):\footnote{The 
 above equation falls under the inductive assumption 
 of Corollary 1 in \cite{alexakis4} because we have increased the value of $\Phi$,
  while keeping all the other parameters fixed. Observe
 that the tensor fields of minimum rank in (\ref{apery'}) will
 have only non-special free indices. Thus
 there is no danger of falling under a ``forbidden case'' of that Corollary.}
For our chosen $z\in Z'_{Max}$ we derive that there is a linear
combination  of $(\mu+1)$-tensor fields with a refined double
character $\vec{\kappa}^z_{ref-doub}$ (indexed in $P$ below) so
that:

\begin{equation}
\label{gtian}
\begin{split}
&\sum_{l\in L^z_\mu} a_l C^{l,i_1\dots i_\mu}_{g}(\Omega_1,\dots ,
\Omega_p,\phi_1,\dots ,\phi_u)\nabla_{i_1}\phi_{u+1}
\nabla_{i_2}\upsilon\dots\nabla_{i_\mu}\upsilon-
\\&\sum_{p\in P} a_P Xdiv_{i_{\mu+1}}C^{p,i_1\dots i_{\mu+1}}_{g}(\Omega_1,\dots ,
\Omega_p,\phi_1,\dots ,\phi_u)\nabla_{i_1}\phi_{u+1}
\nabla_{i_2}\upsilon\dots\nabla_{i_\mu}\upsilon=
\\&\sum_{j\in J} a_j C^j_{g}(\Omega_1,\dots ,
\Omega_p,\phi_1,\dots ,\phi_{u+1},\upsilon^{\mu-1});
\end{split}
\end{equation}
here the contractions indexed in $J$ are simply
 subsequent to $\vec{\kappa}_{simp}$. Now, setting
$\phi_{u+1}=\upsilon$, we derive our claim in this case. $\Box$
\newline

{\it Proof of (\ref{apery}):} This equation just follows by
considering the equation \\$Im^{1,\beta}_{\phi_{u+1}}[L_g]=0$ (see (\ref{antrikos}) and
then replacing $\nabla\omega$ by an $Xdiv$, by virtue of the last 
Lemma in the Appendix of \cite{alexakis1}. $\Box$

\subsection{Proof of Lemma \ref{pskovb} in the subcase $M\ge 2,B_1=1, \pi>1$:}

\par This subcase follows by a very similar reasoning.
We arbitrarily pick out some $z\in Z'_{Max}$ and we will show our
claim for the tensor fields indexed in $L^z$. If we can do this
then by induction we can derive Lemma \ref{pskovb} in this subcase.
 We recall that for each $l\in L^z, z\in Z'_{Max}$, $C^{l,i_1\dots
i_\mu}_{g}$ has one factor with $M>1$ free indices and all other
$\pi>1$  factors that contain free indices will each contain only one free index.
 Therefore, by the definition of {\it maximal}
 refined double character (see the beginning of this section)
for any {\it non-maximal} $C^{l,i_1\dots i_\mu}$, any given factor
will contain at most $M-1$ free indices. For notational
convenience we assume wlog that for each $l\in L^z$, the indices
${}_{i_1}, \dots ,{}_{i_M}$ in $C^{l,i_1\dots i_\mu}_{g}$ belong
to the same factor.

\par We will prove our claim in this case by considering the equation
$Im^{1,\beta}_{\phi_{u+1}}[L_{g}]=0$ (i.e.~(\ref{antrikos})),
where we now set $\omega=\phi_{u+2}$. In order to analyze this
equation and derive our claim, we will introduce some notation:

{\it Notation:} We denote by $\vec{\kappa}^{++}_z$ the refined
$(u+2,\mu -2)$-double character that arises from $\vec{L^z}$ as
 follows: Consider all
the refined double characters of the $(\mu-2)$-tensor fields
$C^{l,i_1\dots i_\mu}_{g}\nabla_{i_\alpha}\phi_{u+1}
\nabla_{i_\beta}\phi_{u+2}$, $\alpha,\beta>M$. Let
$\vec{\kappa}^{++}_z$ be the {\it maximal} such
$(u+2,\mu-2)$-refined double character (if there are many such
refined double characters we pick out one arbitrarily). We will
write $\vec{\kappa}^{++}$ instead of $\vec{\kappa}^{++}_z$ for
brevity.

We assume with no loss of generality (and only for notational
convenience) that $C^{l,i_1\dots
i_\mu}_{g}\nabla_{i_{\mu-1}}\phi_{u+1}\nabla_{i_\mu}\phi_{u+2}$
 has a maximal refined double character $\vec{\kappa}^{++}$.
 We denote by
$$\Sum_{l\in L_{\vec{\kappa}^{++}}} a_l C^{l,i_1\dots
i_\mu}_{g}(\Omega_1,\dots , \Omega_p, \phi_1,\dots
,\phi_u)\nabla_{i_\alpha}\phi_{u+1}\nabla_{i_\beta} \phi_{u+2}$$
 the sublinear combination in
$$\Sum_{\alpha,\beta>M}
\Sum_{l\in L^z} a_l C^{l,i_1\dots i_\mu}_{g}(\Omega_1,\dots ,
\Omega_p, \phi_1,\dots
,\phi_u)\nabla_{i_\alpha}\phi_{u+1}\nabla_{i_\beta} \phi_{u+2}$$
that consists of complete contractions with a refined double
character $\vec{\kappa}^{++}$. By definition, it follows that
there is a {\it nonzero, universal} combinatorial constant
$(Const)$,\footnote{By ``universal'' we mean that it depends only
on the refined double character $\vec{L}^z$.} for which:

\begin{equation}
\label{bastata}
\begin{split}
&\Sum_{l\in L_{\vec{\kappa}^{++}}} a_l C^{l,i_1\dots
i_\mu}_{g}(\Omega_1,\dots , \Omega_p, \phi_1,\dots
,\phi_u)\nabla_{i_1}\upsilon\dots\nabla_{i_\mu}\upsilon=
\\&(Const)\Sum_{l\in L^z} a_l \cdot C^{l,i_1\dots
i_\mu}_{g}(\Omega_1,\dots , \Omega_p, \phi_1,\dots
,\phi_u)\nabla_{i_1}\upsilon\dots\nabla_{i_\mu}\upsilon.
\end{split}
\end{equation}

\par Now, refer to the grand conclusion; we observe that for each $l\in L^z$,
 any $(\mu -2)$-tensor field
$C^{l,i_1\dots i_\mu}_{g}(\Omega_1,\dots , \Omega_p, \phi_1,\dots
,\phi_u)\nabla_{i_\alpha}\phi_{u+1}\nabla_{i_\beta}\phi_{u+2}$
with either $\alpha\le M$ or $\beta\le M$, has a weak character
that is different from $Weak(\vec{\kappa}^{++})$. In addition, we
observe by the definition of refined double characters  that for each $l\in L_\mu\setminus L^z$ all the
$(\mu-2)$-tensor fields $C^{l,i_1\dots i_\mu}_{g}(\Omega_1,\dots ,
\Omega_p, \phi_1,\dots ,\phi_u)\nabla_{i_\alpha}\phi_{u+1}
\nabla_{i_\beta}\phi_{u+2}$ with a weak character
$Weak(\vec{\kappa}^{++})$ have a $(u+2)$-simple character
$Simp(\vec{\kappa}^{++})$ and a refined double character that is
either subsequent to, or equipolent to $\vec{\kappa}^{++}$.
Therefore, with these conventions if we pick out the sublinear
combination with weak character $Weak(\vec{\kappa}^{++})$ (recall
that this sublinear combination must vanish separately) in the
equation $Im^{1,\beta}_{\phi_{u+1}}[L_{g}]=0$ we derive a new equation:

\begin{equation}
\label{reputations2}
\begin{split}
&\Sum_{l\in L_{\vec{\kappa}^{++}}} a_l
Xdiv_{i_1}\dots\hat{Xdiv}_{i_\alpha}\dots
\hat{Xdiv}_{i_\beta}\dots  Xdiv_{i_\mu}C^{l,i_1\dots
i_\mu}_{g}(\Omega_1,\dots , \Omega_p, \phi_1,\dots
,\phi_u)\\&\nabla_{i_\alpha}\phi_{u+1}\nabla_{i_\beta}\phi_{u+2}+
\Sum_{l\in L'} a_l Xdiv_{i_1}\dots\hat{Xdiv}_{i_\alpha}\dots
\hat{Xdiv}_{i_\beta}\dots  Xdiv_{i_\mu} \\&C^{l,i_1\dots
i_\mu}_{g}(\Omega_1,\dots , \Omega_p, \phi_1,\dots
,\phi_u)\nabla_{i_\alpha}\phi_{u+1}\nabla_{i_\beta}\phi_{u+2}+
\\&\Sum_{h\in H} a_h Xdiv_{i_1}\dots\hat{Xdiv}_{i_\alpha}\\&\dots
\hat{Xdiv}_{i_\beta}\dots Xdiv_{i_{\mu+1}}C^{h,i_1\dots
i_{\mu+1}}_{g}(\Omega_1,\dots , \Omega_p, \phi_1,\dots
,\phi_u)\nabla_{i_\alpha}\phi_{u+1}\nabla_{i_\beta}\phi_{u+2}+
\\&\Sum_{j\in J} a_j C^{j,i_{*}i_{**}}_{g}(\Omega_1,\dots , \Omega_p,
\phi_1,\dots ,\phi_u)\nabla_{i_{*}}\phi_{u+1}\nabla_{i_{**}}.
\phi_{u+2}=0.
\end{split}
\end{equation}
Here the tensor fields indexed in $L'$ are all acceptable and have
a $(u+2,\mu-2)$-refined double character that is either subsequent
  or equipolent to $\vec{\kappa}^{++}$ (refer to \cite{alexakis4} for the strict definition of this notion,
   or to the introduction of this paper for a rough description.). The complete contractions
 indexed in $J$ are simply subsequent to $\vec{\kappa}_{simp}$. The tensor fields indexed in
 $H$ have have a $u$-simple character $\vec{\kappa}_{simp}$ and rank $\ge\mu-1$,
  but they may have one or two unacceptable factor(s)
  $\nabla\Omega_x$ (and furthermore any such factors must be
 contracting against $\nabla\phi_{u+1},\nabla\phi_{u+2}$),
 and possibly one or both of
 the factors $\nabla\phi_{u+1},\nabla\phi_{u+2}$ may be contracting
 against an internal index in some factor $\nabla^{(m)}R_{ijkl}$
 or an index ${}_k,{}_l$ in $S_{*}\nabla^{(\nu)}R_{ijkl}$.

\par Now, by repeating {\it exactly} the same argument (just formally replacing any
 expression $\nabla_i\Omega_h\nabla^i\Omega_{h'}$ by $\nabla_i\Omega_h\nabla^i\phi_{u+1}\nabla_j\Omega_{h'}\nabla^j\phi_{u+2}$)
  that showed that we may ``get rid''
of the sublinear combination $\sum_{h\in H^{\S\S}}\dots$ (modulo
introducing correction terms in the form $\sum_{h\in
H^{\S}}\dots+\sum_{j\in J}\dots$\footnote{That argument,
 did not depend on the
``forbidden cases''.}) in the ``grand conclusion'',\footnote{See the the subsection 
after the ``grand conclusion''.} we may also
assume that all tensor fields indexed in $H$ in
(\ref{reputations2}) have at most one factor $\nabla\Omega_h$.

\par Then, applying Lemma 4.1 in \cite{alexakis4}
(or Lemma 4.2 there if $\sigma=3$) we may assume wlog that all tensor
fields indexed in $H$ in (\ref{reputations2}) have no factors
$\nabla\Omega_h$.\footnote{Since $\mu\ge 4$, by
 weight considerations  there is no danger of  ``forbidden cases''.} 
Under that additional assumption, we may apply
the generalized version of  Lemma 4.10 in \cite{alexakis4} 
if necessary, and additionally assume  that in
(\ref{reputations2}) no tensor fields indexed in $H$ have a factor
$\nabla\phi_{u+1}$ or $\nabla\phi_{u+2}$ contracting against a
special index. We are then in a position to apply Proposition
\ref{giade}  to (\ref{reputations2});
we derive that there exists a linear combination of
acceptable $(\mu-1)$-tensor fields (indexed in $P$ below) with a
$(u+2)$-simple character $Simp(\vec{\kappa}^{++})$ so that:

\begin{equation}
\label{reputations2'}
\begin{split}
&\Sum_{l\in L_{\vec{\kappa}^{++}}} a_l C^{l,i_1\dots
i_\mu}_{g}(\Omega_1,\dots , \Omega_p, \phi_1,\dots
,\phi_u)\nabla_{i_1}\upsilon\dots\nabla_{i_\alpha}\phi_{u+1}
\dots\nabla_{i_\beta}\phi_{u+2}\dots
\nabla_{i_\mu}\upsilon+
\\&\Sum_{p\in P} a_p Xdiv_{i_{\mu+1}}C^{h,i_1\dots
i_{\mu+1}}_{g}(\Omega_1,\dots , \Omega_p, \phi_1,\dots
,\phi_u)\nabla_{i_1}\upsilon\dots\nabla_{i_\alpha}\phi_{u+1}
\dots\nabla_{i_\beta}\phi_{u+2}\\&\dots
\nabla_{i_\mu}\upsilon
+\Sum_{j\in J} a_j C^{j,i_{*}i_{**}}_{g}(\Omega_1,\dots , \Omega_p,
\phi_1,\dots
,\phi_u,\upsilon^{\mu-2})\nabla_{i_{*}}\phi_{u+1}\nabla_{i_{**}}
\phi_{u+2}=0;
\end{split}
\end{equation}
(here the complete contractions indexed in $J$ on the RHS are
simply subsequent to the $u$-simple character
$\vec{\kappa}_{simp}$.

\par Setting $\phi_{u+1}=\phi_{u+2}=\upsilon$ in the above we
derive case B of Lemma \ref{pskovb} in this subcase.

\subsection{Proof of  Lemma \ref{pskovb} in the subcases
$M=\mu\ge 3$ and $M=\mu-1\ge2$.}

\par These are more challenging subcases. We will first consider the case where $M=\mu(\ge 3)$.
Our claim in the case $M=\mu-1\ge 2$ will follow by a simplification of
the same ideas and will be discussed at the end of this proof.

\par The case $M=\mu$ corresponds to the setting where the $\mu$-tensor fields
$C^{l,i_1\dots i_\mu}_{g}$ of maximal refined double character
have {\it all} the free indices ${}_{i_1},\dots {}_{i_\mu}$ belonging to the same factor.
Therefore, the sets $L^z\subset L_\mu,z\in Z_{Max}$ that index the tensor
 fields of maximal refined double character can be easily described in this setting: Let us
list by $F_1,\dots F_c$ all the non-generic factors in
$\vec{\kappa}_{simp}$.\footnote{Recall the ``generic'' factors in
$\vec{\kappa}_{simp}$ are those factors in the form
$\nabla^{(m)}R_{ijkl}$ that are {\it not} contracting against any
factors $\nabla\phi_h$.}
 Let us denote by $L_{c+1}$ the set of all generic factors
on $\vec{\kappa}_{simp}$. Then we may number the maximal refined double
 characters of the $\mu$-tensor fields appearing in (\ref{hypothese2})
 as $\vec{L}^z, z\in \{1,\dots ,c,c+1\}$: The refined double character
$\vec{L}^z, z\le c$ stands for the refined double character that arises from the
simple character $\vec{\kappa}_{simp}$ by assigning $\mu$ free indices to the
 factor $F_z$ (all these free indices must be non-special).\footnote{Recall
that free indices are called ``special'' if they are
 internal indices in some factor $\nabla^{(m)}R_{ijkl}$ or indices ${}_k,{}_l$
 in some factor $S_{*}\nabla^{(\nu)}R_{ijkl}$.} The refined double character
$\vec{L}^{c+1}$ stands for the refined double character that
arises from the simple character $\vec{\kappa}_{simp}$ by
assigning $\mu$ free indices to one of the
 generic factors $\nabla^{(m)}R_{ijkl}$.

\par So, in this case we observe that the index set $L^z, z\in
Z'_{Max}$ corresponds to {\it one} of the index sets
$L^z, c\in \{1,\dots ,c+1\}$ introduced above. So we prove 
Lemma \ref{pskovb} for any chosen index set
$L^z,z\in \{1,\dots ,c+1\}$. So from now on $z$ is a fixed number
from the set $\{1,\dots ,c+1\}$.

\par  Now, we denote by $T_z$ the factor in each $C^{l,i_1\dots
i_\mu}_{g}$, $l\in L^z$  that contains the $M$ free indices
(recall that we have called this factor the {\it critical
factor}). We recall that $\vec{L}^z$ stands for the
$(u,\mu)$-refined double character of the tensor fields
$C^{l,i_1\dots i_\mu}_{g}$, $l\in L^z$.

\par We introduce some notation that will be useful for our proof:

\begin{definition}
\label{fefferm} Consider the refined double characters that is
formally constructed out
 of $\vec{L}^z$ by erasing one of the $M$ free indices from the
  critical factor and then adding a free derivative
index on some other
   factor $S_{*}\nabla^{(\nu)}R_{ijkl}$,
    or $\nabla^{(m)}R_{ijkl}$ or
$\nabla^{(p)}\Omega_x$.
 We call those the {\it linked} refined double
characters. We denote the list of those {\it linked} refined
double characters by $\{\vec{\kappa}_1,\dots
,\vec{\kappa}_\gamma\}$.\footnote{Slightly abusing language we
 will say that $\vec{\kappa}_h$ arises from $\vec{L^z}$ by
 hitting the factor $F_h$ by a derivative $\nabla_{i_{*}}$.}

We then define $L_{\vec{\kappa}_1},\dots ,
L_{\vec{\kappa}_\gamma}\subset L_\mu$ to be the index sets of the
tensor fields $C^{l,i_1\dots i_\mu}$ in (\ref{hypothese2})
with a refined double
character $\vec{\kappa}_1,\dots , \vec{\kappa}_\gamma$
respectively. For notational convenience, we assume that for each
$C^{l,i_1\dots i_\mu}_{g}$ with a linked double character, the
indices ${}_{i_1},\dots ,{}_{i_{\mu-1}}$ belong to the critical factor $T_z$
(and thus the free index ${}_{i_\mu}$ belongs to $F_h$ which we will
 call the second critical factor).
\end{definition}

{\it Note:} Let us consider any $C^{l,i_1\dots i_\mu}_{g}$
  with a refined double character $\vec{\kappa}_h$,
 $1\le h\le \gamma$. Then, by the hypothesis $L_\mu^+=\emptyset$ of our Lemma
\ref{pskovb}, we can assume without
loss of generality that the free index ${}_{i_\mu}$ will
 be a derivative index.

  We now define the second linked refined double
 characters of each $\vec{L^z}$:

 \begin{definition}
\label{fefferm2} Consider the refined double characters that are
 formally constructed out of $\vec{L}^z$ by erasing {\it two}
free indices from the critical factor $T_z$ and adding
 two free indices onto some other factor $\nabla^{(m)}R_{ijkl}$,
 $S_{*}\nabla^{(\nu)}R_{ijkl}$ (necessarily non-simple),
  $\nabla^{(p)}\Omega_h$. These
 (formally constructed) refined double characters will be the
 second linked refined double characters.
\end{definition}
Observe that there is an obvious correspondence between the linked
and the second linked refined double characters of $\vec{L}^z$
based on the non-critical factor  that we hit with one or two free
 derivative indices.   Therefore, we denote the second linked
 refined double characters by
$\vec{\kappa}'_1,\dots ,\vec{\kappa}'_\gamma$, so that each 
$\vec{\kappa}'_h$ corresponds to $\vec{\kappa}_h$.

{\it Technical remark concerning Definitions \ref{fefferm}, \ref{fefferm2}:} If
either $F_h$ or $T_z$ are
 not generic factors $\nabla^{(m)}R_{ijkl}$,  we will then have that
 $\vec{\kappa}_h\ne \vec{\kappa}'_h$ for every $h=1,\dots
 ,\gamma$. If both $F_z$ and $T_z$ are generic factors
 $\nabla^{(m)}R_{ijkl}$ and $M\ge 4$ this is still true. In the case
 where both $F_h, T_z$ are generic factors $\nabla^{(m)}R_{ijkl}$ not
contracting against $\nabla\phi$'s and
 $M=3$ we have that $\vec{\kappa}_h=\vec{\kappa}'_h$.
 Furthermore,
 if $M=4$ and both $F_h,T_z$ are generic we see that in
 $\vec{\kappa}'_h$ there is no well-defined defined critical
 factor. In that case, abusing language, when we refer to the
 critical factor we will in fact be counting
 $C^{l,i_1i_2i_3i_4}_{g}$ twice:  Once with the the factor
to which ${}_{i_1},{}_{i_2}$ belong being the critical factor and once
with the factor to which ${}_{i_3},{}_{i_4}$ belong being the critical factor.
Moreover, if $M=3$ and $F_h,T_z$ are generic then
$\vec{\kappa}_h=\vec{\kappa}'_h$. Therefore, when we speak of the
complete contractions indexed in $L_{\vec{\kappa}_h}$ and
$L_{\vec{\kappa}'_h}$ we are counting the same tensor fields
twice. This, however, will not affect our conclusions
 further down.
\newline

 Now, some more notation: We denote by
$\vec{\kappa}^{+}_z$ the refined double character of the
$(\mu-1)$-tensor fields $C^{l,i_1\dots i_\mu}_{g}
\nabla_{i_1}\phi_{u+1}$, $l\in L^z$. For each $l\in
L_{\vec{\kappa}_h}$, $h=1,\dots ,\gamma$ (which is
 linked to $\vec{L^z}$), we denote by
$\dot{C}^{l,i_1 i_2\dots \hat{i}_\mu}_{g} (\Omega_1,\dots
,\Omega_p,\phi_1,\dots,\phi_u)$ the $(\mu-1)$-tensor field that
arises from $C^{l,i_1\dots i_\mu}_{g}$ by erasing the index
${}_{i_\mu}$. We also denote by $\dot{C}^{l,i_1 i_2\dots
\hat{i}_\mu|i_{*}}_{g} (\Omega_1,\dots
,\Omega_p,\phi_1,\dots,\phi_u)$ the $(\mu-1)$-tensor field that
arises from $\dot{C}^{l,i_1 i_2\dots \hat{i}_\mu}_{g}
(\Omega_1,\dots ,\Omega_p,\phi_1,\dots,\phi_u)$ by adding a
derivative index ${\nabla}_{i_{*}}$ onto the critical factor with $M-1$ free
indices.

\par Now, we apply the grand conclusion to
$L_{g}(\Omega_1,\dots ,\Omega_p,\phi_1,\dots ,\phi_u)$, making
$T_z$ the {\it selected factor}. We derive an equation:

\begin{equation}
\label{thousands}
\begin{split}
& Q_z\cdot \Sum_{l\in L^z} a_l Xdiv_{i_2}\dots
Xdiv_{i_\mu}C^{l,i_1\dots i_\mu}_{g}(\Omega_1,\dots
,\Omega_p,\phi_1,\dots ,\phi_u)\nabla_{i_1}\phi_{u+1}+
\\&\Sum_{y=1}^{\gamma} \overline{2}_y\Sum_{l\in L_{\vec{\kappa}_y}} a_l
Xdiv_{i_1}\dots Xdiv_{i_{\mu-1}}[\dot{C}^{l,i_1 i_2\dots
\hat{i}_\mu,i_{*}}_{g}(\Omega_1,\dots ,\Omega_p,\phi_1,\dots
,\phi_u)\nabla_{i_{*}}\phi_{u+1}]
\\&+\Sum_{l\in L'} a_l Xdiv_{i_1}\dots Xdiv_{i_{\mu-1}}
C^{l,i_1\dots i_\mu}_{g}(\Omega_1,\dots ,\Omega_p, \phi_1,\dots
,\phi_u)\nabla_{i_\mu}\phi_{u+1}+
\\&\Sum_{h\in H} a_h Xdiv_{i_1}\dots
Xdiv_{i_\mu} C^{h,i_1\dots i_\mu,i_{*}}_{g}(\Omega_1,\dots
,\Omega_p,\phi_1,\dots ,\phi_u)\nabla_{i_1}\phi_{u+1}+
\\&\Sum_{f\in F} a_f C^{f,i_{*}}_{g}(\Omega_1,\dots , \Omega_p,
\phi_1,\dots ,\phi_u)\nabla_{i_{*}}\phi_{u+1}=0,
\end{split}
\end{equation}
modulo complete contractions of length $\ge\sigma+u+2$. Here $Q_z$
stands for the first coefficient in the grand conclusion (it
depends on the form of the factor $T_z$) and $|I_1|=\mu$.
Also recall that $\overline{2}_y$ stands for 2 if $F_y$ is
 of the form $\nabla^{(m)}R_{ijkl}$ or $S_{*}\nabla^{(\nu)}R_{ijkl}$ and
 1 if $F_y$ is of the form $\nabla^{(p)}\Omega_h$. The tensor
 fields indexed in $L'$ are generic, acceptable $(\mu-1)$-tensor
fields with a simple character $Simp(\vec{\kappa}^{+}_z)$ but a
refined double character that is either doubly subsequent or equipolent
to $\vec{\kappa}^{+}_z$. The tensor
 fields indexed in $H$ are contributors (see Definition \ref{contributeur}). 

\par Now, by applying Lemma \ref{funny} 
to the above equation,
we may assume with no loss of generality that all the tensor
fields indexed in $H$ are acceptable and have a $(u+1)$-simple
character $\vec{\kappa}_z^{+}$.

 Then, applying Corollary 1 in \cite{alexakis4}\footnote{There
are no special free indices in the tensor fields of minimum rank,
hence there is no danger of falling under ``forbidden cases''.}
 to the above and
 picking out the sublinear combination of terms where all $\nabla\upsilon$'s
 are contracting against the factor $T_z$, we derive that
there is a linear combination of  acceptable $\mu$-tensor
fields (indexed in $P$ below) with $(u+1,\mu-1)$-refined double
character $\vec{\kappa}^{+}_z$ so that:

\begin{equation}
\label{skretas1}
\begin{split}
&Q_z\cdot \Sum_{l\in L^z} a_l C^{l,i_1\dots
i_\mu}_{g}(\Omega_1,\dots ,\Omega_p,\phi_1,\dots
,\phi_u)\nabla_{i_1}\phi_{u+1}
\nabla_{i_2}\upsilon\dots\nabla_{i_\mu}\upsilon+
\\&\Sum_{y=1}^{\gamma} \overline{2}_y\Sum_{l\in L_{\vec{\kappa}_y}}
 a_l [\dot{C}^{l,i_1 i_2\dots
\hat{i}_\mu,i_{*}}_{g}(\Omega_1,\dots ,\Omega_p,\phi_1,\dots
,\phi_u)\nabla_{i_{*}}\phi_{u+1}]\nabla_{i_1}\upsilon\dots\nabla_{i_{\mu-1}}\upsilon
\\&-\Sum_{p\in P} a_p
Xdiv_{i_{\mu+1}} C^{p,i_1\dots i_\mu,i_{\mu+1}}_{g}(\Omega_1,\dots
,\Omega_p,\phi_1,\dots ,\phi_u)\nabla_{i_1}\phi_{u+1}
\nabla_{i_2}\upsilon\dots\nabla_{i_\mu}\upsilon
\\&=\Sum_{j\in J} a_j C^{j,i_1\dots i_\mu}_{g}(\Omega_1,
\dots,\Omega_p,\phi_1,\dots ,\phi_u)\nabla_{i_1}\phi_{u+1}
\nabla_{i_2}\upsilon\dots\nabla_{i_\mu}\upsilon,
\end{split}
\end{equation}
where the tensor fields indexed in $J$ are $u$-subsequent to
$\vec{\kappa}_{simp}$.
\newline

 Since $\mu\ge 3$ we may assume with no loss of generality that $\nabla\phi_{u+1}$
is contracting against a derivative index. By applying
$Erase_{\phi_{u+1}}$  we obtain an equation:

\begin{equation}
\label{elmar}
\begin{split}
&Q_z\cdot \Sum_{l\in L^z} a_l C^{l,i_2\dots
i_\mu}_{g}(\Omega_1,\dots ,\Omega_p,\phi_1,\dots ,\phi_u)
\nabla_{i_2}\upsilon\dots\nabla_{i_\mu}\upsilon+
\\&\Sum_{y=1}^{\gamma} \overline{2}_y\Sum_{l\in L_{\vec{\kappa}_y}}
 a_l [\dot{C}^{l,i_1 i_2\dots
\hat{i}_\mu}_{g}(\Omega_1,\dots ,\Omega_p,\phi_1,\dots
,\phi_u)]\nabla_{i_1}\upsilon\dots\nabla_{i_{\mu-1}}\upsilon
\\&-\Sum_{p\in P} a_p
Xdiv_{i_{\mu+1}} C^{p,i_2\dots i_\mu,i_{\mu+1}}_{g}(\Omega_1,\dots
,\Omega_p,\phi_1,\dots ,\phi_u)
\nabla_{i_2}\upsilon\dots\nabla_{i_\mu}\upsilon=
\\&\Sum_{j\in J} a_j C^{j,i_2\dots i_\mu}_{g}(\Omega_1,
\dots,\Omega_p,\phi_1,\dots ,\phi_u)
\nabla_{i_2}\upsilon\dots\nabla_{i_\mu}\upsilon,
\end{split}
\end{equation}
where the tensor fields indexed in $J$ are $u$-subsequent to
$\vec{\kappa}_{simp}$.
\newline

{\it Remark 1:} Thus the above equation involves the $\mu$-tensor fields indexed in 
$L^z$, which our Lemma \ref{pskovb} deals with in this subcase, but it {\it also}
involves the $\mu$-tensor fields in the second 
line. We can therefore {\it not} derive our
 Lemma \ref{pskovb} in this subcase from the above equation alone. 
We seek to derive two more equations in order to obtain a closed system of 
three equations in three different sublinear combinations. 
 In order to formulate and derive our next two
 equations we will need to introduce some more notation:
\newline

{\bf Notation:} For each $h, 1\le h\le c+1$ we denote by
$\dot{C}^{l,\hat{i}_1i_2\dots i_\mu, i_{*}\rightarrow F_h}_{g}$
 the
tensor field that formally  arises from $C^{l,i_1\dots i_\mu}_{g}$, $l\in
L^z$ by erasing a derivative index ${}_{i_\mu}$ and hitting the
(one of the) factor(s) $F_h$ by a derivative $\nabla_{i_{*}}$ (and
adding over all tensor fields we thus obtain if $h=c+1$). Also,
for each $l\in L_{\vec{\kappa}'_v}$, $v=1\dots ,\gamma$, we denote
by $C^{l,i_1\dots i_{\mu-1} \hat{i}_\mu|i_{*}}_{g}$ the tensor
field that arises from $C^{l,i_1\dots i_\mu}_{g}$ by
 erasing the free index ${}_{i_\mu}$,\footnote{Recall that we are assuming ${}_{i_{\mu-1}},{}_{i_\mu}$
 belong to the second critical factor.} (we may always assume it
 is a derivative index) and adding a derivative index ${\nabla}_{i_{*}}$
 onto the critical factor $T_z$. So these tensor fields have $M-1$
  free indices in $T_z$.

\par Now, we apply the grand conclusion to $L_{g}$ making
$F_h$ the selected factor. In order to describe the equation we
then obtain, we introduce some notation: Consider $(\mu-1)$-tensor
fields with a factor $\nabla\phi_{u+1}$ contracting against
 a non-special index in $F_h$,\footnote{Recall that a special index
 is an internal index in some factor $\nabla^{(m)}R_{ijkl}$ or an index
  ${}_k,{}_l$ in some factor $S_{*}\nabla^{(\nu)}R_{ijkl}$.}  while all the other $\mu-1$ free
 indices belong to the factor $T_z$. We denote by $\vec{\kappa}^{+}_h$
 the refined $(u+1,\mu-1)$-double character that corresponds to these tensor fields.
  In this setting we will
denote by $\Sum_{l\in L'} a_l C^{l,i_1\dots i_\mu}_{g}
\nabla_{i_\mu}\phi_{u+1}$ a generic linear combination of
 tensor fields with a $(u+1)$-simple character
$Simp(\vec{\kappa}^{+}_h)$ but a $(u+1,\mu-1)$-refined double character
 that is either doubly subsequent or equipolent to $\vec{\kappa}^{+}_h$.
We derive:

\begin{equation}
\label{podarakia}
\begin{split}
&(\overline{2}_zM-{{M}\choose{2}})\Sum_{l\in L^z} a_l
Xdiv_{i_2}\dots Xdiv_{i_\mu} \dot{C}^{l,\hat{i}_1i_2\dots
i_\mu,i_{*}\rightarrow F_h}_{g}(\Omega_1,\dots
,\Omega_p,\phi_1,\dots ,\phi_u)\\&\nabla_{i_{*}}\phi_{u+1}
+Q_h\cdot \Sum_{l\in
L_{\vec{\kappa}_h}} a_l Xdiv_{i_1}\dots Xdiv_{i_{\mu
-1}}C^{l,i_1\dots i_\mu}_{g}(\Omega_1,\dots ,\Omega_p,\phi_1,\dots
,\phi_u)\nabla_{i_\mu}\phi_{u+1} 
\\&+\Sum_{l\in L_{\vec{\kappa}'_h}}
a_l Xdiv_{i_1}\dots Xdiv_{i_{\mu -2}}Xdiv_{i_{*}}
 C^{l,i_1\dots i_{\mu-1}\hat{i}_\mu|
i_{*}}_{g}(\Omega_1,\dots ,\Omega_p,\phi_1,\dots
,\phi_u)\\&\nabla_{i_{\mu-1}}\phi_{u+1}+
\\&\Sum_{l\in L'} a_l Xdiv_{i_1}\dots Xdiv_{i_{\mu-1}}
C^{l,i_1\dots i_\mu}_{g}(\Omega_1,\dots ,\Omega_p, \phi_1,\dots
,\phi_u)\nabla_{i_\mu}\phi_{u+1}
\\&+ \Sum_{h\in H} a_h Xdiv_{i_1}\dots Xdiv_{i_a}
C^{h,i_1 \dots i_a}_{g}(\Omega_1, \dots ,\Omega_p, \phi_1,\dots
,\phi_{u+1})=
\\&\Sum_{j\in J} a_j C^j_{g}(\Omega_1,
\dots ,\Omega_p, \phi_1,\dots ,\phi_{u+1}),
\end{split}
\end{equation}
where $Q_h$ stands for  the coefficient in the grand conclusion
with selected factor $F_h$ and $|I_1|=1$. Recall $\overline{2}_z$
stands for $2$ if the factor $T_z$ is of the form
$\nabla^{(m)}R_{ijkl}$ or $S_{*}\nabla^{(\nu)}R_{ijkl}$. If $T_z$
is of the form $\nabla^{(p)}\Omega_j$ then it stands for 1 if
$T_z$ is not contracting against a factor $\nabla\phi_h$ and 0
otherwise. The complete contractions in $J$ are $u$-simply
subsequent to $\vec{\kappa}_{simp}$.
\newline

{\it Remark 2:}  Now, we observe that the equation above involves the 
sublinear combinations indexed in $L^z$, $L_{\vec{\kappa}_h}$, but it also involves 
the sublinear combination indexed in $L_{\vec{\kappa}'_h}$. Thus, we seek
 to derive a third equation in order to obtain a closed system. 
\newline

\par Therefore, we invoke the equation
$Im^{1,\beta}_{\phi_{u+1}}[L_{g}]=0$. i.e.~equation (\ref{antrikos}). 
We then define an operation
$SimpOp$ that acts on the terms in
$Im^{1,\beta}_{\phi_{u+1}}[L_{g}]$ by replacing the factor
$\nabla\omega$ by an $Xdiv$. It follows that
$SimpOp\{Im^{1,\beta}_{\phi_{u+1}}[L_{g}]\}=0$ (by virtue of the 
last Lemma  in the Appendix of \cite{alexakis1}).

We will again pick any $h\le c+1$ and focus on the
terms in $Im^{1,\beta}_{\phi_{u+1}}[L_g]=0$ with the 
 factor $\nabla\phi_{u+1}$ contracting against the factor
$F_h$, and $\nabla\omega$ contracting against $T_z$. 
 (So we again focus on the
 $(u+1,\mu-1)$-refined double character $\vec{\kappa}^{+}_h$
 as before). 
Picking out this sublinear combination (which vanishes separately), 
and acting on it by $SimpOp[\dots]$, we derive:

\begin{equation}
\label{podarakia2}
\begin{split}
&{{M}\choose{2}}\Sum_{l\in L^z} a_l Xdiv_{i_2}\dots
 Xdiv_{i_\mu}
\dot{C}^{l,\hat{i}_1i_2\dots i_\mu,i_{*}}_{g}(\Omega_1,\dots
,\Omega_p,\phi_1,\dots ,\phi_u)\nabla_{i_{*}}\phi_{u+1}+
\\&(M-1)\Sum_{l\in L_{\vec{\kappa}_h}} a_l Xdiv_{i_1}\dots
Xdiv_{i_{\mu -1}}C^{l,i_1\dots i_\mu}_{g}(\Omega_1,\dots
,\Omega_p,\phi_1,\dots ,\phi_u)\nabla_{i_\mu}\phi_{u+1}+
\\&\Sum_{l\in L_{\vec{\kappa}'_h}} a_l
 Xdiv_{i_1}\dots Xdiv_{i_{\mu -1}}
 \hat{Xdiv}_{i_\mu}Xdiv_{i_{*}}
C^{l,i_1\dots i_{\mu-1}|i_{*}}_{g}(\Omega_1,\dots ,\Omega_p,
\phi_1,\dots ,\phi_u)\\&\nabla_{i_{\mu-1}}\phi_{u+1}+
\Sum_{l\in L'} a_l Xdiv_{i_1}\dots Xdiv_{i_{\mu-1}}
C^{l,i_1\dots i_\mu}_{g}(\Omega_1,\dots ,\Omega_p, \phi_1, \dots
,\phi_u)\nabla_{i_\mu}\phi_{u+1}+
\\& \Sum_{h\in H} a_h Xdiv_{i_1}\dots Xdiv_{i_a}
C^{h,i_1 \dots i_a}_{g}(\Omega_1, \dots ,\Omega_p, \phi_1,
\dots,\phi_{u+1})=
\\&\Sum_{j\in J} a_j C^j_{g}(\Omega_1,
\dots ,\Omega_p, \phi_1,\dots ,\phi_{u+1});
\end{split}
\end{equation}
(the sublinear combinations in $L',H,J$ stand for {\it generic}
linear combinations as in (\ref{podarakia})).

\par Now, we may first apply Lemma \ref{funny} (to ensure
 all tensor fields in $H$ are acceptable and have 
the same $(u+1)$-simple character $\vec{\kappa}^+_h$), 
and then Corollary 1 in \cite{alexakis4}\footnote{There will be $\mu-1\ge 2$ non-special
free indices among the tensor fields of minimum rank, hence no
danger of falling under a ``forbidden case''.} to the two
equations (\ref{podarakia}), (\ref{podarakia2}), and pick out the sublinear combinations where there are
$(\mu-1)$ factors $\nabla\upsilon$ contracting against $T_z$ and
$\nabla\phi_{u+1}$ is contracting against $F_h$. These sublinear
combinations will vanish separately, and we thus derive two
new equations:

\begin{equation}
\label{podarakiacor}
\begin{split}
&(\overline{2}_zM-{{M}\choose{2}})\Sum_{l\in L^z} a_l
 \dot{C}^{l,\hat{i}_1i_2\dots
i_\mu,i_{*}\rightarrow F_h}_{g}(\Omega_1,\dots
,\Omega_p,\phi_1,\dots
,\phi_u)\nabla_{i_{*}}\phi_{u+1}\nabla_{i_2}\upsilon\dots\\&\nabla_{i_\mu}\upsilon+
Q_h\cdot \Sum_{l\in
L_{\vec{\kappa}_h}} a_l C^{l,i_1\dots i_\mu}_{g}(\Omega_1,\dots
,\Omega_p,\phi_1,\dots
,\phi_u)\nabla_{i_\mu}\phi_{u+1}\nabla_{i_1}
\upsilon\dots\nabla_{i_{\mu-1}}\upsilon+
\\&\Sum_{l\in L_{\vec{\kappa}'_h}}
a_l  C^{l,i_1\dots i_{\mu-1}\hat{i}_\mu| i_{*}}_{g}(\Omega_1,\dots
,\Omega_p,\phi_1,\dots
,\phi_u)\nabla_{i_{\mu-1}}\phi_{u+1}\nabla_{i_1}
\upsilon\dots\nabla_{i_{\mu-2}}\upsilon\nabla_{i_{*}}\upsilon+
\\& \Sum_{h\in H_1} a_h Xdiv_{i_{\mu+1}}
C^{h,i_1 \dots i_{\mu+1}}_{g}(\Omega_1, \dots ,\Omega_p,
\phi_1,\dots ,\phi_u)\nabla_{i_1}\phi_{u+1}
\nabla_{i_2}\upsilon\dots\nabla_{i_\mu}\upsilon=
\\&\Sum_{j\in J} a_j C^j_{g}(\Omega_1,
\dots ,\Omega_p, \phi_1,\dots ,\phi_{u+1},\upsilon^{\mu-1}),
\end{split}
\end{equation}

\begin{equation}
\label{podarakia2cor}
\begin{split}
&{{M}\choose{2}}\Sum_{l\in L^z} a_l \dot{C}^{l,\hat{i}_1i_2\dots
i_\mu,i_{*}}_{g}(\Omega_1,\dots ,\Omega_p,\phi_1,\dots
,\phi_u)\nabla_{i_{*}}\phi_{u+1}\nabla_{i_2}\upsilon\dots\nabla_{i_\mu}\upsilon+
\\&(M-1)\Sum_{l\in L_{\vec{\kappa}_h}} a_l C^{l,i_1\dots i_\mu}_{g}(\Omega_1,\dots
,\Omega_p,\phi_1,\dots
,\phi_u)\nabla_{i_\mu}\phi_{u+1}\nabla_{i_1}
\upsilon\dots\nabla_{i_{\mu-1}}\upsilon+
\\&\Sum_{l\in L_{\vec{\kappa}'_h}} a_l
C^{l,i_1\dots i_{\mu-1}|i_{*}}_{g}(\Omega_1,\dots ,\Omega_p,
\phi_1,\dots ,\phi_u)\nabla_{i_{\mu-1}}\phi_{u+1}\nabla_{i_1}
\upsilon\dots\nabla_{i_{\mu-1}}\upsilon+
\\& \Sum_{h\in H_2} a_h Xdiv_{i_{\mu+1}}
C^{h,i_1 \dots i_{\mu+1}}_{g}(\Omega_1, \dots ,\Omega_p,
\phi_1,\dots ,\phi_u)\nabla_{i_1}\phi_{u+1}
\nabla_{i_2}\upsilon\dots\nabla_{i_\mu}\upsilon=
\\&\Sum_{j\in J} a_j C^j_{g}(\Omega_1,
\dots ,\Omega_p, \phi_1,\dots ,\phi_{u+1},\upsilon^{\mu-1}).
\end{split}
\end{equation}

\par In the above two equations, all the tensor fields in the
first three lines are acceptable, and have a $(u+1,\mu-1)$-refined
double character $\vec{\kappa}^{+}_h$. The tensor fields indexed
in $H_1$ are acceptable contributors with a $(u+1)$-simple character $\vec{\kappa}^+_h$ 
(see Definition \ref{contributeur}). The complete contractions in $J$ are in
both cases simply subsequent to $Simp(\vec{\kappa}^{+}_h)$, and
hence also to $\vec{\kappa}_{simp}$.

\par Subtracting (\ref{podarakia2cor}) from (\ref{podarakiacor})
we derive a new equation:

\begin{equation}
\label{podarakiacor''}
\begin{split}
&(\overline{2}_zM-2{{M}\choose{2}})\Sum_{l\in L^z} a_l
 \dot{C}^{l,\hat{i}_1i_2\dots
i_\mu,i_{*}\rightarrow F_h}_{g}(\Omega_1,\dots
,\Omega_p,\phi_1,\dots
,\phi_u)\nabla_{i_{*}}\phi_{u+1}\nabla_{i_2}\upsilon\dots\\&\nabla_{i_\mu}\upsilon+
(Q_h-(M-1))\cdot \Sum_{l\in
L_{\vec{\kappa}_h}} a_l C^{l,i_1\dots i_\mu}_{g}(\Omega_1,\dots
,\Omega_p,\phi_1,\dots
,\phi_u)\nabla_{i_\mu}\phi_{u+1}\nabla_{i_1}
\upsilon\dots\\&\nabla_{i_{\mu-1}}\upsilon+
 \Sum_{h\in H'} a_h Xdiv_{i_{\mu+1}}
C^{h,i_1 \dots i_{\mu+1}}_{g}(\Omega_1, \dots ,\Omega_p,
\phi_1,\dots ,\phi_u)\nabla_{i_1}\phi_{u+1}
\nabla_{i_2}\upsilon\dots\nabla_{i_\mu}\upsilon
\\&=\Sum_{j\in J} a_j C^j_{g}(\Omega_1,
\dots ,\Omega_p, \phi_1,\dots ,\phi_{u+1},\upsilon^{\mu-1}),
\end{split}
\end{equation}
where the tensor fields and complete contractions indexed in
$H',J$ have the same general properties as the ones indexed in
$H_1,J$ in (\ref{podarakiacor}).

\par We act on (\ref{podarakiacor''})
 by the operation $Erase_{\phi_{u+1}}$ (this operation is
 formally well-defined and produces acceptable tensor fields) and divide by $Q_h-(M-1)$. We
 thus derive an equation:

\begin{equation}
\label{skretas5}
\begin{split}
&\Sum_{l\in L_{\vec{\kappa}_h}} a_l C^{l,i_1\dots i_{\mu-1}}_{g}
(\Omega_1,\dots ,\Omega_p,\phi_1,\dots ,\phi_u)
\nabla_{i_1}\upsilon\dots\nabla_{i_{\mu-1}}\upsilon =
\\&\frac{(2{{M}\choose{2}}-2_zM)}{Q_h-(M-1)}
\Sum_{l\in L^z} a_l \dot{C}^{l,\hat{i}_1i_2\dots
i_\mu}_{g}(\Omega_1,\dots , \Omega_p,
\phi_1,\dots,\phi_u)\nabla_{i_2}
\upsilon\dots\nabla_{i_\mu}\upsilon +
\\& \Sum_{p\in P} a_p  Xdiv_{i_{\mu}}
C^{p,i_1 \dots i_{\mu}}_{g}(\Omega_1, \dots ,\Omega_p,
\phi_1,\dots,\phi_{u+1})\nabla_{i_1}\upsilon\dots
\nabla_{i_{\mu-1}}\upsilon+
\\&\Sum_{j\in J} a_j C^{j,i_1\dots i_{\mu-1}}_{g}
(\Omega_1,\dots ,\Omega_p, \phi_1,\dots ,\phi_{u+1})
\nabla_{i_1}\upsilon\dots \nabla_{i_{\mu-1}}\upsilon.
\end{split}
\end{equation}

Notice the coefficient
$\frac{(2{{M}\choose{2}}-2_zM)}{Q_h-(M-1)}$
is {\it non-positive} (because $M\ge 3$ and $Q_h<M-1$ by
inspection). We denote
 this coefficient by $c_h$.
\newline

\par Now, replacing the above into (\ref{elmar}), and since
$(Q_z+\Sum_{h=1}^\gamma 2_h\cdot c_h)<0$, we obtain:

\begin{equation}
\label{elmar2}
\begin{split}
& \Sum_{l\in L^z} a_l C^{l,i_2\dots
i_\mu}_{g}(\Omega_1,\dots ,\Omega_p,\phi_1,\dots ,\phi_u)
\nabla_{i_2}\upsilon\dots\nabla_{i_\mu}\upsilon
\\&-\Sum_{p\in P} a_p
Xdiv_{i_{\mu+1}} C^{p,i_2\dots i_\mu,i_{\mu+1}}_{g}(\Omega_1,\dots
,\Omega_p,\phi_1,\dots ,\phi_u)
\nabla_{i_2}\upsilon\dots\nabla_{i_\mu}\upsilon=
\\&\Sum_{j\in J} a_j C^{j,i_2\dots i_\mu}_{g}(\Omega_1,
\dots,\Omega_p,\phi_1,\dots ,\phi_u)
\nabla_{i_2}\upsilon\dots\nabla_{i_\mu}\upsilon,
\end{split}
\end{equation}
modulo complete contractions of length $\ge\sigma+u$.

 Now, define an operation $Op$ that acts on the above by adding
  a derivative $\nabla_i$ onto the factor $T_z$
which is contracting
 against the $\mu-1$ factors $\nabla\upsilon$, 
snd contracting $\nabla_i$ against a factor $\nabla\upsilon$. Since (\ref{elmar2}) holds
   formally, $Op$ produces a true equation.
Applying Lemma \ref{funny} to this equation, we derive
 our desired conclusion for the case $M=\mu\ge 3$.
\newline

\par We now prove case B of Lemma \ref{pskovb} in the 
case, $M=\mu-1\ge 2$.
 Slightly abusing notation, we carry over
the notation from the previous discussion, only now
$L^z=\emptyset$. In this convention, the $\mu$-tensor fields 
 in (\ref{hypothese2}) of maximal refined double character
will correspond to the index sets $L_{\vec{\kappa}'_h}$,
$h=1,\dots ,\gamma$.

 We can again make use of equations (\ref{podarakia})
and (\ref{podarakia2}), where we now have that $L^z=\emptyset$. We
then subtract these two equations, obtaining a new true equation.
 Since
$Q_h-(M-1)<0$ as before, we may divide by that constant and set
$\phi_{u+1}=\upsilon$. That gives us our
 desired conclusion. $\Box$

\subsection{Proof of Lemma \ref{pskovb}  in the subcase $M=\mu=2$.}

\par As before, the two main ingredients from the earlier
discussion that we will be using here are the grand conclusion,
and also the equation (\ref{antrikos}),
$Im^{1,\beta}_{\phi_{u+1}}[L_{g}]=0$.
\newline

\par This case presents certain particular difficulties.
 We first divide $L_\mu=L_2$ into two subsets: We say
$l\in L^{1,1}$ if the two free indices ${}_{i_1},{}_{i_2}$ belong to
 different factors.  We say $l\in L^{2,0}$ if the two free
indices ${}_{i_1},{}_{i_2}$ belong to the same factor. We recall
our Lemma hypothesis $L_\mu^{*}=\emptyset$ which ensures that for
each $l\in L^{2,0}$ we cannot have ${}_{i_1},{}_{i_2}$ belonging
to a factor $\nabla^{(2)}\Omega_h$.

\par Firstly let us understand what the refined double characters
associated to $\vec{\kappa}_{simp}$ are. This is easy to do, in
this case: We denote by $F_1,\dots ,F_\tau$ the list of factors in
$\vec{\kappa}_{simp}$ which are contracting against a factor
$\nabla\phi_h$ or are of the form $\nabla^{(p)}\Omega_y$. The rest of the factors $F_l$ in
$\vec{\kappa}_{simp}$ that do not belong to this list will be
generic factors of the form $\nabla^{(m)}R_{ijkl}$.\footnote{Recall
that a generic factor is a factor of the form
$\nabla^{(m)}R_{ijkl}$ that is not contracting against any
$\nabla\phi_h$.} Now, for each $h\le\tau$ we may unambiguously
speak of {\it the} factor $F_h$.

\par For each $h\le \tau$, we denote by $L^{2,0|h}\subset L^{2,0}$ the index set of
the 2-tensor fields $C^{l,i_1i_2}_{g}$ with the two free indices
${}_{i_1},{}_{i_2}$ belonging to the factor $F_h$. On the other hand, we
index in $L^{2,0|\tau+1}\subset L^{2,0}$ all the 2-tensor fields
which have the two indices ${}_{i_1},{}_{i_2}$ belonging to a generic factor
$\nabla^{(m)}R_{ijkl}$. We observe that:

$$\Sum_{l\in L^{2,0}}
a_l C^{l,i_1i_2}_{g}=\Sum_{h=1}^{\tau+1}\Sum_{l\in L^{2,0|h}}a_l
C^{l,i_1i_2}_{g}.$$ Notice that (in this subcase) 
the 2-tensor fields of {\it
maximal} refined double character (in (\ref{hypothese2})) are the ones indexed in
$L^{2,0}$. Observe that the index set $\bigcup_{z\in Z'_{Max}}
L^z$ will correspond to {\it one} index set $L^{2,0|h}$. We denote
the $h$ that corresponds to this {\it one} index set by $\alpha$.
We recall that the factor $F_\alpha$ is called the {\it critical
factor}.

\par Now, for notational convenience we will be assuming that the
free index ${}_{i_2}$ in each $C^{l,i_1i_2}_{g}, l\in L^{2,0}$ is a
derivative index. This can be done with no loss of generality by
virtue of the assumptions  
 of Lemma \ref{pskovb}  (no special free
indices in any $C^{l,i_1i_2}_{g}$ and also $L_\mu^+=\emptyset$).

\par For any $a,b\le \tau$ we now denote by $L^{1,1|a,b}$ the index set of the 2-tensor
fields indexed in $L^{1,1}$ with the additional feature that one
free index belongs to $F_a$ and the other to $F_b$. On the other
hand for any $a\le \tau$, we denote by $L^{1,1|a,\tau+1}$ the
index set of the 2-tensor fields where one index belongs to the
factor $F_a$ and the other belongs to a generic  factor
$\nabla^{(m)}R_{ijkl}$. Finally, we denote by
$L^{1,1|\tau+1,\tau+1}$ the index of the 2-tensor fields where
both free indices belong to (different) generic factors
$\nabla^{(m)}R_{ijkl}$. Recall that we may assume wlog that if
$F_\alpha$ or $F_b$ are simple factors of the form $S_{*}\nabla^{(\rho)}R_{ijkl}$,
then $L^{1,1|\alpha,b}=\emptyset$.

\par  We distinguish two subcases:
Either the critical factor $F_\alpha$ is in one of the forms $\nabla^{(A)}\Omega_h$,
$S_{*}\nabla^{(\nu)} R_{ijkl}$, or it is in the form
$\nabla^{(m)}R_{ijkl}$. We first consider the first subcase.
\newline

{\it Proof of case $M=\mu=2$ when the critical factor is in one of
the forms $\nabla^{(A)}\Omega_h$, $S_{*}\nabla^{(\nu)} R_{ijkl}$:}

 We will firstly  show that for each pair $(\alpha,b)$,
 (where $b\ne \alpha$ if $\alpha\le \tau+1$ and neither
 $F_\alpha,F_b$ is a simple factor of the from
 $S_{*}\nabla^{(\nu)}R_{ijkl}$):

\begin{equation}
\label{kubernhseis}
\begin{split}
&\Sum_{l\in L^{1,1|\alpha,b}} a_l C^{l,i_1i_2}_{g}(\Omega_1,\dots
,\Omega_p,\phi_1,\dots
,\phi_u)\nabla_{i_1}\upsilon\nabla_{i_2}\upsilon+
\\& 2\Sum_{l\in
L^{2,0|\alpha}}a_l \dot{C}^{l,i_1|i_{*}\rightarrow b}_{g}
(\Omega_1,\dots ,\Omega_p,\phi_1,\dots
,\phi_u)\nabla_{i_1}\upsilon\nabla_{i_{*}}\upsilon =
\\& Xdiv_{i_3}\Sum_{h\in H} a_h C^{h,i_1i_2i_3}_{g}(\Omega_1,\dots ,\Omega_p,\phi_1,\dots
,\phi_u)\nabla_{i_1}\upsilon\nabla_{i_2}\upsilon+
\\&\Sum_{j\in J}
a_j C^{j,i_1i_2}_{g}(\Omega_1,\dots ,\Omega_p,\phi_1,\dots
,\phi_u)\nabla_{i_1}\upsilon\nabla_{i_2}\upsilon;
\end{split}
\end{equation}
here $\dot{C}^{l,i_1|i_{*}\rightarrow b}_{g}$ stands for the
tensor field that arises from $C^{l,i_1i_2}_{g}$ by erasing
 the (derivative) index ${}_{i_2}$ and adding a free derivative
 index onto $F_b$ (and $S_{*}$-symmetrizing if needed).
$\sum_{h\in H}\dots$ stands for a generic linear combination of
acceptable
 3-tensor fields with a $u$-simple character $\vec{\kappa}_{simp}$.
 $\sum_{j\in J}\dots$ stands for  a linear combination of contractions
 that are simply subsequent to $\vec{\kappa}_{simp}$.

{\it Note regarding the notation in (\ref{kubernhseis}):}
 For each $b\le\tau$ the tensor field $\dot{C}^{l,i_1|i_{*}\rightarrow
b}_{g} (\Omega_1,\dots ,\Omega_p,\phi_1,\dots ,\phi_u)$ is
well-defined. On the other hand, if $b=\tau+1$ 
$\dot{C}^{l,i_1|i_{*}\rightarrow \tau+1}_{g} (\Omega_1,\dots
,\Omega_p,\phi_1,\dots ,\phi_u)$ will be the sublinear
combination in $\nabla_{i_{*}}[\dot{C}^{l,i_1}_{g} (\Omega_1,\dots
,\Omega_p,\phi_1,\dots ,\phi_u)]$ where $\nabla_{i_{*}}$ is only
allowed to hit one of the factors $\nabla^{(m)}R_{ijkl}$ which are
not contracting against a factor $\nabla\phi_h$. Furthermore, if
$\alpha=\tau+1$ we additionally require that $\nabla_{i_{*}}$ can
not hit the factor $\nabla^{(m)}R_{ijkl}$ to which ${}_{i_1}$
belongs.

\par We will show (\ref{kubernhseis}) below. For now, we 
show how (\ref{kubernhseis}) implies our claim. 
\newline

{\it Proof that (\ref{kubernhseis}) implies our claim (when the
critical factor is in one of the forms $\nabla^{(A)}\Omega_h$,
$S_{*}\nabla^{(\nu)} R_{ijkl}$):} We make note of a
straightforward corollary of (\ref{kubernhseis}): Making the
$\nabla\upsilon$'s into $Xdiv$'s (using the last Lemma in the Appendix of \cite{alexakis1})
 and substituting this into our Lemma hypothesis, $L_{g}=0$ ,
we derive:

\begin{equation}
\label{paparares}
\begin{split}
&\Sum_{l\in L^{2,0|\alpha}} a_l
Xdiv_{i_1}Xdiv_{i_2}C^{l,i_1i_2}_{g}(\Omega_1,\dots
,\Omega_p,\phi_1,\dots ,\phi_u)+
\\& \Sum_{f\ne \alpha}\Sum_{l\in L^{2,0|f}} a_l
Xdiv_{i_1}Xdiv_{i_2}C^{l,i_1i_2}_{g}(\Omega_1,\dots
,\Omega_p,\phi_1,\dots ,\phi_u)+
\\&\Sum_{l\in
L^{2,0|a}}a_l Xdiv_{i_1}Xdiv_{i_{*}}\Sum_{b\le \tau+1,b\ne
\alpha}\dot{C}^{l,i_1|i_{*}\rightarrow b}_{g} (\Omega_1,\dots
,\Omega_p,\phi_1,\dots ,\phi_u)+
\\&\Sum_{1\le c<d\le\tau+1,c,d\ne \alpha}
\Sum_{l\in L^{1,1|c,d}} a_l Xdiv_{i_1}Xdiv_{i_2}
C^{l,i_1i_2}_{g} (\Omega_1,\dots ,\Omega_p,\phi_1,\dots ,\phi_u)+
\\&\Sum_{h\in H} a_h Xdiv_{i_1}\dots
Xdiv_{i_t}C^{h,i_1\dots i_t}_{g}(\Omega_1,\dots
,\Omega_p,\phi_1,\dots ,\phi_u)\\&+\Sum_{j\in J} a_j
C^j_{g}(\Omega_1,\dots ,\Omega_p,\phi_1,\dots ,\phi_u)=0.
\end{split}
\end{equation}
In particular, we have succeeded in replacing the old sublinear
combinations $\Sum_{l\in L^{1,1|\alpha,c}}\dots$of 2-tensor fields
 with the two free indices belonging to different factors with new
 ones, which in particular arise from
$\Sum_{\in L^{2,0|\alpha}}$ in the precise way outlined above.

\par Now, we will be applying the grand conclusion to
(\ref{paparares}), making $F_\alpha$ the selected factor.
We will be interested in the 2-tensor fields in the grand
conclusion that will have the factor $\nabla\phi_{u+1}$ and the
free index ${}_{i_1}$ contracting against/ belonging to the crucial
factor $F_\alpha$.

\par  Now, in order to apply the grand conclusion we
must check that the extra claims are satisfied in this setting. We
indeed have that $L_2^{*}=\emptyset$ by hypothesis and also 
 that $L_2^{+}=\emptyset$, by inspection in (\ref{paparares}). 
Then, we apply
 Lemma \ref{funny} to obtain an equation
that is the same as (\ref{paparares}) only with the additional
restriction that the index set $H$ is replaced by an index set
$H'$ which indexes terms which are contributors 
which are acceptable and have $u$-simple character 
$\vec{\kappa}_{simp}$ and $\nabla\phi_{u+1}$ not contracting against 
a special index. 
 In view of this, we may now apply the grand conclusion to
(\ref{paparares}) and we obtain an equation:

\begin{equation}
\label{umnhswmen}
\begin{split}
&(2q_a+2\sigma_{a,*})\Sum_{l\in L^{2,0|a}} a_l
Xdiv_{i_2}C^{l,i_1i_2}_{g}(\Omega_1,\dots ,\Omega_p,\phi_1,\dots
,\phi_u)\nabla_{i_1}\phi_{u+1}+
\\&\Sum_{l\in \tilde{L}} a_l
Xdiv_{i_2}C^{l,i_1i_2}_{g}(\Omega_1,\dots ,\Omega_p,\phi_1,\dots
,\phi_u)\nabla_{i_1}\phi_{u+1}+
\\&\Sum_{h\in H} a_h Xdiv_{i_2}\dots
Xdiv_{i_t}C^{h,i_1\dots i_t}_{g}(\Omega_1,\dots
,\Omega_p,\phi_1,\dots ,\phi_u)\nabla_{i_1}\phi_{u+1}+
\\&\Sum_{j\in
J} a_j C^j_{g}(\Omega_1,\dots ,\Omega_p,\phi_1,\dots
,\phi_{u+1})=0.
\end{split}
\end{equation}
Here $q_\alpha$ stands for the coefficient between parentheses in the
grand conclusion (see the discussion {\it after} the ``grand 
conclusion''--recall that the coefficient $q_\alpha$ depends on the {\it form} of the selected factor 
$F_\alpha$) with $|I_1|=2$ for the factor $F_\alpha$, while
$\sigma_{\alpha,*}=-\Sum_{v=1,v\ne \alpha}^\sigma 2_v$ (recall the
definition of $2_v$ from the discussion above the ``grand conclusion''). Here all the
1-tensor fields are acceptable, and also have the factor
$\nabla\phi_{u+1}$ contracting against the factor $F_\alpha$.
Moreover, the 2-tensor fields indexed in $\tilde{L}$ have the free
index ${}_{i_2}$ {\it not} belonging to the selected factor $F_\alpha$.
(It will in fact be a non-dangerous index in some factor other than
$F_\alpha$). Terms indexed in $H$ are contributors. Now, applying Lemma \ref{funny} if
necessary,\footnote{There is no danger of falling under a
``forbidden case'', by inspection.} we may also assume that
$\nabla\phi_{u+1}$ is not contracting against a special index in
$F_\alpha$ for any of the tensor fields indexed in $H$.\footnote{Recall that a special index
is an internal index in some factor $\nabla^{(m)}R_{ijkl}$ or one
of the indices ${}_k,{}_l$ in some factor
$S_{*}\nabla^{(\nu)}R_{ijkl}$.}

Observe that $2(q_\alpha+\sigma_{\alpha,*})< 0$ ($q_\alpha<0$
because $|I_1|=2$).

\par But then, just applying Corollary 1 in 
\cite{alexakis4}\footnote{Notice that for the tensor fields of minimum rank
there is a non-special free index.} to the above (and picking out
the sublinear combination where $\nabla\phi_{u+1},\nabla\upsilon$
are contracting against the same factor), we derive an equation:

\begin{equation}
\label{gotnorythm}
\begin{split}
&(2q_\alpha+2\sigma_{\alpha,*})\Sum_{l\in L^{2,0|\alpha}} a_l
C^{l,i_1i_2}_{g}(\Omega_1,\dots ,\Omega_p,\phi_1,\dots
,\phi_u)\nabla_{i_1}\phi_{u+1}\nabla_{i_2}\upsilon+
\\&\Sum_{h\in H} a_h Xdiv_{i_3}\dots
Xdiv_{i_t}C^{h,i_1\dots i_t}_{g}(\Omega_1,\dots
,\Omega_p,\phi_1,\dots
,\phi_u)\nabla_{i_1}\phi_{u+1}\nabla_{i_2}\upsilon+
\\&\Sum_{j\in J}
a_j C^j_{g}(\Omega_1,\dots ,\Omega_p,\phi_1,\dots
,\phi_{u+1},\upsilon)=0.
\end{split}
\end{equation}
 Setting $\phi_{u+1}=\upsilon$ we derive our
claim in this case $M=\mu=2$.
\newline

{\bf Proof of  (\ref{kubernhseis}):} We pick out some $b\ne
\alpha$ and will show our claim for the index set
$L^{1,1|(\alpha,b)}$. We distinguish
 two cases for the proof, based on the form of the factor
$F_\alpha$. Either $F_\alpha$ is of the form
$\nabla^{(p)}\Omega_k$ or of
 the form $S_{*}\nabla^{(\nu)}R_{ijkl}$. We begin with the first case,
  which is the easiest.
\newline

 {\it First case: the factor $F_\alpha$ is in the form $\nabla^{(p)}\Omega_h$:} 
 We consider the equation $L_{g}(\Omega_1,\dots ,
  \Omega_{k-1},\Omega_k\cdot\phi_{u+1},\dots\Omega_p,\phi_1,
\dots,\phi_u)=0$ and we pick out the sublinear combination
  $T_{g}(\nabla\phi_{u+1})$ of complete contractions of
 length
  $\sigma+u+1$ with a factor $\nabla\phi_{u+1}$
  contracting against $F_b$. Clearly, it follows that
  $T_{g}(\nabla\phi_{u+1})=0$, modulo complete contraction
 of length $\ge\sigma +u+2$. Moreover, we calculate:

  \begin{equation}
  \label{hilbert}
  \begin{split}
&(0=)T_{g}(\nabla\phi_{u+1})=\Sum_{l\in L^{1,1|\alpha,b}}a_l
Xdiv_{i_1}C^{l,i_1i_2}_{g}(\Omega_1,\dots ,\Omega_p,\phi_1,\dots
\phi_u)\nabla_{i_2}\phi_{u+1}
\\& +2\Sum_{l\in L^{2,0|\alpha}} a_l
  Xdiv_{i_1} \dot{C}^{l,i_1|i_{*}\rightarrow b}_{g}(\Omega_1,\dots ,\Omega_p,\phi_1,\dots
  \phi_u)\nabla_{i_{*}}\phi_{u+1}
  \\&+\Sum_{l\in L'} a_l
  Xdiv_{i_2} C^{l,i_1i_2}_{g}(\Omega_1,\dots ,\Omega_p,\phi_1,\dots
  \phi_u)\nabla_{i_1}\phi_{u+1}+
  \\&\Sum_{h\in H} a_h Xdiv_{i_2}\dots Xdiv_{i_t}C^{h,i_1\dots
  i_t}_{g}(\Omega_1,\dots ,\Omega_p,\phi_1,\dots
  \phi_u)\nabla_{i_1}\phi_{u+1}
\\&+\Sum_{j\in J} a_j C^j_{g}(\Omega_1,\dots ,\Omega_p,\phi_1,\dots
  \phi_{u+1})=0;
  \end{split}
  \end{equation}
here the tensor fields indexed in $L'$ have the free index ${}_{i_2}$
{\it not} belonging to the factor $F_\alpha$ and also have a $u$-simple character
 $\vec{\kappa}_{simp}$ and $\nabla\phi_{u+1}$ is
 not contracting against a special index in $F_\alpha$.
 The tensor fields in $H$
 have a $u$-simple character $\vec{\kappa}_{simp}$ and  have rank $t\ge 3$.
The factor $\nabla\phi_{u+1}$ may be contracting against a
 special index in $F_\alpha$, and in that case if $t=3$ then
 the other two free indices must be non-special.
  Also, we note that
the tensor fields indexed in $L', H$ potentially have one
non-acceptable factor $\nabla\Omega_k$ (with only one derivative),
and in that case if $t=3$ then the two free indices are not
special. We call the tensor fields with a factor $\nabla\Omega_k$
or with the factor $\nabla\phi_{u+1}$ contracting against a
special index ``bad'' tensor fields.

\par In the case where $\nabla\Omega_k$ is not contracting against a
factor $\nabla\phi$,
  we may ``get rid'' of the bad tensor fields  in (\ref{hilbert})
(modulo introducing tensor fields in the general form $\sum_{l\in
L'}\dots$, $\sum_{h\in H}\dots$ which are not bad) via  Corollary
2 in \cite{alexakis4}\footnote{This can be applied since the terms
with minimum rank do not have special free indices, hence there is
no danger of falling under a ``forbidden case''.} or Lemma
4.7 in that paper. In the case where $\nabla\Omega_k$ is contracting
against a factor $\nabla\phi$, we may ``get rid'' of the bad
tensor fields  in (\ref{hilbert}) (modulo introducing tensor
fields in the general form $\sum_{l\in L'}\dots$, $\sum_{h\in
H}\dots$ which are not bad) via  Lemma
4.6 in \cite{alexakis4}.\footnote{This Lemma can be applied since the
terms with minimum rank do not have special free indices, hence
there is no danger of falling under a ``forbidden case''.}

\par Therefore, we may assume wlog that
all the tensor fields in (\ref{hilbert}) are
acceptable and $\nabla\phi_{u+1}$ is not contracting against
a special index. Therefore, this modified equation
 (\ref{hilbert}) shows (\ref{kubernhseis})
 in the case $F_\alpha=\nabla^{(A)}\Omega_x$.
\newline

{\it Proof of (\ref{kubernhseis}) when the crucial factor $F_\alpha$ is in the form
$S_{*}\nabla^{(\nu)} R_{ijkl}$:}

 We define an operation $Link_{\alpha b}$
 that acts on the complete contractions in $L_{g}$ 
(recall that $L_g$ stands for the LHS of the hypothesis of Lemma \ref{pskovb}) by re-writing
 them in dimension $n+2$ and then
 hitting the factor $F_\alpha$ by a derivative $\nabla_c$ and the
factor $F_b$ by a derivative $\nabla^c$. (As long as the weight of the complete contractions is
$-n-2$, we will be considering the re-writing of our complete
contractions in dimension $n+2$).

\par We now consider the  $Image^{1,\beta}_{\phi_{u+1}}[L_{g}]$.\footnote{Recall the 
definition of this sublinear combination from \cite{alexakis6}.}
 Recall the sublinear combinations
$Image^{1,\beta,\sigma+u}_{\phi_{u+1}}[L_{g}]$,
$Image^{1,\beta,\sigma+u+1}_{\phi_{u+1}}[L_{g}]$.
 We recall that in the sublinear combination of length
$\sigma+u+1$ in $Image^{1,\beta}_{\phi_{u+1}}[L_{g}]$ we can still
identify
 the factors $F_\alpha,F_b$. On the other hand, in the sublinear
 combination of length $\sigma+u$ in
$Image^{1,\beta}_{\phi_{u+1}}[L_{g}]$ we will now define the
factors $F_\alpha$, $F_b$: Recall that any complete contraction
$C^t_{g}(\Omega_1,\dots ,\Omega_p,\phi_{u+1},\phi_1,\dots,\phi_u)$
in this sublinear combination has arisen by replacing a (possibly
symmetrized) curvature factor $\nabla^{(f)} R_{ijkl}$ by an
expression $\nabla^{(f+2)}\phi_{u+1}\otimes g$, provided
 the two
indices in $g$ then contract against two indices in the same
factor. Now, if the curvature factor that was replaced was {\it
not} $F_\alpha$ or $F_b$, then we can straightforwardly  identify
$F_\alpha$ or $F_b$ in $C^t_{g}$. If the curvature factor that was
replaced was $F_a$, we will now define this new term
$\nabla^{(p+2)}\phi_{u+1}$ to be the factor $F_a$. We use the same
convention if the factor $F_b$ was the curvature term that was
replaced.

Thus, we can define the operation $Link_{a,b}[\dots]$ on the
complete contractions in the sublinear combination
$Image^{1,\beta}_{\phi_{u+1}}[L_{g}]$.

\par In view of the equation $Image^{1,\beta}_{\phi_{u+1}}[L_{g}]=0$
 we derive an equation:

\begin{equation}
\label{proxrouts}
\begin{split}
&Link_{a,b}\{ Image^{1,\beta}_{\phi_{u+1}}[L_{g}]\}-
Image^{1,\beta}_{\phi_{u+1}}[Link_{a,b}\{L_{g}\}]= \Sum_{w\in W}
a_w C^w_{g}(\phi_{u+1}),
\end{split}
\end{equation}
which holds  modulo complete contractions of length $\ge\sigma +u +2$.
The complete contractions $C^w$ have length $\sigma+u+1$ and
 a factor $\nabla^{(p)}\phi_{u+1}$, $p\ge 2$.
We observe that the left hand side of the above consists
 of complete contractions with length $\ge\sigma +u+1$.
Thus, we can derive:
\begin{equation}
\label{xrouts}
\begin{split}
&Link_{a,b}\{ Image^{1,\beta}_{\phi_{u+1}}[L_{g}]\}-
Image^{1,\beta}_{\phi_{u+1}}[Link_{a,b}\{L_{g}\}]=0,
\end{split}
\end{equation}
 modulo complete contractions of length $\ge\sigma +u+2$.

 Moreover,
  any complete contraction of length $\sigma+u+1$ in the left
   hand side of the above will have a factor
$\nabla\phi_{u+1}$
    and an internal contraction involving a derivative index.
 We denote the left hand side by $Diff[L_{g}]$. We then
 define
$Diff^{*}[L_{g}]$ the sublinear combination of complete
contractions where $\nabla\phi_{u+1}$ is contracting
 against the factor $F_\alpha$ and the internal contraction belongs to the
factor $F_b$. We then derive:

\begin{equation}
\label{xrouts}
\begin{split}
Diff^{*}[L_{g}]=0.
\end{split}
\end{equation}
It is actually quite easy to understand how 
$Diff^{*}[L_{g}]$ arises from $L_g$:

{\it Description of the linear combination $Diff^{*}[L_{g}]$:}
 Let us write out $L_{g}$ as a linear
combination of complete contractions:

$$L_{g}=\Sum_{t\in T} a_t C^t_{g}(\Omega_1,\dots ,
\Omega_p,\phi_1,\dots ,\phi_u)$$ Then, for each $t\in T$ we denote
by $Z^t$ the set of
 particular contractions $({}_x,{}_y)$ in $C^t_{g}$ for which ${}_x$
 belongs to $F_\alpha$ and ${}_y$ belongs to $F_b$. We formally define
$Rep_{(x,y)}[C^t_{g}]$ to stand for the complete contraction that
 arises from $C^t_{g}$ be erasing the contraction $g^{xy}$
 and making ${}_x$ contract against a factor
$\nabla^x\phi_{u+1}$ and then adding a derivative index
$\nabla^y$ onto the factor $F_b$ to which the index ${}_y$
belongs (so now ${}_y$ is contracting against $\nabla^y$ and we
have obtained an internal contraction).

 \par We then calculate:

\begin{equation}
\label{heat} (0=) Diff^{*}[L_{g}]=\Sum_{t\in T} a_t \Sum_{(x,y)\in
Z^t} Rep_{(x,y)}[C^t_{g}],
\end{equation}
modulo complete contractions of length $\ge\sigma +u+2$. We also
define $Rep^{*}_{(x,y)}[\dots]$ to stand for the operation that
replaces the internal contraction $({}^y,{}_y)$ by an expression
$(\nabla^y\omega,{}_y)$. By virtue of the operation $Sub_\omega$
(defined in the Appendix of \cite{alexakis1})
applied to (\ref{heat}) we derive:

\begin{equation}
\label{hea2} 0= \Sum_{t\in T} a_t \Sum_{(x,y)\in Z^l}
Rep^{*}_{(x,y)}[C^t_{g}],
\end{equation}
modulo complete contractions of length $\ge\sigma +u+3$.

\par Now in order to prove our claim (\ref{kubernhseis}), we define the
 operation $Hit^{\nabla\phi_{u+1}}_{F_\alpha}$ that acts on
 complete contractions and tensor fields by hitting the factor
 $F_\alpha$ by a derivative index $\nabla_i$ which we then contract
 against some factor $\nabla^i\phi_{u+1}$ (and in addition since
  $F_\alpha$ is in the form $S_{*}\nabla^{(\nu)}R_{ijkl}$ we
$S_{*}$-symmetrize.

\par Furthermore, by construction we now observe that we can
then write the right hand side of the above as:

\begin{equation}
\label{kryo}
\begin{split}
&(0=)\Sum_{t\in T} a_t \Sum_{(x,y)\in Z^t}
Rep^{*}_{(x,y)}[C^t_{g}]\\&= \Sum_{l\in L^{1,1|a,b}} a_l Xdiv_{i_1}
Hit^{\nabla\phi_{u+1}}_{F_a}C^{l,i_1i_2}_{g} (\Omega_1,\dots
,\Omega_p,\phi_1,\dots ,\phi_u)\nabla_{i_2} \omega
\\&+2\Sum_{l\in L^{2,0|a}} a_l Xdiv_{i_1}
Hit^{\nabla\phi_{u+1}}_{F_a}[\dot{C}]^{l,i_1|i_{*} \rightarrow
b}_{g}(\Omega_1,\dots ,\Omega_p,\phi_1,\dots ,
\phi_u)\nabla_{i_{*}}\omega+
\\&\Sum_{l\in \tilde{L}} a_l Xdiv_{i_3}
C^{l,i_1i_2i_3}_{g}(\Omega_1,\dots ,\Omega_p,\phi_1,\dots ,
\phi_u)\nabla_{i_1}\phi_{u+1}\nabla_{i_2}\omega+
\\&\Sum_{h\in H} a_h Xdiv_{i_3}\dots Xdiv_{i_w}
C^{h,i_1\dots i_w}_{g} (\Omega_1,\dots ,\Omega_p,\phi_1,\dots ,
\phi_u)\nabla_{i_1}\phi_{u+1}\nabla_{i_2}\omega+
\\&\Sum_{j\in J} a_j C^{j,i_1i_2}_{g}
(\Omega_1,\dots ,\Omega_p,\phi_1,\dots ,
\phi_u)\nabla_{i_1}\phi_{u+1}\nabla_{i_2}\omega;
\end{split}
\end{equation}
here the tensor fields indexed in $\tilde{L}$ are all acceptable
and have the feature that the free index ${}_{i_3}$ does {\it not}
belong to $F_\alpha$. Moreover, for each tensor field
indexed in
 $\tilde{L}$ both the factors
$\nabla\phi_{u+1}$, $\nabla\omega$ are {\it not} contracting
against dangerous indices.\footnote{Recall that a dangerous index
is either an internal index in some $\nabla^{(m)}R_{ijkl}$ or an
index ${}_k,{}_l$ in some $S_{*}\nabla^{(\nu)}R_{ijkl}$, or an
index ${}_j$ in some $S_{*}R_{ijkl}$.}

\par The tensor fields indexed in $H$ each have $w\ge 4$ and
are necessarily acceptable. Moreover, either one or both of the factors
$\nabla\phi_{u+1}$, $\nabla\omega$ may be
 contracting against special indices (if they
are contracting against curvature factors).\footnote{This follows by construction, 
since these tensor fields arise from the tensor fields of minimum 
rank $2$ in (\ref{hypothese2}), all of whose free indices are non-special. } Furthermore,
  the tensor fields indexed in $J$ are $u$-subsequent to
$\vec{\kappa}_{simp}$.

\par Now, we break up the index set $H$: We index in $H_{*}$ the
tensor fields that have the factor $\nabla\phi_{u+1}$  contracting
against a special index in the factor
$F_\alpha=S_{*}\nabla^{(\nu)} R_{ijkl}$. We want to derive an
equation that will be precisely like (\ref{kryo}), only with the
extra restriction that $H_{*}=\emptyset$.
\newline

{\it Proof that we may assume wlog that
$H_{*}=\emptyset$ in (\ref{kryo}):} We may assume with no loss of 
generality (by just switching the last two indices in a curvature factor) that for each
$h\in H_{*}$, $\nabla\phi_{u+1}$ is contracting against the index
${}_k$ in $F_\alpha$. We then denote by
$$C^{h,i_2\dots i_w}_{g} (\Omega_1,\dots
,\Omega_p,Y,\phi_1,\hat{\phi}_{y_1},\dots\,
\hat{\phi}_{y_{t-1}},\dots , \phi_u)$$ the tensor fields that
arises from $C^{h,i_1\dots i_t}_{g}$ by replacing the expression

$$S_{*}\nabla^{(\nu)}_{r_1\dots
r_\nu}R_{ijkl}\nabla^i\tilde{\phi}_h\nabla^k\phi_{u+1}\nabla^{r_1}\phi_{y_1}\dots
\nabla^{r_t}\phi_{y_t}$$ by $\nabla^{(\nu+2-(t-1))}_{r_t\dots
r_\nu jl}Y$ (recall that $t>0$).
 (Notice that we are obtaining complete contractions of
weight $-n', n'\le n$, and moreover with $\sigma_1+\sigma_2=s-1$).
We observe that all the tensor fields thus constructed will have
the same $(u-(t-1))$-simple character, which we denote by
$Cut(\vec{\kappa}_{simp})$. Then,  we
derive an equation:

\begin{equation}
\label{soula}
\begin{split}
&\Sum_{h\in H_{*}} a_h Xdiv_{i_3}\dots Xdiv_{i_t} C^{h,i_2\dots
i_t}_{g} (\Omega_1,\dots ,\Omega_p,Y,\phi_1,\dots
,\hat{\phi}_{y_1},\dots ,\hat{\phi}_{y_{t-1}},\dots
\phi_u)\nabla_{i_2}\omega
\\&+\Sum_{j\in J} a_j C^{j,i_2}_{g}
(\Omega_1,\dots ,\Omega_p,Y,\phi_1,\dots ,\hat{\phi}_{y_1},\dots ,\hat{\phi}_{y_{t-1}},\dots
\phi_u)\nabla_{i_2}\omega=0,
\end{split}
\end{equation}
where each $C^{j,i_2}_{g}$ is simply subsequent to
$\vec{\kappa}_{simp}$.

\par Let $\tau\ge 2$ be the minimum rank of the tensor fields appearing above,
 and suppose they are indexed in $H_{*,\tau}$.
Thus, (except for some special cases which we will treat momentarily),
 we  apply our inductive
assumption of Corollary 1 in \cite{alexakis4} to 
the above;\footnote{{\it Except} when there are 
``forbidden cases'' appearing in the above--we will treat that case 
below. Note that if $\tau=2$ (\ref{soula}) can not fall
 under a forbidden case, since all free indices in the terms 
of minimum rank will be non-special.} and derive
that for some linear combination of acceptable tensor fields
(indexed in $P$ below) with a simple character
$Cut(\vec{\kappa}_{simp})$, so that:

\begin{equation}
\label{soula2}
\begin{split}
&\Sum_{h\in H_{*,\tau}} a_h C^{h,i_2\dots i_{\tau-1}}_{g}
(\Omega_1,\dots
,\Omega_p,Y,\phi_1,\hat{\phi}_{y_1},\dots\,\hat{\phi}_{y_{t-1}},\dots
,\phi_u)\nabla_{i_2}\omega\nabla_{i_3}\upsilon\dots\nabla_{i_{\tau-1}}\upsilon-
\\&\Sum_{p\in P} a_p Xdiv_{i_\tau}C^{p,i_2\dots i_{\tau}}_{g}
(\Omega_1,\dots
,\Omega_p,Y,\phi_1,\hat{\phi}_{y_1},\dots\,\hat{\phi}_{y_{t-1}},\dots
,\phi_u)\nabla_{i_2}\omega\nabla_{i_3}\upsilon\dots\nabla_{i_{\tau-1}}\upsilon
\\&=\Sum_{j\in J} a_j C^{j,i_2\dots i_{\tau-1}}_{g} (\Omega_1,\dots
,\Omega_p,Y,\phi_1,\hat{\phi}_{y_1},\dots\,\hat{\phi}_{y_{t-1}},\dots
,\phi_u)\nabla_{i_2}\omega\nabla_{i_3}\upsilon\dots\nabla_{i_\tau-1}\upsilon.
\end{split}
\end{equation}
\par Now, we act on the above by an operation $Op$ that formally
replaces the expression $\nabla^A_{r_{*}r_1\dots
r_A}\phi_{u+1}\nabla^{r_{*}}\phi_{y_t}$ by an expression
$$S_{*}\nabla^{A+t-1}_{r_1\dots r_{t-1}r_{*}\dots
r_{A-2}}R_{ir_{A-1}kr_A}\nabla^i\tilde{\phi}_x\nabla^k\phi_{u+1}\nabla^{r_1}\phi_{y_1}\dots
\nabla^{y_{t-1}}\phi_{y_{t-1}}\nabla^{r_{*}}\phi_{y_t}.$$
Thus, as explained in the proof that Lemma 3.1 in \cite{alexakis4}  implies
Proposition 2.1 in \cite{alexakis4}, we derive a new equation:

\begin{equation}
\label{soula3}
\begin{split}
&\Sum_{h\in H_{*,\tau}} a_h C^{h,i_1\dots i_{\tau-1}}_{g}
(\Omega_1,\dots ,\Omega_p,\phi_1,\dots
,\phi_u)\nabla_{i_1}\phi_{u+1}\nabla_{i_2}\omega\nabla_{i_3}\upsilon\dots
\nabla_{i_{\tau-1}}\upsilon=
\\&\Sum_{h\in H'_{\tau}} a_h C^{h,i_1\dots i_{\tau-1}}_{g}
(\Omega_1,\dots ,\Omega_p,\phi_1,\dots
,\phi_u)\nabla_{i_1}\phi_{u+1}\nabla_{i_2}\omega
\nabla_{i_3}\upsilon\dots\nabla_{i_{\tau-1}}\upsilon+
\\&\Sum_{p\in P} a_p Xdiv_{i_\tau}C^{p,i_2\dots i_{\tau}}_{g}
(\Omega_1,\dots ,\Omega_p,\phi_1,\dots ,
\phi_u)\nabla_{i_1}\phi_{u+1}\nabla_{i_2}\omega
\nabla_{i_3}\upsilon\dots\nabla_{i_{\tau-1}}\upsilon+
\\&\Sum_{j\in J} a_j C^{j,i_2\dots i_{\tau-1}}_{g} (\Omega_1,\dots
,\Omega_p,\phi_1,\dots , \phi_u)\nabla_{i_1}\phi_{u+1}\nabla_{i_2}
\omega\nabla_{i_3}\upsilon\dots\nabla_{i_\tau-1}\upsilon,
\end{split}
\end{equation}
where the tensor fields indexed in $H'_\tau$ are acceptable, have
a $u$-simple character $\vec{\kappa}_{simp}$ and the factor
$\nabla\phi_{u+1}$ is contracting against $F_\alpha$, but not
against a special index. The tensor fields indexed in $J$ are
simply subsequent to $\vec{\kappa}_{simp}$.

\par Now, replace the $\nabla\upsilon$'s by $Xdiv$s, and then replace
into (\ref{kryo}); iterating this step gives our desired
conclusion (which is an equation in the form of (\ref{kryo}) with
$H_{*}=\emptyset$). If at that last step of this process we 
fall under a ``forbidden case'' of Corollary 1, we apply instead the
 ``weaker version'' of Corollary 1 from the Appendix in \cite{alexakis4}. 
As we have noted above, $\tau>2$ is we fall under ``forbidden cases'', hence 
the ``weaker version'' of Corollary 1 gives us our claim. 

\par Therefore, we may now assume that each of the tensor fields
indexed in $H$ have the factor $\nabla\phi_{u+1}$ not contracting
against a special index in $F_\alpha$. $\Box$
\newline

\par Now our claim is only one step away: We apply the operation
$Erase_{\phi_{u+1}}$ in (\ref{kryo}) (since
$H_{*}=\emptyset$ and $F_\alpha$ is a {\it non-simple} factor
of the form $S_{*}\nabla^{(\nu)}R_{ijkl}$) this is well-defined). This is precisely
our claim for (\ref{kubernhseis}).
\newline

{\it Proof of our claim when the critical factor is of the form
$\nabla^{(m)}R_{ijkl}$:}
 An important observation: This hypothesis
means that there are {\it no} 2-tensor fields in $L_2$ in our
Lemma hypothesis with two free indices belonging to the same
factor of the form $\nabla^{(B)}\Omega_h$ or $S_{*}\nabla^{(\nu)}
R_{ijkl}$. This follows by the definition of the critical factor.

\par We denote by $F_1,\dots F_{\tau_1}$ the non-generic factors $\nabla^{(m)}R_{ijkl}$
in $\vec{\kappa}_{simp}$. We will then have that
$L^{2,0}_2=\bigcup_{f=1}^{\tau_1} L_2^{2,0|f}\bigcup
L_2^{2,0|\tau+1}$. Again denote by $L^{2,0|\alpha}_2$ the index
set that corresponds to $\bigcup_{z\in Z'_{Max}}L^z$. Thus the
critical factor will again be denoted by $F_\alpha$ (it will now
be in the form $\nabla^{(m)} R_{ijkl}$). We first consider all the
sets $L^{1,1|\alpha,b}$, $b=\tau_1+1,\dots b=\tau$ for which $F_b$
is {\it not} in the form $\nabla^{(m)}R_{ijkl}$ (and is also {\it
not} a simple factor of the form $S_{*}\nabla^{(\nu)}R_{ijkl}$). In that case, 
we  claim that we can write:

\begin{equation}
\label{mep}
\begin{split}
&Xdiv_{i_1}Xdiv_{i_2}\Sum_{l\in L^{1,1|\alpha,b}}a_l
C^{l,i_1i_2}_{g}(\Omega_1,\dots ,\Omega_p,\phi_1,\dots ,\phi_u)=
\\&-Xdiv_{i_1}Xdiv_{i_2}\Sum_{l\in L^{2,0|\alpha}}a_l \dot{C}^{l,i_1|i_{*}\rightarrow b}_{g}
(\Omega_1,\dots ,\Omega_p,\phi_1,\dots ,\phi_u)+
\\&\Sum_{h\in H} a_h Xdiv_{i_1}\dots Xdiv_{i_t} C^{h,i_1\dots
i_t}_{g} (\Omega_1,\dots ,\Omega_p,\phi_1,\dots ,
\phi_u)\\&+\Sum_{j\in J} a_j C^j_{g}(\Omega_1,\dots
,\Omega_p,\phi_1,\dots , \phi_u),
\end{split}
\end{equation}
where the tensor fields indexed in $H$ $u$-simple character
$\vec{\kappa}_{simp}$ and are acceptable and have $t\ge 3$. The
complete contractions $C^j$ are $u$-subsequent to
$\vec{\kappa}_{simp}$. Note that if we can show 
the above, then by applying Lemma \ref{oui2},
we may also assume that they satisfy all the extra hypotheses of
Lemma \ref{pskovb} pertaining to the extra claims.

\par Secondly, we consider the tensor fields indexed
in $L^{1,1|\alpha,b}$ where $F_b$ is in the form
$\nabla^{(m)}R_{ijkl}$, where if $\alpha\in \{1,\dots, \tau_1\}$
then $b\ne\alpha$ (if $\alpha=\tau_1+1$, there are no
restrictions). We claim that we can write:

\begin{equation}
\label{mep2}
\begin{split}
&Xdiv_{i_1}Xdiv_{i_2}\Sum_{l\in L^{1,1|\alpha,b}}a_l
C^{l,i_1i_2}_{g}(\Omega_1,\dots ,\Omega_p,\phi_1,\dots ,\phi_u)=
\\&-Xdiv_{i_1}Xdiv_{i_2}\Sum_{l\in L^{2,0|\alpha}}a_l \dot{C}^{l,i_1|i_{*}\rightarrow b}_{g}
(\Omega_1,\dots ,\Omega_p,\phi_1,\dots ,\phi_u)
\\&-Xdiv_{i_1}Xdiv_{i_2}\Sum_{l\in
L^{2,0|b}}a_l \dot{C}^{l,i_1|i_{*}\rightarrow a}_{g}
(\Omega_1,\dots ,\Omega_p,\phi_1,\dots ,\phi_u)+
\\&\Sum_{h\in H} a_h Xdiv_{i_1}\dots Xdiv_{i_t} C^{h,i_1\dots
i_t}_{g} (\Omega_1,\dots ,\Omega_p,\phi_1,\dots ,
\phi_u)\\&+\Sum_{j\in J} a_j C^j_{g}(\Omega_1,\dots
,\Omega_p,\phi_1,\dots , \phi_u),
\end{split}
\end{equation}
with the same notational conventions as above.

{\it Proof of (\ref{mep}) and (\ref{mep2}):} Both equations are
easy to derive: We just consider the equation
$Im^{1,\beta}_{\phi_{u+1}}[L_{g}]=0$ (see (\ref{antrikos}) and we pick out the sublinear
combination where $\nabla\phi_{u+1}$ is contracting against
$F_\alpha$ and $\nabla\omega$ against $F_b$. We then change both
$\nabla\phi_{u+1}$ and $\nabla\omega$ into $Xdiv$s. $\Box$
\newline

{\it Derivation of Lemma \ref{pskovb} in the case $M=\mu=2$
 from the equations (\ref{mep}), (\ref{mep2}) when 
the critical factor is of the form $\nabla^{(m)}R_{ijkl}$:}
 Some notation: For each $\alpha\le\tau_1$ and each
$b\le\tau_1$ with $b\ne \alpha$ or $b=\tau+1$, we denote by
$C^{l,\{i_1i_2\}\rightarrow F_b}_{g}$ the tensor field that
formally arises from $C^{l,i_1i_2}_{g}$, $l\in L^{2,0|\alpha}$ by
erasing the two free indices $\nabla_{i_1},\nabla_{i_2}$ from
$F_\alpha$ and adding them onto the factor $F_b$ (if $b=\tau+1$ we
add them onto any generic factor $\nabla^{(m)}R_{ijkl}$ and then
add over all the tensor fields we can thus obtain). For
$\alpha=\tau+1$ and $b\le\tau_1$ we define
$C^{l,\{i_1i_2\}\rightarrow F_b}_{g}$ in the same way. If
$\alpha=\tau+1,b=\tau+1$, we impose the extra restriction that
 the free indices are erased from the factor $F_{*}$ to which they belong and then we add them to
 any other generic factor $\nabla^{(m)}R_{ijkl}$ other than $F_{*}$;
  we then add over all the tensor fields that we have obtained.

\par Now, for each $a\le \tau_1$ and also for $a=\tau+1$ we
consider the grand conclusion with $F_a$ being the selected
factor. We denote by $\Sum_{t\in T_a} a_t C^{t,i_1i_2}_{g}$ a
generic linear combination of acceptable vector fields for which
the free index ${}_{i_2}$ does not belong to $F_a$. We then derive an equation
 for each such $a$:

\begin{equation}
\label{britz}
\begin{split}
&(-2q_a-2(\sigma_1-2+\sigma_2)) \Sum_{l\in L^{2,0|a}_2} a_l
Xdiv_{i_2}C^{l,i_1i_2}_{g} (\Omega_1,\dots ,\Omega_p,\phi_1,\dots
,\phi_u)\nabla_{i_1}\phi_{u+1}
\\&-2\Sum_{b\in \{i,\dots ,\tau_1,\tau+1\}}
\Sum_{l\in L^{2,0|b}_2}a_l Xdiv_{i_2} C^{l,\{i_1i_2\} \rightarrow
F_a}_{g}(\Omega_1,\dots ,\Omega_p,\phi_1,\dots
,\phi_u)\nabla_{i_1}\phi_{u+1}
\\&+\Sum_{t\in T_a} a_t Xdiv_{i_2}C^{t,i_1i_2}_{g}(\Omega_1,\dots ,\Omega_p,\phi_1,\dots
,\phi_u)\nabla_{i_1}\phi_{u+1}+
\\&\Sum_{h\in H} a_h Xdiv_{i_2}\dots Xdiv_{i_t}
C^{h,i_1\dots i_t}_{g}(\Omega_1,\dots ,\Omega_p,\phi_1,\dots
,\phi_u)\nabla_{i_1}\phi_{u+1}+
\\&\Sum_{j\in J} a_j C^{j,i_1}_{g}(\Omega_1,\dots ,\Omega_p,\phi_1,\dots
,\phi_u)\nabla_{i_1}\phi_{u+1},
\end{split}
\end{equation}
where each $C^j$ is simply subsequent to $\vec{\kappa}_{simp}$ and
each of the tensor fields indexed in $H$ have $t\ge 3$ and are
otherwise as in the conclusion of Lemma \ref{pskovb}. The constant
$q_a$ is equal to $n-2u-\mu-2$ and is strictly positive
 (therefore observe that $-2q_a-2(\sigma_1-2+\sigma_2)<0$).
Therefore, by applying the inductive assumption of Corollary
1 in \cite{alexakis4}\footnote{The terms of minimum rank above all involve
terms with non-special free indices, hence there is no danger of falling under
a ``forbidden case''. Also, we fall under the inductive assumption of that Lemma
because we increased the number of factors $\nabla\phi_h$, 
while keeping all the other parameters of 
the induction fixed.} we derive an equation, for each $F_a$ in
the form $\nabla^{(m)}R_{ijkl}$:

\begin{equation}
\label{britz2}
\begin{split}
&\Sum_{l\in L^{2,0|a}_2} a_l C^{l,i_1i_2}_{g}(\Omega_1,\dots
,\Omega_p,\phi_1,\dots
,\phi_u)\nabla_{i_1}\phi_{u+1}\nabla_{i_2}\upsilon+
\\&\frac{1}{q_a+2(\sigma_1-2+\sigma_2)}
\Sum_{b\in \{i,\dots ,\tau_1,\tau+1\}} \Sum_{l\in L^{2,0|b}_2}a_l
C^{l,\{i_1i_2\}\rightarrow F_a}_{g}(\Omega_1,\dots
,\Omega_p,\phi_1,\dots
,\phi_u)\\&\nabla_{i_1}\phi_{u+1}\nabla_{i_2}\upsilon
+\Sum_{h\in H} a_h Xdiv_{i_3} C^{h,i_1\dots i_t}_{g}(\Omega_1,\dots ,\Omega_p,\phi_1,\dots
,\phi_u)\nabla_{i_1}\phi_{u+1}\nabla_{i_2}\upsilon
\\&=\Sum_{j\in J} a_j C^{j,i_1}_{g}(\Omega_1,\dots ,\Omega_p,\phi_1,\dots
,\phi_u)\nabla_{i_1}\phi_{u+1}.
\end{split}
\end{equation}
We are therefore reduced to proving our claim under the assumption
that the sublinear combination $\Sum_{l\in L^{2,0}_2} a_l
C^{l,i_1i_2}_{g}$ in our Lemma hypothesis can be expressed
 in the form:

$$\nabla^{2,spread}_{i_1i_2}[\Sum_{b\in B} a_b C^b_{g}(\Omega_1,\dots ,\Omega_p,\phi_1,\dots
,\phi_u)],$$ where the complete contractions $C^b_g$ are in the form (\ref{form2}), have weight
$-n+4$, are acceptable of simple character
$\vec{\kappa}_{simp}$.
The symbol $\nabla^{2,spread}_{i_1i_2}$ means that we may hit any
factor  of the form $\nabla^{(m)}R_{ijkl}$, and then we add over all
these factors we have hit.

\par Therefore, applying (\ref{mep}) and (\ref{mep2}) to this
setting we may in addition assume that the sublinear combination
of 2-tensor fields in (\ref{hypothese2}) where one free index
${}_{i_1}$ belongs to a factor $\nabla^{(m)} R_{ijkl}$ and one to a factor
$S_{*}\nabla^{(\nu)} R_{ijkl}$ is:

$$-\nabla^{2|i_1\rightarrow \nabla^{(m)}R_{ijkl},i_2\rightarrow
S_{*}\nabla^{(\nu)} R_{ijkl}} [\Sum_{b\in B} a_b
C^b_{g}(\Omega_1,\dots ,\Omega_p,\phi_1,\dots ,\phi_u)],$$ where
the symbol outside brackets means that we are considering the
sublinear combination in $\nabla^{(2)}_{i_1i_2}$ where ${}_{i_1}$
is forced to hit a factor $\nabla^{(m)}R_{ijkl}$ and ${}_{i_2}$ is
forced to hit a factor $S_{*}\nabla^{(\nu)} R_{ijkl}$ (and we can
analogously obtain a new true equation for the sublinear
combination of terms where one free index belongs to a  factor
$\nabla^{(m)}R_{ijkl}$ and the other to a simple factor of the
form $\nabla^{(y)}\Omega_h$).

Moreover, applying (\ref{mep2}) and making $\nabla
\phi_{u+1}$, $\nabla\upsilon$ into $Xdiv$s, 
(by virtue of the last Lemma in the Appendix of \cite{alexakis1}) 
we may assume that the sublinear combination of 2-tensor
fields in (\ref{hypothese2}) where both free indices ${}_{i_1},{}_{i_2}$
belong to (different) factors $\nabla^{(m)} R_{ijkl}$ is:

$$-2\nabla^{2|i_1\rightarrow \nabla^{(m)}R_{ijkl},i_2\rightarrow
\nabla^{(m)} R_{ijkl}} [\Sum_{b\in B} a_b C^b_{g}(\Omega_1,\dots
,\Omega_p,\phi_1,\dots ,\phi_u)],$$ where the symbol outside
brackets means that we are considering the sublinear combination
in $\nabla^{i_1i_2}$ where ${}_{i_1},{}_{i_2}$ are forced to hit two
different factors of the form $\nabla^{(m)}R_{ijkl}$.

\par Now pick any $A\in \{1,\dots,\tau_1+1\}$. 
Applying the grand conclusion with $F_A$ being the
selected factor, we derive:

\begin{equation}
\label{epitelous}
\begin{split}
&(-2q_1-2(\sigma_1-1)-2\sigma_2) \Sum_{b\in B} a_b
Xdiv_{i_2}\nabla^{2,\{i_1,i_2\}\rightarrow F_1}_{i_1i_2}
C^b_{g}(\Omega_1,\dots ,\Omega_p,\phi_1,\dots
,\phi_u)\\&\nabla_{i_1}\phi_{u+1}+
\Sum_{t\in T} a_tXdiv_{i_2} C^{t,i_1i_2}_{g}(\Omega_1,\dots
,\Omega_p,\phi_1,\dots ,\phi_u)\nabla_{i_1}\phi_{u+1}+
\\&\Sum_{h\in H} a_h Xdiv_{i_2}\dots Xdiv_{i_t} C^{h,i_1\dots i_t}_{g}
(\Omega_1,\dots ,\Omega_p,\phi_1,\dots
,\phi_u)\nabla_{i_1}\phi_{u+1}+
\\&\Sum_{j\in J} a_j C^j_{g}(\Omega_1,\dots ,\Omega_p,\phi_1,\dots
,\phi_{u+1})=0.
\end{split}
\end{equation}
The symbol in the first line means that we hit the factor $F_A$
with two derivatives $\nabla_{i_1 i_2}$. Here the tensor fields
indexed in $T$ have ${}_{i_2}$ {\it not} belonging to $F_A$.
 Since the coefficient in the first line is non-zero, we may apply our inductive assumption of
 Corollary 1 in \cite{alexakis4},\footnote{The terms of minimum rank above 
contain no special free indices,
hence there is no danger of falling under
a ``forbidden case''.}
 pick out the sublinear combination where $\nabla\upsilon$ is contracting against $F_1$
 and then set $\phi_{u+1}=\upsilon$. This is our desired
 conclusion. $\Box$

\subsection{Proof of Lemma \ref{pskovb} in the subcase $\mu=1$.}
\label{mu=1}

\par We proceed to show the remaining case for Lemma
\ref{pskovb}, the case $\mu=1$. Recall that by Lemma \ref{thenewlemma} we may assume that 
there are to $1$-tensor fields in our Lemma hypothesis for which 
the free index belongs to a factor $S_{*}\nabla^{(\nu)}R_{ijkl}$.  

{\it Special cases:} We single out certain special cases which will
 be treated in a Mini-Appendix at the end of this paper. 
The ``special subcases'' are when $\sigma_2>0$ in 
$\vec{\kappa}_{simp}$,\footnote{We recall that $\sigma_2$ stands for the 
number of factors $S_*\nabla^{(\nu)}R_{ijkl}$ in $\vec{\kappa}_{simp}$.} 
and the terms of maximal 
refined double character in (\ref{hypothese2}) 
{\it either} have no removable index,\footnote{In particular, all their 
factors must be in the form $R_{ijkl},S_*R_{ijkl}$ 
without derivatives, {\it or} in the form $\nabla^{(2)}\Omega_h$.} {\it or} 
the refined double characters correspond to the form:
\begin{equation}
\label{tokolpo}
\begin{split} 
&contr(\nabla_{(free)}R_{\sharp\sharp\sharp\sharp}\otimes R_{\sharp\sharp\sharp\sharp}
\otimes\dots\otimes R_{\sharp\sharp\sharp\sharp}\otimes 
\\&S_*R_{i\sharp\sharp\sharp}\otimes\dots
\otimes S_*R_{i\sharp\sharp\sharp}\otimes 
\nabla^{(2)}_{y\sharp}\Omega_1\otimes\dots\otimes\nabla^{(2)}_{y\sharp}\Omega_p\otimes\nabla\phi_1
\otimes\dots\otimes\nabla\phi_u). 
\end{split}
\end{equation}
(In the above, each index ${}_\sharp$ must contract against 
another index in the form ${}_\sharp$; the indices ${}_y$ 
 are either contracting against indices ${}_\sharp,{}_y$ or 
contract against a factor $\nabla\phi_h$).
We remark that in the rest of this proof we will be assuming 
that the terms of maximal refined double character in (\ref{hypothese2}) are {\it not} 
in the form (\ref{tokolpo}) with $\sigma_2>0$.\footnote{We are again 
applying the convention that we do not write out any derivative indices 
in a factor $\nabla^{(m)}R_{ijkl}$ that contract against factors $\nabla\phi$).} 
\newline

In this case the different refined
double characters of the vector fields $C^{l,i_1}_{g}, l\in L_1$
are fully characterized by specifying the factor in
$\vec{\kappa}_{simp}$ to which ${}_{i_1}$ belongs. Recall the
discussion on the index sets $L^z,z\in Z'_{Max}$.  Observe that in this case
the index set $Z'_{Max}$ will consist of one element, and the sublinear 
combination stands for the sublinear combination of 1-tensor fields
in the LHS of the Lemma hypothesis for which the free 
index ${}_{i_1}$ belongs to the {\it critical(=crucial) factor}.\footnote{Recall 
that given the simple character $\vec{\kappa}_{simp}$ of the 
tensor fields appearing in the hypothesis of our Lemma, the crucial factor is either
a well-defined factor in one of the forms $\nabla^{(B)}\Omega_h$, 
$\nabla^{(m)}R_{ijkl}$, $S_{*}\nabla^{(\nu)}R_{ijkl}$ 
{\it or} the generic set of factors $\nabla^{(m)}R_{ijkl}$ 
that are not contracting against any factor $\nabla\phi_h$.} Denote by
$F_\alpha$ the crucial factor.

Again, we start with the case where the critical factor $F_\alpha$ is of
the form $\nabla^{(B)}\Omega_h$.
\newline

 {\bf Proof of the claim when the critical factor 
is in the form $\nabla^{(B)}\Omega_h$:} We denote the index set of the 1-tensor
fields in the Lemma hypothesis where ${}_{i_1}$ belongs to the
crucial factor by $L_1^\alpha$. (In other words 
$\bigcup_{z\in Z'_{Max}} L^z=L_1^\alpha$).

 We will then initially be using an analogue of
(\ref{kubernhseis}). We introduce some notation to serve our purposes.

{\it Notation:} For each $l\in L_1$, we denote by
$C^l_{g}(\Omega_1,\dots ,\Omega_p,\phi_1,\dots ,\phi_u)$
 the complete contraction (of weight $-n+2$) that arises from
$C^{l,i_1}_{g}(\Omega_1,\dots ,\Omega_p,\phi_1,\dots , \phi_u)$ by
erasing the (derivative) index ${}_{i_1}$. For each $x\le\sigma_1$, we then define
\\$Pass_x^{i_1}[C^l_{g}(\Omega_1,\dots ,\Omega_p,\phi_1, \dots
,\phi_u)]$ as follows: If $x\le \sigma_1$, it
 will stand for the vector field that arises from
$C^l_{g}$ by hitting the
 $x^{th}$ factor $\nabla^{(m)}R_{ijkl}$  by a derivative index
$\nabla_{i_1}$.

\par For each $x\le\sigma_1$, we denote by $L^x_1\subset L_1$
the index set of the vector fields in our
 Lemma hypothesis for which ${}_{i_1}$ belongs to the $x^{th}$
 factor. We denote by $\vec{\kappa}_x$ the $(u+1)$-simple
  character that arises by contracting the free index ${}_{i_1}$
against a factor $\nabla\phi_{u+1}$.
\newline

\par We will prove that for each $x, 1\le x\le\sigma_1$:

\begin{equation}
\label{plagiarism}
\begin{split}
& \Sum_{l\in L_1^{x}} a_l C^{l,i_1}_{g} (\Omega_1,\dots
,\Omega_p,\phi_1,\dots ,\phi_u) \nabla_{i_1}\phi_{u+1}
\\& +\Sum_{l\in L_1^{\alpha}} a_l Pass_x^{i_1} C^{l,i_1}_{g}
(\Omega_1,\dots ,\Omega_p,\phi_1,\dots ,\phi_u)
\nabla_{i_1}\phi_{u+1}+
\\& \Sum_{h\in H^x} a_h Xdiv_{i_2}
C^{h,i_1i_2}_{g}(\Omega_1,\dots ,\Omega_p,\phi_1,\dots ,
\phi_u)\nabla_{i_1}\phi_{u+1}+
\\& \Sum_{j\in J} a_j C^{j,i_1}_{g}
(\Omega_1,\dots ,\Omega_p,\phi_1,\dots ,\phi_u)
\nabla_{i_1}\phi_{u+1}=0,
\end{split}
\end{equation}
modulo complete contractions of length $\ge\sigma+u+2$.
   Each $C^{j,i_1}_{g}$ is
$u$-subsequent to $\vec{\kappa}_{simp}$, each vector field
$C^{h,i_1i_2}_{g}$ has a $u$-simple character
$\vec{\kappa}_{simp}$, $(u+1)$-weak character
$Weak(\vec{\kappa}_x)$
 and is either acceptable or has one unacceptable factor $\nabla\Omega_h$.
  We will discuss how (\ref{plagiarism}) follows below. Let us now check
how it applies to show our Lemma.
\newline

{\it Proof that (\ref{plagiarism}) implies case B of Lemma
\ref{pskovb} in this case $\mu=1$ when the critical factor is of
 the form $\nabla^{(p)}\Omega_h$:}
\newline

We make the factor $\nabla\phi_{u+1}$ into an $Xdiv$ in
 each of the above equations (appealing 
 to the last Lemma in the Appendix of \cite{alexakis1}),
  and then we replace the sublinear combination
$$\Sum_{x=1}^{\sigma_1}
Xdiv_{i_1}\Sum_{l\in L_1^{x}} a_l C^{l,i_1}_{g} (\Omega_1,\dots
,\Omega_p,\phi_1,\dots ,\phi_u)$$ in our Lemma hypothesis by

\begin{equation}
\label{ouxaa}
\begin{split}
&-\Sum_{x=1}^{\sigma_1} \Sum_{l\in
L_1^{\alpha}} a_l Xdiv_{i_1} Pass_x^{i_1} C^{l,i_1}_{g}
(\Omega_1,\dots ,\Omega_p,\phi_1,\dots ,\phi_u)+
\\& \Sum_{h\in H} a_h Xdiv_{i_1}Xdiv_{i_2}
C^{h,i_1i_2}_{g}(\Omega_1,\dots ,\Omega_p,\phi_1,\dots , \phi_u)+
\\& \Sum_{j\in J} a_j C^{j}_{g}
(\Omega_1,\dots ,\Omega_p,\phi_1,\dots ,\phi_u).
\end{split}
\end{equation}

\par Thus, we obtain an equation:

\begin{equation}
\label{ouxaa2}
\begin{split}
&\Sum_{l\in L_1^{\alpha}} a_l Xdiv_{i_1}C^{l,i_1}_{g} (\Omega_1,\dots
,\Omega_p,\phi_1,\dots ,\phi_u)
\\&-\Sum_{x=1}^{\sigma_1}
\Sum_{l\in L_1^{\alpha}} a_l  Xdiv_{i_1} Pass_x^{i_1} C^{l,i_1}_{g}
(\Omega_1,\dots ,\Omega_p,\phi_1,\dots ,\phi_u)+
\\& \Sum_{l\in \tilde{L}_1} a_l Xdiv_{i_1}C^{l,i_1}_{g}
(\Omega_1,\dots ,\Omega_p,\phi_1,\dots ,\phi_u)+
\\& \Sum_{h\in H} a_h Xdiv_{i_1}Xdiv_{i_2}
C^{h,i_1i_2}_{g}(\Omega_1,\dots ,\Omega_p,\phi_1,\dots , \phi_u)+
\\& \Sum_{j\in J} a_j C^{j}_{g}
(\Omega_1,\dots ,\Omega_p,\phi_1,\dots ,\phi_u).
\end{split}
\end{equation}

\par The vector fields indexed in $\tilde{L}_1$
 are acceptable but ${}_{i_1}$ belongs to some non-crucial factor
in the form $\nabla^{(B)}\Omega_h$.

\par Again by inspection we have that $L^{+}_1=\emptyset$. By 
applying Lemma \ref{funny}, we may assume wlog that 
 the tensor fields indexed in $H$ are all acceptable.
We also apply Lemmas \ref{oui}, \ref{oui2} if necessary to ensure 
that the above fulfills the requirements to apply the ``grand conclusion''. 
 We may then apply the grand
conclusion to the above, making $F_1$ the selected factor. We
derive an equation:

\begin{equation}
\label{ouxaa2}
\begin{split}
&\Sum_{l\in L_1^{\alpha}} a_l (-q_1- \Sum_{x=1,x\ne
\alpha}^{\sigma_1}\overline{2}_x)C^{l,i_1}_{g}
(\Omega_1,\dots ,\Omega_p,\phi_1,\dots ,\phi_u)
\nabla_{i_1}\phi_{u+1}+
\\& \Sum_{h\in H} a_h Xdiv_{i_2}
C^{h,i_1i_2}_{g}(\Omega_1,\dots ,\Omega_p,\phi_1,\dots ,
\phi_u)\nabla_{i_1}\phi_{u+1}+
\\& \Sum_{j\in J} a_j C^{j}_{g}
(\Omega_1,\dots ,\Omega_p,\phi_1,\dots ,\phi_{u+1});
\end{split}
\end{equation}
here the tensor fields indexed in $H$ have a $u$-simple
 character $\vec{\kappa}_{simp}$, a weak $(u+1)$-character
 $Weak(\vec{\kappa}^{+}_{simp})$ and possibly one
un-acceptable
 factor $\nabla\Omega_h$ ($q_1$ stands for the coefficient
inside parentheses in the first line of the grand conclusion, and
it depends on the form of the selected factor). Using Lemma
\ref{funny}  if necessary we may assume that all terms in $H$ are
acceptable with a $(u+1)$-simple character
$\vec{\kappa}^{+}_{simp}$.

\par Now, under these assumptions, notice that if  $(-q_1-
\Sum_{x=1,x\ne a}^{\sigma_1}2_x)\ne 0$ then dividing by
this number we derive our claim. Let us now derive our claim in the case where 
$(-q_1-\Sum_{x=1,x\ne a}^{\sigma_1}2_x)= 0$.

{\it The remaining case:} Notice that the only case in which the number in parentheses is
zero is when we have $\sigma_1=\sigma_2=0$ and all factors
$\nabla^{(B)}\Omega_h$ have all their non-free indices contracting
against factors $\nabla\phi_h$.
 We also observe that furthermore there can only be one
 free index among all the tensor fields in our induction
 hypothesis (in other words $L_{>1}=\emptyset$). In that case we just
read off our claim directly from
 (\ref{hypothese2}) as follows: We break $L_1$ into sets
 $L^h_1,h=1,\dots ,p$ depending on which factor
$\nabla^{(B)}\Omega_h$ contains the free index. Then, for each
pair $1\le h_1<h_2\le p$ we pick out the
 sublinear combination in (\ref{hypothese2}) where
the one contraction not involving factors $\nabla\phi_x$ is
 between $\nabla^{(B)}\Omega_{h_1},\nabla^{(C)}\Omega_{h_2}$.
 This sublinear combination must vanish separately. We then erase
 the two (derivative) indices involved in these complete
 contractions and we derive:

\begin{equation}
\label{turqui} \{\Sum_{l\in L^{h_1}_1} a_l+\Sum_{l\in L^{h_2}_2}
a_l\}C^l_{g}
\end{equation}
\par This then implies that $\Sum_{l\in L^h_1} a_l=0$ for every $h=1,\dots ,p$, since $p=\sigma\ge 3$.
\newline

\par So, we only have to show (\ref{plagiarism}). But this follows
 by  the exact same argument we used to prove equation (\ref{kubernhseis}) in the case $\mu=2$: 
 We replace factor crucial factor $\nabla^{(B)}\Omega_h$ by 
 $\nabla^{(B)}(\Omega_h\cdot\omega)$  in our Lemma hypothesis 
 and then pick out the sublinear combination with a factor $\nabla\omega$
  contracting against the $x^{th}$ factor $\nabla^{(m)}R_{ijkl}$. 
This sublinear combination must vanish separately and this is precisely 
  our claim, (\ref{plagiarism}). $\Box$
\newline

{\it Proof of the claim when the crucial factor is of the form
$\nabla^{(m)}R_{ijkl}$:} Again, by the definition of the critical(=crucial) 
factor, it follows that no vector field indexed in $L_1$ will have
its free index belonging to a factor $\nabla^{(B)}\Omega_h$ or
$S_{*}\nabla^{(\nu)} R_{ijkl}$.

\par We then repeat the same argument as in the previous case: Let
us denote by $L_1^1,\dots L_1^{\tau_1}, L_1^{\tau_1+1}$ all the
index sets that correspond to the various factors $F_1,\dots
,F_{\tau_1},F_{\tau+1}$\footnote{Recall that $F_1,\dots
F_{\tau_1}$ stand for the factors $\nabla^{(m)}R_{ijkl}$ in
$\vec{\kappa}'_{simp}$ that are contracting against some factor
$\nabla\phi_h$. $F_{\tau+1}$ stands for the set of factors
$\nabla^{(m)}R_{ijkl}$ in $\vec{\kappa}_{simp}$ that are not
contracting against any factors $\nabla\phi_h$.} to which the free
index may belong. Recall that $q_1$ will stand for the coefficient
in the first line of the grand conclusion; recall $q_1\ge 0$. We
distinguish two cases: Either $q_1=0$ or $q_1>0$ (recall $q_1$ is
the coefficient between parentheses in the grand conclusion, in
this case $q_1=n-2u- \mu-3$).  We first prove our claim in the first (very
special) case.

{\it Subcase $q_1=0$:} In that case we clearly have
 that $\sigma_1=1$ and all the indices in each factor
$\nabla^{(B)}\Omega_h$ are contracting against a factor
$\nabla\phi_y$, and also all the indices ${}_{r_1},\dots
,{}_{r_\nu},{}_j$ in each factor $S_{*}\nabla^{(\nu)} R_{ijkl}$ are
contracting against some factor $\nabla\phi'_h$. Furthermore, the
crucial factor must be in the form $\nabla^{(m)}_{r_1\dots
r_m}\nabla_{i_1}R_{ijkl}$ (where ${}_{r_1}, \dots ,{}_{r_m}$ are
contracting against factors $\nabla\phi_x$).
 So we can write:
$$\Sum_{l\in L'_1} a_l C^{l,i_1}_{g}=\nabla_{i_1}^{cruc}[\Sum_{b\in B} a_b C^b_{g}]$$
modulo longer complete contractions, where each $C^b_{g}$ is an
acceptable complete contraction of weight $-n+2$ with simple
character $\vec{\kappa}_{simp}$ and $\nabla_{i_1}^{cruc}$ means
that the derivative $\nabla_{i_1}$ is forced to hit the (unique)
crucial factor $\nabla^{(m)}R_{ijkl}$.

\par We then pick any factor $F_c\ne F_1$ ($F_1$ is the
 crucial factor) and we apply the grand conclusion to
$L_{g}$ with $F_c$ being the selected factor. We obtain an
equation:

\begin{equation}
\label{saoud}
\begin{split}
&\Sum_{b\in B} a_b 2\nabla^{F_c}_{i_1}[C^b_{g}(\Omega_1, \dots
,\Omega_p,\phi_1,\dots ,\phi_u)] \nabla^{i_1}\phi_{u+1}+
\\&\Sum_{h\in H} a_h Xdiv_{i_2}\dots Xdiv_{i_a}
C^{h,i_1\dots i_a}_{g}(\Omega_1, \dots ,\Omega_p,\phi_1,\dots
,\phi_u)\nabla^{i_1} \phi_{u+1}+
\\&\Sum_{j\in J} a_j C^j_{g}(\Omega_1,
\dots ,\Omega_p,\phi_1,\dots ,\phi_u)=0,
\end{split}
\end{equation}
modulo complete contractions of length $\ge\sigma+u+2$. Here the
tensor fields indexed in $H$ have $a\ge 2$ and have
 a $u$-simple character $\vec{\kappa}_{simp}$ and the factor $\nabla\phi_{u+1}$
 contracting against $F_c$. They may also have one
 unacceptable factor $\nabla\Omega_x$.

\par Now, by applying Lemma \ref{funny},\footnote{By inspection,
 the above does not fall under a forbidden case of 
Lemma \ref{funny}.} we may
assume wlog that for each tensor field indexed in $H$ above there
are no unacceptable factors and the factor $\nabla\phi_{u+1}$ are
contracting against a non-special index.

\par We then apply $Erase_{\phi_{u+1}}$ to the above and then
 add a derivative index $\nabla_{i_{*}}$ onto the
crucial factor $F_1$, which we then contract against a factor
$\nabla_{i_{*}}\phi_{u+1}$. This gives us an equation:

\begin{equation}
\label{saoud2}
\begin{split}
&\Sum_{b\in B} a_b 2\nabla^{cruc}_{i_1}[C^b_{g}(\Omega_1, \dots
,\Omega_p,\phi_1,\dots ,\phi_u)]\nabla^{i_1}\phi_{u+1}+
\\&\Sum_{h\in H'} a_h Xdiv_{i_2}\dots Xdiv_{i_a}
C^{h,i_1\dots i_a}_{g}(\Omega_1, \dots ,\Omega_p,\phi_1,\dots
,\phi_u)\nabla^{i_1} \phi_{u+1}+
\\&\Sum_{j\in J} a_j C^j_{g}(\Omega_1,
\dots ,\Omega_p,\phi_1,\dots ,\phi_u)=0,
\end{split}
\end{equation}
where the tensor fields indexed in $H'$ have a $u$-simple
 character $\vec{\kappa}_{simp}$ and are acceptable,
  and $\nabla\phi_{u+1}$ is contracting against a
non-special index in the crucial factor, by construction. This is
our desired conclusion in this case.
\newline

{\it Subcase $q_1>0$:} 
Observe that we will either have $q_1>2(\sigma_1-1)$ or
$q_1=2(\sigma_1-1)$.\footnote{Notice (by a counting argument)
that if the 1-tensor fields of maximal refined double character 
in (\ref{hypothese2}) are in the form (\ref{tokolpo}), then $q_1=2(\sigma_1-1)$.} 
For both cases the starting point will be the
same:

\par We pick any  $c\in
\{1,\dots,\tau_1,\tau_1+1\}$ and we consider the grand conclusion
where we make $F_c$ the selected factor.\footnote{Recall that $F_c$
is always a factor $\nabla^{(m)}R_{ijkl}$.} We derive:

\begin{equation}
\label{outside}
\begin{split}
&\Sum_{l\in L_1^c} a_l (-q_1-2) C^{l,i_1}_{g}(\Omega_1,\dots
,\Omega_p,\phi_1,\dots ,\phi_u)\nabla_{i_1}\phi_{u+1}+
\\&2\Sum_{b\in \{1,\dots ,\tau_1,\tau_1+1\}} \Sum_{l\in L^b_1}
a_l C^{l,i_1\rightarrow F_c}_{g}(\Omega_1,\dots
,\Omega_p,\phi_1,\dots ,\phi_u)\nabla_{i_1}\phi_{u+1}+
\\&(\sum_{b\in B'} a_b C^{b,i_1}_g(\Omega_1,\dots,\Omega_p,\phi_1,\dots,\phi_u)\nabla_{i_1}\phi_{u+1})
\\&+ \Sum_{h\in H}
a_l Xdiv_{i_2}\dots Xdiv_{i_y} C^{h,i_1\dots
i_y}_{g}(\Omega_1,\dots ,\Omega_p,\phi_1,\dots
,\phi_u)\nabla_{i_1}\phi_{u+1}=
\\&\Sum_{j\in J} a_j C^j_{g}(\Omega_1,\dots ,\Omega_p,\phi_1,\dots
,\phi_{u+1}),
\end{split}
\end{equation}
where the tensor fields indexed in $H$ and complete contractions
indexed in $J$ are as in the conclusion of the grand conclusion.
The linear combination $\sum_{b\in B'}\dots$ (defined in 
Definition \ref{generalni}) arises {\it only} when the 1-tensor fields 
of maximal refined double character in (\ref{hypothese2}) 
are in the form (\ref{tokolpo}). 
Applying Lemma \ref{funny} if necessary,\footnote{The Lemma can be applied, 
since we are assuming that we did not start out with 
the ``special cases'', described  in the beginning of our subsection.} we may
 assume the tensor fields in $H$  are acceptable.

\par Now, we divide by $(-q_1-2)$, and we make
$\nabla\phi_{u+1}$ into an $Xdiv$ and then replace into our Lemma
hypothesis, (\ref{hypothese2}). We see that we are reduced to
proving our Lemma \ref{pskovb} in this case with 
the sublinear combination $\Sum_{l\in L_1} a_l\dots$ replaced by a new
sublinear combination $\Sum_{l\in \underline{L}_1} a_l\dots$ with
certain special features:

{\it Special features:} Let us divide $\underline{L}_1$ as before
into $\underline{L}_1^c$,
 $c\in \{ 1,\dots ,\tau_1,\tau_1+1\}$. Then the sublinear 
combinations $\sum_{l\in \underline{L}_1^c}\dots$ are related as follows:
There exists
 a linear combination of acceptable complete contractions
  $$\Sum_{t\in T} a_t C^t_{g}(\Omega_1,\dots ,\Omega_p,\phi_1,\dots ,\phi_u),$$
  with length $\sigma+u$,
 weight $-n+2$ and $u$-simple character $\vec{\kappa}_{simp}$ so that for each
 $c\in \{ 1,\dots ,\tau_1,\tau+1\}$:

\begin{equation}
\label{thought}
\begin{split}
&\Sum_{l\in \underline{L}^c_1} a_l C^{l,i_1}_{g}(\Omega_1,\dots
,\Omega_p,\phi_1,\dots ,\phi_u)= \nabla^{i_1}_{F^c} \Sum_{t\in T}
a_t C^{t}_{g} (\Omega_1,\dots ,\Omega_p,\phi_1,\dots ,\phi_u),
\end{split}
\end{equation}
where $\nabla^{i_1}_{F^c}$ means that $\nabla^{i_1}$ is forced to
hit the factor $F^c$.

\par Then, applying the grand conclusion with $F^c$ being the
selected factor, we derive:

\begin{equation}
\label{outside2}
\begin{split}
&(-q_1+2(\sigma_1-1))\Sum_{t\in T} a_t
\nabla^{i_1}_{F^c}[C^t_{g}(\Omega_1,\dots ,\Omega_p,\phi_1,\dots
,\phi_u)]\nabla_{i_1}\phi_{u+1}+
\\& \Sum_{h\in H}
a_l Xdiv_{i_2}\dots Xdiv_{i_y} C^{h,i_1\dots
i_y}_{g}(\Omega_1,\dots ,\Omega_p,\phi_1,\dots
,\phi_u)\nabla_{i_1}\phi_{u+1}=
\\&\Sum_{j\in J} a_j C^j_{g}(\Omega_1,\dots ,\Omega_p,\phi_1,\dots
,\phi_{u+1}).
\end{split}
\end{equation}
Thus, if $q_1>2(\sigma_1-1)$,
 since the quantity in parentheses is non-zero,
we only have to apply Lemma \ref{funny} to ensure that all tensor
fields indexed in $H$ are acceptable and we are done.
\newline

\par If $q_1=2(\sigma_1-1)$, we distinguish two further subcases:
Either $\sigma_2+p>0$ or $\sigma_2=p=0$.

\par In the first case, it again follows 
(from the definition of $q_1$ and a counting argument) that for each vector
field in our Lemma hypothesis, the factors
$S_{*}\nabla^{(\nu)}_{r_1\dots r_\nu}R_{ijkl}$ must have all their
indices ${}_{r_1},\dots,{}_{r_\nu},{}_j$ contracting against
factors $\nabla\phi'_h$ and all factors $\nabla^{(B)}\Omega_x$
must have at least two of their indices contracting against
factors $\nabla\phi_h$.

\par Let us consider the first subcase: Making a factor $F_d\ne\nabla^{(m)} R_{ijkl}$
into the selected factor and applying the grand conclusion we
derive an equation:

\begin{equation}
\label{outside3}
\begin{split}
&2\sigma_1\nabla^{i_1}_{F^d}[\Sum_{t\in T} a_t
C^t_{g}(\Omega_1,\dots ,\Omega_p,\phi_1,\dots
,\phi_u)]\nabla_{i_1}\phi_{u+1}+
\\& \Sum_{h\in H}
a_l Xdiv_{i_2}\dots Xdiv_{i_y} C^{h,i_1\dots
i_y}_{g}(\Omega_1,\dots ,\Omega_p,\phi_1,\dots
,\phi_u)\nabla_{i_1}\phi_{u+1}=
\\&\Sum_{j\in J} a_j C^j_{g}(\Omega_1,\dots ,\Omega_p,\phi_1,\dots
,\phi_{u+1}),
\end{split}
\end{equation}
where the tensor fields indexed in $H$ have $\nabla\phi_{u+1}$
contracting against the selected factor $F_c$. As usual, they may
have one unacceptable factor $\nabla\Omega_x$ but then
$\nabla\phi_{u+1}$ will not be contracting against a special
index.

\par Now, by applying Lemma \ref{funny} if necessary,  we may
assume that all tensor fields indexed in $H$ are acceptable and
that $\nabla\phi_{u+1}$ is not contracting against a special index
in any of the tensor fields in (\ref{outside3}).

\par Under all the assumptions above, we apply the eraser to $\nabla\phi_{u+1}$
 to the above (notice that this can be done and produces
 acceptable tensor fields, by virtue of our assumptions)
  and then adding a contracted derivative index
$\nabla^{i_1}$ 
onto the crucial factor $F_1=\nabla^{(m)}R_{ijkl}$, and then
 contacting against a factor $\nabla_{i_1}\phi_{u+1}$, we obtain:

\begin{equation}
\label{outside4}
\begin{split}
& 2\Sum_{t\in T} a_t \nabla^{i_1}_{F^1}[C^t_{g}(\Omega_1,\dots
,\Omega_p,\phi_1,\dots ,\phi_u)]\nabla_{i_1}\phi_{u+1}+
\\& \Sum_{h\in H'}
a_l Xdiv_{i_2}\dots Xdiv_{i_y} C^{h,i_1\rightarrow F_1,i_2 \dots
i_y}_{g}(\Omega_1,\dots ,\Omega_p,\phi_1,\dots
,\phi_u)\nabla_{i_1}\phi_{u+1}=
\\&\Sum_{j\in J} a_j C^j_{g}(\Omega_1,\dots ,\Omega_p,\phi_1,\dots
,\phi_{u+1}).
\end{split}
\end{equation}
This is our desired conclusion, in the case $\sigma_2+p>0$.
\newline

\par Finally, the case where $\sigma_2=p=0$, $q_1=2(\sigma_1-1)$.
 Notice that necessarily
then, in our Lemma hypothesis there are no complete contractions
$C^j_{g}$. Recall (\ref{thought}), which still holds in this case.

Since $q_1=2(\sigma_1-1)$ it follows that each $C^t_{g}$ has the
property that in {\it each} factor $\nabla^{(m)}R_{ijkl}$ {\it
all} the derivative indices must be contracting against a factor
$\nabla\phi_h$. Therefore, each complete contraction in our Lemma
hypothesis\footnote{In other words, we momentarily think of each $Xdiv$ 
in our Lemma hypothesis as a sum of complete contractions.} must have precisely two derivative indices among all
the factors $\nabla^{(m)}R_{ijkl}$ that are not contracting
against a factor $\nabla\phi_x$ (this follows by weight
considerations).

 Also, the only tensor fields appearing in
(\ref{hypothese2}), other than the ones indexed in $L_1$
 can have rank 2 (recall that this sublinear combination of these
2-tensor fields is denoted by $\Sum_{l\in L_2} a_l
C^{l,i_1i_2}_{g}$), and they must be in one of two special forms:

 \begin{enumerate}
\item{ {\it Either} $C^{l,i_1i_2}_{g}$ will have two factors
$\nabla^{(m)} R_{ijkl}$ with the indices ${}_i$ being free and all
derivative indices contracting against some factor $\nabla\phi_x$,
and furthermore all the other $\sigma-2$ factors
$\nabla^{(m)}R_{ijkl}$ must have no free indices and each
derivative index contracting against some factor $\nabla\phi_h$.}
\item{ {\it Or} $C^{l,i_1i_2}_{g}$ will have both free indices
being internal non-antisymmetric indices in some factor
$\nabla^{(m)}R_{ijkl}$, and all derivative indices in all 
the other $\sigma-1$ factors
$\nabla^{(m)}R_{ijkl}$ are contracting against factor
$\nabla\phi_x$. }
\end{enumerate}

 Let us prove the claim in this case: Denote by $L_{2,diff}\subset L_2$ the index set of
 2-tensor fields where the two free indices ${}_{i_1},{}_{i_2}$ belong to different factors. 
 Let $C^{l,i_1i_2}_g(\phi_1,\dots,\phi_u)   g_{i_1i_2}$ be the complete contraction that arises from 
 $C^{l,i_1i_2}_g$ by contracting the two indices ${}_{i_1},{}_{i_2}$ against each other.

 Now, given a complete contraction $C^t_g(\phi_1,\dots,\phi_u)$ in the form (\ref{form2}) 
 with a $u$-simple character $\vec{\kappa}_{simp}$, we denote by $\tau[t]$ 
 the number of pairs of factors $(T_a,T_b)$ with at least one particular contraction between them. 
 
 We claim:
   
 \begin{equation}
 \label{mainclaim}
\sum_{l\in L_{2,diff}} a_l C^{l,i_1i_2}_g(\phi_1,\dots,\phi_u)   g_{i_1i_2}=
\sum_{t\in T} a_t \tau[t] C^t_g(\phi_1,\dots,\phi_u).   
\end{equation}

Let us check how the above implies our claim: 
\newline

{\it Proof that (\ref{mainclaim}) implies our claim:} We will prove this inductively: 
Let $M$ be the maximum value of $\tau[t], t\in T$. Clearly $M>0$ (simply because the exist
a pair of indices in two different factors that contract against each other, by construction). 
Denote by $T_M\subset T$
 the corresponding index set.  We will then prove:

\begin{equation}
\label{the.step}
\sum_{t\in T_M} a_t \nabla_{i_1}C^t_g\nabla^{i_1}\phi_{u+1}=
\sum_{h\in H} a_h Xdiv_{i_2} C^{h,i_1i_2}_g\nabla_{i_1}\phi_{u+1}
+\sum_{t\in T'} a_t \nabla_{i_1}C^t_g\nabla^{i_1}\phi_{u+1};
\end{equation}
here the terms indexed in $T'$ are generic complete contractions in the
 form (\ref{form2}) with a $u$-simple character $\vec{\kappa}_{simp}$
 and weight $-n+2$ and moreover $\tau[t]<M$. The 2-tensor fields in  $\sum_{h\in H} a_h \dots$
 are as described in the claim of Lemma \ref{pskovb}. If we can prove this, then by 
 replacing the $\nabla\phi_{u+1}$ by an $Xdiv$,\footnote{(Using the 
 last Lemma in the Appendix of  \cite{alexakis1}).} replacing back 
 into our Lemma hypothesis and then iterating this step at most four times, 
 we derive our claim. 
 
 Thus, matters are reduced to proving (\ref{the.step}):

{\it Proof of (\ref{the.step}):} Refer to equation (\ref{proolaxreiazontai2}).
 In view of (\ref{mainclaim}), we derive an equation: 

\begin{equation}
\label{beig}
\sum_{t\in T_M} M\cdot a_t \nabla_{i_1}C^t_g\nabla^{i_1}\phi_{u+1}
=\sum_{h\in H} a_h Xdiv_{i_2} C^{h,i_1i_2}_g
-\sum_{t\in T\setminus T_M} a_t\cdot \tau[t] \nabla_{i_1}C^t_g\nabla^{i_1}\phi_{u+1}.
\end{equation}
Thus, dividing by $M$ in the above, we derive (\ref{the.step}). $\Box$
\newline

{\it Proof of (\ref{mainclaim}):} For each $C^{l,i_1i_2}_g, l\in L_{2,diff}$, 
let us denote by $T_{i_1}, T_{i_2}$ the two 
factors $\nabla^{(m)}R_{ijkl}, \nabla^{(m'} R_{i'j'k'l'}$ 
to which the indices ${}_{i_1},{}_{i_2}$ belong. 
\begin{definition}
 \label{definition}
We assume wlog that they occupy the positions ${}_i,{}_{i'}$. 
Let us denote by $L_{2,diff}^\alpha$, 
$L_{2,diff}^\beta$, $L_{2,diff}^\gamma$, $L_{2,diff}^\delta$ the index sets
 of tensor fields for which the 
two factors $T_{i_1}, T_{i_2}$ in $C^{l,i_1i_2}_g$ have three, two, one and no particular
 contractions between them, respectively. Now, given any 
$C^{l,i_1i_2}_g$, $l\in L_{2,diff}$, 
 we denote by $\hat{C}^l_g$ the complete contraction that arises by
 contracting ${}_{i_1}$ against ${}_{i_2}$ and then hitting 
 the factors $T_{i_1}, T_{i_2}$ by derivatives $\nabla_{i_*}, \nabla^{i_*}$ 
 (which contract against each other). 
\end{definition} 
 
 Now, given a complete contraction $C^t_g$, $t\in T$, let us denote by $Pair_t^\alpha$,
 $Pair_t^\beta$, $Pair_t^\gamma$, $Pair_t^\delta$ the set  of (ordered) pairs of factors $(T_a,T_b)$ in $C^t_g$ 
 that have four, three, two and one particular contractions between them, respectively. 
 Given any $C^t_g, t\in T$ and any pair of factors $(T_a,T_b)\in Pair_t^\alpha$ etc,
  we denote by $Hit_{T_a,T_b}[C^t_g]$ the complete 
contraction in $Xdiv_{i_1}\nabla_{i_1}C^t_g$ where the derivative
 $\nabla_{i_1}$ is forced to hit $T_a$ and then the derivative 
 $\nabla^{i_1}$ in $Xdiv_{i_1}$ is forced to hit $T_b$. 

  We will then prove that:
 
 \begin{equation}
 \label{the.technicalalpha}
 \sum_{l\in L_{2,diff}^\alpha} a_l \hat{C}^l_g=\sum_{t\in T} a_t 
\sum_{(T_a,T_b)\in Pair_t^\alpha} Hit_{T_a,T_b}[C^t_g], 
 \end{equation}

 \begin{equation}
 \label{the.technicalbeta}
 \sum_{l\in L_{2,diff}^\beta} a_l \hat{C}^l_g=\sum_{t\in T} 
a_t \sum_{(T_a,T_b)\in Pair_t^\beta} Hit_{T_a,T_b}[C^t_g], 
 \end{equation}

 \begin{equation}
 \label{the.technicalgamma}
 \sum_{l\in L_{2,diff}^\gamma} a_l \hat{C}^l_g=\sum_{t\in T} 
a_t \sum_{(T_a,T_b)\in Pair_t^\gamma} Hit_{T_a,T_b}[C^t_g], 
 \end{equation}
 
 \begin{equation}
 \label{the.technicaldelta}
 \sum_{l\in L_{2,diff}^\delta} a_l \hat{C}^l_g=\sum_{t\in T} a_t \sum_{(T_a,T_b)\in Pair_t^\delta} Hit_{T_a,T_b}[C^t_g]. 
 \end{equation}
 
 We note that if we can prove the above, then by just applying the eraser to the pair of 
 contracting  derivative indices $\nabla^{i_*},\nabla_{i_*}$\footnote{Note that by weight 
considerations there is exactly one such pair.} in 
 all four of the above equations and then adding, 
 we derive (\ref{mainclaim}).

Let us first prove  (\ref{the.technicalalpha}), which contains (in a 
simple form) the idea for the proof of the other equations. We will then prove
(\ref{the.technicalbeta}), and explain how the other two 
equations follow by the same argument. 
 \newline

{\it Proof of (\ref{the.technicalalpha}):} Pick out the sublinear combination in (\ref{hypothese2}) with two 
factors $\nabla_a R_{ijkl}, \nabla_{a'} R_{i'j'k'l'}$,\footnote{We are again not writing out the indices 
in those factors that contract against factors $\nabla\phi_h$). }  with 
{\it five} particular contractions between them. It is clear
 that this sublinear combination must vanish separately, and 
that the only terms in $\sum_{t\in T} a_t Xdiv_{i_1} C^{t,i_1}_g$ 
that can contribute to this sublinear combination 
from are precisely the terms $\sum_{t\in T} a_t \sum_{(T_a,T_b)\in Pair_t^\alpha} Hit_{T_a,T_b}[C^t_g]$.
 From $\sum_{l\in L_{2,diff}} a_l Xdiv_{i_1}Xdiv_{i_2} C^{l,i_1i_2}_g$, 
 the only terms that can contribute to this are the 
 tensor fields indexed in $ L_{2,diff}^\alpha$, 
 if we {\it force} the derivative $\nabla^{i_1}$
 to hit $T_{i_2}$ and $\nabla^{i_2}$ to hit $T_{i_1}$.  
Our claim then just follows by the second Bianchi identity. $\Box$
\newline

{\it Proof of (\ref{the.technicalbeta}), (\ref{the.technicalgamma}),
 (\ref{the.technicaldelta}):} Pick out the sublinear combination
 in (\ref{hypothese2}) with two factors $\nabla_aR_{ijkl}, \nabla_{a'}R_{i'j'k'l'}$
 with {\it four} particular contractions between them. Denote this 
 sublinear combination (which clearly must vanish separately) by $Y_g=0$. 
 Let us also observe  (by virtue of the Bianchi identities and the 
 anti-symmetries of the curvature tensor) that the four particular contractions 
 can either be according to the pattern $\nabla_aR_{ijkl}\nabla_{a'}R^{ijkl}$,
 {\it or} according to the pattern $\nabla_aR_{ijkl}\nabla^a{R_{i'}}^{jkl}$; denote these 
 two sublinear combination by $Y^1_g, Y^2_g$. It follows (using the variation
  of the curvature tensor by a symmetric 2-tensor $v_{ij}$), that $Y^1_g=0, Y^2_g=0$. 

Now, notice that the only terms in $\sum_{t\in T} a_t Xdiv_{i_1} C^{t,i_1}_g$ 
that can contribute to $Y^2_g$ are precisely the terms 
$\sum_{t\in T} a_t \sum_{(T_a,T_b)\in Pair_t^\alpha} Hit_{T_a,T_b}[C^t_g]$.
\\ From $\sum_{l\in L_{2,diff}} a_l Xdiv_{i_1}Xdiv_{i_2} C^{l,i_1i_2}_g$, 
 the only terms that can contribute to this are the 
 tensor fields indexed in $ L_{2,diff}^\alpha$, 
 if we {\it force} the derivative $\nabla^{i_1}$
 to hit $T_{i_2}$ and $\nabla^{i_2}$ to hit $T_{i_1}$, and then apply the second Bianchi identity to 
 replace any expression $\nabla^{i'}R_{ijkl}\nabla^i{{R^{i'}}_{j'}}^{k'l'}$   
  by $\nabla^sR_{ijkl}\nabla_s{{R^i}_j}^{kl}$ (the {\it other term}
that arises in the second Bianchi identity falls into $Y^1_g$), or 
  $\nabla^{i'}R_{ijkl}\nabla^i{{R_{i'}}^{jk}}_{l'}$ (no correction term in this case). 

(\ref{the.technicalgamma}), (\ref{the.technicaldelta}) follow by 
the same reasoning. Thus we derive our claim. $\Box$

We have shown Lemma \ref{pskovb} when $\mu=1$, except for the special 
cases described at the beginning of this subsection. We derive Lemma \ref{pskovb}
in those special settings in the Mini-Appendix below:

\subsection{Mini-Appendix: Proof of Lemma \ref{pskovb} when 
$\mu=1$, in the ``special subcases''. }

We now prove our claim in the special subcases, defined at the beginning
 of subsection (\ref{mu=1}); we recall that  $\sigma_2>0$ in 
the special subcases. We recall 
that by virtue of the assumption $L^+_\mu=\emptyset$ in our Lemma assumption, 
 we may assume wlog that  
all tensor fields of minimum rank 1 in (\ref{hypothese2}) contain no free index in any 
factor of the form $S_*\nabla^{(\nu)}R_{ijkl}$. We first consider the 
special subcases where the `-tensor fields of maximal refined double character 
 in (\ref{hypothese2}) correspond to (\ref{tokolpo}):

{\it The second special case,  (\ref{tokolpo}):}  We 
recall also that no $2$-tensor field in (\ref{hypothese2}) 
can contain more than one free index in any simple  factor 
$S_*\nabla^{(\nu)}R_{ijkl}$, by weight 
considerations. Again by weight considerations, we derive that 
the maximum rank of tensor fields appearing in (\ref{hypothese2}) is 2, 
and moreover those tensor fields can have 
no removable free index. (In particular any given 
simple factor $S_*\nabla^{(\nu)}R_{ijkl}$ must have $\nu=0$). Moreover, 
for any $\mu$-tensor field in (\ref{hypothese2}) and any given 
factor $S_*\nabla^{(\nu)}R_{ijkl}\nabla^i\tilde{\phi}_h$ can  
have $\nu=0$ or $\nu=1$. 

Now, denote by $L_1^x\subset L_1$ the index set of $1$-tensor fields in 
(\ref{hypothese2}) with $\nu=0$ on the factor 
$S_*\nabla^{(\nu)}R_{ijkl}\nabla^i\tilde{\phi}_1$. 
We denote by $L^q_1\subset L_1$  the index set of $\mu$-tensor fields in 
(\ref{hypothese2}) with $\nu=1$ in the factor 
$S_*\nabla^{(\nu)}R_{ijkl}\nabla^i\tilde{\phi}_1$.
By our remark above, $L_1^x\bigcup L^q_1=L_1$. 

We will prove the following: First, that 
there exists a linear combination of $2$-tensor
 fields (indexed in $H'$ below), as required by Lemma \ref{pskovb}, such that:
\begin{equation}
\label{riberi}
\begin{split}
&\sum_{l\in L^x_1} a_l C_g^{l,i_1}\nabla_{i_1}\upsilon-Xdiv_{i_2}\sum_{h\in H'} a_h 
C^{h,i_1i_2}_g\nabla_{i_1}\upsilon=
\\&\sum_{l\in L_1'} a_l C^{l,i_1}_g\nabla_{i_1}\upsilon
+\sum_{j\in J} a_j C^{j,i_1}_g\nabla_{i_1}\upsilon.
\end{split}
\end{equation}
Here $\sum_{l\in L_1'} a_l C^{l,i_1}_g$ stands for a linear combination of acceptable 
$\mu$-tensor fields in the  form (\ref{form2}) with 
a $u$-simple character $\vec{\kappa}_{simp}$, and in addition have $\nu=1$  
in the factor $S_*\nabla^{(\nu)}R_{ijkl}\nabla^i\tilde{\phi}_1$. 
The terms indexed in $J$ are simply subsequent to $\vec{\kappa}_{simp}$. 
The above holds modulo terms of length $\ge\sigma+u+2$. 
Thus, making the $\nabla\upsilon$ into an $Xdiv$ and replacing 
into (\ref{hypothese2}), we may assume wlog that $L^x_\mu=\emptyset$.

We then will prove that:  

\begin{equation}
\label{riberi2}
\begin{split}
&\sum_{h\in H'\bigcup L_{>1}} a_h 
C^{h,i_1i_2}_g\nabla_{i_1}\upsilon\nabla_{i_2}\upsilon=0.
\end{split}
\end{equation}
(Here the terms indexed in $L_{>1}$ are {\it the same} terms 
appearing in (\ref{hypothese2}); the terms indexed in $H'$ are the same that appear in (\ref{riberi}). 
Thus, making the $\nabla\upsilon$ into an $Xdi$ and replacing 
into (\ref{hypothese2}), we may additionaly assume wlog that 
$H\bigcup H'=\emptyset$. 
Finally, under that additional assumption we will show that: 

\begin{equation}
\label{riberi3}
\sum_{l\in L^x_1} a_l C^{l,i_1}\nabla_{i_1}\upsilon=0.
\end{equation}

Clearly, if we can show the above three equations, then our claim will follow. 

The proof of the above equation just relies on the trick that we introduced 
to prove (in \cite{alexakis5}) the   Lemma 5.1 in \cite{alexakis4}: 

Pick out the sublinear combination in (\ref{hypothese2}) with an 
(undifferentiated) factor $S_*R_{ijkl}\nabla^i\tilde{\phi}_1$. 
This sublinear combination must vanish separately, thus we derive:

\begin{equation}
 \label{step}
\sum_{l\in L^x_1} a_l X_*div_{i_1}C^{l,i_1}
-X_*div_{i_1}X_*div_{i_2}\sum_{h\in H} a_h C^{h,i_1i_2}_g+\sum_{j\in J} a_j C^j_g=0.
\end{equation}
Now, we denote by $H_a\subset L_{>1}$ the index set of 2-tensor fields 
in the above that contain a free index, say the index ${}_{i_2}$ wlog,  in the 
factor $S_*R_{ijkl}\nabla^i\tilde{\phi}_1$. We also 
formally replace the factor $S_*R_{ijkl}$ by an 
expression $\nabla_j\omega\nabla_k\omega\nabla_l y-\nabla_j\omega\nabla_l\omega\nabla_k y$, 
thus deriving a new true equation. Denote the resulting 
tensor fields by $\tilde{C}^{l,i_1},\tilde{C}^{h,i_1i_2}_g$. 
We then apply 
Lemma 5.1 in \cite{alexakis4} to the resulting equation. We derive  that:

\begin{equation}
\label{ileana}\begin{split}
&\sum_{l\in L^x_1} a_l \tilde{C}^{l,i_1}\nabla_{i_1}\upsilon
+X_*div_{i_2}\sum_{h\in H_a} a_h \tilde{C}^{h,i_1i_2}_g
\nabla_{i_1}\upsilon=X_*div_{i_2}\sum_{h\in H'} a_h 
C^{h,i_1i_2}_g\nabla_{i_1}\upsilon
\\&+\sum_{j\in J} a_j C^j_g.
\end{split}
\end{equation}
Here the terms indexed in $H'$ are of the same form as the terms 
in the LHS, but have rank 2; furthermore, the free indices
 ${}_{i_1},{}_{i_2}$ {\it do not} belong to any of 
the factors $\nabla\omega,\nabla y$. Now, by just formally replacing the expression 
$\nabla_a\omega\nabla_b\omega\nabla_cy$ by $S_*R_{i(ab)c}\nabla^i\tilde{\phi}_1$
 (this gives a new true equation), we derive (\ref{riberi}). (\ref{riberi2}) 
follows by the exact same argument. (\ref{riberi3}) follows by the same argument, 
the only difference being that we now pick out the terms with a 
factor $\nabla_dS_*R_{ijkl}\nabla^i\tilde{\phi}_1$ and replace them by an 
expression 
$\nabla_d\omega\nabla_j\omega\nabla_k\omega\nabla_ly-\nabla_d\omega\nabla_j\omega\nabla_l\omega\nabla_ky$.
We then apply the  Lemma 5.1 in \cite{alexakis4},\footnote{See 
the Appendix of \cite{alexakis4}},
 and in the end formally replace each expression 
$\nabla_d\omega\nabla_j\omega\nabla_k\omega\nabla_ly$ by a factor 
$\nabla_{(d}S_*{R^i}_{jk)l}\nabla_i\tilde{\phi}_1$. This proves (\ref{riberi3}). $\Box$ 
\newline

We now consider the special cases where the 1-tensor fields of
maximal refined double character in (\ref{hypothese2}) have no removable free indices: 

{\it The first special case:}

We start by observing that (by weight considerations) 
there can be no tensor fields of rank higher than $1$ in (\ref{hypothese2}), and that $\sigma_2>0$. 

 In this case, the claim of Lemma \ref{pskovb} is a straightforward 
application of the ``generalized form'' of Lemma 5.1 in \cite{alexakis4}:
 Consider a factor $S_{*}R_{ijkl}\nabla^i\tilde{\phi}_x$ in $\vec{\kappa}_{simp}$. 
Wlog we asume that $x=1$. By virtue of the assumption $L_\mu^+=\emptyset$, 
we know that this factor does not contain a free index 
for any of the tensor fields in (\ref{hypothese2}).  
 We pick out the sublinear combination 
of terms with an (undifferentiated)  factor $S_*R_{ijkl}\nabla^i\tilde{\phi}_1$ 
in (\ref{hypothese2}). This sublinear combination must clearly vanish 
separately, thus we derive that: 

\begin{equation}
\label{ovidios}\begin{split}
&\sum_{l\in L_\mu} a_l X_*div_{i_1}\dots X_*div_{i_\mu} 
C^{l,i_1\dots i_\mu}_g(\Omega_1,\dots,\Omega_p,\phi_1,\dots,\phi_u)\\&+
\sum_{j\in J} a_j C^j_g(\Omega_1,\dots,\Omega_p,\phi_1,\dots,\phi_u)=0. 
\end{split}\end{equation}
(Here $X_*div_{i_s}[\dots]$ stands for the sublinear combination 
in $Xdiv_{i_s}[\dots]$ where the derivative $\nabla^{i_s}$ is not 
allowed to hit the factor $S_{*}R_{ijkl}\nabla^i\tilde{\phi}_1$; 
notice that none of the free indices belong to 
$S_*R_{ijkl}\nabla^i\tilde{\phi}_1$).

Now, let $\tilde{C}^{l,i_1\dots i_\mu}_g$ 
stand for the tensor field that formally 
arises from each $C^{l,i_1\dots i_\mu}_g$ 
by replacing the factor $S_*R_{ijkl}\nabla^i\tilde{\phi}_1$ 
by an expression $\nabla_j\omega\nabla_k\omega\nabla_l y-\nabla_j\omega\nabla_l\omega\nabla_jy$. 
We perform this formal replacement to the terms in (\ref{ovidios}), 
and then we apply Lemma 5.1 in \cite{alexakis4} to the resulting equation. We derive that:
$$ \sum_{l\in L_\mu} a_l  
\tilde{C}^{l,i_1\dots i_\mu}_g(\Omega_1,\dots,\Omega_p,\phi_1,\dots,\phi_u)\nabla_{i_1}\upsilon\dots\nabla_{i_\mu}\upsilon=0.$$
Thus, formally replacing each expression 
$\nabla_j\omega\nabla_k\omega\nabla_l y$ by $S_*R_{i(jk)l}\nabla^i\tilde{\phi}_1$,
 we derive our claim in this remaining special case. $\Box$

\section{Appendix: Proof of Lemmas 3.3, 3.4 in \cite{alexakis4}.}
\label{appendixc}

{\bf Definitions related to Lemmas 3.3, 3.4 in \cite{alexakis4}:}
 In case A of Lemma \ref{pskovb} we denote by
$L_\mu^{*}\subset L_\mu$ the index set of those tensor fields
$C^{l,i_1\dots i_\mu}_{g}$ in (\ref{hypothese2}) for which some
 factor $\nabla^{(A)}_{r_1\dots r_A}\Omega_x$ (for a single $x$,
which we are free to define) has $A=2$ and both
 indices ${}_{r_1},{}_{r_2}$ are free indices.

\par Also,  we define $L_\mu^{+}\subset L_\mu$ to
 stand for the index
set of those $\mu$-tensor fields that have a free
 index (${}_{i_\mu}$ say) belonging to a
factor $S_{*}R_{ijkl}\nabla^i \tilde{\phi}_h$ (without
derivatives).

\par We now consider the sublinear combination indexed in
$L_{>\mu}$ in (\ref{hypothese2}). We  define
$L''_{+}\subset L_{>\mu}$ to stand for the index set of $(\mu+1)$-tensor
fields with a
factor $S_*R_{ijkl}\nabla^i\tilde{\phi}_h$ with {\it both} indices
${}_j,{}_k$ free. 
\newline

{\bf The Lemmas proven in the present paper:} 
After all these definitions, we are prepared to re-state 
Lemmas 3.3, 3.4 in \cite{alexakis4}:

\begin{lemma}
\label{oui} Assume (\ref{hypothese2}), where the terms in the LHS of that equation
have weigh $-n$, real length $\sigma$, $\Phi$ factors 
$\nabla\phi, \nabla\phi',\nabla\tilde{\phi}$ and $\sigma_1+\sigma_2$ curvature 
factors $\nabla^{(m)}R_{ijkl},S_{*}\nabla^{(\nu)}R_{ijkl}$;\footnote{See the discussion
on the {\it induction} in the introduction.}  assume also  that no
$\mu$-tensor field there has any special free indices.
 We claim that there is a linear combination of
acceptable $(\mu+1)$-tensor fields,  $\Sum_{p\in P} a_p
C^{p,i_1\dots i_{\mu +1}}_{g}(\Omega_1,\dots ,\Omega_p,
\phi_1,\dots ,\phi_u)$ with a simple character
$\vec{\kappa}_{simp}$ so that:

\begin{equation}
\label{dog}
\begin{split}
&\Sum_{l\in L_\mu^{*}\bigcup L_\mu^{+}} a_l C^{l,i_1\dots
i_\mu}_{g}(\Omega_1,\dots ,\Omega_p, \phi_1,\dots
,\phi_u)\nabla_{i_1}\upsilon\dots \nabla_{i_\mu}\upsilon+
\\& \Sum_{p\in P} a_p Xdiv_{i_{\mu +1}}
C^{p,i_1\dots i_{\mu +1}}_{g}(\Omega_1,\dots ,\Omega_p,
\phi_1,\dots ,\phi_u)\nabla_{i_1}\upsilon\dots
\nabla_{i_\mu}\upsilon=
\\& \Sum_{j\in J} a_j C^{j,i_1\dots i_\mu}_{g}(\Omega_1,\dots ,\Omega_p,
\phi_1,\dots ,\phi_u)\nabla_{i_1}\upsilon\dots
\nabla_{i_\mu}\upsilon +
\\& \Sum_{l\in \tilde{L}} a_l C^{l,i_1\dots
i_\mu}_{g} (\Omega_1,\dots ,\Omega_p,\phi_1,\dots ,\phi_u)
\nabla_{i_1}\upsilon\dots\nabla_{i_\mu}\upsilon,
\end{split}
\end{equation}
modulo complete contractions of length $\ge\sigma +u+\mu +1$. The
tensor fields indexed in $J$ on the right hand side are simply
subsequent to $\vec{\kappa}_{simp}$. The terms indexed in $\tilde{L}'$ 
in the RHS are acceptable terms in the form (\ref{form2})
with a simple character $\vec{\kappa}_{simp}$. They are not in the 
forms corresponding to tensor fields indexed in $L^*_\mu\bigcup L^+_\mu$. 
\end{lemma}

\begin{lemma}
\label{oui2} Assume (\ref{hypothese2})  with weight $-n$, real length $\sigma$,
 $\Phi$ factors $\nabla\phi,\nabla\phi',\nabla\tilde{\phi}$ 
and $\sigma_1+\sigma_2$ factors
$\nabla^{(m)}R_{ijkl},S_{*}\nabla^{(\nu)}R_{ijkl}$; assume also that 
 none of the $\mu$-tensor fields have special free indices,
and  that
 $L_\mu^{*}\bigcup L_\mu^{+}=\emptyset$. 
 We claim that there exists a linear combination of
acceptable $(\mu+2)$-tensor fields, \\$\Sum_{p\in P} a_p C^{p,i_1\dots i_{\mu
+2}}_{g}(\Omega_1,\dots ,\Omega_p, \phi_1,\dots ,\phi_u)$ with simple
 character $\vec{\kappa}_{simp}$, so that:

\begin{equation}
\label{dog2}
\begin{split}
&\Sum_{l\in L''_{+}} a_l Xdiv_{i_1}\dots
Xdiv_{i_{\mu+1}} C^{l,i_1\dots i_{\mu+1}}_{g}(\Omega_1,\dots
,\Omega_p, \phi_1,\dots ,\phi_u)+
\\& \Sum_{p\in P} a_p Xdiv_{i_1}\dots Xdiv_{i_{\mu+2}}
C^{p,i_1\dots i_{\mu +2}}_{g}(\Omega_1,\dots ,\Omega_p,
\phi_1,\dots ,\phi_u)=
\\& \Sum_{j\in J} a_j C^{j}_{g}(\Omega_1,\dots ,\Omega_p,
\phi_1,\dots ,\phi_u)\\&+ \Sum_{l\in \tilde{L}'} a_l Xdiv_{i_1}\dots
Xdiv_{i_{\mu+1}} C^{l,i_1\dots i_{\mu+1}}_{g} (\Omega_1,\dots
,\Omega_p, \phi_1,\dots ,\phi_u),
\end{split}
\end{equation}
modulo complete contractions of length $\ge\sigma +u+1$.
$\Sum_{j\in J}\dots$ stands for  a linear combination of complete
contractions that are simply subsequent to $\vec{\kappa}_{simp}$. 
The terms indexed in $\tilde{L}'$ in the RHS are acceptable terms in the form (\ref{form2})
with a simple character $\vec{\kappa}_{simp}$. They are not in the 
forms corresponding to tensor fields indexed in $ L''_+$. 
\end{lemma}

{\it Proof of Lemma \ref{oui}:} We prove Lemma \ref{oui} by breaking it into
 two Lemmas. In Lemma \ref{ouia} below, we aim to ``get rid''
 of the tensor fields indexed in $L^{*}_\mu$. (Recall that
 $L_\mu^{*}\subset L_\mu$ stands for the index set of factors with
two free indices belonging to a factor $\nabla^{(2)}\Omega_x$.

With no loss of generality
(up to
 re-labelling the functions $\Omega_h$) we may assume that
 $x=1$.

\begin{lemma}
\label{ouia} We claim that there is a linear combination of
acceptable $(\mu+1)$-tensor fields,  $\Sum_{p\in P} a_p
C^{p,i_1\dots i_{\mu +1}}_{g}(\Omega_1,\dots ,\Omega_p,
\phi_1,\dots ,\phi_u)$ so that:

\begin{equation}
\label{doga}
\begin{split}
&\Sum_{l\in L_\mu^{*}} a_l C^{l,i_1\dots i_\mu}_{g}(\Omega_1,\dots
,\Omega_p, \phi_1,\dots ,\phi_u)\nabla_{i_1}\upsilon\dots
\nabla_{i_\mu}\upsilon+
\\& \Sum_{p\in P} a_p Xdiv_{i_{\mu +1}}
C^{p,i_1\dots i_{\mu +1}}_{g}(\Omega_1,\dots ,\Omega_p,
\phi_1,\dots ,\phi_u)\nabla_{i_1}\upsilon\dots
\nabla_{i_\mu}\upsilon=
\\& \Sum_{j\in J} a_j C^{j,i_1\dots i_\mu}_{g}(\Omega_1,\dots ,\Omega_p,
\phi_1,\dots
,\phi_u)\nabla_{i_1}\upsilon\dots\nabla_{i_\mu}\upsilon
\\&+ \Sum_{l\in \tilde{L}} a_l
C^{l,i_1\dots i_\mu}_{g} (\Omega_1,\dots ,\Omega_p,\phi_1,\dots
,\phi_u) \nabla_{i_1}\upsilon\dots\nabla_{i_\mu}\upsilon.
\end{split}
\end{equation}
where each $C^{l,i_1\dots i_\mu}_g$ on the RHS is acceptable in the form (\ref{form2}) 
with a $u$-simple character $\vec{\kappa}_{simp}$ has a factor $\nabla^{(A)}
\Omega_1$ with $A\ge 3$. The complete contractions indexed in $J$
are simply subsequent to $\vec{\kappa}_{simp}$. The above holds
modulo complete contractions of length $\ge\sigma +u+\mu +1$.
\end{lemma}

If we can show the above, then  we will be reduced to
showing Lemma \ref{oui} under the assumption that
$L_\mu^{*}=\emptyset$.

\par Our next claim, which will
 ``get rid'' of the sublinear combination indexed in $L^{+}_\mu$:
In fact, we make a stronger claim:

\begin{lemma}
\label{ouic} Assume (\ref{hypothese2}) where no $\mu$-tensor field there has 
special free indices;\footnote{Recall that a 
free index is called ``special'' if it occupies the position ${}_k,{}_l$ in a 
factor $S_{*}\nabla^{(\nu)}R_{ijkl}$ or the position 
${}_i,{}_j,{}_k,{}_l$ in a factor $\nabla^{(m)}R_{ijkl}$. } consider a given simple 
factor $S_{*}\nabla^{(\nu)}R_{ijkl}$.\footnote{Recall that a
 factor $S_{*}\nabla^{(\nu)}R_{ijkl}$ is called ``simple'' when 
it is contracting against no factor $\nabla\phi'_h$ 
in the simple character $\vec{\kappa}_{simp}$.}
 Let $\tilde{L}^+_\mu\subset L^+_\mu$
stand for the index set of tensor fields which contain exactly 
one free index in the selected factor $S_{*}\nabla^{(\nu)}R_{ijkl}$. 
 We claim that
there is a linear combination of acceptable $(\mu+1)$-tensor
fields,  $\Sum_{p\in P} a_p C^{p,i_1\dots i_{\mu
+1}}_{g}(\Omega_1,\dots ,\Omega_p, \phi_1,\dots ,\phi_u)$ so that:

\begin{equation}
\label{dogb}
\begin{split}
&\Sum_{l\in L_\mu^{+}} a_l C^{l,i_1\dots i_\mu}_{g}
(\Omega_1,\dots ,\Omega_p, \phi_1,\dots,\phi_u)
\nabla_{i_1}\upsilon\dots \nabla_{i_\mu}\upsilon+
\\& \Sum_{p\in P} a_p Xdiv_{i_{\mu +1}}
C^{p,i_1\dots i_{\mu +1}}_{g}(\Omega_1,\dots ,\Omega_p,
\phi_1,\dots ,\phi_u)\nabla_{i_1}\upsilon\dots
\nabla_{i_\mu}\upsilon=
\\& \Sum_{j\in J} a_j C^j_{g}(\Omega_1,\dots ,\Omega_p,
\phi_1,\dots ,\phi_u,\upsilon^\mu).
\end{split}
\end{equation}
\end{lemma}

\par Clearly, if we can show the two Lemmas above, then Lemma
\ref{oui}  will follow. Now, Lemma \ref{oui2} will also follow
 by two claims. Recall the definitions of the sets $L''_{+}$ from
  the previous subsection. We then claim:

\begin{lemma}
\label{oui2a} In the notation above, and under the assumption that
 $ L_\mu^{+}=\emptyset$,
 we claim that there exists an acceptable  linear combination of
$(\mu+2)$-tensor fields, $\Sum_{p\in P} a_p C^{p,i_1\dots
i_{\mu+2}}_{g}(\Omega_1,\dots ,\Omega_p, \phi_1,\dots ,\phi_u)$ so
that:

\begin{equation}
\label{dog2a}
\begin{split}
&\Sum_{l\in L''_{+}} a_l Xdiv_{i_1}\dots Xdiv_{i_{\mu+1}}
C^{l,i_1\dots i_{\mu+1}}_{g}(\Omega_1,\dots ,\Omega_p,
\phi_1,\dots ,\phi_u)+
\\& \Sum_{p\in P} a_p Xdiv_{i_1}\dots Xdiv_{i_{\mu+2}}
C^{p,i_1\dots i_{\mu +1}}_{g}(\Omega_1,\dots ,\Omega_p,
\phi_1,\dots ,\phi_u)=
\\& \Sum_{j\in J} a_j C^j_{g}(\Omega_1,\dots ,\Omega_p,
\phi_1,\dots ,\phi_u)\\&+ \Sum_{l\in \tilde{L}'} a_l Xdiv_{i_1}\dots
Xdiv_{i_{\mu+1}} C^{l,i_1\dots i_{\mu+1}}_{g} (\Omega_1,\dots
,\Omega_p, \phi_1,\dots ,\phi_u),
\end{split}
\end{equation}
modulo complete contractions of length $\ge\sigma +u+1$. Here
 each $C^{l,i_1\dots i_{\mu+1}}$ in the 
RHS is acceptable in the form (\ref{form2}), 
with a simple character $\vec{\kappa}_{simp}$, 
 and has no factors $S_{*}R_{ijkl}$ with two free 
 indices.
\end{lemma}

Clearly, if we can show Lemma \ref{oui2a}, then
Lemma \ref{oui2} will follow.

\subsection{Proof of Lemma \ref{ouia}.}

 Observe that by hypothesis $\mu\ge 2$ in (\ref{hypothese2}) in this setting. 
We now prove our claim in all cases except for a ``special subcase'', 
which will be treated at the end of this subsection. The 
``special subcase'' is when $\mu=2$ and the terms indexed in 
$L_\mu^*$ have the two free indices in a factor $\nabla^{(2)}\Omega_h$,
{\it and all other factors have no removable free 
indices}.\footnote{See definition 4.1 in \cite{alexakis4}. 
Notice that by weight considerations, 
if this is true of one term in $L_\mu^*$ 
then it will be true of all of them.} We now prove our claim, assming that 
(\ref{hypothese2}) does not fall under the ``special subcase''. 
\newline

 Wlog, by just re-naming factors, let us assume that $h=1$. 
(So the two free indices belong to the factor $\nabla^{(2)}\Omega_1$). 

\par We then let
$\Sum_{q\in Q} a_q C^q_{g}(\Omega_1,\dots ,\Omega_p,
\phi_1,\dots,\phi_u,\phi_{u+1})$ stand for a generic linear
combination of complete contractions (not necessarily acceptable)
with length $\sigma +u+1$ and a factor $\nabla\phi_{u+1}$, 
contracting against a factor
$\nabla^{(A)}\Omega_1$ with $A\ge 2$.

\par We also denote by
$\Sum_{c\in C} a_c C^{c,i_1\dots i_a}_{g}(\Omega_1, \dots
,\Omega_p,\phi_1,\dots,\phi_u,\phi_{u+1})$ a
 generic linear combination of tensor fields with 
 a $u$-simple character $\vec{\kappa}_{simp}$,
 with precisely one  un-acceptable factor $\nabla\Omega_x$, 
 contracting against a factor $\nabla\phi_{u+1}$ and
$a\ge \mu$. Finally, we let $\sum_{j\in J}
a_jC^j_{g}(\Omega_1,\dots ,\Omega_p,
\phi_1,\dots,\phi_u,\phi_{u+1})$ stand for a linear combination
 of contractions with length $\sigma+u+1$ and simply subsequent to $\vec{\kappa}_{simp}$.

\par We make $F_1=\nabla^{(A)}\Omega_1$ the selected factor
and consider the equation \\$Image^{1,+}_{\phi_{u+1}}[L_{g}]=0$, see 
equation (\ref{proolaxreiazontai3}).
{\it Note: By the discussion regarding the derivation of 
the super divergece formula above, this
equation can always be applied to the equation $L_{g}=0$--we do
not need any assumptions regarding the sublinear
 combinations indexed in  $L_\mu^{*},L^{+}_\mu$, etc.} By our
choice of selected factor, all contractions in
$Image^{1,+}_{\phi_{u+1}}[L_{g}]$ will contain at least
$\sigma+u+1$ factors.
  For each $l\in L_\mu^{*}$ we denote by $F_2,\dots F_\sigma$ the real
factors other than $F_1=\nabla^{(A)}\Omega_1$.  We thus derive:

\begin{equation}
\label{garding}
\begin{split}
&\Sum_{l\in L^{*}_\mu} a_l Xdiv_{i_3}\dots Xdiv_{i_\mu}
Xdiv_{i_{*}} \Sum_{k=2}^{\sigma} C^{l,i_3\dots
i_{\mu},k(i_{*})}_{g} (\Omega_2,\dots
,\Omega_p,\phi_1,\dots,\phi_u)
\\& \nabla^{i_1}\Omega_1\nabla_{i_1}\phi_{u+1}+\Sum_{c\in C} a_c Xdiv_{i_1}\dots Xdiv_{i_a}
C^{c,i_1\dots i_a}_{g}(\Omega_1,\dots ,\Omega_p,\phi_1,
\dots,\phi_u,\phi_{u+1})
\\&+\Sum_{q\in Q} a_q C^q_{g}(\Omega_1,\dots ,\Omega_p,
\phi_1,\dots,\phi_u,\phi_{u+1})
\\&+\Sum_{j\in J} a_j C^j_{g}(\Omega_1,\dots ,\Omega_p,
\phi_1,\dots,\phi_u,\phi_{u+1})=0,
\end{split}
\end{equation}
modulo complete contractions of length $\ge\sigma +u+2$. Here each
tensor field $C^{l,i_3\dots i_{\mu},k(i_{*})}_{g}$
 arises from $C^{l,i_1\dots i_\mu}_{g}$ by
erasing the factor $F_1=\nabla^{(2)}_{i_1i_2}\Omega_1$ and hitting
the $k^{th}$ factor $F_k$
 by a derivative index $\nabla_{i_{*}}$. In the above
equation we pick out the sublinear combination of
contractions containing an expression $\nabla_i\Omega_1\nabla^i\phi_{u+1}$.
This sublinear combination must vanish separately, since (\ref{garding}) holds formally.
 Denote the resulting true equation by $S_{g}=0$. Now, we formally erase the
expression $\nabla_i\Omega_1\nabla^i\phi_{u+1}$ in $S_{g}$, and we
obtain a new true equation:

\begin{equation}
\label{garding8}
\begin{split}
&\Sum_{l\in L^{*}_\mu} a_l Xdiv_{i_3}\dots Xdiv_{i_\mu}
Xdiv_{i_{*}} \Sum_{k=1}^{\sigma-1} C^{l,i_1\dots
i_{\mu},k(i_{*})}_{g} (\Omega_2,\dots
,\Omega_p,\phi_1,\dots,\phi_u)
\\&+ \Sum_{c\in C} a_c Xdiv_{i_1}\dots Xdiv_{i_a}
C^{c,i_1\dots i_a}_{g}(\Omega_2,\dots ,\Omega_p,\phi_1,
\dots,\phi_u)+
\\&\Sum_{j\in J} a_j C^j_{g}(\Omega_2,\dots ,\Omega_p,
\phi_1,\dots,\phi_u)=0,
\end{split}
\end{equation}
which holds modulo complete contractions of length $\ge\sigma +u$.
Here $C^{c\dots}$ has arisen from the previous equation by just
erasing the expression $\nabla_i\Omega_1\nabla^i\phi_{u+1}$. We
denote by $\vec{\kappa}_{simp}'$ the $(u-1)$-simple
character of those complete
 contractions.  The contractions
$C^j$ in this setting have at least one factor $\nabla\phi_h, h\in
Def(\vec{\kappa}_{simp})$ contracting
 against a derivative index.
\newline

\par Now, we prove Lemma \ref{ouia} in pieces.
 Consider the $\mu$-tensor fields in $L^{*}_\mu$.
We subdivide $L^{*}_\mu$ into subsets
 $L^{*,z}_\mu, z\in Z$ so that two tensor fields
$C^{l_1,i_1\dots i_\mu}, C^{l_2,i_1\dots i_\mu}$ with the same
refined double character $\vec{\kappa}^z_{ref-doub}$ will be
indexed in the same index set $L^{*,z}_\mu$.
 We let $Z_{Max}\subset Z$
 stand for the index sets corresponding to the maximal refined double characters. Suppose
 $M(\ge 1)$ is the maximum number of free indices that can appear in a given factor $F_d\ne F_1$
  among all the tensor fields indexed in $L^{*}_\mu$. Then, by the definition
on the maximal refined double character in \cite{alexakis4}, $M$
  will also be the maximum number of free indices that can appear in a given factor
among all the tensor fields indexed in $\bigcup_{z\in
Z_{Max}}L^{*,z}_\mu$.
  We denote by $cut(\vec{\kappa})^{z}_{ref-doub}$ the refined double
  character that formally arises from $\vec{\kappa}^{z}_{ref-doub}$
  by erasing the entry that corresponds to
  $\nabla^{(2)}_{i_1i_2}\Omega_1$.
 Now, among all the factors $F_d,d\ne 1$ in all the tensor fields indexed in
 $\bigcup_{z\in Z_{Max}}L^{*,z}_\mu$ that have $M$ free indices, we pick out one
 (or a category of generic factors $\nabla^{(m)}R_{ijkl}$)
 canonically, using the same method that was used to choose
 the critical factor in \cite{alexakis4}. We call that (set of) factors the
(set of) $\alpha$-factor(s). We index in $Z'_{Max}\subset Z_{Max}$
 the set of maximal refined double characters $\vec{\kappa}^z_{ref-doub}, z\in Z_{Max}$
that have $M$ free indices in the (an) $\alpha$-factor.
   Now, for each $z\in Z'_{Max}$, we denote by $\tilde{\kappa}^{z}_{ref-doub}$
 the refined double character that arises from
$cut(\vec{L})^{*,z}_\mu$ by formally adding a (derivative) free index $\nabla_{i_{*}}$
onto the (one of the) $\alpha$-factor(s).

 \par We consider (\ref{garding8}) and we observe that
the maximal refined double characters among the tensor fields in Erase[(\ref{garding8})]
 will be $\vec{\kappa}^z_{ref-doub}$.
Now,  assume with no loss of generality that
 the $\alpha$-factor(s) is (are) $F_2$
 ($F_2,\dots ,F_d$). Then, in
 (\ref{garding8}) the $(\mu-1)$-tensor fields with
 $(M+1)$ free indices on the (one of the) $\alpha$-factor(s)
 will be precisely the sublinear combination
 in the first line with $k=2$ (or $k=2,\dots, d$).
    Therefore, applying our inductive assumption
 of Proposition \ref{giade},\footnote{Since all the tensor fields of minimum rank in (\ref{garding8})
 have all free indices being  non-special, there is no danger of falling under
 a ``forbidden case'' of that Proposition.} we derive that for
 each $z\in Z'_{Max}$ there is a linear combination of acceptable
 $\mu$-tensor fields with a $(\mu-1)$-refined double
 character $\tilde{L}^{*,z}_\mu$, indexed in $P$ below so that:

 \begin{equation}
\label{garding9}
\begin{split}
&\Sum_{l\in L^{*}_\mu} a_l \Sum_{k=1}^d C^{l,i_1\dots
i_{\mu},k(i_{*})}_{g} (\Omega_2,\dots
,\Omega_p,\phi_1,\dots,\phi_u)\nabla_{i_3}
\upsilon\dots\nabla_{i_{*}}\upsilon+
\\& \Sum_{p\in P} a_p Xdiv_{i_\mu}
C^{p,i_1\dots i_{\mu-1}i_\mu}_{g}(\Omega_2,\dots ,\Omega_p,
\phi_1,\dots,\phi_u)\nabla_{i_1}\upsilon\dots
\nabla_{i_{\mu-1}}\upsilon=
\\&\Sum_{j\in J} a_j C^j_{g}(\Omega_2,\dots ,\Omega_p,
\phi_1,\dots,\phi_u,\upsilon^{\mu-1}).
\end{split}
\end{equation}

Now, in the case $\mu>2$, we prove our claim by our standard formal manipulation of the above equation:
In that case, we identify in each of the
 contractions above the one factor that is  contracting against
 the most factors $\nabla\upsilon$-it will be
the $\alpha$-factor which will be contracting against $M+1\ge 2$ factors
 $\nabla\upsilon$. Then, we erase one of the factors $\nabla\upsilon$ that is
contracting against the $\alpha$-factor (thus obtaining a new true
equation) and we multiply the new true equation by
$\nabla_{ij}\Omega_1\nabla^i \upsilon\nabla^j\upsilon$).\footnote{Since $\mu>2$ 
it follows that $M+1>1$; thus since all factors $\nabla\upsilon$ contract 
against non-special indices, at least one of the $M+1$ factors
 $\nabla\upsilon$ will be a derivative index, hence the Eraser can be applied.} This
further true equation is precisely the claim of Lemma \ref{ouia}
(when $\mu>2$). In the case, $\mu=2$, we first apply 
Lemma 4.6 or Corollary 2 or Corollary 3 from 
\cite{alexakis4},\footnote{The fact that we are not 
dealing with the ``special subcase'' ensures that (\ref{garding9})
does not fall under a ``forbidden case'' of any 
of those Lemmas, by weight considerations.} 
to ensure that for the terms indexed in $P$ above, 
the (unique) factor $\nabla\upsilon$ contracts 
against a derivative index, and if it contracts against a 
factor $\nabla^{(B)}\Omega_x$ then $B\ge 3$. With this extra 
restriction, we repeat the argument above and derive our claim. $\Box$
\newline

{\it Proof of Lemma \ref{ouia}  in the ``special subcases'':} 
Let us write out:
$$\sum_{l\in L_\mu^*} a_l C^{l,i_1i_2}_g=\sum_{y\in Y} a_y C^y_g\cdot \nabla^{(2)}_{i_1i_2}\Omega_1.$$
(In other words, we ``factor out'' the term $\nabla^{(2)}_{i_1i_2}\Omega_1$ 
which contains the two free indices); we are 
then left with a complete contraction. 

Now, we apply the ``inverse integration by parts'' technique 
which was introduced in section 3 in \cite{alexakis5}, and then apply the 
``silly divergence formula'', obtained by integrating by parts 
with respect to the function $\Omega_1$.\footnote{Refer 
to section 3 in \cite{alexakis5} 
for  a detailed description of these operations.} 
Pick out the sublinear combination of terms with length $\sigma+u$, with 
{\it no} internal contractions and with $u$ factors $\nabla\phi_h$.
The resulting equation will be in the form: 

$$\Omega_1\cdot\sum_{y\in Y} a_y Spread^{\nabla^t,\nabla_t}\{Spread^{\nabla^s,\nabla_s}[C^t_g] \} =0.$$ 
(Recall that $Spread^{\nabla^t,\nabla_t}$ 
stands for a formal operation that acts on 
complete contractions in the form (\ref{form2}) by hitting two 
different factors by derivatives $\nabla^t,\nabla_t$ 
that contract against each other, and then adding over all the terms we can thus obtain).

Now, since the above holds formally, we derive that:

$$\sum_{y\in Y} a_y C^t_g =0.$$ 
Multiplying the above by $\nabla^{(2)}_{i_1i_2}\Omega_1\nabla^{i_1}\upsilon\nabla^{i_2}\upsilon$,
we derive our claim in the special subcase. $\Box$

\subsection{Proof of Lemmas \ref{ouic} and \ref{oui2a}.}

{\it Notation:} Firstly, we denote by $L^\nu\subset L_\mu$ the
index set of the $\mu$-tensor fields indexed in $L_\mu$, for which
the special factor $S_{*}\nabla^{(r)}R_{ijkl}$ has the index
${}_i$ contracting against a factor $\nabla^i\tilde{\phi}_\nu$ and
contains exactly one (non-special) free index.  We will also
denote by $\overline{L}^\nu$ the index set of the $(\mu+1)$-tensor
fields in (\ref{hypothese2}) which have two free indices in the
expression $S_{*}\nabla^{(\rho)}R_{ijkl}\nabla^i\tilde{\phi}_\nu$,
one of which is special--we assume wlog that the special free
index is ${}_k={}_{i_{\mu+1}}$.

\begin{lemma}
\label{thenewlemma} In the notation above, we claim that there
exists a linear combination of $(\mu+1)$-tensor fields, with a
$u$-simple character $\vec{\kappa}_{simp}$ and with certain
additional properties explained below (\ref{sobarwn}) so that:

\begin{equation}
\label{sobarwn}
\begin{split}
&\sum_{l\in L^\nu} a_l C^{l,i_1\dots
i_\mu}_g(\Omega_1,\dots,\Omega_p,\phi_1,\dots,\phi_u)
\nabla_{i_1}\upsilon\dots\nabla_{i_\mu}\upsilon+
\\& \sum_{l\in \overline{L}^\nu} a_l Xdiv_{i_{\mu+1}}C^{l,i_1\dots
i_{\mu+1}}_g(\Omega_1,\dots,\Omega_p,\phi_1,\dots,\phi_u)
\nabla_{i_1}\upsilon\dots\nabla_{i_\mu}\upsilon-
\\&Xdiv_{i_{\mu+1}}\sum_{p\in P} a_p C^{p,i_1\dots i_{\mu+1}}_g
(\Omega_1,\dots,\Omega_p,\phi_1,\dots,\phi_u)
\nabla_{i_1}\upsilon\dots\nabla_{i_\mu}\upsilon=
\\&\sum_{j\in J} a_j C^{j,i_1\dots i_\mu}_g
(\Omega_1,\dots,\Omega_p,\phi_1,\dots,\phi_u)
\nabla_{i_1}\upsilon\dots\nabla_{i_\mu}\upsilon.
\end{split}
\end{equation}
 The additional
properties of the tensor fields indexed in $P$ are as follows:
Firstly only the index ${}_{i_\mu}$ among the above free
indices ${}_{i_1},\dots,{}_{i_\mu}$ belongs to the (special) factor
$S_{*}\nabla^{(\rho)}R_{ijkl}\nabla^i\tilde{\phi}_\nu$, and
secondly {\it if} ${}_{i_{\mu+1}}$ does belong to the above factor
then $\rho>0$.
\end{lemma}

\par We observe that if we can prove the above then by making the
$\nabla\upsilon$s into $Xdiv$s (by virtue of the last Lemma 
in the Appendix of \cite{alexakis1}) we can derive both Lemma
\ref{ouic}, and Lemma \ref{oui2a} in case B.
\newline

{\it Proof of Lemma \ref{thenewlemma}:}

\begin{definition}
\label{skarta} We denote by $Cut(\vec{\kappa}_{simp})$ the
$(u-1)$-simple character that formally arises from
$\vec{\kappa}_{simp}$ by replacing the expression
$S_{*}\nabla^{(\rho)}R_{ijkl}\nabla^i\tilde{\phi}_\nu$ by a factor 
$\nabla^{(\rho+2)}Y$ ($Y$ is treated as a function
$\Omega_{p+1}$).

\par We then denote by $C^{l,i_1\dots i_\mu|A}_g$ the tensor fields
 that arises from $C^{l,i_1\dots i_\mu}_g$ by
replacing the expression $S_{*}\nabla^{(\rho)}_{r_1\dots
r_\rho}R_{ijkl}\nabla^i\tilde{\phi}_\nu$ by a factor
$\nabla^{(\rho+2)}_{r_1\dots r_\rho jk}Y\nabla_l\phi_\nu$. We also
denote by $C^{l,i_1\dots i_\mu|B}_g$ the tensor field that arises
from $C^{l,i_1\dots i_\mu}_g$ by formally replacing the expession
$S_{*}\nabla^{(\rho)}_{r_1\dots
r_\rho}R_{ijkl}\nabla^i\tilde{\phi}_\nu$ by
$-\nabla^{(\rho+2)}_{r_1\dots r_\rho jl}Y\nabla_k\phi_\nu$.

\par Analogously, for each $l\in \overline{L}^\nu$ we denote by
$C^{l,i_1\dots i_{\mu+1}|A}_g$  the tensor field that arises from
by by replacing the expression $S_{*}\nabla^{(\rho)}_{r_1\dots
r_\rho}R_{ijkl}\nabla^i\tilde{\phi}_\nu$ by
$\nabla^{(\rho+2)}_{r_1\dots r_\rho jk}Y\nabla_l\phi_\nu$ (${}_l$
is not a free index). We also denote by $C^{l,i_1\dots i_\mu|B}_g$
the tensor field that arises from $C^{l,i_1\dots i_\mu}_g$ by
replacing the factor $S_{*}\nabla^{(\rho)}_{r_1\dots
r_\rho}R_{ijkl}\nabla^i\tilde{\phi}_\nu$ by
$-\nabla^{(\rho+2)}_{r_1\dots r_\rho jl}Y\nabla_k\phi_\nu$ (now
${}_k$ is the free index ${}_{i_{\mu+1}}$).
\end{definition}

A note is in order: When we refer to the tensor field
$C^{l,i_1\dots i_{\mu+1}|A}_g$ below and we write
$\tilde{Xdiv}_{i_{\mu+1}}C^{l,i_1\dots i_{\mu+1}|A}_g$,
$\tilde{Xdiv}_{i_{\mu+1}}$  will stand for the sublinear
combination in $Xdiv_{i_{\mu+1}}$ where $\nabla^{i_{\mu+1}}$ is
{\it not} allowed to hit the factor $\nabla^{(B)}Y$. Furthermore,
when we write $\tilde{Xdiv}_{i_{\mu+1}}C^{l,i_1\dots
i_{\mu+1}|B}_g$ below, $\tilde{Xdiv}_{i_{\mu+1}}$ will stand for
the regular $Xdiv_{i_{\mu+1}}$ but we will ``forget'' this
structure--i.e. we will treat as a sum of $\mu$-tensor fields.

\par We will now denote by

$$\Sum_{u\in U} a_u C^{u,i_1\dots i_a}_{g}(\Omega_1,\dots ,\Omega_p,Y,\phi_1,\dots
,\phi_u)$$ a generic linear combination of $a$-tensor fields
($a\ge \mu+1$) with length $\sigma +u$, with the factor
$\nabla\phi_\nu$ {\it not} contracting against the factor
$\nabla^{(A)}Y$ and {\it not} containing a free index.

\par Now, considering the sublinear combination in
$Image^1_Y[L_g](=0)$ which consists of terms where the factor
$S_{*}\nabla^{(\rho)}R_{ijkl}\nabla^i\tilde{\phi}_\nu$, is
replaced by $\nabla^{(\rho+2)}Y$ and $\nabla\phi_\nu$ is {\it not}
contracting against $\nabla^{(\rho+2)}Y$, we derive a new true
equation:

\begin{equation}
\label{yusupov}
\begin{split}
&\sum_{l\in L^\nu\bigcup (L_\mu\setminus L^\nu)} a_l
Xdiv_{i_1}\dots Xdiv_{i_\mu}\{ C^{l,i_1\dots
i_\mu|A}_g(\Omega_1,\dots ,\Omega_p,Y,\phi_1,\dots ,\phi_u)
\\&+C^{l,i_1\dots i_\mu|B}_g(\Omega_1,\dots ,\Omega_p,Y,\phi_1,\dots
,\phi_u) \}+
\\&\Sum_{l\in \overline{L}^\nu} a_l Xdiv_{i_1}\dots Xdiv_{i_\mu}\tilde{Xdiv_{i_{\mu+1}}}
\{ C^{l,i_1\dots i_{\mu+1}|A}_{g}(\Omega_1,\dots
,\Omega_p,Y,\phi_1,\dots ,\phi_u)
\\&+C^{l,i_1\dots
i_{\mu+1}|B}_{g}(\Omega_1,\dots ,\Omega_p,Y,\phi_1,\dots
,\phi_u)\}+
\\&\Sum_{u\in U} a_u Xdiv_{i_1}\dots Xdiv_{i_a}C^{u,i_1\dots i_a}_{g}(\Omega_1,\dots
,\Omega_p,Y,\phi_1,\dots ,\phi_u)+
\\&\sum_{j\in J} a_j C^j_g(\Omega_1,\dots
,\Omega_p,Y,\phi_1,\dots ,\phi_u)=0;
\end{split}
\end{equation}
here the terms indexed in $J$ are simply subsequent to the simple
character $Cut(\vec{\kappa}_{simp})$.

\par We then apply the inductive assumption of Lemma
4.10 in \cite{alexakis4}\footnote{Notice that some tensor fields of minimum
rank $\mu$ in (\ref{yusupov}), i.e. the ones indexed in $L^\nu$, 
will have only non-special free indices,
therefore there is no danger of falling under a ``forbidden case''
of that Lemma.} to the above, and pick out the sublinear
combination of terms where one factor $\nabla\upsilon$ is
contracting against the factor $\nabla^{(B)}Y$ and the other
$\mu-1$ factors $\nabla\upsilon$ are contracting against other
 factors. This sublinear combination must
vanish separately, thus we derive a new equation:

\begin{equation}
\label{yusupovcor}
\begin{split}
&\sum_{l\in L^\nu} a_l \{C^{l,i_1\dots i_\mu|A}_g(\Omega_1,\dots
,\Omega_p,Y,\phi_1,\dots ,\phi_u)
\\&+C^{l,i_1\dots
i_\mu|B}_g(\Omega_1,\dots ,\Omega_p,Y,\phi_1,\dots
,\phi_u)\}\nabla_{i_1}\upsilon\dots\nabla_{i_\mu}\upsilon+
\\&\Sum_{l\in \overline{L}^\nu_{\mu+1}} a_l \tilde{Xdiv_{i_{\mu+1}}}
\{ C^{l,i_1\dots i_{\mu+1}|A}_{g}(\Omega_1,\dots
,\Omega_p,Y,\phi_1,\dots ,\phi_u)
\\&+C^{l,i_1\dots
i_{\mu+1}|B}_{g}(\Omega_1,\dots ,\Omega_p,Y,\phi_1,\dots
,\phi_u)\}\nabla_{i_1}\upsilon\dots\nabla_{i_\mu}\upsilon+
\\&\Sum_{u\in U} a_u Xdiv_{i_{\mu+1}}C^{u,i_1\dots i_{\mu+1}}_{g}(\Omega_1,\dots
,\Omega_p,Y,\phi_1,\dots
,\phi_u)\nabla_{i_1}\upsilon\dots\nabla_{i_\mu}\upsilon+
\\&\sum_{j\in J} a_j C^j_g(\Omega_1,\dots
,\Omega_p,Y,\phi_1,\dots ,\phi_u,\upsilon^\mu)=0.
\end{split}
\end{equation}
(Here the tensor fields indexed in $U$ are generic acceptable
 $(\mu+1)$-tensor fields with a $u$-simple character
 $Cut(\vec{\kappa}_{simp})$--the free index ${}_{i_{\mu+1}}$
 {\it does not} belong to the factor $\nabla\phi_\nu$).

\par Now, we define an operation $Op[\dots]$ which acts on the
tensor fields above by replacing the expression
$\nabla^{(B)}_{r_1\dots r_B}Y\nabla^{r_B}\upsilon\nabla_a\phi_\nu$
by an expression $\nabla^{(B+1)}_{r_1\dots
r_{B-1}ar_B}\phi_\nu\nabla^{r_B}\upsilon$ (denote the
$(u-1)$-simple character that we thus construct by
$Cut'(\vec{\kappa}_{simp})$--the factor $\nabla^{(A)}\phi_{u+1}$
is treated as a factor $\nabla^{(A)}\Omega_{p+1}$). Since
(\ref{yusupovcor}) holds formally, we derive:

\begin{equation}
\label{yusupovcor'}
\begin{split}
&\sum_{l\in L^\nu} a_l \{Op[C]^{l,i_1\dots
i_\mu|A}_g(\Omega_1,\dots ,\Omega_p,Y,\phi_1,\dots ,\phi_u)
\\&+Op[C]^{l,i_1\dots i_\mu|B}_g(\Omega_1,\dots
,\Omega_p,Y,\phi_1,\dots ,\phi_u)
\}\nabla_{i_1}\upsilon\dots\nabla_{i_\mu}\upsilon+
\\&\Sum_{l\in \overline{L}^\nu} a_l \tilde{Xdiv_{i_{\mu+1}}}
\{ Op[C]^{l,i_1\dots i_{\mu+1}|A}_{g}(\Omega_1,\dots
,\Omega_p,Y,\phi_1,\dots ,\phi_u)
\\&+Op[C]^{l,i_1\dots
i_{\mu+1}|B}_{g}(\Omega_1,\dots ,\Omega_p,Y,\phi_1,\dots
,\phi_u)\}\nabla_{i_1}\upsilon\dots\nabla_{i_\mu}\upsilon+
\\&\Sum_{u\in U'} a_u Xdiv_{i_{\mu+1}}Op[C]^{u,i_1\dots i_{\mu+1}}_{g}(\Omega_1,\dots
,\Omega_p,Y,\phi_1,\dots
,\phi_u)\nabla_{i_1}\upsilon\dots\nabla_{i_\mu}\upsilon+
\\&\sum_{j\in J} a_j Op[C]^j_g(\Omega_1,\dots
,\Omega_p,Y,\phi_1,\dots ,\phi_u,\upsilon^\mu)=
\\&\sum_{z\in Z} a_z C^z_g(\Omega_1,\dots
,\Omega_p,Y,\phi_1,\dots ,\phi_u,\upsilon^\mu);
\end{split}
\end{equation}
here the tensor fields indexed in $Z$ on the RHS have length
$\sigma+u+\mu$ (as opposed to all the terms in the LHS which have
length $\sigma+u-1+\mu$), and in addition have a factor
$\nabla^{(A)}\phi_{u+1}$ with $A\ge 2$. The correction terms arise by 
repeating the formal permutations by which the LHS is 
made formally zero by the LHS of (\ref{yusupovcor'}). 
The claim $A\ge 2$ follows because the rightmost two 
indices in each factor will {\it not} be permuted.

\par Now, we observe that for each $l\in L^\nu$:

\begin{equation}
\label{reesoula}
\begin{split}
& \{Op[C]^{l,i_1\dots i_\mu|A}_g(\Omega_1,\dots
,\Omega_p,Y,\phi_1,\dots
,\phi_u)+Op[C]^{l,i_1\dotsi_\mu|B}_g(\Omega_1,\dots
,\Omega_p,\\&Y,\phi_1,\dots
,\phi_u)\}\nabla_{i_1}\upsilon\dots\nabla_{i_\mu}\upsilon
=C^{l,i_1\dots i_\mu}_g(\Omega_1,\dots
,\Omega_p,\phi_1,\dots
,\phi_u)\nabla_{i_1}\upsilon\dots\nabla_{i_\mu}\upsilon+
\\&\sum_{z\in Z} a_z C^z_g(\Omega_1,\dots
,\Omega_p,Y,\phi_1,\dots ,\phi_u,\upsilon^\mu),
\end{split}
\end{equation}
(where the terms indexed in $Z$ are generic complete contractions
as defined above). Analogously, we derive that for each $l\in
\overline{L}^\nu$:

\begin{equation}
\label{reesoula}
\begin{split}
& \tilde{Xdiv}_{i_{\mu+1}}\{ Op[C]^{l,i_1\dots
i_{\mu+1}|A}_{g}(\Omega_1,\dots ,\Omega_p,Y,\phi_1,\dots
,\phi_u)\nabla_{i_1}\dots\nabla_{i_\mu}\upsilon+
\\&Op[C]^{l,i_1\dots
i_{\mu+1}|B}_{g}(\Omega_1,\dots ,\Omega_p,Y,\phi_1,\dots
,\phi_u)\}\nabla_{i_1}\dots\nabla_{i_\mu}\upsilon=
\\&Xdiv_{i_{\mu+1}}C^{l,i_1\dots i_{\mu+1}}_g(\Omega_1,\dots
,\Omega_p,\phi_1,\dots
,\phi_u)\nabla_{i_1}\upsilon\dots\nabla_{i_\mu}\upsilon+
\\&\sum_{z\in Z} a_z C^z_g(\Omega_1,\dots
,\Omega_p,Y,\phi_1,\dots ,\phi_u,\upsilon^\mu).
\end{split}
\end{equation}

\par We then substitute  the above two equations into
(\ref{yusupovcor}) and we obtain a new equation:

\begin{equation}
\label{beyogl}
\begin{split}
&\Sum_{u\in U} a_u Xdiv_{i_{\mu+1}}Op[C]^{u,i_1\dots
i_{\mu+1}}_{g}(\Omega_1,\dots ,\Omega_p,Y,\phi_1,\dots
,\phi_u)\nabla_{i_1}\upsilon\dots\nabla_{i_\mu}\upsilon+
\\&\sum_{j\in J} a_j Op[C]^j_g(\Omega_1,\dots
,\Omega_p,Y,\phi_1,\dots ,\phi_u,\upsilon^\mu)+
\\&\sum_{l\in L^\nu} a_l C^{l,i_1\dots i_\mu}_g(\Omega_1,\dots
,\Omega_p,\phi_1,\dots
,\phi_u)\nabla_{i_1}\upsilon\dots\nabla_{i_\mu}\upsilon+
\\&\Sum_{l\in \overline{L}^\nu_{\mu+1}} a_l Xdiv_{i_{\mu+1}}
C^{l,i_1\dots i_{\mu+1}}_g(\Omega_1,\dots ,\Omega_p,\phi_1,\dots
,\phi_u)\nabla_{i_1}\upsilon\dots\nabla_{i_\mu}\upsilon =
\\&\sum_{z\in Z} a_z C^z_g(\Omega_1,\dots
,\Omega_p,Y,\phi_1,\dots ,\phi_u,\upsilon^\mu),
\end{split}
\end{equation}
(denote the $(u-1+\mu)$-simple character o the tensor fields in
the LHS of the above by $Ext[\vec{\kappa}_{simp}]$ (the factors
$\nabla\upsilon$ are now treated as factors $\nabla\phi_h$).

\par We now derive our claim from (\ref{beyogl})
via an induction: We inductively assume an equation:

\begin{equation}
\label{beyogl'}
\begin{split}
&\Sum_{u\in U^\delta} a_u Xdiv_{i_{\mu+1}}\dots
Xdiv_{i_{\mu+\delta}}\tilde{C}^{u,i_1\dots
i_{\mu+1}}_{g}(\Omega_1,\dots
,\Omega_p,\phi_{\nu},\phi_1,\dots,\hat{\phi}_\nu,\dots
,\phi_u)\nabla_{i_1}\upsilon\\&\dots\nabla_{i_\mu}\upsilon
+\sum_{j\in J} a_j \tilde{C}^j_g(\Omega_1,\dots
,\Omega_p,\phi_\nu,\phi_1,\dots ,\hat{\phi}_\nu,\dots
,\phi_u,\upsilon^\mu)+
\\&\sum_{l\in L^\nu} a_l C^{l,i_1\dots i_\mu}_g(\Omega_1,\dots
,\Omega_p,\phi_1,\dots
,\phi_u)\nabla_{i_1}\upsilon\dots\nabla_{i_\mu}\upsilon+
\\&\Sum_{l\in \overline{L}^\nu} a_u Xdiv_{i_{\mu+1}}
C^{l,i_1\dots i_{\mu+1}}_{g}(\Omega_1,\dots ,\Omega_p,\phi_1,\dots
,\phi_u)\nabla_{i_1}\upsilon\dots\nabla_{i_\mu}\upsilon+
\\&\sum_{p\in P} a_p Xdiv_{i_{\mu+1}}C^{p,i_1\dots i_{\mu+1}}_g
(\Omega_1,\dots ,\Omega_p,\phi_1,\dots
,\phi_u)\nabla_{i_1}\upsilon\dots\nabla_{i_\mu}\upsilon=
\\&\sum_{z\in Z} a_z C^z_g(\Omega_1,\dots
,\Omega_p,Y,\phi_1,\dots ,\phi_u,\upsilon^\mu).
\end{split}
\end{equation}
Here the tensor fields indexed in $U^\delta$ are like the ones
indexed in $U$ in (\ref{beyogl}) (in particular they have a factor
$\nabla^{(A)}_{y_1\dots y_A}Y\nabla^{y_A}\upsilon$ with $A\ge 3$)
but in addition have rank $\delta>0$ (thus (\ref{beyogl}) is a
special case of (\ref{beyogl'}) with $\delta=1$). Furthermore, the tensor fields
indexed in $P$ are as described in the claim of Lemma \ref{thenewlemma}.

Using the generic notation introduced above, we will then show
(using generic notation in the first line below) that we can
write:

\begin{equation}
\label{beyogl''}
\begin{split}
&\Sum_{u\in U^{\delta+1}} a_u Xdiv_{i_{\mu+1}}\dots
Xdiv_{i_{\mu+\delta+1}}\tilde{C}^{u,i_1\dots
i_{\mu+\delta+1}}_{g}(\Omega_1,\dots
,\Omega_p,\phi_{\nu},\phi_1,\dots,\hat{\phi}_\nu,\dots
,\phi_u)\\&\nabla_{i_1}\upsilon\dots\nabla_{i_\mu}\upsilon
+\sum_{j\in J} a_j \tilde{C}^j_g(\Omega_1,\dots
,\Omega_p,\phi_\nu,\phi_1,\dots ,\hat{\phi}_\nu,\dots
,\phi_u,\upsilon^\mu)+
\\&\sum_{l\in L^\nu} a_l C^{l,i_1\dots i_\mu}_g(\Omega_1,\dots
,\Omega_p,\phi_1,\dots
,\phi_u)\nabla_{i_1}\upsilon\dots\nabla_{i_\mu}\upsilon+
\\&\Sum_{l\in \overline{L}^\nu} a_u \tilde{Xdiv}_{i_{\mu+1}}
C^{l,i_1\dots i_{\mu+1}}_{g}(\Omega_1,\dots
,\Omega_p,\phi_1,\dots
,\phi_u)\nabla_{i_1}\upsilon\dots\nabla_{i_\mu}\upsilon+
\\&\sum_{p\in P} a_p Xdiv_{i_{\mu+1}}C^{p,i_1\dots i_{\mu+1}}_g
(\Omega_1,\dots ,\Omega_p,\phi_1,\dots
,\phi_u)\nabla_{i_1}\upsilon\dots\nabla_{i_\mu}\upsilon=
\\&\sum_{z\in Z} a_z C^z_g(\Omega_1,\dots
,\Omega_p,Y,\phi_1,\dots ,\phi_u,\upsilon^\mu).
\end{split}
\end{equation}
If we can prove the above, then by iterative repetition we derive
our claim.
\newline

{\it Proof of (\ref{beyogl''}):} We treat the factors
$\nabla\upsilon$ as factors $\nabla\phi_h,h>u$ (this can be done
easily by a simple polarization). We then notice that
(\ref{beyogl'}) immediately implies:

\begin{equation}
\label{thskw} \begin{split}
&\Sum_{u\in U^\delta} a_u
Xdiv_{i_{\mu+1}}\dots Xdiv_{i_{\mu+\delta}}\tilde{C}^{u,i_1\dots
i_{\mu+1}}_{g}(\Omega_1,\dots
,\Omega_p,\phi_{\nu},\phi_1,\dots,\hat{\phi}_\nu,\dots
,\phi_u)\\&\nabla_{i_1}\upsilon\dots\nabla_{i_\mu}\upsilon
+\sum_{j\in J} a_j \tilde{C}^j_g(\Omega_1,\dots
,\Omega_p,\phi_\nu,\phi_1,\dots ,\hat{\phi}_\nu,\dots
,\phi_u,\upsilon^\mu)=0,
\end{split}
\end{equation}
 modulo complete contractions of length
$\ge\sigma+u+\mu$.

\par Therefore, we  apply our inductive assumption of Corollary
1 in \cite{alexakis4} to 
the above, or if the above falls under a forbidden
case of Corollary 1, we then apply the ``weak substitute''
of that Corollary from the Appenix in \cite{alexakis4}. 
(Notice that if the terms in $U^\delta$ contain ``forbidden tensor fields'' for Corollary 1, 
then necessarily by construction $\delta>1$). We derive that there exists a
linear combination of acceptable tensor fields with a
$(u-1+\mu)$-simple character $Ext[\vec{\kappa}_{simp}]$ and with
rank $\delta+1$ (indexed in $U^{\delta+1}$ below) so that modulo
complete contractions of length $\sigma+u+\mu+\delta$:

\begin{equation}
\label{barely}
\begin{split}
&\Sum_{u\in U^\delta} a_u \tilde{C}^{u,i_1\dots
i_{\mu+\delta}}_{g}(\Omega_1,\dots
,\Omega_p,\phi_{\nu},\phi_1,\dots,\hat{\phi}_\nu,\dots
,\phi_u)\nabla_{i_1}\upsilon\dots\nabla_{i_\mu}\upsilon\\&\nabla_{i_{\mu+1}}\omega\dots
\nabla_{i_{\mu+\delta}}\omega-
Xdiv_{i_{\mu+\delta+1}} \Sum_{u\in U^{\delta+1}} a_u \tilde{C}^{u,i_1\dots
i_{\mu+\delta+1}}_{g}(\Omega_1,\dots
,\Omega_p,\phi_{\nu},\phi_1,\dots,\hat{\phi}_\nu,\dots
,\phi_u)\\&\nabla_{i_1}\upsilon\dots\nabla_{i_\mu}\upsilon\nabla_{i_{\mu+1}}\omega\dots
\nabla_{i_{\mu+\delta}}\omega
=\sum_{j\in J} a_j C^{u,i_1\dots
i_{\mu+\delta}}_{g}(\Omega_1,\dots
,\Omega_p,\phi_{\nu},\phi_1,\dots,\\&\hat{\phi}_\nu,\dots
,\phi_u)\nabla_{i_1}\upsilon\dots\nabla_{i_\mu}\upsilon\nabla_{i_{\mu+1}}\omega\dots
\nabla_{i_{\mu+\delta}}\omega.
\end{split}
\end{equation}

\par Now, since the above holds formally at the linearized level,
it follows that the correction terms of length will be in the form
$$\sum_{p\in P} a_p C^{u,i_1\dots
i_{\mu+\delta}}_{g}\nabla_{i_1}\upsilon\dots\nabla_{i_\mu}\upsilon\nabla_{i_{\mu+1}}\omega\dots
\nabla_{i_{\mu+\delta}}\omega+\sum_{j\in J}\dots +\sum_{z\in
Z}\dots$$ Then, making the factors $\nabla\omega$ into $Xdiv$'s, 
by virtue of the last Lemma in 
the Appendix of \cite{alexakis1}, we
derive our claim. $\Box$

\end{document}